\documentclass[10pt, a4paper]{article}
\usepackage[T1]{fontenc}
\usepackage{a4wide, url}
\usepackage[english]{babel}
\usepackage{graphicx}
\usepackage{array}
\usepackage{natbib}
\usepackage[textfont=sc,labelfont=bf]{caption}
\usepackage{setspace}
\usepackage{multirow}
\usepackage{fancyhdr}
\usepackage{amsthm}
\usepackage{amsmath}
\usepackage{amsfonts}
\usepackage{amssymb}
\usepackage{bigints}
\usepackage{pgf} 
\usepackage{bm}
\usepackage{paralist}
\usepackage{hyperref}
\usepackage{pdflscape}
\usepackage{xfrac}
\usepackage[linesnumbered,ruled,vlined]{algorithm2e}
\SetKwInput{KwInput}{Input}                
\SetKwInput{KwOutput}{Output}              
\usepackage{enumitem,kantlipsum}
\usepackage{xr}
\usepackage{yfonts}
\usepackage{algorithmic}
\usepackage{xcolor}
\hypersetup{
    colorlinks,
    linkcolor={red!50!black},
    citecolor={blue!50!black},
    urlcolor={blue!80!black}
}

\usepackage{comment}
\excludecomment{ext}\excludecomment{comment}

\usepackage[top=1in,bottom=1in,left=.75in,right=.75in]{geometry}

\usepackage{tabu}

\newcolumntype{L}[1]{>{\raggedright\let\newline\\\arraybackslash\hspace{0pt}}m{#1}}
\newcolumntype{C}[1]{>{\centering\let\newline\\\arraybackslash\hspace{0pt}}m{#1}}
\newcolumntype{R}[1]{>{\raggedleft\let\newline\\\arraybackslash\hspace{0pt}}m{#1}}

\DeclareMathAlphabet\mathbfcal{OMS}{cmsy}{b}{n}

\newcommand{\mbf}{\mathbf}
\newcommand{\beq}{\begin{equation}}
\newcommand{\eeq}{\end{equation}}
\newcommand{\bea}{\begin{eqnarray}}
\newcommand{\eea}{\end{eqnarray}}
\newcommand{\ba}{\begin{array}}
\newcommand{\ea}{\end{array}}
\newcommand{\al}[1]{\begin{align}#1\end{align}}
\newcommand{\bit}{\begin{itemize}}
\newcommand{\eit}{\end{itemize}}
\newcommand{\ben}{\begin{enumerate}} 
\newcommand{\een}{\end{enumerate}}
\newcommand{\bpm}{\begin{pmatrix}}
\newcommand{\epm}{\end{pmatrix}}
\newcommand{\bbm}{\begin{bmatrix}}
\newcommand{\ebm}{\end{bmatrix}}

\renewcommand{\l}{\left}
\renewcommand{\r}{\right}

\newcommand{\E}[0]{\mathbb{E}}
\newcommand{\Var}[0]{\mathbb{V}\mathrm{ar}}
\newcommand{\Cov}[0]{\mathbb{C}\mathrm{ov}}

\newcommand{\ra}{\rightarrow}
\newcommand{\nn}{\nonumber}
\newcommand{\wh}{\widehat}
\newcommand{\wt}{\widetilde}

\newcommand\bianco[1]{\textcolor{white}{#1}}

\newtheorem{ass}{Assumption}

\newtheorem{prop}{Proposition}
\newtheorem{rem}{Remark}
\newtheorem{lem}{Lemma}
\newtheorem{cor}{Corollary}

\graphicspath{{Figures/Submission/}}
\title{\textsc{\large 	
Quasi Maximum Likelihood Estimation and Inference\\ \vskip .2cm of Large Approximate Dynamic Factor Models\\ \vskip .2cm
 via the EM algorithm\\ \vskip -.2cm}}
\date{ }

\begin{document}
\maketitle

\begin{center}\vspace{-1.5cm}
\begin{tabular}{cp{1cm}c}
 Matteo Barigozzi &&  Matteo Luciani\\
\footnotesize Universit\`a di Bologna &&\footnotesize Federal Reserve Board\\[-.2cm]
\scriptsize  matteo.barigozzi@unibo.it &&\scriptsize matteo.luciani@frb.gov\\[.5cm]
\end{tabular}

\small This version: \today\\
\footnotesize{First draft: November 8, 2017$^*$}\\
\end{center}

\begin{abstract}
We study estimation of large Dynamic Factor models implemented through the Expectation Maximization (EM) algorithm, jointly with the Kalman smoother.~We prove that as both the cross-sectional dimension, $n$, and the sample size, $T$, diverge to infinity: (i) the estimated loadings are $\sqrt T$-consistent, asymptotically normal and equivalent to their Quasi Maximum Likelihood estimates; (ii) the estimated factors are $\sqrt n$-consistent, asymptotically normal and equivalent to their Weighted Least Squares estimates. Moreover, the estimated loadings are asymptotically as efficient as those obtained by Principal Components analysis, while the estimated factors are more efficient if the idiosyncratic covariance is sparse enough.
\vspace{0.5cm}

\noindent \textit{Keywords:} 
Approximate Dynamic Factor Model; Expectation Maximization Algorithm; Kalman Smoother; Quasi Maximum Likelihood.
\end{abstract}

\renewcommand{\thefootnote}{$\ast$} 
\thispagestyle{empty}

\footnotetext{A previous version of some results in this paper appeared in the paper ``Common factors, trends, and cycles in large datasets'',  available at \href{https://arxiv.org/abs/1709.01445}{arXiv:1709.01445} or at \href{ https://doi.org/10.17016/FEDS.2017.111}{ FEDS.2017.111}.\\

\noindent
M. Barigozzi gratefully acknowledges financial support from MIUR (PRIN2020, Grant 2020N9YFFE).\\[.03cm]

\noindent
Disclaimer: the views expressed in this paper are those of the authors and do not necessarily reflect the views and policies of the Board of Governors or the Federal Reserve System.\\[.03cm]

\noindent
We thank for helpful comments: Majid Al-Sadoon, Leopoldo Catania, Giuseppe Cavaliere, Manfred Deistler, Massimo Franchi, Ivan Petrella, Esther Ruiz, Haihan Tang, Lorenzo Trapani, and Luca Trapin. } 

\renewcommand{\thefootnote}{\arabic{footnote}}

%
%
\section{Introduction}

Factor analysis can be considered a pioneering technique in unsupervised statistical learning \citep{GZ96}. It originally gained popularity in the early decades of the twentieth century as a dimension-reduction technique used in psychometrics \citep{spearman04}. Since then, it has become a classical method used for the statistical analysis of complex datasets in many human, natural, and social sciences (see, e.g., \citealp[Chapter 1, and references therein]{lawleymaxwell71}). In the last thirty years, factor analysis has seen significant success in financial and macroeconometrics because it allows to analyze and predict economic activity by summarizing large panels of economic time series in a simple and effective way (see, e.g., the survey by \citealp{stockwatson16} and references therein).

An $r$-factor model is defined by 
\beq\label{eq:Fsimple}
x_{it} =\mu_i+ \bm\lambda_i^\prime \mbf F_t + \xi_{it}, \quad i=1,\ldots, n.,\quad t=1,\ldots,T,
\eeq
where $x_{it}$ is the observation for the $i$th cross-section at time $t$, $\mu_i$ is a constant, and $\mbf F_t$ and $\bm\lambda_i$ are $r$-dimensional latent column vectors of \textit{factors} and factor \textit{loadings}, with $r\ll n$. We call $\bm\lambda_i^\prime \mbf F_t$ the \textit{common} component and $\xi_{it}$ the \textit{idiosyncratic} component. Throughout, we consider the standard case in which all $\{x_{it}\}$ are zero-mean weakly stationary processes or are the result of a transformation to stationarity.

Furthermore, in the case of time series, the factors are likely to be autocorrelated. For example, we can assume simple first order autoregressive dynamics: 
\beq\label{eq:VARsimple}
\mbf F_{t}= \mbf A\mbf F_{t-1}+\mbf v_t, \quad t=1,\ldots ,T,
\eeq
with $\mbf v_{t}$ being an $r$-dimensional vector of innovations. Likewise, the idiosyncratic components might be autocorrelated. The measurement equation \eqref{eq:Fsimple} and the state equation \eqref{eq:VARsimple} form a state-space model, or, equivalently, a Dynamic Factor Model (DFM) (this is a restricted version of the more general model by \citealp{FHLR00}, where factors can be loaded also with lags). Thanks to its simplicity and empirical success, the DFM is the most common approach to factor analysis of high-dimensional time series. 

In large dimensional macroeconomic and financial datasets, the idiosyncratic components are likely to be also cross-correlated. Indeed, although macroeconomic or financial market dynamics are the main drivers of the comovement in these datasets, sectoral and local comovements are non-negligible sources of fluctuations. 
In the case of correlated idiosyncratic components the factor model is called \textit{approximate} as opposed to an \textit{exact} factor model having uncorrelated idiosyncratic components. 

In an exact factor model a small number of variables is enough to estimate the loadings by Quasi Maximum Likelihood (QML), but we cannot consistently estimate the factors \citep{lawleymaxwell71}. 
In an approximate factor model, we can disentangle the common and idiosyncratic components only in the extreme case when $n\to\infty$ \citep{chamberlainrothschild83}---in other words, we do not suffer the typical ``curse of dimensionality'' but rather benefit from the ``blessing of dimensionality''. 
In particular, when $n\to\infty$, we can consistently estimate the factors by some form of linear projection onto the estimated loadings. However, QML estimation of the loadings is now unfeasible since, in principle, it requires to also jointly estimate all the $T$ idiosyncratic (auto)covariance matrices, each of size $n\times n$,
as well as the $T$ (auto)covariance matrices of the factors.  Thus, we need to explore alternative approaches. 

There are three main solutions to this problem. The first is Principal Component (PC) analysis, which delivers the optimal non-parametric estimator of a large approximate factor model, see, e.g., \citet{stockwatson02JASA} and \citet{Bai03}. The second is QML estimation based on a mis-specified exact model with no autocorrelations in the idiosyncratic components and the factors, and a diagonal or sparse idiosyncratic covariance matrix, see, e.g., \citet{baili16} and \citet{bailiao16}. These two solutions only consider equation \eqref{eq:Fsimple}, thus effectively estimating a static factor model, not a DFM.
There is a third solution, which is the focus of this paper, that considers joint estimation of approximate DFM defined in \eqref{eq:Fsimple}-\eqref{eq:VARsimple}: the Expectation Maximization (EM) algorithm (\citealp{quahsargent93}, \citealp{DGRqml}).

The EM algorithm is fully parametric and based on an iteration of two steps: (E-step) given the loadings and all the model's parameters, the factors and their second moments are estimated via the Kalman smoother and are used to compute the expected log-likelihood conditional on the observed data; (M-step) the expected log-likelihood is maximized to obtain a new estimate of the loadings and all the model's parameters. These steps are always implemented by considering a mis-specified likelihood with uncorrelated  idiosyncratic components, thus making estimation feasible and producing estimators with closed-form expression 
As such, the EM algorithm should be regarded as an approximation of the QML estimation method because it maximizes a mis-specified likelihood of an exact DFM using an iterative procedure. 

This paper focuses on the theoretical properties of the EM and Kalman smoother estimators. In a couple of breakthrough studies, \citet{DGRfilter,DGRqml} provided the first fundamental theoretical treatment of these estimators, which quickly became popular in empirical macroeconomic research. Indeed, this approach allows the user to easily deal with data irregularities and missing values and impose restrictions that reflect any prior knowledge about the data on the model.\footnote{PC analysis in presence of missing data and parameter constraints has been studied by, e.g., \citet{baing19}, \citet{fan2022we}, \citet{XP20}, among others.} 

We make three main contributions. First, we prove that, as $n,T\to\infty$, the estimator of the loadings obtained via the EM algorithm converges asymptotically to a unique maximum of the likelihood. It is well known that in the Gaussian quasi-likelihood case, the EM algorithm produces approximate QML estimators \citep{wu83,BWB17}. In this paper, we refine this result by showing that, as $n,T\to\infty$, the approximation error not only depends on the number of iterations but also becomes negligible. Moreover, the EM estimator of the loadings is asymptotically equivalent to the unfeasible Ordinary Least Squares estimator we would obtain if we had observed the factors. This result holds for any reasonable, but not necessarily consistent, pre-estimator of the loadings used to initialize the EM algorithm. We derive similar results for all the estimated parameters, namely, the idiosyncratic variances, the VAR coefficients, and the covariance matrix of the VAR residuals in equation \eqref{eq:VARsimple}. 

Second, we prove that, as $n,T\to\infty$, the estimator of the factors obtained via the Kalman smoother computed using the parameters estimated via the EM algorithm is equivalent to the unfeasible Weighted Least Squares estimator we would obtain if we had observed the loadings and idiosyncratic variances. As a by-product of this result, we also show that the Kalman smoother and filter estimators are asymptotically equivalent.

Third, we show that the EM estimator of the loadings is asymptotically equivalent to the PC estimator; hence, it has the same consistency rate ($\min(n,\sqrt T)$), is asymptotically normal and equally efficient. Likewise, the Kalman smoother has the same consistency rate ($\min(\sqrt n,T)$) as the PC estimator, it is asymptotically normal, and if the idiosyncratic covariance matrix is sufficiently sparse, it is more efficient than the PC estimator.

This paper is the first to fully characterize the asymptotic properties of the EM and Kalman smoother estimators. Other papers provide results that are close to ours, using more restrictive approaches or deriving only partial asymptotic results. \citet{baili16} considered QML estimation of the loadings for the static factor model \eqref{eq:Fsimple} only and did not study the convergence of the employed maximization algorithm. \citet{DGRfilter} considered the Kalman smoother obtained using the PC estimator of the loadings but did not derive its asymptotic distribution and obtained a slower consistency rate.  Last,  \citet{DGRqml} proved consistency of the Kalman smoother obtained from the EM algorithm but derived a slower rate and did not prove its asymptotic normality, nor did they prove consistency of the EM estimator of the loadings.

Our results lay the theoretical foundations for the wide empirical success of the EM algorithm for estimating large dimensional DFMs (see the next section for a list of applications) and answer two long-standing critiques. First, by showing the equivalence of the consistency rates of the EM and PC estimators, we reverse the belief that PC is a superior approach. 
Second, by providing the asymptotic distributions, we answer the call by \citet{TW87} and \citet{geweke93}, who advocated a Bayesian approach based on the Gibbs sampler because the EM algorithm provides only point estimates.

The paper is organized as follows. In Section \ref{sec:lit}, we briefly review the main applications of the EM algorithm and KS in factor analysis and alternative methods proposed to estimate DFM. In Section \ref{sec:est_nutshell}, we describe the estimation and give a guide for implementing it. All assumptions are in Section \ref{sec:sdfm}. 
The asymptotic results are in Section \ref{sec:asymp}. In Section \ref{sec:PCeff}, we discuss efficiency of the EM estimator and the KS and compare them with the PC estimator. In Section \ref{sec:test}, we propose estimators of the asymptotic covariance matrices.  In Section \ref{sec:mc}, we present an extensive MonteCarlo study, and in Section \ref{sec:emp}, we apply the EM algorithm to US macroeconomic data. Section \ref{sec:conclusions} concludes. The proofs of all theoretical results are in the Appendix.

%
%
\section{Related literature}\label{sec:lit}
The EM approach is arguably the most popular for conducting QML estimation of high-dimensional DFMs. This approach dates back to the 1970s when it was introduced in a low-dimensional setting by, e.g., \citet{SS77}, \citet{shumwaystoffer82}, \citet{watsonengle83}, and \citet{harveypeters90}, while its use in a high-dimensional setting was first suggested by \citet{quahsargent93} and then formalized by \citet{DGRqml}.

The EM approach for high-dimensional DFMs has been extensively employed by empirical macroeconomic researchers, particularly those in central banks. Its most successful applications include (see also the survey by \citealp{PRM20}):
\begin{inparaenum} 
\item [(i)] counterfactual analysis \citep{harvey96,GRS06,GLR19};
\item [(ii)] conditional forecasts \citep{banburagiannonelenza15};
\item [(iii)] nowcasting (\citealp{Nowcasting,NowcastingReview,BGMR13,kimswanson18,CGLM2023});
\item [(iv)] dealing with data irregularly spaced (\citealp{marianomurasawa03,JKVW2011,banburamodugno14,marcellino16});
\item [(v)] imposing constraints on the loadings to account for smooth cross-sectional dependence in the case of ordered units (\citealp{koopman13,JKVW2014}) or a block-specific factor structure (\citealp{CGM16,altavilla2017});
\item [(vi)] building indicators of economic activity (\citealp{reiswatson10,OGAP,ng2023constructing,AL});
\item [(vi)] impulse response analysis (\citealp{juvenalpetrella2015,smokinggun});
\item [(vii)] modeling international stock market dynamics (\citealp{Linton21});
\item [(viii)] extract trends from micro-panels (\citealp{BCGM21}).
\end{inparaenum}

In addition to the EM approach, the literature has proposed several multi-step approaches to estimate the DFM in \eqref{eq:Fsimple}-\eqref{eq:VARsimple}. 
\begin{inparaenum} [(a)]
\item \citet{baing07} and \citet{FGLR09} employ PC analysis followed by VAR estimation. 
\item \citet{DGRfilter} employ PC analysis followed by VAR estimation and the Kalman smoother.  
\item \citet{ng15} consider QML estimation of the loadings based on a matrix decomposition technique that allows the use of the Newton-Raphson method, then followed by the Kalman filter. 
\item \citet{baili16} consider QML estimation of the loadings based on the EM algorithm by \citet{RT82}, followed by Weighted Least Squares to estimate the factors, VAR estimation, and, finally, the Kalman smoother. 
\item \citet{JK15} consider QML estimation of the loadings for a low-dimensional projection of the data based on the prediction error likelihood obtained from the Kalman filter.
\item \citet{LM02} propose an alternating minimization algorithm based on a penalized loss accounting for cross-autocorrelation in the idiosyncratic components, followed by Generalized Least Squares to estimate the factors.
\item \citet{KM09} consider estimation using sub-space methods. 
\item \citet{MCG24} propose a modified version of the EM algorithm used in this paper where the M-step allows for sparsity in the loadings.
\end{inparaenum}

Approaches (a)-(d) consider estimation of the loadings, either via PC or QML, based only on \eqref{eq:Fsimple}, while estimation of \eqref{eq:VARsimple} is in a second step. Approaches (e)-(h) consider joint estimation of 
 \eqref{eq:Fsimple}-\eqref{eq:VARsimple}.
 However, none of these  fully develops an asymptotic theory for the proposed estimators. Moreover, (f) and (h) do not study the convergence of their numerical algorithms.


Last, an alternative to the EM algorithm is represented by Bayesian estimation of large DFMs using Gibbs sampling, see, e.g., by \citet{KOW03}, \citet{lucianiricci}, \citet{BW15},  \citet{DGLM16}, and \citet{koopman2017empirical}, among many others. 

To conclude, classical references for QML estimation of an exact factor model with no autocorrelations are, e.g., \citet{AR56} and \citet{AFP87}, while \citet{tippingbishop99} suggest to simplify the maximization problem by considering a mis-specified likelihood with homoskedastic idiosyncratic components. \citet{baili12} extend classical QML estimation to the high-dimensional case. Alternatively, \citet{stockwatson89}  combine QML estimation with the Kalman filter by using the prediction error likelihood. A review of these methods is in \citet{MBQML}.



\subsection*{Notation}
An $m\times m$ identity matrix is denoted as $\mbf I_m$. Vectors are always considered as one-column matrices.  An $m$ dimensional vector of ones is denoted as $\bm\iota_m$. An $m$ dimensional vector of zeros is denoted as $\mbf 0_m$, an $m\times p$ matrix of zeros is denoted as $\mbf 0_{m\times p}$.

The generic $(i,j)$ entry of a matrix $\mbf A$ is denoted as $[\mbf A]_{ij}$. Unless otherwise specified, we denote as $\nu^{(k)}(\mbf A)$ the $k$-th largest eigenvalue of a generic squared matrix $\mbf A$. 
The spectral norm for a real $p\times m$ matrix $\mbf A$ is defined by $\Vert \mbf A\Vert=(\nu^{(1)}(\mbf A\mbf A^\prime))^{1/2}$. The Frobenius norm is defined by $\Vert \mbf A\Vert_F=(\mbox{tr}(\mbf A\mbf A^\prime))^{1/2}$. For a generic $p$-dimensional vector $\bm v=(v_1\cdots v_p)^\prime$, we consider the norms: $\Vert \bm v\Vert=(\sum_{j=1}^p v_j^2)^{1/2}$, and $\Vert \bm v\Vert_{\max}=\max_{j=1,\ldots, p} \vert v_j\vert$. 

The indicator function on the event $A$ is denoted as $\mathbb I(A)$, i.e, $\mathbb I(A)=1$ if $A$ is true and 0 otherwise.

All random variables (scalars, vectors and matrices) are assumed to belong to $L_2((\Omega,\mathcal F, \mathrm P))$, where $(\Omega,\mathcal F, \mathrm P)$ is a common probability space. For a generic $p$-dimensional process $\{\mbf y_t\}$ we adopt the following definitions.
\begin{inparaenum}
\item [(i)] Expectations are computed using the true values of the underlying distribution unless otherwise indicated, so we write $\E[\mbf y_t]= \int_{\mathbb R^p} \bm y \mathrm d F_{\mbf y_t}(\bm y,\bm\varphi_n)$, where  $F_{\mbf y_t}(\bm y,\bm\varphi_n)$ is the cumulative distribution function of $\mbf y_t$ computed when using as parameters the true ones $\bm\varphi_n$, and we write $\E_{\wh{\varphi}_n}[\mbf y_t]$ when using as parameters $\wh{\bm\varphi}_n$ to compute the cdf. 
\item [(ii)]  For any $t\in\mathbb Z$, given the $pT$-dimensional vector $\bm Y_{t}=(\mbf y_1^\prime\cdots\mbf y_t^\prime)^\prime$ we denote conditioning on $\bm Y_t$ as an abbreviation for conditioning on the $\sigma$-algebra generated by $\{\mbf y_{t-k}, k\ge 0\}$.
\end{inparaenum}

Limits are always taken as $\min(n,T)\to\infty$ unless otherwise specified. We adopt the Landau $O(\cdot)$ and $o(\cdot)$ notation and the ``in probability'' $O_p(\cdot)$ and $o_p(\cdot)$ analogues. We denote convergence in probability and in distribution by $\stackrel{p}{\to}$ and $\stackrel{d}{\to}$, respectively.

For all quantities having a dimension growing with $n$ and/or $T$, we highlight such dependence. The true scalars, vectors or matrices are denoted as, e.g., $\sigma_i^2$, $\bm\Lambda_n$, $\bm\phi_n$, $\bm\theta$, $\bm F_T$,. The corresponding scalars, vectors or matrices containing generic values of parameters are underlined, so: $\underline{\sigma}_i^2$, $\underline{\bm\Lambda}_n$,  $\underline{\bm\phi}_n$, $\underline{\bm\theta}$, $\underline{\bm F}_T$. 

For a generic process $\{\mbf y_t\}$ and any $T\in\mathbb N$, letting $\mathcal F^0_{-\infty}$ and $\mathcal F_T^\infty$ be the $\sigma$-algebras generated by $\{\mathbf y_t, \, t\le 0\}$ and $\{\mathbf y_t, \, t\ge T\}$, respectively, we define the strong mixing coefficients of  $\{\mbf y_t\}$ as $\alpha_y(T)=\sup_{A\in\mathcal F^0_{-\infty},\, B\in  \mathcal F_T^\infty}\vert \mathrm P(A)\mathrm P(B)-\mathrm P(AB)\vert$.

%
%

\section{Estimation via the EM algorithm}\label{sec:est_nutshell}

In this section, we present the EM algorithm described by \citet{shumwaystoffer82} and implemented in the codes by \citet{DGRqml}, which are available from theirs or ours webpages.\footnote{See the replication codes available at:
\href{https://sites.uw.edu/dgiannon/domenico-giannone-s-homepage/}{https://sites.uw.edu/dgiannon/domenico-giannone-s-homepage/}, or\\ 
\href{http://www.barigozzi.eu/codes.html}{http://www.barigozzi.eu/codes.html}, or
\href{https://sites.google.com/site/lucianimatteo/matlab-codes}{https://sites.google.com/site/lucianimatteo/matlab-codes}} This is the approach typically followed by applied researchers.

Let us assume to observe an $n$-dimensional stochastic process $\mbf x_{nt}=(x_{1t}\cdots x_{nt})^\prime$ over $T$ periods. The DFM \eqref{eq:Fsimple}-\eqref{eq:VARsimple} reads:
\begin{align}
\mbf x_{nt}&=\bm\mu_n+\bm\Lambda_n \mbf F_t +\bm  \xi_{nt}, \label{eq:SDFM1R}\\
\mbf F_t &=\sum_{h=1}^{p_F} \mbf A_h \mbf F_{t-h}
+ \mbf v_t,\label{eq:SDFM2R}
\end{align}
where $\bm\Lambda_n=(\bm\lambda_1\cdots\bm\lambda_n)^\prime$ is the $n\times r$ matrix of factor loadings, $\mbf v_t$ is an $r$-dimensional vector of factor innovations with covariance $\bm\Gamma^v=\E[\mbf v_t\mbf v_t^\prime]$, and $p_F$ is a finite integer such that $p_F\ge 1$. Throughout, we assume that the $r$-dimensional process of factors, $\{\mbf F_t\}$, and the $n$-dimensional process of idiosyncratic components $\{\bm\xi_{nt}\}$ are zero-mean covariance stationary processes, so that $\{\mbf x_{nt}\}$ is also covariance stationary and $\E[\mbf x_{nt}]=\bm\mu_n$, with $\bm\mu_n=(\mu_1\cdots\mu_n)^\prime$.

Let us introduce the following notation. Define the  $nT$-dimensional vectors $\bm X_{nT}=(\mbf x_{n1}^\prime\cdots\mbf x_{nT}^\prime)^\prime$ and $\bm \Xi_{nT}=(\bm\xi_{n1}^\prime\cdots\bm\xi_{nT}^\prime)^\prime$, and the $rT$-dimensional vector $\bm F_T=(\mbf F_1^\prime\cdots\mbf F_T^\prime)^\prime$. Moreover, let $\bm\Sigma_{n}^\xi$ be the diagonal matrix containing the $n$ diagonal terms of $\E[\bm\xi_{nt}\bm\xi_{nt}^\prime]$, which we denote as $\sigma_i^2$, $i=1,\ldots, n$, and ${\bm\Omega}^F_T=\E[\bm F_T\bm F_T^\prime]$. 

In principle, to achieve QML, we should estimate (\textit{i}) $nr$ loadings, (\textit{ii}) $\simeq n^2T^2$ elements of the covariance matrix of $\bm \Xi_{nT}$, (\textit{iii}) $\simeq r^2T^2$ elements of the covariance matrix of $\bm F_T$, and (iv) the $n$ constants in $\bm\mu_n$. Clearly the QML estimator of $\bm\mu_n$ is the sample mean $\bar{\mbf x}_n=T^{-1}\sum_{t=1}^T\mbf x_{nt}$ and we can then work with centered data: $\mbf x_{nt}-\bar{\mbf x}_n$ \citep[p.440]{baili12}. However,
all other parameters depend on each other, so they must be estimated jointly---see, e.g., the first-order conditions in \citet{baili12,baili16} when the autocorrelation of the  factors is not modeled. 
 This is an unfeasible task since we have only $nT$ data points; therefore, we need some regularization to reduce the number of parameters to estimate.  

 To this end, first, we adopt the extreme form of regularization possible for the full-covariance matrix of $\bm\Xi_{nT}$ by replacing it with the diagonal matrix $\mbf I_T\otimes \bm\Sigma_n^\xi$. Second, we assume that the parametric model \eqref{eq:SDFM2R} describes the dynamics of the factors---thus, we write ${\bm\Omega}^F_T\equiv {\bm\Omega}^F_T(\bm{\mathcal A},\bm\Gamma^v)$, with $\bm{\mathcal A}=(\mbf A_1\cdots \mbf A_{p_F})$, in order to highlight its dependence on the VAR parameters. 
 Thanks to these assumptions, we reduce the number of unknown parameters that need to be estimated to $Q_n=nr+n+r^2p_F+r(r+1)/2$, the elements of the vector $\bm{\varphi}_n=(\bm \phi_n^\prime\;\bm\theta^\prime)^\prime$, where $\bm\phi_n=(\text{vec}(\bm\Lambda_n)^\prime\; {\sigma}^2_{1}\cdots {\sigma}^2_{n})^\prime$ and $\bm \theta=(\text{vec}(\bm{\mathcal A})^\prime\; \text{vech}(\bm\Gamma^v)^\prime)^\prime$.
 
Our starting point is then the following log-likelihood, computed in the generic values of the parameters, $\underline{\bm\phi}_n$ and $\underline{\bm\theta}$:
\begin{align}
\ell\big(\bm X_{nT};\underline{\bm\phi}_n,\underline{\bm\theta}\big)=&\,-\frac 12 \log\det\Big(\{\mbf I_T\otimes\underline{\bm\Lambda}_n\}{\bm\Omega}_T^F(\underline{\bm{\mathcal A}}, \underline{\bm\Gamma}^v) \{\mbf I_T\otimes\underline{\bm\Lambda}_n^\prime\}+ \{\mbf I_T\otimes \underline{\bm\Sigma}_n^\xi\}\Big)\label{eq:LL0true}\\
&-\frac 12 \l[ (\bm X_{nT}-\bm\iota_T\otimes\bar{\mbf x}_n)^\prime \Big(
\{\mbf I_T\otimes\underline{\bm\Lambda}_n\}{\bm\Omega}_T^F(\underline{\bm{\mathcal A}}, \underline{\bm\Gamma}^v) \{\mbf I_T\otimes\underline{\bm\Lambda}_n\}^\prime  + \{\mbf I_T\otimes \underline{\bm\Sigma}_n^\xi\}\Big)^{-1}(\bm X_{nT}-\bm\iota_T\otimes\bar{\mbf x}_n) \r],\nn
\end{align} 
where we removed the constant terms to simplify the notation. Since the assumptions in Section \ref{sec:sdfm} allow for correlations among idiosyncratic components, the expression in \eqref{eq:LL0true} is a mis-specified or quasi log-likelihood. Such mis-specification is appealing because it coincides with the classical factor analysis under the exact factor structure, and is standard in DFM estimation \citep{DGRqml,baili16}. Thus, the maximization of \eqref{eq:LL0true} is QML estimation rather than ML estimation.
As long as the idiosyncratic components are weakly correlated, the mis-specification introduced by maximizing \eqref{eq:LL0true} has no effect on consistency but only on  efficiency of the estimators (see Section \ref{sec:PCeff}).
Hereafter, we denote the vector of QML estimators, which are the maximizers of \eqref{eq:LL0true}, as $\wh{\bm\varphi}_n^*=(\wh{\bm\phi}_n^{*\prime}\;\wh{\bm\theta}^{*\prime})^\prime$. 

Despite reducing the number of parameters to be estimated by introducing mis-specifications, direct maximization of \eqref{eq:LL0true} is still unfeasible because $\bm\Omega_T^F$ is a full matrix, and to estimate its entries, we need to estimate the factors as well. Therefore, we still face a curse of dimensionality problem because we need to jointly estimate the $rT$ values of the factors and the $Q_n$ parameters using just the $nT$ observations of $\bm X_{nT}$.

There are three solutions to this problem. The first consists of rewriting the log-likelihood \eqref{eq:LL0true} using its prediction error formulation, where the prediction errors and their covariance are obtained via the Kalman filter (see, e.g., \citealp[Chapter 3.4]{harvey90} or \citealp[Chapter 7]{DK01}). This is the standard practice in low-dimensional DFMs \citep{stockwatson89}. However, since there is no closed form solution for the QML estimator of the parameters obtained in this way, numerical maximization is required and this approach becomes quickly unfeasible even for moderate values of $n$. In high dimensions, \citet{JK15} propose to follow this approach subject to a preliminary step in which the data are projected onto a lower dimensional space but do not derive the theoretical properties of this estimator. In particular, it is not clear what are the effects of this first step on the asymptotic properties of the final estimator.

The second solution, proposed by \citet{baili16}, further mis-specify the model by treating the factors as if they were serially uncorrelated. Thus, by imposing the standard identifying constraint $\E[\mbf F_t\mbf F_t^\prime]=\mbf I_r$ and by replacing in \eqref{eq:LL0true} the full-matrix ${\bm\Omega}_T^F(\underline{\bm{\mathcal A}}, \underline{\bm\Gamma}^v)$ with just $\mbf I_T\otimes \mbf I_r$, the log-likelihood is considerably simplified, and its maximization becomes feasible. However, because it does not exist a closed form solution, we still need a numerical approach to estimate the model, e.g., the iterative algorithm proposed by \citet{RT82} and reintroduced by \citet{baili12}. To our knowledge, the convergence of the \citet{RT82} algorithm to the QML estimator has never been formally proved.

 \begin{table}[ht!]
 \footnotesize
	\caption{}
	\label{alg111}
	\algsetup{indent=1.5em}
	\noindent\rule[0.25ex]{\linewidth}{.7pt}
	\textsc{Expectation Maximization algorithm}
	\vskip -.1cm
	\noindent\rule[0.25ex]{\linewidth}{.7pt}
	\textbf{Input:}\\
	\hspace*{\algorithmicindent} 		$n$ dimensional vector of data $\{\mbf x_{nt}\}_{t=1}^T$\\
	\hspace*{\algorithmicindent}	number of common factors $r$\\
	\hspace*{\algorithmicindent}	 VAR order for the factors $p_F$\\
	\hspace*{\algorithmicindent}	 maximum number of iterations $k_{\max}$\\
	\hspace*{\algorithmicindent}	 threshold for convergence $\varepsilon$\\
	\hspace*{\algorithmicindent} 	initial estimators $\{\wh{\bm\lambda}_i^{(0)}\}_{i=1}^n$, $\{\wh{\sigma}_{i}^{2(0)}\}_{i=1}^n$, $\{\wh{\mbf A}^{(0)}_j\}_{j=1}^{p_F}$, $\wh{\bm\Gamma}^{v(0)}$, see Appendix \ref{app:prest}
		\vskip .2cm
		 \textbf{Ouput:}\\
		\hspace*{\algorithmicindent} estimated loadings $\{\wh{\bm\lambda}_i\}_{i=1}^n$\\
		\hspace*{\algorithmicindent} estimated factors $\{\wh{\mbf F}_t\}_{t=1}^T$
	\begin{algorithmic}[0]		
		
		\FOR{$k=0$ \TO $k=k_{\max}$} 
		\STATE
		\vskip .1cm
		 \COMMENT{\underline{E-STEP}}	\\
		\STATE 	run Kalman Smoother with $\{\mbf x_{nt}\}_{t=1}^T$, $\{\wh{\bm\lambda}_i^{(k)}\}_{i=1}^n$, $\{\wh{\sigma}_{i}^{2(k)}\}_{i=1}^n$, $\{\wh{\mbf A}^{(k)}_j\}_{j=1}^{p_F}$, $\wh{\bm\Gamma}^{v(k)}$\\
		$\rightarrow \{\mbf F_{t|T}^{(k)}\}_{t=1}^T, \{\mbf P_{t|T}^{(k)}\}_{t=1}^T, \{\{\mbf C_{t,t-j|T}^{(k)}\}_{t=j+1}^T\}_{j=1}^{p_F}$, see Appendix \ref{app:KFKS}\\
		\vskip .2cm
		\STATE 
		compute expected log-likelihood and sufficient statistics\\
		$\rightarrow\mathcal Q\l(\{\underline{\bm\lambda}_i\}_{i=1}^n,\{\underline{\sigma}_{i}^{2}\}_{i=1}^n,\{\underline{\mbf A}_j\}_{j=1}^{p_F},\underline{\bm \Gamma}^v;\{\wh{\bm\lambda}^{(k)}_i\}_{i=1}^n,\{\wh{\sigma}_{i}^{2(k)}\}_{i=1}^n,\{\wh{\mbf A}_j^{(k)}\}_{j=1}^{p_F},\wh{\bm\Gamma}^{v(k)}\r)$,\\ 
		$\rightarrow \l\{\l(\sum_{t=1}^T \mbf F_{t|T}^{(k)} x_{it}\r)\r\}_{i=1}^n, \l(\sum_{t=1}^T \mbf F_{t|T}^{(k)}\mbf F_{t|T}^{(k)\prime}+\mbf P_{t|T}^{(k)} \r), \l\{\l(\sum_{t=j+1}^T \mbf F_{t|T}^{(k)}\mbf F_{t-j|T}^{(k)\prime}+\mbf C_{t,t-j|T}^{(k)}\r)\r\}_{j=1}^{p_F}$,\\
		see \eqref{eq:LLbayes_exp_bis}, \eqref{eq:LLbayes_exp_bis1}, \eqref{eq:LLbayes_exp_bis2},  \eqref{eq:suffstat0}, \eqref{eq:suffstatAPP}
		\vskip .2cm
		
		 \COMMENT{\underline{M-STEP}}\\
		\STATE maximize expected log-likelihood\\
		$\rightarrow \{\wh{\bm\lambda}_i^{(k+1)}\}_{i=1}^n$, $\{\wh{\sigma}_{i}^{2(k+1)}\}_{i=1}^n$, $\{\wh{\mbf A}_j^{(k+1)}\}_{j=1}^{p_F}$, $\wh{\bm\Gamma}^{v(k+1)}$,\\
		see \eqref{eq:param1}, \eqref{eq:param3}, \eqref{eq:param4}, \eqref{eq:paramGv}
		\vskip .2cm
		
		\COMMENT{\underline{CONVERGENCE}}
		\IF{$k< k_{\max}$ AND $\Delta\ell_{k}<\varepsilon$, see \eqref{eq:convEM}}  
		\STATE $k^*\leftarrow k$
		\STATE run Kalman Smoother with $\{\mbf x_{nt}\}_{t=1}^T$, $\{\wh{\bm\lambda}_i^{(k^*+1)}\}_{i=1}^n$, 
		$\{\wh{\sigma}_{i}^{2(k^*+1)}\}_{i=1}^n$, $\{\wh{\mbf A}_j^{(k^*+1)}\}_{j=1}^{p_F}$, $\wh{\bm\Gamma}^{v(k^*+1)}$\\  
		$\rightarrow \{{\mbf F}_{t|T}^{(k^*+1)}\}_{t=1}^T$, see Appendix \ref{app:KFKS}
		\STATE $\wh{\bm\lambda}_i \leftarrow \wh{\bm\lambda}_i^{(k^*+1)}$, for all $i=1,\ldots, n$
		\STATE $\wh{\mbf F}_t \leftarrow {\mbf F}_{t|T}^{(k^*+1)}$, for all $t=1,\ldots, T$
		\RETURN $\{\wh{\bm\lambda}_i\}_{i=1}^n$, $\{\wh{\mbf F}_t\}_{t=1}^T$
		\STATE {\bf break}	
		\ENDIF
		\vskip .1cm
		\IF {$k=k_{\max}$}
		\PRINT ``algorithm did not converge''
		\STATE {\bf break}
		\ENDIF
		\vskip .1cm
		\STATE $k+1\leftarrow k$
		\vskip .1cm
		\ENDFOR
	\end{algorithmic} 	
	\vskip -.2cm
	\noindent\rule[0.25ex]{\linewidth}{.7pt}
\end{table}

A third solution consists of computing an approximation of the QML estimator using the EM algorithm, as formalized by \citet{DGRqml}. We consider this approach in this paper, and we refer to Table \ref{alg111} for its implementation.

The EM algorithm is an iterative procedure which allows for QML estimation in presence of missing data \citep{DLR77}. In a nutshell, consider a given iteration $k\ge 0$ and assume to have an estimate of the parameters $\wh{\bm\varphi}_n^{(k)}$. By taking expectations of the log-likelihood \eqref{eq:LL0true} with respect to the conditional distribution of $\bm F_T$ given $\bm X_{nT}$ and computed using $\wh{\bm\varphi}_n^{(k)}$, we get:
\begin{align}
\ell(\bm X_{nT};\underline{\bm\varphi}_n)&=\E_{\wh{\varphi}_{n}^{(k)}}[\ell(\bm X_{nT},\bm F_T;\underline{\bm\varphi}_n)|\bm X_{nT}]
-\E_{\wh{\varphi}_{n}^{(k)}}[\ell(\bm F_T|\bm X_{nT};\underline{\bm\varphi}_n)|\bm X_{nT}] \nn\\
&=\mathcal Q(\underline{\bm\varphi}_n,\wh{\bm\varphi}_{n}^{(k)})-\mathcal H(\underline{\bm\varphi}_n,\wh{\bm\varphi}_{n}^{(k)}), \; \text{say}.\label{eq:LLbayes_exp}
\end{align}
A maximum of the log-likelihood is then a maximum of the right hand side of \eqref{eq:LLbayes_exp}. Now, by definition of Kullback-Leibler divergence,  for any $k\ge0$ it holds that
\beq\label{eq:timmagini}
{\mathcal H}(\wh{\bm\varphi}_{n}^{(k+1)};\wh{\bm\varphi}_n^{(k)})\leq {\mathcal H}(\wh{\bm\varphi}_n^{(k)};\wh{\bm\varphi}_n^{(k)}),
\eeq
i.e., $\mathcal H(\underline{\bm\varphi}_n,\wh{\bm\varphi}_{n}^{(k)})$ is maximum at $\underline{\bm\varphi}_n=\wh{\bm\varphi}_{n}^{(k)}$. It is then enough to look for the maximum only of the expected full-information log-likelihood $\mathcal Q(\underline{\bm\varphi}_n,\wh{\bm\varphi}_{n}^{(k)})$. 
This is accomplished in two steps: in the first step, for given estimated parameters $\wh{\bm\varphi}_{n}^{(k)}$, we compute $\mathcal Q(\underline{\bm\varphi}_n,\wh{\bm\varphi}_{n}^{(k)})$ using an estimate of the factors with their associated MSE; in the second step, for a given estimate of the factors, we maximize such log-likelihood to compute a new estimate of the parameters $\wh{\bm\varphi}_{n}^{(k+1)}$. As shown below, this approach solves the curse of dimensionality problem, and its computational burden is minimal because all estimates have an explicit expression. Below we detail the main features of the two steps.

\subsection{E-step and Kalman smoother}
We obtain an initial estimate of the loadings and the idiosyncratic variances using the PC estimator, and of the VAR parameters by fitting a VAR on the PC estimator of the factors (see Appendix \ref{app:prest}). Then, for any iteration $k\ge 0$, in the E-step, we compute the expected full-information log-likelihood, which we can decomposed as:
\begin{align}\label{eq:LLbayes_exp_bis}
\mathcal Q(\underline{\bm\varphi}_n,\wh{\bm\varphi}_{n}^{(k)})
=&\,\E_{\wh{\varphi}_{n}^{(k)}}[\ell(\bm X_{nT}|\bm F_T;\underline{\bm\phi}_n)|\bm X_{nT}]+\E_{\wh{\varphi}_{n}^{(k)}}[\ell(\bm F_T;\underline{\bm\theta})|\bm X_{nT}].
\end{align}
Consistently with the mis-specified log-likelihood \eqref{eq:LL0true}, the first term on the right-hand side of \eqref{eq:LLbayes_exp_bis} is:
\begin{align}
\ell(\bm X_{nT}|\bm F_T;\underline{\bm \phi}_{n})&= -\frac T2 \log\det(\underline{\bm\Sigma}_n^\xi) -\frac 12\sum_{t=1}^T (\mbf x_{nt}-\bar{\mbf x}_n-\underline{\bm\Lambda}_n{\mbf F}_t)^\prime (\underline{\bm\Sigma}_n^\xi)^{-1}(\mbf x_{nt}-\bar{\mbf x}_n-\underline{\bm\Lambda}_n{\mbf F}_t)\nn\\
&=\sum_{i=1}^n\l\{-\frac T2 \log(\underline{\sigma}^2_{i})-\frac 12\sum_{t=1}^T\frac{(x_{it}-\bar x_i-\underline{\bm\lambda}_i^\prime\mbf F_t)^2}{\underline{\sigma}^2_{i}}
\r\},
\label{eq:LLbayes_exp_bis1}
\end{align}
which depends only on $\underline{\bm\phi}_n$ and is well defined as long as all the idiosyncratic components have finite positive variances. 

As for the second term on the right-hand side of \eqref{eq:LLbayes_exp_bis}, we assume that $\mbf F_{t}=\mbf 0_r$ for $t\le 0$ and consider
$\ell(\bm F_T;\underline{\bm\theta})= \sum_{t=1}^T \ell(\mbf F_t|\mbf F_{t-1},\ldots, \mbf F_{t-p_F};\underline{\bm\theta})$,
which depends only on $\underline{\bm\theta}$. 
Then we have
\begin{align}
\ell(\bm F_T;\underline{\bm\theta})=&\,-\frac {T}2 \log\det (\underline {\bm\Gamma}^v)-\frac 12\!\sum_{t=1}^T\l(\mbf F_t-\sum_{h=1}^{p_F}\underline{\mbf A}_h\mbf F_{t-h}\r)^\prime(\underline {\bm\Gamma}^v)^{-1}\l(\mbf F_t-\sum_{h=1}^{p_F}\underline{\mbf A}_h\mbf F_{t-h}\r)\label{eq:LLbayes_exp_bis2},
\end{align}
which is well defined provided that the VAR innovations $\{\mbf v_t\}$ have a finite full-rank covariance matrix. 

Given \eqref{eq:LLbayes_exp_bis1} and \eqref{eq:LLbayes_exp_bis2}, in order to compute the expected log-likelihood \eqref{eq:LLbayes_exp_bis} we need to compute the sufficient statistics:
\begin{align}
\E_{\wh{\varphi}_{n}^{(k)}}[{\mbf F}_t|\bm X_{nT}], \qquad 
\E_{\wh{\varphi}_{n}^{(k)}}[{\mbf F}_t{\mbf F}_t^\prime|\bm X_{nT}], 
\qquad 
\E_{\wh{\varphi}_{n}^{(k)}}[{\mbf F}_t{\mbf F}_{t-h}^\prime|\bm X_{nT}], \quad h=1,\ldots, p_F.\label{eq:suffstat0}
\end{align}
Although exact expressions for the quantities in \eqref{eq:suffstat0} might be hard to compute, we can approximate them by using the output of the Kalman smoother, which gives the linear projection ${\mbf F}_{t|T}^{(k)}=\mathrm {Proj}_{\wh{\varphi}_{n}^{(k)}}[{\mbf F}_t|\bm X_{nT}]$ and the associated conditional covariance and lag-$h$ autocovariance matrices ${\mbf P}^{(k)}_{t|T}$ and ${\mbf C}^{(k)}_{t,t-h|T}$, respectively (see Appendix \ref{app:KFKS}, for the explicit expressions). This approximation does not affect consistency of our estimators (see Proposition \ref{prop:plainvanilla}). 

Hence, in the EM algorithm we let: 
\begin{align}
&\E_{\wh{\varphi}_{n}^{(k)}}[{\mbf F}_t|\bm X_{nT}]= {\mbf F}_{t|T}^{(k)},\qquad \qquad \E_{\wh{\varphi}_{n}^{(k)}}[{\mbf F}_t{\mbf F}_t^\prime|\bm X_{nT}] = {\mbf F}_{t|T}^{(k)}{\mbf F}_{t|T}^{(k)\prime} +{\mbf P}_{t|T}^{(k)},\label{eq:suffstatAPP}\\
&\E_{\wh{\varphi}_{n}^{(k)}}[{\mbf F}_t{\mbf F}_{t-h}^\prime|\bm X_{nT}]={\mbf F}_{t|T}^{(k)}{\mbf F}_{t-h|T}^{(k)\prime} +{\mbf C}_{t,t-h|T}^{(k)},\qquad h=1,\ldots, p_F.\nn
\end{align}

Summing up, in the E-step, we compute $rT$ values of the factors for a given value, $\wh{\varphi}_{n}^{(k)}$, of the parameters. In general, this step is feasible because we have $T(n-r)\gg 0$ degrees of freedom. Moreover, since the Kalman filter and smoother are linear procedures entailing the inversion of a positive definite $r\times r$ matrix, we just have to compute $\simeq T$ recursions.

\begin{rem}\label{rem:PtTMSE}
\upshape{
The Kalman smoother requires first running the forward iterations of the Kalman filter, which in turn requires either inverting the $n\times n$ full covariance matrix of the data or inverting the full $n\times n$ idiosyncratic covariance. This task might be challenging, if not impossible, in a high-dimensional setting. To overcome this problem, we implement the Kalman filter using an estimator of the diagonal matrix $\bm\Sigma_n^\xi$ instead of the full idiosyncratic covariance matrix, ensuring also positive definiteness of the data covariance matrix, and thus making also its inversion feasible. 
This simplification is consistent with the mis-specified log-likelihood \eqref{eq:LLbayes_exp_bis1} because, to guarantee that \eqref{eq:timmagini} holds, we must take expectations with respect to the same distribution as the one used to compute the log-likelihood.
}
\end{rem}

\subsection{M-step}\label{sec:Mstep}
In the M-step, we have to maximize \eqref{eq:LLbayes_exp_bis} with respect to $\underline{\bm\varphi}_n$ to obtain a new estimate of the parameters $\wh{\bm\varphi}_n^{(k+1)}$. This maximization has a closed form solution for all elements of $\wh{\bm\varphi}_n^{(k+1)}$. Specifically, at a given iteration $k\ge 0$, by using \eqref{eq:suffstatAPP}, we obtain the loadings estimators as 
\begin{align}
\wh{\bm\lambda}_{i}^{(k+1)} &= 
\l\{\sum_{t=1}^T \l({\mbf F}_{t|T}^{(k)}{\mbf F}_{t|T}^{(k)\prime}+{\mbf P}^{(k)}_{t|T}\r)\r\}^{-1}\l(\sum_{t=1}^T {\mbf F}_{t|T}^{(k)}\, (x_{it}-\bar x_i)\r), \qquad \text{for } i=1,\ldots, n\label{eq:param1}
\end{align}
and $\wh{\bm\Lambda}_n^{(k+1)}=(\wh{\bm\lambda}_{1}^{(k+1)} \cdots \wh{\bm\lambda}_{n}^{(k+1)} )^\prime$. Similarly, the estimator of the idiosyncratic variances is:
\begin{align}
\wh{\sigma}_i^{2(k+1)}
=&\;\frac 1T\sum_{t=1}^T\l\{x_{it}^2+ \wh{\bm\lambda}_i^{(k+1)\prime}\l({\mbf F}_{t|T}^{(k)}{\mbf F}_{t|T}^{(k)\prime}+{\mbf P}^{(k)}_{t|T}\r) \wh{\bm\lambda}_i^{(k+1)} -2x_{it}{\mbf F}_{t|T}^{(k)\prime}\wh{\bm\lambda}_i^{(k+1)}\r\} \qquad \text{for } i=1,\ldots, n.
\label{eq:param3}
\end{align}
Because we consider a mis-specified log-likelihood, we do not estimate the out-of-diagonal terms of the idiosyncratic covariance matrix $\bm\Gamma_n^\xi$, and we act as if those terms are equal to  zero.

For simplicity, let $p_F=1$ and denote $\mbf A\equiv \mbf A_1$. The estimator of $\mbf A$ is then given by:
\begin{align}
\wh{\mbf A}^{(k+1)}
&=  \l\{\sum_{t=2}^T\l({\mbf F}_{t|T}^{(k)}{\mbf F}_{t-1|T}^{(k)\prime}+{\mbf C}^{(k)}_{t,t-1|T} \r)\r\}
\l\{\sum_{t=2}^T \l({\mbf F}_{t-1|T}^{(k)}{\mbf F}_{t-1|T}^{(k)\prime}+{\mbf P}^{(k)}_{t-1|T}\r)\r\}^{-1}.\label{eq:param4}
\end{align}
For $p_F>1$ we can simply write the VAR in companion form and derive the analogous of \eqref{eq:param4}.

Finally, the estimator of the covariance matrix of the VAR innovations $\{\mbf v_t\}$ is:
\begin{align}
\wh{\bm\Gamma}^{v (k+1)}
=\frac 1T\sum_{t=2}^T&
\l\{
{\mbf F}_{t|T}^{(k)}{\mbf F}_{t|T}^{(k)\prime}+{\mbf P}^{(k)}_{t|T} 
-\wh{\mbf A}^{(k+1)}\l({\mbf F}_{t-1|T}^{(k)}{\mbf F}_{t-1|T}^{(k)\prime}+{\mbf P}^{(k)}_{t-1|T}\r)\wh{\mbf A}^{(k+1)\prime}\r.\label{eq:paramGv}\\
&\l. \phantom{\{}- \l({\mbf F}_{t|T}^{(k)}{\mbf F}_{t-1|T}^{(k)\prime}+{\mbf C}^{(k)}_{t,t-1|T}\r)\wh{\mbf A}^{(k+1)\prime}
-\wh{\mbf A}^{(k+1)} \l({\mbf F}_{t-1|T}^{(k)}{\mbf F}_{t|T}^{(k)\prime}+{\mbf C}^{(k)\prime}_{t,t-1|T}\r)
\r\}.\nn
\end{align}

Summing up, in the M-step, we need to compute $Q_n$ values of the parameters for a given estimator of the factors, $\mbf F_{t|T}^{(k)}$, $t=1,\ldots, T$. This step is feasible because we used the mis-specified log-likelihood \eqref{eq:LLbayes_exp_bis1} to estimate $\bm\phi_n$. Therefore, we decomposed our estimation problem into $n$ separate maximizations, each requiring estimating $r+1$ parameters using $T$ observations. In this way, estimating the high-dimensional parameter vector $\bm\phi_n$ becomes straightforward, and estimating $\bm\theta$ poses no problem because it is a low-dimensional problem with a closed-form solution.

\subsection{Convergence of the EM algorithm and final estimators}\label{sec:contra}

Given the Gaussian quasi-likelihoods \eqref{eq:LLbayes_exp_bis1} and \eqref{eq:LLbayes_exp_bis2}, it is easy show that, for any fixed $n$, there exists an $\omega>0$ such that for any $k\ge0$
\al{
\l\{\ell(\bm X_{nT};\wh{\bm\varphi}_n^{(k+1)})-\ell(\bm X_{nT};\wh{\bm\varphi}_n^{(k)})\r\}&\ge \l\{ \mathcal Q(\wh{\bm\varphi}_n^{(k+1)},\wh{\bm\varphi}_n^{(k)})-\mathcal Q(\wh{\bm\varphi}_n^{(k)},\wh{\bm\varphi}_n^{(k)})\r\}\ge  \omega \Vert\wh{\bm\varphi}_n^{(k+1)}-\wh{\bm\varphi}_n^{(k)} \Vert^2,\label{eq:contra}
}
where the first inequality follows from \citet[Lemma 1]{DLR77} and the second is due to strong concavity of  $\mathcal Q(\cdot,\underline{\bm\varphi}_n)$ for any $\underline{\bm\varphi}_n$ \citep[Condition 1]{wu83}.  Moreover, the left-hand-side of \eqref{eq:contra} tends to zero as $k\to\infty$ \citep[Theorem 3]{wu83}. Therefore, the EM algorithm defines a contractive map. Consequently, the sequence $\{\wh{\bm\varphi}_n^{(k)}\}$ will converge to a maximum of the log-likelihood, as $k\to\infty$ (see Lemma \ref{lem:localmax} for a formal proof when $n\to\infty$).

In practice, we stop the EM algorithm when the log-likelihood shows no further appreciable increase. This is ensured according to a standard convergence rule (see Appendix \ref{rem:kstar}), depending on a pre-specified threshold $\varepsilon$.  

We denote the last iteration of the EM algorithm as $k^*$. Upon convergence, the EM estimator of the parameters is $\wh{\bm\varphi}_n\equiv \wh{\bm\varphi}_n^{(k^*+1)}$. 
By running the Kalman smoother one last time using $\wh{\bm\varphi}_n$, we have the estimator of the factors $\wh{\mbf F}_t=\mbf F_{t|T}^{(k^*+1)}$, $t=1,\ldots, T$. Finally, we estimate the common components as $\wh{\chi}_{it}=\wh{\bm\lambda}_i^\prime\wh{\mbf F}_t$, where $\wh{\bm\lambda}_i\equiv \wh{\bm\lambda}_i^{(k^*+1)}$,  $i=1,\ldots, n$.

%
%
\section{The Large Approximate Dynamic Factor Model}\label{sec:sdfm}

\subsection{Main assumptions}
This section presents the assumptions under which we can consistently estimate the DFM given in  \eqref{eq:SDFM1R}-\eqref{eq:SDFM2R}. These assumptions are stated for an infinite-dimensional stochastic process $\{ x_{it},\, i\in\mathbb N,\,t\in\mathbb Z\}$ of which $\{\mbf x_{nt}=(x_{1t}\cdots x_{nt})^\prime,\,t\in\mathbb Z\}$ is an $n$-dimensional subprocess and the $T\times n$ matrix $(\mbf x_{n1}\cdots\mbf x_{nT})^\prime$ is an observed realization.
Likewise, $\{\bm\xi_{nt}=(\xi_{1t}\cdots \xi_{nt})^\prime,\,t\in\mathbb Z\}$ is an $n$-dimensional sub-process of the infinite-dimensional stochastic process of idiosyncratic components $\{\xi_{it},\, i\in\mathbb N,\,t\in\mathbb Z\}$, and $\{\mbf F_t=(F_{1t}\cdots F_{rt})^\prime,\, t\in\mathbb Z\}$ and $\{\mbf v_t=(v_{1t}\cdots v_{rt})^\prime,\,t\in\mathbb Z\}$ are the $r$-dimensional processes of common factors and the corresponding VAR innovations, respectively. The $n\times r$ matrix of factor loadings $\bm\Lambda_n=(\bm\lambda_1\cdots\bm\lambda_n)^\prime$ forms a nested sequence as $n$ increases.
Finally, as already noticed, the QML estimator of $\bm\mu_n$ is immediately obtained as the sample mean $\bar{\mbf x}_n$. Thus, hereafter, for simplicity, we consider the DFM for pre-centered data, or, equivalently, we set $\bm\mu_n=\mbf 0_n$.


\begin{ass}[\textsc{loadings and factors}]\label{ass:common} $\,$
\begin{compactenum}[(a)]
	\item  There exists an integer $N_0$ such that for all $n> N_0$, $\Vert n^{-1}\bm\Lambda_n^\prime\bm\Lambda_n-\bm\Sigma_{\Lambda}\Vert=0$, where $\bm\Sigma_{\Lambda}$ is $r\times r$ and positive definite; moreover, for all $n\in\mathbb N$, $m_\lambda\le\max_{i=1,\ldots,n}\Vert\bm\lambda_i\Vert\le M_\lambda$ for some finite positive reals $M_\lambda$ and $m_\lambda$ 	independent of $n$.  

	\item For all $t\in\mathbb Z$, $\bm\Gamma^F=\E_{}[\mbf F_t\mbf F_t^\prime]$ is $r\times r$ and positive definite, and $\Vert\bm\Gamma^F\Vert\le M_F$ for some finite positive real $M_F$.

	\item There exists an integer $N_1$ such that for all $n> N_1$, $r$ is a finite positive integer, independent of $n$, and such that $r\le N_1$.

	\item $\mbf {A}(z)=\sum_{k=1}^{p_F} \mbf{A}_{k}z^{k-1}$, such that $p_F$ is a finite positive integer,
	${\mbf A}_{k}$ are $r\times r$,  and $\det(\mbf I_r-{\mbf A}(z))\ne 0$ for all $z\in\mathbb C$ such that $|z|\le M_A$ for some finite positive real $M_A< 1$. 

 	\item For all $t\in\mathbb Z$, $\E_{}[\mbf v_t]=\mbf 0_r$, $\bm\Gamma^v=\E_{}[\mbf v_t\mbf v_t^\prime]$ is $r\times r$ positive definite, and $\Vert\bm\Gamma^v\Vert\le M_v$ for some finite positive real $M_v$.
 
 	\item For all $t\in\mathbb Z$ and all $k\in\mathbb Z$ with $k\ne 0$, $\mbf v_{t}$ and $\mbf v_{t-k}$ are independent.
 
 	\item [(g)] For all $j_1,j_2,j_3,j_4=1,\ldots, r$, all $t=1,\ldots, T$, and all $T\in\mathbb N$, 
	$$
	\frac 1T\sum_{s_1,s_2=1}^T \vert \E[v_{j_1 s_1}v_{j_2 t} v_{j_3 s_2} v_{j_4 t}]\vert \le \mathrm K_v, \
	\quad 
	\frac 1T\sum_{s_1,s_2=1}^T \vert \E[v_{j_1 s_1}v_{j_2 t}]\vert \vert \E[ v_{j_3 s_2} v_{j_4 t}]\vert \le \mathrm K_v, 
	$$
	 for some finite positive real $\mathrm K_v$ independent of $j_1,j_2,j_3,j_4$, and $T$. 
 
 	\item [(h)] For all $t\in\mathbb Z$, $\mbf v_t$ has pdf $f_{\mbf v_t}(\bm u)$ such that $\int_{\mathbb R^r} \vert f_{\mbf v_t}(\bm u+\bm v)-f_{\mbf v_t}(\bm v)\vert \mathrm d\bm u\le C_f \Vert \bm v\Vert$ for any $\bm v\in\mathbb R^r$ and for some finite positive real $C_f$ independent of $t$;
	\item [(i)] For all $t\le 0$, $\mbf v_t = \mbf 0_r$.
\end{compactenum}
\end{ass}
 
Parts (a) and (b) imply that the loadings matrix has asymptotically maximum column rank $r$ (part (a)), and the factors have a finite full-rank covariance matrix (part (b))---these assumptions are similar to the requirements in \citet[Assumptions A and B]{baili16} and \citet[Assumptions A and B]{Bai03}.  Moreover, because of part (a), for any given $n\in\mathbb N$, all the factors have a finite contribution to each series (upper bound on $\max_{i=1,\ldots, n}\Vert\bm\lambda_i\Vert$), and there is at least one factor that contributes to at least one series (lower bound on $\max_{i=1,\ldots, n}\Vert\bm\lambda_i\Vert$). While the former condition is common, the latter is less standard but very mild as it simply guarantees that at least one loading is non-zero for any fixed $n\in\mathbb N$. 

Part (c) implies the existence of a finite number of common factors $r$. In particular, $r$ is identified only for $n\to\infty$ (see also the next section). Hereafter, in parts (a) and (c), we can assume $N_0=N_1=N$, say, without loss of generality. 

The remaining conditions of Assumption \ref{ass:common} characterize the VAR for the factors in \eqref{eq:SDFM2R}.
Part (d) implies that $\{\mbf F_t\}$ is a weakly stationary process with a causal autoregressive representation. 
And in parts (e), (f), and (g), we assume that $\{\mbf v_t\}$ is a zero-mean $r$-dimensional independent process with finite positive definite covariance matrix and  finite summable 4th order cumulants. Parts (d) and (e) imply also that $\bm\Gamma^F$ is finite, as required in part (b).

Part (h) is an integral Lipschitz condition, which is satisfied by most continuous densities. This assumption guarantees that $\{\mbf F_t\}$ is a strong mixing, or equivalently $\alpha$-mixing, process with mixing coefficients $\alpha_F(T)\le \exp\l\{-c_F T^{\gamma_F}\r\}$, for all $T\in\mathbb N$ and some finite positive reals $c_F$ and $\gamma_F$ independent of $T$ (\citealp[Theorem 3.1]{PT85}).\footnote{Independence in part (f) is not strictly necessary for having $\{\mbf F_t\}$ strong mixing, as we could allow for GARCH effects by assuming geometric ergodicity of $\{\mbf v_t\}$ instead (\citealp{FZ06}). Indeed, geometric ergodicity implies $\beta$-mixing, which implies strong mixing. } 
Strongly mixing factors with exponentially decaying  mixing coefficients are directly assumed by \citet[Assumption 2c]{FLM13}.

Part (i) implies $\mbf F_t=\mbf 0_r$ for $t\le 0$. 
This assumption is standard and it fixes the initial conditions for the solution of the VAR in \eqref{eq:SDFM2R}. 


\begin{ass}[\textsc{idiosyncratic component}]\label{ass:idio}
  $\,$
\begin{compactenum}[(a)]
	\item For all $i\in\mathbb N$ and all $t\in\mathbb Z$, $\E[\xi_{it}]=0$ and $\sigma_i^2=\E[\xi_{it}^2]$ is such that $C_\xi^{-1}\le \sigma_i^2\le C_\xi$, for some finite positive real $C_\xi$ independent of $i$.

	\item For all $i,j\in\mathbb N$, all $t\in\mathbb Z$, and all $k\in\mathbb Z$, $\vert \E[\xi_{it}\xi_{j,t-k}]\vert\le \rho^{\vert k\vert} M_{ij}$, where $\rho$ and $M_{ij}$ are finite positive reals independent of $t$ such that $0\le \rho <1$, $M_{ii}=\sigma_i^2$, $\sum_{j=1, j\ne i}^n M_{ij}\le M_\xi$, and $\sum_{i=1, i\ne j}^n M_{ij}\le M_\xi$ for some finite positive real $M_{\xi}$ independent of $n$. 

	\item For all $i\in\mathbb N$,  $\{\xi_{it}\}$ is a strong mixing process with mixing coefficients such that $\alpha_{\xi_i}(T)\le \exp(-c_{\xi} T^{\gamma_{\xi}})$, for all $T\in\mathbb N$, and for some finite positive reals $c_\xi$ and $\gamma_\xi$ independent of $T$ and $i$. 

	\item [(d)] 
	For all $j=1,\ldots, n$ and all $n,T\in\mathbb N$, 
	\al{
	&\frac 1{nT}\sum_{t,s=1}^T\sum_{i_1,i_2=1}^n \vert \E[ \xi_{i_1t}\xi_{jt}\xi_{i_2s}\xi_{js}]\vert \le \mathrm K_\xi, 
	\quad \frac 1{nT}\sum_{t,s=1}^T\sum_{i_1,i_2=1}^n \vert \E[ \xi_{i_1t}\xi_{jt}]\vert \vert\E[\xi_{i_2s}\xi_{js}]\vert \le \mathrm K_\xi,
	\nn
	}
	and for all $t=1,\ldots, T$ and all $n,T\in\mathbb N$, 
	\al{
	&\frac 1{nT}\sum_{s_1,s_2=1}^T\sum_{i,j=1}^n \vert \E[ \xi_{is_1}\xi_{it}\xi_{js_2}\xi_{jt}]\vert \le \mathrm K_\xi, 
	\quad \frac 1{nT}\sum_{s_1,s_2=1}^T\sum_{i,j=1}^n \vert \E[ \xi_{is_1}\xi_{it}]\vert \vert \E[\xi_{js_2}\xi_{jt}]\vert \le \mathrm K_\xi, \nn
	}
	for some finite positive real $\mathrm K_\xi$ independent of $n$, and $T$.

	\item [(e)] For all $t\in\mathbb Z$, as $n\to\infty$,
	\[
	\frac 1{\sqrt n}\sum_{i=1}^n
	 \frac{\bm\lambda_i  \xi_{it}}{\sigma_i^2}
	\stackrel{d}{\to}
	\mathcal N\l(\mbf 0_r, \lim_{n\to\infty}\frac 1 n\sum_{i,j=1}^n 
	\frac{\bm\lambda_i\bm\lambda_j^\prime\E_{}[\xi_{it}\xi_{jt}]}{\sigma_i^2\sigma_j^2}
	\r).
	\]
	\item [(f)] For all $n\in\mathbb N$,  $\bm\Gamma_n^\xi= \E[\bm\xi_{nt}\bm\xi_{nt}^\prime]$ is such that
	$\nu^{(n)}(\bm\Gamma_n^\xi)\ge L_\xi$ for some finite positive real $L_\xi$ independent of $n$.

 \end{compactenum}
 \end{ass}

Part (a) imposes that the idiosyncratic components have zero mean, and the idiosyncratic variances are positive and finite. We strengthen this assumption in part (f) by requiring that the whole idiosyncratic covariance matrix is positive definite.

Part (b) has a twofold purposes. First, it limits the degree of serial correlation of the idiosyncratic components by assuming standard geometric decay of the autocovariances. Second, it limits the degree of cross-sectional correlation between idiosyncratic components, a standard assumption for approximate DFMs. This assumption implies the usual conditions required by \citet[Assumptions C2, C3, and C4]{Bai03}, \citet[Assumption 2b]{FLM13}, and \citet[Assumptions C3, C4, and E1]{baili16} (see Lemma \ref{lem:Gxi}). 

In part (c), we assume that each idiosyncratic component is strongly mixing with exponentially decaying coefficients. This requirement is quite general since it allows the idiosyncratic to be non linear processes---\citet[Assumption 2c]{FLM13} make the same assumption. 

In part (d), we require finite summable fourth-order cumulants---a standard requirement found in the literature (see, e.g., \citealp[Assumption F1]{Bai03} for the first condition, and \citealp[Assumption E2]{baili16} for the second one---which, jointly with the mixing assumption in part (c), allows for consistent estimation of  covariances by means of the sample covariances 
\citep[Theorem 6, Chapter IV, p. 210]{hannan}.

In part (e), we assume a Central Limit Theorem, a standard assumption in the literature (e.g., \citealp[Assumption F3]{Bai03}, and \citealp[Assumption F3]{baili16}). This is a high-level requirement and in order to properly derive it, we should introduce some notion of dependence for random fields, which typically requires some ordering of the cross-sectional units. Now, in most applications, there is not a natural ordering of the variables. For this reason, we avoid to spell out primitive conditions that guarantee part (e) to hold.

We then add the natural requirement  of independence between common shocks and idiosyncratic components.

\begin{ass}[\textsc{Independence between common shocks and idiosyncratic components}]\label{ass:ind} The processes
$\{\xi_{it},\, i\in\mathbb N,\, t\in\mathbb Z\}$ and $\{v_{jt},\, j=1,\ldots, r,\, t\in\mathbb Z\}$ are mutually independent. 
\end{ass}

Assumption \ref{ass:ind} implies that the factors and the common components are independent of the idiosyncratic components at all leads and lags and across all units. This assumption is compatible with the idea that the structural macroeconomic shocks driving the common component are independent of the idiosyncratic components representing measurement errors or local dynamics.\footnote{In principle, we could relax this requirement to allow for weak dependence as in, e.g., \citet[Assumption D]{Bai03}.}

Last,  Assumption \ref{ass:ind} jointly with Assumptions \ref{ass:common}(f), \ref{ass:common}(g), \ref{ass:common}(h), \ref{ass:idio}(c), and \ref{ass:idio}(d), implies that the process $\{\mbf F_t\xi_{it}\}$ is strongly mixing with finite fourth moments \citep[Theorem 5.1.a]{bradley05}. Then, from \citet[Theorem 1.7]{ibra62}, we have the following Central Limit Theorem, for all $i\in\mathbb N$, as $T\to\infty$,
\beq
\frac 1{\sqrt T}\sum_{t=1}^T\mbf F_t\xi_{it}\stackrel{d}{\ra} \mathcal N\l(\mbf 0_r,\lim_{T\to\infty}\frac 1 T\sum_{s,t=1}^T\E_{} \l[\mbf F_t\mbf F_s^\prime \xi_{it}\xi_{is}\r]
\r).\nn
\eeq
This last result is typically assumed in the literature---see, e.g., \citealp[Assumption F4]{Bai03}, and \citealp[Assumption F1]{baili16}.

\subsection{Additional assumptions}
We now state two more assumptions. These are needed only to derive some of our asymptotic results.

\begin{ass}[\textsc{Linearity}]\label{ass:linear}
For all $j=1,\ldots, r$, $t=1,\ldots, T$, and $n,T\in\mathbb N$, $\E[F_{jt}|\bm X_{nT}]$ is a linear function of $\bm X_{nT}$.
\end{ass}

Assumption \ref{ass:linear} would hold if we directly assume joint Gaussianity of $\bm X_{nT}$ and $\bm F_T$. However, assuming Gaussianity might be too stringent, while Assumption \ref{ass:linear} is more general. For example, consider the case $r=1$, let $\mathrm f_t$ and $\bm{\mathrm X}_{nT}$ be realizations of the factor and the data, and let  $f:\mathbb R^{nT+1}\to\mathbb R$ be the joint pdf of $F_t$ and $\bm X_{nT}$. Then, from \citet{steyn1960regression} and \citet[Chapter 44.3]{kotz2004continuous}, we see that if $f$ belongs to the family of multivariate Pearson-type distributions, and $\lim_{\mathrm  f_t\to\pm \infty}\mathrm  f_t^2 f(\mathrm f_t,\bm{\mathrm X}_{nT})=0$, then $\E[F_t|\bm X_{nT}=\bm{\mathrm X}_{nT}]$ is linear in $\bm{\mathrm X}_{nT}$. \citet[p.292]{quahsargent93} and \citet[Assumption R]{DGRqml} made similar assumptions.

\begin{ass}[\textsc{Tails}]\label{ass:tails}
$\,$
\begin{compactenum}
 	\item [(a)] For all $t\in\mathbb Z$, all $j=1,\ldots, r$, and all $s>0$, $\text {\upshape P}(|v_{jt}|\ge s)\le \exp\l\{-K_v s^{\delta_v}\r\}$ for some finite positive reals $K_v$ and $\delta_v\le 2$ independent of $t$ and $j$.

	\item [(b)] For all $t\in\mathbb Z$, all $n\in\mathbb N$, all $s>0$, 
	\beq
	\mathrm P\l(\l\Vert\frac 1{\sqrt n} \sum_{i=1}^n \frac{\bm\lambda_i \xi_{it}}{\sigma_i^2}\r\Vert \ge s\r)\le 
	r 
	\exp\Big\{- \kappa_{1} s^2 \Big\}+
	rn \exp\bigg\{-\kappa_{2} \Big(s \sqrt n\Big)^{\alpha}\bigg\}
	,\nn
	\eeq
	for some finite positive reals $\kappa_1$, $\kappa_2$, and $\alpha\le 2$ independent of $t$ and $n$.
	\end{compactenum}
\end{ass}

In part (a), we assume an exponential-type tail inequality for the common shocks, which implies that for all $s>0$, $t\in\mathbb Z$, and $j=1,\ldots,r$,
$\mathrm P(|F_{jt}|>s)\le \exp\big\{-K_F s^{\delta_v}\big\}$, for some finite positive reals $K_F$ and $\delta_v\le 2$ independent of $t$ and $j$ (see Lemma \ref{lem:tail} and \citealp[Corollary 4]{maleki}). The same applies to part (b), which, by setting $n=1$, implies that for all $s>0$, $t\in\mathbb Z$, and $i\in\mathbb N$, 
$\mathrm P(\vert \xi_{it}\vert \ge s)\le \exp\big\{-K_\xi s^{\delta_\xi}\big\}$, for some finite positive real $K_\xi$ and $\delta_\xi \le 2$ independent of $t$ and $j$. This is because 
$\Vert \bm\lambda_i\Vert\le M_\lambda$ and $\sigma_i^2\ge C_\xi^{-1}$, for all $i\in\mathbb N$, by Assumptions \ref{ass:common}(a) and \ref{ass:idio}(a), respectively.  \citet[Assumption 2c]{FLM13} also assume the factors and idiosyncratic components to belong to a distribution having exponentially decaying tails. 

Depending on the values of $\delta_v$ and $\delta_\xi$, we are able to consider not only distributions with sub-Gaussian ($\delta_v,\delta_\xi=2$) or sub-exponential tails ($\delta_v,\delta_\xi=1$), which include the Laplace and the Generalized Error distribution \citep[Chapter 2]{vershynin18}, but also distributions with sub-Weibull tails ($\delta_v,\delta_\xi<1$), which can mimic a heavy tail behavior even if all moments exist \citep{MN98,KC18,vladimirova2020sub}.  This is clarified in Remark \ref{rem:weibull} below.

In general, part (b) is a Bernstein-type inequality implying that the weighted sums of the idiosyncratic components have Gaussian tails, as $n\to\infty$. This is a high-level requirement and, as for Assumption \ref{ass:idio}(e) in order to properly derive it, we should introduce some notion of dependence for random fields. Here, instead, we just notice that in the simplest case of cross-sectionally independent idiosyncratic components, this condition would be a direct consequence of \citet[Corollary 4]{maleki} (see also \citealp[Corollary 1]{vladimirova2020sub} for a similar result, and \citealp[Theorems 2.8.2 and 2.6.3, for the cases $\delta_\xi=1$ and $\delta_\xi=2$, respectively]{vershynin18}).

Finally, as a consequence of parts (a) and (b), jointly with Assumptions \ref{ass:common}(f), \ref{ass:common}(h), \ref{ass:idio}(c), and \ref{ass:ind},
we can show that not only $\{\mbf F_t\xi_{it}\}$ is a strongly mixing process but its components have also exponentially decaying tails. Therefore, for all $T\in\mathbb N$ and all $i\in\mathbb N$,  the following Bernstein-type inequality holds:
\beq
\mathrm P\l(\l\Vert\frac 1 {\sqrt T}\sum_{t=1}^T \mbf F_{t}\xi_{it} \r\Vert\ge s \r)\le 
r \exp\Big\{- \kappa_{3} s^2\Big\}+
rT \exp\bigg\{- \kappa_{4}\Big(s\sqrt  T\Big)^{\beta}\bigg\}
,\nn
\eeq
for some finite positive reals $\kappa_{3}$, $\kappa_{4}$, and $\beta<1$ independent of $i$ and $T$ (\citealp[Theorem 1]{MPR11}, \citealp[Theorem 1.4, p.31]{bosq12}, and Lemma \ref{lem:tail}).

\begin{rem}\label{rem:weibull}
\upshape{
Following \citet{KC18} and \citet{kuchibhotla2021}, we say that a given random variable $y$ is sub-Weibull with exponent $\delta$ if $\mathrm P(\vert y \vert \ge s)\le \exp\l\{-K_y s^{\delta}\r\}$ for any $s>0$ and some positive real $K_y$. This is equivalent to requiring the following Cram\'er type condition to hold: $\sup_{m\ge 1} r^{-1/\delta} (\E[\vert y \vert^m])^{1/m}\le M_y$ for some positive real $M_y$.\footnote{Another necessary and sufficient condition for a random variable $y$ to be sub-Weibull is given by the following condition on its Orlicz norm: $\Vert y\Vert_{\psi_{\delta}}=\inf \{\eta>0 : \E [\exp((\vert y \vert/\eta)^{\delta})] \le 2\}\le M_y^\prime$ for some finite positive real $M_y^\prime$.
} This shows that although all moments exist, they can be quite large for small values of $\delta$. For example, the $m$th moment of a Weibull is given by $\E[y^m]=\Gamma(1+\frac{m}{\delta})$, which is rapidly increasing as $\delta$ decreases
(see, e.g., \citealp{lehman1963}, for tabulated values of the first four moments as functions of $\delta$). 
}
\end{rem}

\subsection{Identification conditions}
To identify the DFM, we need to address four issues: first, we need to identify the number of factors; second, we need to identify the true parameters of the model; third, we need to ensure that the linear system is identified; and fourth, we need to guarantee the existence of the maxima of the log-likelihood.

\subsubsection{Identification of the number of factors and of the common component}
Starting with the number of factors $r$, let the covariance matrix of $\{\bm \chi_{nt}\}$ be $\bm\Gamma_n^{\chi}=\bm\Lambda_n\bm\Gamma^F\bm\Lambda_n^\prime$ and denote as $\mu_{jn}^\chi$ the $j$-th largest eigenvalue of $\bm\Gamma_n^{\chi}$, then Assumptions \ref{ass:common}(a) and \ref{ass:common}(b) imply that, for any $j=1,\ldots, r$,
\beq\label{eq:diveval}
\underline C_j\le \lim\inf_{n\to\infty} \frac{ \mu_{jn}^\chi}n\le\lim\sup_{n\to\infty} \frac{\mu_{jn}^\chi}n\le \overline C_j, 
\eeq
for some finite positive reals $\underline C_j$ and $\overline C_j$.
Furthermore, from Assumption \ref{ass:idio}(b) it follows that the largest eigenvalue of the idiosyncratic covariance matrix, $\bm\Gamma_n^\xi$, denoted as $\mu_{1n}^\xi$, is such that
\beq\label{eq:diveval_xi}
\mu_{1n}^\xi =\Vert\bm\Gamma_n^\xi\Vert \le M_\xi.
\eeq
From conditions \eqref{eq:diveval} and \eqref{eq:diveval_xi} and Weyl's 
inequality, it follows that the $r$ largest eigenvalues of the covariance matrix of $\{\mbf x_{nt}\}$ diverge linearly in $n$, whereas all remaining eigenvalues stay bounded for all $n\in\mathbb N$. Lemma \ref{lem:Gxi} proves this result, condition \eqref{eq:diveval}, and condition \eqref{eq:diveval_xi}---these conditions are directly imposed by \citet[Assumptions A1 and A2, respectively]{DGRqml}.


The asymptotic behavior of the eigenvalues of the covariance matrix allows for the identification of the common and idiosyncratic component when $n\to\infty$ and it is the basis for all  existing methods for determining $r$ (see, e.g., \citealp{baing02}).

\subsubsection{Identification of the loadings, the factors, and the VAR parameters}
All the structures equivalent to \eqref{eq:SDFM1R}-\eqref{eq:SDFM2R} can be obtained through an $r\times r$ invertible matrix $\mbf R$, as follows: 
\[
\mbf F_t^o=\mbf R^{-1}\mbf F_t,\;\; \bm{\Lambda}_n^{o}=\bm{\Lambda}_n\mbf R, \;\; \mbf A^o(L)=\mbf R^{-1}\mbf A(L)\mbf R, \;\;
\bm\Gamma^{vo} = \mbf R^{-1}\bm\Gamma^v,\;\;
 \bm\Gamma^{\xi o}_n=\bm\Gamma^{\xi}_n.
\]
Under such relationships, using only first- and second-moment information we cannot distinguish the model specified by $\bm{\Lambda}_n^{o}$, $\mbf A^o(L)$, $\bm\Gamma^{vo}$, and $\bm\Gamma^{\xi o}$, from the one given by $\bm{\Lambda}_n$, $\mbf A(L)$, $\bm\Gamma^v$, and $\bm\Gamma^{\xi}$. Once the loadings and the factors are identified, then $\mbf A(L)$ and $\bm\Gamma^v$ are also identified; $\bm\Gamma_n^\xi$ is always identified.

To identify the model, we need enough a priori structure to preclude any but the trivial transformation $\mbf R=\mbf I_r$. This can be achieved by imposing additional $r^2$ identifying constraints.  Let $\mbf M_n^\chi$ be the $r\times r$ diagonal matrix with as elements the eigenvalues $\mu_{jn}^\chi$, $j=1,\ldots, r$, of the covariance matrix of the common component, $\bm\Gamma_n^\chi$, sorted in descending order. Let  $\mbf V_n^\chi$ be the $n\times r$ matrix with the corresponding normalized eigenvectors as columns. Then, we assume the following identifying constraints.

\begin{ass}[\textsc{Identification}] \label{ass:ident}
$\,$
\begin{compactenum}[(a)]
	\item The eigenvalues of $\bm \Sigma_\Lambda\bm\Gamma^F$ are distinct.
	\item $\bm\Sigma_\Lambda$ is diagonal and $\bm\Gamma^F=\mbf I_r$. 
	\item For all $j=1,\ldots, r$, $[\bm \Lambda_n]_{1j}\ge 0$.
\end{compactenum}
\end{ass}

Part (a) is standard. Since the eigenvalues of $\bm\Sigma_\Lambda\bm\Gamma^F$ are equal to the $r$ non-zero eigenvalues of $\lim_{n\to\infty}n^{-1}\bm\Gamma_n^\chi$, given by $\lim_{n\to\infty}n^{-1}\mbf M_n^\chi$, it implies that in \eqref{eq:diveval} we have $\overline C_{j}<\underline C_{j-1}$ for any $j=2,\ldots, r$, and, thus, it avoids the uninteresting difficulties related with asymptotically multiple eigenvalues, which would require more restrictions to identify the space spanned by the columns of $\bm\Lambda_n$. 

Part (b) is similar to what is usually imposed in PC estimation (see Remark \ref{rem:identificL} below). It implies that $\bm\Gamma_n^\chi=\bm\Lambda_n\bm\Lambda_n^\prime=\mbf V_n^\chi\mbf M_n^\chi\mbf V_n^{\chi\prime}$.  Hence, the $r$ non-zero eigenvalues of $\lim_{n\to\infty}n^{-1}\bm\Gamma_n^\chi$, given by $\lim_{n\to\infty} n^{-1}\mbf M_n^\chi$, coincide with the diagonal entries of $\bm\Sigma_\Lambda$, which are then distinct because of part (a). Since part (b) concerns second moments and sums of squares, it allows us to identify $\bm\Lambda_n$ only up to a sign. 
 To achieve global identification we must fix also the column sign of $\bm\Lambda_n$ which is done through part (c) (\citealp[Remark 1]{baing13}). 
 
Summing up, in part (b), we are imposing $r^2$ restrictions: $r(r + 1)/2$ by requiring orthonormality of the factors, and $r(r-1)/2$ by requiring that $\bm\Sigma_\Lambda$ is diagonal. 
Consistently with the fact that in the typical empirical applications the focus is on the common component only, these assumed identification conditions do not provide economic meaning to the factors; in this sense ours is an exploratory rather than confirmatory factor analysis. Moreover, and most importantly, the restrictions are imposed only in the limit $n,T\to\infty$. Thus, in our setting, the model is only asymptotically identified. This is enough to derive our asymptotic theory.



\begin{rem}\label{rem:identificL}
\upshape{Typically, in PC estimation it is required that: 
(i) for all $n\in\mathbb N$, $n^{-1}\bm\Lambda_n^\prime\bm\Lambda_n$ is a diagonal matrix with finite distinct elements, and (ii) for all $T\in\mathbb N$, $T^{-1}\sum_{t=1}^T \mbf F_t\mbf F_t^\prime=\mbf I_r$ (\citealp{FGLR09}, \citealp{DGRfilter,DGRqml}, and \citealp[]{baing13}). In classical QML estimation constraint (ii) is the same as in PC estimation, while (i) is replaced with the requirement that for all $n\in\mathbb N$, $n^{-1}\bm\Lambda_n^\prime(\bm\Sigma_n^\xi)^{-1}\bm\Lambda_n$ is a diagonal matrix with finite distinct elements (\citealp[constraint IC3 therein]{baili12,baili16}). Differently from the present setting these constraints are assumed to hold for any given $n,T\in\mathbb N$.}
\end{rem}

\begin{rem}\upshape{
By letting $\bm{\mathcal S}$ be a diagonal $r\times r$ matrix with entries $\mathbb I([\mbf V_n^\chi]_{1j}\ge 0)-\mathbb I([\mbf V_{n}^\chi]_{1j}<0)$, $j=1,\ldots,r$, we can show that, as $n\to\infty$, $\bm\lambda_i^\prime$ coincides with $\mbf v_i^{\chi\prime}\bm{\mathcal S}(\mbf M_n^{\chi})^{1/2}$, where $\mbf v_i^{\chi\prime}$ is the $i$th row of $\mbf V_n^{\chi}$ (see Lemma \ref{lem:lambda}).
It also follows that, by linear projection of $\bm\chi_{nt}$ onto $\bm\Lambda_n$ at each given $t$, the true factors $\mbf F_t$ are identified, as $n\to\infty$, as the first $r$ normalized PCs of the common component given by $(\mbf M_n^{\chi})^{-1/2}\bm{\mathcal S}\mbf V_n^{\chi\prime} \bm\chi_{nt}$, which are clearly orthonormal, indeed, $\E[(\mbf M_n^{\chi})^{-1/2}\bm{\mathcal S}\mbf V_n^{\chi\prime} \bm\chi_{nt}\bm\chi_{nt}^\prime\mbf V_n^{\chi} \bm{\mathcal S}(\mbf M_n^{\chi})^{-1/2}]=\mbf I_r=\bm\Gamma^F$ as requested. 
}
\end{rem}


\subsubsection{Identification of the linear system}\label{sec:AD}
From our Assumptions \ref{ass:common}(a), \ref{ass:common}(d), and \ref{ass:common}(e),  we can show that the state space formulation \eqref{eq:SDFM1R}-\eqref{eq:SDFM2R} is minimal and stable, for all $n>N$, where $N= N_0$ in Assumption \ref{ass:common}(a). It is convenient to consider here the factorization $\bm \Gamma^{v}=\mbf H\mbf H^\prime$ for some $r\times r$ matrix $\mbf H$ having full rank. Consider again for simplicity the case $p_F=1$ with $\mbf A \equiv \mbf A_1$. Then, stability, which requires $\vert \nu^{(1)}(\mbf A)\vert <1$, is a direct consequence of stationarity in Assumption \ref{ass:common}(d). While minimality holds because the couple $(\mbf A, \mbf H)$ is controllable due to Assumption \ref{ass:common}(e) and the couple  $(\mbf A, \bm\Lambda_n)$ is observable for all $n>N$ because Assumption \ref{ass:ident}(b) implies that $\text{rk}(\bm\Lambda_n)=\text{rk}(\mbf V_n^\chi)=r$, for all $n>N$.\footnote{A linear system with $p_F=1$, is controllable if and only if $\text{rk}[\mbf H\; (\mbf A\mbf H)\cdots (\mbf A^{(r-1)}\mbf H)]=r$ and it is observable if and only if $\text{rk}[\bm\Lambda_n^\prime\; (\bm\Lambda_n\mbf A)^\prime\cdots (\bm\Lambda_n\mbf A^{r-1})^\prime ]=r$ (\citealp[Appendix C, pp. 341-342]{AM79}).}  This implies that, for all $n>N$, the linear system satisfies the mini-phase condition:
\[
\text{rk}\l(\ba{cc}
\mbf I_r-\mbf Az & \;-\mbf H\\
\bm\Lambda_n & \;\mbf 0_{n\times r}
\ea\r)=2r, \;\text{ for }\; |z|\ge1.
\]
In other words, the linear system \eqref{eq:SDFM1R}-\eqref{eq:SDFM2R} is the minimal state-space representation of a DFM 
having as McMillan degree the number of factors $r$  \citep[Section II]{andersondeistler08}. 

This result, together with Assumption \ref{ass:ident}, guarantees that the transfer function matrix $\bm W(z)=\bm\Lambda_n(\mbf I_r-\mbf Az)^{-1} \mbf H$ is identified for all $z\in\mathbb C$, with the exception of a zero-measure set (\citealp[see also][Section 4.3, for a proof]{LDA}, and \citealp{heatonsolo04}, for similar results). This implies generic identifiability of the linear system \eqref{eq:SDFM1R}-\eqref{eq:SDFM2R}.

\subsubsection{Existence of the EM and QML estimators}\label{enzuccio}
Finally, we consider the issue of identification of the maxima of the log-likelihood. For any given $i\in\mathbb N$, define 
\begin{align}
\mathcal O_{\lambda_i}&=\{\underline{\bm\lambda}_i\in\mathbb R^{r},\,  \underline{\bm\lambda}_i\in[-M_\lambda,M_\lambda]^r\}\; \text{ and }\; \mathcal O_{\sigma_i^2}=\{ \underline{\sigma}_i^2\in\mathbb R,\,  \underline{\sigma}_i^2\in[C_\xi^{-1},C_\xi]\},\nn
\end{align}
where $M_\lambda$ and $C_\xi$ are defined in Assumptions \ref{ass:common}(a) and \ref{ass:idio}(a), respectively. Likewise, define
\begin{align}
\mathcal O_{\mathcal A}&=\{\text{vec}(\underline{\mbf A})\in \mathbb R^{r^2},\, \text{vec}(\underline{\mbf A})\in [-M_A,M_A]^{r^2p_F}\},\nn\\
\mathcal O_{\Gamma^v}&=\{\text{vech}(\underline{\bm \Gamma}^v)\in \mathbb R^{r(r+1)/2},\, \nu^{(j)}(\underline{\bm \Gamma}^v)\in[M_v^{-1},M_v], j=1,\ldots,r\},\nn
\end{align}
where $M_A$ and $M_v$ are defined in Assumptions \ref{ass:common}(d) and \ref{ass:common}(e), respectively. Moreover, for any given $n\in\mathbb N$, let also
\begin{align}
\mathcal E_{\Lambda_n}&=\{\text {vec}(\underline {\bm\Lambda}_n)^\prime \in\mathbb R^{nr},\, \underline C_r\le \nu^{(r)}(n^{-1}\underline {\bm\Lambda}_n^\prime \underline {\bm\Lambda}_n)\le\nu^{(1)}(n^{-1}\underline {\bm\Lambda}_n^\prime \underline {\bm\Lambda}_n)\le \overline C_1,  \underline {\bm\Lambda}_n^\prime\underline {\bm\Lambda}_n \text{diagonal}  \},\nn\\
\mathcal E_{\Gamma^\xi_n}&=\{\text{vech}(\underline{\bm\Gamma}_n^\xi)^\prime\in\mathbb R^{n(n+1)/2},\,L_\xi \le \nu^{(n)}(\underline{\bm\Gamma}_n^\xi)\le \nu^{(1)}(\underline{\bm\Gamma}_n^\xi)\le M_\xi \},\nn
\end{align}
where $\underline C_r$ and $\overline C_1$ are defined in \eqref{eq:diveval}, while  $M_\xi$ and $L_\xi$ are defined in Assumptions \ref{ass:idio}(b) and \ref{ass:idio}(f), respectively. 
Then, the search for the maximum of the expected log-likelihood in the EM algorithm and the one of the log-likelihood takes place on the set $\mathcal O_n=\{\mathcal O_{\lambda_i}^n\cap \mathcal E_{\Lambda_n}\}\times\{\mathcal O_{\sigma_i^2}^n\cap\mathcal E_{\Gamma^\xi_n}\} \times \mathcal O_{\mathcal A}\times \mathcal O_{\Gamma^v}$, which has dimension $Q_n=n(r+1)+r^2p_F+r(r+1)/2$ growing with $n$. 

Because of Assumptions \ref{ass:common}(a) and \ref{ass:idio}(a), the rows of the loadings matrix $\bm\lambda_i$ and the idiosyncratic variances $\sigma_i^2$ belong to $\mathcal O_{\lambda_i}\times\mathcal O_{\sigma_i^2} \subset \mathbb R^{r+1}$, which is a compact set for any given $i\in\mathbb N$. Similarly, because of Assumptions \ref{ass:common}(d) and \ref{ass:common}(e), the entries of $\mbf A_k$, $k=1,\ldots, p_F$, and $\mbf H$ belong to $\mathcal O_{\mathcal A}\times \mathcal O_{\Gamma^v}\subset \mathbb R^{r^2p_F+r(r+1)/2}$, which is also a compact set. These properties are crucial as they ensure the existence of a maximum which is a solution of the EM algorithm.
Indeed, for any iteration $k\ge 0$ the expected log-likelihood $\mathcal Q(\underline{\bm\varphi}_n,\wh{\bm\varphi}_{n}^{(k)})$ in \eqref{eq:LLbayes_exp_bis}, which we maximize in the M-step, is made of two terms: a term  in \eqref{eq:LLbayes_exp_bis2} associated to the VAR for the factors \eqref{eq:SDFM2R}, which is defined on the finite-dimensional set $\mathcal O_{\mathcal A}\times \mathcal O_{\Gamma^v}$, and a term \eqref{eq:LLbayes_exp_bis1} associated to the factor equation \eqref{eq:SDFM1R} and thus defined over a set $\mathcal O_{\lambda_i}^n\times\mathcal O_{\sigma_i^2}^n$ of dimension growing with $n$.
Now, the former term poses no problem because we can use compactness of $\mathcal O_{\mathcal A}\times \mathcal O_{\Gamma^v}$ to show that the maxima $\wh{\bm\theta}^{(k+1)}$ exist. Regarding the latter term, in order to prove the existence of $\wh{\bm\phi}_n^{(k+1)}$ we can still use the same compactness argument by noticing that maximizing this term amounts to separately maximizing of $n$ terms, each depending only of $\underline{\bm\lambda}_i$ and $\underline{\sigma}_i^2$ for given $i\in\mathbb N$, and thus defined on the finite-dimensional compact set $\mathcal O_{\lambda_i}\times\mathcal O_{\sigma_i^2}$. It is also easy to show that the elements of $\wh{\bm\varphi}_n^{(k+1)}$ are unique because the M-step gives a closed form solution for each of them.
This reasoning guarantees the existence and uniqueness of the EM estimators (see Lemma \ref{lem:existEM}).

Let us turn to the QML estimator that maximizes the full log-likelihood \eqref{eq:LL0true}. The existence of the QML estimator $\wh{\bm\theta}^*$ poses no problem because being a finite-dimensional vector, we can use the compactness argument. Direct proof of the existence of the QML estimator $\wh{\bm\phi}_n^*$ is instead more challenging, as in this case, we cannot rely on compactness because the full log-likelihood is defined on a set of increasing dimension $n$. 

Nevertheless, we give an indirect proof of the existence of $\wh{\bm\phi}_n^*$ by noticing that, under our identification Assumptions \ref{ass:ident}(a)-\ref{ass:ident}(c), as $n,T\to\infty$, the elements of $\wh{\bm\phi}_n^*$ are asymptotically equivalent to the unfeasible Ordinary Least Squares estimators we would obtain if the factors were observed (see Lemma \ref{lem:starols}). Therefore, this argument ensures, at least asymptotically, both the existence and also the uniqueness of the QML estimators.  We also refer to the next section for more details. 

\begin{rem}\label{rem:miao}
\upshape{Typically, this literature assumes that the QML estimators $\wh{\sigma}_i^{2*}$ belong to a compact set (\citealp[Assumption D]{baili16}). In our set-up, this assumption is implied  (\citealp[Theorem 3.1]{GGM21}, and \citealp[Theorem 1]{mao2024statistical}), and not needed for proving our results \citep{MBPCAQML}. 
}
\end{rem}

\section{Asymptotic properties}\label{sec:asymp}
This section presents the asymptotic properties of the EM estimator of the parameters---i.e., of the factor loadings, idiosyncratic variances, VAR coefficients, and the covariance matrix of the VAR innovations---and of the Kalman smoother estimator of the factors. We assume that $r$, the number of common factors, is known; without loss of generality, we fix the VAR order in \eqref{eq:SDFM2R} to $p_F=1$, and we let $\mbf A\equiv \mbf A_1$ so that $\mbf A(L)\equiv\mbf A L$.  We briefly discuss the estimation of $r$ and $p_F$ at the end of the section.

\subsection{Consistency under basic assumptions}\label{sec:cons00}
We start by proving the consistency of the EM algorithm and the Kalman smoother under the most general case where we neither impose linearity of the conditional mean (Assumption \ref{ass:linear}) nor exponentially decaying tails (Assumption \ref{ass:tails}).

\begin{prop}\label{prop:plainvanilla}
Consider the EM estimators of the parameters $\wh{\bm\Lambda}_n=(\wh{\bm \lambda}_1\cdots \wh{\bm \lambda}_n)^\prime$ with $\wh{\bm \lambda}_i \equiv \wh{\bm \lambda}_i^{(k^*+1)}$, $\wh{\sigma}_i^2\equiv \wh{\sigma}_i^{2(k^*+1)}$ , $i=1,\ldots, n$, 
$\wh{\mbf A}\equiv\wh{\mbf A}^{(k^*+1)}$, and $\wh{\bm\Gamma}^v\equiv\wh{\bm\Gamma}^{v(k^*+1)}$, and the Kalman smoother estimator of the factors, $\wh{\mbf F}_t\equiv \mbf F_{t|T}^{(k^*+1)}$, $t=1,\ldots, T$, $k^*\ge 0$. Then, under Assumptions \ref{ass:common}, \ref{ass:idio}, \ref{ass:ind}, and \ref{ass:ident}:
\begin{compactenum}
\item [(a)] for all $\epsilon>0$, there exist a positive real $\eta(\epsilon)$, and integers $n^*(\epsilon)$ and $T^*(\epsilon)$, all independent of $i$, such that, for all $n\ge n^*(\epsilon)$ and all $T\ge T^*(\epsilon)$,
	\al{
	(a.1)&\quad\mathrm P\l(\min( n,\sqrt {T})\, \Vert\wh{\bm \lambda}_{i}-\bm\lambda_i\Vert\ge \eta(\epsilon)\r) < \epsilon,\;\text{ for any given $i=1,\ldots, n$,}\nn\\
	(a.2)&\quad\mathrm P\l(\min( n,\sqrt {T})\, n^{-1/2}\Vert\wh{\bm \Lambda}_{n}-\bm\Lambda_n\Vert\ge \eta(\epsilon)\r) < \epsilon,\nn\\
	(a.3)&\quad \mathrm P\l(\min( n,\sqrt {T})\, \vert\wh{\sigma}_{i}^2-\sigma_i^2\vert\ge \eta(\epsilon)\r) < \epsilon,\;\text{ for any given $i=1,\ldots, n$,}\nn\\
	(a.4)&\quad\mathrm P\l(\min( n,\sqrt {T})\, \Vert\wh{\mbf A}-\mbf A\Vert\ge \eta(\epsilon)\r) < \epsilon,\nn\\
	(a.5)&\quad \mathrm P\l(\min( n,\sqrt {T})\, \Vert\wh{\bm \Gamma}^v-\bm\Gamma^v\Vert\ge \eta(\epsilon)\r) < \epsilon;\nn
	}
\item[(b)] for all $\epsilon>0$, there exist a positive real $\eta(\epsilon)$, and integers $n^{**}(\epsilon)$ and $T^{**}(\epsilon)$, all independent of $t$, such that, for all $n\ge n^{**}(\epsilon)$ and all $T\ge T^{**}(\epsilon)$,
	\[
	\mathrm P\l(\min( \sqrt n,\sqrt {T})\, \Vert\wh{\mbf F}_{t}-\mbf F_t\Vert\ge \eta(\epsilon)\r) < \epsilon,
	\]
	for any given $t=1,\ldots, T$.
\end{compactenum}
\end{prop}

The consistency rate of the estimated parameter is the same as the PC estimator \citep[Theorem 2]{Bai03}. This is because by initializing the algorithm with the consistent PC estimator, we can consider the EM estimator a ``one-step'' estimator (\citealp[Theorem 4.3]{LC06}). As for the the estimated factors, they converge at a slower rate than the PC estimator, due to the $\sqrt T$ term \citep[Theorem 1]{Bai03}. 

\begin{rem}
\upshape{
The estimation error for the factors and the one for the loadings both depend on the estimation error of the diagonal idiosyncratic covariance matrix $\Vert \wh{\bm\Sigma}_n^\xi-\bm\Sigma_n^\xi\Vert =\max_{i=1,\ldots, n}\vert \wh{\sigma}_i^{2}-\sigma_i^2\vert$. 
The latter produces a term in the estimation error of the factors which is $O_p(T^{-1/2})$ (term $D.2$ in the proof of Lemma \ref{lem:aprile24}) and a term in the estimation error of the loadings which is also $O_p(T^{-1/2})$ (term $III_d$ in the proof  of Proposition \ref{prop:plainvanilla}). These two terms have non-standard asymptotic distributions and are  non-negligible. Thus, we cannot prove asymptotic normality of our estimators without additional assumptions.
}
\end{rem}

\begin{rem}
\upshape{The results of Proposition \ref{prop:plainvanilla}  apply to the case in which the observed data, $\mbf y_{nt}=(y_{1t}\cdots y_{nt})^\prime$ , are such that $y_{it}=y_{i0}+\mu_i t + \bm\lambda_i^\prime \mbf G_t + z_{it}$. In this case, if we write $\mbf F_t =\Delta \mbf G_t$ and $\xi_{it}=\Delta z_{it}$, then Proposition \ref{prop:plainvanilla} holds for $x_{it}=\Delta y_{it}$. This strategy always works for the loadings \citep[see][for details]{BLL2} so part (a) still holds, but estimation of the factors must be modified if $\mbf G_t$ is a cointegrated vector. First, we must model $\mbf G_t$ as a VAR in levels. Second, we estimate $\mbf G_t$ running the Kalman smoother in levels, and whenever $z_{it}\sim I(1)$ for some $i$, we add a latent state. If  $z_{it}\sim I(0)$ for all $i$, then part (b) stands; if $z_{it}\sim I(1)$ for some $i$, we conjecture that part (b) would remain unchanged. We leave the derivation of the asymptotic properties of the estimated factors in this last case for further research.
}
\end{rem}



\subsection{Consistency and asymptotic normality}
If we also assume that Assumptions \ref{ass:linear} and \ref{ass:tails} hold, we can refine the previous result because we can now guarantee that the EM algorithm converges to the QML estimator. 

\begin{prop}[\textsc{Loadings}]\label{prop:load} 
Consider the EM estimators of the loadings 
$\wh{\bm\Lambda}_n=(\wh{\bm \lambda}_1\cdots \wh{\bm \lambda}_n)^\prime$ with $\wh{\bm \lambda}_i \equiv \wh{\bm \lambda}_i^{(k^*+1)}$, $k^*\ge 0$. 
Then, under Assumptions \ref{ass:common}, \ref{ass:idio}, \ref{ass:ind}, \ref{ass:linear}, \ref{ass:tails}, and \ref{ass:ident}:
\begin{compactenum}	
	\item [(a)] for all $\epsilon>0$, there exist a positive real $\eta(\epsilon)$, and integers $n^*(\epsilon)$ and $T^*(\epsilon)$, all independent of $i$, such that, for all $n\ge n^*(\epsilon)$ and all $T\ge T^*(\epsilon)$, and some $0<\delta_v\le 2$,
	\al{
	(a.1)&\quad\mathrm P\l(\min( n/ {\log^{2/\delta_v} T},\sqrt {T})\, \Vert\wh{\bm \lambda}_{i}-\bm\lambda_i\Vert\ge \eta(\epsilon)\r) < \epsilon,\;\text{ for any given $i=1,\ldots, n$,}\nn\\
	(a.2)&\quad\mathrm P\l(\min( n/ {\log^{2/\delta_v} T},\sqrt {T})\, n^{-1/2}\Vert\wh{\bm \Lambda}_{n}-\bm\Lambda_n\Vert\ge \eta(\epsilon)\r) < \epsilon;\nn
	}
	\item [(b)] for any given $i=1,\ldots, n$, as $n,T\to\infty$, if $n^{-1}\sqrt {T}\log^{2/\delta_v}T\to 0$,
	\al{
	\sqrt T (\wh{\bm\lambda}_{i}-\bm\lambda_i)\stackrel{d}{\to}\mathcal N(\mbf 0_r,\bm{\mathcal V}_i),\nn
	}
	 where  $$\bm{\mathcal V}_{i}=(\bm\Gamma^F)^{-1}\l(\lim_{T\to\infty} \frac 1T\sum_{t=1}^T\sum_{s=1}^T \E[\xi_{it}\xi_{is}]\E[\mbf F_t\mbf F_s^\prime]\r)(\bm\Gamma^F)^{-1},
	 $$ with $\bm\Gamma^F=\lim_{T\to\infty}T^{-1} \sum_{t=1}^T \mbf F_t\mbf F_t^\prime=\mbf I_r$, because of Assumption \ref{ass:ident}(b);
	 \item [(c)] for any given $i=1,\ldots, n$, if $\E[\xi_{it}\xi_{is}]=0$ for all $t,s=1,\ldots, T$ with $t\ne s$, then, $\bm{\mathcal V}_{i}=\sigma_i^2(\bm\Gamma^F)^{-1} $, with $\bm\Gamma^F=\mbf I_r$, because of Assumption \ref{ass:ident}(b).
\end{compactenum}
\end{prop}

The rate of consistency of the estimated loadings, $\min(n/{\log^{2/\delta_v} T},\sqrt T)$, given in Proposition \ref{prop:load}, is new to the EM literature. This rate is the same, up to a logarithmic factor, as the one of the PC estimator (\citealp[Theorem 2]{Bai03}), which, in turn, is equivalent to the unfeasible Ordinary Least Squares (OLS) we would obtain if the factors were observed.  The EM estimator is also asymptotically equivalent to the QML estimator considered by \citet[Theorem 1]{baili16} when imposing no autocorrelation for the factors.
For an explanation of the logarithmic term we refer to Remark \ref{rem:logs} below.  Efficiency is discussed in Section \ref{sec:PCeff}.

It is important to stress that Proposition \ref{prop:load} requires not only $T\to\infty$, as in classical QML estimation theory, but also $n\to\infty$ otherwise no consistency can be proved. In particular, as $n\to\infty$ the factors can be treated as observed, therefore, there is no more an issue of missing information and the QML estimator of the loadings must coincide with the unfeasible OLS. This is a manifestation of the blessing of dimensionality which is a fundamental feature of approximate factor models.

The proof of Proposition \ref{prop:load} is based on the following decomposition of the estimation error into four terms:
\begin{align}\label{eq:MLpropL}
\sqrt T (\wh{\bm\lambda}_i-\bm\lambda_i)&=
\underbrace{\sqrt T({\bm\lambda}_{i}^{\text{\tiny{OLS}}}-{\bm\lambda}_{i})}_{O_p(1)}+
\underbrace{\sqrt T(\wh{\bm\lambda}_{i}^{*}-{\bm\lambda}_{i}^{\text{\tiny OLS}})}_{O_p(n^{-1}\sqrt T   \log^{2/\delta_v}T)}+
\underbrace{\sqrt T(\wh{\bm\lambda}_i^{**}-\wh{\bm\lambda}_i^{*})}_{O_p(n^{-1}\sqrt T  \log^{2/\delta_v}T) }+
\underbrace{\sqrt T(\wh{\bm\lambda}_i-\wh{\bm\lambda}_i^{**})}_{o_p( n^{-1}\sqrt T \log^{2/\delta_v}T)},
\end{align}
which shows that the EM estimator $\wh{\bm\lambda}_i$ is asymptotically equivalent to the OLS estimator we would obtain had we observed the factors.

To prove our result, we first show that the QML estimator of the loadings, $\wh{\bm\lambda}_i^*$ is asymptotically equivalent to the PC estimator, which, in turn, is asymptotically equivalent to the unfeasible $\sqrt T$-consistent OLS estimator ${\bm\lambda}_{i}^{\text{\tiny OLS}}$ (second term on the rhs of \eqref{eq:MLpropL} proved in Lemma \ref{lem:starols}, see also \citealp[Theorem 3 and Corollary 1]{MBPCAQML}). In particular, both approximation errors are $o_p(T^{-1/2})$, whenever $n^{-1}\sqrt {T}\log^{2/\delta_v} T\to 0$. This result extends to the DFM the result by \citet[Theorem 1]{baili16} obtained for QML estimation of a static factor model 
i.e., when we replace the full-matrix ${\bm\Omega}_T^F(\underline{\bm{\mathcal A}}, \underline{\mbf H})$ with just $\mbf I_{rT}$ in the log-likelihood \eqref{eq:LL0true}. 
Notice that, differently from the proofs in \citet{baili16} and, as mentioned in Remark \ref{rem:miao}, this result does not depend on the QML estimator $\wh{\sigma}_i^{*2}$.
 
Next, we show that the EM estimator converges to a global maximum of the likelihood (third term on the rhs of \eqref{eq:MLpropL}). As we discussed in Section \ref{sec:contra}, we know that the sequence of estimators $\{\wh{\bm\lambda}_i^{(k)},\, k\ge 0\}$ converges to a local maximum of the likelihood, say $\wh{\bm\lambda}_i^{**}$, as $k\to\infty$ (see Lemma \ref{lem:localmax}). However, in general, the likelihood might have many maxima due to the identification indeterminacy of the loadings. Nevertheless, once we make the identifying Assumptions \ref{ass:ident}(b) and \ref{ass:ident}(c), there is only a unique maximum, $\wh{\bm\lambda}_i^*$, which is asymptotically equivalent to the unique OLS estimator, whenever $n^{-1}\sqrt {T}\log^{2/\delta_v} T\to 0$ (see Lemma \ref{lem:localglobal} and \citealp[Section 4]{ruud91}, for a similar result in the case of one-to-one mapping from the factors to the data, corresponding to the case of no idiosyncratic component). This is also clear from the asymptotic expansions of $\wh{\bm\lambda}_i^*-\bm\lambda_i$ obtained by \citet{baili12,baili16} for the static model (i.e., the factors have no dynamics) and using identification schemes different than the one we use here.

Third, we show that asymptotically the error coming from running the EM a finite number of times vanishes (fourth term on the rhs of \eqref{eq:MLpropL}). Indeed, due to the finite number of iterations, $k^*$, the EM algorithm delivers an estimator $\wh{\bm\lambda}_i\equiv\wh{\bm\lambda}_i^{(k^*+1)}$, which is just an approximation of the local maximum $\wh{\bm\lambda}_i^{**}$ that we would attained after an infinite number of iterations. We show that the error entailed by such approximation depends on the ratio of the Hessians of the complete and incomplete log-likelihoods, i.e., on how much information is missing because the factors are not observed (\citealp{MR94}, \citealp[Chapter 3.9]{MLT07}, and \citealp[Chapter 8]{sundberg2019statistical}). In this case the approximation error is $o_p(T^{-1/2})$, provided $n^{-1}\sqrt {T}\log^{2/\delta_v}T\to 0$ (see Lemma \ref{lem:convEM1}). This last result is a refinement of the results by \citet[Theorem 2]{BWB17} on the convergence of the EM algorithm. We refer to Section \ref{sec:contra} below for more details.

Last, ${\bm\lambda}_{i}^{\text{\tiny{OLS}}}$ is a $\sqrt T$-consistent estimator of the factor loadings $\bm\lambda_{i}$ (first term on the rhs of \eqref{eq:MLpropL}).


\begin{prop}[\textsc{Factors}]\label{prop:factors}
Consider the Kalman smoother estimator of the factors 
$\wh{\bm{\mathcal F}}_T=(\wh{\mbf F}_1\cdots \wh{\mbf F}_T)^\prime$, with $\wh{\mbf F}_t \equiv {\mbf F}_{t|T}^{(k^*+1)}$, $t=1,\ldots, T$, $k^*\ge 0$.
Then, under Assumptions \ref{ass:common}, \ref{ass:idio}, \ref{ass:ind}, \ref{ass:linear}, \ref{ass:tails}, and \ref{ass:ident}:
\begin{compactenum}
	\item [(a)]  for all $\epsilon>0$, there exist a positive real $\eta(\epsilon)$, and integers $n^{**}(\epsilon)$ and $T^{**}(\epsilon)$, all independent of $t$, such that, for all $n\ge n^{**}(\epsilon)$ and all $T\ge T^{**}(\epsilon)$,
	\al{
	(a.1)&\quad\mathrm P\l(\min( \sqrt n, {T}/\sqrt{\log n})\, \Vert\wh{\mbf F}_{t}-\mbf F_t\Vert\ge \eta(\epsilon)\r) < \epsilon,\;\text{ for any given $t=1,\ldots, T$,}\nn\\
	(a.2)&\quad\mathrm P\l(\min( \sqrt n, {T}/\sqrt{\log n})\, T^{-1/2}\Vert\wh{\bm{\mathcal F}}_{T}-\bm{\mathcal F}_T\Vert\ge \eta(\epsilon)\r) < \epsilon;\nn
	}
	\item [(b)] as $n,T\to\infty$, if $T^{-1}\sqrt {n\log n}\to 0$, 
	\[
	\sqrt n (\wh{\mbf F}_{t}-\mbf F_t)\stackrel{d}{\to}\mathcal N(\mbf 0_r,\bm{\mathcal W}_t),
	\]
	for any given $t=1,\ldots, T$, where 
	$$
	\bm{\mathcal W}_t=(\bm\Sigma_{\Lambda\Sigma\Lambda})^{-1} \l(\lim_{n\to\infty} \frac 1 n\sum_{i=1}^n\sum_{j=1}^n \frac{\E[\xi_{it}\xi_{jt}] \bm\lambda_i\bm\lambda_j  }{\sigma_i^2\sigma_j^2} \r)(\bm\Sigma_{\Lambda\Sigma\Lambda})^{-1},
	$$ 
	with $\bm\Sigma_{\Lambda\Sigma\Lambda}=\lim_{n\to\infty} n^{-1}\sum_{i=1}^n \bm\lambda_i (\sigma_i^2)^{-1}\bm\lambda_i^\prime$;
	 \item [(c)] for any given $t=1,\ldots, T$, if $\E[\xi_{it}\xi_{jt}]=0$ for all $i,j=1,\ldots, n$ with $i\ne j$, then, $\bm{\mathcal W}_{t}=(\bm\Sigma_{\Lambda\Sigma\Lambda})^{-1}$.
\end{compactenum}
\end{prop}
 
The rate of consistency of the estimated factors given in Proposition \ref{prop:factors} is faster than the rate originally derived by \citet[Proposition 1]{DGRqml} for the same estimator---$\min(\sqrt n,T/\sqrt{\log n})$ vs. $\min(\sqrt n,T^{1/4}/\sqrt{\log n})$. Moreover, this consistency rate is the same (up to a logarithmic factor) as that of the PC estimator (\citealp[Theorem 1]{Bai03}). However, while the PC estimator is equivalent to the unfeasible Ordinary Least Squares (OLS) we would obtain if the loadings were observed, the Kalman smoother estimator we are considering is equivalent to the unfeasible Weighted Least Squares (WLS) we would obtain if the loadings were observed and we knew the idiosyncratic variances. As such, the Kalman smoother is also equivalent to the feasible WLS studied by \citet[Theorem 2]{baili16} and computed using the QML estimator of the loadings for a static factor model.
For an explanation of the logarithmic term we refer to Remark \ref{rem:logs} below.  
Efficiency is discussed in Section \ref{sec:PCeff}.

The proof of Proposition \ref{prop:factors} is based on the following decomposition of the estimation error into four terms
\begin{align}\label{eq:MLpropF}
\sqrt n (\wh{\mbf F}_t-\mbf F_t)&=
\underbrace{\sqrt n(\mbf F_{t}^{\text{\tiny WLS}}-{\mbf F}_{t})}_{O_p(1)}+
\underbrace{\sqrt n(\wh{\mbf F}_{t}^{\text{\tiny WLS}}-{\mbf F}_{t}^{\text{\tiny WLS}})}_{O_p(T^{-1}\sqrt {n\log n} )}+
\underbrace{\sqrt n(\wh{\mbf F}_{t|t}-\wh{\mbf F}_{t}^{\text{\tiny WLS}})}_{O_p(n^{-1/2})}+
\underbrace{\sqrt n(\wh{\mbf F}_{t}-\wh{\mbf F}_{t|t})}_{O_p(n^{-1/2} )}.
\end{align}
which shows that the Kalman smoother estimator $\wh{\mbf F}_t$ is asymptotically equivalent to the WLS estimator we would obtain had we observed the loadings and had we known the idiosyncratic variances. 

To prove our result, we first show that our estimator of the factors, which is obtained via the Kalman smoother computed using the EM estimators of the parameters, is asymptotically equivalent to the Kalman filter, $\wh{\mbf F}_{t|t}$, which in turn is asymptotically equivalent to the WLS estimator $\wh{\mbf F}_{t}^{\text{\tiny WLS}}$ (fourth and third term on the rhs of \eqref{eq:MLpropF}, respectively). Both approximation errors  are $O_p(n^{-1})$ (see Lemma \ref{lem:FFOnhat} and \citealp[Section 2.3]{PR22}, for the one-factor case). 

Then, we take into account the estimation error of the parameters given in Proposition \ref{prop:load}, which implies that the WLS estimator $\wh{\mbf F}_{t}^{\text{\tiny WLS}}$ converges to its unfeasible counterpart computed using the true value of the parameters, ${\mbf F}_{t}^{\text{\tiny WLS}}$, with a rate which is $o_p(n^{-1/2})$, provided $T^{-1}\sqrt {n\log n}\to 0$ (second term on the rhs of \eqref{eq:MLpropF}). This result refines the result of Proposition \ref{prop:plainvanilla} because, thanks to Assumption \ref{ass:tails}, we are now able to derive a tighter bound for $\Vert \wh{\bm\Sigma}_n^\xi-\bm\Sigma_n^\xi \Vert=\max_{i=1,\ldots, n}\vert \wh{\sigma}_i^2-\sigma_i^2\vert$, as shown in Proposition \ref{prop:altri} below.

Last, ${\mbf F}_{t}^{\text{\tiny WLS}}$  is a $\sqrt n$-consistent estimator of the realizations of the factors $\mbf F_t$ (first term on the rhs of \eqref{eq:MLpropF}).


\begin{prop}[\textsc{Common component}]\label{th:chi} 
Consider the EM plus the Kalman smoother estimator of the common component $\wh{\chi}_{it}  \equiv \wh{\bm\lambda}_i^{(k^*+1)\prime}{\mbf F}_{t|T}^{(k^*+1)}$, 
$i=1,\ldots, n$, $t=1,\ldots, T$, with $k^*\ge 0$.
Then, under Assumptions \ref{ass:common}, \ref{ass:idio}, \ref{ass:ind}, \ref{ass:linear}, \ref{ass:tails}, and \ref{ass:ident}:
\begin{compactenum}
	\item [(a)] for all $\epsilon>0$, there exist a positive real $\eta(\epsilon)$, and integers $n^{\circ}(\epsilon)$ and $T^{\circ}(\epsilon)$, all independent of $i$ and $t$, such that, for all $n\ge n^{\circ}(\epsilon)$ and $T\ge T^{\circ}(\epsilon)$,
	\[
	\mathrm P\l(\min( \sqrt n, \sqrt T)\, \vert\wh{\chi}_{it}-\chi_{it}\vert\ge \eta(\epsilon)\r) < \epsilon,
	\]	
	for any given $i=1,\ldots, n$, $t=1,\ldots, T$;
	\item [(b)]  as $n,T\to\infty$,
	\[
	(T^{-1}\mathcal C^\lambda_{it}+n^{-1}\mathcal C^F_{it})^{-1/2}(\wh{\chi}_{it}-\chi_{it})\stackrel{d}{\to}\mathcal N(0,1),
	\]
	for any given $i=1,\ldots,n$ and $t=1,\ldots, T$,	where $\mathcal C^\lambda_{it}=\mbf F_t^\prime\bm{\mathcal V}_i\mbf F_t$ and $\mathcal C^F_{it}=\bm\lambda_i^\prime\bm{\mathcal W}_t\bm\lambda_i$, with $\bm{\mathcal V}_i$ defined in Proposition \ref{prop:load}(b), and $\bm{\mathcal W}_t$ in Proposition \ref{prop:factors}(b).
\end{compactenum}
\end{prop}

Proposition \ref{th:chi} does not require a limit for $T/n$ or $n/T$, so it holds without any constraint between the rates of divergence of $n$ and $T$. That said, Proposition \ref{th:chi} has two special cases: (a) if $n/T\to 0$, then $\sqrt n(\wh{\chi}_{it}-\chi_{it})\stackrel{d}{\to}\mathcal N(0,\mathcal C^F_{it})$; and, (b) if $T/n\to 0$, then $\sqrt T(\wh{\chi}_{it}-\chi_{it})\stackrel{d}{\to}\mathcal N(0,\mathcal C^\lambda_{it})$. This is the same rate of consistency we obtain for the PC estimator of the common component (\citealp[Theorem 3]{Bai03}).

All other estimated parameters are also consistently estimated.
\begin{prop}[\textsc{Idiosyncratic variances and VAR parameters}]\label{prop:altri} 
Consider the EM estimators of the parameters $\wh{\bm\Sigma}_n^\xi=\text{\upshape diag}(\wh{\sigma}_1^2\cdots \wh{\sigma}_n^2)$,
with  $\wh{\sigma}_i^2\equiv \wh{\sigma}_i^{2(k^*+1)}$, $i=1,\ldots, n$, 
$\wh{\mbf A}\equiv\wh{\mbf A}^{(k^*+1)}$, and $\wh{\bm\Gamma}^v\equiv\wh{\bm\Gamma}^{v(k^*+1)}$, $k^*\ge 0$. 
Then, under Assumptions \ref{ass:common}, \ref{ass:idio}, \ref{ass:ind}, \ref{ass:linear}, \ref{ass:tails}, and \ref{ass:ident}:
\begin{compactenum}	
	\item [(a)] for all $\epsilon>0$, there exist a positive real $\eta(\epsilon)$, and integers $n^*(\epsilon)$ and $T^*(\epsilon)$, all independent of $i$, such that, for all $n\ge n^*(\epsilon)$ and all $T\ge T^*(\epsilon)$, and some $0<\delta_v\le 2$,
	\al{
	(a.1)&\quad \mathrm P\l(\min( n/ {\log^{2/\delta_v} T},\sqrt {T})\, \vert\wh{\sigma}_{i}^2-\sigma_i^2\vert\ge \eta(\epsilon)\r) < \epsilon,\;\text{ for any given $t=i,\ldots, n$,}\nn\\
	(a.2)&\quad \mathrm P\l(\min( n/ {\log^{2/\delta_v} T},\sqrt {T/\log n})\, \Vert\wh{\bm \Sigma}_{n}^\xi-\bm\Sigma_n^\xi\Vert\ge \eta(\epsilon)\r) < \epsilon,\nn\\
	(a.3)&\quad \mathrm P\l(\min( n/ {\log^{2/\delta_v} T},\sqrt {T})\, \Vert\wh{\mbf A}-\mbf A\Vert\ge \eta(\epsilon)\r) < \epsilon,\nn\\
	(a.4)&\quad \mathrm P\l(\min( n/ {\log^{2/\delta_v} T},\sqrt {T})\, \Vert\wh{\bm \Gamma}^v-\bm\Gamma^v\Vert\ge \eta(\epsilon)\r) < \epsilon;\nn
	}
	\item [(b)] as $n,T\to\infty$, if $n^{-1}\sqrt {T}\log^{2/\delta_v}T\to 0$, 
	\al{
	(b.1)&\quad\sqrt T (\wh{\sigma}_{i}^2-\sigma_i^2)\stackrel{d}{\to}\mathcal N(0,\sigma_i^4(\kappa_i+2)),\;\text{ for any given $i=1,\ldots, n$,}\nn\\
	(b.2)&\quad\sqrt T (\text{\upshape vec}(\wh{\mbf A})-\text{\upshape vec}(\mbf A))\stackrel{d}{\to}\mathcal N(\mbf 0_{r^2},\bm\Gamma^v\otimes (\bm\Gamma^F)^{-1}),\nn\\
	(b.3)&\quad\sqrt T (\text{\upshape vech}(\wh{\bm\Gamma}^v)-\text{\upshape vech}(\bm\Gamma^v))\stackrel{d}{\to}\mathcal N(\mbf 0_{r(r+1)/2},2\mbf D^\dag(\bm\Gamma^v\otimes\bm\Gamma^v)(\mbf D^\dag)^\prime),\nn
	}
	with $\bm\Gamma^F=\mbf I_r$, because of Assumption \ref{ass:ident}(b), and 
	where $\kappa_i= \E[\xi_{it}^4]/\sigma_i^4-3$ and $\mbf D^\dag=(\mbf D^\prime\mbf D)^{-1}\mbf D$ with $\mbf D$ being $r^2\times r(r+1)/2$ such that $\mbf D\text{\upshape vech}(\bm\Gamma^v)=\text{\upshape vec}(\bm\Gamma^v).$
\end{compactenum}
\end{prop}

We conclude with a series of general remarks.

\begin{rem}\label{rem:ident}\upshape{
From Propositions \ref{prop:load} and \ref{prop:factors}, we immediately have that, as $n,T\to\infty$,
\al{
&\l\Vert\frac {\wh{\bm\Lambda}_n^\prime\wh{\bm\Lambda}_n}n -\bm \Sigma_\Lambda\r\Vert
=o_p(1),\qquad
\l\Vert \frac{\wh{\bm{\mathcal F}}_T^\prime\wh{\bm{\mathcal F}}_T} T - \bm\Gamma^F\r\Vert 
=o_p(1),\nn
}
by Assumption \ref{ass:common}(a), which defines $\bm\Sigma_\Lambda$, and because the sample covariance matrix $T^{-1}{\bm{\mathcal F}}_T^\prime{\bm{\mathcal F}}_T$ is a consistent estimator of $\bm\Gamma^F$ (see Lemma \ref{lem:consistCOV}). Moreover, by Assumption \ref{ass:ident}(b), $\bm \Sigma_\Lambda=\lim_{n\to\infty} n^{-1}\mbf M_n^\chi$, which is a positive definite diagonal matrix and $\bm\Gamma^F=\mbf I_r$. Therefore, the EM estimator of the loadings and the related Kalman Smoother estimator of the factors satisfy the identifying constraints asymptotically, as $n,T\to\infty$. We verify this result numerically in Section \ref{sec:mc}. This is in agreement with Assumption \ref{ass:ident}(b) which imposes the identifying constraints only in the limit $n,T\to\infty$. This approach differs from other works on factor models (see, e.g., \citealp{baili12,baili16}, or \citealp{baing13}) where the identifying constraints are assumed to hold for any given $n$ and $T$ (see Remark \ref{rem:identificL}).

In principle, we could obtain estimators satisfying the identifying constraints in finite samples by imposing them ex-post in an additional step, as \citet[Section 8]{baili12} suggested in QML estimation of a static factor model. However, empirical works using the EM algorithm rarely apply this additional step.
}
\end{rem}

\begin{rem}\label{rem:logs}
\upshape{
The constraints $T^{-1}\sqrt {n}\log n\to 0$ and $n^{-1}\sqrt {T}\log^{2/\delta_v} T\to 0$ are common  (up to the presence of logarithmic terms) in the factor model literature (see, e.g., \citealp[Theorems 1 and 2]{Bai03}, for PC estimation) and are compatible. Indeed, they are simultaneously fulfilled if we assume that there exist some finite positive reals $\underline{\gamma}>1/2$ and $\overline{\gamma}<2$ such that $T^{\underline{\gamma}}<n<T^{\overline{\gamma}}$, as $T\to\infty$. 
When $n$ and $T$ have the same order of magnitude, as in many macroeconomic and financial datasets, these assumptions on the relative rates of divergence of $n$ and $T$ are very mild. 

In particular, the logarithmic term in the consistency rates of Propositions \ref{prop:load} and \ref{prop:factors} comes from Assumption \ref{ass:tails} of sub-Weibull tails, which is slightly more general than the typical assumption of sub-Gaussianity made when studying estimators of high-dimensional models. This is, however, a modest price to pay. In particular, under Gaussianity, $\delta_v=2$, while under distributions with sub-exponential tails, $\delta_v=1$. These logarithmic terms come essentially from two errors. The first, which is $O_p(n^{-1}\log^{2/\delta_v} T)$, is due to the necessity of finding a uniform bound over $t$ for the difference between the log-likelihood of the DFM in \eqref{eq:LL0true} and the static factor model log-likelihood considered in \citet{baili16}. This terms involves the sum of squared factors (see Lemma \ref{lem:2max}). The second one, which is $O_p(T^{-1/2}\sqrt{\log n})$, is due to the uniform bound for $\max_{i=1,\ldots, n}\vert \wh{\sigma}_i^2-\sigma_i^2\vert$ obtained when estimating the factors (see Lemma \ref{lem:sigmaunifhat} and also \citealp[Lemmas A3 and B1]{FLM11}).}
\end{rem}

\begin{rem}
\upshape{
The asymptotic properties of the estimators are unaffected if the number of factors $r$ is estimated. Indeed, consider a consistent estimator $\wh r$ i.e., such that $\mathrm P(\wh r = r)\to 1$, as $n,T\to\infty$, as for example the one in \citet{baing02}. Then, for any $z\in\mathbb R$ and any $i=1,\ldots, n$ and $t=1,\ldots,T$, it is easy to prove that $\mathrm P(\wh{\chi}_{it}\le z)=\mathrm P(\{\wh{\chi}_{it}\le z\}\,\vert\, \{\wh r=r\})+o_p(1)$ (see  \citealp[footnote 5]{Bai03}).

Similarly, the asymptotic properties of the estimators are unaffected if the order $p_F$ of the the VAR for the factors is estimated. This approach is asymptotically equivalent to computing the BIC using the true factors $\mbf F_t$ because we could estimate $p_F$ through the consistent PC estimator of the factors $\wt{\mbf F}_t$  (see Lemma \ref{lem:PCAF} and \citealp[Theorem 1]{Bai03}). And, in turn, the BIC is known to select the true lag order consistently \citep[Theorem 1]{hannan80}. Therefore, $\mathrm P(\wh p_F = p_F)\to 1$, as $n,T\to\infty$. Following the same reasoning of the previous remark, it is easy to show that, for any $z\in\mathbb R$ and any $i=1,\ldots, n$ and $t=1,\ldots,T$, $\mathrm P(\wh{\chi}_{it}\le z)=\mathrm P(\{\wh{\chi}_{it}\le z\}\,\vert\, \{\wh p_F= p_F\})+o_p(1)$. 	
}
\end{rem}

\begin{rem}
\upshape{
Given the asymptotic equivalence of Kalman filter, smoother, and WLS estimators of the factors, we might expect that the MSE obtained from either the Kalman filter or the smoother, i.e.,
${\mbf P}_{t|t}$ or ${\mbf P}_{t|T}$, respectively, asymptotically coincide (inflated by $n$) with the asymptotic covariance matrix $\bm{\mathcal W}_t$ of $\wh{\mbf F}_t$ defined in Proposition \ref{prop:factors}.
 However, this is not the case since we estimate a mis-specified model. Indeed, as $n\to\infty$, we can shown that  both $n{\mbf P}_{t|t}$ and $n{\mbf P}_{t|T}$ are asymptotically equivalent to $(\bm\Sigma_{\Lambda\Sigma\Lambda})^{-1}$ (see Lemma \ref{lem:KFGLS}), which, as shown in Proposition \ref{prop:factors}(c) is the asymptotic covariance matrix of $\wh{\mbf F}_t$ only if the model is correctly specified. In other words, the mis-specified Kalman filter and smoother do not estimate the true MSE. Although this has no effect on our asymptotic results (see Remark \ref{rem:PtTMSE}), still we cannot use the estimated MSEs ${\mbf P}_{t|t}^{(k^*+1)}$ or ${\mbf P}_{t|T}^{(k^*+1)}$ for making inference.
The true Kalman filter MSE, accounting for the model mis-specification, is derived by \citet[Section 2.1]{harvey2009computing} and, for any $t=1,\ldots, T$, it is given by the recursions
\al{
{\bm \Pi}_{t|t} =&\, \bm \Pi_{t|t-1} +  \mbf P_{t|t-1} \bm\Lambda_n^\prime ( \bm\Lambda_n\mbf P_{t|t-1}\bm\Lambda_n^\prime+\bm\Sigma_n^\xi)^{-1}( \bm\Lambda_n\bm\Pi_{t|t-1}\bm\Lambda_n^\prime+\bm\Gamma_n^\xi)( \bm\Lambda_n\mbf P_{t|t-1}\bm\Lambda_n^\prime+\bm\Sigma_n^\xi)^{-1}\bm\Lambda_n \mbf P_{t|t-1}\nn\\
&- \mbf P_{t|t-1} \bm\Lambda_n^\prime ( \bm\Lambda_n\mbf P_{t|t-1}\bm\Lambda_n^\prime+\bm\Sigma_n^\xi)^{-1}\bm\Lambda_n \bm \Pi_{t|t-1}
- \bm \Pi_{t|t-1}\bm\Lambda_n^\prime ( \bm\Lambda_n\mbf P_{t|t-1}\bm\Lambda_n^\prime+\bm\Sigma_n^\xi)^{-1}\bm\Lambda_n\mbf P_{t|t-1},\label{eq:dellefoche}\\
\bm \Pi_{t|t-1}=&\, \mbf A{\bm \Pi}_{t-1|t-1}  \mbf A^\prime+\bm\Gamma^v,\nn
}
where $\mbf P_{t|t-1}$ is the one-step-ahead Kalman filter MSE (see \eqref{eq:pred2}).
%
%
%
As expected, $n\bm{\Pi}_{t|t}$ is asymptotically equivalent, as $n\to\infty$, to the asymptotic covariance matrix $\bm{\mathcal W}_t$ of $\wh{\mbf F}_t$ (see Lemma \ref{lem:foche}). However,  since both $\mbf P_{t|t-1}$ and $\bm \Pi_{t|t-1}$ depend on $\mbf A$ and $\bm\Gamma^v$, for finite $n$ this MSE accounts explicitly also for the autocorrelation on the factors.
}
\end{rem}

\begin{rem}
\upshape{ The results of Propositions \ref{prop:plainvanilla}, \ref{prop:load}, \ref{prop:factors}, \ref{th:chi}, and \ref{prop:altri} would also hold if we had allowed for serial heteroskedasticity of the idiosyncratic components, i.e., for time-varying second moments so that $\E[\bm\xi_{nt}\bm\xi_{ns}^\prime]=\bm\Gamma_{n,ts}^\xi$. For this case, which we do not consider explicitly, we refer to \citet{baili16}, who show that the estimators of the idiosyncratic variances, $\wh{\sigma}^2_i$, $i=1,\ldots,n$, have to be considered as estimators of the average variances $\bar{\sigma}_i^2=T^{-1}\sum_{t=1}^T\sigma_{i,t}^2$; hence, in all above results, we should replace $\sigma_i^2$ with $\bar{\sigma}_i^2$. This approach amounts to maximizing a log-likelihood that has an additional degree of mis-specification because we use the time average idiosyncratic variances rather than the true time-varying variances. In this case, the asymptotic covariance matrix $\bm{\mathcal W}_t$ of the estimated factors in Proposition \ref{prop:factors} becomes effectively a time-varying matrix.}
\end{rem}

\subsection{Convergence of the EM algorithm under generic initialization}\label{sec:contra}
In this section, we discuss how our asymptotic results change if we initialize the EM algorithm with a generic initial estimator of the parameters, say $\check{\bm\varphi}_n^{(0)}=(\text{vec}(\check{\bm\Lambda}_n^{(0)})^\prime\; \check{\sigma}^{2(0)}_{1}\cdots \check{\sigma}^{2(0)}_{n}\; \text{vec}(\check{\mbf A}^{(0)})^\prime\;\text{vech}(\bm\Gamma^{v(0)})^\prime)^\prime$, having still elements belonging to $\mathcal O_n$ as defined in Section \ref{enzuccio}, and, thus, satisfying Assumptions \ref{ass:common}, \ref{ass:idio}, and \ref{ass:ident}.

For fixed  $n$ and in a general setting, \citet[]{BWB17} prove that the EM algorithm defines a contraction path towards a local maximum of the likelihood, $\wh{\bm\varphi}_n^{**}$. To give more details and understand the relation with our results, we need to introduce some general definitions.  First, consider any initial estimator belonging to a {closed} neighborhood of the local maximum of given Euclidean radius $\varrho>0$, i.e.,  
$\check{\bm\varphi}_n^{(0)}\in\mathcal B(\varrho; \wh{\bm\varphi}_n^{**})\subset \mathbb R^Q$. In our setting, we can think of $\mathcal B(\varrho; \wh{\bm\varphi}_n^{**})\equiv \mathcal O_n$. Then, define the EM operator $\bm M_T:\mathbb R^{Q}\to\mathbb R^{Q}$ such that $\bm M_T(\wh{\bm\varphi}_n^{(k)})=\wh{\bm\varphi}_n^{(k+1)}$. We have $\Vert \E[M_T(\underline{\bm\varphi}_n)]-\wh{\bm\varphi}_n^{**}\Vert \le \beta\Vert \underline{\bm\varphi}_n-\wh{\bm\varphi}_n^{**}\Vert $ for some $\beta\in(0,1)$ and all  $\underline{\bm\varphi}_n\in \mathcal B(\varrho; \wh{\bm\varphi}_n^{**})$
(\citealp[Theorem 1]{BWB17}). Second, for any given $T$ and $\delta\in(0,1)$, let $\varepsilon_{T,\delta}$ be the smallest scalar such that
$\mathrm P(\sup_{\underline{\bm\varphi}_n\in \mathcal B(\varrho;\wh{\bm\varphi}_n^{**})} \Vert M_T(\underline{\bm\varphi}_n)- \E[M_T(\underline{\bm\varphi}_n)]\Vert \le \varepsilon_{T,\delta})\ge 1-\delta$. 

It follows that, if $T$ is large enough such that $\varepsilon_{T,\delta} \le (1-\beta)\varrho$, then, for any $k\ge 0$, the EM operator defines a contraction towards the maximum of the likelihood with high-probability (\citealp[Theorem 2]{BWB17}):
\beq
\mathrm P\l( \Vert \wh{\bm\varphi}_n^{(k+1)}-\wh{\bm\varphi}_n^{**}\Vert \le \beta^{k+1} \Vert \check{\bm\varphi}_n^{(0)}-\wh{\bm\varphi}_n^{**}\Vert +\frac{\varepsilon_{T,\delta}}{1-\beta} \r)\ge 1-\delta.\label{eq:BALA}
\eeq

Now, under our mixing Assumptions \ref{ass:common}(d)-\ref{ass:common}(h) and \ref{ass:idio}(c), $\varepsilon_{T,\delta}\to 0$, as $T\to\infty$. This, jointly with \eqref{eq:BALA}, has two implications, as $T\to\infty$:
\ben
\item [(a)] for all $k\ge 0$, $\Vert \wh{\bm\varphi}_n^{(k+1)}-\wh{\bm\varphi}_n^{**}\Vert\le \beta^{k+1} \Vert \check{\bm\varphi}_n^{(0)}-\wh{\bm\varphi}_n^{**}\Vert +o_p(1)$; 
\item [(b)] if $k \ge \log_{1/\beta} \varepsilon^{-1}_{T,\delta}\equiv k_T$, then $\Vert \wh{\bm\varphi}_n^{(k+1)}-\wh{\bm\varphi}_n^{**}\Vert=o_p(1)$.
\een

These two results apply if $n$ is fixed. But when $n\to\infty$, we might conjecture that those results will hold for each component of $\bm\varphi_n$ separately. And this is, indeed, what we verify in this paper. Consider the loadings $\bm\lambda_i$. As mentioned in the previous section, as $n\to\infty$ the likelihood has a unique global maximum, $\wh{\bm\lambda}_i^{*}$, which is consistent because it is the QML estimator. Therefore, $\Vert \wh{\bm\lambda}_i^{**}-{\bm\lambda}_i\Vert\le \Vert \wh{\bm\lambda}_i^{**}-\wh{\bm\lambda}_i^*\Vert+\Vert \wh{\bm\lambda}_i^{*}-{\bm\lambda}_i\Vert=o_p(1)$, as $n,T\to\infty$. If we initialize the EM algorithm with the consistent PC estimator $\wh{\bm\lambda}_i^{(0)}$, we prove that, as $n,T\to\infty$, result (a) above still applies, i.e., there exists a $\beta_\lambda\in(0,1)$ such that, for all $k\ge 0$,
\beq
	\Vert\wh{\bm\lambda}_i^{(k+1)}-\wh{\bm\lambda}_i^{**} \Vert \le \l\{\Vert\wh {\bm\lambda}_i^{(0)}-{\bm\lambda}_i \Vert+ \Vert\wh {\bm\lambda}_i^{**}-{\bm\lambda}_i \Vert\r\} \beta^{k+1}_{\lambda}+ o_p(1)=o_p(1).\label{eq:BALABL}
\eeq
Actually, we are  also able to show that  the contraction factor is such that, as $n,T\to\infty$:
\beq
	\beta_{\lambda}=\l\Vert \mbf I_r-\l(\sum_{t=1}^T \l\{\mbf F_{t|T}^*\mbf F_{t|T}^{*\prime}+\mbf P^{*}_{t|T}\r\}\r)^{-1}\l(\sum_{t=1}^T \mbf F_t\mbf F_t^\prime\r)\r\Vert +o_p(1)=o_p(1), \label{eq:BALABLbeta}
\eeq
with $\mbf F_{t|T}^{*}$ and $\mbf P^{*}_{t|T}$ being the factors and their associated MSEs obtained from the Kalman smoother when using the QML estimator of the parameters---we refer to Lemma \ref{lem:convEM1} for a proof of \eqref{eq:BALABL} and \eqref{eq:BALABLbeta}. It follows that the convergence rate in \eqref{eq:BALABL} is faster than the one of the initial PC estimator, and we can treat the EM estimator as if it were the QML estimator. Moreover, \eqref{eq:BALABLbeta} means that \eqref{eq:BALABL} holds for any initial estimator (see Lemma \ref{lem:convEM1muzzo}). The intuition for this result is that, provided $n^{-1}\check{\bm\Lambda}^{(0)\prime}\check{\bm\Lambda}^{(0)}$ has full-rank, any cross-sectional averaging is likely to recover in high-dimensions a factors' space not too different from the true one, an argument often used when considering estimation methods based on aggregations schemes alternative to PCs \citep{WU15,fan2022learning}.
	 
For all other parameters, a relation like \eqref{eq:BALABL} holds as well, with possibly different contraction factors. Hence, result (a) still applies (see Lemma \ref{lem:convEMiAH}). However, we could not derive a result like \eqref{eq:BALABLbeta} in this case because computing analytic expressions of all those contraction factors is difficult. Thus, when we initialize the algorithm with any generic initial estimator, we have to rely on result (b); that is, we can guarantee convergence of the EM algorithm to the QML estimator, as $n,T\to\infty$, only if we let the algorithm run for a number of iterations $k\ge k_T$, where the larger $k$ is, the faster is the convergence rate.

%
%
\section{Efficiency and comparison with PC analysis}\label{sec:PCeff}
	
In this section, we compare the asymptotic covariances of the EM and Kalman smoother with those of the PC estimators, which are the optimal non-parametric estimators. 

From  Propositions \ref{prop:load}(b)   and \ref{prop:factors}(b), we see that consistency is not affected by estimating a mis-specified model with uncorrelated idiosyncratic components, but there is an efficiency loss due to this mis-specification, as shown by the sandwich forms of the asymptotic covariance matrices. In the case of uncorrelated idiosyncratic components, the log-likelihood we maximize is correctly specified. If the model is correctly specified, the EM estimator is the most efficient one because it is asymptotically equivalent to the Maximum Likelihood estimator. Thus, its asymptotic covariance attains the classical lower bound of the OLS estimator (see Propositions \ref{prop:load}(c)). Likewise, the asymptotic covariance of the factors attains the WLS lower bound (see Propositions \ref{prop:factors}(c)).

Although, in general, the EM algorithm and the Kalman smoother do not provide the most efficient estimators, they can provide advantages with respect to the PC estimators. 

\begin{prop}[\textsc{Efficiency}]\label{prop:eff} 
Let $\bm{\mathcal V}_i^{\text{\tiny{\upshape PC}}}$ and $\bm{\mathcal W}_t^{\text{\tiny\upshape PC}}$ be the asymptotic covariance matrices of the loadings and factors estimated via PC analysis, then, under Assumptions \ref{ass:common}, \ref{ass:idio}, \ref{ass:ind}, \ref{ass:linear}, \ref{ass:tails}, and \ref{ass:ident}:
\begin{compactenum}
\item [(a)] if $\sqrt T\log^{1/\delta_v}T/n\to 0$, as $n,T\to\infty$, then $\bm{\mathcal V}_i^{\text{\tiny{\upshape PC}}}=\bm{\mathcal V}_i$ for any $i=1,\ldots, n$;
\item [(b)] if $\sqrt {n\log n}/T\to 0$ and  $n^{-1}\sum_{\substack{i,j=1,i\ne j}}^n\vert\E[\xi_{it}\xi_{jt}]\vert\to 0$,  as $n,T\to\infty$, then $(\bm{\mathcal W}^{\text{\tiny\upshape PC}}_t-\bm{\mathcal W}_t)$ is a positive definite matrix. 
\end{compactenum}
\end{prop}

Part (a) follows immediately once we impose the identifying Assumption \ref{ass:ident} to the results about PC estimation of the loadings
(see \citealp[Theorem 1]{MBPCAQML}, and \citealp[Theorem 2]{Bai03}).  Therefore, although we estimate a mis-specified model, the EM estimator of the loadings is as efficient as the PC estimator. 

Turning to part (b), by imposing the identifying Assumption \ref{ass:ident}, the asymptotic covariance of the PC estimator of the factors is
$
\bm{\mathcal W}^{\text{\tiny PC}}_t=(\bm\Sigma_\Lambda)^{-1}\{\lim_{n\to\infty} n^{-1} \bm\Lambda_n^\prime \bm\Gamma_n^\xi \bm\Lambda_n\}(\bm\Sigma_\Lambda)^{-1}
$ 
(see Lemma \ref{lem:PCAF}, and \citealp[Theorem 1]{Bai03}), and if the true model were an exact factor model, i.e., $\E[\xi_{it}\xi_{jt}]=0$ if $i\ne j$ so that $\bm\Gamma_n^\xi=\bm\Sigma_n^\xi$, then we would have 
$
\bm{\mathcal W}^{\text{\tiny PC}}_t=(\bm\Sigma_\Lambda)^{-1}\{\lim_{n\to\infty} n^{-1} \bm\Lambda_n^\prime \bm\Sigma_n^\xi \bm\Lambda_n\}(\bm\Sigma_\Lambda)^{-1}
$, which is the unfeasible OLS asymptotic covariance matrix in presence of heteroskedastic errors. For the Kalman smoother, from Proposition \ref{prop:factors} we have
$
\bm{\mathcal W}_t =(\bm\Sigma_{\Lambda\Sigma\Lambda})^{-1} 
\{\lim_{n\to\infty} n^{-1} \bm\Lambda_n^\prime(\bm\Sigma_n^{\xi})^{-1} \bm\Gamma_n^\xi (\bm\Sigma_n^{\xi})^{-1}\bm\Lambda_n\}
(\bm\Sigma_{\Lambda\Sigma\Lambda})^{-1}
$, and if the true model were an exact factor model, this would reduce to
$
\bm{\mathcal W}_t=(\bm\Sigma_{\Lambda\Sigma\Lambda})^{-1} 
$, which is the unfeasible WLS asymptotic covariance matrix. In this case, the Kalman smoother estimator is more efficient because it takes into account heteroskedasticity, whereas the PC estimator ignores the possibility of individual-specific variances. 
	
Here, we go one step further and show that if  the idiosyncratic covariance matrix ${\bm\Gamma}_n^\xi$ is sparse enough, we can still have efficiency gains compared to PC analysis. Indeed, if the total contribution of the off-diagonal elements ${\bm\Gamma}_n^\xi$ is negligible compared to $n$, we can expect the asymptotic covariance of the estimated factors to be close to the one we would have for an exact factor model, which we know is smaller than the PC asymptotic covariance. The sparsity condition we assume is the same as the one on \citet[Assumption 3.1]{bailiao16}. Although this sparsity condition is hard to verify in practice and is seldom exactly satisfied by the data, in our MonteCarlo results in Section \ref{sec:mc}, we show that the Kalman smoother estimator tends to perform better than the PC estimator even under more general idiosyncratic covariance structures like banded matrices.

\begin{rem}
\upshape{If we assume uncorrelated and homoskedastic idiosyncratic components, i.e., such that $\bm\Gamma_n^\xi=\psi\mbf I_n$ for some $\psi>0$, then it is easy to see that $\bm{\mathcal V}_i = \bm{\mathcal V}_i^{\text{\tiny PC}}=\psi(\bm\Gamma^F)^{-1}=\psi\mbf I_r$, by Assumption \ref{ass:ident}(b), and $\bm{\mathcal W}_t=\bm{\mathcal W}_t^{\text{\tiny PC}}=\psi(\bm\Sigma_\Lambda)^{-1}$. In this case, the EM and Kalman smoother estimators are as efficient as the PC estimators.}
\end{rem}

\begin{rem}
\upshape{There are other estimators that could be more efficient. 
First, \citet{bailiao16}, \citet{wang2019penalized}, and \citet{poignard2020statistical} proposed penalized QML-type estimators of the loadings and of the idiosyncratic covariance matrix. These estimators are used to build a GLS estimator of the factors. This approach addresses cross-sectional idiosyncratic correlations and heteroskedasticity, but not 
serial idiosyncratic correlations. 
Second, \citet{BT11} propose a GLS estimator of the loadings, based on the classical \citet{CU49} approach, and a WLS estimator for the factors based on that loadings estimator and estimates of the idiosyncratic variances. This approach addresses cross-sectional idiosyncratic heteroskedasticity and serial idiosyncratic correlations, but not cross-sectional idiosyncratic correlations. Finally, \citet{LM02} address all idiosyncratic cross-autocorrelations by 
assuming a sparse VAR model for the idiosyncratic components. To this end, they apply the \citet{CU49} approach and embed a penalty into an alternating minimization algorithm. They recover the factors via GLS. Their theoretical analysis applies only to the finite-dimensional case.
None of these three approaches models the factors' dynamics.
}
\end{rem}

%
%

\section{Inference}\label{sec:test}
To conduct inference, we need asymptotic covariances of the loadings and factors matrices and their
estimators. 
\begin{cor} \label{cor:matriciAN}
Under the same assumptions of Propositions \ref{prop:load}  and \ref{prop:factors}, as $n,T\to\infty$:
\begin{compactenum}
	\item [(a.1)]  for any finite $\bar n$ and any given sequence of integers $\{s(1)\ldots,s(\bar n)\} \subset \{1,\ldots,n\}$,
	 let $\text{\upshape vec}(\wh{\bm\Lambda}_{\bar n})=(\wh{\bm\lambda}_{s(1)}^\prime\cdots \wh{\bm\lambda}_{s(\bar n)}^\prime)^\prime$ 
	 and $\text{\upshape vec}({\bm\Lambda}_{\bar n})=({\bm\lambda}_{s(1)}^\prime\cdots {\bm\lambda}_{s(\bar n)}^\prime)^\prime$,  
	 as $n,T\to\infty$, if $n^{-1}\sqrt {T}\log^{2/\delta_v}T\to 0$,
	 \al{
	\sqrt T \big\{\text{\upshape vec}(\wh{\bm\Lambda}_{\bar n})-\text{\upshape vec}({\bm\Lambda}_{\bar n})\big\}\stackrel{d}{\to}\mathcal N(\mbf 0_{r\bar n},\bm{\mathcal V}_{\bar n}),\nn
	}
	where  
	$$
	\bm{\mathcal V}_{\bar n}=(\mbf I_{\bar n}\otimes \bm\Gamma^F)^{-1}\l(\lim_{T\to\infty} \frac 1T\sum_{t=1}^T\sum_{s=1}^T\E[\bm\xi_{\bar n t}\bm\xi_{\bar ns}] \otimes\E[\mbf F_t\mbf F_s^\prime]\r)(\mbf I_{\bar n}\otimes \bm\Gamma^F)^{-1},
	$$
	with $\bm\xi_{\bar n t}=(\xi_{s(1)t}\cdots \xi_{s(\bar n)t})^\prime$, with $\bm\Gamma^F=\mbf I_r$, because of Assumption \ref{ass:ident}(b); 
	\item [(a.2)] if $\E[(\bm\xi_{n1}^\prime\cdots\bm\xi_{nT}^\prime)^\prime(\bm\xi_{n1}^\prime\cdots\bm\xi_{nT}^\prime)]=\mbf I_T\otimes \bm\Sigma_n^\xi$ for all $n,T\in\mathbb N$, then $\bm{\mathcal V}_{\bar n}={\bm\Sigma}_{\bar n}^\xi\otimes\, (\bm\Gamma^F)^{-1}$, with $\bm\Gamma^F=\mbf I_r$, because of Assumption \ref{ass:ident}(b); 
		\item [(b.1)] for any finite $\bar T$ and any given sequence of integers $\{s(1)\ldots,s(\bar T)\} \subset \{1,\ldots,T\}$,
	 let $\text{\upshape vec}(\wh{\bm{\mathcal F}}_{\bar T})=(\wh{\mbf F}_{s(1)}^\prime\cdots \wh{\mbf F}_{s(\bar T)}^\prime)^\prime$ 
	 and $\text{\upshape vec}({\bm {\mathcal F}}_{\bar T})=({\mbf F}_{s(1)}^\prime\cdots {\mbf F}_{s(\bar T)}^\prime)^\prime$, then, 
	 as $n,T\to\infty$, if $T^{-1}\sqrt {n\log n}\to 0$,
	 \[
	 \sqrt n \big\{\text{\upshape vec}(\wh{\bm{\mathcal F}}_{\bar T})-\text{\upshape vec}({\bm{\mathcal F}}_{\bar T})\big\}\stackrel{d}{\to}\mathcal N(\mbf 0_{r\bar T},\bm{\mathcal W}_{\bar T}),
	 \]
	 where 
	 $$
	\bm{\mathcal W}_{\bar T}=(\mbf I_{\bar T}\otimes \bm\Sigma_{\Lambda\Sigma\Lambda})^{-1} \l(\lim_{n\to\infty} \frac 1 n\sum_{i=1}^n\sum_{j=1}^n \frac{\E[\bm\zeta_{i\bar T}\bm\zeta_{j\bar T}] \otimes \bm\lambda_i\bm\lambda_j  }{\sigma_i^2\sigma_j^2} \r)(\mbf I_{\bar T}\otimes\bm\Sigma_{\Lambda\Sigma\Lambda})^{-1},
	$$
	with $\bm \zeta_{i\bar T}=(\xi_{is(1)}\cdots\xi_{is(\bar T)})^\prime$;
	\item [(b.2)] if $\E[(\bm\zeta_{1T}^\prime\cdots \bm\zeta_{nT}^\prime)^\prime(\bm\zeta_{1T}^\prime\cdots \bm\zeta_{nT}^\prime)]=\bm\Sigma_n^\xi\otimes\mbf I_T$ for all $n,T\in\mathbb N$, then $\bm{\mathcal W}_{\bar T}=\mbf I_{\bar T}\otimes(\bm\Sigma_{\Lambda\Sigma\Lambda})^{-1}$.
\end{compactenum}
\end{cor}


For any given $i,j=1,\ldots,n$, the most general estimator of the asymptotic covariance between $\wh{\bm\lambda}_i$ and $\wh{\bm\lambda}_j$ is given by:\footnote{As discussed in Remark \ref{rem:ident}, the sample covariance of $\wh{\mbf F}_t$ is equal to the $r$-dimensional identity matrix only asymptotically; hence, it must be included in the estimator of the asymptotic covariance.}
\begin{align}
\wh{\bm{ \mathcal V}}_{i,j}^{\text{\tiny(HAC)}}
&=\l(\frac 1 T\sum_{t=1}^T \wh {\mbf F}_t\wh{\mbf F}_t^\prime\r)^{-1}\l( \frac 1T \sum_{t=1}^T\sum_{s=1}^T \mathrm K(t,s) \l\{ \wh {\mbf F}_t\wh {\mbf F}_s^\prime \wh{\xi}_{it}\wh{\xi}_{js} \r\} \r)\l(\frac 1 T\sum_{t=1}^T \wh {\mbf F}_t\wh{\mbf F}_t^\prime\r)^{-1},\label{eq:robcovloadiest}
\end{align}
where $\wh{\xi}_{it}=x_{it}-\wh{\chi}_{it}$ is the estimated idiosyncratic component of the $i$th variable at time $t$. $\mathrm K(t,s) = 1-\frac{|t-s|}{M_T+1}$, if $|t-s|\le M_T$ and zero otherwise, with $M_T$ is such that $M_T\to\infty$ and $M_T/T\to 0$, as $T\to\infty$.
Consistency of \eqref{eq:robcovloadiest}, as $n,T\to\infty$, follows from \citealp[Theorem 6]{Bai03}  combined with Propositions \ref{prop:factors} and \ref{th:chi}. 

For any given $t=1,\ldots, T$ and $k=0,\ldots, M_T$, with $M_T$ defined above, the most general estimator of the asymptotic covariance between $\wh{\mbf F}_t$ and $\wh{\mbf F}_{t-k}$ is given by:
\begin{align}
\wh{\bm{\mathcal W}}_{t,t-k}^{\text{\tiny(HAC)}} &= \l(\frac 1 n\sum_{i=1}^n\frac{\wh{\bm{\lambda}}_i\wh{\bm{\lambda}}_i^\prime}{\wh{\sigma}_i^2}\r)^{-1}
\!\!\l(
\frac 1 n \sum_{i=1}^n\sum_{j=1}^n
\mathrm K(i,j)
\l\{\frac{\wh{\bm{\lambda}}_i\wh{\bm{\lambda}}_j^\prime}{\wh{\sigma}_i^2\wh{\sigma}_j^2} \wh{\gamma}^\xi_{ij, k}
\r\}
\r)
 \l(\frac 1 n\sum_{i=1}^n\frac{\wh{\bm{\lambda}}_i\wh{\bm{\lambda}}_i^\prime}{\wh{\sigma}_i^2}\r)^{-1},\label{eq:robcovfacttest}
\end{align}
where $\wh{\gamma}^\xi_{ij, k} = T^{-1}\sum_{t=k+1}^T \wh{\xi}_{it}\wh{\xi}_{j,t-k}$ and $\mathrm K(i,j)=1$ if $1\le i,j\le m_{n,T}$ and zero otherwise, with $m_{n,T}\to\infty$ and $m_{n,T}/\min(n,T)\to 0$, as $n,T\to\infty$. 
Consistency of \eqref{eq:robcovfacttest}, as $n,T\to\infty$, follows from \citet[Theorem 4]{baing06} combined with Propositions \ref{prop:load}, \ref{th:chi}, and \ref{prop:altri}. For larger values of $k$, $\E[\xi_{it}\xi_{j,t-k}]$ is likely to be small due to our Assumption \ref{ass:idio}(b), and, thus, we can consider $\wh{\mbf F}_t$ and $\wh{\mbf F}_{t-k}$ as asymptotically uncorrelated. An alternative kernel function of the correlation based distance between $i$ and $j$ is considered by \citet{kim2022robust} who considers also averages of  \eqref{eq:robcovfacttest} when computed choosing different random permutations of the selected $m_{n,T}$ cross-sectional units. Alternatively, rather that smoothing via the use of a kernel, \citet{BR24} consider an estimator based on thresholding of the idiosyncratic sample covariance matrix (see also \citealp{FLM13}).

  
Finally, we can obtain an estimator of $\bm{\mathcal W}_t$  that accounts also for the autocorrelation of the factors from the true MSE of the Kalman filter given in \eqref{eq:dellefoche} when using the estimated parameters. Once again this requires an estimator of the idiosyncratic covariance matrix like the estimators discussed above.
The asymptotic covariance between $\wh{\mbf F}_t$ and $\wh{\mbf F}_{t-k}$ can be obtained by considering the Kalman filter with the augmented state vector $(\mbf F_t^\prime\cdots \mbf F_{t-k}^\prime)^\prime$ and by then considering the corresponding $r\times r$ off-diagonal block of the resulting $r(k+1)\times r(k+1)$ MSE.

%

Having the estimated loadings $\wh{\bm\lambda}_i$, the estimated factors $\wh{\mbf F}_t$, and any of the above estimators of $\bm{\mathcal V}_{i,j}$ and $\bm{\mathcal W}_{t,t-k}$, we can estimate the variance of the estimated common component by plugging these quantities into the expression in Proposition \ref{th:chi}(b). 

To conclude, the covariance estimators defined in this section can be used, together with the asymptotic distributions derived in Propositions \ref{prop:load} and \ref{prop:factors}, for inferential purposes. Examples are in Section \ref{sec:emp}.

\begin{rem}
\upshape{
From parts (a.2) and (b.2) we see that, under a correctly specified model, each estimated row $i$ of the loadings matrix is asymptotically uncorrelated with the other rows, and, likewise, each estimated realization of the factors at a given point in time $t$ is asymptotically uncorrelated with the other time periods. Estimators of the asymptotic covariances are easily built in this case as:
$\wh{\bm{ \mathcal V}}_{i}^{(0)} =\widehat{\sigma}_i^2(T^{-1}\sum_{t=1}^T \wh {\mbf F}_t\wh{\mbf F}_t^\prime)^{-1}$, 
and $\wh{\bm{ \mathcal W}}_{t}^{(0)}= (n^{-1}\sum_{i=1}^n (\wh{\sigma}_i^2)^{-1} {\wh{\bm{\lambda}}_i\wh{\bm{\lambda}}_i^\prime})^{-1}$,  
while $\wh{\bm{ \mathcal V}}_{i,j}^{(0)} =\mbf 0_{r\times r}$  if  $i\ne j$, and $\wh{\bm{ \mathcal W}}_{t,t-k}^{(0)}= \mbf 0_{r\times r}$ if $k\ne 0$. Consistency of these two estimators follows directly from Propositions \ref{prop:load}, \ref{prop:factors}, \ref{th:chi}, and \ref{prop:altri}.
}
\end{rem}

%
%
\section{Monte Carlo study}\label{sec:mc}
Throughout, we consider a model with $r=4$ factors, and we simulate the data according to 
\begin{align}
x_{it}&
={\bm\ell}_{i}^\prime {\bm f}_t+ \phi_i \xi_{it}, \quad 
{\bm f}_t ={\bm A}  {\bm f}_{t-1} + \mbf u_t,\quad
\xi_{it}=\alpha_i\xi_{it-1}+e_{it},\label{eq:Sim4}
\end{align}
where ${\bm\ell}_{i}$ has entries ${\ell}_{ij}\stackrel{iid}{\sim}\mathcal{N}(1,1)$, $i=1,\ldots, n$, $j=1,\ldots, r$; 
${\bm A}=\mu\check{\bm A} \{\nu^{(1)}({\check{\bm A}})\}^{-1}$, where $[\check{\bm A}]_{jj}\sim U[0.5,0.8]$, $[\check{\bm A}]_{jk}\sim U[0,0.3]$, $j,k=1,\ldots,r$, and $\mu=0.7$;
$u_{jt}\stackrel{iid}{\sim}(0,1)$, $j=1,\ldots, r$ following either a Gaussian, an Asymmetric Laplace, or a Skew-$t$ distribution, and with $\Cov(u_{it},u_{jt})=0$, for $i\ne j$;
$\alpha_i=\{0,\delta_i\}$, $i=1,\ldots, n$, where $\delta_i\stackrel{iid}{\sim}\mathcal{U}(0,\delta)$, and $\delta\in\{0,0.5\}$;
$e_{it}\stackrel{iid}{\sim}(0,\sigma_{ei}^2)$, $i=1,\ldots, n$, following either a Gaussian distribution with $\sigma_{ei}^2\sim U[0.5, 1.5]$, or an Asymmetric Laplace distribution with $\sigma_{ei}^2=1$, or a Skew-$t$ distribution with $\sigma_{ei}^2=1$;
%
%
%
 $\Cov(e_{it},e_{jt})=\tau^{\vert i-j\vert }$, $i,j=1,\ldots, n$, with $\tau\in\{0,0.5\}$ if $\vert i-j\vert \le 10$, and $\Cov(e_{it},e_{jt})=0$ otherwise;  and, last,
    $\phi_i=\sqrt{\theta_i{\wh{\Var}(\chi_{it})}(\wh{\Var}(\xi_{it}))^{-1}}$, $i=1,\ldots, n$, where $\wh{\Var}(\cdot)$ denotes the sample variance, and 
    $\theta_i\stackrel{iid}{\sim}\mathcal{U}(\bar{\theta}-0.25,\bar{\theta})$, and $\bar{\theta}=0.5$.
  The parameters $\mu$, $\tau$, $\delta$, and $\bar\theta$ are crucial and control: the persistence of the factors, the degrees of
cross-sectional and serial idiosyncratic correlation in the idiosyncratic components, and the noise-to-signal ratio, respectively.\footnote{In the case of the Asymmetric Laplace distribution, all the innovations have location 0, scale index $\lambda$, and asymmetry index $\kappa$, with $\kappa \sim \mathcal{U}(.9,1.1)$ and $\lambda=\sqrt{(1+\kappa^4)\kappa^{-2}}$, so that all the shocks have variance 1. In the case of the Skew-$t$ distribution, all the shocks have location 0, dispersion 1, skewness index $\gamma$, and tail index $\nu$, with $\nu_u\sim {U}(4,12)$, $\gamma_u \sim {U}(-.1,.1)$, $\nu_e\sim {U}(3,13)$, and $\gamma_e \sim {U}(-.15,.15)$.} 

Finally, in order to satisfy Assumptions \ref{ass:ident}(b) and \ref{ass:ident}(c), after we have generated the common component $\bm\chi_{nt}$ as in \eqref{eq:Sim4}, we construct the factors as $\mbf F_t=(\mbf M_n^{\chi})^{-1/2}\mbf V_n^{\chi\prime} \bm\chi_{nt}$ and the loadings as $\bm \Lambda_n=\mbf V_n^{\chi} (\mbf M_n^{\chi})^{1/2}$, where $\mbf M_n^\chi$ is the $r\times r$ diagonal matrix containing the eigenvalues of $T^{-1}\sum_{t=1}^T \bm\chi_{nt}\bm\chi_{nt}^\prime$, and $\mbf V_n^\chi$ is the $n\times r$ matrix having as columns the corresponding normalized eigenvectors and with sign fixed such that it has non-negative entries in the first row. 

We consider $n\in\{100,200,300, 500\}$, $T\in\{100,200,300, 500\}$ and $B=1000$ replications. At each replication $b$, we run the EM algorithm as described in Algorithm \ref{alg111}, thus obtaining an estimate of the loadings $\wh{\bm\lambda}_i^{(b)}$, the factors $\wh{\mbf F}_t^{(b)}$, and the common component $\wh{\chi}_{it}^{(b)}=\wh{\bm\lambda}_i^{(b)\prime}\wh{\mbf F}_t^{(b)}$. We initialize the EM algorithm either through the PC estimators as explained in Appendix \ref{app:prest}, or through a contaminated version of the PC estimator obtained using contaminated eigenvectors of the data. In particular, 
let $\wh{\mbf V}_n^x$ be the $n\times r$ matrix of the $r$ leading eigenvectors of $T^{-1}\sum_{t=1}^T \mbf x_{nt}\mbf x_{nt}^\prime$, then the contaminated eigenvector is $\wh{\mbf V}_n^{x,c}=\wh{\mbf V}_n^x+\bm {\mathcal Z}\bm \Upsilon^{1/2}$, where $[\bm {\mathcal Z}]_{ij}\stackrel{iid}{\sim} \mathcal N(0,1)$, $i=1,\ldots, n$, $j=1,\ldots,r$, and $[\bm \Upsilon]_{ij}=\iota^{\vert i-j\vert}$, $i,j=1,\ldots, r$, with $\iota \in (0,1)$---the bigger $\iota$ is, the stronger is the contamination.

\begin{figure}[t!]\caption{Simulation results - Convergence of the EM algorithm}\label{fig:EM}
\centering 
\setlength{\tabcolsep}{.0\textwidth}
\begin{tabular}{@{}C{.32\textwidth}C{.02\textwidth}C{.32\textwidth}C{.02\textwidth}C{.32\textwidth}@{}} 
\footnotesize Gaussian &&\footnotesize Asymmetric Laplace &&\footnotesize Skew-$t$\\
\includegraphics[width=.31\textwidth]{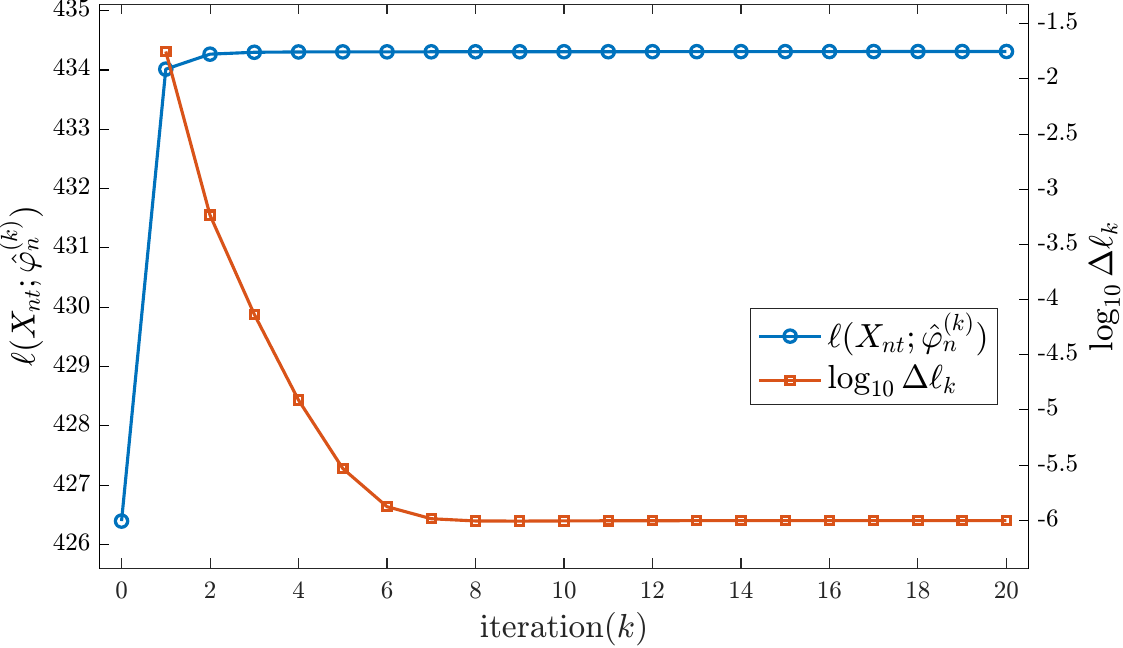} &&
\includegraphics[width=.31\textwidth]{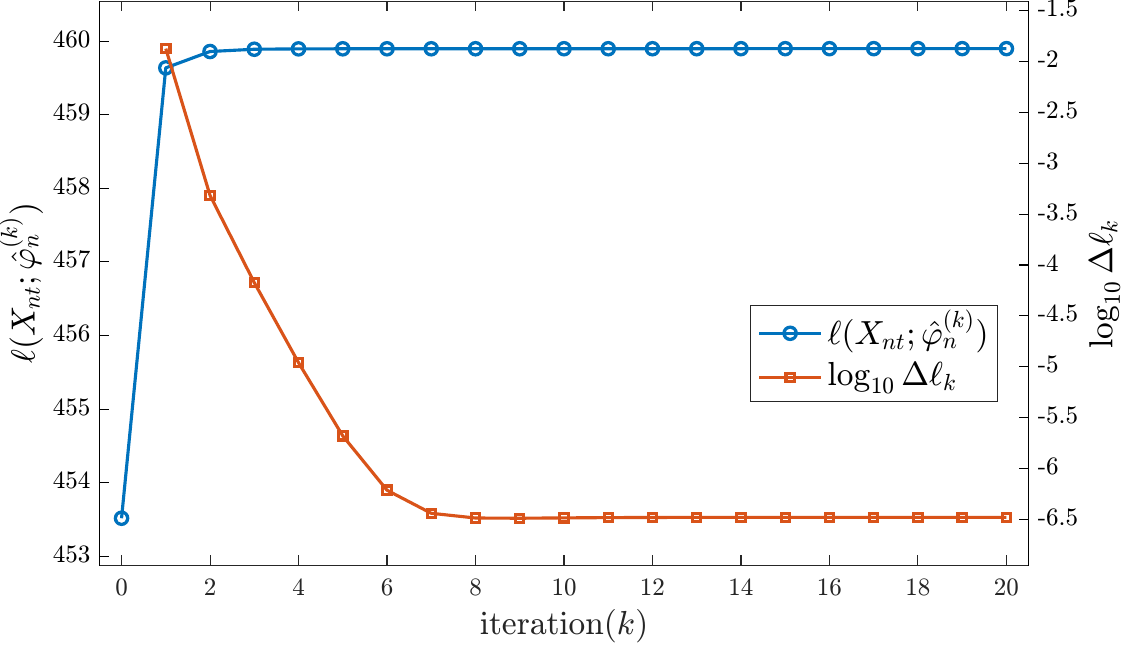}&&
\includegraphics[width=.31\textwidth]{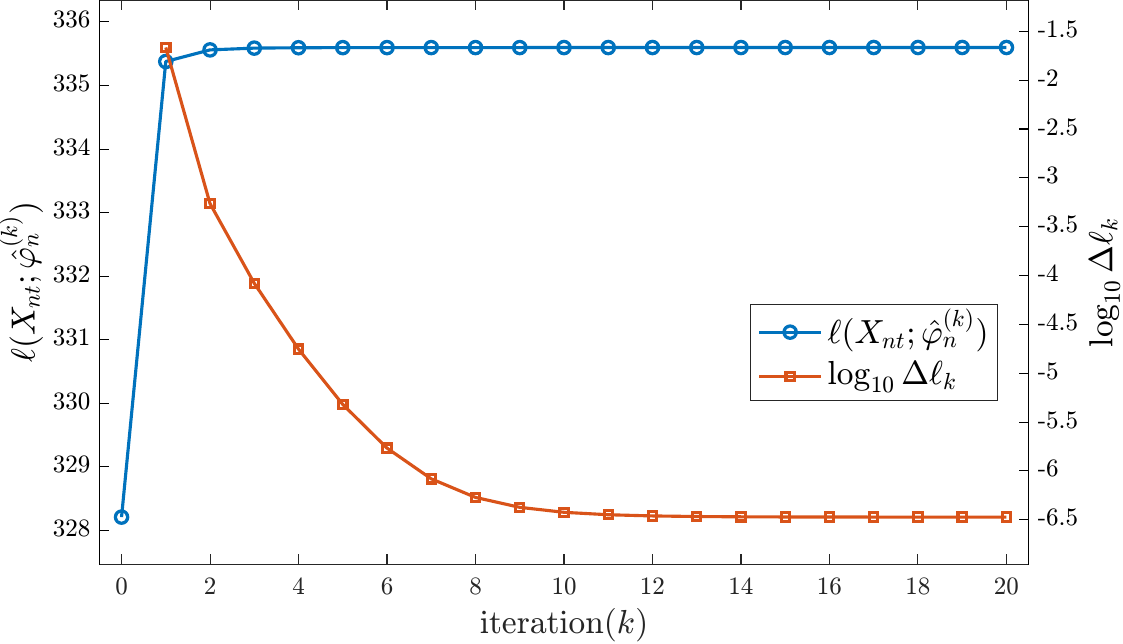}
\\
\footnotesize Gaussian &&\footnotesize Asymmetric Laplace &&\footnotesize Skew-$t$\\
\includegraphics[width=.275\textwidth,height=.175\textwidth]{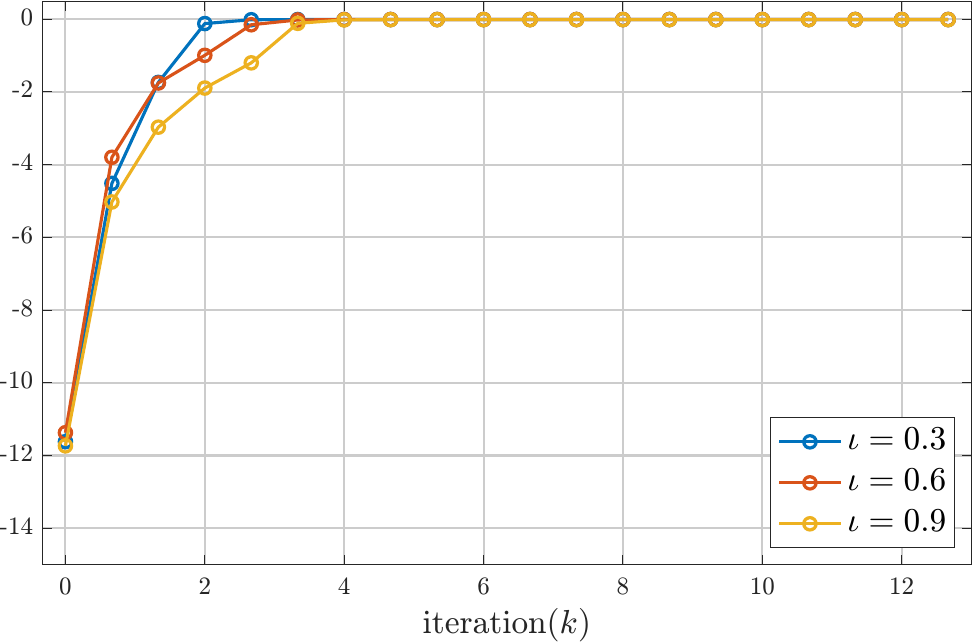} &&
\includegraphics[width=.275\textwidth,height=.175\textwidth]{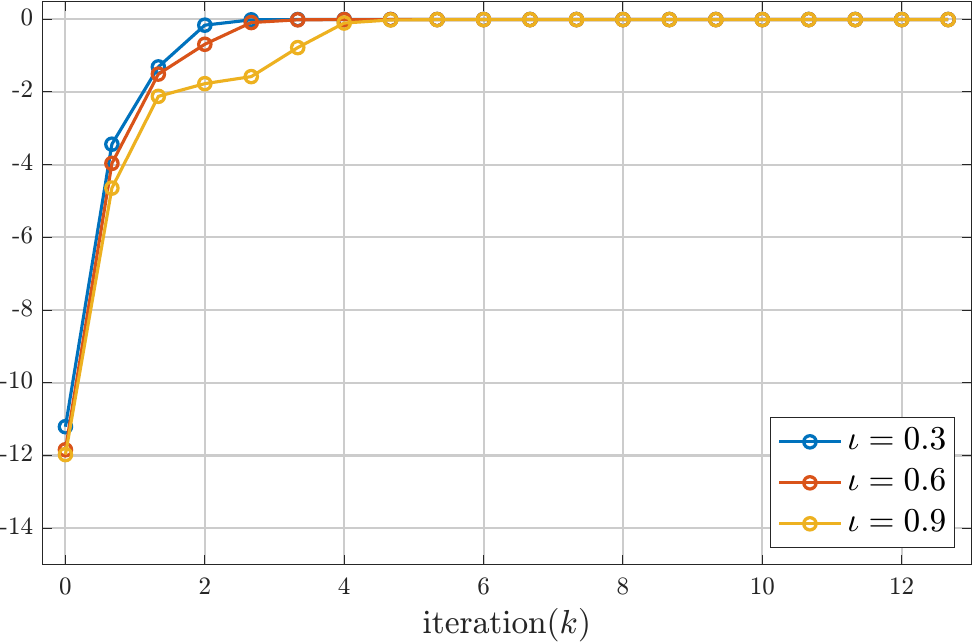}&&
\includegraphics[width=.275\textwidth,height=.175\textwidth]{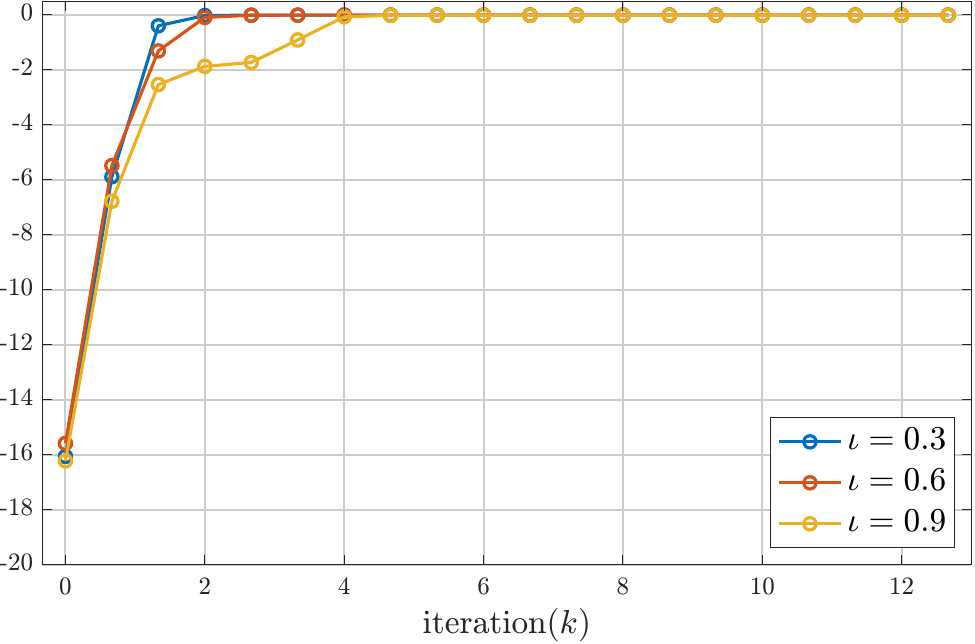}
\\
\end{tabular}

\begin{tabular}{p{\textwidth}}\scriptsize
The upper plots show the log-likelihood $\ell(\bm X_{nT};\wh{\bm{\varphi}}_n^{(k)})$ (blue line, left scale), and the convergence criterion $\Delta\ell_k$ defined in \eqref{eq:convEM} (red line, right scale) when we initialized the EM algorithm with the PC estimator. The lower plots show the percentage deviation of the log-likelihood of the algorithm initialized with the contaminated estimator ($\ell^c_k$) from the one initialized with the PC estimator ($\ell_k$), that is ${\ell^c_k}/{\ell_k}-1$. The bigger $\iota$ is, the stronger is the contamination. The log-likelihoods in this figure were obtained from a single simulation where $T=100$, $n=100$, $\mu=0.7$, $\delta=0.5$, $\tau=0.5$, and $\theta=0.5$.
\end{tabular}

\end{figure}

The upper plots of Figure \ref{fig:EM} show the log-likelihood $\ell(\bm X_{nT};\wh{\bm{\varphi}}_n^{(k)})$ (blue line, left scale) as a function of the iteration $k$ of the EM algorithm, and the convergence criterion $\Delta\ell_k$ defined in \eqref{eq:convEM} (red line, right scale) when we initialize the EM algorithm with the PC estimator. The results in Figure \ref{fig:EM} show that the log-likelihood is an increasing function in the number of iterations $k$ (as it should be). Moreover, the algorithm converges very fast: within two iterations $\Delta\ell_k\le 10^{-3}$. This is to be expected since we initialize the algorithm with the consistent PC estimator. 

The lower plots of Figure \ref{fig:EM} show the percentage deviation of the log-likelihood of the algorithm initialized with the contaminated estimator ($\ell^c_k$) from the one initialized with the PC estimator ($\ell_k$), that is $\ell^c_k/\ell_k-1$. 
Contaminating the initialization implies starting from a much lower likelihood. However, with just a few iterations, the log-likelihood initialized with the contaminated estimator is the same as the one properly initialized. This result shows that initializing the EM with a non-consistent estimator is fine, as it just requires running the algorithm for a few more iterations. Therefore, hereafter, we focus on the case in which we inizialize the EM algorithm with the PC estimator. 

\begin{figure}[t!]\caption{Simulation results - Common components}\label{fig:chi_rmse}
\centering

\footnotesize
{Root Mean Squared Errors}\\

\setlength{\tabcolsep}{.0\textwidth}
\begin{tabular}{@{}C{.5\textwidth}C{.5\textwidth}@{}}
\\[-4pt]
{Absolute} &  {Relative to PC}\\
\includegraphics[width=.4\textwidth]{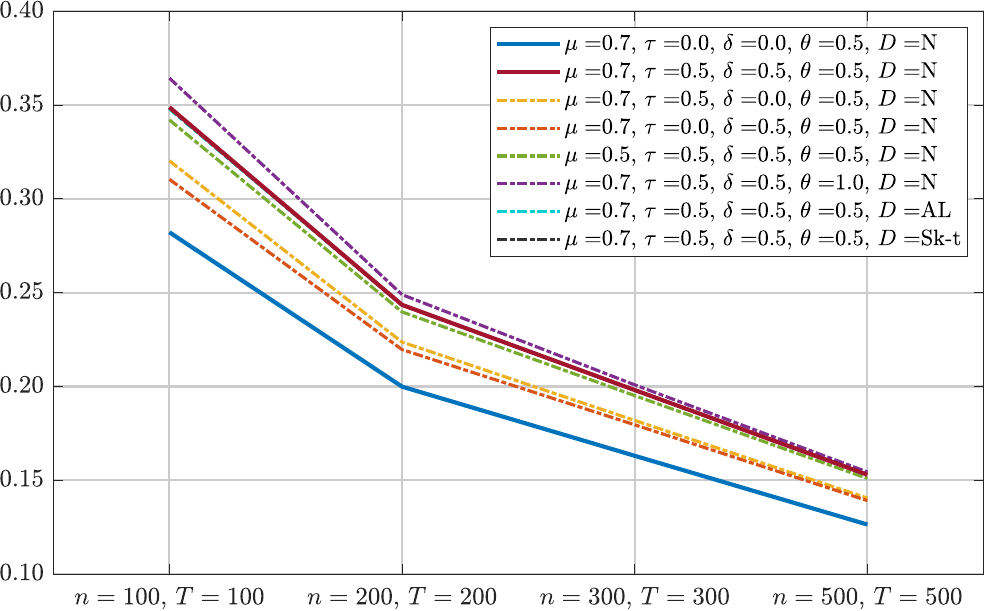}&
\includegraphics[width=.4\textwidth]{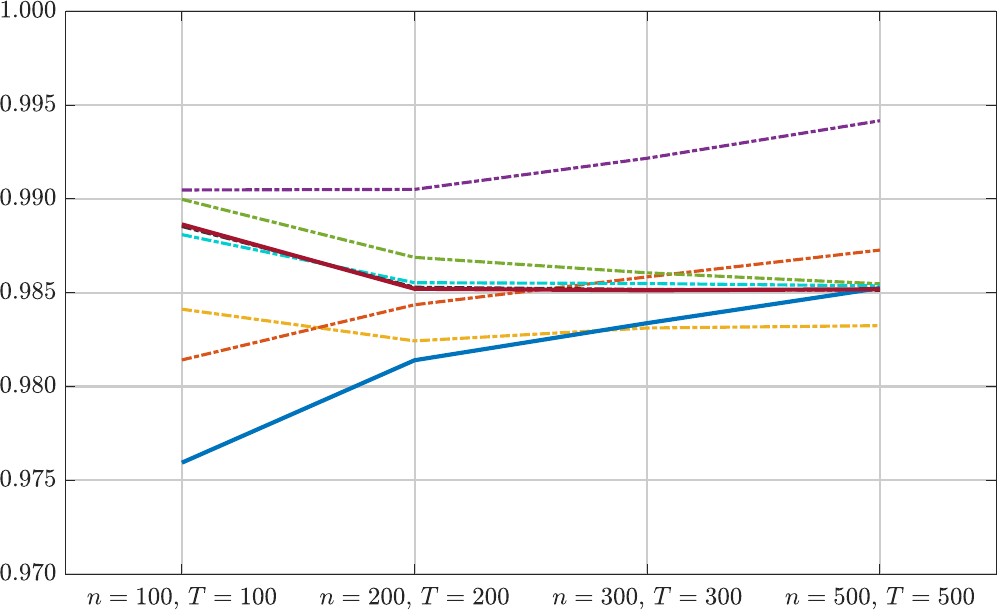}\\
\end{tabular}
\end{figure}

Moving to the performance of our estimators, we begin with Proposition \ref{th:chi} part (a), which gives consistency and rates for the common component's estimator. The left plot in Figure \ref{fig:chi_rmse} shows the root mean squared error:
\begin{align}
\mathrm{RMSE}=\sqrt{\frac{1}{nTB}\sum_{i=1}^n\sum_{t=1}^T\sum_{b=1}^B{(\wh{\chi}_{it}^{(b)}-\chi_{it}^{(b)})}^2}.\nn
\end{align}
Instead, the right plot shows the relative RMSE of our estimator over the RMSE of the PC estimator---values smaller than one indicate a better performance of our estimators.  We show results for several DGPs, which allow us to disentangle the effects of each single mis-specification on the performance of our estimator.

Two main results emerge from Figure \ref{fig:chi_rmse}. First, as $n$ and $T$ grow, the RMSE of all DGPs converge toward zero, thus indicating that the mis-specification introduced by estimating a model with uncorrelated and possibly non-Gaussian idiosyncratic components, even when that is not the case, does not affect our estimator. 
In particular, between serial correlation (the orange line) and cross-sectional correlation (the yellow line), the mis-specification that hurts the most is cross-sectional correlation. 
This is good news for practitioners because the cross-sectional correlation between the idiosyncratic components can somehow be limited by avoiding to include in the dataset variables that are too similar with one another (see, e.g., the discussions in \citealp{BN06} and \citealp{lucioADFM}). Moreover,  the model is consistently estimated when the shocks are asymmetric and have heavy tails even when they come from a distribution that does not meet Assumption \ref{ass:tails} of exponentially decaying tails of the distribution.\footnote{The dark gray line for the Skew-t distribution cannot be seen in the left plot because it is underneath, thus it coincides with, the red line.} This is also good news for practitioners because this result tells us that we can use this model in settings likely to be non-Gaussian, as is the case in datasets of disaggregated inflation rates (see, e.g., \citealp{reiswatson10}, and \citealp{AL}), which are notoriously skewed and fat-tailed.

Second, overall, our estimator behaves very similarly to but slightly better than the PC estimator, despite the latter being non-parametric and, thus, not affected by mis-specifications. This result confirms the conjectures based on extensive numerical studies made by \citet[]{DGRqml} and \citet[]{baili16}, showing that, in a high-dimensional setting, the EM estimator is ``as good as'' the PC estimator. 

To evaluate the estimates of the factors and the loadings, at each replication $b$, we consider a multivariate version of the $R^2$ (see also \citealp{DGRqml}): 
\begin{align}
\mathrm{TR}_F^{(b)}=\frac{\mbox{tr}\left(({\bm{\mathcal F}}_T^{(b)\prime}\wh{{\bm{\mathcal F}}}_T^{(b)})(\wh{{\bm{\mathcal F}}}_T^{(b)\prime}\wh{\bm{ \mathcal F}}_T^{(b)})^{-1}(\wh{\bm{\mathcal F}}_T^{(b)\prime}\bm{\mathcal F}_T^{(b)})\right)}{\mbox{tr}\left(\bm{\mathcal F}_T^{(b)\prime}\bm{\mathcal F}_T^{(b)}\right)}, \quad
\mathrm{TR}_\Lambda^{(b)}=\frac{\mbox{tr}\left((\bm \Lambda_n^{(b)\prime}\wh{\bm \Lambda}_n^{(b)})(\wh{\bm \Lambda}_n^{(b)\prime}\wh{\bm \Lambda}_n^{(b)})^{-1}(\wh{\bm \Lambda}_n^{(b)\prime}\bm \Lambda_n^{(b)})\right)}{\mbox{tr}\left(\bm \Lambda_n^{(b)\prime}\bm \Lambda_n^{(b)}\right)},\nonumber
\end{align}
where  $\wh{\bm{ \mathcal F}}_T^{(b)}$  and $\wh{\bm \Lambda}_n^{(b)}$ are the $T\times r$ and $n\times r$ matrices of the estimated factors and loadings, respectively. These trace statistics are smaller than one, and they tend to one when the empirical canonical correlations between the true quantities and their estimates tend to one.

\begin{figure}[t!]\caption{Simulation results - Factors and loadings}\label{fig:FL_tracce}
\centering

\footnotesize Trace statistics relative to PC\\

\setlength{\tabcolsep}{.0\textwidth}
\begin{tabular}{@{}C{.5\textwidth}C{.5\textwidth}@{}}
\\[-4pt]
{Factors} &  {Loadings}\\
\includegraphics[width=.4\textwidth]{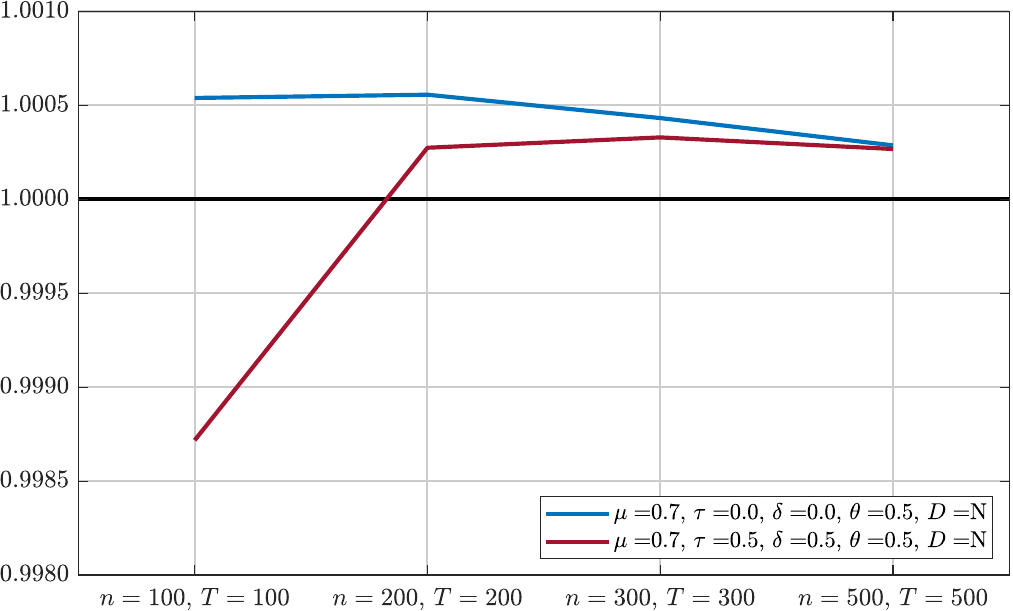}&
\includegraphics[width=.4\textwidth]{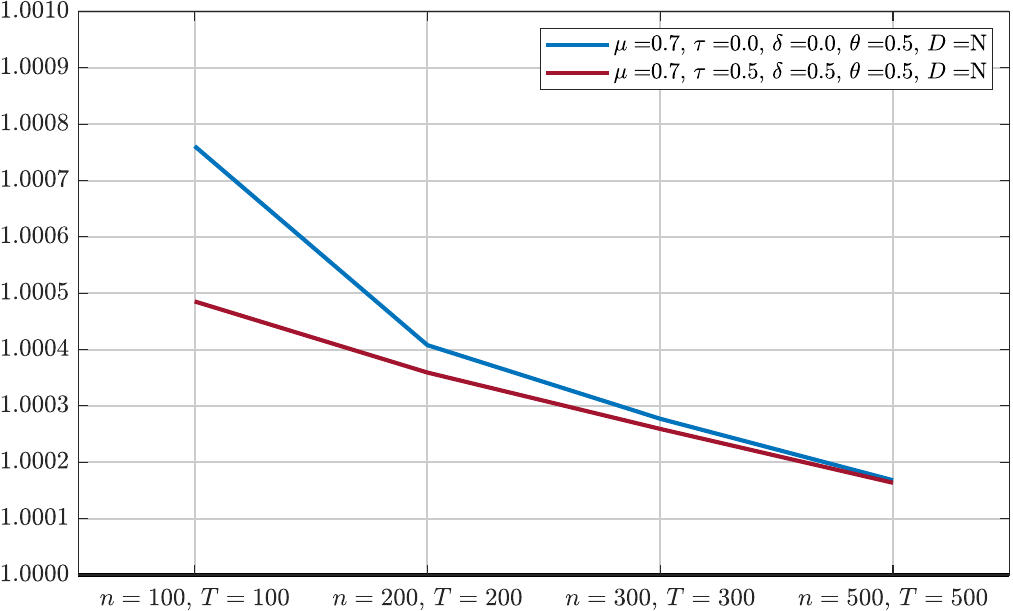}\\
\end{tabular}

\end{figure}

Figure \ref{fig:FL_tracce} reports the values of $\mathrm{TR}_F^{(b)}$ and $\mathrm{TR}_\Lambda^{(b)}$, relative to the same measures computed for the PC estimator, averaged over all $B$ replications (values larger than one  indicate a better performance of our estimators). The results in Figure \ref{fig:FL_tracce} mirror those of Figure \ref{fig:chi_rmse}: our estimator behaves very similarly to but slightly better than the PC estimator.

\begin{figure}[t!]\caption{Simulation results - Factors and loadings}\label{fig:FL_identifica}
\centering

\footnotesize Identification criteria\\

\setlength{\tabcolsep}{.0\textwidth}
\begin{tabular}{@{}C{.5\textwidth}C{.5\textwidth}@{}}
\\[-4pt]
{Factors - $B^{-1}\sum_{b=1}^BI_F^{(b)}$} &  {Loadings - $\log(B^{-1}\sum_{b=1}^BI_\Lambda^{(b)})$}\\
\includegraphics[width=.4\textwidth]{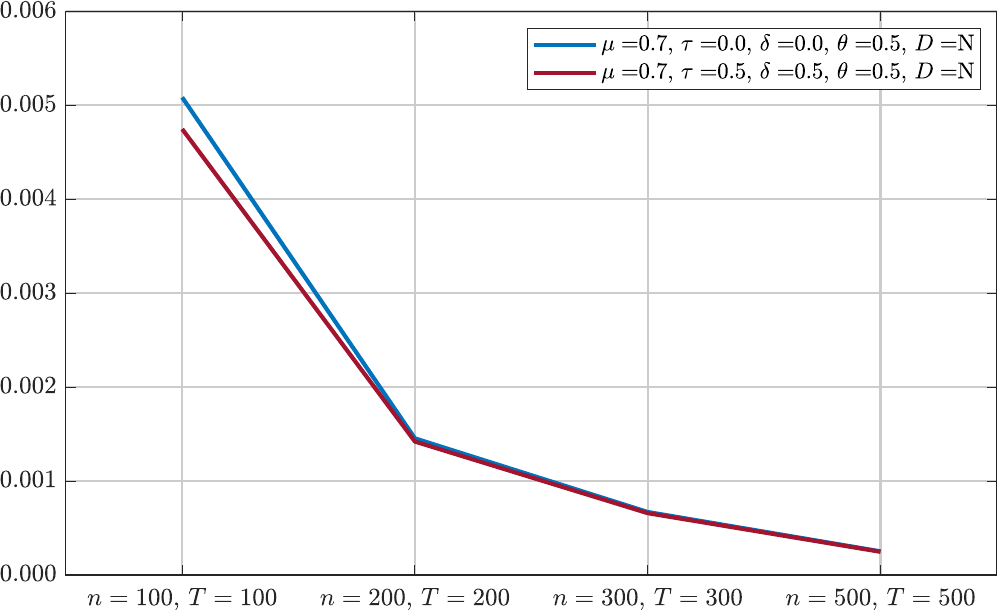}&
\includegraphics[width=.4\textwidth]{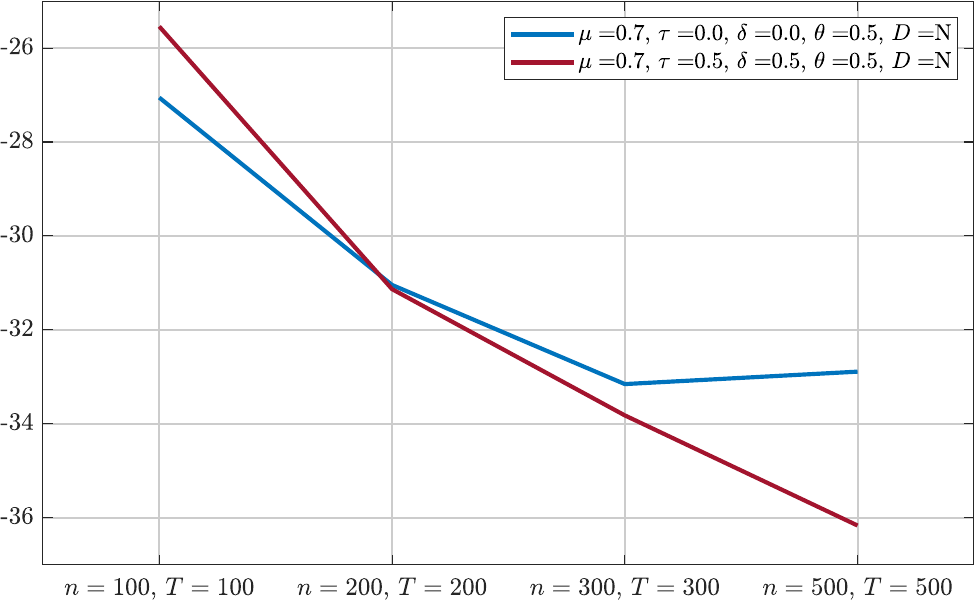}\\
\end{tabular}

\end{figure}

To derive our asymptotic results, we assumed that the true factors and the loadings are asymptotically identified as in Assumption \ref{ass:ident}(b). As explained in Remark \ref{rem:ident} the estimated factors and loading satisfy this assumption asymptotically. To verify that this is the case, in Figure \ref{fig:FL_identifica}, we show the two identification errors:
\begin{align}
I_F^{(b)} =\frac{1}{r}\sum_{j=1}^r \l(\nu^{(j)}(\wh{\bm \Gamma}^{F{(b)}})-1\r)^2, \qquad
I_\Lambda^{(b)}=\frac{1}{r}\sum_{j=1}^r \l(\nu^{(j)}(\wh{\bm \Sigma}_\Lambda^{(b)})-[\wh{\bm \Sigma}_\Lambda^{(b)}]_{jj}\r)^2, \nonumber
\end{align}
where 
$\wh{\bm \Gamma}^{F{(b)}}=T^{-1}\wh{\bm {\mathcal F}}_T^{(b)\prime} \wh{\bm {\mathcal F}}_T^{(b)}$, 
and $\wh{\bm \Sigma}_\Lambda^{(b)}=n^{-1}\wh{\bm \Lambda}_n^{(b)\prime} \wh{\bm \Lambda}_n^{(b)}$.
According to our asymptotic results, both these quantities should tend to zero as $n$ and $T$ increase because, given our simulation design, the true factors are orthonormal and the true loadings are orthogonal in agreement with Assumption \ref{ass:ident}(b). Figure \ref{fig:FL_identifica} shows $I_F^{(b)}$ and $I_\Lambda^{(b)}$ averaged over all $B$ replications. The identifying constraints are more and more precisely satisfied as $n$ and $T$ grow.

Next, we move to the asymptotic distribution of the common component. To this end, for each replication $b$ and any $i,t$, we compute
\begin{align}
Z_{it}^{(b)} &= \l(\frac 1 T\wh{\mbf F}_t^{(b)\prime}\wh{\bm{\mathcal V}}_{i}^{(\text{\tiny(HAC},b\text{\tiny)}}\wh{\mbf F}^{(b)}_t+\frac 1 n\wh{\bm\lambda}_i^{(b)\prime}\wh{\bm{\mathcal W}}_t^{(\text{\tiny(HAC},b\text{\tiny)}}\wh{\bm\lambda}^{(b)}_i\r)^{-1/2} (\wh{\chi}_{it}^{(b)}-{\chi}_{it}),\label{eq:ZitHAC}
\end{align}
where we use
 the robust estimators of the asymptotic covariance matrices defined in \eqref{eq:robcovloadiest} and \eqref{eq:robcovfacttest}, respectively. For comparison we also consider $Z_{it}^{(0,b)}$ defined as in \eqref{eq:ZitHAC}, but when using estimators of the asymptotic covariance matrices in the case of non-correlated idiosyncratic components (henceforth, non-robust covariance matrices).

\begin{figure}[t!]\caption{Simulation results - Histograms of $Z_{it}^{(b)}$}\label{fig:Hist}
\centering \footnotesize

Serially and cross-correlated idiosyncratic components ($\tau=0.5$, $\delta=0.5$) - Robust covariances\\

\begin{tabular}{ccc}
\\[-4pt]
\footnotesize $n=100$, $T=100$    & \footnotesize $n=200$, $T=200$ & \footnotesize $n=300$, $T=300$ \\
\footnotesize Gaussian & \footnotesize Gaussian & \footnotesize Gaussian \\
\includegraphics[width=.22\textwidth]{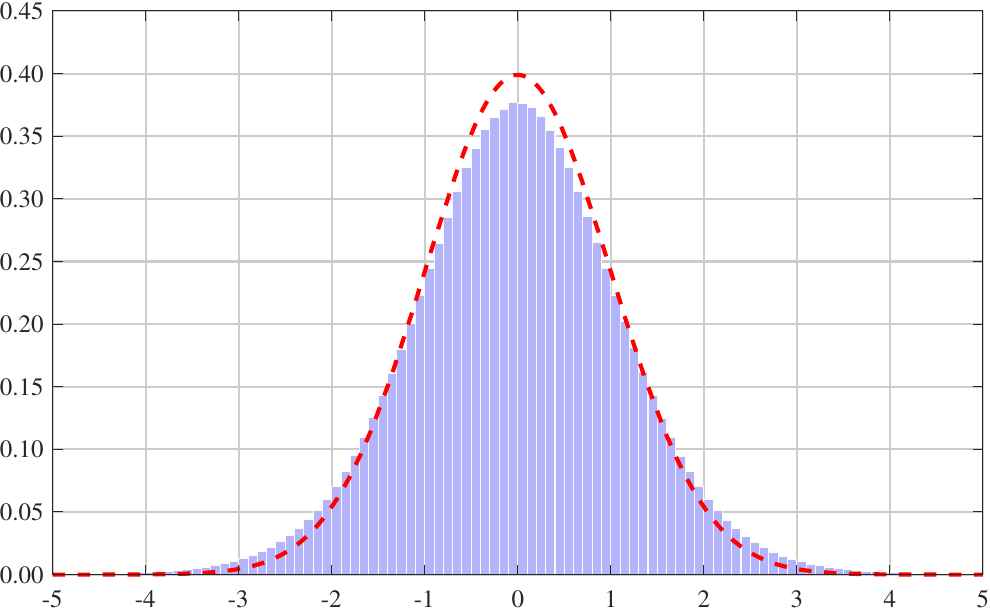}&
\includegraphics[width=.22\textwidth]{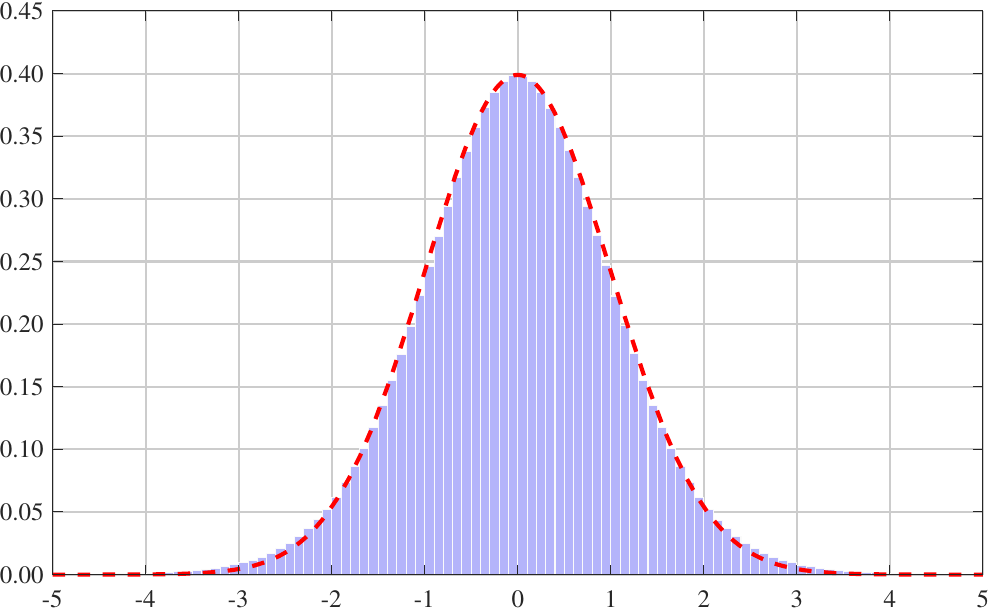}&
\includegraphics[width=.22\textwidth]{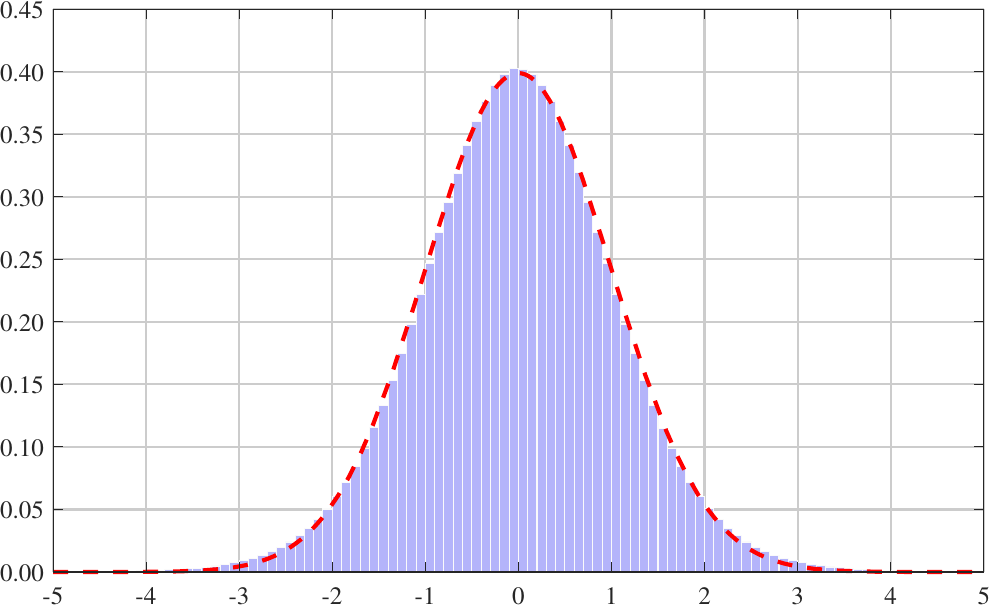}\\
\end{tabular}
\begin{tabular}{cccc}
\\[-4pt]
\footnotesize $n=200$, $T=200$ & \footnotesize $n=200$, $T=200$ & \footnotesize $n=300$, $T=300$ & \footnotesize $n=300$, $T=300$  \\
\footnotesize Asymmetric Laplace & \footnotesize Skew-$t$&\footnotesize Asymmetric Laplace & \footnotesize Skew-$t$\\
\includegraphics[width=.22\textwidth]{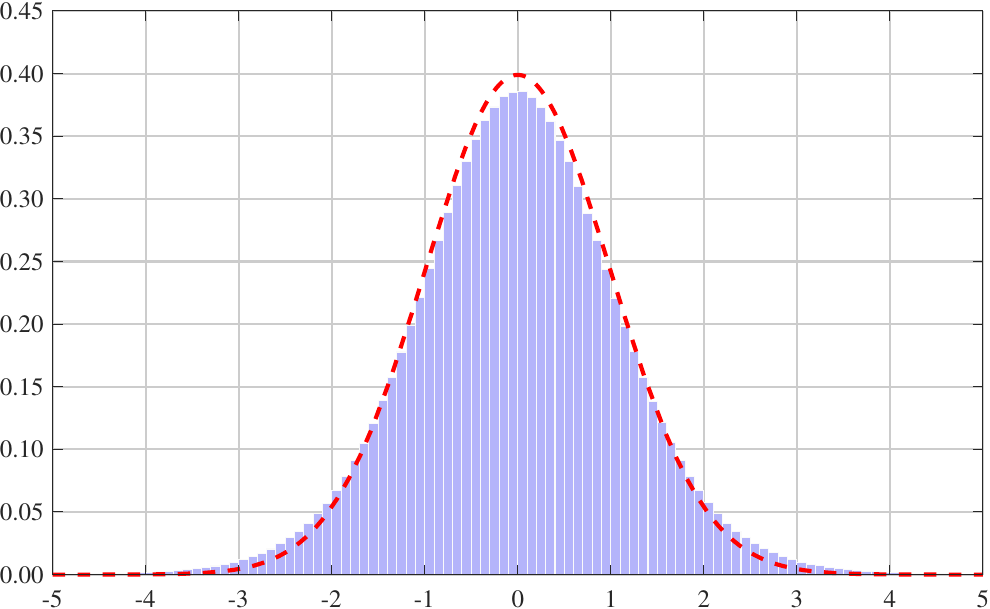}&
\includegraphics[width=.22\textwidth]{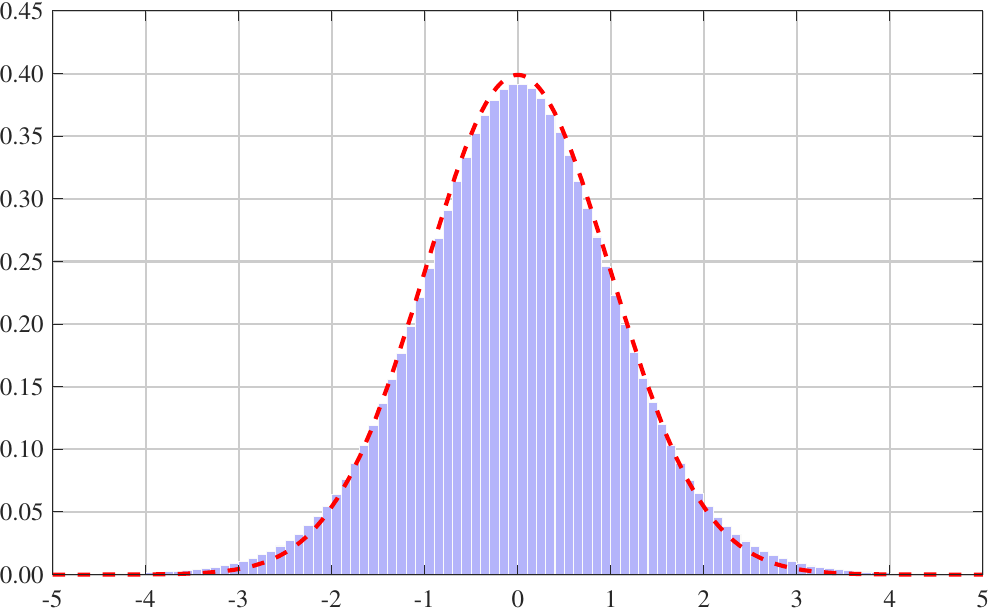}&
\includegraphics[width=.22\textwidth]{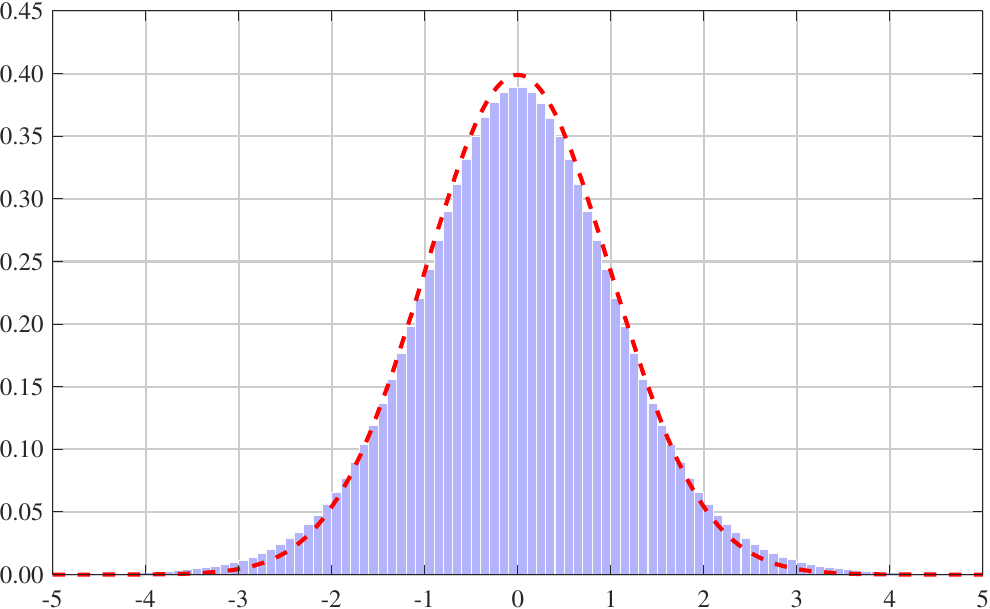}&
\includegraphics[width=.22\textwidth]{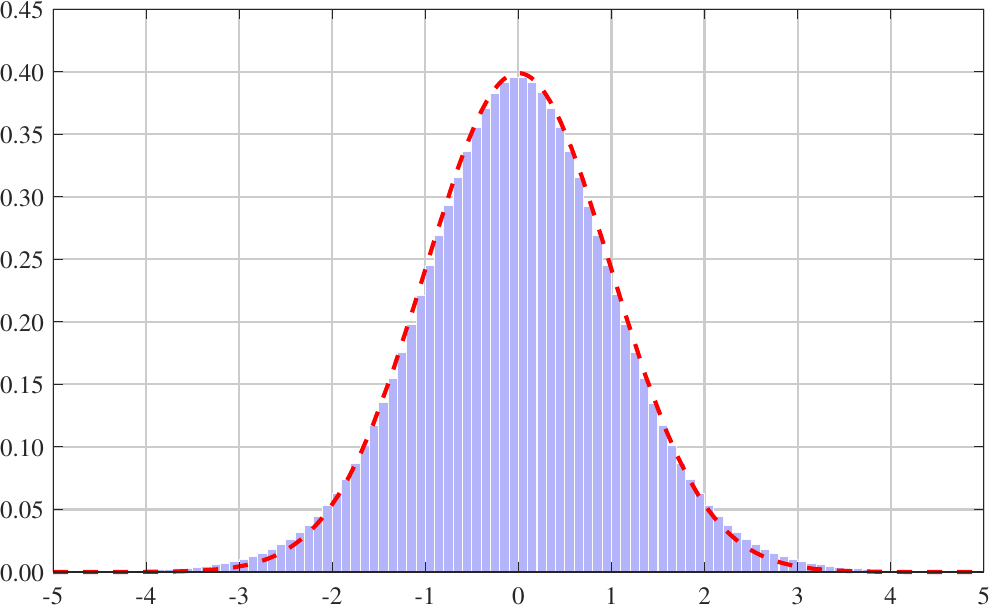}\\
\end{tabular}
\end{figure}

According to Proposition \ref{th:chi}, $Z_{it}^{(b)}\stackrel{d}{\to} \mathcal N(0,1)$ as $n,T\to\infty$. To evaluate the goodness of our theoretical results, we compute the average coverage
\[
C(1-\alpha)=\frac 1{nTB}\sum_{i=1}^n\sum_{t=1}^T\sum_{b=1}^B\, \mathbb I\l(\mathcal Z_{\alpha/2} \le Z_{it}^{(b)}\le \mathcal Z_{1-\alpha/2}\r),
\]
where $\mathcal Z_\alpha$ is the $\alpha$-quantile of the standard normal distribution. In Table \ref{tab:chi_coverage_ALL}, we report the coverage  $C(1-\alpha)$, for selected values of $\alpha\in(0,1)$, while for illustration purposes, in Figure \ref{fig:Hist}, we show histograms of $\{Z_{it}^{(b)}: i=1,\ldots, n,\, t=1,\ldots, T, \, b=1,\ldots, B\}$, for some of the cases considered in Table \ref{tab:chi_coverage_ALL}.  We stress that throughout this exercise the chosen bandwidths for all robust covariance estimators are not data driven, but rather fixed a priori.\footnote{To estimate $\wh{\bm{\mathcal V}}_i^{\text{\tiny(HAC)}}$ we set $M_T=\lfloor T^{1/4}\rfloor$ and to estimate $\wh{\bm{\mathcal W}}^{\text{\tiny(HAC)}}_t$ we set $m=\lfloor n^{4/5}\rfloor$.}  

Results confirm the derived asymptotic distribution. When the idiosyncratic components are uncorrelated, the non-robust estimators of the covariance matrices offers almost perfect coverage, while the robust estimators give a slight over-coverage. In the relevant cases of  serially and cross-correlated idiosyncratic components, the considered robust estimators work very well.
For comparison we show also results for 
$Z_{it}^{(b)}$ when $\chi_{it}$ is estimated with the PC estimator, and the asymptotic covariances are estimated using the estimators in \citet{baing06}. In this case, deviations from Gaussianity seem to lead to serious under-coverage. 

\begin{table}[t!]\caption{Simulation results - Common components - Coverage, $C(1-
\alpha)$}\label{tab:chi_coverage_ALL}
\centering
\setlength{\tabcolsep}{.0\textwidth}

\scriptsize
\begin{tabular}{@{}C{.3\textwidth}|C{.1\textwidth}C{.1\textwidth}|C{.1\textwidth}C{.1\textwidth}|C{.1\textwidth}C{.1\textwidth}} 
\hline\hline
&&&&&\\[-7pt]
							&&&\multicolumn{2}{c|}{$(1-\alpha)=0.90$}&\multicolumn{2}{c}{$(1-\alpha)=0.95$}\\	
							\cline{4-7}
							&&&&&\\[-6pt]
							&$n$ & $T$&EM & PC &EM &PC\\\hline
 							& 100 & 100 & 0.89& 0.88& 0.94& 0.93\\
Gaussian, $\tau=0$, $\delta=0$  	& 200 & 200 & 0.89& 0.89& 0.95& 0.94\\
Non-robust covariances			& 300 & 300 & 0.90& 0.89& 0.95& 0.95\\
            						& 500 & 500 & 0.90& 0.90& 0.95& 0.95\\\hline
					  		&  100   &  100   &  0.91&  0.89&  0.95&  0.94\\
Gaussian, $\tau=0$, $\delta=0$        &  200   &  200   &  0.92 &  0.90&  0.96&  0.95\\
Robust covariances           		&  300   &  300   &  0.92&  0.91&  0.96&  0.95\\
			                                 &  500   &  500    &  0.93 &  0.91 &  0.96 &  0.95\\\hline
		                                         &  100   &  100   &  0.86&  0.84 &  0.92 &  0.90\\
Gaussian, $\tau=0.5$, $\delta=0.5$   &  200   &  200   &  0.88 &  0.86 &  0.93&  0.92\\
Robust covariances 	             		&  300   &  300    &  0.89 &  0.87  &  0.94  &  0.93\\
			                              &  500   &  500   &  0.89&  0.88 &  0.94&  0.93\\\hline
			                              &  100   &  100   &  0.86 &  0.81 &  0.92 &  0.88\\
Asymmetric Laplace, $\tau=0.5$, $\delta=0.5$     &  200   &  200   &  0.88&  0.83 &  0.93 &  0.89\\
Robust covariances 	             		&  300   &  300   &  0.89 &  0.83 &  0.94 &  0.90\\
			                           &  500   &  500    &  0.89  &  0.84  &  0.94  &  0.90\\\hline
			                                                      &  100   &  100   &  0.86 &  0.81 &  0.92 &  0.88\\
Skew-$t$, $\tau=0.5$, $\delta=0.5$     &  200   &  200   &  0.88&  0.83&  0.93 &  0.89\\
Robust covariances 	             &  300   &  300   &  0.89&  0.83 &  0.94&  0.90\\
		                           &  500   &  500   &  0.89 &  0.84 &  0.94 &  0.90\\
		                           \hline\hline
\end{tabular}
\end{table}

Last, we move to Proposition \ref{prop:eff}, which states that under a sparsity condition on the covariance matrix of the idiosyncratic components, the Kalman smoother estimator of the factors is more efficient than the PC estimator. To verify this result, we compare the theoretical asymptotic covariance of the factors, $\bm{\mathcal W}_t$ and $\bm{\mathcal W}^{\text{\tiny\upshape PC}}_t$. Specifically, for each simulation, we look at the sign of the smallest eigenvalue of the matrix $(\bm{\mathcal W}^{\text{\tiny\upshape PC}}_t-\bm{\mathcal W}_t)$, computed using the true simulated parameters, which should always be positive. 

Table \ref{tab:efficiency} reports the percentage of times out of 5000 simulations in which $(\bm{\mathcal W}^{\text{\tiny\upshape PC}}_t-\bm{\mathcal W}_t)$ is positive definite. The DGP used for this exercise is the one described at the beginning of this section, but for the case indicated as $\tau=0.5^\ast$ in which we set $\Cov(e_{it},e_{jt})=0.5^{\vert i-j\vert}$ if $i,j\le \lfloor \sqrt{n}\rfloor$, and $\Cov(e_{it},e_{jt})=0$, otherwise, in order to better proxy the assumed sparsity condition. The results in Table \ref{tab:efficiency} confirm those in Proposition \ref{prop:eff}. When the sparsity condition is verified (this is the case of columns $\tau=0$ and $\tau=0.5^\ast$), the Kalman smoother is more efficient than the PC estimator. That said, even when the sparsity condition is not verified (column $\tau=0.5$), the Kalman smoother tends to be more efficient than the PC estimator. We also note that the two largest eigenvalues (not shown) are always positive, meaning that, in our simulations, the PC estimator is never more efficient than the EM estimator.

\begin{table}[t!]\caption{Percentage of simulations in which $(\bm{\mathcal W}^{\text{\tiny\upshape PC}}_t-\bm{\mathcal W}_t)$ is positive semidefinite }\label{tab:efficiency} \centering \scriptsize
\setlength{\tabcolsep}{.0\textwidth}
\begin{tabular}{@{}C{.15\textwidth}C{.15\textwidth}C{.15\textwidth}C{.15\textwidth}C{.15\textwidth}@{}}\hline\hline
$n$ &   $T$ &   $\tau=0$    &   $\tau=0.5^\ast$ &   $\tau=0.5$  \\ \hline
100 &   100 &   100.00  &   \bianco{1}97.74 &   \bianco{1}50.68 \\
200 &   200 &   100.00  &   \bianco{1}99.96 &   \bianco{1}85.84 \\
300 &   300 &   100.00  &   100.00  &   \bianco{1}96.82 \\
400 &   400 &   100.00  &   100.00  &   \bianco{1}99.10 \\
500 &   500 &   100.00  &   100.00  &   \bianco{1}99.60 \\
1000    &   1000    &   100.00  &   100.00  &   100.00  \\\hline\hline
\end{tabular}
\end{table}

%
%
\section{Empirical application}\label{sec:emp}
\begin{figure}[t!]\caption{US data - Convergence of the EM algorithm}\label{fig:EM_EMP}
\centering 
\footnotesize
\vskip .2cm
\includegraphics[width=.4\textwidth]{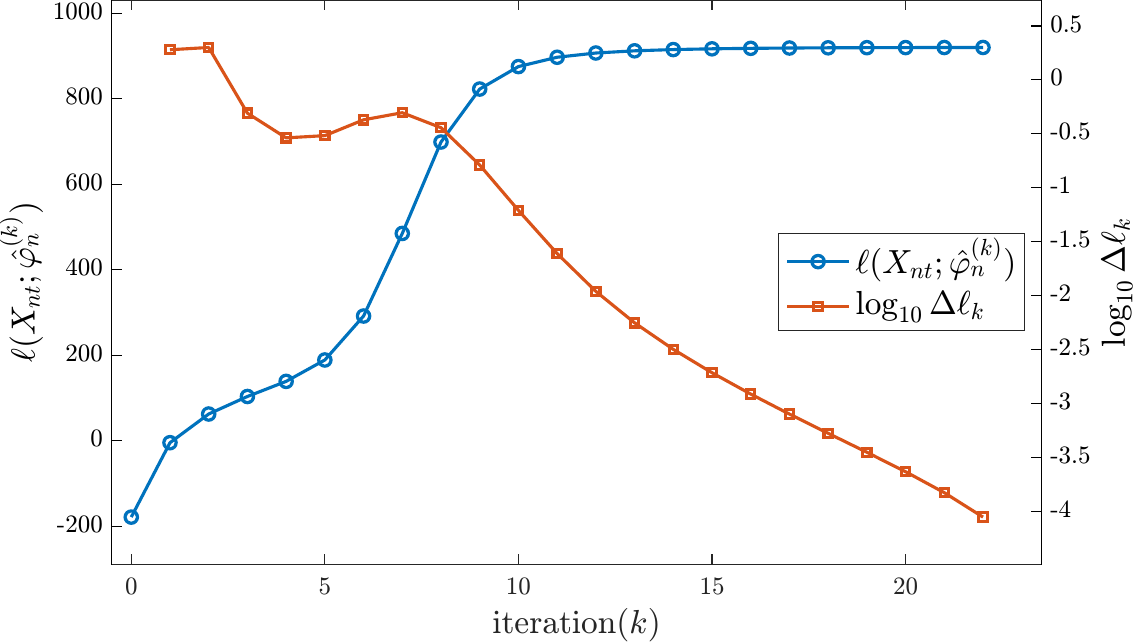}

\begin{tabular}{p{.4\textwidth}}\scriptsize
The plot shows the log-likelihood $\ell(\bm X_{nT};\wh{\bm{\varphi}}_n^{(k)})$ (blue line, left scale), and the convergence criterion $\Delta\ell_k$ defined in \eqref{eq:convEM} (red line, right scale) when we run the EM algorithm on US data and we initialize it with the PC estimator. 
\end{tabular}

\end{figure}

In this section, we consider the dataset used by \citet{OGAP}, a typical panel of $n=103$ US macroeconomic quarterly indicators observed from 1960:Q1 to 2018:Q4, thus $T = 236$. All variables are transformed to stationarity. In particular, we follow the common approach of taking first differences of price indexes and keeping interest rates in levels (see, e.g., \citealp{BBE05}). The information criterion by \citet{baing02} indicates $\wh r=6$ common factors. 
We estimate that the order of the VAR for the factors is $\wh p_F=2$.  
As shown in Figure \ref{fig:EM_EMP}, the  EM algorithm converges in 22 iterations when setting the threshold in \eqref{eq:convEM} to $\varepsilon=10^{-4}$ and the log-likelihood increases monotonically.

As we said in the Introduction, there is extensive literature showing the effectiveness of the EM algorithm in estimating large DFMs. Therefore, the purpose of this section is not to show that this method works or that it is superior to the PC estimator. Rather, we concentrate on the innovations brought about in this paper that have an impact on empirical applications, namely: the confidence bands for the common components and the factors, and a test of hypothesis on the factor loadings. 

Throughout, we compute the asymptotic covariances by using $\wh{\bm{\mathcal V}}_{i}^{\text{\tiny HAC}}$ for the loadings, as given in \eqref{eq:robcovloadiest}, with bandwidth $M_T=\lceil T^{1/4}\rceil$, and $\wh{\bm{\mathcal W}}_t^{\text{\tiny KF}}=n\wh{\bm \Pi}_{t|t}$ for the factors, computed as in \eqref{eq:dellefoche} 
by using the estimated parameters and applying local thresholding to the sample idiosyncratic covariance matrix \citep{FLM13,BR24}.

\begin{figure}[t!]\caption{US data - Estimated common components}\label{fig:chiest}
\centering
\footnotesize
\vskip .2cm
\setlength{\tabcolsep}{.0\textwidth}
\begin{tabular}{@{}C{.025\textwidth}C{.32\textwidth}C{.32\textwidth}C{.32\textwidth}}
& \footnotesize GDP growth rate & \footnotesize Core CPI inflation & \footnotesize Fed funds rate \\
\rotatebox{90}{\hskip .5cm  \footnotesize EM} &
\includegraphics[width=.31\textwidth]{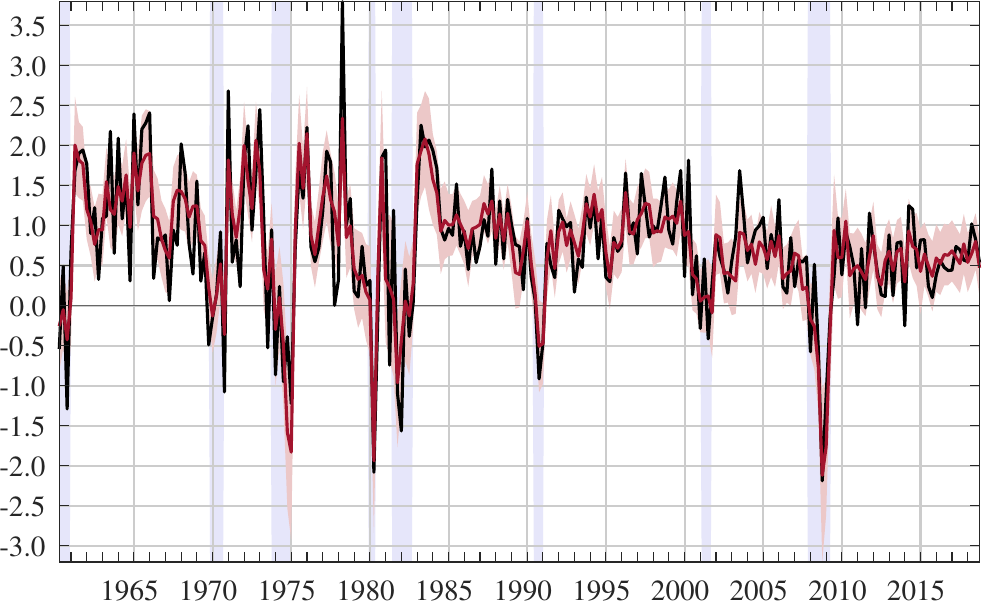}& 
\includegraphics[width=.31\textwidth]{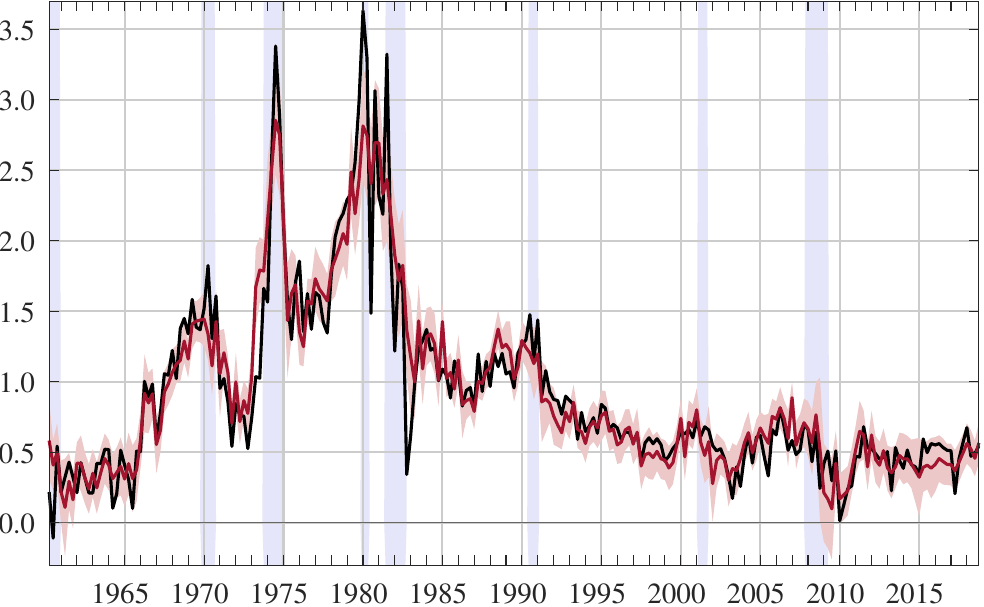}&
\includegraphics[width=.31\textwidth]{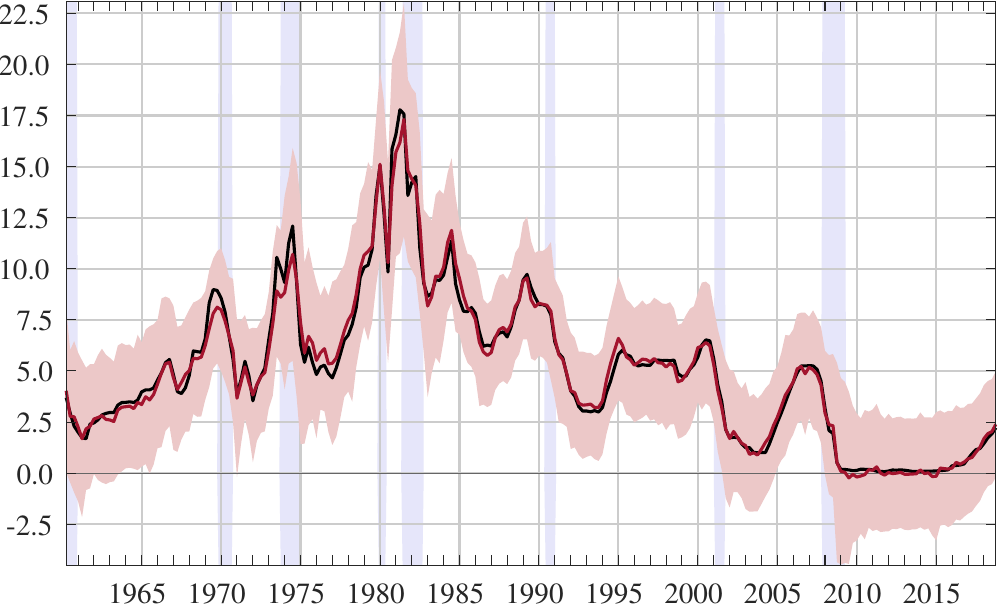}\\
\rotatebox{90}{\hskip .5cm  \footnotesize PC} &
\includegraphics[width=.31\textwidth]{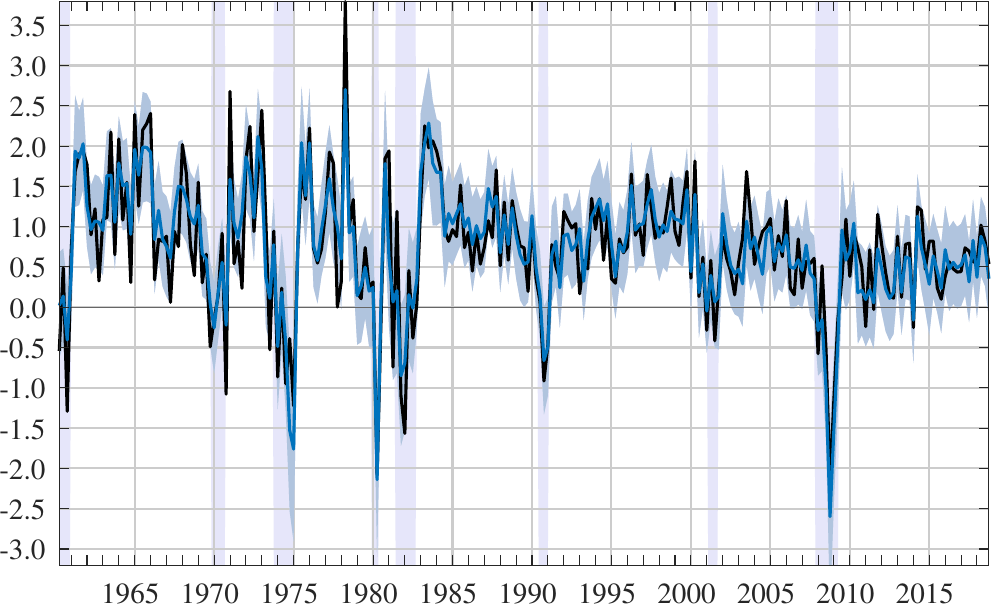}&
\includegraphics[width=.31\textwidth]{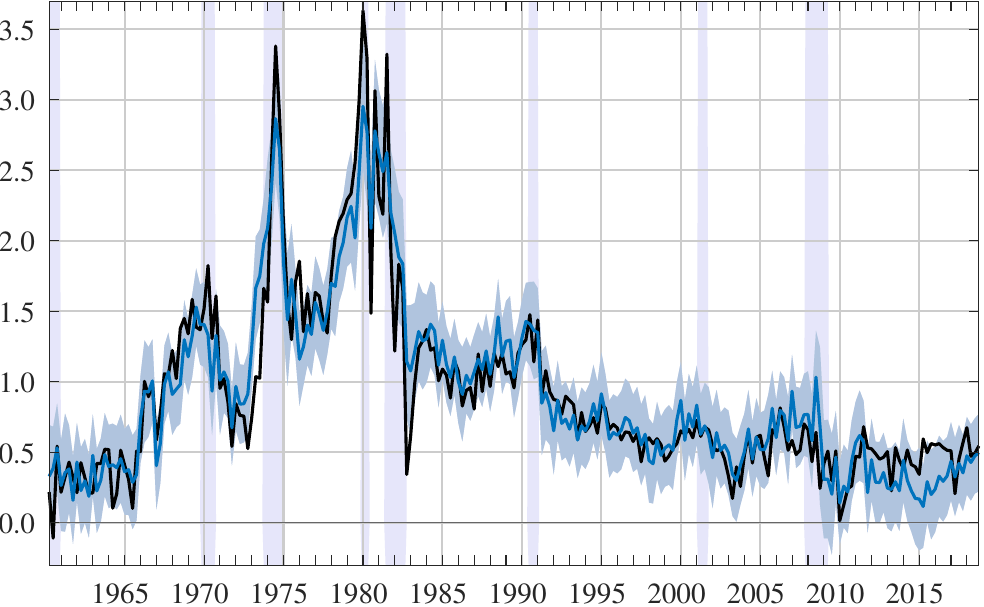}&
\includegraphics[width=.31\textwidth]{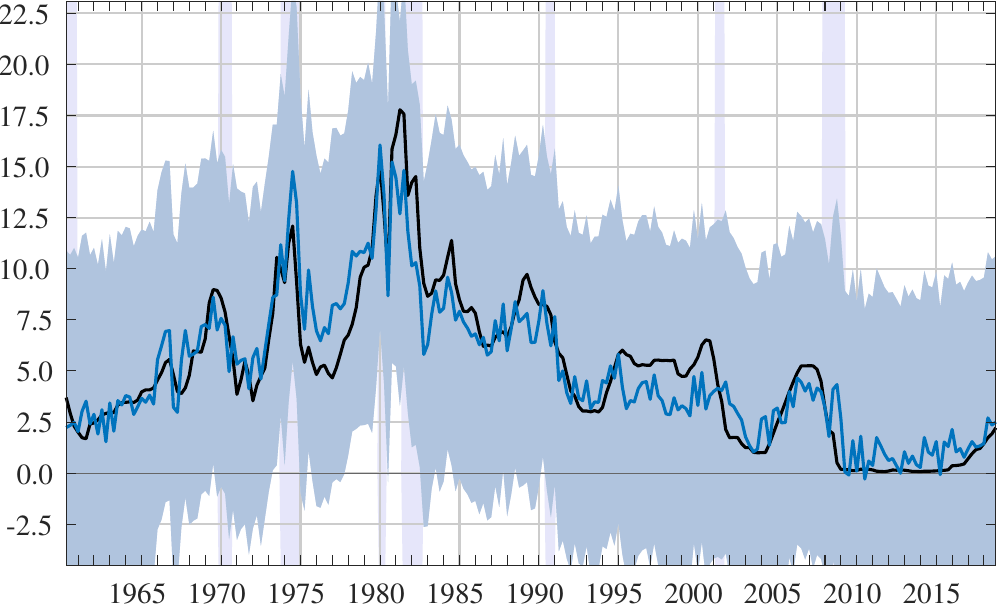}\\
\end{tabular}

\begin{tabular}{p{.95\textwidth}} \scriptsize
The shaded area is the 95\% confidence band.
\end{tabular}
\end{figure}

The first row of Figure \ref{fig:chiest} shows the common components of a few variables of interest estimated with the EM algorithm (the red line) with their 95\% confidence bands, together with the observed data (the black line). Specifically, for any given $i=1,\ldots, n$ and all $t=1,\ldots, T$, a  $(1-\alpha)\%$ confidence interval for the estimated common component is given by
\begin{align}
 {\mathcal I}_{\wh{\chi}_{it}}(\alpha)
=\l[
\wh{\chi}_{it}-z_{(1-\alpha/(2T))}\sqrt{\frac {\wh{\mbf F}_t^\prime\wh{\bm{\mathcal V}}_{i}^{\text{\tiny HAC}}\wh{\mbf F}_t}T
+\frac{\wh{\bm\lambda}_i^\prime\wh{\bm{\mathcal W}}_t^{\text{\tiny KF}}\wh{\bm\lambda}_i}n}
\;,\;
\wh{\chi}_{it}+z_{(1-\alpha/(2T))}\sqrt{\frac {\wh{\mbf F}_t^\prime\wh{\bm{\mathcal V}}_{i}^{\text{\tiny HAC}}\wh{\mbf F}_t}T
+\frac{\wh{\bm\lambda}_i^\prime\wh{\bm{\mathcal W}}_t^{\text{\tiny KF}}\wh{\bm\lambda}_i}n}
\,
\r].\nn
\end{align}
Moreover, in order to have confidence bands valid for all $T$ observations we apply a Bonferroni correction to the critical values. For comparison, in the second row of Figure \ref{fig:chiest}, we report the estimated common components obtained with the PC estimator (the blue line) with their 95\% confidence bands computed using the HAC estimators in \citet{baing06}. 

The common components of core CPI inflation and the Fed funds rate estimated with the EM track the observed series better than those estimated with the PC estimator. This result is possibly due to local departures from stationarity in those series, which create problems for PC analysis but not for the EM, as the Kalman smoother is able to track changes in the dynamics due to its recursive character.
In particular, the confidence band of the common component of the Fed funds rate for the PC estimator is much wider than that of the EM estimator because the variance of the idiosyncratic component obtained with the PC estimator is nearly eight times larger than that of the EM estimator and it is also far more persistent. 

In addition to computing confidence bands for the in-sample estimate of the common component, we can also do so for the unconditional and conditional forecasts obtained from the model. In other words, our results open the possibility of computing the uncertainty around GDP nowcast, which is nothing else than a short-term conditional forecast, and for scenario analysis performed with DFMs. Here we give a couple of simplified examples.

\begin{figure}[t!]\caption{US data - GDP growth forecast}\label{fig:chifore}
\centering
\footnotesize
\vskip .2cm
\setlength{\tabcolsep}{.0\textwidth}
\begin{tabular}{@{}C{.5\textwidth}C{.5\textwidth}}
Forecast of 2018:Q1 & Forecast through 2018:Q4\\
\includegraphics[width=.4\textwidth]{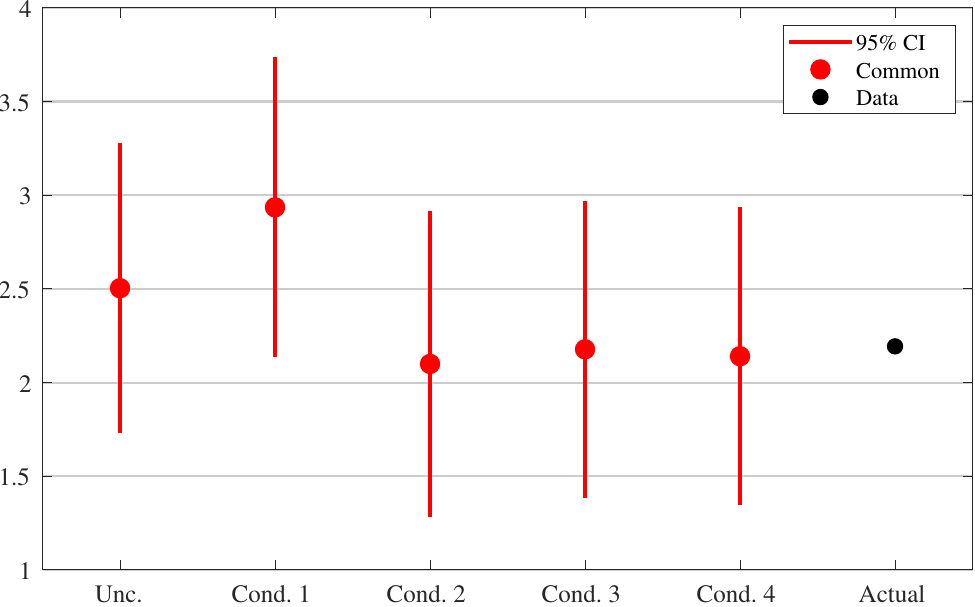}& 
\includegraphics[width=.4\textwidth]{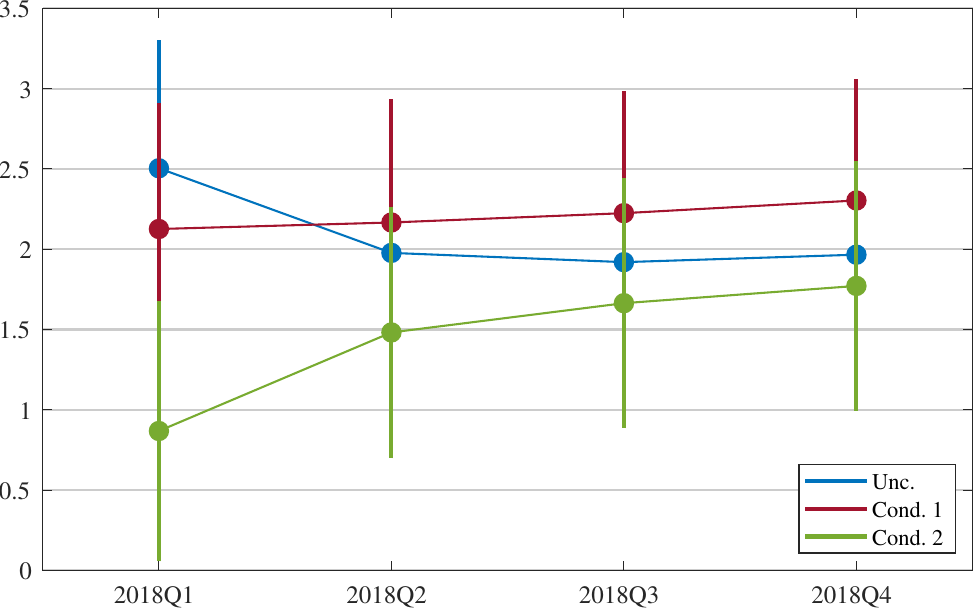}\\
\end{tabular}

\begin{tabular}{p{\textwidth}} \scriptsize
In the left chart ``Unc.'' is the unconditional forecast given data up to 2017:Q4; ``Cond. 1'' is the forecast conditioning on the value in 2018:Q1 of all the high-frequency variables in our dataset (stock prices, oil prices, surveys, and interest rates); ``Cond. 2'' is the forecast conditioning on the value in 2018:Q1 of all the high-frequency indicators the labor market indicators published in the BLS employment report; ``Cond. 3'' is the forecast conditional also on the value in 2018:Q1 of the CPI and PPI data, housing market indicators, and industrial production; ``Cond. 4'' is the common component estimated given all the data up to 2018:Q1; and, lastly, ``Actual'' is the actual value of GDP growth published by the BEA. All these indicators are released prior to GDP, which is released at the end of the following month: high-frequency indicators are available almost in real time, the labor report is published in the first week of the following month, and the CPI, PPI, Industrial Production, and housing indicators, are published about the mid of the following month.  \\\scriptsize
In the right chart ``Unc.'' is the unconditional forecast given data up to 2017:Q4; ``Cond. 1'' is the forecast conditional on payroll employment growing at the same average pace as in 2017; ``Cond. 2'' is the pessimistic scenario assuming that payroll employment grows at half the average pace as in 2017.
\end{tabular}
\end{figure}

The left chart of Figure \ref{fig:chifore} shows the 1-step ahead forecast of 2018:Q1. That is, we estimate the model up to 2017:Q4, and then we produce four different forecasts of GDP growth conditioning on the observations of an increasing number of variables for 2018:Q1. 
The conditional forecasts are obtained using the Kalman filter estimator of the factors, $\mbf F_{t|t}^{(k^*+1)}$, while
the unconditional forecast is obtained using their one-step-ahead prediction, $\mbf F_{t|t-1}^{(k^*+1)}=\wh{\mbf A}\mbf F_{t-1|t-1}^{(k^*+1)}$. This exercise mimics a simplified nowcasting setting---we are omitting the aspect of mixed frequency---as the variables we are conditioning on are published earlier than GDP, and the sequence of conditioning mimics the calendar of data releases. As shown in the left chart of Figure \ref{fig:chifore}, the model adjusts the forecast in the right direction as more hard data becomes available.

In the second exercise, we produce forecasts for GDP growth for each quarter of 2018 (blue line). Then, we adjust them based on different scenarios for payroll employment. We have two scenarios: one where employment grows at the same pace as the previous year (red line)---this is a slightly lower average pace than the model expected; another where employment grows at half the previous year's pace, resulting in a more pessimistic forecast (green line). As shown in the right chart of Figure \ref{fig:chifore}, one and two quarters ahead, the forecasts from these scenarios differ significantly.

Finally, we can  test for linear restrictions on the loadings. Consider testing for $s$ linear restrictions
$\text H_0:\; \bm R^\prime\,\text{\upshape vec}({\bm\Lambda}_{n})=\bm q$, against the alternative $\text H_1:\;\bm R^\prime\,\text{\upshape vec}({\bm\Lambda}_{n})\ne \bm q$,
where $\bm R$ is $nr\times s$ and $\bm q$ is $s\times 1$. Then,  the usual Wald-type test statistic is computed as
\beq\label{eq:wald}
{\mathsf W}_{\wh{\bm\Lambda}_n}=T \; \l( \bm R^\prime\,\text{\upshape vec}(\wh{\bm\Lambda}_{n})-\bm q\r)^\prime
\l(
\bm R^\prime \, 
\wh{\bm{\mathcal V}}_n^{\text{\tiny HAC}}
\bm R
\r)^{-1}
\l( \bm R^\prime\,\text{\upshape vec}(\wh{\bm\Lambda}_{n})-\bm q\r)
,
\eeq
where $\bm R$ selects only $sr$ rows of $\wh{\bm{\mathcal V}}_n^{\text{\tiny HAC}}$.  Under $\text H_0$, from Proposition \ref{prop:load} it follows that ${\mathsf W}_{\wh{\bm\Lambda}_n}\stackrel{d}{\to}\chi^2_{(r)}$, as $n,T\to\infty$. This test is the analogous of the test derived in the case of PC estimation by \citet{baing13}. Testing for equal loadings is equivalent to testing for equal common components for all $t=,1\ldots, T$.


Table \ref{tab:test} shows the result of the test of five different null hypotheses. The first hypothesis we test (column (A)) is whether GDP and GDI, both measures of US aggregate output, have equal loadings and, as a consequence, have the same common component. Our test does not reject the hypothesis of equal loadings. This result supports the idea recently explored in the literature that combining GDP and GDI can better estimate aggregate output \citep[e.g.,][]{ADNSS,BLgdo}.

\begin{table}[t!]
\caption{US data - Testing hypothesis on the loadings}\label{tab:test}
\centering

\scriptsize
\vskip 2pt
\setlength{\tabcolsep}{.0\textwidth}

\begin{tabular}{ L{.1\textwidth}  L{.1\textwidth}  C{.1\textwidth} C{.1\textwidth} C{.1\textwidth} C{.1\textwidth} C{.1\textwidth} C{.1\textwidth}}\hline\hline
	&&&&&\\[-8pt]										
	&& 	(A)	& 	(B)	& 	(C)	& 	(D)	& 	(E)	& 	(F)	\\[-8pt]
	&&&&&\\ \hline												
	&&&&&\\[-8pt] 												
$r=1$ &${\mathsf W}_{\wh{\bm\Lambda}_n}$	& 	0.1947&	 0.0087&	0.0736&	0.1358 &	1.7421&	9.122	\\
	  &$p$-value	& 0.66&	0.93&	0.79&	0.71&	0.19&	0.00		\\[-8pt]
	&&&&&\\				
		\hline					
$r=2$ &${\mathsf W}_{\wh{\bm\Lambda}_n}$	&0.2203&	0.3708&	0.5392&	0.6124&	3.0266&	11.7768 		\\
	  &$p$-value& 0.90&	0.83&	0.76&	0.74&	0.22&	0.00			\\[-8pt]			
	&&&&&\\ \hline
$r=3$ &${\mathsf W}_{\wh{\bm\Lambda}_n}$	&0.2462&	4.6652&	1.5638&	0.8066&	5.3210&	13.6319 		\\
	  &$p$-value& 	0.97&	0.20&	0.67&	0.85&	0.15&	0.00	\\[-8pt]			
	&&&&&\\ \hline
$r=4$ &${\mathsf W}_{\wh{\bm\Lambda}_n}$	&0.3803&	9.040&	7.8699&	1.2168&	5.1899&	110.8475	\\
	  &$p$-value& 	0.98&	0.06&	0.10&	0.88&	0.27&	0.00	\\[-8pt]			
	&&&&&\\ \hline
$r=5$& ${\mathsf W}_{\wh{\bm\Lambda}_n}$	& 1.8934&	14.7642&	9.7824&	1.6959&	5.4872&	128.9513	\\
	  &$p$-value& 0.86&	0.01&	0.08&	0.89&	0.36&	0.00 	\\[-8pt]			
	&&&&&\\ \hline
$r=6$& ${\mathsf W}_{\wh{\bm\Lambda}_n}$	& 	1.2704&	28.6640&	14.4703&	4.7851&	11.9874&	126.7982	\\
	  &$p$-value	& 0.97&	0.00&	0.02&	0.57&	0.06&	0.00		\\[-8pt]			
	&&&&&\\ \hline
	\hline		
\end{tabular}

\begin{tabular}{p{.8\textwidth}} \scriptsize
The null hypotheses are: (A) $\bm\lambda_{\text{GDP}} =\bm \lambda_{\text{GDI}}$;
(B) $\bm\lambda_{\text{CPI}} =\bm \lambda_{\text{PCE}}$;
(C) $\bm\lambda_{\text{CPI}_\text{core}} =\bm \lambda_{\text{PCE}_\text{core}}$;
(D) $\bm\lambda_{\text{CPI}_\text{energy}} =\bm \lambda_{\text{PCE}_\text{energy}}$; 
(E) $\bm\lambda_{\text{CPI}_\text{food}} =\bm \lambda_{\text{PCE}_\text{food}}$; and, 
(F) $\bm\lambda_{\text{GDP}} =\bm \lambda_{\text{Payroll}}$.
\end{tabular}
\end{table}

The second hypothesis (column (B)) that we test is whether CPI inflation and PCE price inflation, which are two alternative measures of consumer price inflation, have equal loadings. These indexes usually differ because they are constructed differently.\footnote{The CPI, which captures the headlines in newspapers, determines the return on Treasury Inflation-Protected Securities, or TIPS, while the inflation objective of the Federal Reserve is specified in terms of PCE price inflation.}
The test does not reject the null whenever $r<6$, which we read as signaling that, indeed, CPI inflation and PCE price inflation respond in the same way to the common factors, and thus, their difference is just idiosyncratic or, perhaps, weak/local factors. Columns (C), (D), and (E) in Table 5 test whether the difference between the PCE and CPI core sub-index, energy sub-index, and food sub-index are just idiosyncratic. The test suggests this to be the case.

Finally, in column (F) of Table \ref{tab:test}, we verify that if we test a non-sense hypothesis, the test rejects it. Specifically, we test whether the loadings of GDP growth and the growth in total non-farm employment are the same, which should not be the case. The test unequivocally reaches the same conclusion.

%
%
\section{Concluding remarks}\label{sec:conclusions}
This paper provides the asymptotic properties of Quasi Maximum Likelihood (QML) estimation for large approximate dynamic factor models, implemented via the Kalman smoother  and the Expectation Maximization (EM) algorithm. Our results provide the statistical foundations of one of the most popular and successful methods for estimating high-dimensional factor models for time series commonly used in many public and private institutions to track and predict economic activity.

From a  technical point of view, we show and prove that the EM approach is feasible even in the high-dimensional case, i.e., when the cross-sectional size $n$ can be much larger than the sample size $T$,  a point also made by \citet{DGRqml}. Then, we show that the EM estimator converges at the same rate as the Principal Components (PC) estimator. Moreover, we show that the EM  estimator of the loadings is always as efficient as the PC estimator, while the Kalman smoother estimator of the factors is more efficient than the PC estimator if the idiosyncratic covariance is sparse enough. 

Compared to the standard PC estimator, the EM approach has the main advantage of allowing the user to easily impose restrictions that reflect any prior knowledge about the data on the model. The user can impose these restrictions because the state-space formulation and the Kalman Smoother allow explicit modeling and estimation of the dynamic evolution of the latent factors and deal with data irregularly spaced in time. Moreover, in the M-step, the user can impose restrictions on the parameters, thus allowing for constrained QML estimation. In contrast, the user cannot use the PC estimator to model the latent factors's dynamics. In an application on a dataset of US macroeconomic time series, we show that the EM algorithm can produce estimates of the common component that track the dynamics of the observed series better than the PC estimator, especially those series displaying periods of high persistence and regime changes, like inflation and interest rates. This result suggests that the Kalman smoother might be more robust to local deviations from stationarity, a feature already highlighted by \citet{kalman60} and \citet{kalmanbucy61}.


\clearpage

\pagestyle{fancy}
\fancyhf{}
\chead{Supplementary material for the paper: QMLE of Large Approximate DFM via the EM algorithm}
\cfoot{Page \thepage}
 \renewcommand{\footrulewidth}{0.2pt}
 \renewcommand{\headrulewidth}{0.2pt}

\setcounter{table}{0}
\setcounter{figure}{0}
\setcounter{page}{1}

\begin{center}
\textit{Supplementary material for the paper} \\[6pt]

\Large{\sc Quasi Maximum Likelihood Estimation and Inference\\ \vskip .2cm of Large Approximate Dynamic Factor Models\\ \vskip .2cm
 via the EM algorithm} \\[12pt]

\begin{tabular}{cp{1cm}c}
 Matteo Barigozzi &&  Matteo Luciani\\
\footnotesize Universit\`a di Bologna &&\footnotesize Federal Reserve Board\\[-.2cm]
\scriptsize  matteo.barigozzi@unibo.it &&\scriptsize matteo.luciani@frb.gov\\[.5cm]
\end{tabular}

\end{center}

\renewcommand{\thefootnote}{ } 

\footnotetext{
M. Barigozzi gratefully acknowledges financial support from MIUR (PRIN2020, Grant 2020N9YFFE).\\[.03cm]

\noindent
Disclaimer: the views expressed in this paper are those of the authors and do not necessarily reflect the views and policies of the Board of Governors or the Federal Reserve System.} 

\gdef\thefootnote{(\roman{footnote})}

\small

\section*{Table of contents}
\contentsline {section}{\numberline {A}Further details on estimation}{2}{appendix.A}%
\contentsline {subsection}{\numberline {A.1}Principal Component estimators}{2}{subsection.A.1}%
\contentsline {subsection}{\numberline {A.2}Kalman filter and smoother}{3}{subsection.A.2}%
\contentsline {subsection}{\numberline {A.3}Stopping rule for the EM algorithm}{4}{subsection.A.3}%
\contentsline {section}{\numberline {B}Proof of main results}{4}{appendix.B}%
\contentsline {subsection}{\numberline {B.1}Proof of Proposition \ref {prop:plainvanilla}}{4}{subsection.B.1}%
\contentsline {subsection}{\numberline {B.2}Proof of Proposition \ref {prop:load}}{8}{subsection.B.2}%
\contentsline {subsection}{\numberline {B.3}Proof of Proposition \ref {prop:factors}}{9}{subsection.B.3}%
\contentsline {subsection}{\numberline {B.4}Proof of Proposition \ref {th:chi}}{14}{subsection.B.4}%
\contentsline {subsection}{\numberline {B.5}Proof of Proposition \ref {prop:altri}}{15}{subsection.B.5}%
\contentsline {subsection}{\numberline {B.6}Proof of Corollary \ref {cor:matriciAN}}{18}{subsection.B.6}%
\contentsline {subsection}{\numberline {B.7}Proof of Proposition \ref {prop:eff}}{19}{subsection.B.7}%
\contentsline {section}{\numberline {C}General lemmas}{22}{appendix.C}%
\contentsline {section}{\numberline {D}Lemmas necessary for proving Proposition \ref {prop:plainvanilla}}{31}{appendix.D}%
\contentsline {section}{\numberline {E}Lemmas necessary for proving Proposition \ref {prop:load}}{47}{appendix.E}%
\contentsline {section}{\numberline {F}Lemmas necessary for proving Proposition \ref {prop:factors}}{78}{appendix.F}%
\contentsline {section}{\numberline {G}Lemmas necessary for proving Proposition \ref {prop:altri}}{80}{appendix.G}%
\contentsline {section}{\numberline {H}Lemmas necessary for proving Proposition \ref {prop:eff}}{84}{appendix.H}%
\contentsline {section}{\numberline {I}Derivation of the Kalman filter MSE}{85}{appendix.I}%

\newpage

\section*{Notation}
\subsection*{Parameters}
\begin{center}
{
\begin{tabular}{ll}
\hline\\[-10pt]
\hline
Parameters\\
\hline
${\bm\varphi}_n$ & true value\\[4pt]

$\underline{\bm\varphi}_n$& generic value\\[4pt]

$\wh{\bm\varphi}_n^{*}$& QML estimator maximizing the log-likelihood \eqref{eq:LL0true} \\[4pt]

$\wh{\bm\varphi}_n^{\dag}$& QML estimator maximizing the log-likelihood in \citet[Eq. 3]{baili16}\\[4pt]

$\wh{\bm\varphi}_n^{(0)}$& PC estimator used in E-step at iteration $0$\\[4pt]

$\wh{\bm\varphi}_n^{(k)}$& estimator used in E-step at iteration $k > 0$ \\[4pt]

$\wh{\bm\varphi}_n^{(k+1)}$ & estimator computed in M-step at iteration $k\ge 0$ \\[4pt]

$\wh{\bm\varphi}_n\equiv\wh{\bm\varphi}_n^{(k^*+1)}$ & final estimator computed in M-step at iteration $k^*$\\
\hline\\[-10pt]
\hline
\end{tabular}
}
\end{center}

An analogous notation is used for the sub-vectors of parameters ${\bm\phi}_n$ and ${\bm\theta}$ and all their elements. \smallskip

\subsection*{Factors}
\begin{center}
{
\begin{tabular}{ll}
\hline\\[-10pt]
\hline
Factors\\
\hline
$\mbf F_t$ & true value\\[4pt]

$\wt{\mbf F}_t$ & PC estimator\\[4pt]

$\mbf F_{t|s}$ and $\mbf P_{t|s}$ ($s=t-1,t,T$)& estimator and its pseudo-MSE computed via KF-KS using ${\bm\varphi}_n$\\[4pt]

$\mbf F_{0,t|s}$ and $\mbf P_{0,t|s}$ ($s=t-1,t,T$)& estimator and its MSE computed via KF-KS using ${\bm\varphi}_n$ but with $\bm\Gamma_n^\xi$\\[4pt]

$\mbf F_{t|s}^*$ and $\mbf P_{t|s}^*$ ($s=t-1,t,T$)& estimator and its pseudo-MSE computed via KF-KS using $\wh{\bm\varphi}_n^*$\\[4pt]

${\mbf F}_{t|s}^{(k)}$ and ${\mbf P}_{t|s}^{(k)}$ ($s=t-1,t,T$) & estimator and its pseudo-MSE computed via KF-KS using $\wh{\bm\varphi}_n^{(k)}$, $k\ge 0$\\[4pt]

$\wh{\mbf F}_{t|t}\equiv {\mbf F}_{t|t}^{(k^*+1)}$ & final estimator computed via KF at iteration $(k^*+1)$ using $\wh{\bm\varphi}_n$\\[4pt]

$\wh{\mbf P}_{t|t}\equiv{\mbf P}_{t|t}^{(k^*+1)}$ & pseudo-MSE computed via KF at iteration $(k^*+1)$ using $\wh{\bm\varphi}_n$\\[4pt]

$\wh{\mbf F}_{t}\equiv  \wh{\mbf F}_{t|T}\equiv {\mbf F}_{t|T}^{(k^*+1)}$ & final estimator computed via KS at iteration $(k^*+1)$ using $\wh{\bm\varphi}_n$\\[4pt]

$\wh{\mbf P}_{t|T}\equiv{\mbf P}_{t|T}^{(k^*+1)}$ & pseudo-MSE computed via KS at iteration $(k^*+1)$ using $\wh{\bm\varphi}_n$\\[4pt]

$\wh{\mbf F}_t^{\text{\tiny WLS}}$ & WLS estimator computed using $\wh{\bm\varphi}_n$\\[4pt]

\hline\\[-10pt]
\hline
\end{tabular}
}
\end{center}

\begin{appendix}
\setcounter{equation}{0}
\numberwithin{equation}{section}

\section{Further details on estimation}\label{app:EM}

Hereafter, we assume without loss of generality that $\bm \mu_n=\mbf 0_n$ and $p_F=1$ with $\mbf A\equiv\mbf A_1$. When $p_F>1$, it is enough to write the VAR in \eqref{eq:SDFM2R} in companion form and to modify the estimation accordingly, using the augmented state vector $(\mbf F_t^\prime\cdots \mbf F_{t-p_F+1}^\prime)^\prime$.
\subsection{Principal Component estimators}\label{app:prest}
Let $\wh{\bm\Gamma}_n^x$ be the sample covariance matrix of the data and denote as $\wh{\mbf M}_n^x$ the diagonal matrix with entries the $r$-largest eigenvalues of $\wh{\bm\Gamma}_n^x$, and as $\wh{\mbf V}_n^x$ the $n\times r$ matrix of the corresponding normalized eigenvectors. Moreover, 
let $\wh{\bm{\mathcal S} }^{(0)}$ be a $r\times r$ diagonal matrix with entries
$\mathbb I([\wh{\mbf V}_n^x]_{1j}\ge 0)- \mathbb I([\wh{\mbf V}_n^x]_{1j}< 0)$, $j=1,\ldots,r$. Then,
\begin{align}
&\wh{\bm\Lambda}_n^{(0)}= \wh{\mbf V}_n^x\wh{\bm{\mathcal S} }^{(0)}(\wh{\mbf M}_n^x)^{1/2}, \nn\\  
&\wt{\mbf F}_t=(\wh{\bm\Lambda}_n^{(0)\prime}\wh{\bm\Lambda}_n^{(0)})^{-1}\wh{\bm\Lambda}_n^{(0)\prime}\mbf x_{nt} = (\wh{\mbf M}_n^x)^{-1/2}\wh{\bm{\mathcal S} }^{(0)}\wh{\mbf V}_n^{x\prime}\mbf x_{nt}, \nn\\
& \wt{\bm\xi}_{nt}=\mbf x_{nt}-\wh{\bm\Lambda}_n^{(0)\prime} \wt{\mbf F}_{t},\quad t=1,\ldots, T\nn\\
&\wh{\mbf A}^{(0)}=  \l(\sum_{t=2}^T\wt{\mbf F}_{t} \wt{\mbf F}_{t-1}^{\prime} \r)
\l(\sum_{t=2}^T \wt{\mbf F}_{t-1}\wt{\mbf F}_{t-1}^{\prime}\r)^{-1},\nn\\
& \wt{\mbf v}_t= \wt{\mbf F}_t-\wh{\mbf A}^{(0)}\wt{\mbf F}_{t-1}, \quad t=1,\ldots, T,\nn\\
& \wh{\bm\Gamma}^{v(0)} = T^{-1}\sum_{t=1}^T \wt{\mbf v}_t\wt{\mbf v}_t^\prime.\nn
\end{align}
Finally, letting $\wh{\bm\lambda}_i^{(0)\prime}$ be the $i$-th row of $\wh{\bm\Lambda}_n^{(0)}$, and $\wt{\xi}_{it}$ the $i$th component of $\wt{\bm\xi}_{nt}$, then
\beq\nn
\wh{\sigma}_i^{2(0)}=
T^{-1}\sum_{t=1}^T
 \wt{\xi}_{it}^{\,2},\qquad i=1,\ldots,n. 
 \eeq 
The vector of initial estimates of parameters is then:
\beq\nn
\wh{\bm\varphi}_n^{(0)}=\l(
\text{vec}(\wh{\bm\Lambda}_n^{(0)})^\prime\; \wh{\sigma}_1^{2(0)}\cdots \wh{\sigma}_n^{2(0)}\; \text{vec}(\wh{\mbf A}^{(0)})^\prime\; \text{vech}(\wh{\bm\Gamma}^{v(0)})^\prime
\r)^\prime,
\eeq
and it is used to run the first iteration of the EM algorithm.

\subsection{Kalman filter and smoother}\label{app:KFKS}


The following iterations are stated for given initial conditions $\mbf F_{0|0}$ and $\mbf P_{0|0}$ and given the true parameters $\bm\varphi_n$. 

\subsubsection*{Forward iterations - Filtering}\label{app:KF0}

The Kalman filter is based on the forward iterations for $t=1,\ldots, T$:
\begin{align}
&\mbf F_{t|t-1} = \mbf A \mbf F_{t-1|t-1},\label{eq:pred1}\\
&\mbf P_{t|t-1} = \mbf A\mbf P_{t-1|t-1} \mbf A^\prime + \bm\Gamma^v,\label{eq:pred2}\\
&\mbf F_{t|t} =\mbf F_{t|t-1}+\mbf P_{t|t-1}\bm\Lambda_n^\prime(\bm\Lambda_n\mbf P_{t|t-1}\bm\Lambda_n^\prime+\bm\Sigma_n^\xi)^{-1}(\mbf x_{nt}-\bm\Lambda_n\mbf F_{t|t-1}),\label{eq:up1}\\
&\mbf P_{t|t} =\mbf P_{t|t-1}-\mbf P_{t|t-1}\bm\Lambda_n^\prime(\bm\Lambda_n\mbf P_{t|t-1}\bm\Lambda_n^\prime+\bm\Sigma_n^\xi)^{-1}\bm\Lambda_n\mbf P_{t|t-1}.\label{eq:up2}
\end{align}
Moreover, by combining \eqref{eq:pred2} and \eqref{eq:up2}, we obtain the Riccati difference equation:
\beq\label{eq:riccati}
\mbf P_{t+1|t}-\mbf A\mbf P_{t|t-1}\mbf A^\prime+\mbf A\mbf P_{t|t-1}\bm\Lambda_n^\prime(\bm\Lambda_n\mbf P_{t|t-1}\bm\Lambda_n^\prime+ \bm\Sigma_n^\xi)^{-1}\bm\Lambda_n\mbf P_{t|t-1}\mbf A^\prime= \bm\Gamma^v.
\eeq

\subsubsection*{Backward iterations - Smoothing}\label{app:KFback}
The Kalman smoother is then based on the backward iterations for $t=T,\ldots, 1$:
\begin{align}
\mbf F_{t|T} &=\mbf F_{t|t}+\mbf P_{t|t}\mbf A^\prime\mbf P_{t+1|t}^{-1}(\mbf F_{t+1|T}-\mbf F_{t+1|t}),\label{eq:KS1}\\
\mbf P_{t|T}&=\mbf P_{t|t} + \mbf P_{t|t} \mbf A^\prime \mbf P_{t+1|t}^{-1}
(\mbf P_{t+1|T}-\mbf P_{t+1|t})\mbf P_{t+1|t}^{-1} \mbf A \mbf P_{t|t}.\label{eq:KS2}
\end{align}
Finally, ${\mbf C}_{t,t-1|T}$ can be obtained from a state space model with an augmented state vector containing both $\mbf F_t$ and $\mbf F_{t-1}$, by taking the $r\times r$ off-diagonal block of the $2r\times 2r$ matrix defined in \eqref{eq:KS2} but for the augmented model. 

An equivalent way of implementing \eqref{eq:KS1}, which does not require matrix inversion is in  \citet[Chapter 4.4, pp.87-91]{DK01}, which is defined by the backward iterations for $t=T,\ldots, 1$
\begin{align}
&\mbf F_{t|T}=\mbf F_{t|t-1}+\mbf P_{t|t-1}\mbf r_{t-1},\label{eq:KS3}\\
&\mbf r_{t-1}=\bm\Lambda_n^\prime(\bm\Lambda_n\mbf P_{t|t-1}\bm\Lambda_n^\prime+\bm\Sigma_n^\xi)^{-1}(\mbf x_t-\bm\Lambda_n\mbf F_{t|t-1})+\mbf L^\prime_t\mbf r_t,\label{eq:KS4}\\
&\mbf P_{t|T}=\mbf P_{t|t-1}(\mbf I_r-\mbf N_{t-1}\mbf P_{t|t-1}),\label{eq:KS5}\\
&\mbf N_{t-1}=\bm\Lambda_n^\prime(\bm\Lambda_n\mbf P_{t|t-1}\bm\Lambda_n^\prime+\bm\Sigma_n^\xi)^{-1}\bm\Lambda_n+\mbf L_t^\prime\mbf N_t\mbf L_t,\label{eq:KS6}\\
&\mbf L_t= \mbf A-\mbf A \mbf P_{t|t-1} \bm\Lambda_n^\prime (\bm\Lambda_n\mbf P_{t|t-1}\bm\Lambda_n^\prime+\bm\Sigma_n^\xi)^{-1} \bm\Lambda_n,\label{eq:KS7}\\
&\mbf C_{t,t+1|T}=\mbf P_{t|t-1}\mbf L_t^\prime (\mbf I_r-\mbf N_{t}\mbf P_{t+1|t}), \qquad \mbf C_{t,t-1|T}=\mbf C_{t,t+1|T}^\prime,\label{eq:KS8}
\end{align}
where $\mbf r_T=\mbf 0_{r}$, $\mbf N_T=\mbf 0_{r}$ and by construction $\mbf A\mbf P_{t|t}=\mbf L_t\mbf P_{t|t-1}$. 

\subsubsection*{Intialization of Kalman filter and smoother}\label{app:initKF}
The Kalman filter is initialized as follows. At the first iteration of the EM algorithm, i.e., when  $k=0$, we set ${\mbf F}^{(0)}_{0|0}=\mbf 0_r$ and ${\mbf P}^{(0)}_{0|0}=\mbf I_r$, consistently with Assumption \ref{ass:ident}(b). Other initializations as ${\mbf P}^{(0)}_{0|0}=\kappa_0 \mbf I_r$ for some finite real $\kappa_0>0$ are also possible.
At any successive iteration of the EM algorithm, i.e, when $k>0$, we set ${\mbf F}^{(k)}_{0|0}={\mbf F}_{0|T}^{(k-1)}$ and ${\mbf P}^{(0)}_{0|0}=\mbf I_r$.

To run the Kalman smoother we start using the last predictions of the Kalman filter. Thus, for any $k\ge 0$, we set ${\mbf F}^{(k)}_{T+1|T}= \wh{\mbf A}^{(k)} {\mbf F}^{(k)}_{T|T}$ where ${\mbf F}^{(k)}_{T|T}$ is obtained from \eqref{eq:up1}, and we set ${\mbf P}^{(k)}_{T+1|T}=\wh{\mbf A}^{(k)}\mbf P_{T|T}^{(k)} \wh{\mbf A}^{(k)\prime} +\wh{\bm\Gamma}^{v(k)}$, where ${\mbf P}^{(k)}_{T|T}$ is obtained from \eqref{eq:up2}.

\subsection{Stopping rule for the EM algorithm}\label{rem:kstar}
To stop the EM algorithm we adopt the following convergence rule. We fix a maximum finite number of iterations $k_{\max}$, and we stop it at the first iteration $k^*\le k_{\max}$ such that:
\beq
\Delta\ell_{k^*}=\frac{\big\vert \ell(\bm X_{nT};\wh{\bm\varphi}_n^{(k^*+1)})-\ell(\bm X_{nT};\wh{\bm\varphi}_n^{(k^*)})\big\vert}
{\frac 12\big\vert \ell(\bm X_{nT};\wh{\bm\varphi}_n^{(k^*+1)})
+\ell(\bm X_{nT};\wh{\bm\varphi}_n^{(k^*)})
\big\vert}<\varepsilon,
\label{eq:convEM}
\eeq
where $\varepsilon$ is a pre-specified tolerance level. In this case, the log-likelihood is computed using its prediction error formulation obtained from the Kalman filter: 
\begin{align}
\ell(\bm X_{nT};\underline{\bm\varphi}_n)=&\,-\frac 12\sum_{t=1}^T \log\det(\underline{\bm\Lambda}_n\underline{\mbf P}_{t|t-1}\underline{\bm\Lambda}_n^\prime+\underline{\bm\Sigma}_n^\xi)-\frac12\sum_{t=1}^T(\mbf x_{nt}-\underline{\bm\Lambda}_n\underline{\mbf F}_{t|t-1})^\prime (\underline{\bm\Lambda}_n\underline{\mbf P}_{t|t-1}\underline{\bm\Lambda}_n^\prime+\underline{\bm\Sigma}_n^\xi)^{-1}(\mbf x_{nt}-\underline{\bm\Lambda}_n\underline{\mbf F}_{t|t-1}),\nn
\end{align}
where $\underline{\mbf F}_{t|t-1}$ and $\underline{\mbf P}_{t|t-1}$ are computed using \eqref{eq:pred1} and \eqref{eq:pred2}, respectively, when using generic values of the parameters.
Similar convergence criteria can be found in \citet{boothhobert99} and \citet[Chapter 4.9]{MLT07}. 

%
%
\setcounter{equation}{0}
\setcounter{lem}{0}
\section{Proof of main results}

Hereafter, we assume without loss of generality that $\bm \mu_n=\mbf 0_n$ and $p_F=1$ with $\mbf A\equiv\mbf A_1$.

\subsection{Proof of Proposition \ref{prop:plainvanilla}}

Consider the EM algorithm initialized using the PC estimators of the parameters as defined in Section \ref{app:prest}.
At $k^*=0$, from \eqref{eq:param1}, we have
\al{
\wh{\bm\lambda}_i^{(1)} = \l(T^{-1}\sum_{t=1}^T\mbf F_{t|T}^{(0)}\mbf F_{t|T}^{(0)\prime}+\mbf P_{t|T}^{(0)} \r)^{-1}\l(T^{-1}\sum_{t=1}^T \mbf F_{t|T}^{(0)} x_{it}\r).\label{eq:vaccazoccola}
}
Now, 
\al{
\l\Vert T^{-1} \sum_{t=1}^T \mbf F_{t|T}^{(0)}\mbf F_{t|T}^{(0)\prime}-T^{-1} \sum_{t=1}^T  \mbf F_{t}\mbf F_{t}^{\prime}\r\Vert \le &\,
2\l\Vert T^{-1} \sum_{t=1}^T (\mbf F_{t|T}^{(0)}-\mbf F_t)\mbf F_{t}^{\prime}\r\Vert+ \l\Vert T^{-1} \sum_{t=1}^T (\mbf F_{t|T}^{(0)}-\mbf F_t)(\mbf F_{t|T}^{(0)}-\mbf F_t)^{\prime}\r\Vert,\label{eq:vaccaboia1}
}
and
\al{
\l\Vert T^{-1} \sum_{t=1}^T \mbf F_{t|T}^{(0)}x_{it}-T^{-1} \sum_{t=1}^T  \mbf F_{t}x_{it}\r\Vert \le &\,
\l\Vert T^{-1} \sum_{t=1}^T (\mbf F_{t|T}^{(0)}-\mbf F_t)x_{it}\r\Vert.\label{eq:vaccaboia2}
}
Throughout, let $\bm y_t=\mbf F_t$ or $\bm y_t = x_{it}$. Then, we have to consider
\al{
\l\Vert T^{-1} \sum_{t=1}^T (\mbf F_{t|T}^{(0)}-\mbf F_{t})\bm y_t^\prime \r\Vert\le&\, 
\l\Vert T^{-1} \sum_{t=1}^T (\mbf F_{t|T}^{(0)}-\mbf F_{t|t}^{(0)})\bm y_t^\prime \r\Vert+\l\Vert T^{-1} \sum_{t=1}^T (\mbf F_{t|t}^{(0)}-\wh{\mbf F}_{t}^{\text{\tiny WLS}(0)})\bm y_t^\prime \r\Vert+\l\Vert T^{-1} \sum_{t=1}^T (\wh{\mbf F}_{t}^{\text{\tiny WLS}(0)}-\mbf F_t)\bm y_t^\prime \r\Vert\nn\\
=&\, I+II+III, \;\text{ say.}\label{eq:vacca}
}
Let us consider each term in \eqref{eq:vacca}. First,
\al{
I\le&\,\max_{t=1,\ldots, T} \Vert{\mbf P}_{t|t}^{(0)}\Vert \, \Vert\wh{\mbf A}^{(0)}\Vert\, \max_{t=1,\ldots, T}\Vert({\mbf P}_{t+1|t}^{(0)})^{-1}\Vert\l\{\l\Vert T^{-1}\sum_{t=1}^T {\mbf F}_{t+1|T}^{(0)}\bm y_t^\prime\r\Vert
+\Vert\wh{\mbf A}^{(0)}\Vert\,\l\Vert T^{-1}\sum_{t=1}^T {\mbf F}_{t+1|t+1}^{(0)}\bm y_t^\prime\r\Vert\r\}\nn\\
=&\, O_p(n^{-1}),\label{eq:I00}
}
by Lemmas \ref{lem:cazzarolahat0}, \ref{lem:cazzarolahat00}, and \ref{lem:FFO1sum}, and since $\Vert\wh{\mbf A}^{(0)}\Vert\le \Vert{\mbf A}\Vert+\Vert\wh{\mbf A}^{(0)}-\mbf A\Vert =O_p(1)$, by Assumption \ref{ass:common}(d) and Lemma \ref{lem:est0VAR}(i).  Second, from \eqref{eq:KFhathat} and \eqref{eq:girandola} in the proof of Lemma \ref{lem:gennaio22bis}
\al{
II\le &\, O_p(n^{-1}) \l\{\l\Vert n^{-1/2}T^{-1}\sum_{t=1}^T \mbf x_{nt}\bm y_t^\prime\r\Vert+\Vert\wh{\mbf A}^{(0)}\Vert\,\l\Vert T^{-1}\sum_{t=1}^T \mbf F_{t-1|t-1}^{(0)}\bm y_t^\prime  \r\Vert\r\}= O_p(n^{-1}),\label{eq:II00}
}
because of Lemma \ref{lem:FFO1sum} and since 
\al{
&\l\Vert n^{-1/2}T^{-1}\sum_{t=1}^T \bm\Lambda_n\mbf F_t\mbf F_t^\prime\bm\lambda_i\r\Vert\le n^{-1/2}\Vert\bm\Lambda_n\Vert\, \l\Vert T^{-1}\sum_{t=1}^T \mbf F_t\mbf F_t^\prime\r\Vert M_\lambda=O_p(1),\label{eq:varibounds}\\
&\l\Vert n^{-1/2}T^{-1}\sum_{t=1}^T \bm\Lambda_n\mbf F_t\xi_{it}\r\Vert
\le n^{-1/2}\Vert\bm\Lambda_n\Vert\, \l\Vert T^{-1}\sum_{t=1}^T \mbf F_t\xi_{it}\r\Vert=O_p(T^{-1/2}),\nn\\
&\l\Vert n^{-1/2}T^{-1}\sum_{t=1}^T \bm\xi_{nt}\mbf F_t^\prime\bm\lambda_i\r\Vert\le\l\Vert n^{-1/2}T^{-1}\sum_{t=1}^T \bm\xi_{nt}\mbf F_t^\prime\r\Vert M_\lambda=O_p(T^{-1/2}),\nn\\
&\l\Vert n^{-1/2}T^{-1}\sum_{t=1}^T \bm\xi_{nt}\xi_{it}\r\Vert
\le \l\Vert n^{-1/2}T^{-1}\sum_{t=1}^T \bm\xi_{nt}\xi_{it}-n^{-1/2}\E[\bm\xi_{nt}\xi_{it}]\r\Vert+  n^{-1/2}\Vert\E[\bm\xi_{nt}\xi_{it}]\Vert=O_p(T^{-1/2})+O(1),\nn
}
by Assumption \ref{ass:common}(a), Lemmas \ref{lem:lambdasqrtn}, \ref{lem:consistCOV}(i), \ref{lem:consistCOV}(ii), \ref{lem:consistCOV}(iii), and \ref{lem:consistCOV}(iv), and because
$n^{-1}\Vert\E[\bm\xi_{nt}\xi_{it}]\Vert^2 \le n^{-1}\sum_{j=1}^n \vert\E[\xi_{jt}^2\xi_{it}^2] \vert\le \mathrm K_\xi$, by Assumption \ref{ass:idio}(d). Note that the first and third relations in \eqref{eq:varibounds} cover also the case $\bm y_t=\mbf F_t^\prime$.

Finally, let us consider the last term in \eqref{eq:vacca}. From \eqref{eq:ultimopezzo} in the proof of Lemma \ref{lem:aprile24}
\al{
III \le &\, \l\Vert
T^{-1}\sum_{t=1}^T (\wh{\bm\Lambda}_n^{(0)\prime}(\wh{\bm\Sigma}_n^{\xi(0)})^{-1}\wh{\bm\Lambda}_n^{(0)})^{-1}\wh{\bm\Lambda}_n^{(0)\prime}(\wh{\bm\Sigma}_n^{\xi(0)})^{-1}(\bm\Lambda_n-\wh{\bm\Lambda}_n^{(0)})\mbf F_t\bm y_t^\prime
\r\Vert\nn\\
&+ \l\Vert T^{-1}\sum_{t=1}^T (\wh{\bm\Lambda}_n^{(0)\prime}(\wh{\bm\Sigma}_n^{\xi(0)})^{-1}\wh{\bm\Lambda}_n^{(0)})^{-1}\wh{\bm\Lambda}_n^{(0)\prime}(\wh{\bm\Sigma}_n^{\xi(0)})^{-1}\bm\xi_{nt}\bm y_t^\prime
\r\Vert\nn\\
\le&\, 
\Vert(\bm\Lambda_n^\prime(\bm\Sigma_n^\xi)^{-1}\bm\Lambda_n)^{-1}\Vert \, \Vert \bm\Lambda_n^\prime(\bm\Sigma_n^\xi)^{-1} (\bm\Lambda_n-\wh{\bm\Lambda}_n^{(0)}) \Vert \l\Vert T^{-1}\sum_{t=1}^T \mbf F_t\bm y_t^\prime\r\Vert\nn\\ 
&+\Vert
n(\wh{\bm\Lambda}_n^{(0)\prime}(\wh{\bm\Sigma}_n^{\xi(0)})^{-1}\wh{\bm\Lambda}_n^{(0)})^{-1}n^{-1/2}\wh{\bm\Lambda}_n^{(0)\prime}(\wh{\bm\Sigma}_n^{\xi(0)})^{-1}-n({\bm\Lambda}_n^{\prime}({\bm\Sigma}_n^{\xi})^{-1}{\bm\Lambda}_n)^{-1}n^{-1/2}{\bm\Lambda}_n^{\prime}({\bm\Sigma}_n^{\xi})^{-1}\Vert \nn\\
&\cdot n^{-1/2}\Vert \bm\Lambda_n-\wh{\bm\Lambda}_n^{(0)}\Vert \, \l\Vert T^{-1}\sum_{t=1}^T \mbf F_t\bm y_t^\prime\r\Vert\nn\\ 
&+n\Vert(\wh{\bm\Lambda}_n^{(0)\prime}(\wh{\bm\Sigma}_n^{\xi(0)})^{-1}\wh{\bm\Lambda}_n^{(0)})^{-1}\Vert \,n^{-1}\l\Vert T^{-1}\sum_{t=1}^T \bm\Lambda_n^\prime(\bm\Sigma_n^\xi)^{-1}\bm\xi_{nt}\bm y_t^\prime \r\Vert\nn\\
&+n \Vert
(\wh{\bm\Lambda}_n^{(0)\prime}(\wh{\bm\Sigma}_n^{\xi(0)})^{-1} \wh{\bm\Lambda}_n^{(0)})^{-1}\Vert\, n^{-1}\l\Vert T^{-1}\sum_{t=1}^T \{\wh{\bm\Lambda}_n^{(0)\prime}(\wh{\bm\Sigma}_n^{\xi(0)})^{-1}-{\bm\Lambda}_n^{\prime}({\bm\Sigma}_n^{\xi})^{-1}\} \bm\xi_{nt}\bm y_t^\prime\r\Vert\nn\\ 
=&\, III_a+III_b+III_c+III_d, \; \text{ say.}\label{eq:III0}
}
Then, 
\beq
III_a= O_p(n^{-1/2}T^{-1/2}),\label{eq:IIIa}
\eeq
by \eqref{eq:varibounds} and \eqref{eq:ultimopezzoAA}-\eqref{eq:ultimopezzoAA0} in the proof of Lemma \ref{lem:aprile24}.
Moreover,  
\beq
III_b=O_p(\max(n^{-2},T^{-1})),\label{eq:IIIb}
\eeq 
by \eqref{eq:varibounds} and Lemmas \ref{lem:est0LOAD}(ii) and \ref{lem:est0_LAST}(v).
Regarding $III_c$, if $\bm y_t=\mbf F_t$, we have
\al{
III_c &= n\Vert(\wh{\bm\Lambda}_n^{(0)\prime}(\wh{\bm\Sigma}_n^{\xi(0)})^{-1}\wh{\bm\Lambda}_n^{(0)})^{-1}\Vert \, n^{-1}\l\Vert T^{-1}\sum_{t=1}^T \bm\Lambda_n^\prime(\bm\Sigma_n^\xi)^{-1}\bm\xi_{nt}\mbf F_t^\prime \r\Vert= O_p(n^{-1/2}T^{-1/2}),\label{eq:IIIc1}
}
by Lemmas \ref{lem:est0_LAST}(iii) and \ref{lem:COVFF0}(iv).  If $\bm y_t=x_{it}$, we have
\al{
III_c \le&\, n\Vert(\wh{\bm\Lambda}_n^{(0)\prime}(\wh{\bm\Sigma}_n^{\xi(0)})^{-1}\wh{\bm\Lambda}_n^{(0)})^{-1}\Vert \, n^{-1}\l\Vert T^{-1}\sum_{t=1}^T \bm\Lambda_n^\prime(\bm\Sigma_n^\xi)^{-1}\bm\xi_{nt}\mbf F_t^\prime\r\Vert \, \Vert\bm\lambda_i \Vert\nn\\
&+n\Vert(\wh{\bm\Lambda}_n^{(0)\prime}(\wh{\bm\Sigma}_n^{\xi(0)})^{-1}\wh{\bm\Lambda}_n^{(0)})^{-1}\Vert\, \l\Vert T^{-1}\sum_{t=1}^T \bm\Lambda_n^\prime(\bm\Sigma_n^\xi)^{-1}\bm\xi_{nt}\xi_{it} \r\Vert\nn\\
=&\, O_p(n^{-1/2}T^{-1/2})),\label{eq:IIIc2}
}
by Assumption \ref{ass:common}(a), and Lemmas \ref{lem:est0_LAST}(iii), \ref{lem:COVFF0}(iv), and \ref{lem:COVFF0}(v). 
Last, 
\beq
III_d=O_p(\max(n^{-1},T^{-1/2})),
\label{eq:IIId}
\eeq 
by  \eqref{eq:varibounds} and Lemma \ref{lem:est0_LAST}(v). From \eqref{eq:III0}, \eqref{eq:IIIa}, \eqref{eq:IIIb}, \eqref{eq:IIIc1}, \eqref{eq:IIIc2}, and \eqref{eq:IIId}
\al{
III =O_p(\max(n^{-1},T^{-1/2})).\label{eq:III00}
}

Combining \eqref{eq:I00}, \eqref{eq:II00}, and \eqref{eq:III00} we have
\al{
\l\Vert T^{-1} \sum_{t=1}^T (\mbf F_{t|T}^{(0)}-\mbf F_t)\mbf F_{t}^{\prime}\r\Vert=O_p(\max(n^{-1},T^{-1/2})),\label{eq:bella1}\\
\l\Vert T^{-1} \sum_{t=1}^T (\mbf F_{t|T}^{(0)}-\mbf F_t)x_{it}\r\Vert=O_p(\max(n^{-1},T^{-1/2})),\label{eq:bella2}
}
which once substituted into \eqref{eq:vaccaboia1} and \eqref{eq:vaccaboia2}, jointly with Lemma \ref{lem:PPOn} give
\al{
\l\Vert T^{-1} \sum_{t=1}^T \{\mbf F_{t|T}^{(0)}\mbf F_{t|T}^{(0)\prime}+\mbf P_{t|T}^{(0)}\}-T^{-1} \sum_{t=1}^T  \mbf F_{t}\mbf F_{t}^{\prime}\r\Vert\le&\,  
\l\Vert T^{-1} \sum_{t=1}^T \mbf F_{t|T}^{(0)}\mbf F_{t|T}^{(0)\prime}-T^{-1} \sum_{t=1}^T  \mbf F_{t}\mbf F_{t}^{\prime}\r\Vert +\max_{t=1,\ldots, T} \Vert \mbf P_{t|T}^{(0)}\Vert \nn\\
=&\, O_p(\max(n^{-1},T^{-1/2})) + O_p(n^{-1}),\nn
}
and
\al{
\l\Vert T^{-1} \sum_{t=1}^T \mbf F_{t|T}^{(0)}x_{it}-T^{-1} \sum_{t=1}^T  \mbf F_{t}x_{it}\r\Vert = O_p(\max(n^{-1},T^{-1/2})).\nn
}
Therefore, from \eqref{eq:vaccazoccola}
\al{
\Vert\wh{\bm\lambda}_i^{(1)} -\bm\lambda_i^{\text{\tiny OLS}}\Vert = O_p(\max(n^{-1},T^{-1/2})),\label{eq:daje1}
}
with $\bm\lambda_i^{\text{\tiny OLS}} = (T^{-1}\sum_{t=1}^T\mbf F_t\mbf F_t^\prime)^{-1}(T^{-1}\sum_{t=1}^T\mbf F_tx_{it})$. And by Lemma \ref{lem:olsT12}(i) we also have
\al{
\Vert\bm\lambda_i^{\text{\tiny OLS}}-\bm\lambda_i\Vert 
= O_p(T^{-1/2}),\label{eq:daje2}
} 
indeed, recalling that $\bm\Gamma^F=\mbf I_r$ by Assumption \ref{ass:ident}(b), by Lemma \ref{lem:consistCOV}(i) and Weyl's inequality \citep[Theorem 1]{MK04} we have $\vert \nu^{(r)}(T^{-1}\sum_{t=1}^T \mbf F_t\mbf F_t^\prime)^{-1}\vert = O_p(T^{-1/2})$ which implies
$\Vert (T^{-1}\sum_{t=1}^T \mbf F_t\mbf F_t^\prime)^{-1}
\Vert = O_p(1).$
From \eqref{eq:daje1} and \eqref{eq:daje2}
\al{
\Vert\wh{\bm\lambda}_i^{(1)} -\bm\lambda_i\Vert \le \Vert\wh{\bm\lambda}_i^{(1)} -\bm\lambda_i^{\text{\tiny OLS}}\Vert + \Vert\bm\lambda_i^{\text{\tiny OLS}}-\bm\lambda_i\Vert= O_p(\max(n^{-1},T^{-1/2})).\label{eq:dajeLOAD}
}
Moreover, by letting $\bm y_t=n^{-1/2} \mbf x_{nt}$ the above proof leads to
\al{
\l\Vert n^{-1/2}T^{-1} \sum_{t=1}^T (\mbf F_{t|T}^{(0)}-\mbf F_t)\mbf x_{nt}^{\prime}\r\Vert=O_p(\max(n^{-1},T^{-1/2})),\label{eq:bella3}
}
and, therefore, using also Lemma \ref{lem:olsT12}(ii), we have
\al{
n^{-1/2} \Vert \wh{\bm\Lambda}_n^{(1)} -\bm\Lambda_n\Vert = O_p(\max(n^{-1},T^{-1/2})).\label{eq:dajeLOADvec}
}
This proves parts (a.1) and (a.2)  when $k^*=0$.

Following the same reasoning leading to \eqref{eq:dajeLOAD} and by Lemma \ref{lem:olsT12}(iv) we can easily prove that
\al{
\Vert\wh{\mbf A}^{(1)} -\mbf A\Vert \le \Vert\wh{\mbf A}^{(1)} -\mbf A^{\text{\tiny OLS}}\Vert + \Vert\mbf A^{\text{\tiny OLS}}-\mbf A\Vert= O_p(\max(n^{-1},T^{-1/2})),\label{eq:dajeAA}
}
where $\wh{\mbf A}^{(1)}$ is defined in \eqref{eq:param4} and $\mbf A^{\text{\tiny OLS}} = (T^{-1}\sum_{t=1}^T \mbf F_t\mbf F_{t-1}^\prime )(T^{-1}\sum_{t=1}^T \mbf F_{t-1}\mbf F_{t-1}^\prime )^{-1}$ (recall that $\mbf F_0=\mbf 0_r$ by Assumption \ref{ass:common}(i)). 
To prove \eqref{eq:dajeAA} we also use the fact that $\max_{t=1,\ldots, T}\Vert \mbf C_{t,t-1|T}^{(0)}\Vert=O(n^{-1})$ since it can be obtained by the upper right block of $\mbf P_{t|T}^{(0)}$ when this is computed from the the Kalman smoother having the augmented state vector $(\mbf F_t^\prime\,\mbf F_{t-1}^\prime)^\prime$. This proves part (a.4)  when $k^*=0$.

Likewise, using \eqref{eq:dajeAA} and the same reasoning leading to \eqref{eq:dajeLOAD}, by Lemma \ref{lem:olsT12}(v), we can easily prove also that 
\beq
\Vert\wh{\bm\Gamma}^{v(1)} -\bm\Gamma^v\Vert\le 
\Vert\wh{\bm\Gamma}^{v(1)} -\bm\Gamma^{v\text{\tiny OLS}}\Vert+
\Vert\bm\Gamma^{v\text{\tiny OLS}} -\bm\Gamma^v\Vert=
O_p(\max(n^{-1},T^{-1/2})).\label{eq:dajeGv}
\eeq
where $\wh{\bm\Gamma}^{v(1)}$ is defined in \eqref{eq:paramGv} and $\bm\Gamma^{v\text{\tiny OLS}}=T^{-1}\sum_{t=1}^T (\mbf F_t-\mbf A^{\text{\tiny OLS}}\mbf F_{t-1}) (\mbf F_t-\mbf A^{\text{\tiny OLS}}\mbf F_{t-1})^\prime$. To prove \eqref{eq:dajeGv} we need to use also the intermediate quantity 
$T^{-1}\sum_{t=1}^T (\mbf F_t-\wh{\mbf A}^{(1)}\mbf F_{t-1}) (\mbf F_t-\wh{\mbf A}^{(1)}\mbf F_{t-1})^\prime$, which is $\min(n,\sqrt T)$-consistent because of \eqref{eq:dajeAA}. This proves part (a.5)  when $k^*=0$.

Finally, using again the same reasoning leading to \eqref{eq:dajeLOAD}, by Lemma \ref{lem:olsT12}(iii), we can prove that
\beq
\vert \wh{\sigma}_i^{2(1)}-\sigma_i^2\vert \le\vert \wh{\sigma}_i^{2(1)}-\sigma_i^{2\text{\tiny OLS}}\vert+\vert \sigma_i^{2\text{\tiny OLS}}-\sigma_i^2\vert=O_p(\max(n^{-1},T^{-1/2})),\label{eq:dajesigma}
\eeq
where $\wh{\sigma}_i^{2(1)}$ is defined in \eqref{eq:param3} and $\sigma_i^{2\text{\tiny OLS}}=T^{-1}\sum_{t=1}^T (x_{it}-{\bm\lambda}_i^{\text{\tiny OLS}\prime}\mbf F_t)^2$. To prove \eqref{eq:dajesigma} we need to use also the intermediate quantity 
$T^{-1}\sum_{t=1}^T (x_{it}-\wh{\bm\lambda}_i^{(1)\prime}\mbf F_t)^2$, which is $\min(n,\sqrt T)$-consistent because of \eqref{eq:dajeLOAD}. This proves part (a.3)  when $k^*=0$.

Now, from \eqref{eq:dajeLOAD} and \eqref{eq:dajesigma} using the same reasoning of the proof of Lemma \ref{lem:est0}(ii), which in turn requires \eqref{eq:bella1} and \eqref{eq:bella2}, again it follows that
\beq
n^{-1}\l\vert\sum_{i=1}^n(\wh{\sigma}_{i}^{(1)2}-\sigma_i^2)\r\vert =O_p(\max(n^{-1},T^{-1/2})). \label{eq:dajesigmamean}
\eeq
And, using the same reasoning  as in the proof of Lemma \ref{lem:est0_LAST} but using now \eqref{eq:dajeLOADvec},  \eqref{eq:dajesigma}, and \eqref{eq:dajesigmamean},
\al{
&n^{-1}\Vert\wh{\bm\Lambda}_n^{(1)\prime}(\wh{\bm \Sigma}_n^{\xi(1)})^{-1}\wh{\bm\Lambda}_n^{(1)}-\bm\Lambda_n^\prime(\bm\Sigma_n^\xi)^{-1}\bm\Lambda_n\Vert = O_p(\max(n^{-1},T^{-1/2})),\nn\\
&n^{-1/2}\Vert\wh{\bm\Lambda}_n^{(1)\prime}(\wh{\bm \Sigma}_n^{\xi(1)})^{-1}-\bm\Lambda_n^\prime(\bm\Sigma_n^\xi)^{-1}\Vert = O_p(\max(n^{-1},T^{-1/2})),\label{eq:varieLSL}\\
&n\Vert(\wh{\bm\Lambda}_n^{(1)\prime}(\wh{\bm \Sigma}_n^{\xi(1)})^{-1}\wh{\bm\Lambda}_n^{(1)})^{-1}\Vert = O_p(1)\nn\\
&n\Vert(\wh{\bm\Lambda}_n^{(1)\prime}(\wh{\bm \Sigma}_n^{\xi(1)})^{-1}\wh{\bm\Lambda}_n^{(1)})^{-1}-(\bm\Lambda_n^\prime(\bm\Sigma_n^\xi)^{-1}\bm\Lambda_n)^{-1}\Vert = O_p(\max(n^{-1},T^{-1/2}))\nn.
}
From  \eqref{eq:dajeLOADvec}, \eqref{eq:dajeAA}, \eqref{eq:dajeGv}, and the relations in \eqref{eq:varieLSL} and following the same steps as in the proofs of Lemmas \ref{lem:cazzarolahat0}, \ref{lem:cazzarolahat00}, and \ref{lem:PPOn}, we get
\al{
&\max_{t=1,\ldots, T}\Vert \mbf P^{(1)}_{t|t-1}\Vert=O_p(1),\quad \max_{t=1,\ldots, T}\Vert (\mbf P^{(1)}_{t|t-1})^{-1}\Vert=O_p(1),\label{eq:varieP}\\
&\max_{t=1,\ldots, T}\Vert \mbf P^{(1)}_{t|t}\Vert=O_p(n^{-1}),\quad \max_{t=1,\ldots, T}\Vert \mbf P^{(1)}_{t|T}\Vert=O_p(n^{-1}).\nn
}
and also $\Vert\mbf F_{t|t}^{(1)}\Vert=O_p(1)$ and $\Vert\mbf F_{t|T}^{(1)}\Vert=O_p(1)$ by the same arguments in Lemma \ref{lem:FFO1}. It follows that we can apply the same steps as in the proofs of Lemmas \ref{lem:FFOn}, \ref{lem:gennaio22bis}, and \ref{lem:aprile24} to get
$$
\Vert {\mbf F}_{t|T}^{(1)}-{\mbf F}_{t}\Vert=O_p(\max(n^{-1/2},T^{-1/2})).
$$
This proves part (b) when $k^*=0$.
 
Then we can show that \eqref{eq:bella1} and \eqref{eq:bella2} still hold when using ${\mbf F}_{t|T}^{(1)}$ in place of ${\mbf F}_{t|T}^{(0)}$ and using also the last of \eqref{eq:varieP} we prove $\min(n,\sqrt T)$-consistency of $\wh{\bm\lambda}_i^{(2)}$. Similarly we can prove $\min(n,\sqrt T)$-consistency of $n^{-1/2}\wh{\bm\Lambda}_n^{(2)}$, $\wh{\mbf A}^{(2)}$, $\wh{\bm\Gamma}^{v(2)}$, and $\sigma_i^{2(2)}$. It is then clear that we can repeat the same reasoning leading to \eqref{eq:dajesigmamean}, \eqref{eq:varieLSL}, and \eqref{eq:varieP} but when $k^*=1$. So these arguments hold for all $k^*\ge 0$. This completes the proof. $\Box$

\subsection{Proof of Proposition \ref{prop:load}}
For part (a.1), for any $k^*\ge 0$, we have (recall that $\wh{\bm\lambda}_{i}\equiv \wh{\bm\lambda}_{i}^{(k^*+1)}$)
\begin{align}
(\wh{\bm\lambda}_{i}-{\bm\lambda}_{i}) =&\,  (\wh{\bm\lambda}_{i}-\wh{\bm\lambda}_{i}^{**})+(\wh{\bm\lambda}_{i}^{**}-\wh{\bm\lambda}_{i}^{*})+(\wh{\bm\lambda}_{i}^{*}-{\bm\lambda}_{i}^{\text{\tiny OLS}})+({\bm\lambda}_{i}^{\text{\tiny OLS}}-\bm\lambda_i) \nn\\
=&\, L.1+L.2+L.3+L.4, \; \text{ say.} 
\label{eq:loadingsconsistenti}
\end{align}
From Lemma \ref{lem:convEM1}
\al{
\Vert L.1\Vert = O_p(\max(
n^{-2}\log^{4/\delta_v}T,
n^{-1}T^{-1}\log^{1/\delta_v}T\sqrt{\log n},
T^{-3/2}\sqrt{\log n}
)).\label{eq:ELLE1}
}
From Lemma \ref{lem:localglobal}(i)
\al{
\Vert L.2\Vert =O_p(\max(n^{-1}\log^{2/\delta_v}T,n^{-1/2}T^{-1/2}\sqrt{\log n},T^{-1})).\label{eq:ELLE2}
}
From Lemma \ref{lem:starols}(i)
\al{
\Vert L.3\Vert=O_p(\max(n^{-1},n^{-1/2}T^{-1/2},T^{-1}) ).\label{eq:ELLE3}
}
From Lemma \ref{lem:olsT12}(i)
\al{
\Vert L.4\Vert = O_p(T^{-1/2}).\label{eq:ELLE4}
}
By using \eqref{eq:ELLE1}, \eqref{eq:ELLE2}, \eqref{eq:ELLE3}, and \eqref{eq:ELLE4} into \eqref{eq:loadingsconsistenti}, we prove part (a.1). 

For part (a.2), for any $k^*\ge 0$, we have (recall that $\wh{\bm\Lambda}_{n}\equiv \wh{\bm\Lambda}_{n}^{(k^*+1)}$)
\begin{align}
(\wh{\bm\Lambda}_{n}-{\bm\Lambda}_{n}) =&\,  (\wh{\bm\Lambda}_{n}-\wh{\bm\Lambda}_{n}^{**})+(\wh{\bm\Lambda}_{n}^{**}-\wh{\bm\Lambda}_{n}^{*})+(\wh{\bm\Lambda}_{n}^{*}-{\bm\Lambda}_{n}^{\text{\tiny OLS}})+({\bm\Lambda}_{n}^{\text{\tiny OLS}}-\bm\Lambda_n),\nn 
\end{align}
and the proof follows from Lemmas  \ref{lem:olsT12}(ii), \ref{lem:starols}(ii), and \ref{lem:convEM1unif}(i), and since
\al{
n^{-1}\Vert \wh{\bm\Lambda}_{n}^{**}-\wh{\bm\Lambda}_{n}^{*}\Vert^2&= n^{-1}\sum_{i=1}^n \Vert \wh{\bm\lambda}_{i}^{**}-\wh{\bm\lambda}_{i}^{*}\Vert^2\le \max_{i=1,\ldots, n} \Vert \wh{\bm\lambda}_{i}^{**}-\wh{\bm\lambda}_{i}^{*}\Vert^2\nn\\
&= O_p(\max(n^{-2}\log^{4/\delta_v}T,n^{-1}T^{-1}{\log n},T^{-2})),\label{eq:medipd}
}
by Lemma \ref{lem:localglobal}(i).

For part (b), from part (a.1) and \eqref{eq:loadingsconsistenti},
 if $n^{-1}\sqrt T\log^{2/\delta_v} T\to 0$, as $n,T\to\infty$, we have 
%
%
%
%
\begin{align}
\sqrt T(\wh{\bm\lambda}_i -\bm\lambda_i)&=\sqrt T({\bm\lambda}_i^{\text{\tiny OLS}} -\bm\lambda_i)+o_p(1)\nn\\
&=\l(T^{-1}\sum_{t=1}^T\mbf F_t\mbf F_t^\prime\r)^{-1}\l(T^{-1/2}\sum_{t=1}^T \mbf F_t \xi_{it} \r)+o_p(1).
\label{eq:baili}
\end{align}
Now, since $\{\mbf F_t \xi_{it}\}$ is strongly mixing with exponentially decaying coefficients by \citet[Theorem 5.1.a]{bradley05}
(see also \eqref{eq:mixFxi} in the proof of Lemma \ref{lem:tail}), and given that by Assumption \ref{ass:tails} the following Cram\'er condition holds 
$$
\sup_{m\ge 1} r^{-1/\delta} (\E[\vert \xi_{it}F_{jt} \vert^m])^{1/m}\le K,
$$ 
for some finite positive reals $\delta\in\l(0,\frac{\delta_v\delta_\xi}{\delta_v+\delta_\xi}\r)$ and $K$  independent of $t$, $i$, and $j$ \citep[Section 2]{KC18}, then the Central Limit Theorem by \citet[Theorem 1.7]{ibra62} applies, i.e.,
\beq
T^{-1/2}\sum_{t=1}^T\mbf F_t\xi_{it}\stackrel{d}{\ra} \mathcal N\l(\mbf 0_r,\lim_{T\to\infty}T^{-1}\sum_{s,t=1}^T\E \l[\mbf F_t\mbf F_s^\prime \xi_{it}\xi_{is}\r]
\r).\label{eq:ibraccio}
\eeq
Therefore, by Lemmas \ref{lem:consistCOV}(i) and \ref{lem:frida}, and Assumption \ref{ass:common}(b)
\al{
\l\Vert \l(T^{-1}\sum_{t=1}^T\mbf F_t\mbf F_t^\prime\r)^{-1}- (\bm\Gamma^F)^{-1}\r\Vert\le &\, 
\l\Vert \l(T^{-1}\sum_{t=1}^T\mbf F_t\mbf F_t^\prime\r)^{-1}\r\Vert \,
\l\Vert T^{-1}\sum_{t=1}^T\mbf F_t\mbf F_t^\prime-\bm\Gamma^F\r\Vert
\,
\Vert (\bm\Gamma^F)^{-1}\Vert\nn\\
=&\,O_p(T^{-1/2}).\label{eq:fridahot2}
}
Thus from \eqref{eq:baili}, \eqref{eq:ibraccio}, and \eqref{eq:fridahot2}, by Slutsky's Theorem, we have
\beq
\sqrt T (\wh{\bm\lambda}_i-\bm\lambda_i) \stackrel{d}{\ra} \mathcal N(\mbf 0_r, \bm{\mathcal V}_{i}), \nn
\eeq
where,
\begin{align}
\bm{\mathcal V}_{i} &
=(\bm\Gamma^F)^{-1}\l\{ \lim_{T\to\infty}T^{-1}\sum_{s,t=1}^T\E \l[\mbf F_t \xi_{it}\xi_{is}\mbf F_s^\prime\r]\r\}(\bm\Gamma^F)^{-1}
=\l\{ T^{-1}\sum_{s,t=1}^T\E[\xi_{it}\xi_{is}]\,\E \!\l[\mbf F_t \mbf F_s^\prime\r]\r\},\nn 
\end{align}
since $\bm\Gamma^F=\mbf I_r$ because of Assumption \ref{ass:ident}(b) and $\{\mbf F_t\}$ and $\{\xi_{it}\}$ are independent processes because of Lemma \ref{lem:fidio}. This proves part (b). Part (c) is straightforward. This completes the proof. $\Box$

\subsection{Proof of Proposition \ref{prop:factors}}

Recall the definitions $\wh{\mbf F}_t\equiv \wh{\mbf F}_{t|T}\equiv \mbf F_{t|T}^{(k^*+1)}$, for any $k^*\ge 0$. From Lemmas \ref{lem:FFOnhat}(ii) and   \ref{lem:FFOnhat}(iii),
\al{
\Vert\wh{\mbf F}_{t|T}-\mbf F_t\Vert&\le \Vert\wh{\mbf F}_{t|T}-\wh{\mbf F}_{t|t}\Vert+\Vert\wh{\mbf F}_{t|t}-\wh{\mbf F}_t^{\text{\tiny WLS}}\Vert + \Vert \wh{\mbf F}_t^{\text{\tiny WLS}}-\mbf F_t\Vert\nn\\
&= \Vert\wh{\mbf F}_t^{\text{\tiny WLS}}-\mbf F_t\Vert + O_p(n^{-1}).\label{eq:ultimopezzoamontehat}
}
where $\wh{\mbf F}_{t}^{\text{\tiny \upshape {WLS}}}=(\wh{\bm\Lambda}_n^{\prime}(\wh{\bm\Sigma}_n^{\xi})^{-1}\wh{\bm\Lambda}_n)^{-1}\wh{\bm\Lambda}_n^{\prime}(\wh{\bm\Sigma}_n^{\xi})^{-1}\mbf x_{nt}.$

Now,
\al{
\Vert \wh{\mbf F}_t^{\text{\tiny WLS}}-\mbf F_t\Vert\le&\, \Vert(\wh{\bm\Lambda}_n^{\prime}(\wh{\bm\Sigma}_n^{\xi})^{-1}\wh{\bm\Lambda}_n)^{-1}\wh{\bm\Lambda}_n^{\prime}(\wh{\bm\Sigma}_n^{\xi})^{-1}(\bm\Lambda_n-\wh{\bm\Lambda}_n)\Vert\, \Vert\mbf F_t\Vert\nn\\
&+ \Vert(\wh{\bm\Lambda}_n^{\prime}(\wh{\bm\Sigma}_n^{\xi})^{-1}\wh{\bm\Lambda}_n)^{-1}\wh{\bm\Lambda}_n^{\prime}(\wh{\bm\Sigma}_n^{\xi})^{-1}\bm\xi_{nt}\Vert
\nn\\
\le&\,
\Vert
(\bm\Lambda_n^{\prime}
({\bm\Sigma}_n^{\xi})^{-1}
{\bm\Lambda}_n)^{-1}
{\bm\Lambda}_n^\prime
({\bm\Sigma}_n^{\xi})^{-1}
(\bm\Lambda_n-\wh{\bm\Lambda}_n)\Vert\, 
\Vert\mbf F_t\Vert\nn\\
&+\Vert
(\wh{\bm\Lambda}_n^{\prime}(\wh{\bm\Sigma}_n^{\xi})^{-1}\wh{\bm\Lambda}_n)^{-1}\wh{\bm\Lambda}_n^{\prime}(\wh{\bm\Sigma}_n^{\xi})^{-1}-({\bm\Lambda}_n^{\prime}({\bm\Sigma}_n^{\xi})^{-1}{\bm\Lambda}_n)^{-1}{\bm\Lambda}_n^{\prime}({\bm\Sigma}_n^{\xi})^{-1}\Vert\nn\\
&\cdot\Vert \bm\Lambda_n-\wh{\bm\Lambda}_n\Vert\, \Vert\mbf F_t\Vert+\Vert ({\bm\Lambda}_n^{\prime}({\bm\Sigma}_n^{\xi})^{-1}{\bm\Lambda}_n)^{-1} {\bm\Lambda}_n^{\prime}({\bm\Sigma}_n^{\xi})^{-1}\bm\xi_{nt}\Vert\nn\\
&+\Vert
(\wh{\bm\Lambda}_n^{\prime}(\wh{\bm\Sigma}_n^{\xi})^{-1}\wh{\bm\Lambda}_n)^{-1}\wh{\bm\Lambda}_n^{\prime}(\wh{\bm\Sigma}_n^{\xi})^{-1}\bm\xi_{nt}-({\bm\Lambda}_n^{\prime}({\bm\Sigma}_n^{\xi})^{-1}{\bm\Lambda}_n)^{-1}{\bm\Lambda}_n^{\prime}({\bm\Sigma}_n^{\xi})^{-1}\bm\xi_{nt}\Vert \nn\\
=&\, \rm A+\rm B+\rm C+\rm D, \;\text{ say.} \label{eq:ultimopezzohat}
}
Let us consider each term in \eqref{eq:ultimopezzohat}. First, consider term $\rm A$ and notice that
\al{
n^{-1/2} \Vert \wh{\bm\Lambda}_n-\bm\Lambda_n^{\text{\tiny OLS}}\Vert\le &n^{-1/2} \Vert \wh{\bm\Lambda}_n-\wh{\bm\Lambda}_n^{**}\Vert
+n^{-1/2} \Vert \wh{\bm\Lambda}_n^{**}-\wh{\bm\Lambda}_n^{*}\Vert
+n^{-1/2} \Vert \wh{\bm\Lambda}_n^*-\bm\Lambda_n^{\text{\tiny OLS}}\Vert\nn\\
=&\,
O_p(\max(
n^{-2}\log^{4/\delta_v}T,
n^{-1}T^{-1}\log^{1/\delta_v}T\sqrt{\log n},
T^{-3/2}{\log n}
))\nn\\
&+O_p(\max(n^{-1}\log^{2/\delta_v}T,n^{-1/2}T^{-1/2}\sqrt{\log n},T^{-1}))\nn\\
&+O_p(\max(n^{-1}\log^{2/\delta_v}T,n^{-1/2}T^{-1/2},T^{-1}))\nn\\
=&\,O_p(\max(n^{-1}\log^{2/\delta_v}T,n^{-1/2}T^{-1/2}\sqrt{\log n},T^{-1})),\label{eq:reggioemilia}
}
by Lemmas  
\ref{lem:starols}(ii), 
\ref{lem:localglobal}(i), and
\ref{lem:convEM1unif}(i) (see also \eqref{eq:medipd} in the proof of Proposition \ref{prop:load}).

Therefore, from \eqref{eq:reggioemilia}
\al{
\rm A\le&\, \Vert
(\bm\Lambda_n^{\prime}
({\bm\Sigma}_n^{\xi})^{-1}
{\bm\Lambda}_n)^{-1}
{\bm\Lambda}_n^\prime
({\bm\Sigma}_n^{\xi})^{-1}
(\bm\Lambda_n-\bm\Lambda_n^{\text{\tiny OLS}})\Vert\, 
\Vert\mbf F_t\Vert\nn\\
&+n \Vert
(\bm\Lambda_n^{\prime}
({\bm\Sigma}_n^{\xi})^{-1}
{\bm\Lambda}_n)^{-1}\Vert\, n^{-1/2}\Vert{\bm\Lambda}_n\Vert\, 
\Vert({\bm\Sigma}_n^{\xi})^{-1}\Vert\,
n^{-1/2} \Vert \wh{\bm\Lambda}_n-\bm\Lambda_n^{\text{\tiny OLS}}\Vert
\, 
\Vert\mbf F_t\Vert\nn\\
=&\, \{\rm  A.1 + \rm A.2\} \Vert\mbf F_t\Vert, \;\text{ say.} \label{eq:ultimopezzoAAhat}
}
Then, 
\al{
\rm A.1 \le &\, n \Vert
(\bm\Lambda_n^{\prime}
({\bm\Sigma}_n^{\xi})^{-1}
{\bm\Lambda}_n)^{-1}\Vert\,n^{-1}\Vert {\bm\Lambda}_n^\prime
({\bm\Sigma}_n^{\xi})^{-1}
(\bm\Lambda_n-\bm\Lambda_n^{\text{\tiny OLS}})\Vert\nn\\
=&\, n \Vert
(\bm\Lambda_n^{\prime}
({\bm\Sigma}_n^{\xi})^{-1}
{\bm\Lambda}_n)^{-1}\Vert\,n^{-1}\l\Vert 
T^{-1}\sum_{t=1}^T {\bm\Lambda}_n^\prime
({\bm\Sigma}_n^{\xi})^{-1}\bm\xi_{nt}\mbf F_t^\prime\r\Vert\,
\l\Vert \l(T^{-1}\sum_{t=1}^T \mbf F_t\mbf F_t^\prime\r)^{-1}
\r\Vert\nn\\
=&\, O_p(n^{-1/2}T^{-1/2}),\label{eq:ultimopezzoAA1hat}
}
by Lemmas \ref{lem:LSL2}(iii), \ref{lem:COVFF0}(iv), and  \ref{lem:frida}.
Moreover, ${\rm A}.2 = O_p(\max(n^{-1}\log^{2/\delta_v}T,n^{-1/2}T^{-1/2}\sqrt{\log n},T^{-1}))$,
because of \eqref{eq:reggioemilia} and Lemmas  \ref{lem:lambdasqrtn}, \ref{lem:LSL2}(iii), and Assumption \ref{ass:idio}(a) which implies $\Vert({\bm\Sigma}_n^{\xi})^{-1}\Vert\le C_\xi$. This, jointly with \eqref{eq:ultimopezzoAAhat} and \eqref{eq:ultimopezzoAA1hat} implies that
\al{
\rm A=&\,O_p(\max(n^{-1}\log^{2/\delta_v}T,n^{-1/2}T^{-1/2}\sqrt{\log n},T^{-1})),  \label{eq:ultimopezzoAA0hat}
}
since $\Vert \mbf F_t\Vert = O_p(1)$ because $\E[F_{jt}^2]=1$, $j=1,\ldots, r$, by Assumption \ref{ass:ident}(b). 

Second, by Proposition \ref{prop:load}(a) and Lemma \ref{lem:eststar_LASTHAT}(v)
\al{
\rm B=&\,  \Vert
n(\wh{\bm\Lambda}_n^{\prime}(\wh{\bm\Sigma}_n^{\xi})^{-1}\wh{\bm\Lambda}_n)^{-1}n^{-1/2}\wh{\bm\Lambda}_n^{\prime}(\wh{\bm\Sigma}_n^{\xi})^{-1}-n({\bm\Lambda}_n^{\prime}({\bm\Sigma}_n^{\xi})^{-1}{\bm\Lambda}_n)^{-1}n^{-1/2}{\bm\Lambda}_n^{\prime}({\bm\Sigma}_n^{\xi})^{-1}\Vert\, n^{-1/2}\Vert \bm\Lambda_n-\wh{\bm\Lambda}_n\Vert\, \Vert\mbf F_t\Vert\nn\\
=&\, O_p(\max(n^{-2}\log^{4/\delta_v}T,T^{-1} \sqrt{\log n}, n^{-1}T^{-1/2} \log^{2/\delta_v}T\sqrt{\log n})),\label{eq:ultimopezzoBBhat}
}
and again since $\Vert \mbf F_t\Vert = O_p(1)$ because $\E[F_{jt}^2]=1$, $j=1,\ldots, r$, by Assumption \ref{ass:ident}(b).

Third,
\al{
\rm C\le&\, n\Vert ({\bm\Lambda}_n^{\prime}({\bm\Sigma}_n^{\xi})^{-1}{\bm\Lambda}_n)^{-1}\Vert\,  n^{-1}\Vert{\bm\Lambda}_n^{\prime}({\bm\Sigma}_n^{\xi})^{-1}\bm\xi_{nt}\Vert = O_p(n^{-1/2}),\label{eq:ultimopezzoCChat}
}
by Lemmas \ref{lem:LSL2}(iii) and \ref{lem:LSXi}(i).

Fourth, and last, 
\al{
\rm D\le&\,n \Vert 
(\wh{\bm\Lambda}_n^{\prime}(\wh{\bm\Sigma}_n^{\xi})^{-1}\wh{\bm\Lambda}_n)^{-1}-
({\bm\Lambda}_n^{\prime}({\bm\Sigma}_n^{\xi})^{-1}{\bm\Lambda}_n)^{-1} \Vert\, n^{-1}\Vert{\bm\Lambda}_n^{\prime}({\bm\Sigma}_n^{\xi})^{-1}\bm\xi_{nt}
\Vert\nn\\
&+ n\Vert ({\bm\Lambda}_n^{\prime}({\bm\Sigma}_n^{\xi})^{-1}{\bm\Lambda}_n)^{-1}\Vert\,n^{-1}\Vert
\{\wh{\bm\Lambda}_n^{\prime}(\wh{\bm\Sigma}_n^{\xi})^{-1}-
{\bm\Lambda}_n^{\prime}({\bm\Sigma}_n^{\xi})^{-1}\}\bm\xi_{nt}\Vert\nn\\
&+ n\Vert 
(\wh{\bm\Lambda}_n^{\prime}(\wh{\bm\Sigma}_n^{\xi})^{-1}\wh{\bm\Lambda}_n)^{-1}-
({\bm\Lambda}_n^{\prime}({\bm\Sigma}_n^{\xi})^{-1}{\bm\Lambda}_n)^{-1} \Vert\, 
n^{-1}\Vert
\{\wh{\bm\Lambda}_n^{\prime}(\wh{\bm\Sigma}_n^{\xi})^{-1}-
{\bm\Lambda}_n^{\prime}({\bm\Sigma}_n^{\xi})^{-1}\}\bm\xi_{nt}
\Vert\nn\\
=&\, \rm D.1+\rm D.2+\rm D.3, \;\text{say.}\label{eq:ultimopezzoDD00ahat}
}
Then, 
${\rm D.1}= O_p(\max(n^{-3/2}\log^{2/\delta_v}T,n^{-1/2}T^{-1/2}\sqrt{\log n}))$, by Lemmas \ref{lem:LSXi}(i) and \ref{lem:eststar_LASTHAT}(iv). 
Moreover,
\al{
\rm D.2=&\, 
n\Vert ({\bm\Lambda}_n^{\prime}({\bm\Sigma}_n^{\xi})^{-1}{\bm\Lambda}_n)^{-1}\Vert\,\nn\\
&\cdot\Big\{
n^{-1}\Vert
(\wh{\bm\Lambda}_n-\bm\Lambda_n)^{\prime}
({\bm\Sigma}_n^{\xi})^{-1}\bm\xi_{nt}\Vert+
n^{-1}\Vert
\bm\Lambda_n^\prime
[(\wh{\bm\Sigma}_n^{\xi})^{-1}-({\bm\Sigma}_n^{\xi})^{-1}
]\bm\xi_{nt}\Vert+
n^{-1}\Vert
(\wh{\bm\Lambda}_n-\bm\Lambda_n)^{\prime}
[(\wh{\bm\Sigma}_n^{\xi})^{-1}-({\bm\Sigma}_n^{\xi})^{-1}]\bm\xi_{nt}\Vert\Big\}
\nn\\
=&\, n\Vert ({\bm\Lambda}_n^{\prime}({\bm\Sigma}_n^{\xi})^{-1}{\bm\Lambda}_n)^{-1}\Vert\,\{\rm D.2.a+ \rm D.2.b+ \rm D.2.c\},\; \text{ say.}\nn
} 
We then have the following results. First,
\al{
\rm D.2.a=&\, 
n^{-1}\l\Vert\sum_{i=1}^n (\wh{\bm\lambda}_i-\bm\lambda_i) (\sigma_i^2)^{-1} \xi_{it}\r\Vert
\le \max_{i=1,\ldots, n}\Vert \wh{\bm\lambda}_i-\bm\lambda_i\Vert n^{-1}\Vert (\bm\Sigma_n^\xi)^{-1}\bm\xi_{nt}\Vert\nn\\
=&\, O_p(\max(n^{-3/2}\log^{2/\delta_v}T,n^{-1/2}T^{-1/2}\sqrt{\log n})),\label{eq:ultimopezzoDD2ahat}
}
by Lemmas \ref{lem:LSXi}(i) and \ref{lem:sigmaunifhat}(i). Second,
\al{
\rm D.2.b=&\,
n^{-1}\l\Vert\sum_{i=1}^n \bm\lambda_i \{(\wh{\sigma}_i^2)^{-1}-({\sigma}_i^2)^{-1}\} \xi_{it}\r\Vert
=n^{-1}\l\Vert\sum_{i=1}^n \bm\lambda_i (\sigma_i^2\wh{\sigma}_i^2)^{-1}\{\wh{\sigma}_i^2-{\sigma}_i^2\} \xi_{it}\r\Vert\nn\\
\le&\, \max_{i=1,\ldots, n}\vert \wh{\sigma}_i^2-{\sigma}_i^2\vert \, 
\l\{n^{-1}\l\Vert\sum_{i=1}^n \bm\lambda_i \{\sigma_i^2(\wh{\sigma}_i^2-\sigma_i^2)\}^{-1}\xi_{it}\r\Vert
+n^{-1}\l\Vert\sum_{i=1}^n \bm\lambda_i (\sigma_i^4)^{-1}\xi_{it}\r\Vert\r\}\nn\\
=&\, O_p(\max(n^{-3/2}\log^{2/\delta_v}T,n^{-1/2}T^{-1/2}\sqrt{\log n})),
\label{eq:ultimopezzoDD2bhat}
}
by Lemmas \ref{lem:LSXi}(i) and \ref{lem:sigmaunifhat}(ii). Third, clearly by \eqref{eq:ultimopezzoDD2ahat} and \eqref{eq:ultimopezzoDD2bhat}
\al{
{\rm D.2.c} = o_p(\max(n^{-3/2}\log^{2/\delta_v}T,n^{-1/2}T^{-1/2}\sqrt{\log n})).
\label{eq:ultimopezzoDD2chat}
}
By \eqref{eq:ultimopezzoDD2ahat}, \eqref{eq:ultimopezzoDD2bhat}, and \eqref{eq:ultimopezzoDD2chat}, and Lemma \ref{lem:LSL2} we have ${\rm D.2}=
O_p(\max(n^{-3/2}\log^{2/\delta_v}T,n^{-1/2}T^{-1/2}\sqrt{\log n})$. Last, 
${\rm D.3}=O_p(\max(n^{-2}\log^{4/\delta_v}T, T^{-1}\log n)){\rm D.2}$, by Lemma \ref{lem:eststar_LASTHAT}(ii), thus  it is dominated by ${\rm D.2}$.
Therefore,
\beq
{\rm D}=  O_p(\max(n^{-3/2}\log^{2/\delta_v}T,n^{-1/2}T^{-1/2}\sqrt{\log n}). \label{eq:ultimopezzoDDhat}
\eeq
By substituting \eqref{eq:ultimopezzoAA0hat}, \eqref{eq:ultimopezzoBBhat}, \eqref{eq:ultimopezzoCChat}, and \eqref{eq:ultimopezzoDDhat}, into \eqref{eq:ultimopezzohat}
we have
\al{
\Vert \wh{\mbf F}_t^{\text{\tiny WLS}}-\mbf F_t\Vert = O_p(\max(n^{-1/2},T^{-1}\sqrt{\log n} )),\label{eq:topol}
}
which, once substituted into \eqref{eq:ultimopezzoamontehat}, proves part (a.1).

For part (a.2), let 
$\wh{\bm{\mathcal F}}_T^{\text{\tiny KF}}=(\wh{\mbf F}_{1|1} \cdots \wh{\mbf F}_{T|T})^\prime$ 
and $\wh{\bm{\mathcal F}}_T^{\text{\tiny WLS}}=(\wh{\mbf F}_{1}^{\text{\tiny WLS}} \cdots \wh{\mbf F}_{T}^{\text{\tiny WLS}})^\prime$, 
and recall that, by definition,
$\wh{\bm{\mathcal F}}_T=(\wh{\mbf F}_{1|T} \cdots \wh{\mbf F}_{T|T})^\prime=(\wh{\mbf F}_{1} \cdots \wh{\mbf F}_{T})^\prime$.
From \eqref{eq:ultimopezzohat}, we have
\al{
T^{-1/2}\Vert\wh{\bm{\mathcal F}}_{T}-\bm{\mathcal F}_T\Vert \le&\,
T^{-1/2}\Vert\wh{\bm{\mathcal F}}_{T}-\wh{\bm{\mathcal F}}_T^{\text{\tiny KF}}\Vert
+T^{-1/2}\Vert\wh{\bm{\mathcal F}}_{T}^{\text{\tiny KF}}-\wh{\bm{\mathcal F}}_T^{\text{\tiny WLS}}\Vert
+T^{-1/2}\Vert\wh{\bm{\mathcal F}}_{T}^{\text{\tiny WLS}}-{\bm{\mathcal F}}_T\Vert.
\label{eq:novella}
} 
By Lemma \ref{lem:FFOnhatunif} the first two terms on the rhs of \eqref{eq:novella} are such that
\al{
T^{-1/2}\Vert\wh{\bm{\mathcal F}}_{T}-\wh{\bm{\mathcal F}}_T^{\text{\tiny KF}}\Vert&\le \max_{t=1,\ldots, T}\Vert\wh{\mbf F}_{t}-\wh{\mbf F}_{t|t}\Vert = O_p(n^{-1}\log^{1/\delta_v}T),\nn\\
T^{-1/2}\Vert\wh{\bm{\mathcal F}}_{T}^{\text{\tiny KF}}-\wh{\bm{\mathcal F}}_T^{\text{\tiny WLS}}\Vert&\le \max_{t=1,\ldots, T}\Vert\wh{\mbf F}_{t|t}-\wh{\mbf F}_{t}^{\text{\tiny WLS}}\Vert= O_p(n^{-1}\log^{1/\delta_v}T).\label{eq:novella12}
}
While, for the last term on the rhs of \eqref{eq:novella},
letting $\bm{\mathcal E}_{nT}=(\bm\xi_{n1}\cdots \bm\xi_{nT})^\prime$,
 we have
\al{
T^{-1/2}\Vert\wh{\bm{\mathcal F}}_{T}^{\text{\tiny WLS}}-{\bm{\mathcal F}}_T\Vert\le &\, 
n\Vert(\wh{\bm\Lambda}^\prime(\wh{\bm\Sigma}_n^\xi)^{-1}\wh{\bm\Lambda}_n)^{-1} \Vert\nn\\
&\cdot n^{-1}T^{-1/2}
\Big\{
\Vert \wh{\bm\Lambda}^\prime(\wh{\bm\Sigma}_n^\xi)^{-1}(\bm\Lambda_n-\wh{\bm\Lambda}_n)\bm{\mathcal F}_T^\prime\Vert
+ \Vert\{\wh{\bm\Lambda}^\prime(\wh{\bm\Sigma}_n^\xi)^{-1}-{\bm\Lambda}^\prime({\bm\Sigma}_n^\xi)^{-1}\} \bm{\mathcal E}_{nT}^\prime\Vert+\Vert
{\bm\Lambda}^\prime({\bm\Sigma}_n^\xi)^{-1}\bm{\mathcal E}_{nT}^\prime
\Vert
\Big\}\nn\\
=&\, n\Vert(\wh{\bm\Lambda}^\prime(\wh{\bm\Sigma}_n^\xi)^{-1}\wh{\bm\Lambda}_n)^{-1} \Vert \l\{\mathcal A+\mathcal B+\mathcal C \r\}, \;\text{say.}\label{eq:novella3}
}

For the first term on the rhs of \eqref{eq:novella3} we have
\al{
\mathcal A\le&\, 
n^{-1}\Vert  {\bm\Lambda}^\prime({\bm\Sigma}_n^\xi)^{-1}(\bm\Lambda_n-\wh{\bm\Lambda}_n)\Vert\, T^{-1/2}\Vert\bm{\mathcal F}_T\Vert
+n^{-1/2}\Vert \wh{\bm\Lambda}^\prime(\wh{\bm\Sigma}_n^\xi)^{-1}-{\bm\Lambda}^\prime({\bm\Sigma}_n^\xi)^{-1}\Vert \,
n^{-1/2}\Vert \wh{\bm\Lambda}_n-\bm\Lambda_n\Vert \, T^{-1/2}\Vert\bm{\mathcal F}_T\Vert\nn\\
=&\, \{\mathcal A_1+\mathcal A_2\}T^{-1/2}\Vert\bm{\mathcal F}_T\Vert, \;\text{say,}\label{eq:novella3A}
}
where
\al{
\mathcal A_1\le &\, n^{-1} \Vert {\bm\Lambda}^\prime({\bm\Sigma}_n^\xi)^{-1}(\bm\Lambda_n-{\bm\Lambda}_n^{\text{\tiny OLS}})\Vert +
n^{-1/2}\Vert {\bm\Lambda}\Vert \, \Vert ({\bm\Sigma}_n^\xi)^{-1}\Vert \, n^{-1/2}\Vert {\bm\Lambda}_n^{\text{\tiny OLS}}-\wh{\bm\Lambda}_n\Vert\nn\\
=&\, O_p(n^{-1/2}T^{-1/2})+O_p(\max(n^{-1}\log^{2/\delta_v}T, n^{-1/2}T^{-1/2}\sqrt{\log n},T^{-1}))\nn\\
=&\,O_p(\max(n^{-1}\log^{2/\delta_v}T, n^{-1/2}T^{-1/2}\sqrt{\log n},T^{-1})),\label{eq:novella3A1}
}
by \eqref{eq:reggioemilia}, \eqref{eq:ultimopezzoAA1hat}, Lemma \ref{lem:lambdasqrtn}, and Assumption \ref{ass:idio}(a),
while 
\al{
\mathcal A_2=O_p(\max(T^{-1/2}\sqrt{\log n},n^{-1}\log^{2/\delta_v}T))O_p(\max(T^{-1/2},n^{-1}\log^{2/\delta_v}T)),\label{eq:novella3A2}
}
 by Proposition \ref{prop:load}(a) and \ref{lem:eststar_LASTHAT}(ii). By substituting \eqref{eq:novella3A1} and  \eqref{eq:novella3A2} into  \eqref{eq:novella3A}
 \al{
 \mathcal A= O_p(\max(n^{-1}\log^{2/\delta_v}T,T^{-1}\sqrt{\log n},n^{-1/2}T^{-1/2}\sqrt{\log n})).\label{eq:novella3A0}
 }
 
 Moving to the second term on the rhs of \eqref{eq:novella3}, we have
 \al{
 \mathcal B\le &\, n^{-1}T^{-1/2}\Vert(\wh{\bm\Lambda}_n-\bm\Lambda_n)^\prime(\bm\Sigma_n^\xi)^{-1}\bm{\mathcal E}_{nT}^\prime \Vert
 +n^{-1}T^{-1/2}\Vert \bm\Lambda_n^\prime \{ (\wh{\bm\Sigma}_n^\xi)^{-1}-(\bm\Sigma_n^\xi)^{-1} \}\bm{\mathcal E}_{nT}^\prime\Vert\nn\\
 &+ n^{-1}T^{-1/2}\Vert (\wh{\bm\Lambda}_n-\bm\Lambda_n)^\prime\{ (\wh{\bm\Sigma}_n^\xi)^{-1}-(\bm\Sigma_n^\xi)^{-1} \}\bm{\mathcal E}_{nT}^\prime\Vert\nn\\
 =&\, \mathcal B_1+\mathcal B_2+\mathcal B_3, \;\text{say.}\label{eq:novella3B}
 }
Then, considering each term on the rhs of \eqref{eq:novella3B},
 \al{
  \mathcal B_1\le &\,n^{-1}T^{-1/2}\Vert({\bm\Lambda}_n^{\text{\tiny OLS}}-\bm\Lambda_n)^\prime(\bm\Sigma_n^\xi)^{-1}\bm{\mathcal E}_{nT}^\prime \Vert
  +n^{-1/2}\Vert\wh{\bm\Lambda}_n-\bm\Lambda_n^{\text{OLS}}\Vert 
\,\Vert  (\bm\Sigma_n^\xi)^{-1}\Vert n^{-1/2}T^{-1/2}\Vert \bm{\mathcal E}_{nT}^\prime \Vert\nn\\
=&\,   \mathcal B_{1.a}+\mathcal B_{1.b}\;\text{say,}\nn
 }
 where
 \al{
  \mathcal B_{1.a}\le&\,  n^{-1}T^{-1/2}\l\Vert T^{-1} \sum_{t=1}^T\mbf F_t\bm{\xi}_{nt}^\prime(\bm\Sigma_n^\xi)^{-1}\bm{\mathcal E}_{nT}^\prime\r \Vert \, 
  \l\Vert \l(T^{-1}\sum_{t=1}^T\mbf F_t\mbf F_t^\prime\r)^{-1}\r\Vert=  O_p(n^{-1/2}T^{-1/2}),\nn
 }
 by Lemmas \ref{lem:COVFF0}(vi) and \ref{lem:frida}. Whereas $\mathcal B_{1.b}=O_p(\max(n^{-1}\log^{2/\delta_v}T, n^{-1/2}T^{-1/2}\sqrt{\log n},T^{-1}))$ by \eqref{eq:reggioemilia} and Lemma \ref{lem:LSXi}(vi).
Therefore,
 \al{
  \mathcal B_1=O_p(\max(\max(n^{-1}\log^{2/\delta_v}T, n^{-1/2}T^{-1/2}\sqrt{\log n},T^{-1})).\label{eq:novella3B1}
  }
Moreover, letting $\bm\zeta_i=(\xi_{i1}\cdots \xi_{iT})^\prime$,
\al{
\mathcal B_2=&\,n^{-1}T^{-1/2}\Vert \bm\Lambda_n^\prime (\wh{\bm\Sigma}_n^\xi)^{-1}\{\bm\Sigma_n^\xi-\wh{\bm\Sigma}_n^\xi\}(\bm\Sigma_n^\xi)^{-1}\bm{\mathcal E}_{nT}\Vert\nn\\
=&\,n^{-1}T^{-1/2}\l\Vert
\sum_{i=1}^n
\bm\lambda_i\bm\zeta_i^\prime (\wh{\sigma}_i^2\sigma_i^2)^{-1}\{\sigma_{i}^2-\wh{\sigma}_i^2\}
\r\Vert\nn\\
\le &\,C_\xi \max_{i=1,\ldots, n} \vert \sigma_{i}^2-\wh{\sigma}_i^2\vert \l\{\min_{i=1,\ldots, n}\wh{\sigma}_i^2\r\}^{-1}\,
 n^{-1}T^{-1/2} \l\Vert
\sum_{i=1}^n
\bm\lambda_i\bm\zeta_i^\prime\r\Vert\nn\\
\le &\,C_\xi \max_{i=1,\ldots, n} \vert \sigma_{i}^2-\wh{\sigma}_i^2\vert \,\Vert (\wh{\bm\Sigma}_n^\xi)^{-1}\Vert
 n^{-1}T^{-1/2} \Vert
\bm\Lambda_n^\prime\bm{\mathcal E}_{nT}^\prime\Vert\nn\\
=&\, O_p(\max(n^{-1}\log^{2/\delta_v}T, T^{-1/2}\sqrt{\log n}))O_p(n^{-1/2}),\label{eq:novella3B2}
}
by Assumption \ref{ass:idio}(a) and Lemmas \ref{lem:sigmaunifhat}(ii), \ref{lem:sigmaunifhat}(ii), and  
Last,
\al{
\mathcal B_3\le&\, n^{-1/2}\Vert \wh{\bm\Lambda}_n-\bm\Lambda_n\Vert\, \Vert (\wh{\bm\Sigma}_n^\xi)^{-1}-(\bm\Sigma_n^\xi)^{-1} \Vert n^{-1/2}T^{-1/2}\Vert\bm{\mathcal E}_{nT}\Vert\nn\\
=&\, O_p(\max(T^{-1/2}, n^{-1}\log^{2/\delta_v}T))\cdot O_p(\max(T^{-1/2}\sqrt{\log n}, n^{-1}\log^{2/\delta_v}T)), \label{eq:novella3B3}
}
by Proposition \ref{prop:load}(a) and Lemmas \ref{lem:LSXi}(vi) and \ref{lem:sigmaunifhat}(v). By substituting 
\eqref{eq:novella3B1}, \eqref{eq:novella3B2}, and \eqref{eq:novella3B3} into \eqref{eq:novella3B}
\al{
\mathcal B=O_p( n^{-1}\log^{2/\delta_v}T, n^{-1/2}T^{-1/2}\sqrt{\log n},T^{-1}).
\label{eq:novella3B0}
}

Finally, for the third term on the rhs of \eqref{eq:novella3}, we have
\al{
\mathcal C=O_p(n^{-1/2}),\label{eq:novella3C0}
}
by Lemma \ref{lem:LSXi}(v). 

By substituting 
\eqref{eq:novella3A0}, \eqref{eq:novella3B0}, and \eqref{eq:novella3C0} into \eqref{eq:novella3}, and since $n\Vert(\wh{\bm\Lambda}^\prime(\wh{\bm\Sigma}_n^\xi)^{-1}\wh{\bm\Lambda}_n)^{-1} \Vert=O_p(1)$ by Lemma \ref{lem:eststar_LASTHAT}(iii),
\al{
T^{-1/2}\Vert\wh{\bm{\mathcal F}}_{T}^{\text{\tiny WLS}}-{\bm{\mathcal F}}_T\Vert =O_p(\max(n^{-1/2}, T^{-1}\sqrt{\log n})),\nn
}
which, once substituted in \eqref{eq:novella} together with \eqref{eq:novella12}, proves part (a.2).

 Turning to part (b), from \eqref{eq:ultimopezzoamontehat},  \eqref{eq:ultimopezzohat}, \eqref{eq:ultimopezzoCChat} and  \eqref{eq:topol},  if $T^{-1}\sqrt {n\log n}\to 0$, as $n,T\to\infty$, we have 
 \al{
\sqrt n( \wh{\mbf F}_t-\mbf F_t) &= \sqrt n( \wh{\mbf F}_t^{\text{\tiny WLS}}-\mbf F_t) +o_p(1) \nn\\
&= n({\bm\Lambda}_n^\prime({\bm\Sigma}_n^{\xi})^{-1}{\bm\Lambda}_n)^{-1} \l\{n^{-1/2}\sum_{i=1}^n\bm\lambda_i\xi_{it} (\sigma_i^2)^{-1}\r\}+o_p(1).\label{eq:sviluppoF}
 }
 Then, by Assumption \ref{ass:idio}(e), as $n\to\infty$, it holds that
\begin{align}
n^{-1/2}\sum_{i=1}^n\bm\lambda_i\xi_{it} (\sigma_i^2)^{-1}\stackrel{d}{\to}\mathcal N\l(\mbf 0_r,\lim_{n\to\infty} n^{-1} \sum_{i,j=1}^n 
\bm\lambda_i\bm\lambda_j^\prime\E_{}[\xi_{it}\xi_{jt}]
(\sigma_i^2\sigma_j^2)^{-1}
\r).\label{eq:numKSCLT}
\end{align}
Moreover, from  Lemma \ref{lem:LSL2}(iii) 
we have that 
$$
\Vert n^{-1} {\bm\Lambda}_n^\prime({\bm\Sigma}_n^{\xi})^{-1}{\bm\Lambda}_n-\bm\Sigma_{\Lambda\Sigma\Lambda}\Vert=o_p(1),
$$
for some finite and positive definite $r\times r$ matrix $\bm\Sigma_{\Lambda\Sigma\Lambda}$, which, jointly with Lemma \ref{lem:LSL2}(v), implies
\beq
\Vert n({\bm\Lambda}_n^\prime({\bm\Sigma}_n^{\xi})^{-1}{\bm\Lambda}_n)^{-1}- (\bm\Sigma_{\Lambda\Sigma\Lambda})^{-1}\Vert = o_p(1).\label{eq:denKSCLT}
\eeq
From \eqref{eq:sviluppoF}, \eqref{eq:numKSCLT}, and  \eqref{eq:denKSCLT}, by Slutsky's Theorem, we have
\[
\sqrt n (\wh{\mbf F}_t-\mbf F_t)\stackrel{d}{\to}\mathcal N(\mbf 0_r,\bm{\mathcal W}_t),
\]
where
\[
\bm{\mathcal W}_t =(\bm\Sigma_{\Lambda\Sigma\Lambda})^{-1}\l\{\lim_{n\to\infty} n^{-1} \sum_{i,j=1}^n 
\bm\lambda_i\bm\lambda_j^\prime\E_{}[\xi_{it}\xi_{jt}]
(\sigma_i^2\sigma_j^2)^{-1}\r\}(\bm\Sigma_{\Lambda\Sigma\Lambda})^{-1}.
\]
This proves part (b). Part (c) is straightforward. This completes the proof.  $\Box$

\subsection{Proof of Proposition \ref{th:chi}}

First notice that 
\beq
\wh{\chi}_{it}-\chi_{it} = (\wh{\mbf F}_t-\mbf F_t)^\prime \bm\lambda_i + \wh{\mbf F}_t^\prime (\wh{\bm\lambda}_i-\bm\lambda_i)
=(\wh{\mbf F}_t-\mbf F_t)^\prime \bm\lambda_i + \mbf F_t^\prime (\wh{\bm\lambda}_i-\bm\lambda_i)+ (\wh{\mbf F}_t-\mbf F_t)^\prime(\wh{\bm\lambda}_i-\bm\lambda_i).\label{eq:chierror}
\eeq
Then, 
\al{
\vert \wh{\chi}_{it}-\chi_{it} \vert\le&\, \Vert \bm\lambda_i\Vert\, \Vert \wh{\mbf F}_t-\mbf F_t\Vert +\Vert\wh{\bm\lambda}_i-\bm\lambda_i \Vert\, \Vert \mbf F_t\Vert+ \Vert \wh{\bm\lambda}_i-\bm\lambda_i\Vert \, \Vert \wh{\mbf F}_t-\mbf F_t\Vert \nn\\
=&\, O_p(\max(n^{-1/2}, T^{-1}\sqrt{\log n}))+ O_p(\max(T^{-1/2}, n^{-1}\log^{2/\delta_v} T))+ o_p(\max(n^{-1/2},T^{-1/2}))\nn\\
=& O_p(\max(n^{-1/2}, T^{-1/2})),\label{eq:chierror2}
}
by Propositions \ref{prop:load}(a) and \ref{prop:factors}(a), Assumption \ref{ass:common}(a), and since $\Vert \mbf F_t\Vert=O_p(1)$ because $\E[F_{jt}^2]=1$, $j=1,\ldots,r$, by Assumption \ref{ass:ident}(b).
This proves part (a).

For part (b), let us denote $\delta_{nT}=\min(\sqrt n,\sqrt T)$, for simplicity of notation. Consider the first term on the rhs of \eqref{eq:chierror}. Define $\bm K_\Lambda=n(\bm\Lambda_n^\prime(\bm\Sigma_n^\xi)^{-1}\bm\Lambda_n)^{-1}$. Then, from \eqref{eq:sviluppoF} in the proof of Proposition \ref{prop:factors}
\begin{align}
\delta_{n,T} \bm\lambda_i^\prime(\wh{\mbf F}_t-\mbf F_t) 
=&\,\delta_{n,T}\, n^{-1} \bm\lambda_i^\prime\bm K_\Lambda  \sum_{j=1}^n \bm\lambda_j^\prime (\sigma_j^2)^{-1} \xi_{jt} +O_p(\delta_{n,T}T^{-1}\sqrt{\log n})\nn\\
=&\,\delta_{nT}\, n^{-1/2} \mathcal A_{it}+o_p(1),\;\text{ say,}\label{eq:KSsemplice}
\end{align}
since $\delta_{nT}T^{-1}\le \delta_{nT}\max(n^{-1},T^{-1})=\delta_{nT}\delta_{nT}^{-2}\to 0$. Similarly, consider the second term on the rhs of \eqref{eq:chierror} and define $\bm K_F=(T^{-1}\sum_{t=1}^T\mbf F_t\mbf F_t^\prime)^{-1}$. Then, from \eqref{eq:baili} in the proof of Proposition \ref{prop:load}
\begin{align}
\delta_{n,T} \mbf F_t^\prime (\wh{\bm\lambda}_i-\bm\lambda_i)=&\,\delta_{n,T} T^{-1} \mbf F_t^\prime \bm K_F\sum_{s=1}^T\mbf F_s \xi_{is}+O_p(\delta_{n,T}n^{-1}\log^{2/\delta_v} T)\nn\\
=&\,\delta_{nT}\, T^{-1/2} \mathcal B_{it}+o_p(1),\;\text{ say,}\label{eq:MLsemplice}
\end{align}
since $\delta_{nT}n^{-1}\le \delta_{nT}\max(n^{-1},T^{-1})=\delta_{nT}\delta_{nT}^{-2}\to 0$. From Propositions  \ref{prop:load}(b) and \ref{prop:factors}(b), as $n,T\to\infty$,
\begin{align}
&\mathcal A_{it} \stackrel{d}{\to}\mathcal N(0, \mathcal C^F_{it})\quad \text{ and }\quad \mathcal B_{it} \stackrel{d}{\to}\mathcal N(0, \mathcal C^\lambda_{it}),\label{eq:CLTB}
\end{align}
where $\mathcal C^F_{it}=\bm\lambda_i^\prime\bm{\mathcal W}_t\bm\lambda_i$ and $\mathcal C^\lambda_{it}=\mbf F_t^\prime \bm{\mathcal V}_i\mbf F_t$.
Moreover, $\mathcal A_{it}$ and $\mathcal B_{it}$ are asymptotically independent, since the former is a cross-sectional sum of random variables, while the latter is the sum of a given time series and under Lemmas \ref{lem:Gxi}(i)-\ref{lem:Gxi}(iii) are weakly serially and cross-sectionally correlated in the same sense as assumed by \citet[Assumption C]{Bai03}.

Define $a_{nT}=\delta_{nT}n^{-1/2}$ and $b_{nT}=\delta_{nT}T^{-1/2}$. Then, substituting \eqref{eq:KSsemplice} and \eqref{eq:MLsemplice} into \eqref{eq:chierror}, we obtain
\begin{align}
\delta_{nT} (\wh{\chi}_{it}-\chi_{it}) 
&= a_{nT}\mathcal A_{it} + b_{nT} \mathcal B_{it} + \delta_{nT}(\wh{\mbf F}_t-\mbf F_t)^\prime(\wh{\bm\lambda}_i-\bm\lambda_i)+
o_p(1)\nn\\
&= a_{nT}\mathcal A_{it} + b_{nT} \mathcal B_{it} + o_p(\delta_{nT}
\max(n^{-1/2},T^{-1/2}))+
o_p(1)\nn\\
&= a_{nT}\mathcal A_{it} + b_{nT} \mathcal B_{it} + o_p(1),\label{ultimissimo}
\end{align}
because of \eqref{eq:chierror2}.
From \eqref{eq:CLTB} and \eqref{ultimissimo}, by Slutsky's theorem and following the same reasoning as in \citet[proof of Theorem 3]{Bai03}, as $n,T\to\infty$, we have 
\[
\delta_{nT}{\l(a_{nT}^2\mathcal C^F_{it} +  b_{nT}^2\mathcal C^\lambda_{it}\r)^{-1/2}} { (\wh{\chi}_{it}-\chi_{it}) }
=(n^{-1}\mathcal C^F_{it} +  T^{-1}\mathcal C^\lambda_{it})^{-1/2}(\wh{\chi}_{it}-\chi_{it})
\stackrel{d}{\ra}\mathcal N(0,1),
\]
which completes the proof. $\Box$

\subsection{Proof of Proposition \ref{prop:altri}}
For part (a.1), for any $k^*\ge 0$, we have (recall that $\wh{\sigma}_{i}^2\equiv \wh{\sigma}_{i}^{2(k^*+1)}$)
\begin{align}
(\wh{\sigma}_{i}^2-{\sigma}_{i}^2) =&\,  (\wh{\sigma}_{i}^2-\wh{\sigma}_{i}^{2**})+(\wh{\sigma}_{i}^{2**}-\wh{\sigma}_{i}^{2*})+(\wh{\sigma}_{i}^{2*}-{\sigma}_{i}^{2\text{\tiny OLS}})+({\sigma}_{i}^{2\text{\tiny OLS}}-\sigma_i^2)\nn\\
=&\,O_p(\max(n^{-1}\log^{2/\delta_v}T,T^{-1/2}))\nn\\
&+O_p(\max(n^{-1}\log^{2/\delta_v}T,n^{-1/2}T^{-1/2}\sqrt{\log n},T^{-1}))\nn\\
&+O_p(\max(n^{-1}\log^{2/\delta_v}T,n^{-1/2}T^{-1/2},T^{-1}))
+O_p(T^{-1/2}), 
\nn
\end{align}
by Lemmas  \ref{lem:olsT12}(iii), \ref{lem:starols}(iii), \ref{lem:localglobal}(ii), and \ref{lem:convEMiAH}(i).

For part (a.2), for any $k^*\ge 0$, we have (recall that $\wh{\bm\Sigma}_{n}^\xi\equiv \wh{\bm\Sigma}_{n}^{\xi(k^*+1)}$)
\begin{align}
(\wh{\bm\Sigma}_{n}^\xi-{\bm\Sigma}_{n}^\xi) =&\,  (\wh{\bm\Sigma}_{n}^\xi-\wh{\bm\Sigma}_{n}^{\xi**})+(\wh{\bm\Sigma}_{n}^{\xi**}-\wh{\bm\Sigma}_{n}^{\xi*})+(\wh{\bm\Sigma}_{n}^{\xi*}-\bm\Sigma_n^\xi),\nn 
\end{align}
and the proof follows from Lemma  \ref{lem:sigmaunif}(iii) and since
\al{
\Vert \wh{\bm\Sigma}_{n}^{\xi**}-\wh{\bm\Sigma}_{n}^{\xi*}\Vert&= \l(\max_{i=1,\ldots, n} \vert \wh{\sigma}_{i}^{2**}-\wh{\sigma}_{i}^{2*}\vert^2\r)^{1/2}\le 
\max_{i=1,\ldots, n} \vert \wh{\sigma}_{i}^{2**}-\wh{\sigma}_{i}^{2*}\vert\nn\\
&= O_p(\max(n^{-1}\log^{2/\delta_v}T,n^{-1/2}T^{-1/2}\sqrt{\log n},T^{-1})),\nn
}
by Lemma \ref{lem:localglobal}(i), and
\al{
\Vert \wh{\bm\Sigma}_{n}^{\xi}-\wh{\bm\Sigma}_{n}^{\xi**}\Vert&= \l(\max_{i=1,\ldots, n} \vert \wh{\sigma}_{i}^{2}-\wh{\sigma}_{i}^{2**}\vert^2\r)^{1/2}\le 
\max_{i=1,\ldots, n} \vert \wh{\sigma}_{i}^{2}-\wh{\sigma}_{i}^{2**}\vert\nn\\
&= O_p(\max(n^{-1}\log^{2/\delta_v}T,T^{-1/2}\sqrt{\log n})),\nn
}
by Lemma \ref{lem:convEM1unif}(ii). 

For part (a.3), for any $k^*\ge 0$, we have (recall that $\wh{\mbf A}\equiv \wh{\mbf A}^{(k^*+1)}$)
\begin{align}
(\wh{\mbf A}-{\mbf A}) =&\,  (\wh{\mbf A}-\wh{\mbf A}^{**})+(\wh{\mbf A}^{**}-\wh{\mbf A}^{*})+(\wh{\mbf A}^{*}-\mbf A^{\text{\tiny OLS}})
+(\mbf A^{\text{\tiny OLS}}-\mbf A)\nn\\
=&\,O_p(\max(n^{-1}\log^{2/\delta_v}T,T^{-1/2}))\nn\\
&+O_p(n^{-1}\log^{2/\delta_v}T)
+O_p(n^{-1}\log^{2/\delta_v}T)
+O_p(T^{-1/2}),\label{eq:tagliacozzo}
\end{align}
by Lemmas
\ref{lem:olsT12}(iv),
\ref{lem:starolsVAR}(i),
\ref{lem:localglobal}(iii),
and
\ref{lem:convEMiAH}(ii). 

For part (a.4), for any $k^*\ge 0$, we have (recall that $\wh{\bm\Gamma}^v\equiv \wh{\bm\Gamma}^{v(k^*+1)}$)
\begin{align}
(\wh{\bm\Gamma}^v-{\bm\Gamma}^v) =&\,  (\wh{\bm\Gamma}^v-\wh{\bm\Gamma}^{v**})+(\wh{\bm\Gamma}^{v**}-\wh{\bm\Gamma}^{v*})+(\wh{\bm\Gamma}^{v*}-\bm\Gamma^{v\text{\tiny OLS}})
+(\bm\Gamma^{v\text{\tiny OLS}}-\bm\Gamma^v)\nn\\
=&\,O_p(\max(n^{-1}\log^{2/\delta_v}T,T^{-1/2}))\nn\\
&+O_p(n^{-1}\log^{2/\delta_v}T)
+O_p(n^{-1}\log^{2/\delta_v}T)
+O_p(T^{-1/2}),\nn
\end{align}
by Lemmas
\ref{lem:olsT12}(v),
\ref{lem:starolsVAR}(ii),
\ref{lem:localglobal}(iv),
and
\ref{lem:convEMiAH}(iii). This completes the proof of part (a). 

To prove part (b), we need sharper rates and, to this end, we use the closed form expressions of the estimated parameters in \eqref{eq:param3}, \eqref{eq:param4}, and \eqref{eq:paramGv}.

Start with part (b.1). From  \eqref{eq:param3}, for any $k^*\ge 0$, we have
\al{
\wh{\sigma}_i^2 = &\, T^{-1}\sum_{t=1}^T (x_{it}-\wh{\bm\lambda}_i^\prime\mbf F_t)^2\nn\\
&+  \wh{\bm\lambda}_i^\prime 
\bigg\{
2 T^{-1} \sum_{t=1}^T(\mbf F_{t|T}^{(k^*)}-\mbf F_t)\mbf F_t^\prime
+T^{-1} \sum_{t=1}^T(\mbf F_{t|T}^{(k^*)}-\mbf F_t)(\mbf F_{t|T}^{(k^*)}-\mbf F_t)^\prime
+T^{-1} \sum_{t=1}^T \mbf P_{t|T}^{(k^*)} 
\bigg\}\wh{\bm\lambda}_i\nn\\
&-2\wh{\bm\lambda}_i^\prime T^{-1}\sum_{t=1}^T (\mbf F_{t|T}^{(k^*)}-\mbf F_t)x_{it}\nn\\
= &\, T^{-1}\sum_{t=1}^T (x_{it}-{\bm\lambda}_i^{\text{\tiny OLS}\prime}\mbf F_t)^2+ 
(\wh{\bm\lambda}_i-{\bm\lambda}_i^{\text{\tiny OLS}})^\prime
\l\{T^{-1}\sum_{t=1}^T \mbf F_t\mbf F_t^\prime\r\}
(\wh{\bm\lambda}_i-{\bm\lambda}_i^{\text{\tiny OLS}})
-2(\wh{\bm\lambda}_i-{\bm\lambda}_i^{\text{\tiny OLS}})^\prime T^{-1}\sum_{t=1}^T \mbf F_t (x_{it}-{\bm\lambda}_i^{\text{\tiny OLS}\prime}\mbf F_t) \nn\\
&+  \wh{\bm\lambda}_i^\prime 
\bigg\{
2 T^{-1} \sum_{t=1}^T(\mbf F_{t|T}^{(k^*)}-\mbf F_t)\mbf F_t^\prime
+T^{-1} \sum_{t=1}^T(\mbf F_{t|T}^{(k^*)}-\mbf F_t)(\mbf F_{t|T}^{(k^*)}-\mbf F_t)^\prime
+T^{-1} \sum_{t=1}^T \mbf P_{t|T}^{(k^*)} 
\bigg\}\wh{\bm\lambda}_i\nn\\
&-2\wh{\bm\lambda}_i^\prime T^{-1}\sum_{t=1}^T (\mbf F_{t|T}^{(k^*)}-\mbf F_t)x_{it}
.\label{eq:esselunga0settembre}
}
Therefore, since by construction $T^{-1}\sum_{t=1}^T \mbf F_t (x_{it}-{\bm\lambda}_i^{\text{\tiny OLS}\prime}\mbf F_t) =0$, from \eqref{eq:esselunga0settembre}
\al{
\l\vert \wh{\sigma}_i^2  - T^{-1}\sum_{t=1}^T (x_{it}-{\bm\lambda}_i^{\text{\tiny OLS}\prime}\mbf F_t)^2\r\vert\le&\,  
\Vert\wh{\bm\lambda}_i-{\bm\lambda}_i^{\text{\tiny OLS}}\Vert^2 \,\l\Vert T^{-1}\sum_{t=1}^T \mbf F_t\mbf F_t^\prime \r\Vert
\nn\\
&+\Vert \wh{\bm\lambda}_i\Vert^2 
\bigg\{2\bigg\Vert T^{-1} \sum_{t=1}^T(\mbf F_{t|T}^{(k^*)}-\mbf F_t)\mbf F_t^\prime\bigg\Vert+\bigg\Vert T^{-1} \sum_{t=1}^T \mbf P_{t|T}^{(k^*)} 
\bigg\Vert\nn\\
&+\bigg\Vert T^{-1} \sum_{t=1}^T(\mbf F_{t|T}^{(k^*)}-\mbf F_t)(\mbf F_{t|T}^{(k^*)}-\mbf F_t)^\prime\bigg\Vert\bigg\}+ 2 \Vert \wh{\bm\lambda}_i\Vert \, \bigg\Vert T^{-1}\sum_{t=1}^T (\mbf F_{t|T}^{(k^*)}-\mbf F_t)x_{it}\bigg\Vert\nn\\
=&\, \Vert\wh{\bm\lambda}_i-{\bm\lambda}_i^{\text{\tiny OLS}}\Vert^2 \,\l\Vert T^{-1}\sum_{t=1}^T \mbf F_t\mbf F_t^\prime \r\Vert+  \Vert \wh{\bm\lambda}_i\Vert^2 \{2\mathcal A+\mathcal B+\mathcal C\} + 2 \Vert \wh{\bm\lambda}_i\Vert \mathcal D, \;\text{say.}\label{eq:esselunga}
}
For term $\mathcal A$  on the rhs of \eqref{eq:esselunga} we have
\al{
\mathcal A = O_p(\max(n^{-1}\log^{2/\delta_v}T, n^{-1/2}T^{-1/2}\sqrt{\log n}, T^{-1}\sqrt{\log n})),\label{eq:esselungaA}
}
by Lemma \ref{lem:ultimo}(i).  Term $\mathcal B$ is dominated by  term $\mathcal A$. For term $\mathcal C$ on the rhs of \eqref{eq:esselunga} we have
\al{
\mathcal C  \le \max_{t=1,\ldots, T} \Vert \mbf P_{t|T}^{(k^*)} \Vert = O_p(n^{-1}), \label{eq:esselungaC}
}
by Lemma \ref{lem:PPOnhat}(iv) when $k^*\ge 1$ and by Lemma \ref{lem:PPOn} when $k^*= 0$.
For term $\mathcal D$  on the rhs of \eqref{eq:esselunga} we have
\al{
\mathcal D \le&\, \bigg\Vert T^{-1}\sum_{t=1}^T (\mbf F_{t|T}^{(k^*)}-\mbf F_t)\mbf F_t\bigg\Vert\, \Vert \bm\lambda_i\Vert + \bigg\Vert T^{-1}\sum_{t=1}^T (\mbf F_{t|T}^{(k^*)}-\mbf F_t)\xi_{it}\bigg\Vert\nn\\
=&\, O_p(\max(n^{-1}\log^{2/\delta_v}T, n^{-1/2}T^{-1/2}\sqrt{\log n}, T^{-1}\sqrt{\log n})),\label{eq:esselungaD}
}
by Lemmas \ref{lem:ultimo}(i) and \ref{lem:ultimo}(ii), and Assumption \ref{ass:common}(a).

Now by substituting \eqref{eq:esselungaA}, \eqref{eq:esselungaC}, and \eqref{eq:esselungaD} into \eqref{eq:esselunga}, we have
\al{
\l\vert \wh{\sigma}_i^2  - T^{-1}\sum_{t=1}^T (x_{it}-{\bm\lambda}_i^{\text{\tiny OLS}\prime}\mbf F_t)^2\r\vert
=&\, \vert \wh{\sigma}_i^2  -\sigma_i^{2\text{\tiny OLS}}\vert
 = \Vert\wh{\bm\lambda}_i-{\bm\lambda}_i^{\text{\tiny OLS}}\Vert^2 \,\l\Vert T^{-1}\sum_{t=1}^T \mbf F_t\mbf F_t^\prime \r\Vert\label{eq:esselunga0}\\ 
&+O_p(\max(n^{-1}\log^{2/\delta_v}T, n^{-1/2}T^{-1/2}\sqrt{\log n}, T^{-1}\sqrt{\log n}))\nn\\
=&\,O_p(\max(n^{-1}\log^{2/\delta_v}T, n^{-1/2}T^{-1/2}\sqrt{\log n}, T^{-1}\sqrt{\log n})) ,\nn
}
by \eqref{eq:loadingsconsistenti}-\eqref{eq:ELLE3} in the proof of Proposition \ref{prop:load}, Lemma \ref{lem:consistCOV}(i) combined with Assumption \ref{ass:ident}(b), and 
since $\Vert \wh{\bm\lambda}_i\Vert\le \Vert \wh{\bm\lambda}_i-\bm\lambda_i\Vert+\Vert {\bm\lambda}_i\Vert=O_p(1)$ by Proposition \ref{prop:load} and Assumption \ref{ass:common}(a).

Therefore, from \eqref{eq:esselunga0} and Lemma \ref{lem:olsT12}(i),
 if $n^{-1}\sqrt T\log^{2/\delta_v} T\to 0$, as $n,T\to\infty$, we have 
 \al{
\sqrt T(\wh{\sigma}_i^2-\sigma_i^2) =&\, \sqrt T({\sigma}_i^{2\text{\tiny OLS}}-\sigma_i^2) +o_p(1)\nn\\
=&\, T^{-1/2}\sum_{t=1}^T \xi_{it}^2+\sqrt T (\bm\lambda_i^{\text{\tiny OLS}}-\bm\lambda_i)^\prime\l\{T^{-1}\sum_{t=1}^T \mbf F_t\mbf F_t^{\prime}\r\}(\bm\lambda_i^{\text{\tiny OLS}}-\bm\lambda_i)-\sqrt T\sigma_i^2+o_p(1)\nn\\
=&\, T^{-1/2}\sum_{t=1}^T \{\xi_{it}^2-\sigma_i^2\}+o_p(1),\label{eq:homelate2}
}
where we used also \eqref{eq:homelate} in the proof of Lemma  \ref{lem:olsT12}. 
Now, by Assumption \ref{ass:idio}(c) and \citet[Theorem 14.1]{davidson}, we have that $\{\xi_{it}^2-\sigma_i^2\}$ is strongly mixing with exponentially decaying coefficients
and such that by Assumption \ref{ass:tails} $\sup_{m\ge 1} r^{-1/\delta_\xi} (\E[\vert \xi_{it}^2-\sigma_i^2\vert^m])^{1/m}\le K_1$ 
for some finite positive real $K_1$  independent of $t$ and $i$ \citep[Section 2]{KC18}. Then,  the Central Limit Theorem by \citet[Theorem 1.7]{ibra62} applies, i.e.,
\al{
T^{-1/2}\sum_{t=1}^T \{\xi_{it}^2-\sigma_i^2\}\stackrel{d}{\ra} \mathcal N(0, \E[\xi_{it}^4]-\sigma_i^4).\label{eq:homelate3}
}
From \eqref{eq:homelate2}, \eqref{eq:homelate3}, and Slutsky's Theorem, and by noticing that the excess kurtosis is given by $\kappa_i = \E[\xi_{it}^4]/\sigma_i^4-3$, so that
$\E[\xi_{it}^4]=\sigma_i^4(\kappa_i+3)$, we prove part (b.1). 

For part (b.2), from  \eqref{eq:param4}, for any $k^*\ge 0$, we have
\begin{align}
\Vert \wh{\mbf A}&\,-\mbf A^{\text{\tiny OLS}}\Vert
\le 
\l\Vert T^{-1} \sum_{t=2}^T\l({\mbf F}_{t|T}^{(k^*)}{\mbf F}_{t-1|T}^{(k^*)\prime}+{\mbf C}^{(k^*)}_{t,t-1|T} \r)-
T^{-1}\sum_{t=2}^T{\mbf F}_{t}{\mbf F}_{t-1}^{\prime}\r\Vert\, 
\l\Vert\l(T^{-1}\sum_{t=2}^T \mbf F_{t-1}\mbf F_{t-1}^\prime \r)^{-1}\r\Vert\nn\\
&+\l\Vert \l\{T^{-1}\sum_{t=2}^T\l({\mbf F}_{t-1|T}^{(k^*)}{\mbf F}_{t-1|T}^{(k^*)\prime}+{\mbf P}^{(k^*)}_{t-1|T} \r)\r\}^{-1}-
\l\{T^{-1}\sum_{t=2}^T{\mbf F}_{t-1}{\mbf F}_{t-1}^{\prime}\r\}^{-1}\r\Vert\, 
\l\Vert T^{-1}\sum_{t=2}^T \mbf F_{t}\mbf F_{t-1}^\prime\r\Vert\nn\\
&+\l\Vert T^{-1}\sum_{t=2}^T\l({\mbf F}_{t|T}^{(k^*)}{\mbf F}_{t-1|T}^{(k^*)\prime}+{\mbf C}^{(k^*)}_{t,t-1|T} \r)-
T^{-1}\sum_{t=2}^T{\mbf F}_{t}{\mbf F}_{t-1}^{\prime}\r\Vert\nn\\
&\cdot \l\Vert \l\{T^{-1}\sum_{t=2}^T\l({\mbf F}_{t-1|T}^{(k^*)}{\mbf F}_{t-1|T}^{(k^*)\prime}+{\mbf P}^{(k^*)}_{t-1|T} \r)\r\}^{-1}-
\l\{T^{-1}\sum_{t=2}^T{\mbf F}_{t-1}{\mbf F}_{t-1}^{\prime}\r\}^{-1}\r\Vert\nn\\
=&\, \mathfrak{A}+ \mathfrak{B}+ \mathfrak{C},\;\text{say.}\label{eq:frakABC}
\end{align}
Now, 
\al{
\mathfrak{A}\le &\, \l\{\l\Vert T^{-1}\sum_{t=2}^T {\mbf F}_{t|T}^{(k^*)}{\mbf F}_{t-1|T}^{(k^*)\prime}-T^{-1}\sum_{t=2}^T{\mbf F}_{t}{\mbf F}_{t-1}^{\prime}\r\Vert\, 
+\max_{t=2,\ldots, T}\Vert {\mbf C}^{(k^*)}_{t,t-1|T}\Vert\r\} \l\Vert\l(T^{-1}\sum_{t=2}^T \mbf F_{t-1}\mbf F_{t-1}^\prime \r)^{-1}\r\Vert\nn\\
=&\,O_p(\max(n^{-1}\log^{2/\delta_v}T, n^{-1/2}T^{-1/2}\sqrt{\log n},T^{-1}\sqrt{\log n})),\label{eq:frakA}
}
and
\al{
\mathfrak B\le&\, 
\l\{\l\Vert T^{-1} \sum_{t=2}^T{\mbf F}_{t-1|T}^{(k^*)}{\mbf F}_{t-1|T}^{(k^*)\prime} -
T^{-1} \sum_{t=2}^T{\mbf F}_{t-1}{\mbf F}_{t-1}^{\prime}\r\Vert + \max_{t=2,\ldots, T} \Vert{\mbf P}^{(k^*)}_{t-1|T}\Vert \r\}\nn\\
&\cdot
\l\Vert\l\{T^{-1} \sum_{t=2}^T\l({\mbf F}_{t-1|T}^{(k^*)}{\mbf F}_{t-1|T}^{(k^*)\prime}+{\mbf P}^{(k^*)}_{t-1|T}\r)
\r\}^{-1}\r\Vert\,
\l\Vert\l(T^{-1}\sum_{t=2}^T{\mbf F}_{t-1}{\mbf F}_{t-1}^{\prime}\r)^{-1}\r\Vert
\,\l\Vert T^{-1}\sum_{t=2}^T \mbf F_{t}\mbf F_{t-1}^\prime\r\Vert\nn\\
=&\,O_p(\max(n^{-1}\log^{2/\delta_v}T, n^{-1/2}T^{-1/2}\sqrt{\log n},T^{-1}\sqrt{\log n})),\label{eq:frakB}
}
where we used: Lemma \ref{lem:ultimo}(i), Lemma \ref{lem:PPOnhat}(iv) when $k^*\ge 1$ or Lemma \ref{lem:PPOn} when $k^*= 0$ which can be both applied to $\mbf C_{t,t-1|T}^{(k^*)}$ can be obtained by the upper right block of $\mbf P_{t|T}^{(k^*)}$ when this is computed from the the Kalman smoother having the augmented state vector $(\mbf F_t^\prime\,\mbf F_{t-1}^\prime)^\prime$,
and we also used  
the fact that  $\Vert(T^{-1}\sum_{t=2}^T \mbf F_{t-1}\mbf F_{t-1}^\prime )^{-1}\Vert=O_p(1)$ by Lemma \eqref{lem:frida}. 
For $\mathfrak B$ we also used Lemma \ref{lem:consistCOV}(i) combined with the fact that $\Vert \bm\Gamma_1^F\Vert\le 1$ by Cauchy-Schwartz inequality and Assumption \ref{ass:ident}(b). Clearly $\mathfrak C$ is dominated by $\mathfrak A$ and $\mathfrak B$.

By substituting \eqref{eq:frakA} and \eqref{eq:frakB} into \eqref{eq:frakABC}, we have
\al{
\Vert \wh{\mbf A}-\mbf A^{\text{\tiny OLS}}\Vert 
=&\,O_p(\max(n^{-1}\log^{2/\delta_v}T, n^{-1/2}T^{-1/2}\sqrt{\log n},T^{-1}\sqrt{\log n})).\nn
}
Therefore,
 if $n^{-1}\sqrt T\log^{2/\delta_v} T\to 0$, as $n,T\to\infty$, we have
\al{
&\sqrt T(\text{\upshape vec}(\wh{\mbf A})-\text{\upshape vec}(\mbf A))=\sqrt T(\text{\upshape vec}({\mbf A}^{\text{\tiny OLS}})-\text{\upshape vec}(\mbf A))+o_p(1),\nn
}
and the proof of part (b.2) follows directly from \citet[Proposition 11.2]{Hamilton} and Slutsky's Theorem.

For part (b.3), following a decomposition analogous to the one in \eqref{eq:esselunga0settembre}, from  \eqref{eq:paramGv}, for any $k^*\ge 0$, we have
\al{
\wh{\bm\Gamma}^v =&\, 
T^{-1}\sum_{t=2}^T (\mbf F_t-\wh{\mbf A} \mbf F_{t-1})(\mbf F_t-{\mbf A}^{\text{\tiny OLS}} \mbf F_{t-1})^\prime+
( \wh{\mbf A}-{\mbf A}^{\text{\tiny OLS}}) \l\{T^{-1}\sum_{t=2}^T \mbf F_{t-1}\mbf F_{t-1}^\prime\r\}( \wh{\mbf A}-{\mbf A}^{\text{\tiny OLS}})^\prime
\nn\\
&-( \wh{\mbf A}-{\mbf A}^{\text{\tiny OLS}}) T^{-1}\sum_{t=2}^T \mbf F_{t-1} (\mbf F_{t}-{\mbf A}^{\text{\tiny OLS}\prime}\mbf F_{t-1})^\prime 
- T^{-1}\sum_{t=2}^T  (\mbf F_{t}-{\mbf A}^{\text{\tiny OLS}\prime}\mbf F_{t-1}) \mbf F_{t-1}^\prime ( \wh{\mbf A}-{\mbf A}^{\text{\tiny OLS}})^\prime
\nn\\
 &+T^{-1}\sum_{t=2}^T ({\mbf F}_{t|T}^{(k^*)}{\mbf F}_{t|T}^{(k^*)\prime}+ {\mbf P}^{(k^*)}_{t|T}) -T^{-1}\sum_{t=2}^T \mbf F_t\mbf F_t^\prime
+\wh{\mbf A} \l\{
T^{-1}\sum_{t=2}^T ({\mbf F}_{t-1|T}^{(k^*)}{\mbf F}_{t-1|T}^{(k^*)\prime}+ {\mbf P}^{(k^*)}_{t-1|T}) -T^{-1}\sum_{t=2}^T \mbf F_{t-1}\mbf F_{t-1}^\prime
\r\}
\wh{\mbf A}^\prime\nn\\
&- \wh{\mbf A} \l\{T^{-1}\sum_{t=2}^T ({\mbf F}_{t-1|T}^{(k^*)}{\mbf F}_{t|T}^{(k^*)\prime}+{\mbf C}^{(k^*)\prime}_{t,t-1|T})-T^{-1}\sum_{t=2}^T {\mbf F}_{t-1}{\mbf F}_{t}^{\prime}\r\}\nn\\
&- \l\{T^{-1}\sum_{t=2}^T ({\mbf F}_{t-1|T}^{(k^*)}{\mbf F}_{t|T}^{(k^*)\prime}+{\mbf C}^{(k^*)}_{t,t-1|T})- T^{-1}\sum_{t=2}^T{\mbf F}_{t}{\mbf F}_{t-1}^{\prime}\r\}\wh{\mbf A}^\prime.\nn
}
Therefore, since by construction $T^{-1}\sum_{t=2}^T \mbf F_{t-1} (\mbf F_{t}-{\mbf A}^{\text{\tiny OLS}\prime}\mbf F_{t-1})^\prime =T^{-1}\sum_{t=2}^T  (\mbf F_{t}-{\mbf A}^{\text{\tiny OLS}\prime}\mbf F_{t-1}) \mbf F_{t-1}^\prime =\mbf 0_{r\times r}$, by using again Lemmas \ref{lem:ultimo}
 and Lemma \ref{lem:PPOnhat}(iv) when $k^*\ge 1$ or Lemma \ref{lem:PPOn} when $k^*= 0$ which can be both applied to $\mbf C_{t,t-1|T}^{(k^*)}$, as argued above, 
and using also \eqref{eq:tagliacozzo}, we obtain
\al{
\Vert \wh{\bm\Gamma}^v-\bm\Gamma^{v\text{\tiny OLS}}\Vert =&\,O_p(\max(n^{-1}\log^{2/\delta_v}T, n^{-1/2}T^{-1/2}\sqrt{\log n},T^{-1}\sqrt{\log n})).\nn
}
Therefore,
 if $n^{-1}\sqrt T\log^{2/\delta_v} T\to 0$, as $n,T\to\infty$, we have
\al{
&\sqrt T(\text{\upshape vech}(\wh{\bm\Gamma}^v)-\text{\upshape vech}(\bm\Gamma^v))=\sqrt T(\text{\upshape vech}({\bm\Gamma}^{v\text{\tiny OLS}})-\text{\upshape vech}(\bm\Gamma^v))+o_p(1),\nn
}
and the proof of part (b.3) follows directly from \citet[Proposition 11.2]{Hamilton} and Slutsky's Theorem. This completes the proof.
$\Box$ 
%
\subsection{Proof of Corollary \ref{cor:matriciAN}}

For part (a), for ease of notation and without loss of generality let $s(j)=j$ for all $j=1,\ldots, \bar n$. 
Then, from \eqref{eq:loadingsconsistenti}, since $\bar n$ is finite,
\begin{align}
\sqrt T (\text{\upshape vec}(\wh{\bm\Lambda}_{\bar n})- \text{\upshape vec}({\bm\Lambda}_{\bar n}))=\l\{\mbf I_{\bar n}\otimes \l(T^{-1}\sum_{t=1}^T\mbf F_t\mbf F_t^\prime\r)^{-1}\r\}
\l(T^{-1/2}\sum_{t=1}^T \text{vec}(\mbf F_t\bm\xi_{\bar n t}^\prime)\r)
+o_p(1).\label{eq:matL1}
\end{align} 
Moreover, since $\bar n$ is finite we can still apply the Central Limit Theorem by \citet[Theorem 1.4]{ibra62} so that, as $T\to\infty$,
\beq
T^{-1/2}\sum_{t=1}^T \text{vec}(\mbf F_t\bm\xi_{\bar n t}^\prime) \stackrel{d}{\to}\mathcal N(\mbf 0_{\bar n r}, \bm{\Sigma}_{\bar n}),\label{eq:matL2}
\eeq
with
\begin{align}
\bm\Sigma_{\bar n}&=\lim_{T\to\infty}\E_{}\l[\l\{
T^{-1/2}\sum_{t=1}^T
\l(
\ba{c}
\mbf F_t \xi_{1t}\\
\vdots\\
\mbf F_t \xi_{\bar nt}
\ea
\r)
\r\}
\l\{
T^{-1/2}\sum_{t=1}^T
\l(
\ba{c}
\mbf F_t \xi_{1t}\\
\vdots\\
\mbf F_t \xi_{\bar nt}
\ea
\r)
\r\}^\prime\, \r]\nn\\
&=\lim_{T\to\infty} T^{-1}\sum_{t,s=1}^T \E\l[\bm\xi_{\bar n t}\bm\xi_{\bar n s}^\prime\otimes \mbf F_t\mbf F_s^\prime \r]=\lim_{T\to\infty} T^{-1}\sum_{t,s=1}^T \E\l[\bm\xi_{\bar n t}\bm\xi_{\bar n s}^\prime\r]\otimes \E\l[\mbf F_t\mbf F_s^\prime \r],\label{eq:squalo}
\end{align}
where we used the fact that $\{\mbf F_t\}$ and $\{\bm\xi_{nt}\}$ are independent processes because of Lemma \ref{lem:fidio}.
The proof of part (a.1) follows from Lemmas \ref{lem:consistCOV}(i) and \ref{lem:frida}, \eqref{eq:matL1}, \eqref{eq:matL2}, and Slutsky's theorem.  For part (a.2) just notice that, if $\E[(\bm\xi_{n1}^\prime\cdots\bm\xi_{nT}^\prime)^\prime(\bm\xi_{n1}^\prime\cdots\bm\xi_{nT}^\prime)]=\mbf I_T\otimes \bm\Sigma_n^\xi$ for all $n,T\in\mathbb N$, then  
$\E\l[\bm\xi_{\bar n t}\bm\xi_{\bar n s}^\prime\r]= \mbf 0_{\bar n\times\bar n}$ if $t\ne s$ while $\E\l[\bm\xi_{\bar n t}\bm\xi_{\bar n t}^\prime\r]= \bm\Sigma_{\bar n}^\xi$. By substituting into \eqref{eq:squalo} we get $\bm\Sigma_{\bar n}=\bm\Sigma_{\bar n}^\xi\otimes \bm\Gamma^F$, and
$\bm{\mathcal V}_{\bar n}=(\mbf I_{\bar n}\otimes\bm\Gamma^F)^{-1}(\bm\Sigma_{\bar n}^\xi\otimes \bm\Gamma^F)(\mbf I_{\bar n}\otimes\bm\Gamma^F)^{-1}=\bm\Sigma_{\bar n}^\xi\otimes(\bm\Gamma^F)^{-1}$, thus proving part (a.2).

Turning to part (b), for ease of notation and without loss of generality let $s(j)=j$ for all $j=1,\ldots, \bar T$. 
Then, from  \eqref{eq:sviluppoF} in the proof of Proposition \ref{prop:factors}, since $\bar T$ is finite
\begin{align}
\sqrt n (\wh{\bm F}_{\bar T}-\bm F_{\bar T})=
\l(\mbf I_{\bar T}\otimes n({\bm\Lambda}_n^\prime({\bm\Sigma}_n^{\xi})^{-1}{\bm\Lambda}_n)^{-1}\r)
\l(\mbf I_{\bar T}\otimes n^{-1/2}{\bm\Lambda}_n^\prime({\bm\Sigma}_n^{\xi})^{-1}\r)\bm \Xi_{n\bar T}+o_p(1).\label{eq:matF1}
\end{align}
where $\bm \Xi_{n\bar T}=(\bm\xi_{n1}^\prime\cdots\bm\xi_{n\bar T}^\prime)^\prime$. Moreover, since $\bar T$ is finite we can still apply the Central Limit Theorem in Assumption \ref{ass:idio}(e), so that as $n\to\infty$
\begin{align}
\l(\mbf I_{\bar T}\otimes n^{-1/2}{\bm\Lambda}_n^\prime({\bm\Sigma}_n^{\xi})^{-1}\r)\bm \Xi_{n\bar T}
=n^{-1/2}\sum_{i=1}^n\text{vec}\l(\bm\lambda_i (\sigma_i^2)^{-1}(\xi_{i1}\cdots\xi_{i\bar T})\r)
\stackrel{d}{\to}\mathcal N(\mbf 0_{\bar Tr}, \bm{\Omega}_{\bar T}),\label{eq:matF2}
\end{align}
with
\begin{align}
\bm{\Omega}_{\bar T}&=\lim_{n\to\infty} n^{-1} \l(\mbf I_{\bar T}\otimes {\bm\Lambda}_n^\prime({\bm\Sigma}_n^{\xi})^{-1}\r)\E_{}[\bm \Xi_{n\bar T}\bm \Xi_{n\bar T}^\prime]
\l(\mbf I_{\bar T}\otimes ({\bm\Sigma}_n^{\xi})^{-1}{\bm\Lambda}_n\r)\nn\\
&=\lim_{n\to\infty}\E_{}\l[
\l\{
n^{-1/2}\sum_{i=1}^n
\l(\ba{c}
\bm\lambda_i(\sigma_i^2)^{-1}\xi_{i1}\\
\vdots\\
\bm\lambda_i(\sigma_i^2)^{-1}\xi_{iT}
\ea
\r)
\r\}
\l\{
n^{-1/2}\sum_{i=1}^n
\l(\ba{c}
\bm\lambda_i(\sigma_i^2)^{-1}\xi_{i1}\\
\vdots\\
\bm\lambda_i(\sigma_i^2)^{-1}\xi_{iT}
\ea
\r)
\r\}^\prime\,
\r]\nn\\
&=\lim_{n\to\infty}\frac 1n\sum_{i,j=1}^n\l\{\E_{}[\bm\zeta_{i\bar T}\bm\zeta_{j\bar T}^\prime]\otimes 
\l({\bm\lambda_i\bm\lambda_j^\prime}(\sigma_i^2\sigma_j^2)^{-1}
\r)\r\},\nn
\end{align}
and recall that $\bm\zeta_{i\bar T}=(\xi_{i1}\cdots\xi_{i\bar T})^\prime$. The proof of part (b.1) follows from \eqref{eq:denKSCLT} in the proof of Proposition \ref{prop:factors}, \eqref{eq:matF1}, \eqref{eq:matF2}, and Slutsky's theorem. For part (b.2) just notice that,  if $\E[(\bm\zeta_{1T}^\prime\cdots \bm\zeta_{nT}^\prime)^\prime(\bm\zeta_{1T}^\prime\cdots \bm\zeta_{nT}^\prime)]=\bm\Sigma_n^\xi\otimes\mbf I_T$ for all $n,T\in\mathbb N$, then  
$\E_{}[\bm\zeta_{i\bar T}\bm\zeta_{j\bar T}^\prime]= \mbf 0_{\bar T\times\bar T}$ if $i\ne j$ while $\E_{}[\bm\zeta_{i\bar T}\bm\zeta_{i\bar T}^\prime]= \sigma_i^2\mbf I_{\bar T}$. By substituting into \eqref{eq:squalo} we get $\bm\Omega_{\bar T}=\mbf I_{\bar T}\otimes \bm\Sigma_{\Lambda\Sigma\Lambda}$, and
$\bm{\mathcal W}_{\bar T}=(\mbf I_{\bar T}\otimes \bm\Sigma_{\Lambda\Sigma\Lambda})^{-1}(\mbf I_{\bar T}\otimes \bm\Sigma_{\Lambda\Sigma\Lambda})(\mbf I_{\bar T}\otimes \bm\Sigma_{\Lambda\Sigma\Lambda})^{-1}=\mbf I_{\bar T}\otimes (\bm\Sigma_{\Lambda\Sigma\Lambda})^{-1}$, thus proving part (b.2). 
This completes the proof. $\Box$

\subsection{Proof of Proposition \ref{prop:eff}}

Under Assumptions \ref{ass:common}, \ref{ass:idio}, \ref{ass:ind}, and \ref{ass:ident}, from \citet[Theorem 1]{MBPCAQML}, the asymptotic covariances of the PC estimator of the loadings is
\al{
\bm{\mathcal V}_i^{\text{\tiny PC}}&=
(\bm\Gamma^F)^{-1}
\l(\lim_{T\to\infty}\frac 1 T\sum_{t,s=1}^T\E_{}[\xi_{it}\xi_{is}]\E[\mbf F_t\mbf F_s^\prime]\r)
(\bm\Gamma^F)^{-1}
,\label{eq:pcavarL}
}
where $\bm\Gamma^F=\mbf I_r$, by Assumption \ref{ass:ident}(b). Therefore, the expression of $\bm{\mathcal V}_i^{\text{\tiny PC}}$ in \eqref{eq:pcavarL}
 coincides with the expression of $\bm{\mathcal V}_i$  for the asymptotic covariances of the EM estimator given  in Proposition \ref{prop:load}(b). This proves part (a).

Turning to part (b). Since, because of Proposition \ref{prop:factors}, $\lim_{n,T\to\infty} \sqrt n\E[(\wh{\mbf F}_t-\mbf F_t)]=\mbf 0_r$, then
\al{
\bm{\mathcal W}_t = \lim_{n,T\to\infty} n\Cov(\wh{\mbf F}_t-\mbf F_t,\wh{\mbf F}_t-\mbf F_t)&=
\lim_{n,T\to\infty}\l\{n\E[(\wh{\mbf F}_t-\mbf F_t)(\wh{\mbf F}_t-\mbf F_t)^\prime]- n \E[(\wh{\mbf F}_t-\mbf F_t)] \E[(\wh{\mbf F}_t-\mbf F_t)^\prime]\r\}\nn\\
&=\lim_{n,T\to\infty}\l\{n\E[(\wh{\mbf F}_t-\mbf F_t)(\wh{\mbf F}_t-\mbf F_t)^\prime]\r\}\label{eq:mcWEM}\\
&=\lim_{n\to\infty} n\{(\bm\Lambda_n^\prime(\bm\Sigma_n^{\xi})^{-1}\bm\Lambda_n)^{-1} \bm\Lambda_n^\prime(\bm\Sigma_n^{\xi})^{-1} \bm\Gamma_n^\xi (\bm\Sigma_n^{\xi})^{-1}\bm\Lambda_n(\bm\Lambda_n^\prime(\bm\Sigma_n^{\xi})^{-1}\bm\Lambda_n)^{-1}\}.\nn
}
%
Similarly, because of Lemma \ref{lem:PCAF},  $\lim_{n,T\to\infty} \sqrt n\E[(\wt{\mbf F}_t-\mbf F_t)]=\mbf 0_r$, which implies
\al{
\bm{\mathcal W}_t^{\text{\tiny PC}} =\lim_{n,T\to\infty} n\Cov(\wt{\mbf F}_t-\mbf F_t,\wt{\mbf F}_t-\mbf F_t)&=
\lim_{n,T\to\infty}\l\{n\E[(\wt{\mbf F}_t-\mbf F_t)(\wt{\mbf F}_t-\mbf F_t)^\prime]-n \E[(\wt{\mbf F}_t-\mbf F_t)] \E[(\wt{\mbf F}_t-\mbf F_t)^\prime]\r\}\nn\\
&=\lim_{n,T\to\infty}\l\{n\E[(\wt{\mbf F}_t-\mbf F_t)(\wt{\mbf F}_t-\mbf F_t)^\prime]\r\}\label{eq:mcWPC}\\
&= \lim_{n\to\infty} n\{(\mbf M_n^\chi)^{-1} \bm\Lambda_n^\prime \bm\Gamma_n^\xi \bm\Lambda_n(\bm\Lambda_n^\prime\bm\Lambda_n)^{-1}\}\nn\\
&= \lim_{n\to\infty} n\{(\bm\Lambda_n^\prime\bm\Lambda_n)^{-1} \bm\Lambda_n^\prime \bm\Gamma_n^\xi \bm\Lambda_n(\bm\Lambda_n^\prime\bm\Lambda_n)^{-1}\},\nn
}
where we used the fact that $\lim_{n\to\infty} n^{-1}\mbf M_n^\chi=\lim_{n\to
\infty} n^{-1}\bm\Lambda_n^\prime\bm\Lambda_n$, because of Assumption \ref{ass:ident}(b).

Moreover,  if $\sqrt{n\log n}/T\to 0$ (which implies also $\sqrt n/T\to0$), Proposition \ref{prop:factors} (see in particular \eqref{eq:sviluppoF} in its proof) and Lemma \ref{lem:PCAF} (see in particular \eqref{eq:sviluppoFPC} in its proof) jointly imply that
\al{
\l(\ba{c}
\sqrt n(\wh{\mbf F}_t-\mbf F_t)\\
\sqrt n(\wt{\mbf F}_t-\mbf F_t)
\ea
\r)&=\l(\ba{cc}
\sqrt n(\bm\Lambda_n^\prime(\bm\Sigma_n^{\xi})^{-1}\bm\Lambda_n)^{-1} \bm\Lambda_n^\prime(\bm\Sigma_n^{\xi})^{-1}\bm \xi_{nt}\\
\sqrt n(\mbf M_n^\chi)^{-1} \bm\Lambda_n^\prime\bm\xi_{nt}
\ea
\r)+o_p(1)\stackrel{p}{\to}\mathcal N\l(
\mbf 0_{2r},\l(
\ba{cc}
\bm{\mathcal W}_t &\bm{\mathcal U}_t \\
\bm{\mathcal U}_t^\prime&\bm{\mathcal W}_t^{\text{\tiny PC}} 
\ea\r)
\r),\nn
}
where
\al{
\bm{\mathcal U}_t=\lim_{n\to\infty} n\Cov(\wh{\mbf F}_t-\mbf F_t,\wt{\mbf F}_t-\mbf F_t)&=
\lim_{n,T\to\infty}\l\{n\E[(\wh{\mbf F}_t-\mbf F_t)(\wt{\mbf F}_t-\mbf F_t)^\prime]-n \E[(\wh{\mbf F}_t-\mbf F_t)] \E[(\wt{\mbf F}_t-\mbf F_t)^\prime]\r\}\nn\\
&=\lim_{n,T\to\infty}\l\{n\E[(\wh{\mbf F}_t-\mbf F_t)(\wt{\mbf F}_t-\mbf F_t)^\prime]\r\}\label{eq:mcWEMWPC}\\
&= \lim_{n\to\infty} n\{(\bm\Lambda_n^\prime(\bm\Sigma_n^\xi)^{-1}\bm\Lambda_n)^{-1} \bm\Lambda_n^\prime (\bm\Sigma_n^\xi)^{-1}\bm\Gamma_n^\xi \bm\Lambda_n(\bm\Lambda_n^\prime\bm\Lambda_n)^{-1}\}.\nn
}

From \eqref{eq:mcWEM} and \eqref{eq:mcWPC} it follows that
\begin{align}
\bm{\mathcal W}^{\text{\tiny PC}}_t- \bm{\mathcal W}_t=&\,\lim_{n,T\to\infty}\l\{n\,\E_{}[(\wt{\mbf F}_t-\wh{\mbf F}_t)(\wt{\mbf F}_t-\wh{\mbf F}_t)^\prime]\r\}\nn\\
&+\lim_{n,T\to\infty}\l\{n\,\E_{}[(\wh{\mbf F}_t-\mbf F_t)(\wt{\mbf F}_t-\wh{\mbf F}_t)^\prime]\r\}
+\lim_{n,T\to\infty}\l\{n\,\E_{}[(\wt{\mbf F}_t-\wh{\mbf F}_t)(\wh{\mbf F}_t-\mbf F_t)^\prime]\r\}\nn\\
=&\,\bm{\mathcal C}_{1}+\bm{\mathcal C}_{2}+\bm{\mathcal C}_{2}^\prime, \; \text{say.}\label{eq:target}
\end{align}
Let us now define $\bm{\mathcal O}_n^\xi=\bm\Gamma_n^\xi-\bm\Sigma_n^\xi$, which is the $n\times n$ matrix of off-diagonal entries of the idiosyncratic covariance. Because of \eqref{eq:mcWEM}, \eqref{eq:mcWEMWPC}, and
Lemma \ref{lem:lambda} we can write
\begin{align}
\bm{\mathcal C}_{2}=&\,\lim_{n,T\to\infty}\l\{n\,\E_{}[(\wh{\mbf F}_t-\mbf F_t)(\wt{\mbf F}_t-\wh{\mbf F}_t)^\prime]\r\}\nn\\
=&\,\lim_{n,T\to\infty}\l\{n\E_{}[(\wh{\mbf F}_t-\mbf F_t)(\wt{\mbf F}_t-\mbf F_t)^\prime]-
n\E_{}[(\wh{\mbf F}_t-\mbf F_t)(\wh{\mbf F}_t-\mbf F_t)^\prime]\r\}= \bm{\mathcal U}_t-\bm{\mathcal W}_t
\label{eq:C2pd}\\
=&\,\lim_{n\to\infty}\l\{ n(\bm\Lambda_n^\prime(\bm\Sigma_n^\xi)^{-1}\bm\Lambda_n)^{-1}
\bm\Lambda_n^\prime(\bm\Sigma_n^\xi)^{-1}\bm\Gamma_n^\xi 
\l[\bm\Lambda_n(\bm\Lambda_n^\prime\bm\Lambda_n)^{-1}-(\bm\Sigma_n^\xi)^{-1}\bm\Lambda_n(\bm\Lambda_n^\prime(\bm\Sigma_n^\xi)^{-1}\bm\Lambda_n)^{-1}\bm\Lambda_n^\prime \r]\r\}\nn\\
=&\,\lim_{n\to\infty}\l\{ n(\bm\Lambda_n^\prime(\bm\Sigma_n^\xi)^{-1}\bm\Lambda_n)^{-1}
\bm\Lambda_n^\prime(\bm\Sigma_n^\xi)^{-1}\bm{\mathcal O}_n^\xi 
\l[\bm\Lambda_n(\bm\Lambda_n^\prime\bm\Lambda_n)^{-1}-(\bm\Sigma_n^\xi)^{-1}\bm\Lambda_n(\bm\Lambda_n^\prime(\bm\Sigma_n^\xi)^{-1}\bm\Lambda_n)^{-1}\bm\Lambda_n^\prime \r]\r\}\nn\\
=&\,(\bm\Sigma_{\Lambda\Sigma\Lambda})^{-1}\bm\Sigma_{\Lambda}^{1/2}
\lim_{n\to\infty}\l\{\mbf V_n^{\chi\prime}(\bm\Sigma_n^\xi)^{-1}\bm{\mathcal O}_n^\xi \mbf V_n^\chi\r\}
\bm\Sigma_\Lambda^{-1/2}\nn\\
&-
(\bm\Sigma_{\Lambda\Sigma\Lambda})^{-1}\bm\Sigma_{\Lambda}^{1/2}
\lim_{n\to\infty}\l\{\mbf V_n^{\chi\prime}(\bm\Sigma_n^\xi)^{-1}\bm{\mathcal O}_n^\xi (\bm\Sigma_n^\xi)^{-1}\mbf V_n^\chi\r\}
\bm\Sigma_\Lambda^{1/2}(\bm\Sigma_{\Lambda\Sigma\Lambda})^{-1}.\nn
\end{align}
Now, let $\bm{\mathfrak V}_n=\mbf V_n^{\chi\prime}(\bm\Sigma_n^\xi)^{-1}\bm{\mathcal O}_n^\xi (\bm\Sigma_n^\xi)^{-1}\mbf V_n^\chi$, then, for any $h,k=1,\ldots,r$,
\begin{align}
[\bm{\mathfrak V}_n]_{hk}&=\!\!\!\sum_{i,j,\ell,m=1}^n 
[\mbf V_n^{\chi\prime}]_{hi}\, 
[(\bm\Sigma_n^\xi)^{-1}]_{ij}\,
[\bm{\mathcal O}_n^\xi]_{j\ell}
[(\bm\Sigma_n^\xi)^{-1}]_{\ell m}\,
[\mbf V_n^{\chi}]_{mk}\nn\\
&=
\sum_{\substack{i,\ell=1\\ i\ne \ell}}^n 
[\mbf V_n^{\chi}]_{ih}\, 
[(\bm\Sigma_n^\xi)^{-1}]_{ii}\,
[\bm{\mathcal O}_n^\xi]_{i\ell}
[(\bm\Sigma_n^\xi)^{-1}]_{\ell \ell}\,
[\mbf V_n^{\chi}]_{\ell k}.\label{eq:VDODV}
\end{align}

%
Now, let $\bm\iota_{ni}$ be the $n$-dimensional vector with $i$th entry equal one and all other equal zero, and let $\bm s_j$ the $r$-dimensional vector with $j$th entry equal one and all other equal zero. Then, since $\mbf V_n^{\chi} =\bm\Gamma_n^\chi\mbf V_n^{\chi}(\mbf M_n^\chi)^{-1}$, there exists a finite positive integer $\bar n$, such that, for all $n \ge \bar n$
\begin{align}
\max_{i=1,\ldots,n}\max_{j=1,\ldots,r}\vert v_{ij}^\chi \vert&=\max_{i=1,\ldots,n}\max_{j=1,\ldots,r} \vert\bm \iota_{ni}^\prime \mbf V_n^\chi \bm s_j \vert\nn\\
&=\max_{i=1,\ldots,n}\max_{j=1,\ldots,r}\vert \bm \iota_{ni}^\prime\bm\Gamma_n^\chi\mbf V_n^{\chi}(\mbf M_n^\chi)^{-1}\bm s_j\vert\nn\\
&\le \max_{i=1,\ldots,n}\Vert\bm\iota_{ni}^\prime\bm\Gamma_n^\chi\Vert\,\Vert \mbf V_n^{\chi}\Vert\,\max_{j=1,\ldots,r} \vert (\mu_{jn}^\chi)^{-1}\vert\nn\\
&=\max_{i=1,\ldots,n}\Vert \bm\lambda_i^\prime \bm\Gamma^F \bm\Lambda_n^\prime\Vert  {n^{-1}\underline C_r^{-1}}\nn\\
&\le \max_{i=1,\ldots,n} \Vert\bm\lambda_i\Vert\, \Vert \bm\Gamma^F\Vert\, \Vert \bm\Lambda_n\Vert{n^{-1}\underline C_r^{-1}}\nn\\
&\le M_\lambda^2 M_F {n^{-1/2}\underline C_r^{-1}},\label{eq:ppevec}
\end{align}
where we used Assumptions \ref{ass:common}(a) and \ref{ass:common}(b), and Lemmas \ref{lem:Gxi}(iv) and \ref{lem:lambdasqrtn}.

Therefore, from \eqref{eq:VDODV} and \eqref{eq:ppevec}, and Assumption \ref{ass:idio}(a),
\begin{align}
\max_{h,k=1,\ldots r}\vert [\bm{\mathfrak V}_n]_{hk}\vert& \le 
\l(\max_{h=1,\ldots r}\max_{i=1,\ldots n}\vert [\mbf V_n^{\chi}]_{ih}\vert\r)^2 
\l(\max_{i=1,\ldots n}(\sigma_i^2)^{-1}\r)^2
\sum_{i,\ell=1}^n\vert [\bm{\mathcal O}_n^\xi]_{i\ell}\vert\nn\\
&\le  {C_v^2 }\l\{n\l(\min_{i=1,\ldots n}\sigma_i^2\r)^2\r\}^{-1} \sum_{i,\ell=1}^n\vert [\bm{\mathcal O}_n^\xi]_{i\ell}\vert
\le n^{-1} {C_v^2 C_\xi^{2}}\sum_{i,\ell=1}^n\vert [\bm{\mathcal O}_n^\xi]_{i\ell}\vert.\label{eq:VDODV2}
\end{align}
Moreover, since by assumption we have
\[
\lim_{n\to\infty}n^{-1} \sum_{i,\ell=1}^n\vert [\bm{\mathcal O}_n^\xi]_{i\ell}\vert=\lim_{n\to\infty}n^{-1} \sum_{\substack{i,\ell=1\\ i\ne \ell}}^n\vert [\bm{\Gamma}_n^\xi]_{i\ell}\vert=0,
\]
then, from \eqref{eq:VDODV2}, we have $\lim_{n\to\infty}\max_{h,k=1,\ldots r}\vert [\bm{\mathfrak V}_n]_{hk}\vert=0$.
Therefore,
\beq
\lim_{n\to\infty}\l\{\mbf V_n^{\chi\prime}(\bm\Sigma_n^\xi)^{-1}\bm{\mathcal O}_n^\xi (\bm\Sigma_n^\xi)^{-1}\mbf V_n^\chi\r\}=\mbf 0_{r\times r}.\label{eq:VDODV3}
\eeq
Following the same reasoning it also holds that
\beq
\lim_{n\to\infty}\l\{\mbf V_n^{\chi\prime}(\bm\Sigma_n^\xi)^{-1}\bm{\mathcal O}_n^\xi \mbf V_n^\chi\r\}=\mbf 0_{r\times r}.\label{eq:VDODV4}
\eeq
By substituting, \eqref{eq:VDODV3} and \eqref{eq:VDODV4} into \eqref{eq:C2pd}, we obtain 
$\bm{\mathcal C}_2=\mbf 0_{r\times r}.$
Finally,
\al{
\bm{\mathcal C}_1=&\, \lim_{n,T\to\infty}\l\{n\,\E_{}[(\wt{\mbf F}_t-\wh{\mbf F}_t)(\wt{\mbf F}_t-\wh{\mbf F}_t)^\prime]\r\}\nn\\
=&\,
\lim_{n,T\to\infty}\l\{n\,\E_{}[(\wh{\mbf F}_t-{\mbf F}_t)(\wh{\mbf F}_t-{\mbf F}_t)^\prime]\r\}+
 \lim_{n,T\to\infty}\l\{n\,\E_{}[(\wt{\mbf F}_t-{\mbf F}_t)(\wt{\mbf F}_t-{\mbf F}_t)^\prime]\r\}\nn\\
&-\lim_{n,T\to\infty}\l\{n\,\E_{}[(\wt{\mbf F}_t-{\mbf F}_t)(\wh{\mbf F}_t-{\mbf F}_t)^\prime]\r\}-
\lim_{n,T\to\infty}\l\{n\,\E_{}[(\wh{\mbf F}_t-{\mbf F}_t)(\wt{\mbf F}_t-{\mbf F}_t)^\prime]\r\}\nn\\
=&\, \bm{\mathcal W}_t+\bm{\mathcal W}_t^{\text{\tiny PC}}-(\bm{\mathcal U}_t+\bm{\mathcal U}_t^\prime).\nn
}
Therefore, $\bm{\mathcal C}_1$ is positive definite since $\bm{\mathcal W}_t$ and $\bm{\mathcal W}_t^{\text{\tiny PC}}$ are positive definite by Assumptions \ref{ass:common}(a), \ref{ass:idio}(a), and \ref{ass:idio}(f),
and, from \eqref{eq:target} we prove part (ii). This completes the proof. $\Box$

%
%

\setcounter{lem}{0}
\numberwithin{lem}{section}
\section{General lemmas}

\begin{lem}\label{lem:Gxi}
Under Assumptions \ref{ass:common} and \ref{ass:idio}:
\begin{compactenum} 
\item [(i)] for all $n\in\mathbb N$ and $T\in\mathbb N$, $(nT)^{-1}\sum_{i,j=1}^n\sum_{t,s=1}^T \vert\E_{}[\xi_{it}\xi_{js}]\vert \le M_1$, for some finite positive real $M_1$ independent of $n$ and $T$;
\item [(ii)] for all $n\in\mathbb N$ and $t\in\mathbb Z$, $n^{-1}\sum_{i,j=1}^n \vert\E_{}[\xi_{it}\xi_{jt}]\vert \le M_2$, for some finite positive real $M_2$ independent of $n$ and $t$;
\item [(iii)] for all $i\in\mathbb N$ and $T\in\mathbb N$, $T^{-1}\sum_{t,s=1}^T \vert\E_{}[\xi_{it}\xi_{is}]\vert \le M_3$, for some finite positive real $M_3$ independent of $i$ and $T$;
\item [(iv)] 
for all $j=1,\ldots,r$, $\underline C_j\!\le \lim\inf_{n\to\infty} n^{-1}{ \mu_{jn}^\chi}\le\lim\sup_{n\to\infty} n^{-1}{\mu_{jn}^\chi}\le\! \overline C_j$, $\!\!$
for some finite positive reals $\underline C_j$ and $\overline C_j$;
\item [(v)] for all $n\in\mathbb N$, $\mu_{1n}^\xi = \Vert\bm\Gamma_n^\xi\Vert \le M_2$, where $M_2$ is defined in part (ii);
\item [(vi)] for all $j=1,\ldots,r$, $\underline C_j\le \lim\inf_{n\to\infty} n^{-1}{ \mu_{jn}^x}\le\lim\sup_{n\to\infty} n^{-1}{\mu_{jn}^\chi}\le \overline C_j$, and for all $n\in\mathbb N$, $\mu_{r+1,n}^x \le M_2$, where $M_2$ is defined in part (ii).
\end{compactenum}
\end{lem}

\noindent\textsc{Proof.} Using Assumptions \ref{ass:idio}(a) and \ref{ass:idio}(b), we have:
\begin{align}
(nT)^{-1}\sum_{i,j=1}^n\sum_{t,s=1}^T \vert\E_{}[\xi_{it}\xi_{js}]\vert &=n^{-1}\sum_{i,j=1}^n\sum_{k=-(T-1)}^{T-1} \l(1-\frac{\vert k\vert}{T}\r) \vert\E_{}[\xi_{it}\xi_{j,t-k}]\vert\nn\\
&\le  n^{-1}\sum_{i,j=1}^n\sum_{k=-\infty}^{\infty} \rho^{\vert k\vert} M_{ij}\nn\\
&=  n^{-1}\sum_{i=1}^n\sum_{k=-\infty}^{\infty} \rho^{\vert k\vert} M_{ii}+n^{-1}\sum_{i,j=1, i\ne j}^n\sum_{k=-\infty}^{\infty} \rho^{\vert k\vert} M_{ij}\nn\\
&=  n^{-1}\sum_{i=1}^n\sum_{k=-\infty}^{\infty} \rho^{\vert k\vert} \sigma_i^2+\max_{i=1,\ldots,n}\sum_{j=1, j\ne i}^n\sum_{k=-\infty}^{\infty} \rho^{\vert k\vert} M_{ij}\nn\\
&\le \frac{C_\xi(1+\rho)}{1-\rho}+\frac{M_\xi(1+\rho)}{1-\rho}.\nn
\end{align}
Similarly, 
\begin{align}
n^{-1}\sum_{i,j=1}^n \vert\E_{}[\xi_{it}\xi_{jt}]\vert &\le n^{-1}\sum_{i=1}^n \sigma_i^2 +\max_{i=1,\ldots,n}\sum_{j=1, j\ne i}^n M_{ij}\le C_\xi+M_\xi,\nn
\end{align}
and 
\begin{align}
T^{-1}\sum_{t,s=1}^T \vert\E_{}[\xi_{it}\xi_{is}]\vert &= \sum_{k=-(T-1)}^{T-1} \l(1-\frac{\vert k\vert}{T}\r) \vert\E_{}[\xi_{it}\xi_{i,t-k}]\vert\nn\\
&\le \sum_{k=-\infty}^{\infty} \rho^{\vert k\vert} M_{ii} \le \frac{1+\rho}{1-\rho} \sigma_{i}^2\le \frac{C_\xi (1+\rho)}{1-\rho}.\nn
\end{align}
Defining, $M_1=\frac{(C_\xi+M_\xi)(1+\rho)}{1-\rho}$, $M_2=C_\xi+M_\xi$, and $M_3= \frac{C_\xi (1+\rho)}{1-\rho}$, we prove parts (i), (ii), and (iii).

For part (iv), first notice that, for  all $n\in\mathbb N$, the $r$ non-zero eigenvalues of $\bm\Gamma_n^\chi$ are also the $r$ eigenvalues of $n^{-1}\bm\Lambda_n^\prime\bm\Lambda_n\bm\Gamma^F$. Thus, by \citet[Theorem 7]{MK04}, for all $j=1,\ldots, r$ and all $n\in\mathbb N$, we have
\beq\nn
n^{-1}{\nu^{(r)}(\bm\Lambda_n^\prime\bm\Lambda_n)}\ \nu^{(j)}(\bm\Gamma^F) \le n^{-1}{\mu_{jn}^\chi} \le n^{-1}{\nu^{(j)}(\bm\Lambda_n^\prime\bm\Lambda_n)}\ \nu^{(1)}(\bm\Gamma^F).
\eeq
The proof then follows from Assumptions \ref{ass:common}(a) and \ref{ass:common}(b). Indeed, by continuity of eigenvalues, Assumption \ref{ass:common}(a) implies that, for all $j=1,\ldots, r$, and all $n>N_0$
\beq
n^{-1}{\nu^{(j)}(\bm\Lambda_n^\prime\bm\Lambda_n)}= \nu^{(j)}(\bm\Sigma_\Lambda),\label{eq:arna}
\eeq
and there exist finite positive reals $m_\lambda$ and $M_\lambda$ such that $0<m_{\lambda}^2\le \nu^{(r)}(\bm\Sigma_\Lambda)\le \nu^{(1)}(\bm\Sigma_\Lambda)\le M_{\lambda}^2<\infty$.
Similarly, Assumption \ref{ass:common}(b), implies that there exist finite positive reals $m_F$ and $M_F$ such that  $0<m_{F}\le \nu^{(r)}(\bm\Gamma^F)\le \nu^{(1)}(\bm\Gamma^F)\le M_F<\infty$. 

For part (v), by Assumptions \ref{ass:idio}(a) and \ref{ass:idio}(b):
\begin{align}
\Vert\bm\Gamma_n^\xi\Vert\le \max_{i=1,\ldots,n}\sum_{j=1}^n \vert \E[\xi_{it}\xi_{jt}]\vert=\max_{i=1,\ldots,n}\sigma_i^2+\max_{i=1,\ldots,n} \sum_{j=1,j\ne i}^n M_{ij} \le  C_\xi+ M_{\xi}=M_2.\nn
\end{align}
Part (vi) follows from parts (iv) and (v) and Weyl's inequality \citep[Theorem 1]{MK04}. This completes the proof. $\Box$


\begin{lem}\label{lem:lambdasqrtn}
Under Assumption \ref{ass:common}, as $n\to\infty$, $n^{-1/2}\Vert\bm\Lambda_n\Vert=O(1)$.
\end{lem}

\noindent
\noindent\textsc{Proof.} 
By definition and using Assumption \ref{ass:common}(a),
$$
\lim_{n\to\infty}  n^{-1/2}\Vert\bm\Lambda_n\Vert = \lim_{n\to\infty}  n^{-1/2} \sqrt{ \nu^{(1)}(\bm\Lambda_n^\prime\bm\Lambda_n)} =\sqrt{\nu^{(1)}(\bm\Sigma_\Lambda)}\le M_\lambda.
$$ 
This completes the proof. $\Box$

\begin{lem}\label{lem:LSL2} 
For any $r\times r$ symmetric and positive definite matrix $\bm P$ with $\Vert\bm P\Vert\le C_P$ for some finite positive real $C_P$, under Assumptions \ref{ass:common} and \ref{ass:idio}, as $n\to\infty$,
\begin{compactenum}[(i)]
\item $n\Vert (\bm\Lambda_n^\prime(\bm\Sigma_n^\xi)^{-1}\bm\Lambda_n+\bm P^{-1})^{-1}\Vert=O(1)$;
\item $n\Vert (\bm\Lambda_n^\prime(\bm\Gamma_n^\xi)^{-1}\bm\Lambda_n+\bm P^{-1})^{-1}\Vert=O(1)$;
\item [(iii)] $n\Vert (\bm\Lambda_n^\prime(\bm\Sigma_n^\xi)^{-1}\bm\Lambda_n)^{-1}\Vert=O(1)$;
\item [(iv)] $n\Vert (\bm\Lambda_n^\prime(\bm\Gamma_n^\xi)^{-1}\bm\Lambda_n)^{-1}\Vert=O(1)$;
\item [(v)] $n^{-1}\Vert \bm\Lambda_n^\prime(\bm\Sigma_n^\xi)^{-1}\bm\Lambda_n\Vert=O(1)$;
\item [(vi)] $n^{-1}\Vert\bm\Lambda_n^\prime(\bm\Gamma_n^\xi)^{-1}\bm\Lambda_n\Vert=O(1)$;
\item [(vii)] $n^{-1/2}\Vert \bm\Lambda_n^\prime(\bm\Sigma_n^\xi)^{-1}\Vert=O(1)$;
\item [(viii)] $n^{-1/2}\Vert \bm\Lambda_n^\prime(\bm\Gamma_n^\xi)^{-1}\Vert=O(1)$.
\end{compactenum}
\end{lem}

\noindent
\textsc{Proof.} By \citet[Theorems 1, which is Weyl's inequality, and 7]{MK04}
\begin{align}
n\Vert (\bm\Lambda_n^\prime(\bm\Sigma_n^\xi)^{-1}\bm\Lambda_n+\bm P^{-1})^{-1}\Vert& = \frac n{\nu^{(r)}(\bm\Lambda_n^\prime(\bm\Sigma_n^\xi)^{-1}\bm\Lambda_n+\bm P^{-1}) }\le \frac n{\nu^{(r)}(\bm\Lambda_n^\prime(\bm\Sigma_n^\xi)^{-1}\bm\Lambda_n)+\nu^{(r)}(\bm P^{-1})}\nn\\
&\le \frac n{\nu^{(r)}(\bm\Lambda_n\bm\Lambda_n^\prime)\nu^{(n)}((\bm\Sigma_n^\xi)^{-1})+[\nu^{(1)}(\bm P)]^{-1}}\nn\\
&\le \frac n{\nu^{(r)}(\bm\Lambda_n^\prime\bm\Lambda_n)[\nu^{(1)}(\bm\Sigma_n^\xi)]^{-1}+[\nu^{(1)}(\bm P)]^{-1}}\nn\\
&\le \frac n{\nu^{(r)}(\bm\Lambda_n^\prime\bm\Lambda_n) C_\xi^{-1}+C_P^{-1}}\nn\\
&= \frac 1{n^{-1}\nu^{(r)}(\bm\Lambda_n^\prime\bm\Lambda_n) C_\xi^{-1}+n^{-1}C_P^{-1}},\label{eq:pulce}
\end{align}
by Assumption \ref{ass:idio}(a) and since $\bm P$ finite by assumption.  Then, as shown in \eqref{eq:arna} in the proof of Lemma \ref{lem:Gxi}(iv),
\beq\label{eq:NLLC}
\underline C_{j}\le \lim\inf_{n\to\infty} n^{-1} \nu^{(j)}(\bm \Lambda_n^\prime\bm \Lambda_n) \le
\lim\sup_{n\to\infty} n^{-1} \nu^{(j)}(\bm\Lambda_n^\prime\bm \Lambda_n) \le \overline C_j, \qquad j=1,\ldots, r.
\eeq
Thus, we have $\lim_{n\to\infty} n^{-1}\nu^{(r)}(\bm\Lambda_n^\prime\bm\Lambda_n) \ge \underline C_r$, which, once substituted in \eqref{eq:pulce}, proves part (i). 

For part (ii) the proof is the same as part (i), but in \eqref{eq:pulce} we use Lemma \ref{lem:Gxi}(v) instead of Assumption \ref{ass:idio}(a), thus replacing $C_\xi$ with $M_2$.
Parts (iii) and (iv) are obvious by just setting $\bm P^{-1}=\mbf 0_{r\times r}$ in the first step of \eqref{eq:pulce}. 

For part (v) we have
\al{
n^{-1}\Vert \bm\Lambda_n^\prime(\bm\Sigma_n^\xi)^{-1}\bm\Lambda_n\Vert&=n^{-1}\nu^{(1)}(\bm\Lambda_n^\prime(\bm\Sigma_n^\xi)^{-1}\bm\Lambda_n)=n^{-1}\nu^{(1)}(\bm\Lambda_n\bm\Lambda_n^\prime(\bm\Sigma_n^\xi)^{-1})\nn\\
&\le n^{-1}\nu^{(1)}(\bm\Lambda_n\bm\Lambda_n^\prime) \nu^{(1)}((\bm\Sigma_n^\xi)^{-1}) =\frac{\nu^{(1)}(\bm\Lambda_n^\prime\bm\Lambda_n)}{n\nu^{(n)}(\bm\Sigma_n^\xi)}\nn\\
&\le \frac{\nu^{(1)}(\bm\Lambda_n^\prime\bm\Lambda_n)}{n C_\xi^{-1}},\label{eq:mostri}
}
by Assumption \ref{ass:idio}(a). Then, by \eqref{eq:NLLC} in the proof of Lemma \ref{lem:LSL}, we have $\lim_{n\to\infty} n^{-1}\nu^{(1)}(\bm\Lambda_n^\prime\bm\Lambda_n) \le \overline C_1$, which, once substituted in \eqref{eq:mostri}, proves part (v).  

For part (vi) the proof is the same as part (v), but in \eqref{eq:mostri} we use Assumption \ref{ass:idio}(f) instead of Assumption \ref{ass:idio}(a), thus replacing $C_\xi^{-1}$ with $L_\xi$. 

For part (vii)
\al{
\Vert \bm\Lambda_n(\bm\Sigma_n^\xi)^{-1}\Vert \le \Vert \bm\Lambda_n\Vert \, \Vert (\bm\Sigma_n^\xi)^{-1}\Vert\le  \Vert \bm\Lambda_n\Vert \,\max_{i=1,\ldots, n}(\sigma_i^2)^{-1} = O(\sqrt n),\nn
}
by Lemma \ref{lem:lambdasqrtn} and Assumption \ref{ass:idio}(a). And for part (viii) the proof is the same as for part (vii) but using Assumption \ref{ass:idio}(f) instead of Assumption \ref{ass:idio}(a). This proves parts (vii) and (viii) and completes the proof. $\Box$

\begin{lem}\label{lem:taylorinv}
Given two invertible matrices $\bm K$ and $\bm H$ the following holds: 
$$
(\bm H+ \bm K)^{-1}=\bm K^{-1}- (\bm H + \bm K)^{-1}\bm H\bm K^{-1}.
$$
\end{lem}

\noindent\textsc{Proof.} 
We have
\begin{align}
(\bm H+ \bm K)^{-1}&= (\bm H + \bm K)^{-1}- \bm K^{-1} + \bm K^{-1}= (\bm H+ \bm K)^{-1}(\bm K - (\bm H + \bm K))\bm K^{-1}+ \bm K^{-1}\nn\\
&= (\bm H + \bm K)^{-1}(-\bm H)\bm K^{-1} + \bm K^{-1}= \bm K^{-1}- (\bm H + \bm K)^{-1}\bm H\bm K^{-1}.
\end{align}
This completes the proof. $\Box$

\begin{lem}\label{lem:denom}
For $m<n$ with $m$ independent of $n$ and given
\begin{compactenum}
\item[(a)] an $m\times m$ matrix $\bm A$ symmetric and positive definite with $\Vert\bm A\Vert \le M_A$; 
\item[(b)] an $n\times n$ matrix $\bm B$ symmetric and positive definite with $\Vert\bm B\Vert\le M_B$;
\item[(c)] an $n\times m$ matrix $\bm C$ such that $\bm C^\prime\bm C$ is positive definite with 
\beq\nn
\underline M_{C_j}\le \lim\inf_{n\to\infty} n^{-1} \nu^{(j)}(\bm C^\prime \bm C) \le
\lim\sup_{n\to\infty} n^{-1} \nu^{(j)}(\bm C^\prime \bm C) \le \overline M_{C_j}, \qquad j=1,\ldots, m;
\eeq
\end{compactenum}
where, $M_A$, $M_B$, $\underline M_{C_j}$ and $\overline M_{C_j}$ are finite positive reals independent of $n$ and $m$, then the following holds
\[
\Vert (\bm A^{-1}+\bm C^\prime\bm B^{-1}\bm C)^{-1}\bm C^\prime\bm B^{-1}\bm C - \mbf I_m\Vert= O(n^{-1}).
\]
\end{lem}

\noindent\textsc{Proof.} In Lemma \ref{lem:taylorinv} set $\bm K=\bm C^\prime\bm B^{-1}\bm C$ and $\bm H=\bm A^{-1}$, then, 
\begin{align}
(\bm A^{-1}+\bm C^\prime\bm B^{-1}\bm C)^{-1}=(\bm C^\prime\bm B^{-1}\bm C)^{-1}-(\bm A^{-1}+\bm C^\prime\bm B^{-1}\bm C)^{-1}\bm A^{-1}(\bm C^\prime\bm B^{-1}\bm C)^{-1}.\label{eq:abcinv00}
\end{align}
which implies
\beq\label{eq:abcinv}
(\bm A^{-1}+\bm C^\prime\bm B^{-1}\bm C)^{-1}\bm C^\prime\bm B^{-1}\bm C - \mbf I_m = (\bm A^{-1}+\bm C^\prime\bm B^{-1}\bm C)^{-1}\bm A^{-1}.
\eeq
Then, 
by Weyl's inequality \citep[Theorem 1]{MK04}:
\beq
\nu^{(m)}(\bm A^{-1}+ \bm C^\prime \bm B^{-1} \bm C) \ge\nu^{(m)}(\bm A^{-1})+\nu^{(m)}(\bm C^\prime \bm B^{-1} \bm C)=\{\nu^{(1)}(\bm A)\}^{-1}+\nu^{(m)}(\bm C^\prime \bm B^{-1} \bm C)
.\label{eq:weyllunga}
\eeq
From \eqref{eq:weyllunga}, we have
\begin{align}
\Vert (\bm A^{-1}+\bm C^\prime\bm B^{-1}\bm C)^{-1}\bm A^{-1}\Vert&\le\Vert (\bm A^{-1}+\bm C^\prime\bm B^{-1}\bm C)^{-1}\Vert\ \Vert\bm A^{-1}\Vert\nn\\
&=\l\{\nu^{(m)}(\bm A^{-1}+ \bm C^\prime \bm B^{-1} \bm C)\r\}^{-1} \ \l\{\nu^{(m)}(\bm A)\r\}^{-1}
\nn\\
&\le \l\{\nu^{(m)}(\bm C^\prime \bm B^{-1} \bm C)+\{\nu^{(1)}(\bm A)\}^{-1}\r\}^{-1}\ \l\{\nu^{(m)}(\bm A)\r\}^{-1}.
\label{eq:abcinv2}
\end{align}
For first term on the rhs of \eqref{eq:abcinv2}, the $m$ eigenvalues of $\bm C'\bm B^{-1}\bm C$ are also the $m$ largest non-zero eigenvalues of $\bm C\bm C' \bm B^{-1}$, and the $m$ largest non-zero  eigenvalues of $\bm C\bm C'$ are also the $m$ eigenvalues of $\bm C'\bm C$. Therefore, because of \citet[Theorem 7]{MK04}:
\[
n^{-1} \nu^{(m)}(\bm C\bm C^\prime \bm B^{-1} )\ge n^{-1} \nu^{(m)}(\bm C\bm C^\prime) \nu^{(m)}(\bm B^{-1})=n^{-1} \nu^{(m)}(\bm C^\prime\bm C)\l\{{\nu^{(1)}(\bm B)}\r\}^{-1}.
\]
Thus, by conditions (b) and (c)
\[
n\l\{\nu^{(m)}(\bm C^\prime \bm B^{-1} \bm C)\r\}^{-1}\le n \l\{\nu^{(m)}(\bm C^\prime\bm C)\r\}^{-1} \nu^{(1)}(\bm B)\le \frac {M_B}{\underline M_{C_m}}.
\]
Moreover, by condition (a), $\nu^{(1)}(\bm A)>0$ and $\nu^{(1)}(\bm A)\le M_A$, i.e., $\{\nu^{(1)}(\bm A)\}^{-1}\ge M_A$, and 
\[
\nu^{(m)}(\bm C^\prime \bm B^{-1} \bm C)+\{\nu^{(1)}(\bm A)\}^{-1}\ge \nu^{(m)}(\bm C^\prime \bm B^{-1} \bm C)\ge \frac {M_B}{\underline M_{C_m}}.
\] 

For the second term on the rhs of \eqref{eq:abcinv2}, by condition (a),
\[
\l\{\nu^{(m)}(\bm A)\r\}^{-1} \leq \frac 1 L_A,
\] 
for some finite positive real $L_A$. Hence, from \eqref{eq:abcinv2} 
\[
n \Vert (\bm A^{-1}+\bm C^\prime\bm B^{-1}\bm C)^{-1}\bm A^{-1}\Vert\le\frac {M_B}{\underline M_{C_m}L_A}. 
\]

 and by using it in \eqref{eq:abcinv} we complete the proof. $\Box$ 

\begin{lem}\label{lem:LSL} For any $r\times r$ symmetric and positive definite matrix $\bm P$ with $\Vert \bm P\Vert\le M_P$ for some finite positive real $M_P$, under Assumptions \ref{ass:common} and \ref{ass:idio}, as $n\to\infty$,
\begin{compactenum}[(i)]
\item
$n\Vert (\bm\Lambda_n^\prime(\bm\Sigma_n^\xi)^{-1}\bm\Lambda_n+\bm P^{-1})^{-1}\bm\Lambda_n^\prime(\bm\Sigma_n^\xi)^{-1}\bm\Lambda_n-\mbf I_r\Vert=O(1)$;
\item
$n\Vert (\bm\Lambda_n^\prime(\bm\Gamma_n^\xi)^{-1}\bm\Lambda_n+\bm P^{-1})^{-1}\bm\Lambda_n^\prime(\bm\Gamma_n^\xi)^{-1}\bm\Lambda_n-\mbf I_r\Vert=O(1)$;
\item 
$n^2\Vert (\bm\Lambda_n^\prime(\bm\Sigma_n^\xi)^{-1}\bm\Lambda_n+\bm P^{-1})^{-1}-(\bm\Lambda_n^\prime(\bm\Sigma_n^\xi)^{-1}\bm\Lambda_n)^{-1}\Vert=O(1)$;
\item 
$n^2\Vert (\bm\Lambda_n^\prime(\bm\Gamma_n^\xi)^{-1}\bm\Lambda_n+\bm P^{-1})^{-1}-(\bm\Lambda_n^\prime(\bm\Gamma_n^\xi)^{-1}\bm\Lambda_n)^{-1}\Vert=O(1)$.
\end{compactenum}
\end{lem}

\noindent
\textsc{Proof.}
Both results follow from Lemma \ref{lem:denom} since: $\bm P$ satisfies condition (a) by assumption, $\bm \Sigma_n^\xi$ satisfies condition (b) because of Assumption \ref{ass:idio}(a) and $\bm \Gamma_n^\xi$ satisfies condition (b) because of Lemma \ref{lem:Gxi}(v) and Assumption \ref{ass:idio}(f), and ${\bm\Lambda}_{n}^\prime{\bm\Lambda}_{n}$ satisfies condition (c) since it is positive definite because of Assumptions \ref{ass:common}(a), and, moreover, its eigenvalues are such that 
$\underline C_{j}\le \lim\inf_{n\to\infty} n^{-1} \nu^{(j)}(\bm \Lambda_n^\prime\bm \Lambda_n) \le
\lim\sup_{n\to\infty} n^{-1} \nu^{(j)}(\bm\Lambda_n^\prime\bm \Lambda_n) \le \overline C_j$, for $j=1,\ldots, r,$ as shown in \eqref{eq:NLLC} in the proof of Lemma \ref{lem:LSL2}. 

Turning to part (iii),
\al{
n^2\Vert &(\bm\Lambda_n^\prime(\bm\Sigma_n^\xi)^{-1}\bm\Lambda_n+\bm P^{-1})^{-1}-(\bm\Lambda_n^\prime(\bm\Sigma_n^\xi)^{-1}\bm\Lambda_n)^{-1}\Vert\nn\\
&\le n^2\Vert (\bm\Lambda_n^\prime(\bm\Sigma_n^\xi)^{-1}\bm\Lambda_n+\bm P^{-1})^{-1}(\bm\Lambda_n^\prime(\bm\Sigma_n^\xi)^{-1}\bm\Lambda_n)-\mbf I_r\Vert\,
\Vert (\bm\Lambda_n^\prime(\bm\Sigma_n^\xi)^{-1}\bm\Lambda_n)^{-1}\Vert = O(1).\nn
}
by part (i) and Lemma \ref{lem:LSL2}(iii). This proves part (iii).

Part (iv) is proved in the same way but using Lemma \ref{lem:LSL2}(iv) instead of Lemma \ref{lem:LSL2}(iv).
This completes the proof. $\Box$


\begin{lem}\label{lem:LSXi} 
Under Assumptions \ref{ass:common} and \ref{ass:idio}, as $n,T\to\infty$, 
\begin{compactenum}
\item [(i)] $n^{-1/2}\Vert \bm\Lambda_n^\prime(\bm\Sigma_n^\xi)^{-1}\bm\xi_{nt}\Vert=O_p(1)$, uniformly in $t$;
\item [(ii)] $n^{-1/2}\Vert \bm\Lambda_n^\prime(\bm\Gamma_n^\xi)^{-1}\bm\xi_{nt}\Vert=O_p(1)$, uniformly in $t$;
\item [(iii)] $n^{-1/2}T^{-1/2}\Vert\sum_{t=1}^T \bm\Lambda_n^\prime(\bm\Sigma_n^\xi)^{-1}\bm\xi_{nt}\Vert=O_p(1)$;
\item [(iv)] $ n^{-1/2}T^{-1/2} \l\Vert
\bm\Lambda_n^\prime\bm{\mathcal E}_{nT}^\prime\r\Vert_F=O_p(1)$;
\item [(v)] $ n^{-1/2}T^{-1/2} \Vert
\bm\Lambda_n^\prime(\bm\Sigma_n^\xi)^{-1}\bm{\mathcal E}_{nT}^\prime\Vert_F=O_p(1)$;
\item [(vi)] $n^{-1/2}T^{-1/2}\Vert \bm{\mathcal E}_{nT}\Vert_F=O_p(1)$;
\end{compactenum}
where $ \bm{\mathcal E}_{nT}=(\bm\xi_{n1}\cdots \bm\xi_{nT})^\prime$.
\end{lem}

\noindent
\textsc{Proof.}  Throughout, let $\lambda_{ij}$ be the $(i,j)$the entry of $\bm\Lambda_n$.
For part (i), we have
\begin{align}
\E\l[\Vert n^{-1/2}\bm\Lambda_n^\prime(\bm\Sigma_n^\xi)^{-1}\bm\xi_{nt}\Vert^2\r]&
=\sum_{j=1}^r n^{-1} \E\l[\l(\sum_{i=1}^n \frac{\lambda_{ij} \xi_{it}}{\sigma_i^2}\r)^2\r]\nn\\
&\le r\max_{j=1,\ldots,r} n^{-2}\sum_{i=1}^n\sum_{k=1}^n \frac{\vert \lambda_{ij}\vert\, \vert \lambda_{kj}\vert }{\sigma_i^2\sigma_k^2} \E\l[  \xi_{it}\xi_{kt}\r]\nn\\
&\le r n^{-1}M_\lambda^2  C_\xi^2\sum_{i=1}^n\sum_{k=1}^n \vert \E[\xi_{it}\xi_{kt}]\vert\le  rM_\lambda^2  C_\xi^2M_2,\label{eq:gigi}
\end{align}
where in the third step we used Assumption \ref{ass:common}(a) (since   $\max_{j=1,\ldots, r}\vert \lambda_{ij}\vert\le \Vert \bm\lambda_i\Vert\le M_\lambda$, for all $i=1,\ldots, n$), and Assumption \ref{ass:idio}(a), and in the last step we used Lemma \ref{lem:Gxi}(ii). 
By Chebychev's inequality and since the constants in \eqref{eq:gigi} do not depend on $t$, we prove part (i). 

For part (ii), we have
\al{
\E\l[\Vert n^{-1/2}\bm\Lambda_n^\prime(\bm\Gamma_n^\xi)^{-1}\bm\xi_{nt} \Vert^2\r]&=
\E\l[n^{-1}\bm\xi_{nt}^\prime(\bm\Gamma_n^\xi)^{-1}\bm\Lambda_n\bm\Lambda_n^\prime(\bm\Gamma_n^\xi)^{-1}\bm\xi_{nt}\r]\nn\\
&=\E\l[n^{-1} \text{tr}\l\{\bm\Lambda_n^\prime(\bm\Gamma_n^\xi)^{-1}\bm\xi_{nt}\bm\xi_{nt}^\prime(\bm\Gamma_n^\xi)^{-1}\bm\Lambda_n\r\} \r]\nn\\
&=n^{-1}  \text{tr}\l\{\bm\Lambda_n^\prime(\bm\Gamma_n^\xi)^{-1}\E\l[\bm\xi_{nt}\bm\xi_{nt}^\prime\r](\bm\Gamma_n^\xi)^{-1}\bm\Lambda_n\r\}\nn\\
&=n^{-1} \text{tr}\l\{\bm\Lambda_n^\prime(\bm\Gamma_n^\xi)^{-1}\bm\Lambda_n\r\}\nn\\
&\le r n^{-1} \nu^{(1)}\l(\bm\Lambda_n^\prime(\bm\Gamma_n^\xi)^{-1}\bm\Lambda_n\r)\nn\\
&= r n^{-1} \Vert \bm\Lambda_n^\prime(\bm\Gamma_n^\xi)^{-1}\bm\Lambda_n\Vert= O(1),\label{eq:domo}
}
by Lemma \ref{lem:LSL2}(vi). By Chebychev's inequality and since the constants in \eqref{eq:domo} do not depend on $t$, we complete the proof of part (ii).
%
%

For part (iii), we have
\begin{align}
\E\l[\l\Vert n^{-1/2}T^{-1/2}\sum_{t=1}^T \bm\Lambda_n^\prime(\bm\Sigma_n^\xi)^{-1}\bm\xi_{nt}\r\Vert^2\r]&
= n^{-1}T^{-1}\sum_{j=1}^r \E\l[\l(\sum_{t=1}^T\sum_{i=1}^n \frac{\lambda_{ij} \xi_{it}}{\sigma_i^2}\r)^2\r]\nn\\
&\le r\max_{j=1,\ldots,r} n^{-1}T^{-1}\sum_{t=1}^T\sum_{s=1}^T\sum_{i=1}^n\sum_{k=1}^n \frac{\vert \lambda_{ij}\vert\, \vert \lambda_{kj}\vert }{\sigma_i^2\sigma_k^2} \E\l[  \xi_{it}\xi_{ks}\r]\nn\\
&\le r n^{-1}T^{-1}M_\lambda^2  C_\xi^2\sum_{t=1}^T\sum_{s=1}^T\sum_{i=1}^n\sum_{k=1}^n \vert \E[\xi_{it}\xi_{ks}]\vert\le rM_\lambda^2  C_\xi^2M_1,\nn
\end{align}
by Assumptions \ref{ass:common}(a), \ref{ass:idio}(a), and Lemma \ref{lem:Gxi}(i). 
By Chebychev's inequality we prove part (iii).  

For part (iv), we have
\al{
\E\l[\l\Vert n^{-1/2}T^{-1/2}\bm\Lambda_n^\prime\bm{\mathcal E}_{nT}^\prime\r\Vert_F^2
\r]=&\, n^{-1}T^{-1}\sum_{k=1}^r\sum_{t=1}^T\E\l[\l(\sum_{i=1}^n \lambda_{ik} \xi_{it}\r)^2\r]\nn\\
\le&\, n^{-1}rM_\lambda^2 \max_{t=1,\ldots, T} \sum_{i,j=1}^n\vert \E[\xi_{it}\xi_{jt}]\vert\le rM_\lambda^2 M_2,\nn
}
by Assumption \ref{ass:common}(a) and Lemma \ref{lem:Gxi}(ii). By Chebychev's inequality we prove part (iv).

Part (v) is proved as parts (iv), but using also Assumption \ref{ass:idio}(a). 

For part (vi), we have
\al{
\E\l[\l\Vert n^{-1/2}T^{-1/2}  \bm{\mathcal E}_{nT}\r\Vert^2_F\r]=&\, n^{-1}T^{-1}\sum_{t=1}^T \sum_{i=1}^n \E[\xi_{it}^2]\nn\\
\le& \max_{t=1,\ldots,T}\max_{i=1,\ldots,n}\E[\xi_{it}^2] = \max_{i=1,\ldots,n}\sigma_i^2 \le C_\xi,\nn
}
by Assumption \ref{ass:idio}(a). By Chebychev's inequality we prove part (vi).
This completes the proof.
$\Box$

\begin{lem}\label{lem:COVFF0}
Under Assumptions \ref{ass:common}, \ref{ass:idio}, \ref{ass:ind}, and \ref{ass:ident}, as $n,T\to\infty$: 
\begin{compactenum}
\item [(i)] $\sqrt{nT} \Vert n^{-1}T^{-1} \sum_{t=1}^T\bm\Lambda_n^\prime \bm\xi_{nt}{\mbf F}_{t}^\prime\Vert_F=O_p(1)$;
\item [(ii)] $\sqrt{nT}\, \Vert n^{-3/2}T^{-1} \sum_{t=1}^T\bm\Lambda_n^\prime \bm\xi_{nt}\bm{\xi}_{nt}^\prime\Vert_F=O_p(1)$;
\item [(iii)] $\sqrt{nT}\, \Vert n^{-2}T^{-1} \sum_{t=1}^T\bm\Lambda_n^\prime \bm\xi_{nt}\bm{\xi}_{nt}^\prime\bm\Lambda_n\Vert_F=O_p(1)$;
\item [(iv)] $\sqrt{nT} \Vert n^{-1}T^{-1} \sum_{t=1}^T\bm\Lambda_n^\prime(\bm\Sigma_n^\xi)^{-1} \bm\xi_{nt}{\mbf F}_{t}^\prime\Vert_F=O_p(1)$;
\item [(v)] $\sqrt{nT} \Vert n^{-3/2}T^{-1} \sum_{t=1}^T\bm\Lambda_n^\prime(\bm\Sigma_n^\xi)^{-1} \bm\xi_{nt}\bm{\xi}_{nt}^\prime\Vert_F=O_p(1)$;
\item [(vi)] $\sqrt{nT} \Vert n^{-1}T^{-3/2} \sum_{t=1}^T \mbf F_t \bm\xi_{nt}^\prime(\bm\Sigma_n^\xi)^{-1} \bm{\mathcal E}_{nT}\Vert_F=O_p(1);$
\end{compactenum}
where $ \bm{\mathcal E}_{nT}=(\bm\xi_{n1}\cdots \bm\xi_{nT})^\prime$.
\end{lem}

\noindent
\textsc{Proof.} 
 Throughout, let $\lambda_{ij}$ be the $(i,j)$the entry of $\bm\Lambda_n$.
For part (i),
\al{
\E\l[\l\Vert n^{-1}T^{-1} \sum_{t=1}^T\bm\Lambda_n^\prime \bm\xi_{nt}{\mbf F}_{t}^\prime\r\Vert^2_F \r]=&\, \,n^{-2}T^{-2} \sum_{t,s=1}^T\sum_{i,j=1}^n\sum_{k,h=1}^r 
\lambda_{ih}\lambda_{jh} \E[\xi_{it}\xi_{js} F_{kt}F_{ks}]\nn\\
\le&\, n^{-2}T^{-2} M_\lambda^2 r^2\max_{k,h=1,\ldots, r} \sum_{t,s=1}^T\sum_{i,j=1}^n\vert\E[\xi_{it}\xi_{js}]\vert\, \vert  \E[F_{kt}F_{ks}]\vert\nn\\
\le&\, n^{-1}T^{-1} M_\lambda^2 r^2M_1,\label{eq:pastiera}
}
by Lemma \ref{lem:Gxi}(i) and because $\vert  \E[F_{kt}F_{ks}]\vert\le 1$ by Cauchy-Schwarz inequality and Assumption \ref{ass:ident}(b). By Chebychev's inequality we prove part (i).

For part (ii),
\al{
\E&\l[\l\Vert n^{-3/2}T^{-1} \sum_{t=1}^T\bm\Lambda_n^\prime \bm\xi_{nt}\bm\xi_{nt}^\prime \r\Vert^2_F \r]
=  \sum_{k=1}^r \sum_{j=1}^n\E\l[
\l\vert
n^{-3/2}T^{-1} \sum_{t=1}^T\sum_{i=1}^n\lambda_{ik} \xi_{it}\xi_{jt}
\r\vert^2
\r]\nn\\
\le&\, rM_\lambda n^{-1}T^{-1} \max_{j=1,\ldots, n}\E\l[
\l\vert
n^{-1/2}T^{-1/2} \sum_{t=1}^T\sum_{i=1}^n  \xi_{it}\xi_{jt}
\r\vert^2
\r]\nn\\
\le&\, rM_\lambda n^{-2}T^{-2} \max_{j=1,\ldots, n}\sum_{t,s=1}^T\sum_{i,\ell=1}^n
\vert\E[
\xi_{it}\xi_{jt}\xi_{\ell s}\xi_{js}
]\vert \nn\\
\le&\, rM_\lambda n^{-1}T^{-1} \mathrm K_\xi,\nn
}
by Assumptions \ref{ass:common}(a) and \ref{ass:idio}(d). By Chebychev's inequality we prove part (ii).

For part (iii), 
following the proof of part (ii),
\al{
\E\l[\l\Vert n^{-2}T^{-1} \sum_{t=1}^T\bm\Lambda_n^\prime \bm\xi_{nt}\bm\xi_{nt}^\prime\bm\Lambda_n \r\Vert^2_F \r]=&\,  \sum_{k,h=1}^r \E\l[
\l\vert
n^{-2}T^{-1} \sum_{t=1}^T\sum_{i,j=1}^n\lambda_{ik}\lambda_{jh} \xi_{it}\xi_{jt}\r\vert^2
\r]\nn\\
\le&\, rM_\lambda n^{-4}T^{-2} \sum_{t,s=1}^T\sum_{i_1,j_1=1}^n\sum_{i_2,j_2=1}^n
\vert \E[\xi_{i_1t}\xi_{j_1t}\xi_{i_2 s}\xi_{j_2s}]\vert
\nn\\
\le&\, rM_\lambda n^{-1}T^{-1} \mathrm K_\xi,\nn
}
by Assumptions \ref{ass:common}(a) and \ref{ass:idio}(d). By Chebychev's inequality we prove part (iii).

Parts (iv) and (v) are proved as parts (i) and (ii), respectively, but using also Assumption \ref{ass:idio}(a). 


For part (vi)
\al{
\E&\l[\l\Vert n^{-1}T^{-3/2} \sum_{t=1}^T \mbf F_t \bm\xi_{nt}^\prime(\bm\Sigma_n^\xi)^{-1} \bm{\mathcal E}_{nT}\r\Vert_F^2\r]
=n^{-2}T^{-3}\sum_{k=1}^r\sum_{s=1}^T \E\l[
\l(\sum_{t=1}^T\sum_{i=1}^n F_{kt} \xi_{it} (\sigma_i^2)^{-1}\xi_{is} \r)^2\r]\nn\\
\le&\,  n^{-2}T^{-3}C_\xi^2 \sum_{k=1}^r\sum_{s=1}^T \sum_{t_1,t_2=1}^T \sum_{i,j=1}^n \E[F_{kt_1}F_{kt_2} \xi_{it_1}\xi_{jt_2}\xi_{is}\xi_{js}]\nn\\
\le&\,  n^{-2}T^{-2} r C_\xi^2 \max_{k=1,\ldots, r} \max_{t_1,t_2=1,\ldots, T}\vert \E[F_{kt_1}F_{kt_2}]\vert  \max_{s=1,\ldots, T} \sum_{t_1,t_2=1}^T \sum_{i,j=1}^n\vert \E[ \xi_{it_1}\xi_{jt_2}\xi_{is}\xi_{js}]\vert\nn\\
\le&\,  n^{-1}T^{-1} r C_\xi^2\max_{k=1,\ldots, r}  \max_{t=1,\ldots,T} \E[F_{kt}^2]\, \mathrm K_\xi\nn\\
\le&\,  n^{-1}T^{-1} r C_\xi^2 \mathrm K_\xi,\nn
}
by Assumptions \ref{ass:idio}(a) and \ref{ass:idio}(d), Lemma \ref{lem:fidio}, and Cauchy-Schwarz inequality jointly with Assumption \ref{ass:ident}(b). By Chebychev's inequality we prove part (vi).
This completes the proof. $\Box$

\begin{lem}\label{lem:lambda}
Under Assumptions \ref{ass:common} and \ref{ass:ident}, for all $n>N_0$: 
$\Vert n^{-1/2} \bm\Lambda_n -n^{-1/2} \mbf V_n^\chi\bm{\mathcal S}(\mbf M_n^{\chi})^{1/2}\Vert = 0$,
for some $r\times r$ positive diagonal matrix $\bm{\mathcal S}$ independent of $n$ with entries $\mathbb I([\mbf V_n^\chi]_{1j}\ge 0)-\mathbb I([\mbf V_{n}^\chi]_{1j}<0)$, $j=1,\ldots,r$, and where $N_0$ is defined in Assumption \ref{ass:common}(a).
\end{lem}
\noindent\textsc{Proof.}  By Assumption \ref{ass:ident}(b), for all $n\in\mathbb N$, 
\beq\label{eq:ritardo}
n^{-1}\bm\Gamma_n^\chi =n^{-1} \bm\Lambda_n\bm\Lambda_n^\prime =  n^{-1}\mbf V_n^\chi\mbf M_n^\chi \mbf V_n^{\chi\prime}.
\eeq
First notice that, the $r$ non-zero eigenvalues of $n^{-1}\bm\Gamma_n^\chi$ are the $r$ eigenvalues of $n^{-1}\bm\Lambda_n^\prime\bm\Lambda_n$ and for all $n>N_0$, $\Vert\bm\Lambda_n^\prime\bm\Lambda_n -\bm\Sigma_\Lambda\Vert =0$, by Assumption \ref{ass:common}(a). Moreover,  $\bm\Sigma_\Lambda$ is diagonal and positive definite by Assumptions \ref{ass:ident}(b) and \ref{ass:common}(a), respectively. Hence, for all $n>N_0$,
\beq
n^{-1}\mbf M_n^\chi =\bm\Sigma_\Lambda\;\text{ and }\; 
n(\mbf M_n^\chi)^{-1} =\bm\Sigma_\Lambda^{-1}.\label{eq:paris}
\eeq

Furthermore, it must be that the columns of $\bm\Lambda_n$ span the same space as the columns of $\mbf V_n^\chi$. Since the eigenvectors are normalized and for all $n>N_0$, $\Vert n^{-1} \mbf M_n^\chi -\bm\Sigma_\Lambda\Vert=0$ by \eqref{eq:paris}, there exist two $r\times r$ matrices $\mbf K_{1n}$ and $\mbf K_{2n}$ such that, for all $n>N_0$,
\beq\label{eq:ritardo3}
n^{-1/2}\bm\Lambda_n =
\mbf V_n^\chi\mbf K_{1n}\; \text{ and } 
\mbf V_n^\chi= 
n^{-1/2} \bm\Lambda_n \mbf K_{2n}.
\eeq
Let 
$$
\mbf K_{1}= \lim_{n\to\infty} n^{-1/2} (\mbf V_n^{\chi\prime}\mbf V_n^{\chi})^{-1}\mbf V_n^{\chi\prime}\bm\Lambda_n= \lim_{n\to\infty} n^{-1/2} \mbf V_n^{\chi\prime}\bm\Lambda_n
$$ 
and 
$$
\mbf K_2= \lim_{n\to\infty}  n(\bm\Lambda_n^\prime\bm\Lambda_n)^{-1} n^{-1/2} \bm\Lambda_n^\prime \mbf V_n^\chi = \lim_{n\to\infty} \sqrt n(\bm\Lambda_n^\prime\bm\Lambda_n)^{-1} \bm\Lambda_n^\prime \mbf V_n^\chi.
$$
Then, by linear projection from \eqref{eq:ritardo3} we have that 
\beq\label{eq:ritardo1}
\lim_{n\to\infty} \mbf K_{1n}  =\mbf K_1,
\eeq
which is positive definite since the columns of $\mbf V_n^\chi$ are linear combinations of the columns of $\bm\Lambda_n$ so, for all $n>N_0$, $\text{rk} (n^{-1/2}\mbf V_n^{\chi\prime}\bm\Lambda_n)=\text{rk} (n^{-1}\bm\Lambda_n^\prime \bm\Lambda_n)=\text{rk}(\bm\Sigma_\Lambda)=r$ by Assumption \ref{ass:common}(a).   Moreover, for all $n>N_0$, $\Vert \mbf K_1\Vert\le n^{-1/2}\Vert \bm\Lambda_n\Vert$, which is finite by Lemma \ref{lem:lambdasqrtn}.

Similarly, from \eqref{eq:ritardo3} we also have that 
\beq\label{eq:ritardo2}
\lim_{n\to\infty}  \mbf K_{2n}  = \mbf K_2.
\eeq
which exists and is positive definite, since $\Vert n(\bm\Lambda_n^\prime\bm\Lambda_n)^{-1}-n\bm\Sigma_\Lambda^{-1}\Vert = 0$, for all $n>N_0$, and $\bm\Sigma_\Lambda$ is finite and positive definite by Assumption \ref{ass:common}(a).  Moreover, for all $n>N_0$, $\Vert \mbf K_2\Vert \le n\Vert(\bm\Lambda_n^\prime \bm\Lambda_n)^{-1}\Vert\,n^{-1/2}\Vert \bm\Lambda_n\Vert= \Vert \bm\Sigma_\Lambda^{-1}\Vert n^{-1/2}\Vert \bm\Lambda_n\Vert$, which is finite by Assumption \ref{ass:common}(a) and Lemma \ref{lem:lambdasqrtn}.

By using \eqref{eq:ritardo3} into the rhs of \eqref{eq:ritardo}, we get, for all $n>N_0$,
$$
n^{-2} \bm\Lambda_n\mbf K_{2n} \mbf M_n^\chi \mbf K_{2n}^\prime \bm\Lambda_n^\prime=
\mbf V_n^\chi\mbf M_n^\chi \mbf V_n^{\chi\prime},
$$
which, since eigenvectors are normalized, implies, that, for all $n>N_0$, we can write 
\beq\label{eq:ritardo4}
n^{-2} \mbf V_n^{\chi\prime}\bm\Lambda_n\mbf K_{2n} \mbf M_n^\chi \mbf K_{2n}^\prime \bm\Lambda_n^\prime\mbf V_n^{\chi}=
n^{-1}\mbf M_n^\chi.
\eeq
From \eqref{eq:ritardo4} we must have 
$\mbf I_r=\lim_{n\to\infty} n^{-1/2} \mbf V_n^{\chi\prime}\bm\Lambda_n\mbf K_{2n} = \mbf K_1\mbf K_2$,
so, as expected $\mbf K_1=\mbf K_2^{-1}$ and $\mbf K_2=\mbf K_1^{-1}$. It follows that, for all $n>N_0$,
\beq\label{eq:tardi}
\mbf V_n^{\chi\prime}\bm\Lambda_n(\bm\Lambda_n^\prime\bm\Lambda_n)^{-1} \bm\Lambda_n^\prime \mbf V_n^\chi=\mbf I_r.
\eeq

Now, from \eqref{eq:ritardo}, we can also write that for all $n>N_0$
\beq\label{eq:b18}
n^{-1/2}\bm\Lambda_n\bm R_n =
n^{-1/2} \mbf V_n^\chi (\mbf M_n^\chi)^{1/2},
\eeq
for some $r\times r$ matrix $\bm R_n$. Let, 
\beq
\bm R = \lim_{n\to\infty} (\bm\Lambda_n^\prime \bm\Lambda_n)^{-1}\bm\Lambda_n^\prime\mbf V_n^\chi(\mbf M_n^\chi)^{1/2}=  \lim_{n\to\infty}  n^{-1/2}\mbf K_{2n} (\mbf M_n^\chi)^{1/2}= \mbf K_{2} \lim_{n\to\infty}  n^{-1/2} (\mbf M_n^\chi)^{1/2}.\label{eq:vrum}
\eeq
Hence, for all $n>N_0$, $\bm R = \mbf K_{2} (\bm\Sigma_\Lambda)^{1/2}$ by \eqref{eq:paris}. So $\bm R$ is finite by Assumption \ref{ass:common}(a) and since $\mbf K_2$ is finite.

Moreover, from \eqref{eq:vrum}
\al{
\bm R^{-1} &=\l\{\lim_{n\to\infty} \sqrt n(\mbf M_n^\chi)^{-1/2}\r\}\mbf K_2^{-1} 
=\l\{\lim_{n\to\infty} \sqrt n(\mbf M_n^\chi)^{-1/2}\r\}\mbf K_1\nn\\
& =\lim_{n\to\infty} \sqrt n(\mbf M_n^\chi)^{-1/2} \mbf K_{1n}= \lim_{n\to\infty}(\mbf M_n^\chi)^{-1/2}\mbf V_n^{\chi\prime}\bm\Lambda_n.\label{eq:vrum2}
}
From \eqref{eq:vrum} and \eqref{eq:vrum2}, and using Assumption \ref{ass:common}(a), we also have
\al{
\bm R^{-1} &= \lim_{n\to\infty}(\mbf M_n^\chi)^{-1/2}\mbf V_n^{\chi\prime}\bm\Lambda_n=
 \lim_{n\to\infty}(\mbf M_n^\chi)^{-1/2}(\mbf M_n^\chi)^{-1/2}(\mbf M_n^\chi)^{1/2} \mbf V_n^{\chi\prime}\bm\Lambda_n ( \bm\Lambda_n^\prime\bm\Lambda_n)^{-1}
 \bm\Lambda_n^\prime\bm\Lambda_n\nn\\
 &= \l\{ \lim_{n\to\infty}n(\mbf M_n^\chi)^{-1}\r\} \bm R^\prime  \bm\Sigma_\Lambda.\label{eq:vrum3}
}
Hence, for all $n>N_0$, by \eqref{eq:paris},
$\bm R^{-1}=(\bm\Sigma_\Lambda)^{-1}\bm R^\prime  \bm\Sigma_\Lambda$. Thus, $\bm R^{-1}$ is finite by Assumption \ref{ass:common}(a) and since $\bm R_2$ is finite.
It follows that $\bm R$ is positive definite.

Moreover, $\bm R$ is orthogonal, indeed, from \eqref{eq:ritardo} and \eqref{eq:vrum}
\al{
\bm R\bm R^\prime&=\lim_{n\to\infty} (\bm\Lambda_n^\prime \bm\Lambda_n)^{-1}\bm\Lambda_n^\prime\mbf V_n^\chi(\mbf M_n^\chi)^{1/2}(\mbf M_n^\chi)^{1/2}\mbf V_n^{\chi\prime}\bm\Lambda_n(\bm\Lambda_n^\prime \bm\Lambda_n)^{-1}\nn\\
&= \lim_{n\to\infty} (\bm\Lambda_n^\prime \bm\Lambda_n)^{-1}\bm\Lambda_n^\prime\bm\Gamma_n^\chi\bm\Lambda_n(\bm\Lambda_n^\prime \bm\Lambda_n)^{-1}\nn\\
&= \lim_{n\to\infty} (\bm\Lambda_n^\prime \bm\Lambda_n)^{-1}\bm\Lambda_n^\prime\bm\Lambda_n\bm\Lambda_n^\prime\bm\Lambda_n(\bm\Lambda_n^\prime \bm\Lambda_n)^{-1}=\mbf I_r.\label{eq:vrum4}
}
By substituting \eqref{eq:vrum4} into \eqref{eq:vrum3}, for all $n>N_0$, by \eqref{eq:paris},
\al{
\bm R^{-1}&= 
\bm\Sigma_\Lambda^{-1}\bm R^{-1}  \bm\Sigma_\Lambda.\label{eq:vrum5}
}
 By right-multiplying \eqref{eq:vrum5} by $\bm R$ and left-multiplying by $\bm\Sigma_\Lambda$ we have
\[
\bm\Sigma_\Lambda=\bm R^{-1}  \bm\Sigma_\Lambda \bm R,
\]
which implies $\bm R=\bm{\mathcal J}$ where $\bm{\mathcal J}$ is an $r\times r$ diagonal matrix with entries $\pm 1$ independent of $n$. Therefore, from \eqref{eq:b18}, for all $n>N_0$,
\[
\bm\Lambda_n
\bm{\mathcal J} =
n^{-1/2} \mbf V_n^\chi (\mbf M_n^\chi)^{1/2}
\]
or, equivalently, for all $n>N_0$,
\[
n^{-1/2}\bm\Lambda_n =
 n^{-1/2} \mbf V_n^\chi \bm{\mathcal J} (\mbf M_n^\chi)^{1/2}.
\]
Finally, by Assumption \ref{ass:ident}(c) it must be that $\bm{\mathcal J}=\bm{\mathcal S}$. This completes the proof. $\Box$

%

\begin{lem}\label{lem:xunif}
Under Assumptions \ref{ass:common}, \ref{ass:idio}, and \ref{ass:ind}, 
as $n\to\infty$,  $n^{-1/2}\Vert\mbf x_{nt}\Vert=O_p(1)$, uniformly in $t$.
\end{lem}

\noindent
\noindent\textsc{Proof.} 
We have,
\begin{align}
\E\l[\Vert n^{-1/2}\mbf x_{nt}\Vert^2\r] &= n^{-1} \sum_{i=1}^n \E[x_{it}^2] =n^{-1} \sum_{i=1}^n \l\{\bm\lambda_i^\prime\bm\Gamma^F\bm\lambda_i+ \sigma_i^2 \r\}
\nn\\
&\le \max_{i=1,\ldots, n}
\bm\lambda_i^\prime\bm\Gamma^F\bm\lambda_i +\max_{i=1,\ldots, n} \sigma_i^2 \le M_\lambda^2 \Vert\bm\Gamma^F\Vert \le M_\lambda^2M_F + C_\xi,\label{eq:disaster}
\end{align}
by Assumption \ref{ass:ind}, Assumptions \ref{ass:common}(a) and \ref{ass:common}(b), and Assumption \ref{ass:idio}(a). The proof of part (i) follows by Chebychev's inequality and by noticing that the bound in \eqref{eq:disaster} does not depend on $t$.
This completes the proof. $\Box$

\begin{lem} \label{lem:fidio} Under Assumptions \ref{ass:common} and \ref{ass:ind}, the processes
$\{\xi_{it},\, i\in\mathbb N,\, t\in\mathbb Z\}$ and $\{F_{jt},\, j=1,\ldots, r,\, t\in\mathbb Z\}$ are mutually independent.
\end{lem}

\noindent\textsc{Proof.} It is enough to notice that $\mbf F_t=\sum_{k=0}^\infty \mbf A^k\mbf H\mbf u_{t-k}$, then, by Assumption \ref{ass:ind}, we complete the proof. $\Box$

\begin{lem}\label{lem:consistCOV}
Under Assumptions \ref{ass:common}, \ref{ass:idio}, and \ref{ass:ind}, as $n,T\to\infty$,
\begin{compactenum}
\item [(i)] for $k=0,1$, $T \E[\Vert T^{-1}\sum_{t=1}^T\mbf F_t\mbf F_{t-k}^\prime-\bm\Gamma^F_k\Vert_F^2]= O(1)$, with $\bm\Gamma^F_k=\E[\mbf F_t\mbf F_{t-k}^\prime];$
\item [(ii)] $T \max_{i=1,\ldots, n}\E[\Vert T^{-1}\sum_{t=1}^T\mbf F_t\xi_{it}\Vert^2]= O(1)$;
\item [(iii)] $T \E[\Vert n^{-1/2} T^{-1}\sum_{t=1}^T\mbf F_t\bm\xi_{nt}\Vert_F^2]= O(1)$;
\item [(iv)] 
$T\max_{i,j=1,\ldots, n}\E[\vert T^{-1}\sum_{t=1}^T \xi_{it}\xi_{jt}-\E[\xi_{it}\xi_{jt}]\vert^2]= O(1)$;
\item [(v)] $T\max_{i,j=1,\ldots, n}\E[\vert T^{-1}\sum_{t=1}^T x_{it}x_{jt}-\E[x_{it}x_{jt}]\vert^2]= O(1)$;
\item [(vi)] $T \E[\Vert  n^{-1}T^{-1}\sum_{t=1}^T\bm\xi_{nt} \bm \xi_{nt}^\prime- n^{-1}\bm\Gamma^\xi_n\Vert_F^2]=O(1)$;
\item [(vii)] $T \E[\Vert  n^{-1}T^{-1}\sum_{t=1}^T\mbf x_{nt} \mbf x_{nt}^\prime- n^{-1}\bm\Gamma^x_n\Vert_F^2]=O(1)$.
\end{compactenum}
\end{lem}

\noindent\textsc{Proof.} 
Part (i) follows since $\{\mbf F_t\}$ is ergodic, because of Assumption \ref{ass:common}(d) which implies that $\{\mbf F_t\}$ has summable autocovariances, and therefore $\{\mbf F_t
\mbf F_{t-k}^\prime\}$ is also ergodic (\citealp[Theorem 3.35]{white01}, and \citealp[pp. 170, 182]{stout1974almost}). In particular, $\bm\Gamma_k^F$ is finite because of Assumptions \ref{ass:common}(b),  \ref{ass:common}(d) and \ref{ass:common}(e), and $\E[\Vert\mbf F_t\Vert^4\Vert]$ is also finite because of Assumptions \ref{ass:common}(d)-\ref{ass:common}(g). See also \citet[Proposition 11.1, pp. 298-299]{Hamilton}, which can be applied using the fact that $\{\mbf v_t\}$ is an independent process by Assumption \ref{ass:common}(f) and thus it is a martingale difference process. This proves part (i).

For part (ii), as $T\to\infty$, 
\begin{align}
\E_{}\l[\l\Vert T^{-1}\sum_{t=1}^T\mbf F_t\xi_{it}\r\Vert^2 
\r]&=T^{-2} \sum_{j=1}^r \E_{}\l[\l(\sum_{t=1}^T F_{jt} \xi_{it}\r)^2\r]=T^{-2} \sum_{j=1}^r\sum_{t,s=1}^T \E_{}\l[ F_{jt} \xi_{it}F_{js} \xi_{is}\r]\label{eq:numlambda}\\
&\le T^{-2} \sum_{j=1}^r\sum_{t,s=1}^T \vert \E_{}\l[ F_{jt} F_{js}\r]\vert \,  \vert\E_{} \l[\xi_{it}\xi_{is}\r]\vert\nn\\
&\le T^{-2} \sum_{j=1}^r\sum_{t,s=1}^T \E[F_{jt}^2] \vert\E_{} \l[\xi_{it}\xi_{is}\r]\vert \le T^{-1} r M_F M_3,\nn
\end{align}
because of Lemma \ref{lem:Gxi}(iii), and where we also used Lemma \ref{lem:fidio}, Cauchy-Schwarz inequality, and the fact that $F_{jt}$ is weakly stationary by Assumptions \ref{ass:common}(b), \ref{ass:common}(d), and \ref{ass:common}(e). By noticing that the constants on the rhs of \eqref{eq:numlambda} do not depend on $i$ we prove part (ii). Part (iii) follows directly from part (ii), indeed
\al{
\E_{}\l[\l\Vert n^{-1/2}T^{-1}\sum_{t=1}^T\mbf F_t\bm\xi_{nt}\r\Vert_F^2\r]&= n^{-1}T^{-2} \sum_{j=1}^r\sum_{i=1}^n \E_{}\l[\l(\sum_{t=1}^T F_{jt} \xi_{it}\r)^2\r]\nn\\
&\le T^{-2} \sum_{j=1}^r\max_{i=1,\ldots,n}\E_{}\l[\l(\sum_{t=1}^T F_{jt} \xi_{it}\r)^2\r]\le T^{-1} r M_F M_3.\nn
}

For part (iv), notice that, since $\{\xi_{it}\}$ is a strongly mixing process with exponentially decaying coefficients, because of Assumption \ref{ass:idio}(c), then it is also ergodic 
(\citealp[Proposition 3.44]{white01}, and \citealp{rosenblatt1972}), so that also $\{\xi_{it}\xi_{jt}\}$ is ergodic (\citealp[Theorem 3.35]{white01}, and \citealp[pp. 170, 182]{stout1974almost}). 
In particular, $\E[\xi_{it}\xi_{jt}]$ is finite by Assumption \ref{ass:idio}(a) and $\E[\vert \xi_{it}\xi_{jt}\xi_{is}\xi_{js}\vert]$ is also finite because by Assumption \ref{ass:idio}(d), and both are bounded by constant independent of $i$ and $j$.

This proves part (iv). Part (v) follows directly from parts (i), (ii), and (iv). Part (vi) follows directly from part (iv), and part (vii) follows directly from part (v).
This completes the proof. $\Box$

\begin{lem}\label{lem:frida}
Under Assumptions \ref{ass:common}, \ref{ass:idio}, and \ref{ass:ind}, as $n,T\to\infty$, $\Vert (T^{-1}\sum_{t=1}^T\mbf F_t\mbf F_t^\prime)^{-1}\Vert = O_p(1)$.
\end{lem}

\noindent
\textsc{Proof.}
From  Lemma \ref{lem:consistCOV}(i), and \citet[Theorem 1]{MK04} which is Weyl's inequality
\al{
\l\vert \nu^{(r)}\l(T^{-1}\sum_{t=1}^T\mbf F_t\mbf F_t^\prime\r) - \nu^{(r)}(\bm\Gamma^F) \r\vert\le \l\Vert T^{-1}\sum_{t=1}^T\mbf F_t\mbf F_t^\prime-\bm\Gamma^F\r\Vert = O_p(T^{-1/2}).\nn
}
This implies (note that $x-y\ge -|x-y|$ for any $x,y\in\mathbb R$)
\al{
\det\l(T^{-1}\sum_{t=1}^T\mbf F_t\mbf F_t^\prime \r) =&\,\prod_{j=1}^r \nu^{(j)}\l(T^{-1}\sum_{t=1}^T\mbf F_t\mbf F_t^\prime\r)\ge \l\{\nu^{(r)}\l(T^{-1}\sum_{t=1}^T\mbf F_t\mbf F_t^\prime\r)\r\}^r\nn\\
\ge&\, \l\{ \nu^{(r)}(\bm\Gamma^F) -\l\vert \nu^{(r)}\l(T^{-1}\sum_{t=1}^T\mbf F_t\mbf F_t^\prime\r) - \nu^{(r)}(\bm\Gamma^F) \r\vert \r\}^r>0,\nn
}
by Assumption \ref{ass:common}(b). Thus, $\Vert (T^{-1}\sum_{t=1}^T\mbf F_t\mbf F_t^\prime)^{-1}\Vert = O_p(1)$. This completes the proof. $\Box$

%

\section{Lemmas necessary for proving Proposition \ref{prop:plainvanilla}}

\begin{lem}\label{lem:est0LOAD}
Consider the initial estimator of the loadings $\wh{\bm\Lambda}_n^{(0)}=(\wh{\bm\lambda}_1^{(0)}\cdots \wh{\bm\lambda}_n^{(0)})^\prime$ defined in Section \ref{app:prest}, then, under Assumptions \ref{ass:common}, \ref{ass:idio}, \ref{ass:ind}, and \ref{ass:ident}, as $n,T\to\infty$: 
\begin{compactenum}
\item [(i)] $\min(n,\sqrt T)\,\Vert\wh{\bm\lambda}_{i}^{(0)}-\bm\lambda_i\Vert = O_p(1)$, uniformly in  $i$; 
\item [(ii)] $\min(n,\sqrt T)\,n^{-1/2}\Vert\wh{\bm\Lambda}_n^{(0)}-\bm\Lambda_n\Vert = O_p(1)$.
\end{compactenum}
\end{lem}

\noindent
\textsc{Proof.} Both results are direct consequences of \citet[Theorem 1]{MBPCAQML}, see also  \citet[Theorem 2]{Bai03} under similar assumptions. This completes the proof. $\Box$

\begin{lem}\label{lem:COVFF}
Consider the initial estimator of the factors $\wt{\mbf F}_t$ defined in Section \ref{app:prest}, then, under Assumptions \ref{ass:common}, \ref{ass:idio}, \ref{ass:ind}, and \ref{ass:ident}, as $n,T\to\infty$: 
\begin{compactenum}
\item [(i)] for $k=0,1$, $\min(n^{-1},T^{-1/2})\,\Vert T^{-1} \sum_{t=k+1}^T (\wt{\mbf F}_t-\mbf F_t){\mbf F}_{t-k}^\prime\Vert= O_p(1)$;
\item [(ii)] $\min(n^{-1},T^{-1/2})\,\Vert n^{-1/2} T^{-1} \sum_{t=1}^T (\wt{\mbf F}_t-\mbf F_t){\bm\xi}_{nt}^\prime\Vert= O_p(1)$;
\item [(iii)] $\min(n^{-1},T^{-1/2})\,\Vert n^{-1} T^{-1} \sum_{t=1}^T (\wt{\mbf F}_t-\mbf F_t){\bm\xi}_{nt}^\prime\bm\Lambda_n\Vert= O_p(1)$;
\item [(iv)] $\min(n^{-1},T^{-1/2})\,\Vert T^{-1} \sum_{t=1}^T (\wt{\mbf F}_t-\mbf F_t){\xi}_{it} \Vert= O_p(1)$, uniformly in $i$;
\item [(v)] $\min(n^{-1},T^{-1/2})\,\Vert T^{-1} \sum_{t=1}^T (\wt{\mbf F}_t-\mbf F_t) \Vert= O_p(1)$.
\end{compactenum}
\end{lem}

\noindent
\textsc{Proof.} For part (i), by definition 
\al{
\wt{\mbf F}_t-\mbf F_t=&\,(\wh{\bm\Lambda}_n^{(0)\prime}\wh{\bm\Lambda}_n^{(0)})^{-1}\wh{\bm\Lambda}_n^{(0)\prime}\mbf x_{nt}-\mbf F_t\nn\\
=&\,\l\{(\wh{\bm\Lambda}_n^{(0)\prime}\wh{\bm\Lambda}_n^{(0)})^{-1}\wh{\bm\Lambda}_n^{(0)\prime}\bm\Lambda_n-\mbf I_r\r\}\mbf F_t+ 
\l\{
(\wh{\bm\Lambda}_n^{(0)\prime}\wh{\bm\Lambda}_n^{(0)})^{-1}\wh{\bm\Lambda}_n^{(0)\prime}-
({\bm\Lambda}_n^{\prime}{\bm\Lambda}_n)^{-1}\bm\Lambda_n^{\prime}\r\}\bm \xi_{nt}\nn\\
&+({\bm\Lambda}_n^{\prime}{\bm\Lambda}_n)^{-1}\bm\Lambda_n^{\prime}\bm\xi_{nt}.\label{pdpd}
}
Using \eqref{pdpd} the first term on the rhs of \eqref{inuefottiti} is such that
\al{
\l\Vert T^{-1} \sum_{t=1}^T (\wt{\mbf F}_t-\mbf F_t){\mbf F}_{t}^\prime\r\Vert\le&\, 
\l\Vert T^{-1} \sum_{t=1}^T \l\{(\wh{\bm\Lambda}_n^{(0)\prime}\wh{\bm\Lambda}_n^{(0)})^{-1}\wh{\bm\Lambda}_n^{(0)\prime}\bm\Lambda_n-\mbf I_r\r\}\mbf F_t{\mbf F}_{t}^\prime\r\Vert\nn\\
&+\l\Vert T^{-1} \sum_{t=1}^T 
\sqrt n \l\{
(\wh{\bm\Lambda}_n^{(0)\prime}\wh{\bm\Lambda}_n^{(0)})^{-1}\wh{\bm\Lambda}_n^{(0)\prime}-
({\bm\Lambda}_n^{\prime}{\bm\Lambda}_n)^{-1}\bm\Lambda_n^{\prime}\r\}n^{-1/2}\bm \xi_{nt}
{\mbf F}_{t}^\prime\r\Vert\nn\\ 
&+\l\Vert T^{-1} \sum_{t=1}^T ({\bm\Lambda}_n^{\prime}{\bm\Lambda}_n)^{-1}\bm\Lambda_n^{\prime}\bm\xi_{nt}{\mbf F}_{t}^\prime\r\Vert\nn\\ 
\le&\, \l\{ \l\Vert T^{-1} \sum_{t=1}^T\mbf F_t{\mbf F}_{t}^\prime\r\Vert + \l\Vert n^{-1/2} T^{-1} \sum_{t=1}^T\bm \xi_{nt}
{\mbf F}_{t}^\prime\r\Vert\r\} O_p(\max(n^{-1},T^{-1/2}))\nn\\ 
&+  \Vert n ({\bm\Lambda}_n^{\prime}{\bm\Lambda}_n)^{-1}\Vert\,
\l\Vert n^{-1}T^{-1} \sum_{t=1}^T\bm\Lambda_n^\prime \bm\xi_{nt}{\mbf F}_{t}^\prime\r\Vert\nn\\
=&\, O_p(\max(n^{-1},T^{-1/2}))+O_p(n^{-1/2}T^{-1/2}),\label{eq:brioche}
}
by Lemma \ref{lem:est0LOAD}(ii), which does not depend on $t$, and Lemmas \ref{lem:lambdasqrtn} (jointly with Assumption \ref{ass:common}(a)), \ref{lem:consistCOV}(i), and \ref{lem:consistCOV}(iii), and part (i).
The case $k=1$ is proved in the same way and this proves part (iii).

For part (ii), as in part (i), we have
\al{
&\l\Vert n^{-1/2}T^{-1} \sum_{t=1}^T (\wt{\mbf F}_t-\mbf F_t){\bm \xi}_{nt}^\prime\r\Vert\nn\\
\le&\, 
 \l\{ \l\Vert n^{-1/2}T^{-1} \sum_{t=1}^T\mbf F_t{\bm \xi}_{nt}^\prime\r\Vert + \l\Vert n^{-1} T^{-1} \sum_{t=1}^T\bm \xi_{nt}
{\bm \xi}_{nt}^\prime\r\Vert\r\} O_p(\max(n^{-1},T^{-1/2}))\nn\\ 
&+  \Vert n ({\bm\Lambda}_n^{\prime}{\bm\Lambda}_n)^{-1}\Vert\,
\l\Vert n^{-3/2}T^{-1} \sum_{t=1}^T\bm\Lambda_n^\prime \bm\xi_{nt}{\bm \xi}_{nt}^\prime\r\Vert\nn\\
=&\, O_p(\max(n^{-1},T^{-1/2}))+ O_p(n^{-1/2}T^{-1/2}),\nn
}
by Lemma \ref{lem:est0LOAD}(ii), which does not depend on $t$, and Lemmas \ref{lem:lambdasqrtn} (jointly with Assumption \ref{ass:common}(a)), \ref{lem:consistCOV}(iii), and \ref{lem:consistCOV}(vi), and Lemma \ref{lem:COVFF0}(ii). This proves part (ii). 

Part (iii) follows directly from part (ii) but using Lemma \ref{lem:COVFF0}(iii).

For part (iv), as in part (i), we have
\al{
\l\Vert T^{-1} \sum_{t=1}^T (\wt{\mbf F}_t-\mbf F_t){\xi}_{it}\r\Vert\le&\, \l\{ \l\Vert T^{-1} \sum_{t=1}^T\mbf F_t{\xi}_{it}\r\Vert + \l\Vert n^{-1/2} T^{-1} \sum_{t=1}^T\bm \xi_{nt}
{\xi}_{it}\r\Vert\r\} O_p(\max(n^{-1},T^{-1/2}))\nn\\ 
&+  \Vert n ({\bm\Lambda}_n^{\prime}{\bm\Lambda}_n)^{-1}\Vert\,
\l\Vert n^{-1}T^{-1} \sum_{t=1}^T\bm\Lambda_n^\prime \bm\xi_{nt}{\xi}_{it}\r\Vert\nn\\
=&\, O_p(\max(n^{-1},T^{-1/2}))+O_p(n^{-1/2}T^{-1/2}),\nn
}
by the same arguments used to prove part (ii) and noticing that now $\vert\xi_{it}\vert =O_p(1)$ by Assumption \ref{ass:idio}(a). Part (v) is proved as part (iv) setting $\xi_{it}=1$.
This completes the proof. $\Box$
\begin{lem}\label{lem:est0VAR}
Consider the initial estimator of the VAR parameters $\wh{\mbf A}^{(0)}$ and $\wh{\bm\Gamma}^{v(0)}$ defined in Section \ref{app:prest}, then, under Assumptions \ref{ass:common}, \ref{ass:idio}, \ref{ass:ind}, and \ref{ass:ident}, as $n,T\to\infty$: 
\begin{compactenum}
\item [(i)] $\min( n,\sqrt T)\,\Vert\wh{\mbf A}^{(0)}-\mbf A\Vert = O_p(1)$;
\item [(ii)] $\min( n,\sqrt T)\,\Vert\wh{\bm\Gamma}^{v(0)}-\bm\Gamma^v\Vert = O_p(1)$.
\end{compactenum}
\end{lem}

\noindent\textsc{Proof.} 
For part (i), in agreement with Assumption \ref{ass:common}(i), we can always set $\wt{\mbf F}_0=\mbf 0_r$, so it follows that, by construction $T^{-1}\sum_{t=2}^T \wt{\mbf F}_{t-1}\wt{\mbf F}_{t-1}^\prime=T^{-1}\sum_{t=1}^T \wt{\mbf F}_{t-1}\wt{\mbf F}_{t-1}^\prime=\mbf I_r$. Then, by Assumption \ref{ass:ident}(b), we have
\al{
\wh{\mbf A}^{(0)}-\mbf A &= T^{-1}\sum_{t=2}^T \wt{\mbf F}_t\wt{\mbf F}_{t-1}^\prime -\bm\Gamma_1^F\nn\\
&=\l\{ T^{-1}\sum_{t=2}^T \wt{\mbf F}_t\wt{\mbf F}_{t-1}^\prime- T^{-1}\sum_{t=2}^T {\mbf F}_t{\mbf F}_{t-1}^\prime\r\}+ \l\{T^{-1}\sum_{t=2}^T {\mbf F}_t{\mbf F}_{t-1}^\prime -\bm\Gamma_1^F\r\}.\label{inueculo}
}
Now,
\al{
T^{-1}\sum_{t=2}^T \wt{\mbf F}_t\wt{\mbf F}_{t-1}^\prime- T^{-1}\sum_{t=2}^T {\mbf F}_t{\mbf F}_{t-1}^\prime=&\, T^{-1} \sum_{t=2}^T (\wt{\mbf F}_t-\mbf F_t){\mbf F}_{t-1}^\prime+T^{-1} \sum_{t=2}^T {\mbf F}_{t}(\wt{\mbf F}_{t-1}-\mbf F_{t-1})^\prime\nn\\
&+T^{-1} \sum_{t=2}^T (\wt{\mbf F}_t-\mbf F_t)(\wt{\mbf F}_{t-1}-\mbf F_{t-1})^\prime.\label{inuefottiti}
}
Now, by Lemma \ref{lem:COVFF}(i) the first and second term in the rhs of \eqref{inuefottiti} are $O_p(\max(n^{-1},T^{-1/2}))$, while the third term on the rhs is dominated by the first two. Hence, for the first term on the rhs of \eqref{inueculo}  we have
\al{
\l\Vert T^{-1}\sum_{t=2}^T \wt{\mbf F}_t\wt{\mbf F}_{t-1}^\prime- T^{-1}\sum_{t=2}^T {\mbf F}_t{\mbf F}_{t-1}^\prime
\r\Vert=O_p(\max(n^{-1},T^{-1/2})).\nn
}
The second term on the rhs of \eqref{inueculo} is $O_p(T^{-1/2})$, by Lemma \ref{lem:consistCOV}(i).
This proves part (i).

Part (ii) follows from the results in \citet[Proposition P]{FGLR09} combined with part (i).  This completes the proof. $\Box$

\begin{lem}\label{lem:est0}
Consider the initial estimator of the idiosyncratic variances $\wh{\sigma}_i^{2(0)}$, $i=1,\ldots, n$, defined in Section \ref{app:prest}, then, under Assumptions \ref{ass:common}, \ref{ass:idio}, \ref{ass:ind}, and \ref{ass:ident}, as $n,T\to\infty$: 
\begin{compactenum}
\item [(i)] $\min(n,\sqrt T)\,\vert\wh{\sigma}_{i}^{(0)2}-\sigma_i^2\vert = O_p(1)$, uniformly in  $i$;
\item [(ii)] $\min(n,\sqrt T)\,n^{-1}\vert\sum_{i=1}^n(\wh{\sigma}_{i}^{(0)2}-\sigma_i^2)\vert = O_p(1)$.
\end{compactenum}
\end{lem}

\noindent\textsc{Proof.} 
Start from
\begin{align}
(\wh{\sigma}_i^{2(0)}- \sigma_i^2)=&\, T^{-1}\sum_{t=1}^T (x_{it}-\wh{\bm\lambda}_i^{(0)\prime}\wt{\mbf F}_t)^2-\E[(x_{it}-\bm\lambda_i^\prime\mbf F_t)^2]\nn\\
=&\, \l\{T^{-1}\sum_{t=1}^T  x_{it}^2-\E[ x_{it}^2] \r\}+ \l\{
\wh{\bm\lambda}_i^{(0)\prime}\wh{\bm\lambda}_i^{(0)}- \bm\lambda_i^\prime\bm\lambda_i\r\}\nn\\
&-2\l\{ T^{-1}\sum_{t=1}^T\bm\lambda_i^\prime\mbf F_t \wt{\mbf F}_t^\prime\wh{\bm\lambda}_i^{(0)}  -\bm\lambda_i^\prime\bm\lambda_i\r\} -2\l\{ T^{-1}\sum_{t=1}^T \xi_{it}\wt{\mbf F}_t^\prime\wh{\bm\lambda}_i^{(0)}  \r\},\label{eq:spesa}
\end{align}
since $\bm\Gamma^F=\mbf I_r$ by Assumption \ref{ass:ident}(b),  $\E[{\mbf F}_t\xi_{it}]=\mbf 0_{r}$ by Lemma \ref{lem:fidio}, and $T^{-1}\sum_{t=1}^T \wt{\mbf F}_{t}\wt{\mbf F}_{t}^\prime=\mbf I_r$ by construction. 
For part (i), from \eqref{eq:spesa} we have
\al{
\vert\wh{\sigma}_i^{2(0)}- \sigma_i^2\vert\le &\, \l\vert T^{-1}\sum_{t=1}^T  x_{it}^2-\E[ x_{it}^2] \r\vert+ \l \vert
\wh{\bm\lambda}_i^{(0)\prime}\wh{\bm\lambda}_i^{(0)}- \bm\lambda_i^\prime\bm\lambda_i\r\vert\nn\\
&+2\l\vert T^{-1}\sum_{t=1}^T\bm\lambda_i^\prime\mbf F_t {\mbf F}_t^\prime\wh{\bm\lambda}_i^{(0)}  -\bm\lambda_i^\prime\bm\lambda_i\r\vert
+2\l\vert T^{-1}\sum_{t=1}^T\bm\lambda_i^\prime\mbf F_t (\wt{\mbf F}_t-\mbf F_t)^\prime\wh{\bm\lambda}_i^{(0)}\r\vert\nn\\
& +2\l\vert T^{-1}\sum_{t=1}^T \xi_{it}{\mbf F}_t^\prime\wh{\bm\lambda}_i^{(0)}  \r\vert 
+2\l\vert T^{-1}\sum_{t=1}^T \xi_{it}(\wt{\mbf F}_t-\mbf F_t)^\prime\wh{\bm\lambda}_i^{(0)}  \r\vert\nn\\
\le &\, \l\vert T^{-1}\sum_{t=1}^T  x_{it}^2-\E[ x_{it}^2] \r\vert+ \l \vert
\wh{\bm\lambda}_i^{(0)\prime}\wh{\bm\lambda}_i^{(0)}- \bm\lambda_i^\prime\bm\lambda_i\r\vert\nn\\
&+2\Vert \bm\lambda_i\Vert\, \Vert \wh{\bm\lambda}_i^{(0)}  -\bm\lambda_i\Vert
+2\Vert \bm\lambda_i\Vert\, \l\Vert T^{-1}\sum_{t=1}^T\mbf F_t {\mbf F}_t^\prime-\mbf I_r\r\Vert\, \{\Vert\wh{\bm\lambda}_i^{(0)}-\bm\lambda_i\Vert+\Vert\bm\lambda_i\Vert \} \nn\\
&+2 \Vert \bm\lambda_i\Vert\, \l\Vert T^{-1}\sum_{t=1}^T \mbf F_t (\wt{\mbf F}_t-\mbf F_t)^\prime\r\Vert\,\{\Vert \wh{\bm\lambda}_i^{(0)}-\bm\lambda_i\Vert +\Vert\bm\lambda_i\Vert\}\nn\\
& +2\l\Vert T^{-1}\sum_{t=1}^T \xi_{it}{\mbf F}_t^\prime\r\Vert\,\Vert{\bm\lambda}_i \Vert
+2\l\Vert T^{-1}\sum_{t=1}^T \xi_{it}{\mbf F}_t^\prime \r\Vert\,\Vert\wh{\bm\lambda}_i^{(0)}-\bm\lambda_i  \Vert\nn\\
& +2\l\Vert T^{-1}\sum_{t=1}^T \xi_{it}(\wt{\mbf F}_t-\mbf F_t)^\prime\r\Vert \, \{\Vert \wh{\bm\lambda}_i^{(0)} -\bm\lambda_i\Vert+\Vert\bm\lambda_i\Vert\}\nn\\
=&\,O_p(\max(n^{-1},T^{-1/2})),\label{eq:laringite}
}
where, we used multiple times Assumption \ref{ass:common}(a) and Lemma \ref{lem:est0LOAD}(i), and, we also used:
Lemma \ref{lem:consistCOV}(v) for the first term,
Lemma \ref{lem:consistCOV}(i) for the  fourth term,
Lemma \ref{lem:COVFF}(i) for the fifth term,
Lemma \ref{lem:consistCOV}(ii) for the sixth and seventh term, and Lemma \ref{lem:COVFF}(iv) for the last term. This proves part (i).

As for part (ii),  from \eqref{eq:spesa} we have
\al{
n^{-1}\l\vert\sum_{i=1}^n (\wh{\sigma}_i^{2(0)}-\sigma_i^2)\r\vert \le&\, 
n^{-1}\l\vert 
\sum_{i=1}^n 
\l\{T^{-1}\sum_{t=1}^T x_{it}^2-\E[x_{it}^2]\r\}
\r\vert +
n^{-1} \l\vert \sum_{i=1}^n \{ \wh{\bm\lambda}_i^{(0)\prime}\wh{\bm\lambda}_i^{(0)}-\bm\lambda_i^\prime\bm\lambda_i\}\r\vert\nn\\
&+ 2n^{-1} \l\vert \sum_{i=1}^n \l\{T^{-1}\sum_{t=1}^T\bm\lambda_i^\prime \mbf F_t\mbf F_t^\prime \wh{\bm\lambda}_i^{(0)} -\bm\lambda_i^\prime\bm\lambda_i \r\} 
\r\vert\nn\\
&+2n^{-1}\l\vert \sum_{i=1}^n\l\{  T^{-1}\sum_{t=1}^T\bm\lambda_i^\prime \mbf F_t(\wt{\mbf F}_t-\mbf F_t)^\prime \wh{\bm\lambda}_i^{(0)}\r\}\r\vert\nn\\
&+2n^{-1}\l\vert \sum_{i=1}^n \l\{T^{-1}\sum_{t=1}^T \xi_{it} \mbf F_t^\prime\r\}\wh{\bm\lambda}_i^{(0)}\r\vert 
+2n^{-1}\l\vert \sum_{i=1}^n\l\{T^{-1}\sum_{t=1}^T \xi_{it} (\wt{\mbf F}_t-\mbf F_t)^\prime\r\}\wh{\bm\lambda}_i^{(0)}\r\vert\nn\\
=&\, 
n^{-1}\l\vert \sum_{i=1}^n \l\{T^{-1}\sum_{t=1}^T x_{it}^2-\E[x_{it}^2]\r\}\r\vert +n^{-1} \l\vert \sum_{i=1}^n \{ \wh{\bm\lambda}_i^{(0)\prime}\wh{\bm\lambda}_i^{(0)}-\bm\lambda_i^\prime\bm\lambda_i\}\r\vert\nn\\
&+ 2 n^{-1/2} \Vert\bm\Lambda_n\Vert\, 
n^{-1/2} \l\Vert\wh{\bm\Lambda}_n^{(0)} -\bm\Lambda_n\r\Vert+ 2n^{-1}\Vert\bm\Lambda_n^\prime\bm\Lambda_n\Vert\,  \l\Vert T^{-1}\sum_{t=1}^T \mbf F_t\mbf F_t^\prime-\mbf I_r \r\Vert \nn\\
&+ 2 n^{-1/2} \Vert\bm\Lambda_n\Vert\, 
n^{-1/2} \l\Vert\wh{\bm\Lambda}_n^{(0)} -\bm\Lambda_n\r\Vert\, \l\Vert T^{-1}\sum_{t=1}^T \mbf F_t\mbf F_t^\prime-\mbf I_r \r\Vert \nn\\
&+  2 n^{-1}\Vert {\bm\Lambda}_n^{\prime}{\bm\Lambda}_n\Vert\, \l\Vert T^{-1}\sum_{t=1}^T \mbf F_t (\wt{\mbf F}_t-\mbf F_t)^\prime\r\Vert\nn\\
&+  2 n^{-1/2}\Vert\bm\Lambda_n\Vert\, n^{-1/2}\Vert \wh{\bm\Lambda}_n^{(0)}-{\bm\Lambda}_n\Vert\, \l\Vert T^{-1}\sum_{t=1}^T \mbf F_t (\wt{\mbf F}_t-\mbf F_t)^\prime\r\Vert\nn\\
&+2n^{-1}\l\vert T^{-1}\sum_{t=1}^T  \mbf F_t^\prime{\bm\Lambda}_n^\prime \bm\xi_{nt}\r\vert 
+2n^{-1}\l\vert T^{-1}\sum_{t=1}^T  \mbf F_t^\prime(\wh{\bm\Lambda}_n^{(0)}-\bm\Lambda_n)^\prime\bm\xi_{nt}\r\vert \nn\\
&+2n^{-1}\l\vert T^{-1}\sum_{t=1}^T  (\wt{\mbf F}_t-\mbf F_t)^\prime {\bm\Lambda}_n^\prime\bm\xi_{nt}\r\vert\nn\\
& +2n^{-1}\l\vert T^{-1}\sum_{t=1}^T (\wt{\mbf F}_t-\mbf F_t)^\prime (\wh{\bm\Lambda}_n^{(0)}-\bm\Lambda_n)^\prime \bm\xi_{nt} \r\vert\nn\\
=&\, O_p(\max(n^{-1},T^{-1/2})). \label{eq:ceramica}
}
The result in \eqref{eq:ceramica} follows from repeated use of Lemmas \ref{lem:lambdasqrtn}, \ref{lem:consistCOV}(i), \ref{lem:est0LOAD}(ii), and \ref{lem:COVFF}(i) as well as the following results. First,
 \al{
 n^{-1}\l\vert 
\sum_{i=1}^n 
\l\{T^{-1}\sum_{t=1}^T x_{it}^2-\E[x_{it}^2]\r\}\r\vert 
\le n^{-1}\l\Vert T^{-1}\sum_{t=1}^T \mbf x_{nt}\mbf x_{nt}^\prime -\E[\mbf x_{nt}\mbf x_{nt}^\prime]\r\Vert_F = O_p(T^{-1/2}),\nn
 } 
 by Lemma \ref{lem:consistCOV}(vii). Second,
 \al{
n^{-1} \l\vert T^{-1}\sum_{t=1}^T \mbf F_t^\prime{\bm\Lambda}_n^{\prime}\bm\xi_{nt}\r\vert 
&\le n^{-1}\l\Vert T^{-1}\sum_{t=1}^T{\bm\Lambda}_n^{\prime}\bm\xi_{nt} \mbf F_t^\prime\r\Vert_F 
= O_p(n^{-1/2}T^{-1/2}),\nn
}
 by Lemma \ref{lem:COVFF0}(i). Third, 
 \al{
n^{-1}\l\vert T^{-1}\sum_{t=1}^T  (\wt{\mbf F}_t-\mbf F_t)^\prime {\bm\Lambda}_n^\prime\bm\xi_{nt}\r\vert
&\le n^{-1}\l\Vert  T^{-1}\sum_{t=1}^T{\bm\Lambda}_n^\prime\bm\xi_{nt}(\wt{\mbf F}_t-\mbf F_t)'
\r\Vert_F\nn\\
&\le \sqrt r n^{-1}\l\Vert  T^{-1}\sum_{t=1}^T{\bm\Lambda}_n^\prime\bm\xi_{nt}(\wt{\mbf F}_t-\mbf F_t)'
\r\Vert\nn\\
&= O_p(\max(n^{-1},T^{-1/2})),\nn
 }
 by Lemma \ref{lem:COVFF}(iii). Fourth,
 \al{
 n^{-1}\l\vert T^{-1}\sum_{t=1}^T  \mbf F_t^\prime(\wh{\bm\Lambda}_n^{(0)}-\bm\Lambda_n)^\prime\bm\xi_{nt}\r\vert &\le 
 n^{-1}\l\Vert T^{-1}\sum_{t=1}^T  (\wh{\bm\Lambda}_n^{(0)}-\bm\Lambda_n)^\prime\bm\xi_{nt}\mbf F_t^\prime\r\Vert_F\nn\\
& \le \sqrt r n^{-1}\l\Vert T^{-1}\sum_{t=1}^T  (\wh{\bm\Lambda}_n^{(0)}-\bm\Lambda_n)^\prime\bm\xi_{nt}\mbf F_t^\prime\r\Vert\nn\\
&\le \sqrt r n^{-1/2}\Vert \wh{\bm\Lambda}_n^{(0)}-\bm\Lambda_n\Vert \, n^{-1/2}\l\Vert T^{-1}\sum_{t=1}^T\bm\xi_{nt}\mbf F_t^\prime\r\Vert\nn\\
&= O_p(\max(n^{-1},T^{-1/2}))+ O_p(T^{-1/2}),\nn
 }
 by Lemmas \ref{lem:consistCOV}(iii) and \ref{lem:est0LOAD}(ii). And, last
 \al{
 n^{-1}\l\vert T^{-1}\sum_{t=1}^T (\wt{\mbf F}_t-\mbf F_t)^\prime (\wh{\bm\Lambda}_n^{(0)}-\bm\Lambda_n)^\prime \bm\xi_{nt} \r\vert&\le 
 n^{-1}\l\Vert T^{-1}\sum_{t=1}^T  (\wh{\bm\Lambda}_n^{(0)}-\bm\Lambda_n)^\prime \bm\xi_{nt} (\wt{\mbf F}_t-\mbf F_t)^\prime \r\Vert_F\nn\\
 &\le \sqrt r  n^{-1}\l\Vert T^{-1}\sum_{t=1}^T  (\wh{\bm\Lambda}_n^{(0)}-\bm\Lambda_n)^\prime \bm\xi_{nt} (\wt{\mbf F}_t-\mbf F_t)^\prime \r\Vert\nn\\
&\le \sqrt r  n^{-1/2}\Vert \wh{\bm\Lambda}_n^{(0)}-\bm\Lambda_n\Vert \, n^{-1/2}\l\Vert T^{-1}\sum_{t=1}^T \bm\xi_{nt} (\wt{\mbf F}_t-\mbf F_t)^\prime \r\Vert\nn\\
&= O_p(\max(n^{-1},T^{-1/2})),\nn
 }
 by Lemmas \ref{lem:consistCOV}(iii) and \ref{lem:COVFF}(ii). This completes the proof. $\Box$

\begin{lem}\label{lem:est0_LAST}
Under Assumptions \ref{ass:common}, \ref{ass:idio}, \ref{ass:ind}, and \ref{ass:ident}, as $n,T\to\infty$:
\begin{compactenum}
\item [(i)] $\min(n,\sqrt T)\,n^{-1}\Vert\wh{\bm\Lambda}_n^{(0)\prime}(\wh{\bm \Sigma}_n^{\xi(0)})^{-1}\wh{\bm\Lambda}_n^{(0)}-\bm\Lambda_n^\prime(\bm\Sigma_n^\xi)^{-1}\bm\Lambda_n\Vert = O_p(1)$;
\item [(ii)] $\min(n,\sqrt T)\,n^{-1/2}\Vert\wh{\bm\Lambda}_n^{(0)\prime}(\wh{\bm \Sigma}_n^{\xi(0)})^{-1}-\bm\Lambda_n^\prime(\bm\Sigma_n^\xi)^{-1}\Vert = O_p(1)$;
\item [(iii)] $n\Vert(\wh{\bm\Lambda}_n^{(0)\prime}(\wh{\bm \Sigma}_n^{\xi(0)})^{-1}\wh{\bm\Lambda}_n^{(0)})^{-1}\Vert = O_p(1)$;
\item [(iv)] $\min(n,\sqrt T)\,n\Vert(\wh{\bm\Lambda}_n^{(0)\prime}(\wh{\bm \Sigma}_n^{\xi(0)})^{-1}\wh{\bm\Lambda}_n^{(0)})^{-1}-(\bm\Lambda_n^\prime(\bm\Sigma_n^\xi)^{-1}\bm\Lambda_n)^{-1}\Vert = O_p(1)$;
\item [(v)] $\min(n,\sqrt T)\,\sqrt n\Vert(\wh{\bm\Lambda}_n^{(0)\prime}(\wh{\bm \Sigma}_n^{\xi(0)})^{-1}\wh{\bm\Lambda}_n^{(0)})^{-1}
\wh{\bm\Lambda}_n^{(0)\prime}(\wh{\bm \Sigma}_n^{\xi(0)})^{-1}
-(\bm\Lambda_n^\prime(\bm\Sigma_n^\xi)^{-1}\bm\Lambda_n)^{-1}\bm\Lambda_n^\prime(\bm\Sigma_n^\xi)^{-1}\Vert = O_p(1)$.
\end{compactenum}
\end{lem}

\noindent\textsc{Proof.} 
Start with
\al{
n^{-1}\Vert\wh{\bm\Lambda}_n^{(0)\prime}(\wh{\bm \Sigma}_n^{\xi(0)})^{-1}\wh{\bm\Lambda}_n^{(0)}&\,-\bm\Lambda_n^\prime(\bm\Sigma_n^\xi)^{-1}\bm\Lambda_n\Vert 
\le 
2 n^{-1}\Vert \{\wh{\bm\Lambda}_n^{(0)}-\bm\Lambda_n\}^\prime({\bm\Sigma}_n^\xi)^{-1}{\bm\Lambda}_n\Vert\nn\\
  &+n^{-1}\Vert  {\bm\Lambda}_n^\prime \{(\wh{\bm\Sigma}_n^{\xi(0)})^{-1} -({\bm\Sigma}_n^\xi)^{-1} \} {\bm\Lambda}_n\Vert\nn\\
 &+2n^{-1}\Vert
 \{\wh{\bm\Lambda}_n^{(0)}-\bm\Lambda_n\}^\prime \{(\wh{\bm\Sigma}_n^{\xi(0)})^{-1} -({\bm\Sigma}_n^\xi)^{-1} \} {\bm\Lambda}_n
 \Vert\nn\\
  &+n^{-1}\Vert
 \{\wh{\bm\Lambda}_n^{(0)}-\bm\Lambda_n\}^\prime \{(\wh{\bm\Sigma}_n^{\xi(0)})^{-1} -({\bm\Sigma}_n^\xi)^{-1} \} \{\wh{\bm\Lambda}_n^{(0)}-\bm\Lambda_n\}\nn
 \Vert\nn\\
 \le &\, 2 n^{-1/2}\Vert \wh{\bm\Lambda}_n^{(0)}-\bm\Lambda_n\Vert \, \Vert({\bm\Sigma}_n^\xi)^{-1}\Vert \, n^{-1/2}\Vert {\bm\Lambda}_n\Vert\nn\\
  &+n^{-1}\Vert  {\bm\Lambda}_n^\prime \{(\wh{\bm\Sigma}_n^{\xi(0)})^{-1} -({\bm\Sigma}_n^\xi)^{-1} \} {\bm\Lambda}_n\Vert\nn\\
&+2n^{-1/2}\Vert
 \wh{\bm\Lambda}_n^{(0)}-\bm\Lambda_n \Vert n^{-1/2} \Vert \{(\wh{\bm\Sigma}_n^{\xi(0)})^{-1} -({\bm\Sigma}_n^\xi)^{-1} \} {\bm\Lambda}_n
 \Vert\nn\\
  &+n^{-1}\Vert \wh{\bm\Lambda}_n^{(0)}- {\bm\Lambda}_n\Vert^2 \Vert (\wh{\bm\Sigma}_n^{\xi(0)})^{-1} -({\bm\Sigma}_n^\xi)^{-1}\Vert.
\label{eq:renollobrutto}
}
Consider each term on the rhs of \eqref{eq:renollobrutto}. The first term is
\al{
n^{-1/2}\Vert \wh{\bm\Lambda}_n^{(0)}-\bm\Lambda_n\Vert \, \Vert({\bm\Sigma}_n^\xi)^{-1}\Vert \, n^{-1/2}\Vert {\bm\Lambda}_n\Vert= O_p(\max(n^{-1},T^{-1/2})),\label{eq:RB1}
}
by Lemma \ref{lem:est0LOAD}(ii), and also Assumption \ref{ass:idio}(a) for which $\Vert({\bm\Sigma}_n^\xi)^{-1}\Vert\le C_\xi$, and Lemma \ref{lem:lambdasqrtn} for which $n^{-1/2}\Vert {\bm\Lambda}_n\Vert=O(1)$. The third term is
\al{
n^{-1/2}\Vert
 \wh{\bm\Lambda}_n^{(0)}&-\bm\Lambda_n \Vert n^{-1/2} \Vert \{(\wh{\bm\Sigma}_n^{\xi(0)})^{-1} -({\bm\Sigma}_n^\xi)^{-1} \} {\bm\Lambda}_n
 \Vert\nn\\
 &\le n^{-1/2}\Vert
 \wh{\bm\Lambda}_n^{(0)}-\bm\Lambda_n \Vert \l\{ \Vert (\wh{\bm\Sigma}_n^{\xi(0)})^{-1}\Vert + \Vert ({\bm\Sigma}_n^\xi)^{-1}\Vert \r\} n^{-1/2}\Vert  {\bm\Lambda}_n\Vert\nn\\
& = O_p(\max(n^{-1},T^{-1/2})),\label{eq:RB3}
}
by Lemmas  \ref{lem:lambdasqrtn} and \ref{lem:est0LOAD}(ii), Assumption \ref{ass:idio}(a), and since, by Lemma \ref{lem:est0}(i), for all $j=1,\ldots, n$,
\al{
\Vert (\wh{\bm\Sigma}_n^{\xi(0)})^{-1}\Vert &= \l\{ \nu^{(n)}(\wh{\bm\Sigma}_n^{\xi(0)}) \r\}^{-1} = \l\{\min_{i=1,\ldots, n} \wh{\sigma}_i^{2(0)}\r\}^{-1}
= \l\{\min_{i=1,\ldots, n} {\sigma}_i^{2}+\wh{\sigma}_i^{2(0)}-{\sigma}_i^{2}\r\}^{-1}\nn\\
&\le \l\{\min_{i=1,\ldots, n} {\sigma}_i^{2}+\min_{i=1,\ldots, n} (\wh{\sigma}_i^{2(0)}-{\sigma}_i^{2})\r\}^{-1}\nn\\
&\le \l\{\min_{i=1,\ldots, n} {\sigma}_i^{2} - \min_{i=1,\ldots, n} \vert \wh{\sigma}_i^{2(0)}-\sigma_i^2 \vert \r\}^{-1}\nn\\
&\le \l\{C_\xi^{-1} - \vert \wh{\sigma}_j^{2(0)}-\sigma_j^2 \vert \r\}^{-1} \le C_\xi  + O_p(\max(n^{-1},T^{-1/2})).\label{eq:SHATPD}
}
By the same arguments, the fourth term is
\al{
n^{-1}\Vert
 \wh{\bm\Lambda}_n^{(0)}&-\bm\Lambda_n \Vert^2  \l\{ \Vert (\wh{\bm\Sigma}_n^{\xi(0)})^{-1}\Vert + \Vert ({\bm\Sigma}_n^\xi)^{-1}\Vert \r\} = o_p(\max(n^{-1},T^{-1/2})). \label{eq:RB4}
}
Finally, the second term on the rhs of \eqref{eq:renollobrutto} is
\al{
n^{-1}\Vert  {\bm\Lambda}_n^\prime \{(\wh{\bm\Sigma}_n^{\xi(0)})^{-1} -({\bm\Sigma}_n^\xi)^{-1} \} {\bm\Lambda}_n\Vert =&\, n^{-1}\l\Vert
\sum_{i=1}^n \bm\lambda_i\bm\lambda_i^\prime \{(\wh{\sigma}_i^{2(0)})^{-1}-(\sigma_i^2)^{-1}\}
\r\Vert\nn\\
\le &\,n^{-1} \l(\sum_{k,h=1}^r  \l[\sum_{i=1}^n \lambda_{ih}\lambda_{ik} \{(\wh{\sigma}_i^{2(0)})^{-1}-(\sigma_i^2)^{-1}\}\r]^2  \r)^{1/2}\nn\\
\le&\, r M_\lambda n^{-1}\l\vert\sum_{i=1}^n 
\{(\wh{\sigma}_i^{2(0)})^{-1}-(\sigma_i^2)^{-1}\}
\r\vert\nn\\
=&\, r M_\lambda n^{-1}\l\vert\sum_{i=1}^n (\wh{\sigma}_i^{2(0)})^{-1}(\sigma_i^2)^{-1}
\{\wh{\sigma}_i^{2(0)}-\sigma_i^2\}
\r\vert\nn\\
=&\, r M_\lambda C_\xi n^{-1}\l\vert
\l\{\min_{i=1,\ldots, n} \wh{\sigma}_i^{2(0)}\r\}^{-1}
\sum_{i=1}^n 
\{\wh{\sigma}_i^{2(0)}-\sigma_i^2\}
\r\vert
\nn\\
=&\, r M_\lambda C_\xi^2 n^{-1}\l\vert
\sum_{i=1}^n 
\{\wh{\sigma}_i^{2(0)}-\sigma_i^2\}
\r\vert
\nn\\
=&\,  O_p(\max(n^{-1},T^{-1/2})),  \label{eq:renollobrutto2}
}
by Lemma \ref{lem:est0}(ii), \eqref{eq:SHATPD}, and Assumptions \ref{ass:common}(a) and \ref{ass:idio}(a). By substituting \eqref{eq:RB1}, \eqref{eq:RB3}, \eqref{eq:RB4}, and  \eqref{eq:renollobrutto2} into \eqref{eq:renollobrutto}, we prove part (i).

For part (ii), we have
\al{
n^{-1}&\Vert\wh{\bm\Lambda}_n^{(0)}(\wh{\bm\Sigma}_n^{\xi(0)})^{-1}-{\bm\Lambda}_n({\bm\Sigma}_n^\xi)^{-1} \Vert^2 = n^{-1}\l\Vert\sum_{i=1}^n \wh{\bm\lambda}_i^{(0)} (\wh{\sigma}_i^{2(0)})^{-1}-\sum_{i=1}^n \bm\lambda_i(\sigma_i^2)^{-1}\r\Vert^2\nn\\
=&\, n^{-1}\l\Vert\sum_{i=1}^n \wh{\bm\lambda}_i^{(0)} (\wh{\sigma}_i^{2(0)}-\sigma_i^2+\sigma_i^2)^{-1}-\sum_{i=1}^n \bm\lambda_i(\sigma_i^2)^{-1}\r\Vert^2\nn\\
=&\, n^{-1}\l\Vert\sum_{i=1}^n \wh{\bm\lambda}_i^{(0)}\sigma_i^2 \{(\wh{\sigma}_i^{2(0)}-\sigma_i^2)/\sigma_i^2+1\}^{-1}-\sum_{i=1}^n \bm\lambda_i(\sigma_i^2)^{-1}\r\Vert^2\nn\\
\le&\, n^{-1}\l\Vert \l\{\min_{i=1,\ldots, n}(\wh{\sigma}_i^{2(0)}-\sigma_i^2)/\sigma_i^2+1\r\}^{-1} \sum_{i=1}^n \wh{\bm\lambda}_i^{(0)}(\sigma_i^2)^{-1} -\sum_{i=1}^n \bm\lambda_i(\sigma_i^2)^{-1}\r\Vert^2\nn\\
=&\, n^{-1}\l\Vert \l\{1-\min_{i=1,\ldots, n}(\wh{\sigma}_i^{2(0)}-\sigma_i^2)/\sigma_i^2+o\l(\min_{i=1,\ldots, n}(\wh{\sigma}_i^{2(0)}-\sigma_i^2)/\sigma_i^2\r)\r\} \sum_{i=1}^n \wh{\bm\lambda}_i^{(0)}(\sigma_i^2)^{-1} - \sum_{i=1}^n \bm\lambda_i(\sigma_i^2)^{-1}\r\Vert^2\nn\\
\le &\, n^{-1}\l\Vert \wh{\bm\Lambda}_n^{(0)}({\bm\Sigma}_n^{\xi})^{-1}-{\bm\Lambda}_n({\bm\Sigma}_n^\xi)^{-1}\r\Vert^2
+\l\{ \min_{i=1,\ldots, n}\l\vert(\wh{\sigma}_i^{2(0)}-\sigma_i^2)/\sigma_i^2\r\vert\r\}^2 \l\Vert \wh{\bm\Lambda}_n^{(0)}({\bm\Sigma}_n^{\xi})^{-1}\r\Vert^2\nn\\
\le &\, n^{-1}\l\Vert \wh{\bm\Lambda}_n^{(0)}({\bm\Sigma}_n^{\xi})^{-1}-{\bm\Lambda}_n({\bm\Sigma}_n^\xi)^{-1}\r\Vert^2
+ \l\{\min_{i=1,\ldots, n}\vert \wh{\sigma}_i^{2(0)}-\sigma_i^2\vert \r\}^2 L_\xi^2\, n^{-1}\Vert {\bm\Lambda}_n\Vert^2 \Vert({\bm\Sigma}_n^{\xi})^{-1}\Vert^2\nn\\
&+ \l\{\min_{i=1,\ldots, n}\l\vert\wh{\sigma}_i^{2(0)}-\sigma_i^2\r\vert\r\}^2 L_\xi^2\, n^{-1}\Vert \wh{\bm\Lambda}_n^{(0)}-\bm\Lambda_n\Vert^2 \Vert({\bm\Sigma}_n^{\xi})^{-1}\Vert^2\nn\\
=&\, O_p(\max(n^{-2},T^{-1})),\nn
}
by Lemmas \ref{lem:lambdasqrtn}, \ref{lem:est0LOAD}(ii), \ref{lem:est0}(ii), and Assumptions \ref{ass:idio}(a) and \ref{ass:idio}(f). This proves part (ii).

For part (iii), by part (ii) and \citet[Theorem 1]{MK04} which is Weyl's inequality, we have
\al{
n^{-1}\vert \nu^{(r)}
(\wh{\bm\Lambda}_n^{(0)\prime}(\wh{\bm \Sigma}_n^{\xi(0)})^{-1}\wh{\bm\Lambda}_n^{(0)})
-
\nu^{(r)}({\bm\Lambda}_n^{\prime}({\bm \Sigma}_n^{\xi})^{-1}){\bm\Lambda}_n)\vert\le 
&\, n^{-1}\Vert \wh{\bm\Lambda}_n^{(0)\prime}(\wh{\bm \Sigma}_n^{\xi(0)})^{-1}\wh{\bm\Lambda}_n^{(0)}
- {\bm\Lambda}_n^{\prime}({\bm \Sigma}_n^{\xi})^{-1}{\bm\Lambda}_n\Vert \nn\\
=&\, O_p(\max(n^{-1},T^{-1/2})).\label{eq:CRUCIAL0}
}
Moreover (note that $x-y\ge -|x-y|$ for any $x,y\in\mathbb R$),
\al{
\det ( \wh{\bm\Lambda}_n^{(0)\prime}(\wh{\bm \Sigma}_n^{\xi(0)})^{-1}\wh{\bm\Lambda}_n^{(0)}) =&\,\prod_{j=1}^r \nu^{(j)}
(\wh{\bm\Lambda}_n^{(0)\prime}(\wh{\bm \Sigma}_n^{\xi(0)})^{-1}\wh{\bm\Lambda}_n^{(0)})\ge \l\{ \nu^{(r)}
(\wh{\bm\Lambda}_n^{(0)\prime}(\wh{\bm \Sigma}_n^{\xi(0)})^{-1}\wh{\bm\Lambda}_n^{(0)}) \r\}^r\nn\\
\ge&\, \l\{ \nu^{(r)}
({\bm\Lambda}_n^{\prime}({\bm \Sigma}_n^{\xi})^{-1}{\bm\Lambda}_n) -
\vert\nu^{(r)}
(\wh{\bm\Lambda}_n^{(0)\prime}(\wh{\bm \Sigma}_n^{\xi(0)})^{-1}\wh{\bm\Lambda}_n^{(0)})
-
\nu^{(r)}({\bm\Lambda}_n^{\prime}({\bm \Sigma}_n^{\xi})^{-1}){\bm\Lambda}_n)\vert \r\}^r,\nn
}
thus, by Lemma \ref{lem:LSL2}(iv), which implies $\lim_{n\to\infty}n^{-1}\nu^{(r)}
({\bm\Lambda}_n^{\prime}({\bm \Sigma}_n^{\xi})^{-1}{\bm\Lambda}_n)>0$,
 and \eqref{eq:CRUCIAL0}, it follows that, with probability tending to one as $n,T\to\infty$, we have $\det (n^{-1} \wh{\bm\Lambda}_n^{(0)\prime}(\wh{\bm \Sigma}_n^{\xi(0)})^{-1}\wh{\bm\Lambda}_n^{(0)})>0$, or, equivalently $n^{-1} \wh{\bm\Lambda}_n^{(0)\prime}(\wh{\bm \Sigma}_n^{\xi(0)})^{-1}\wh{\bm\Lambda}_n^{(0)}$ is positive definite, i.e.,
 $n\Vert(\wh{\bm\Lambda}_n^{(0)\prime}(\wh{\bm \Sigma}_n^{\xi(0)})^{-1}\wh{\bm\Lambda}_n^{(0)})^{-1}\Vert=O_p(1)$.  This proves part (iii).

For part (iv), we have
\al{
n\Vert&(\wh{\bm\Lambda}_n^{(0)\prime}(\wh{\bm \Sigma}_n^{\xi(0)})^{-1}\wh{\bm\Lambda}_n^{(0)})^{-1}-(\bm\Lambda_n^\prime(\bm\Sigma_n^\xi)^{-1}\bm\Lambda_n)^{-1}\Vert\nn\\
&\le
n\Vert(\wh{\bm\Lambda}_n^{(0)\prime}(\wh{\bm \Sigma}_n^{\xi(0)})^{-1}\wh{\bm\Lambda}_n^{(0)})^{-1}\Vert\,
n^{-1}\Vert \wh{\bm\Lambda}_n^{(0)\prime}(\wh{\bm \Sigma}_n^{\xi(0)})^{-1}\wh{\bm\Lambda}_n^{(0)}-\bm\Lambda_n^\prime(\bm\Sigma_n^\xi)^{-1}\bm\Lambda_n\Vert\,
n\Vert (\bm\Lambda_n^\prime(\bm\Sigma_n^\xi)^{-1}\bm\Lambda_n)^{-1}\Vert\nn\\
&=O_p(\max(n^{-1},T^{-1/2})),\nn
}
because of parts (i) and (iii) and Lemma \ref{lem:LSL2}(iii). 

Part (v) follows directly from parts (ii) and (iv). This completes the proof. $\Box$


\begin{lem} \label{lem:detP} 
For all $T\in\mathbb N$,
\begin{compactenum}[(i)]
\item $\mbf P_{t|t-1}$, $\mbf P_{t|t}$, and $\mbf P_{t|T}$ are deterministic $r\times r$ matrices,  for all $t=1,\ldots, T$;
\item  $\Vert \mbf P_{t+1|t}\Vert\le \Vert \mbf P_{t|t-1}\Vert$,  for all $t=1,\ldots, T-1$;
\item $\mbf P_{0,t|t-1}$, $\mbf P_{0,t|t}$, and $\mbf P_{0,t|T}$ are deterministic $r\times r$ matrices,   for all $t=1,\ldots, T$;
\item $\Vert \mbf P_{0,t+1|t}\Vert\le \Vert \mbf P_{0,t|t-1}\Vert$, for all $t=1,\ldots, T-1$. 
\end{compactenum}
\end{lem}

\noindent\textsc{Proof.} Since $\mbf P_{0|0}=\mbf I_r$ is deterministic, then, for all $t=1,\ldots, T$, $\mbf P_{t|t-1}$, $\mbf P_{t|t}$, and $\mbf P_{t|T}$ do not depend on the actual observations  because of \eqref{eq:pred2} and \eqref{eq:up2}. This proves part (i).

As for part (ii), since $\mbf F_{t|t-1}$ is based on less information than $\mbf F_{t+1|t}$ and since ${\mbf P}_{t|t-1}$ and $\mbf P_{t+1|t}$ are deterministic, then $({\mbf P}_{t|t-1} - {\mbf P}_{t+1|t})$ is a positive definite matrix for all $t=1,\ldots, T-1$ (see, e.g., \citealp[Chapter 3.3, p. 123]{harvey90}). As consequence, $\Vert \mbf P_{t+1|t}\Vert\le \Vert \mbf P_{t|t-1}\Vert$ for all $t=1,\ldots, T-1$ (see, e.g., \citealp[Proposition L1, p. 360]{MOA11}). 

The proof of parts (iii) and (iv) is identical to parts (i) and (ii), respectively, since also in this case $\mbf P_{0,0|0}=\mbf I_r$. This completes the proof.
$\Box$


\begin{lem}\label{lem:cazzarola}
 Under Assumptions \ref{ass:common} and \ref{ass:idio}, for all $T\in\mathbb N$,
\begin{compactenum}[(i)] 
\item $\max_{t=1,\ldots, T}\Vert \mbf P_{t|t-1}\Vert\le M_P$ for some finite positive real $M_P$; 
\item $\min_{t=1,\ldots, T}\nu^{(r)}(\mbf P_{t|t-1})\ge \underline M_P$ for some finite positive real $\underline M_P$.
\item $\max_{t=1,\ldots, T}\Vert \mbf P_{0,t|t-1}\Vert\le M_P$ for some finite positive real $M_P$; 
\item $\min_{t=1,\ldots, T}\nu^{(r)}(\mbf P_{0,t|t-1})\ge \underline M_P$ for some finite positive real $\underline M_P$.
\end{compactenum}
\end{lem}

\noindent\textsc{Proof.} Given that $\mbf P_{0|0}=\mbf I_r$ is obviously positive definite, by \eqref{eq:pred2} and Weyl's inequality \citep[Theorem 1]{MK04} it follows that
\beq\label{eq:tossa}
\nu^{(r)}(\mbf P_{1|0})\ge \nu^{(r)}(\mbf A\mbf A^\prime) +\nu^{(r)}(\bm\Gamma^v)\ge M_v^{-1},
\eeq
for some finite positive real $M_v^{-1}$, since $\bm\Gamma^v$ has full rank by Assumption \ref{ass:common}(e) and $\nu^{(r)}(\mbf A\mbf A^\prime)$ is real and such that $\nu^{(r)}(\mbf A\mbf A^\prime)\ge (\nu^{(r)}(\mbf A))^2\ge 0$. From $\mbf P_{0|0}=\mbf I_r$  and \eqref{eq:pred2} it follows also that
\beq\label{eq:mamma}
\Vert \mbf P_{1|0}\Vert \le \Vert\mbf A\Vert^2 \Vert \mbf P_{0|0}\Vert +\Vert\bm\Gamma^v\Vert^2\le M_A^2+M_v^2=M_P,\;\text{ say.}
\eeq
By Lemma \ref{lem:detP}(ii) and \eqref{eq:mamma}, we have
\[
\max_{t=2,\ldots, T}\Vert\mbf P_{t|t-1}\Vert\le \Vert\mbf P_{1|0}\Vert\le M_P,
\]
since $M_P$ is independent of $t$ and $\mbf P_{t|t-1}$ is deterministic because of Lemma \ref{lem:detP}(i). This proves part (i).

For part (ii), by \citet[Theorems 1 and 7]{MK04} and \eqref{eq:pred2}
\beq\label{eq:tossa2}
\nu^{(r)}(\mbf P_{t|t-1})\ge \nu^{(r)}(\mbf A\mbf A^\prime)\nu^{(r)}(\mbf P_{t|t}) +\nu^{(r)}(\bm\Gamma^v)\ge  M_v^{-1},
\eeq
by the same arguments leading to \eqref{eq:tossa} and since $\nu_{\min}(\mbf P_{t|t})\ge 0$ because $\mbf P_{t|t}$ is at least positive semidefinite by construction. By letting $\underline M_P= M_v^{-1}$, we prove part (ii). 
Parts (iii) and (iv) are proved exactly as parts (i) and (ii), respectively, since also in this case $\mbf P_{0,0|0}=\mbf I_r$. 
This completes the proof.
$\Box$


\begin{lem}\label{lem:cazzarolahat0}
 Under Assumptions \ref{ass:common}, \ref{ass:idio}, \ref{ass:ind}, and \ref{ass:ident}, as $n,T\to\infty$:
\begin{compactenum}[(i)] 
\item $\max_{t=1,\ldots, T}\Vert \mbf P^{(0)}_{t|t-1}\Vert=O_p(1)$; 
\item $\max_{t=1,\ldots, T}\Vert(\mbf P^{(0)}_{t|t-1})^{-1}\Vert =O_p(1)$.
\end{compactenum}
\end{lem}

\noindent
\textsc{Proof.} For part (i), 
\al{
\max_{t=1,\ldots, T}\Vert {\mbf P}^{(0)}_{t|t-1}\Vert&\le \max_{t=1,\ldots, T}\Vert{\mbf P}_{t|t-1}\Vert+ \max_{t=1,\ldots, T}\Vert{\mbf P}_{t|t-1}^{(0)}-{\mbf P}_{t|t-1}\Vert\nn\\
&= O_p(1) + O_p(\max(n^{-1},T^{-1/2})),\nn
}
by Lemma \ref{lem:cazzarola}(i) and since the second term on the rhs depends only on the estimation error of $\wh{\mbf A}^{(0)}$, $\wh{\bm\Gamma}^{v(0)}$, $n^{-1/2}\wh{\bm\Lambda}_n^{(0)}$, 
$n^{-1/2}\wh{\bm\Lambda}_n^{(0)\prime}(\wh{\bm\Sigma}_n^{\xi(0)})^{-1}$
and $n^{-1}(\wh{\bm\Lambda}_n^{(0)\prime}
(\wh{\bm\Sigma}_n^{\xi(0)})^{-1}
\wh{\bm\Lambda}_n^{(0)})^{-1}$, which are all bounded by Lemmas \ref{lem:est0LOAD}(ii), \ref{lem:est0VAR}, \ref{lem:est0_LAST}(ii), and \ref{lem:est0_LAST}(iv). This proves part (i).

Part (ii) is proved in the same way as part (i) but using Lemma  \ref{lem:cazzarola}(ii). This completes the proof. $\Box$


\begin{lem}\label{lem:searle}
For $m<n$ with $m$ independent of $n$ and given
\begin{compactenum}
\item[(a)] an $n\times n$ matrix $\bm A$ symmetric and positive definite with $\Vert\bm A\Vert \le M_A$; 
\item[(a)] an $m\times m$ matrix $\bm B$ symmetric and positive definite with $\Vert\bm B\Vert \le M_B$; 
\item[(c)] an $n\times m$ matrix $\bm U$ such that $\Vert n^{-1}\bm U^\prime \bm U\Vert\le M_U$ and $\text{rk}(\bm U)=m$;
\item[(d)] an $m\times n$ matrix $\bm V$ such that $\Vert n^{-1}\bm V \bm V^\prime\Vert\le M_V$ and $\text{rk}(\bm V)=m$;
\end{compactenum}
where, $M_A$, $M_B$, $M_U$, and $M_{V}$ are finite positive reals independent of $n$ and $m$, then the following holds
\begin{compactenum}
\item[(i)] $(\bm A+\bm U\bm B\bm V)^{-1}=\bm A^{-1}-\bm A^{-1}\bm U\bm B(\mbf I_m+\bm V\bm A^{-1}\bm U\bm B)^{-1}\bm V\bm A^{-1}$;
\item[(ii)] $(\bm A+\bm U\bm B\bm V)^{-1}=\bm A^{-1}-\bm A^{-1}\bm U\bm B\bm V(\mbf I_m+\bm V\bm A^{-1}\bm U\bm B\bm V)^{-1}\bm A^{-1}$.
\end{compactenum}
\end{lem}

\noindent
\textsc{Proof.} Both results are proved by \citet[eq. (24) and eq. (25), respectively]{HS81}. $\Box$


\begin{lem}\label{lem:pioveadir8}
For any $r\times r$ symmetric and positive definite matrix $\bm P$ with $\Vert \bm P\Vert\le M_P$ for some finite positive real $M_P$, under Assumptions \ref{ass:common} and \ref{ass:idio}
\al{
\bm P\bm\Lambda_n^\prime (\bm\Lambda_n\bm P \bm\Lambda_n^\prime+\bm\Sigma_n^\xi)^{-1}
\bm\Lambda_n = \bm P(({\bm\Lambda}_n^{\prime}({\bm\Sigma}_n^{\xi})^{-1}{\bm\Lambda}_n)^{-1}+\bm P)^{-1}.\nn
}
\end{lem}

\noindent
\textsc{Proof.} In Lemma \ref{lem:searle}(i) set $\bm A={\bm\Sigma}_n^{\xi}$, $\bm B=\mbf I_r$, $\bm U={\bm\Lambda}_n\bm P$, and $\bm V={\bm\Lambda}_n^{\prime}$. Then, by noticing that the assumptions of Lemma \ref{lem:searle} are satisfied because of Assumptions \ref{ass:common}(a) and \ref{ass:idio}(a),  it follows that:
\[
 (\bm\Lambda_n\bm P \bm\Lambda_n^\prime+\bm\Sigma_n^\xi)^{-1} = ({\bm\Sigma}_n^{\xi})^{-1}-({\bm\Sigma}_n^{\xi})^{-1}\bm\Lambda_n\bm P(\mbf I_r+ \bm\Lambda_n^\prime({\bm\Sigma}_n^{\xi})^{-1}\bm\Lambda_n\bm P)^{-1}\bm\Lambda_n^\prime({\bm\Sigma}_n^{\xi})^{-1}.
\]
Therefore,
\al{
\bm P\bm\Lambda_n^\prime (\bm\Lambda_n\bm P \bm\Lambda_n^\prime+\bm\Sigma_n^\xi)^{-1}
\bm\Lambda_n &=\bm P\bm\Lambda_n^\prime\{
({\bm\Sigma}_n^{\xi})^{-1}-({\bm\Sigma}_n^{\xi})^{-1}\bm\Lambda_n\bm P(\mbf I_r+ \bm\Lambda_n^\prime({\bm\Sigma}_n^{\xi})^{-1}\bm\Lambda_n\bm P)^{-1}\bm\Lambda_n^\prime({\bm\Sigma}_n^{\xi})^{-1}
\}\bm\Lambda_n\nn\\
&= \bm P
\{
\bm\Lambda_n^\prime({\bm\Sigma}_n^{\xi})^{-1}\bm\Lambda_n
-\bm\Lambda_n^\prime({\bm\Sigma}_n^{\xi})^{-1}\bm\Lambda_n\bm P 
(\mbf I_r+\bm\Lambda_n^\prime({\bm\Sigma}_n^{\xi})^{-1}\bm\Lambda_n\bm P)\bm\Lambda_n^\prime({\bm\Sigma}_n^{\xi})^{-1}\bm\Lambda_n
\}.\label{eq:searle1}
}
Now, in Lemma \ref{lem:searle}(ii) set $\bm A=(\bm\Lambda_n^\prime({\bm\Sigma}_n^{\xi})^{-1}\bm\Lambda_n)^{-1}$, $\bm B=\bm P$, $\bm U=\mbf I_r$, and $\bm V=\mbf I_r
$, and notice that the assumptions therein are satisfied because of Lemmas \ref{lem:LSL2}(iii) and \ref{lem:LSL2}(v).
Then, for the last line of \eqref{eq:searle1} we have
\al{
\bm\Lambda_n^\prime({\bm\Sigma}_n^{\xi})^{-1}\bm\Lambda_n
-\bm\Lambda_n^\prime({\bm\Sigma}_n^{\xi})^{-1}\bm\Lambda_n\bm P 
(\mbf I_r+\bm\Lambda_n^\prime({\bm\Sigma}_n^{\xi})^{-1}\bm\Lambda_n\bm P)\bm\Lambda_n^\prime({\bm\Sigma}_n^{\xi})^{-1}\bm\Lambda_n= ((\bm\Lambda_n^\prime({\bm\Sigma}_n^{\xi})^{-1}\bm\Lambda_n)^{-1}+\bm P)^{-1}. \label{eq:searle2}
}
By substituting \eqref{eq:searle2} into \eqref{eq:searle1} we complete the proof. $\Box$


\begin{lem}\label{lem:cazzarolahat00}
 Under Assumptions \ref{ass:common}, \ref{ass:idio}, \ref{ass:ind}, and \ref{ass:ident}, as $n,T\to\infty$,
$\max_{t=1,\ldots, T}n\Vert \mbf P^{(0)}_{t|t}\Vert=O_p(1)$.
\end{lem}

\noindent
\textsc{Proof.} From \eqref{eq:up2} by using Lemma \ref{lem:pioveadir8}, but with  $\wh{\bm\Lambda}_n^{(0)}$ and $\wh{\bm\Sigma}_n^{\xi(0)}$ in place of $\bm\Lambda_n$ and $\bm\Sigma_n^\xi$,
it holds that:
\al{
{\mbf P}_{t|t}^{(0)}=&\, {\mbf P}_{t|t-1}^{(0)}-{\mbf P}_{t|t-1}^{(0)}\wh{\bm\Lambda}_n^{(0)\prime}
(\wh{\bm\Lambda}_n^{(0)}{\mbf P}_{t|t-1}^{(0)}\wh{\bm\Lambda}_n^{(0)\prime}+\wh{\bm\Sigma}_n^{\xi(0)})^{-1}\wh{\bm\Lambda}_n^{(0)}{\mbf P}_{t|t-1}^{(0)}\nn\\
=&\,\l\{ \mbf I_r-{\mbf P}_{t|t-1}^{(0)}\wh{\bm\Lambda}_n^{(0)\prime}
(\wh{\bm\Lambda}_n^{(0)}{\mbf P}_{t|t-1}^{(0)}\wh{\bm\Lambda}_n^{(0)\prime}+\wh{\bm\Sigma}_n^{\xi(0)})^{-1}\wh{\bm\Lambda}_n^{(0)}\r\}{\mbf P}_{t|t-1}^{(0)}\nn\\
=&\,\l\{ \mbf I_r-\mbf P_{t|t-1}^{(0)}((\wh{\bm\Lambda}_n^{(0)\prime}(\wh{\bm\Sigma}_n^{\xi(0)})^{-1}\wh{\bm\Lambda}_n^{(0)})^{-1}+\mbf P_{t|t-1}^{(0)})^{-1}\r\}{\mbf P}_{t|t-1}^{(0)}\label{eq:searle3}.
}
Then, by setting in Lemma \ref{lem:taylorinv} $\bm K= {\mbf P}_{t|t-1}^{(0)}$ and $\bm H=(\wh{\bm\Lambda}_n^{(0)\prime}(\wh{\bm\Sigma}_n^{\xi(0)})^{-1}\wh{\bm\Lambda}_n^{(0)})^{-1}$, for the last line of \eqref{eq:searle3} we have
\al{
((\wh{\bm\Lambda}_n^{(0)\prime}(\wh{\bm\Sigma}_n^{\xi(0)})^{-1}\wh{\bm\Lambda}_n^{(0)})^{-1}+\mbf P_{t|t-1}^{(0)})^{-1} =&\,(\mbf P_{t|t-1}^{(0)})^{-1}\label{eq:searle4}\\
&-
((\wh{\bm\Lambda}_n^{(0)\prime}(\wh{\bm\Sigma}_n^{\xi(0)})^{-1}\wh{\bm\Lambda}_n^{(0)})^{-1}+\mbf P_{t|t-1}^{(0)})^{-1} (\wh{\bm\Lambda}_n^{(0)\prime}(\wh{\bm\Sigma}_n^{\xi(0)})^{-1}\wh{\bm\Lambda}_n^{(0)})^{-1} (\mbf P_{t|t-1}^{(0)})^{-1}.\nn
}
By substituting \eqref{eq:searle4} into \eqref{eq:searle3} we get
\begin{align}
{\mbf P}_{t|t}^{(0)}
=&\,{\mbf P}_{t|t-1}^{(0)} ((\wh{\bm\Lambda}_n^{(0)\prime}(\wh{\bm\Sigma}_n^{\xi(0)})^{-1}\wh{\bm\Lambda}_n^{(0)})^{-1}+ {\mbf P}_{t|t-1}^{(0)})^{-1} (\wh{\bm\Lambda}_n^{(0)\prime}(\wh{\bm\Sigma}_n^{\xi(0)})^{-1}\wh{\bm\Lambda}_n^{(0)})^{-1}.\label{eq:searle5}
\end{align}
Finally, by using again \eqref{eq:searle4} into \eqref{eq:searle5}
\begin{align}
{\mbf P}_{t|t}^{(0)}
=&\,{\mbf P}_{t|t-1}^{(0)}\l\{
({\mbf P}_{t|t-1}^{(0)})^{-1}-((\wh{\bm\Lambda}_n^{(0)\prime}(\wh{\bm\Sigma}_n^{\xi(0)})^{-1}\wh{\bm\Lambda}_n^{(0)})^{-1}+ {\mbf P}_{t|t-1}^{(0)})^{-1}
(\wh{\bm\Lambda}_n^{(0)\prime}(\wh{\bm\Sigma}_n^{\xi(0)})^{-1}\wh{\bm\Lambda}_n^{(0)})^{-1}({\mbf P}_{t|t-1}^{(0)})^{-1}
\r\}\nn\\
&\cdot(\wh{\bm\Lambda}_n^{(0)\prime}(\wh{\bm\Sigma}_n^{\xi(0)})^{-1}\wh{\bm\Lambda}_n^{(0)})^{-1}\nn\\
=&\,\l\{
\mbf I_r-{\mbf P}_{t|t-1}^{(0)}((\wh{\bm\Lambda}_n^{(0)\prime}(\wh{\bm\Sigma}_n^{\xi(0)})^{-1}\wh{\bm\Lambda}_n^{(0)})^{-1}+ {\mbf P}_{t|t-1}^{(0)})^{-1}
(\wh{\bm\Lambda}_n^{(0)\prime}(\wh{\bm\Sigma}_n^{\xi(0)})^{-1}\wh{\bm\Lambda}_n^{(0)})^{-1}({\mbf P}_{t|t-1}^{(0)})^{-1}
\r\}\nn\\
&\cdot(\wh{\bm\Lambda}_n^{(0)\prime}(\wh{\bm\Sigma}_n^{\xi(0)})^{-1}\wh{\bm\Lambda}_n^{(0)})^{-1}\nn\\
=&\,(\wh{\bm\Lambda}_n^{(0)\prime}(\wh{\bm\Sigma}_n^{\xi(0)})^{-1}\wh{\bm\Lambda}_n^{(0)})^{-1}\nn\\
&- {\mbf P}_{t|t-1}^{(0)}
((\wh{\bm\Lambda}_n^{(0)\prime}(\wh{\bm\Sigma}_n^{\xi(0)})^{-1}\wh{\bm\Lambda}_n^{(0)})^{-1}+ {\mbf P}_{t|t-1}^{(0)})^{-1}
(\wh{\bm\Lambda}_n^{(0)\prime}(\wh{\bm\Sigma}_n^{\xi(0)})^{-1}\wh{\bm\Lambda}_n^{(0)})^{-1}\nn\\
&\cdot({\mbf P}_{t|t-1}^{(0)})^{-1}(\wh{\bm\Lambda}_n^{(0)\prime}(\wh{\bm\Sigma}_n^{\xi(0)})^{-1}\wh{\bm\Lambda}_n^{(0)})^{-1}.
\label{PtThathat}
\end{align}
Notice that we could use Lemmas \ref{lem:taylorinv} and \ref{lem:pioveadir8} to derive \eqref{PtThathat}  since all inverses used
are well defined because of Lemmas \ref{lem:est0_LAST}(iii), \ref{lem:cazzarolahat0}(i),  and \ref{lem:cazzarolahat0}(ii) and \eqref{eq:SHATPD} in the proof of Lemma \ref{lem:est0}.

Therefore, from \eqref{PtThathat}
\al{
\max_{t=1,\ldots, T} n\Vert {\mbf P}_{t|t}^{(0)}\Vert \le &\, 
 n \Vert (\wh{\bm\Lambda}_n^{(0)\prime}(\wh{\bm\Sigma}_n^{\xi(0)})^{-1}\wh{\bm\Lambda}_n^{(0)})^{-1}\Vert\nn\\
&+ \max_{t=1,\ldots, T} \Vert {\mbf P}_{t|t-1}^{(0)}\Vert \, \max_{t=1,\ldots, T} \Vert ({\mbf P}_{t|t-1}^{(0)})^{-1}\Vert\,
n\Vert (\wh{\bm\Lambda}_n^{(0)\prime}(\wh{\bm\Sigma}_n^{\xi(0)})^{-1}\wh{\bm\Lambda}_n^{(0)})^{-1}\Vert^2\nn\\
&\cdot \Vert ((\wh{\bm\Lambda}_n^{(0)\prime}(\wh{\bm\Sigma}_n^{\xi(0)})^{-1}\wh{\bm\Lambda}_n^{(0)})^{-1}+ {\mbf P}_{t|t-1}^{(0)})^{-1}\Vert\nn\\
=&\, O_p(1) + O_p(n^{-1}),\nn
}
 because of Lemmas \ref{lem:est0_LAST}(iii), \ref{lem:cazzarolahat0}(i), and \ref{lem:cazzarolahat0}(ii), and  since, by \citet[Theorem 1]{MK04} which is Weyl's inequality,
 \al{
 \Vert &((\wh{\bm\Lambda}_n^{(0)\prime}(\wh{\bm\Sigma}_n^{\xi(0)})^{-1}\wh{\bm\Lambda}_n^{(0)})^{-1}+ {\mbf P}_{t|t-1}^{(0)})^{-1}\Vert =
 \l\{\nu^{(r)}( (\wh{\bm\Lambda}_n^{(0)\prime}(\wh{\bm\Sigma}_n^{\xi(0)})^{-1}\wh{\bm\Lambda}_n^{(0)})^{-1}+ {\mbf P}_{t|t-1}^{(0)})  \r\}^{-1}\nn\\
 \le &\,  \l\{\nu^{(r)}( (\wh{\bm\Lambda}_n^{(0)\prime}(\wh{\bm\Sigma}_n^{\xi(0)})^{-1}\wh{\bm\Lambda}_n^{(0)})^{-1})+ \nu^{(r)}( {\mbf P}_{t|t-1}^{(0)})  \r\}^{-1}\nn\\
 = &\,  \l\{
 \l[\nu^{(1)} (\wh{\bm\Lambda}_n^{(0)\prime}(\wh{\bm\Sigma}_n^{\xi(0)})^{-1}\wh{\bm\Lambda}_n^{(0)})\r]^{-1}
 + \nu^{(r)}( {\mbf P}_{t|t-1}^{(0)})  \r\}^{-1}\nn\\
  = &\,  \l\{
 \l[\nu^{(1)} (\wh{\bm\Lambda}_n^{(0)\prime}(\wh{\bm\Sigma}_n^{\xi(0)})^{-1}\wh{\bm\Lambda}_n^{(0)})\nu^{(r)}( {\mbf P}_{t|t-1}^{(0)})\r]^{-1}
 + 1  \r\}^{-1} \l\{\nu^{(r)}( {\mbf P}_{t|t-1}^{(0)})\r\}^{-1}\nn\\
 =&\, \l\{1-\l[\nu^{(1)} (\wh{\bm\Lambda}_n^{(0)\prime}(\wh{\bm\Sigma}_n^{\xi(0)})^{-1}\wh{\bm\Lambda}_n^{(0)})\nu^{(r)}( {\mbf P}_{t|t-1}^{(0)})\r]^{-1}\r\}\l\{\nu^{(r)}( {\mbf P}_{t|t-1}^{(0)})\r\}^{-1}+O_p(n^{-2})\nn\\
 =&\, O_p(1),\label{eq:bochum}
 }
again by Lemmas \ref{lem:est0_LAST}(iii), \ref{lem:cazzarolahat0}(i), and \ref{lem:cazzarolahat0}(ii). This completes the proof. $\Box$

\begin{lem}\label{lem:PPOn}
 Under Assumptions \ref{ass:common}, \ref{ass:idio}, \ref{ass:ind}, and \ref{ass:ident}, as $n,T\to\infty$, 
$\max_{t=1,\ldots, T}n\Vert \mbf P^{(0)}_{t|T}\Vert=O_p(1)$.
\end{lem}

\noindent
\textsc{Proof.} From \eqref{eq:KS2}, we get
\al{
\Vert \mbf P_{t|T}^{(0)}-\mbf P_{t|t}^{(0)}\Vert\le&\, \Vert \mbf P_{t|t}^{(0)}\Vert^2\, \Vert \wh{\mbf A}^{(0)}\Vert^2\,  
\Vert(\mbf P_{t+1|t}^{(0)})^{-1}\Vert^2
\{\Vert\mbf P_{t+1|T}^{(0)}\Vert+\Vert\mbf P_{t+1|t}^{(0)}\Vert\}.\label{eq:piva}
}
Start with $t=T-1$, then from \eqref{eq:piva},
\al{
\Vert \mbf P_{T-1|T}^{(0)}-\mbf P_{T-1|T-1}^{(0)}\Vert\le&\, \Vert \mbf P_{T-1|T-1}^{(0)}\Vert^2\, \Vert \wh{\mbf A}^{(0)}\Vert^2\,  
\Vert(\mbf P_{T|T-1}^{(0)})^{-1}\Vert^2
\{\Vert\mbf P_{T|T}^{(0)}\Vert+\Vert\mbf P_{T|T-1}^{(0)}\Vert\}\nn\\
 =&\, O_p(n^{-2}).\label{eq:piva2}
}
by Lemmas \ref{lem:cazzarolahat0}(i), \ref{lem:cazzarolahat0}(ii), and \ref{lem:cazzarolahat00},  and since $\Vert\wh{\mbf A}^{(0)}\Vert\le \Vert{\mbf A}\Vert+\Vert\wh{\mbf A}^{(0)}-\mbf A\Vert =O_p(1)$, by Assumption \ref{ass:common}(d) and Lemma \ref{lem:est0VAR}(i). From \eqref{eq:piva2} 
it follows that 
\al{
\Vert \mbf P^{(0)}_{T-1|T}\Vert\le \Vert \mbf P^{(0)}_{T-1|T-1}\Vert+\Vert \mbf P^{(0)}_{T-1|T}-\mbf P^{(0)}_{T-1|T-1}\Vert = O_p(n^{-1})+O_p(n^{-2}).\label{eq:piva3}
}
Thus, at $t=T-2$, from \eqref{eq:piva} and \eqref{eq:piva3}, 
\al{
\Vert \mbf P_{T-2|T}^{(0)}-\mbf P_{T-2|T-2}^{(0)}\Vert\le&\, \Vert \mbf P_{T-2|T-2}^{(0)}\Vert^2\, \Vert \wh{\mbf A}^{(0)}\Vert^2\,  
\Vert(\mbf P_{T-1|T-2}^{(0)})^{-1}\Vert^2
\{\Vert\mbf P_{T-1|T}^{(0)}\Vert+\Vert\mbf P_{T-1|T-2}^{(0)}\Vert\}\nn\\
 =&\, O_p(n^{-2}).\label{eq:piva4}
}
From \eqref{eq:piva4} it follows that 
\al{
\Vert \mbf P^{(0)}_{T-2|T}\Vert\le \Vert \mbf P^{(0)}_{T-2|T-2}\Vert+\Vert \mbf P^{(0)}_{T-2|T}-\mbf P^{(0)}_{T-2|T-2}\Vert = O_p(n^{-1})+O_p(n^{-2}).\label{eq:piva5}
}
Since all the bounds in \eqref{eq:piva2}-\eqref{eq:piva5} are the same for all $t$, from Lemma \ref{lem:cazzarolahat00} and \eqref{eq:piva} we have
\al{
\max_{t=1,\ldots, T}n\Vert \mbf P^{(0)}_{t|T}\Vert\le \max_{t=1,\ldots, T}n\Vert \mbf P^{(0)}_{t|t}\Vert+\max_{t=1,\ldots, T}n\Vert \mbf P^{(0)}_{t|T}-\mbf P^{(0)}_{t|t}\Vert = O_p(1)+O_p(n^{-1}).\nn
}
This completes the proof. $\Box$

\begin{lem}\label{lem:wood}
For $m< n$, and given symmetric positive definite matrices $\bm A$ of dimension $m\times m$ and $\bm B$ of dimension $n\times n$, and for $\bm C$ of dimension $n\times m$ with $\text{rk}(\bm C)=m$, the following holds
\beq\label{statement}
\bm A \bm C^\prime (\bm C\bm A\bm C^\prime+\bm B)^{-1} = (\bm A^{-1}+\bm C^\prime\bm B^{-1}\bm C)^{-1}\bm C^\prime\bm B^{-1}.
\eeq
\end{lem}

\noindent\textsc{Proof.} 
Recall the Woodbury forumla
\begin{align}\label{woodbury}
(\bm C\bm A\bm C^\prime+\bm B)^{-1}=\bm B^{-1}-\bm B^{-1}\bm C(\bm A^{-1}+\bm C^\prime\bm B^{-1}\bm C)^{-1}\bm C^\prime\bm B^{-1}.
\end{align}
Denote $\bm D=(\bm A^{-1}+\bm C^\prime\bm B^{-1}\bm C)^{-1}$ then 
from \eqref{woodbury} the lhs of \eqref{statement} is equivalent to
\begin{align}
\bm A\bm C^\prime\l[\bm B^{-1}-\bm B^{-1}\bm C\bm D\bm C'\bm B^{-1}\r] = \bm A\l[\bm C'\bm B^{-1}-\bm C^\prime\bm B^{-1}\bm C\bm D\bm C^\prime\bm B^{-1}\r]=\bm A\l[\bm I-\bm C^\prime\bm B^{-1}\bm C\bm D\r]\bm C^\prime\bm B^{-1}.\nn
\end{align}
Then, \eqref{statement} becomes
\[
\bm A\l[\bm I-\bm C^\prime\bm B^{-1}\bm C\bm D\r]\bm C^\prime\bm B^{-1}=\bm D\bm C^\prime\bm B^{-1},
\]
or equivalently multiplying both sides on the right by $\bm B\bm C(\bm C^\prime\bm C)^{-1}$ 
\beq\label{statement2}
\bm A\l[\bm I-\bm C^\prime\bm B^{-1}\bm C\bm D\r]=\bm D.
\eeq
Now multiplying \eqref{statement2} on the left by $\bm A^{-1}$ and on the right by $\bm D^{-1}$
\[
\l[\bm D^{-1}-\bm C^\prime\bm B^{-1}\bm C\r]=\bm A^{-1},
\]
which is equivalent to
\[
\bm A^{-1}+\bm C^\prime\bm B^{-1}\bm C-\bm C^\prime\bm B^{-1}\bm C-\bm A^{-1}=\bm 0_{m\times m},
\]
which is always true. $\Box$

\begin{lem}\label{lem:FFO1}
 Under Assumptions \ref{ass:common}, \ref{ass:idio}, \ref{ass:ind}, and \ref{ass:ident}, as $n,T\to\infty$, for all $s=0,\ldots, T$,
$\Vert \mbf F^{(0)}_{t|s}\Vert=O_p(1)$, uniformly in $t\le s$.
\end{lem}

\noindent
\textsc{Proof.} 
Let
$\wh{\bm\Omega}^{F(0)}_s=\E_{\wh{\varphi}^{(0)}}[\bm F_s\bm F_s^\prime]$, which is $rs\times rs$ having the $r\times r$ generic $(t_1,t_2)$ block denoted by $[\wh{\bm\Omega}^{F(0)}_s]_{t_1,t_2}$ and such that $[\wh{\bm\Omega}^{F(0)}_s]_{t_2,t_1}=[\wh{\bm\Omega}^{F(0)}_s]_{t_1,t_2}^\prime$ and
\[
\text{vec}\l([\wh{\bm\Omega}^{F(0)}_s]_{t_1,t_2}\r) = 
(\mbf I_r\otimes \wh{\mbf A}^{(0)})^{|t_1-t_2|}
(\mbf I_{r^2}-\{\wh{\mbf A}^{(0)}\otimes\wh{\mbf A}^{(0)}\})^{-1}\text{vec}(\wh{\bm\Gamma}^{v(0)}), \quad t_1,t_2= 1,\ldots, s.
\] 
Notice that although $\wh{\bm\Omega}^{F(0)}_s$ depends on $\wh{\mbf A}^{(0)}$ and $\wh{\bm\Gamma}^{v(0)}$, for simplicity of notation, hereafter, we omit such dependence. 
Let also ${\bm\Omega}^{F}_s=\E[\bm F_s\bm F_s^\prime]$, clearly ${\bm\Omega}^{F}_s$ is positive definite, since by Assumptions \ref{ass:common}(b) and \ref{ass:common}(d), 
$\Vert 
[{\bm\Omega}^{F}_s]_{t_1,t_2}
\Vert \le \Vert \bm\Gamma^F\Vert$ for all $t_1\ne t_2$. Moreover, recall that ${\bm\Omega}^{F}_s$ is a block-Toeplitz matrix and define the corresponding circulant matrix as ${\bm\Phi}^{F}_s$, then (see, e.g., \citealp[Lemma 4.3 and Section 3.1]{gray2006toeplitz}).
$s^{-1} \vert \nu^{(1)}({\bm\Omega}^{F}_s) -\nu^{(1)}({\bm\Psi}^{F}_s)\vert = O(s^{-1/2})$,
and
$\nu^{(1)}({\bm\Psi}^{F}_s)  = O(s)$.
Thus, 
$$
\Vert {\bm\Omega}^{F}_s\Vert= \nu^{(1)}({\bm\Omega}^{F}_s) \le \nu^{(1)}({\bm\Psi}^{F}_s) + \vert \nu^{(1)}({\bm\Omega}^{F}_s) -\nu^{(1)}({\bm\Psi}^{F}_s)\vert = O(s)+O(\sqrt s),
$$ 
which implies  
\beq
s\Vert ({\bm\Omega}^{F}_s)^{-1}\Vert = O(1).\label{eq:puffinrockinv}
\eeq
Now the $r\times r$ generic $(t_1,t_2)$ block of $({\bm\Omega}^{F})_s^{-1}$ is an analytic function of $\mbf A$ and $\bm\Gamma^v$, which, in the case $r=1$, is given by (see \citealp{akaike73} for the case $r>1$)
\[
[({\bm\Omega}^{F}_s)^{-1}]_{t_1,t_2} = \E[v_t^2](1-A^2)^{-2} \cdot 
\l\{
\ba{ccccc}
1 & \text{if}& t_1=t_2=1 & \text{and}& t_1=t_2=s,\\ 
1+A^2& \text{if}& t_1=t_2 & \text{and}& 1<t_1,t_2<s,\\
-A &\text{if}& \vert t_1-t_2\vert=1 ,\\
0 &\text{otherwise}.
\ea
\r.
\]
Then, because of Lemma \ref{lem:est0VAR} and \eqref{eq:puffinrockinv}, we have that, as $n,s\to\infty$, $s(\wh{\bm\Omega}^{F(0)}_s)^{-1}$ is positive definite with probability tending to one, i.e.,
\al{
s \Vert (\wh{\bm\Omega}^{F(0)}_s)^{-1}\Vert = O_p (1).\label{eq:puffinrock}
}

Let now
\al{
\wh{\bm K}_{ns}^{(0)} &= \wh{\bm\Omega}^{F(0)}_s (\mbf I_s\otimes \wh{\bm\Lambda}_n^{(0)\prime}) \l\{(\mbf I_s\otimes \wh{\bm\Lambda}_n^{(0)})\wh{\bm\Omega}^{F(0)}_s
(\mbf I_s\otimes \wh{\bm\Lambda}_n^{(0)\prime})+(\mbf I_s\otimes \wh{\bm\Sigma}_n^{\xi(0)})
\r\}^{-1}\nn\\
&=\l\{(\mbf I_s\otimes \wh{\bm\Lambda}_n^{(0)\prime})
(\mbf I_s\otimes \wh{\bm\Sigma}_n^{\xi(0)})^{-1}
(\mbf I_s\otimes \wh{\bm\Lambda}_n^{(0)}) + (\wh{\bm\Omega}^{F(0)}_s)^{-1} \r\}^{-1}
\l\{(\mbf I_s\otimes \wh{\bm\Lambda}_n^{(0)\prime})
(\mbf I_s\otimes \wh{\bm\Sigma}_n^{\xi(0)})^{-1}\r\}\nn\\
&=\l\{\mbf I_s\otimes \wh{\bm\Lambda}_n^{(0)\prime}(\wh{\bm\Sigma}_n^{\xi(0)})^{-1}\wh{\bm\Lambda}_n^{(0)} + (\wh{\bm\Omega}^{F(0)}_s)^{-1} \r\}^{-1}
\l\{\mbf I_s\otimes \wh{\bm\Lambda}_n^{(0)\prime}(\wh{\bm\Sigma}_n^{\xi(0)})^{-1}\r\},\label{eq:andnow?}
}
which is $rs\times ns$ and where in the second line we used Lemma \ref{lem:wood}. Notice that all the inverses in \eqref{eq:andnow?} are well defined by \eqref{eq:puffinrock},
Lemma \ref{lem:est0_LAST}(iii),  and \eqref{eq:SHATPD} in the proof of Lemma \ref{lem:est0}.

Then, by definition of linear projection, we have
\al{
\mbf F_{t|s}^{(0)} =&\, \mathrm {Proj}_{\wh{\varphi}_{n}^{(0)}}[{\mbf F}_t|\bm X_{ns}]= (\bm\iota_{t,s}^\prime \otimes \mbf I_r)\{\wh{\bm K}_{ns}^{(0)}\}\bm X_{ns}\nn\\
=&\, (\bm\iota_{t,s}^\prime \otimes \mbf I_r)\l\{\mbf I_s\otimes \wh{\bm\Lambda}_n^{(0)\prime}(\wh{\bm\Sigma}_n^{\xi(0)})^{-1}\wh{\bm\Lambda}_n^{(0)} \r\}^{-1} \l\{\mbf I_s\otimes \wh{\bm\Lambda}_n^{(0)\prime}(\wh{\bm\Sigma}_n^{\xi(0)})^{-1}\r\}\bm X_{ns}\label{eq:lindt}\\
&+(\bm\iota_{t,s}^\prime \otimes \mbf I_r) \l[
\{\wh{\bm K}_{ns}^{(0)}\}^{-1} - \l\{\mbf I_s\otimes (\wh{\bm\Lambda}_n^{(0)\prime}(\wh{\bm\Sigma}_n^{\xi(0)})^{-1}\wh{\bm\Lambda}_n^{(0)} )^{-1} (\mbf I_s\otimes \wh{\bm\Lambda}_n^{(0)\prime}(\wh{\bm\Sigma}_n^{\xi(0)})^{-1} \r\}
\r]   \bm X_{ns},\nn
}
where $\bm\iota_{t,s}^\prime$ is the $t$th row of $\mbf I_s$ and  $\bm X_{ns}=(\mbf x_{n1}^\prime\cdots\mbf x_{ns}^\prime)^\prime$ is an  $ns$-dimensional vector.

Then, 
by 
\eqref{eq:abcinv00} in the proof of Lemma \ref{lem:denom}, 
we have 
\al{
n^{3/2}s &\l\Vert\{\wh{\bm K}_{ns}^{(0)}\}^{-1} -
\l\{\mbf I_s\otimes (\wh{\bm\Lambda}_n^{(0)\prime}(\wh{\bm\Sigma}_n^{\xi(0)})^{-1}\wh{\bm\Lambda}_n^{(0)})^{-1}
(\wh{\bm\Lambda}_n^{(0)\prime}(\wh{\bm\Sigma}_n^{\xi(0)})^{-1})
 \r\}\r\Vert \nn\\
\le&\, n^{2} s\l\Vert
\l\{\mbf I_s\otimes (\wh{\bm\Lambda}_n^{(0)\prime}(\wh{\bm\Sigma}_n^{\xi(0)})^{-1}\wh{\bm\Lambda}_n^{(0)})^{-1}+(\wh{\bm\Omega}_s^{F(0)})^{-1}\r\}  -
\l\{\mbf I_s\otimes (\wh{\bm\Lambda}_n^{(0)\prime}(\wh{\bm\Sigma}_n^{\xi(0)})^{-1}\wh{\bm\Lambda}_n^{(0)})^{-1}\r\} 
 \r\Vert\nn\\
 &\cdot n^{-1/2} \l\Vert
 \mbf I_s\otimes (\wh{\bm\Lambda}_n^{(0)\prime}(\wh{\bm\Sigma}_n^{\xi(0)})^{-1})
 \r\Vert\nn\\
 =&\, n^{2} s\l\Vert
\l\{\mbf I_s\otimes (\wh{\bm\Lambda}_n^{(0)\prime}(\wh{\bm\Sigma}_n^{\xi(0)})^{-1}\wh{\bm\Lambda}_n^{(0)})^{-1}+(\wh{\bm\Omega}_s^{F(0)})^{-1}\r\} (\wh{\bm\Omega}_s^{F(0)})^{-1}
\l\{\mbf I_s\otimes (\wh{\bm\Lambda}_n^{(0)\prime}(\wh{\bm\Sigma}_n^{\xi(0)})^{-1}\wh{\bm\Lambda}_n^{(0)})^{-1}\r\} 
 \r\Vert\nn\\
 &\cdot n^{-1/2} \l\Vert
 \mbf I_s\otimes (\wh{\bm\Lambda}_n^{(0)\prime}(\wh{\bm\Sigma}_n^{\xi(0)})^{-1})
 \r\Vert= O_p(1). \label{eq:bomba}
}
Indeed, for the first term on the rhs of \eqref{eq:bomba}, by Lemmas \ref{lem:LSL2}(i), \ref{lem:LSL2}(v), and \ref{lem:est0_LAST}(i), and by \eqref{eq:puffinrock}, we have
\al{
n^{2} s& \l\Vert
\l\{\mbf I_s\otimes (\wh{\bm\Lambda}_n^{(0)\prime}(\wh{\bm\Sigma}_n^{\xi(0)})^{-1}\wh{\bm\Lambda}_n^{(0)})^{-1}+(\wh{\bm\Omega}_s^{F(0)})^{-1}\r\} (\wh{\bm\Omega}_s^{F(0)})^{-1}
\l\{\mbf I_s\otimes (\wh{\bm\Lambda}_n^{(0)\prime}(\wh{\bm\Sigma}_n^{\xi(0)})^{-1}\wh{\bm\Lambda}_n^{(0)})^{-1}\r\} 
 \r\Vert\nn\\
 \le &\, n \l\Vert
\l\{\mbf I_s\otimes (\wh{\bm\Lambda}_n^{(0)\prime}(\wh{\bm\Sigma}_n^{\xi(0)})^{-1}\wh{\bm\Lambda}_n^{(0)})^{-1}+(\wh{\bm\Omega}_s^{F(0)})^{-1}\r\} \r\Vert\, s \Vert (\wh{\bm\Omega}_s^{F(0)})^{-1}\Vert \, n\Vert \mbf I_s\otimes (\wh{\bm\Lambda}_n^{(0)\prime}(\wh{\bm\Sigma}_n^{\xi(0)})^{-1}\wh{\bm\Lambda}_n^{(0)})^{-1}\Vert \nn\\
=&\, O_p(1).\label{eq:bomba1}
}
Alternatively to bound the first term on the rhs of \eqref{eq:bomba1} we can use directly Lemma  \ref{lem:est0_LAST}(iii) and \ref{eq:puffinrock}.

And, for the second term on the rhs of \eqref{eq:bomba}, by Lemmas  \ref{lem:LSL2}(vii) and \ref{lem:est0_LAST}(ii), we have
\al{
n^{-1/2} \Vert
 \mbf I_s\otimes (\wh{\bm\Lambda}_n^{(0)\prime}(\wh{\bm\Sigma}_n^{\xi(0)})^{-1})
 \Vert = O_p(1).\label{eq:bomba2}
 }

Now, let
$
 \wh{\bm \Pi}_s^{(0)}=\l[
\{\wh{\bm K}_{ns}^{(0)}\}^{-1} - \l\{\mbf I_s\otimes (\wh{\bm\Lambda}_n^{(0)\prime}(\wh{\bm\Sigma}_n^{\xi(0)})^{-1}\wh{\bm\Lambda}_n^{(0)} )^{-1} (\mbf I_s\otimes \wh{\bm\Lambda}_n^{(0)\prime}(\wh{\bm\Sigma}_n^{\xi(0)})^{-1} \r\}
\r],$
and similarly define $\bm \Pi_s$ when using the true parameters. By Lemmas \ref{lem:est0LOAD}(ii), \ref{lem:est0VAR}, \ref{lem:est0_LAST}(ii), and \ref{lem:est0_LAST}(iv), and \eqref{eq:bomba}, we have
\al{
n^{3/2}s\Vert  \wh{\bm \Pi}_s^{(0)}-\bm \Pi_s\Vert = O_p(\max(n^{-1},T^{-1/2})).\label{eq:bicocca}
}
Moreover, $ \bm \Pi_s $ is $rs\times ns$ such that
\al{
(\bm\iota_{t,s}^\prime \otimes \mbf I_r) \bm \Pi_s\bm X_{ns} &=
(\bm\iota_{t,s}^\prime \otimes \mbf I_r)\l(
\ba{c}
 \sum_{\tau=1}^s \bm \pi_{1\tau} \mbf x_{n\tau}\\
\vdots\\
 \sum_{\tau=1}^s \bm \pi_{s\tau} \mbf x_{n\tau}
\ea
\r)= \l(\sum_{\tau=1}^s \bm \pi_{1\tau} \mbf x_{n\tau}\cdots   \sum_{\tau=1}^s \bm \pi_{s\tau} \mbf x_{n\tau}\r)\bm\iota_{t,s}\nn\\
&= \sum_{\tau=1}^s \bm \pi_{t\tau} \mbf x_{n\tau},\label{eq:blocchi}
}
with $\bm \pi_{t\tau}$, $t,\tau=1,\ldots, n$ being $r\times n$ sub-block of $\bm\Pi_s$ and since $\text{vec}(\bm A\bm B\bm C)=(\bm C^\prime \otimes \bm A)\text{vec}(\bm B)$. Furthermore, since clearly from \eqref{eq:puffinrockinv} 
\al{
&\max_{t=1,\ldots, s}\sqrt s\Vert (\bm\iota_{t,s}^\prime \otimes \mbf I_r) [({\bm\Omega}_s)^{-1}]\Vert= O(1),
&\max_{\tau=1,\ldots, s}\sqrt s\Vert  [({\bm\Omega}_s)^{-1}](\bm\iota_{\tau,s} \otimes \mbf I_r)\Vert= O(1),
}
 using the same reasoning as in \eqref{eq:bomba1} and \eqref{eq:bomba2}, but when considering the true parameters, we also have that \linebreak
$\max_{t,\tau=1,\ldots, s} \Vert \bm\pi_{t\tau}\Vert = O(n^{-3/2}s^{-1/2})$. Thus, from \eqref{eq:blocchi} we have
\al{
\Vert (\bm\iota_{t,s}^\prime \otimes \mbf I_r) \bm \Pi_s\bm X_{ns} \Vert \le \max_{t,\tau=1,\ldots, s} \Vert \bm \pi_{t\tau} \Vert\,  \l\Vert \sum_{\tau=1}^s \mbf x_{n\tau}\r\Vert = O_p(n^{-1}), \label{eq:dore'}
}
since, by Assumption \ref{ass:common}(d) and Lemma \ref{lem:Gxi}(iii)
\al{
\E\l[\l\Vert \sum_{\tau=1}^s \mbf x_{n\tau} \r\Vert^2\r] 
&=\sum_{i=1}^n \sum_{\tau_1=1}^s\sum_{\tau_2=1}^s \E[x_{i\tau_1}x_{i\tau_2}]=\sum_{i=1}^n \sum_{\tau_1=1}^s\sum_{\tau_2=1}^s \bm\lambda_i^\prime \E[\mbf F_{\tau_1}\mbf F_{\tau_2}]\bm\lambda_i +\E[\xi_{i\tau_1}\xi_{i\tau_2}]\nn\\
&\le \sum_{i=1}^n \sum_{\tau_1=1}^s\sum_{\tau_2=1}^s \E[x_{i\tau_1}x_{i\tau_2}]=\sum_{i=1}^n \sum_{\tau_1=1}^s\sum_{\tau_2=1}^s \vert \bm\lambda_i^\prime \mbf A^{\vert \tau_1-\tau_2\vert}\bm\lambda_i \vert +\vert \E[\xi_{i\tau_1}\xi_{i\tau_2}]\vert\nn \\
&\le nM_\lambda^2  \l\{\sum_{\tau_1=1}^s\sum_{\tau_2=1}^s \Vert \mbf A^{|\tau_1-\tau_2|}\Vert +\sum_{\tau_1=1}^s\sum_{\tau_2=1}^s \vert \E[\xi_{i\tau_1}\xi_{i\tau_2}]\vert \r\}\nn\\
&\le ns M_\lambda^2 \sum_{k=-(s-1)}^{s-1} \l(1-s^{-1}|k|\r) M_A^{k}  + ns M_\lambda^2  M_3\nn\\
&\le ns M_\lambda^2 2 \sum_{k=0}^{s-1} M_A^{k}  + ns M_\lambda^2  M_3\nn\\
&\le ns M_\lambda^2 2 (1-M_A)^{-1} + ns M_\lambda^2  M_3.\nn
}

Thus, from \eqref{eq:lindt}, \eqref{eq:bomba}, \eqref{eq:bicocca}, and \eqref{eq:dore'}, and since $\text{vec}(\bm A\bm B\bm C)=(\bm C^\prime \otimes \bm A)\text{vec}(\bm B)$ and $\Vert\bm\iota_{t,s}^\prime \otimes \mbf I_r\Vert =1$,
\al{
\Vert& \mbf F_{t|s}^{(0)} \Vert \le 
\Vert 
( \wh{\bm\Lambda}_n^{(0)\prime}(\wh{\bm\Sigma}_n^{\xi(0)})^{-1}\wh{\bm\Lambda}_n^{(0)}  )^{-1}
  \wh{\bm\Lambda}_n^{(0)\prime}(\wh{\bm\Sigma}_n^{\xi(0)})^{-1} \mbf x_t
\Vert
\nn\\
&+ 
\l\Vert (\bm\iota_{t,s}^\prime \otimes \mbf I_r) \bm \Pi_s\bm X_{ns} \r\Vert + \Vert  \bm\iota_{t,s}^\prime \otimes \mbf I_r\Vert \, \Vert  \wh{\bm \Pi}_s^{(0)}-\bm \Pi_s\Vert 
\,\Vert \bm X_{ns}\Vert\nn\\
\le&\, n \Vert 
( \wh{\bm\Lambda}_n^{(0)\prime}(\wh{\bm\Sigma}_n^{\xi(0)})^{-1}\wh{\bm\Lambda}_n^{(0)}  )^{-1}\Vert \, n^{-1/2}\Vert \wh{\bm\Lambda}_n^{(0)\prime}(\wh{\bm\Sigma}_n^{\xi(0)})^{-1}\Vert\, n^{-1/2}\Vert \mbf x_t\Vert + O_p(n^{-1})+o_p(n^{-1}s^{-1/2})\nn\\
=&\, O_p(1),\label{eq:ponti}
}
where we used also Lemmas 
\ref{lem:LSL2}(v), \ref{lem:LSL2}(vii), \ref{lem:xunif}, \ref{lem:est0_LAST}(i), and \ref{lem:est0_LAST}(ii), and the facts that  $\Vert \bm X_{ns}\Vert=O(\sqrt{ns})$ by the same reasoning as in Lemma \ref{lem:xunif}, and $ \Vert  \bm\iota_{t,s}^\prime \otimes \mbf I_r\Vert=1$.
This completes the proof. $\Box$

\begin{lem}\label{lem:FFOn}
 Under Assumptions \ref{ass:common}, \ref{ass:idio}, \ref{ass:ind}, and \ref{ass:ident}, as $n,T\to\infty$,
$n\Vert \mbf F^{(0)}_{t|T}-\mbf F^{(0)}_{t|t}\Vert=O_p(1)$, uniformly in $t$.
\end{lem}

\noindent
\textsc{Proof.}  From \eqref{eq:KS1} and \eqref{eq:pred1}
\begin{align}
\Vert{\mbf F}_{t|T}^{(0)}- {\mbf F}_{t|t}^{(0)}\Vert&\le \Vert{\mbf P}_{t|t}^{(0)}\Vert \, \Vert\wh{\mbf A}^{(0)}\Vert\, \Vert({\mbf P}_{t+1|t}^{(0)})^{-1}\Vert\, \{\Vert{\mbf F}_{t+1|T}^{(0)}\Vert+\Vert{\mbf F}_{t+1|t}^{(0)}\Vert\}\nn\\
&\le \Vert{\mbf P}_{t|t}^{(0)}\Vert \, \Vert\wh{\mbf A}^{(0)}\Vert\, \Vert({\mbf P}_{t+1|t}^{(0)})^{-1}\Vert\, \{\Vert{\mbf F}_{t+1|T}^{(0)}\Vert+\Vert\wh{\mbf A}^{(0)}\Vert\,\Vert{\mbf F}_{t|t}^{(0)}\Vert\}\nn\\
&= O_p(n^{-1}),\label{eq:ponti2}
\end{align}
by Lemmas 
\ref{lem:cazzarolahat0}(ii),
\ref{lem:cazzarolahat00}, and
\ref{lem:FFO1} (when $s=T$ and $s=t$), and since $\Vert\wh{\mbf A}^{(0)}\Vert\le \Vert{\mbf A}\Vert+\Vert\wh{\mbf A}^{(0)}-\mbf A\Vert =O_p(1)$, by Assumption \ref{ass:common}(d) and Lemma \ref{lem:est0VAR}(i). This completes the proof.

\begin{lem}\label{lem:gennaio22bis}
 Under Assumptions \ref{ass:common}, \ref{ass:idio}, \ref{ass:ind}, and \ref{ass:ident}, as $n,T\to\infty$,
$n \Vert {\mbf F}_{t|t}^{(0)}-{\mbf F}_{t}^{\text{\tiny \upshape {WLS}}(0)}\Vert=O_p(1)$, uniformly in  $t$,
where ${\mbf F}_{t}^{\text{\tiny \upshape {WLS}}(0)}=(\wh{\bm\Lambda}_n^{(0)\prime}(\wh{\bm\Sigma}_n^{\xi(0)})^{-1}\wh{\bm\Lambda}_n^{(0)})^{-1}\wh{\bm\Lambda}_n^{(0)\prime}(\wh{\bm\Sigma}_n^{\xi(0)})^{-1}\mbf x_{nt}.$
\end{lem}

\noindent
\noindent\textsc{Proof.} From \eqref{eq:up1} and \eqref{eq:pred1}, by Lemma \ref{lem:wood} 
\begin{align}
{\mbf F}_{t|t}^{(0)}&={\mbf F}_{t|t-1}^{(0)}+ {\mbf P}_{t|t-1}^{(0)}\wh{\bm\Lambda}_n^{(0)\prime}(\wh{\bm\Lambda}_n^{(0)}{\mbf P}_{t|t-1}^{(0)}\wh{\bm\Lambda}_n^{(0)\prime}+\wh{\bm\Sigma}_n^{\xi(0)})^{-1}(\mbf x_{nt}-\wh{\bm\Lambda}_n^{(0)}{\mbf F}_{t|t-1}^{(0)})\nn\\
=&\,  (\wh{\bm\Lambda}_n^{(0)\prime}(\wh{\bm\Sigma}_n^{\xi(0)})^{-1}\wh{\bm\Lambda}_n^{(0)}+({\mbf P}_{t|t-1}^{(0)})^{-1})^{-1}\wh{\bm\Lambda}_n^{(0)\prime}(\wh{\bm\Sigma}_n^{\xi(0)})^{-1}\mbf x_{nt}\nn\\
&+\l\{\mbf I_r-
(\wh{\bm\Lambda}_n^{(0)\prime}(\wh{\bm\Sigma}_n^{\xi(0)})^{-1}\wh{\bm\Lambda}_n^{(0)}+(\wh{\mbf P}_{t|t-1}^{(0)})^{-1})^{-1}\wh{\bm\Lambda}_n^{(0)\prime}(\wh{\bm\Sigma}_n^{\xi(0)})^{-1}\wh{\bm\Lambda}_n^{(0)}\r\}\wh{\mbf A}^{(0)}{\mbf F}_{t-1|t-1}^{(0)}\nn\\
=&\,  
(\wh{\bm\Lambda}_n^{(0)\prime}(\wh{\bm\Sigma}_n^{\xi(0)})^{-1}\wh{\bm\Lambda}_n^{(0)})^{-1}
\wh{\bm\Lambda}_n^{(0)\prime}(\wh{\bm\Sigma}_n^{\xi(0)})^{-1}\mbf x_{nt}\nn\\
&+\l\{(\wh{\bm\Lambda}_n^{(0)\prime}(\wh{\bm\Sigma}_n^{\xi(0)})^{-1}\wh{\bm\Lambda}_n^{(0)}+({\mbf P}_{t|t-1}^{(0)})^{-1})^{-1}
-
(\wh{\bm\Lambda}_n^{(0)\prime}(\wh{\bm\Sigma}_n^{\xi(0)})^{-1}\wh{\bm\Lambda}_n^{(0)})^{-1}
\r\}
\wh{\bm\Lambda}_n^{(0)\prime}(\wh{\bm\Sigma}_n^{\xi(0)})^{-1}\mbf x_{nt}\nn\\
&+\l\{\mbf I_r-
(\wh{\bm\Lambda}_n^{(0)\prime}(\wh{\bm\Sigma}_n^{\xi(0)})^{-1}\wh{\bm\Lambda}_n^{(0)}+(\wh{\mbf P}_{t|t-1}^{(0)})^{-1})^{-1}\wh{\bm\Lambda}_n^{(0)\prime}(\wh{\bm\Sigma}_n^{\xi(0)})^{-1}\wh{\bm\Lambda}_n^{(0)}\r\}\wh{\mbf A}^{(0)}{\mbf F}_{t-1|t-1}^{(0)}.
\label{eq:KFhathat}
\end{align}
Notice that the inverses in \eqref{eq:KFhathat} are all well defined by Lemmas \ref{lem:est0_LAST}(iii) and \ref{lem:cazzarolahat0}(ii) and \eqref{eq:SHATPD} in the proof of Lemma \ref{lem:est0}.

Now, by Lemma \ref{lem:denom}
\al{
\Vert & (\wh{\bm\Lambda}_n^{(0)\prime}(\wh{\bm\Sigma}_n^{\xi(0)})^{-1}\wh{\bm\Lambda}_n^{(0)}+({\mbf P}_{t|t-1}^{(0)})^{-1})^{-1}\wh{\bm\Lambda}_n^{(0)\prime}(\wh{\bm\Sigma}_n^{\xi(0)})^{-1}\wh{\bm\Lambda}_n^{(0)}-\mbf I_r\Vert
= O_p(n^{-1}).\label{eq:otite}
}
Furthermore, by Lemmas \ref{lem:LSL}(iii) and \ref{lem:est0_LAST}(iv)
\al{
\Vert (\wh{\bm\Lambda}_n^{(0)\prime}(\wh{\bm\Sigma}_n^{\xi(0)})^{-1}\wh{\bm\Lambda}_n^{(0)}+({\mbf P}_{t|t-1}^{(0)})^{-1})^{-1}
-
(\wh{\bm\Lambda}_n^{(0)\prime}(\wh{\bm\Sigma}_n^{\xi(0)})^{-1}\wh{\bm\Lambda}_n^{(0)})^{-1}\Vert =O(n^{-2}),\label{eq:ennemeno2}
}
and by Lemmas \ref{lem:LSL2}(vii) and \ref{lem:est0_LAST}(ii), 
\al{
\Vert \wh{\bm\Lambda}_n^{(0)\prime}(\wh{\bm\Sigma}_n^{\xi(0)})^{-1}\Vert = O_p(\sqrt n).\label{eq:ennemeno05}
}
Indeed, we can apply Lemmas \ref{lem:denom} and \ref{lem:LSL}(iii), since $\Vert({\mbf P}_{t|t-1}^{(0)})^{-1}\Vert=O_p(1)$ by Lemma \ref{lem:cazzarolahat0}(ii),
$\Vert(\wh{\bm\Sigma}_{n}^{\xi(0)})^{-1}\Vert=O_p(1)$
by \eqref{eq:SHATPD} in the proof of Lemma \ref{lem:est0}, and, by Lemmas \ref{lem:lambdasqrtn} and \ref{lem:est0LOAD}(ii)  we have
\al{
n^{-1}\Vert \wh{\bm\Lambda}_n^{(0)\prime}\wh{\bm \Lambda}_n^{(0)}-{\bm\Lambda}_n^\prime{\bm \Lambda}_n\Vert\le&\, 2n^{-1} \Vert{\bm \Lambda}_n^\prime(\wh{\bm \Lambda}_n^{(0)}-{\bm \Lambda}_n)\Vert + n^{-1}\Vert 
(\wh{\bm \Lambda}_n^{(0)}-{\bm \Lambda}_n)^\prime(\wh{\bm \Lambda}_n^{(0)}-{\bm \Lambda}_n)\Vert\nn\\
\le &\, 2n^{-1/2}\Vert  {\bm \Lambda}_n\Vert \, n^{-1/2} \Vert \wh{\bm \Lambda}_n^{(0)}-{\bm \Lambda}_n\Vert + n^{-1}\Vert \wh{\bm \Lambda}_n^{(0)}-{\bm \Lambda}_n\Vert^2\nn\\
=&\,O_p(\max(n^{-1},{T}^{-1/2})).\nn
}
which, by Weyl's inequality \citep[Theorem 1]{MK04}, implies 
\al{
n^{-1}\vert \nu^{(j)}(\wh{\bm \Lambda}_n^{(0)\prime}\wh{\bm \Lambda}_n^{(0)})- \nu^{(j)}({\bm \Lambda}_n^\prime{\bm \Lambda}_n)\vert \le&\, n^{-1}\Vert \wh{\bm\Lambda}_n^{(0)\prime}\wh{\bm \Lambda}_n^{(0)}-{\bm\Lambda}_n^\prime{\bm \Lambda}_n\Vert=O_p(\max(n^{-1},{T}^{-1/2})),\nn
}
and, therefore, for $j=1,\ldots, r,$,
$$
\underline C_{j}\le \text{p-}\lim\!\!\!\!\!\!\inf_{n,T\to\infty} n^{-1} \nu^{(j)}(\wh{\bm \Lambda}_n^{(0)\prime}\wh{\bm \Lambda}_n^{(0)}) \le
\text{p-}\lim\!\!\!\!\!\sup_{n,T\to\infty} n^{-1} \nu^{(j)}(\wh{\bm\Lambda}_n^{(0)\prime}\wh{\bm \Lambda}_n^{(0)}) \le \overline C_j.
$$

By using \eqref{eq:otite}, \eqref{eq:ennemeno2}, \eqref{eq:ennemeno05} into \eqref{eq:KFhathat}:
\al{
\Vert& {\mbf F}_{t|t}^{(0)}-(\wh{\bm\Lambda}_n^{(0)\prime}(\wh{\bm\Sigma}_n^{\xi(0)})^{-1}\wh{\bm\Lambda}_n^{(0)})^{-1}\wh{\bm\Lambda}_n^{(0)\prime}(\wh{\bm\Sigma}_n^{\xi(0)})^{-1}\mbf x_{nt}\Vert\nn\\
 \le &\, n
\Vert(\wh{\bm\Lambda}_n^{(0)\prime}(\wh{\bm\Sigma}_n^{\xi(0)})^{-1}\wh{\bm\Lambda}_n^{(0)}+({\mbf P}_{t|t-1}^{(0)})^{-1})^{-1}
-
(\wh{\bm\Lambda}_n^{(0)\prime}(\wh{\bm\Sigma}_n^{\xi(0)})^{-1}\wh{\bm\Lambda}_n^{(0)})^{-1}
\Vert\,n^{-1/2}\Vert
\wh{\bm\Lambda}_n^{(0)\prime}(\wh{\bm\Sigma}_n^{\xi(0)})^{-1}\Vert\,n^{-1/2}\Vert \mbf x_{nt}\Vert \nn\\
&+ \Vert \mbf I_r-
(\wh{\bm\Lambda}_n^{(0)\prime}(\wh{\bm\Sigma}_n^{\xi(0)})^{-1}\wh{\bm\Lambda}_n^{(0)}+(\wh{\mbf P}_{t|t-1}^{(0)})^{-1})^{-1}\wh{\bm\Lambda}_n^{(0)\prime}(\wh{\bm\Sigma}_n^{\xi(0)})^{-1}\wh{\bm\Lambda}_n^{(0)}\Vert
\,\Vert \wh{\mbf A}^{(0)}\Vert\, \Vert {\mbf F}_{t-1|t-1}^{(0)}\Vert\nn\\
=&\, O_p(n^{-1}),\label{eq:girandola}
}
by Lemmas \ref{lem:xunif} and \ref{lem:FFO1} (when $s=t-1$), and since $\Vert\wh{\mbf A}^{(0)}\Vert\le \Vert{\mbf A}\Vert+\Vert\wh{\mbf A}^{(0)}-\mbf A\Vert =O_p(1)$, by Assumption \ref{ass:common}(d) and Lemma \ref{lem:est0VAR}(i). This completes the proof. $\Box$

\begin{lem}\label{lem:aprile24}
 Under Assumptions \ref{ass:common}, \ref{ass:idio}, \ref{ass:ind}, and \ref{ass:ident}, as $n,T\to\infty$,
$\min(\sqrt n,\sqrt T) \Vert {\mbf F}_{t|T}^{(0)}-{\mbf F}_{t}\Vert=O_p(1)$, uniformly in  $t$.
\end{lem}

\noindent
\textsc{Proof.} From Lemmas \ref{lem:FFOn} and \ref{lem:gennaio22bis}
\al{
\Vert\mbf F_{t|T}^{(0)}-\mbf F_t\Vert&\le \Vert\mbf F_{t|T}^{(0)}-\mbf F_{t|t}^{(0)}\Vert+\Vert\mbf F_{t|t}^{(0)}-\wh{\mbf F}_t^{\text{\tiny WLS}(0)}\Vert + \Vert \wh{\mbf F}_t^{\text{\tiny WLS}(0)}-\mbf F_t\Vert\nn\\
&= \Vert\wh{\mbf F}_t^{\text{\tiny WLS}(0)}-\mbf F_t\Vert + O_p(n^{-1}).\label{eq:ultimopezzoamonte}
}
where ${\mbf F}_{t}^{\text{\tiny \upshape {WLS}}(0)}=(\wh{\bm\Lambda}_n^{(0)\prime}(\wh{\bm\Sigma}_n^{\xi(0)})^{-1}\wh{\bm\Lambda}_n^{(0)})^{-1}\wh{\bm\Lambda}_n^{(0)\prime}(\wh{\bm\Sigma}_n^{\xi(0)})^{-1}\mbf x_{nt}.$

Now,
\al{
\Vert \wh{\mbf F}_t^{\text{\tiny WLS}(0)}-\mbf F_t\Vert\le&\, \Vert(\wh{\bm\Lambda}_n^{(0)\prime}(\wh{\bm\Sigma}_n^{\xi(0)})^{-1}\wh{\bm\Lambda}_n^{(0)})^{-1}\wh{\bm\Lambda}_n^{(0)\prime}(\wh{\bm\Sigma}_n^{\xi(0)})^{-1}(\bm\Lambda_n-\wh{\bm\Lambda}_n^{(0)})\Vert\, \Vert\mbf F_t\Vert\nn\\
&+ \Vert(\wh{\bm\Lambda}_n^{(0)\prime}(\wh{\bm\Sigma}_n^{\xi(0)})^{-1}\wh{\bm\Lambda}_n^{(0)})^{-1}\wh{\bm\Lambda}_n^{(0)\prime}(\wh{\bm\Sigma}_n^{\xi(0)})^{-1}\bm\xi_{nt}\Vert
\nn\\
\le&\,
\Vert
(\bm\Lambda_n^{\prime}
({\bm\Sigma}_n^{\xi})^{-1}
{\bm\Lambda}_n)^{-1}
{\bm\Lambda}_n^\prime
({\bm\Sigma}_n^{\xi})^{-1}
(\bm\Lambda_n-\wh{\bm\Lambda}_n^{(0)})\Vert\, 
\Vert\mbf F_t\Vert\nn\\
&+\Vert
(\wh{\bm\Lambda}_n^{(0)\prime}(\wh{\bm\Sigma}_n^{\xi(0)})^{-1}\wh{\bm\Lambda}_n^{(0)})^{-1}\wh{\bm\Lambda}_n^{(0)\prime}(\wh{\bm\Sigma}_n^{\xi(0)})^{-1}-({\bm\Lambda}_n^{\prime}({\bm\Sigma}_n^{\xi})^{-1}{\bm\Lambda}_n)^{-1}{\bm\Lambda}_n^{\prime}({\bm\Sigma}_n^{\xi})^{-1}\Vert\nn\\
&\cdot\Vert \bm\Lambda_n-\wh{\bm\Lambda}_n^{(0)}\Vert\, \Vert\mbf F_t\Vert+\Vert ({\bm\Lambda}_n^{\prime}({\bm\Sigma}_n^{\xi})^{-1}{\bm\Lambda}_n)^{-1} {\bm\Lambda}_n^{\prime}({\bm\Sigma}_n^{\xi})^{-1}\bm\xi_{nt}\Vert\nn\\
&+\Vert
(\wh{\bm\Lambda}_n^{(0)\prime}(\wh{\bm\Sigma}_n^{\xi(0)})^{-1}\wh{\bm\Lambda}_n^{(0)})^{-1}\wh{\bm\Lambda}_n^{(0)\prime}(\wh{\bm\Sigma}_n^{\xi(0)})^{-1}\bm\xi_{nt}-({\bm\Lambda}_n^{\prime}({\bm\Sigma}_n^{\xi})^{-1}{\bm\Lambda}_n)^{-1}{\bm\Lambda}_n^{\prime}({\bm\Sigma}_n^{\xi})^{-1}\bm\xi_{nt}\Vert \nn\\
=&\, A+B+C+D, \;\text{ say.} \label{eq:ultimopezzo}
}
Let us consider each term in \eqref{eq:ultimopezzo}. First, consider term $A$ and let
$\bm\Lambda_n^{\text{\tiny OLS}}=(\sum_{t=1}^T \mbf x_{nt}\mbf F_t^\prime) (\sum_{t=1}^T \mbf F_t\mbf F_t^\prime)^{-1}$. Then, from \citet[Corollary 1]{MBPCAQML}
\al{
n^{-1/2} \Vert \wh{\bm\Lambda}_n^{(0)}-\bm\Lambda_n^{\text{\tiny OLS}}\Vert=O_p(\max(n^{-1},n^{-1/2}T^{-1/2})).\label{eq:PCAOLS}
}
Therefore, from \eqref{eq:PCAOLS}
\al{
A\le&\, \Vert
(\bm\Lambda_n^{\prime}
({\bm\Sigma}_n^{\xi})^{-1}
{\bm\Lambda}_n)^{-1}
{\bm\Lambda}_n^\prime
({\bm\Sigma}_n^{\xi})^{-1}
(\bm\Lambda_n-\bm\Lambda_n^{\text{\tiny OLS}})\Vert\, 
\Vert\mbf F_t\Vert\nn\\
&+n \Vert
(\bm\Lambda_n^{\prime}
({\bm\Sigma}_n^{\xi})^{-1}
{\bm\Lambda}_n)^{-1}\Vert\, n^{-1/2}\Vert{\bm\Lambda}_n\Vert\, 
\Vert({\bm\Sigma}_n^{\xi})^{-1}\Vert\,
n^{-1/2} \Vert \wh{\bm\Lambda}_n^{(0)}-\bm\Lambda_n^{\text{\tiny OLS}}\Vert
\, 
\Vert\mbf F_t\Vert\nn\\
=&\, \{A.1 + A.2\} \Vert\mbf F_t\Vert, \;\text{ say.} \label{eq:ultimopezzoAA}
}
Then, 
\al{
A.1 \le &\, n \Vert
(\bm\Lambda_n^{\prime}
({\bm\Sigma}_n^{\xi})^{-1}
{\bm\Lambda}_n)^{-1}\Vert\,n^{-1}\Vert {\bm\Lambda}_n^\prime
({\bm\Sigma}_n^{\xi})^{-1}
(\bm\Lambda_n-\bm\Lambda_n^{\text{\tiny OLS}})\Vert\nn\\
=&\, n \Vert
(\bm\Lambda_n^{\prime}
({\bm\Sigma}_n^{\xi})^{-1}
{\bm\Lambda}_n)^{-1}\Vert\,n^{-1}\l\Vert 
T^{-1}\sum_{t=1}^T {\bm\Lambda}_n^\prime
({\bm\Sigma}_n^{\xi})^{-1}\bm\xi_{nt}\mbf F_t^\prime\r\Vert\,
\l\Vert \l(T^{-1}\sum_{t=1}^T \mbf F_t\mbf F_t^\prime\r)^{-1}
\r\Vert\nn\\
=&\, O_p(n^{-1/2}T^{-1/2}),\label{eq:ultimopezzoAA1}
}
by Lemmas \ref{lem:LSL2}(iii) and \ref{lem:COVFF0}(iv),
and also, recalling that $\bm\Gamma^F=\mbf I_r$ by Assumption \ref{ass:ident}(b), by Lemma \ref{lem:consistCOV}(i) and Weyl's inequality \citep[Theorem 1]{MK04} we have $\vert \nu^{(r)}(T^{-1}\sum_{t=1}^T \mbf F_t\mbf F_t^\prime)-1\vert = O_p(T^{-1/2})$ which implies
$\Vert (T^{-1}\sum_{t=1}^T \mbf F_t\mbf F_t^\prime)^{-1}
\Vert = O_p(1).$

Moreover, $A.2 = O_p(\max(n^{-1},n^{-1/2}T^{-1/2}))$,
because of \eqref{eq:PCAOLS} and Lemmas  \ref{lem:lambdasqrtn}, \ref{lem:LSL2}(iii), and Assumption \ref{ass:idio}(a) which implies $\Vert({\bm\Sigma}_n^{\xi})^{-1}\Vert\le C_\xi$. This, jointly with \eqref{eq:ultimopezzoAA} and \eqref{eq:ultimopezzoAA1} implies that
\al{
A=&\,O_p(\max(n^{-1},n^{-1/2}T^{-1/2})),  \label{eq:ultimopezzoAA0}
}
since $\Vert \mbf F_t\Vert = O_p(1)$ because $\E[F_{jt}^2]=1$, $j=1,\ldots, r$, by Assumption \ref{ass:ident}(b).

Second, by Lemmas \ref{lem:est0LOAD}(ii) and \ref{lem:est0_LAST}(v),
\al{
B=&\,  \Vert
n(\wh{\bm\Lambda}_n^{(0)\prime}(\wh{\bm\Sigma}_n^{\xi(0)})^{-1}\wh{\bm\Lambda}_n^{(0)})^{-1}n^{-1/2}\wh{\bm\Lambda}_n^{(0)\prime}(\wh{\bm\Sigma}_n^{\xi(0)})^{-1}-n({\bm\Lambda}_n^{\prime}({\bm\Sigma}_n^{\xi})^{-1}{\bm\Lambda}_n)^{-1}n^{-1/2}{\bm\Lambda}_n^{\prime}({\bm\Sigma}_n^{\xi})^{-1}\Vert\nn\\
&\cdot n^{-1/2}\Vert \bm\Lambda_n-\wh{\bm\Lambda}_n^{(0)}\Vert\, \Vert\mbf F_t\Vert\nn\\
=&\, O_p(\max(n^{-2},T^{-1})),\label{eq:ultimopezzoBB}
}
and since $\Vert \mbf F_t\Vert = O_p(1)$ because $\E[F_{jt}^2]=1$, $j=1,\ldots, r$, by Assumption \ref{ass:ident}(b).
Third,
\al{
C\le&\, n\Vert ({\bm\Lambda}_n^{\prime}({\bm\Sigma}_n^{\xi})^{-1}{\bm\Lambda}_n)^{-1}\Vert\,  n^{-1}\Vert{\bm\Lambda}_n^{\prime}({\bm\Sigma}_n^{\xi})^{-1}\bm\xi_{nt}\Vert = O_p(n^{-1/2}),\label{eq:ultimopezzoCC}
}
by Lemmas \ref{lem:LSL2}(iii) and \ref{lem:LSXi}(i).
Fourth, and last, 
\al{
D\le&\,n \Vert 
(\wh{\bm\Lambda}_n^{(0)\prime}(\wh{\bm\Sigma}_n^{\xi(0)})^{-1}\wh{\bm\Lambda}_n^{(0)})^{-1}-
({\bm\Lambda}_n^{\prime}({\bm\Sigma}_n^{\xi})^{-1}{\bm\Lambda}_n)^{-1} \Vert\, n^{-1}\Vert{\bm\Lambda}_n^{\prime}({\bm\Sigma}_n^{\xi})^{-1}\bm\xi_{nt}
\Vert\nn\\
&+ n\Vert ({\bm\Lambda}_n^{\prime}({\bm\Sigma}_n^{\xi})^{-1}{\bm\Lambda}_n)^{-1}\Vert\,n^{-1/2}\Vert
\wh{\bm\Lambda}_n^{(0)\prime}(\wh{\bm\Sigma}_n^{\xi(0)})^{-1}-
{\bm\Lambda}_n^{\prime}({\bm\Sigma}_n^{\xi})^{-1}
\Vert\, n^{-1/2}\Vert\bm\xi_{nt}\Vert\nn\\
&+ n\Vert 
(\wh{\bm\Lambda}_n^{(0)\prime}(\wh{\bm\Sigma}_n^{\xi(0)})^{-1}\wh{\bm\Lambda}_n^{(0)})^{-1}-
({\bm\Lambda}_n^{\prime}({\bm\Sigma}_n^{\xi})^{-1}{\bm\Lambda}_n)^{-1} \Vert\, n^{-1/2}\Vert
\wh{\bm\Lambda}_n^{(0)\prime}(\wh{\bm\Sigma}_n^{\xi(0)})^{-1}-
{\bm\Lambda}_n^{\prime}({\bm\Sigma}_n^{\xi})^{-1}
\Vert\, n^{-1/2}\Vert\bm\xi_{nt}\Vert\nn\\
=&\, D.1+D.2+D.3, \;\text{say.}\nn
}
Then, $D.1= O_p(\max(n^{-3/2},n^{-1/2}T^{-1/2}))$, by Lemmas \ref{lem:LSXi}(i) and \ref{lem:est0_LAST}(iv),
while we have $D.2= O_p(\max(n^{-1}, T^{-1/2}))$, because of Lemmas \ref{lem:LSL2}(iii) and \ref{lem:est0_LAST}(ii) and since $\Vert\bm\xi_{nt}\Vert= O_p(\sqrt n)$ because $\E[\xi_{it}^2]=\sigma_i^2\le C_\xi$ by Assumption \ref{ass:idio}(a). Last $D.3= O_p(n^{-2},T^{-1})$ by Lemmas \ref{lem:est0_LAST}(ii) and \ref{lem:est0_LAST}(iv) and since $\Vert\bm\xi_{nt}\Vert= O_p(\sqrt n)$.
Therefore,
\beq
D=  O_p(\max(n^{-1}, T^{-1/2})). \label{eq:ultimopezzoDD}
\eeq
By substituting \eqref{eq:ultimopezzoAA0}, \eqref{eq:ultimopezzoBB}, \eqref{eq:ultimopezzoCC}, and \eqref{eq:ultimopezzoDD}, into \eqref{eq:ultimopezzo} and then into \eqref{eq:ultimopezzoamonte}, we complete the proof. $\Box$

\begin{lem}\label{lem:FFO1sum}
 Under Assumptions \ref{ass:common}, \ref{ass:idio}, \ref{ass:ind}, and \ref{ass:ident}, as $n,T\to\infty$, for $s=t$ and $s=T$:
\begin{compactenum}
\item [(i)] $\Vert T^{-1}\sum_{t=1}^T  \mbf F^{(0)}_{t|s}\mbf F_t^\prime\Vert=O_p(1)$;
\item [(i)] $\Vert T^{-1}\sum_{t=1}^T  \mbf F^{(0)}_{t|s}\mbf F_t^\prime\bm\lambda_i\Vert=O_p(1)$, uniformly in $i$;
\item [(ii)] $\Vert T^{-1}\sum_{t=1}^T  \mbf F^{(0)}_{t|s}\xi_{it}\Vert=O_p(1)$, uniformly in $i$.
\end{compactenum}
\end{lem}

\noindent
\textsc{Proof.} 
First notice that, for all $k=t-T,\ldots, t-1$,
\al{
\l\Vert n^{-1/2} T^{-1}\sum_{t=1}^T  \mbf x_{n,t-k}\mbf F_t^\prime\r\Vert=O_p(1), \label{eq:campane}
}
by Lemma \ref{lem:consistCOV}. The proof of part (i) follows by iterating either forward or backwards since both
$\Vert T^{-1}\sum_{t=1}^T  \mbf F_{t|t}^{(0)} \mbf F_t^\prime\Vert$ and $\Vert  T^{-1}\sum_{t=1}^T  \mbf F_{t|T}^{(0)} \mbf F_t^\prime\Vert$
are functions of \eqref{eq:campane}, because of Lemmas \ref{lem:FFOn} and \ref{lem:gennaio22bis}.

Part (ii) follows from part (i) and Assumption \ref{ass:common}(a). Part (iii) follows by substituting $\mbf F_t$ with $\xi_{it}$ in \eqref{eq:campanestar} and then by applying Lemma \ref{lem:consistCOV}(ii). This completes the proof. $\Box$

\begin{lem}\label{lem:olsT12} 
Let ${\bm\phi}_n^{\text{\tiny \upshape OLS}}=(\text{\upshape vec}({\bm\Lambda}_n^{\text{\tiny\upshape OLS}})^\prime\; {\sigma}^{2\text{\tiny\upshape OLS}}_{1}\cdots \wh{\sigma}^{2\text{\tiny\upshape OLS}}_{n})^\prime$ be the vector of OLS estimators of the entries of ${\bm\phi}_n$ obtained when $\mbf F_t$ is known, that is whose entries maximize $\ell(\bm X_{nT}|\bm F_T;\underline{\bm\phi}_n)$. Let also ${\bm\theta}^{\text{\tiny \upshape OLS}}=(\text{\upshape vec}({\mbf A}^{\text{\tiny\upshape OLS}})^\prime\; \text{\upshape vech}(\bm\Gamma^{v\text{\tiny \upshape OLS}}))^\prime$ be the vector of OLS estimators of the entries of ${\bm\theta}$ obtained when $\mbf F_t$ is known, that is whose entries maximize $\ell(\bm F_T;\underline{\bm\theta})$. Then, under Assumptions \ref{ass:common}, \ref{ass:idio}, \ref{ass:ind}, and \ref{ass:ident}, as $n,T\to\infty$,
\begin{compactenum}
\item [(i)] 
$\sqrt T\, \Vert \bm\lambda_i^{\text{\tiny \upshape OLS}}-\bm\lambda_i\Vert=O_p(1)$, uniformly in $i$;
\item [(ii)] $\sqrt T\,
n^{-1/2} \Vert \bm\Lambda_n^{\text{\tiny \upshape OLS}}-\bm\Lambda_n\Vert= O_p(1)$;
\item [(iii)] $\sqrt T \vert\sigma_i^{2\text{\tiny \upshape OLS}}-\sigma_i^2\vert=O_p(1)$;
\item [(iv)] $\sqrt T \Vert\wh{\mbf A}^{\text{\tiny \upshape OLS}}-\mbf A \Vert=O_p(1)$;
\item [(v)]  $\sqrt T \Vert\wh{\bm\Gamma}^{v \text{\tiny \upshape OLS}}-\bm\Gamma^v \Vert=O_p(1)$.
\end{compactenum}
\end{lem}

\noindent
\textsc{Proof.} For part (i),
\al{
\Vert \bm\lambda_i^{\text{\tiny OLS}}-\bm\lambda_i\Vert\le \l\Vert \l(T^{-1}\sum_{t=1}^T \mbf F_t\mbf F_t^\prime\r)^{-1}\r\Vert \, \l\Vert T^{-1}\sum_{t=1}^T \mbf F_t x_{it} \r\Vert=O_p(T^{-1/2}),
}
where for the numerator we used Lemma 
 \ref{lem:consistCOV}(ii) while for the denominator we used Lemma \ref{lem:frida}.

Part (ii) is proved in the same way but using Lemma \ref{lem:consistCOV}(iii)  instead of Lemma \ref{lem:consistCOV}(ii).

For part (iii), first notice that
\al{
{\sigma}_i^{2\text{\tiny OLS}}=&\, T^{-1}\sum_{t=1}^T{\xi}_{it}^{\text{\tiny OLS}2}
=
T^{-1}\sum_{t=1}^T (x_{it}-\bm\lambda_i^{\text{\tiny OLS}\prime}\mbf F_t)^2\nn\\
=&\, T^{-1}\sum_{t=1}^T \{x_{it}^2 
+ \bm\lambda_i^{\prime}\mbf F_t\mbf F_t^\prime\bm\lambda_i
-2x_{it}\bm\lambda_i^\prime \mbf F_t\}\nn\\
&+T^{-1}\sum_{t=1}^T
\{ 2\bm\lambda_i^{\text{\tiny OLS}\prime}\mbf F_t\mbf F_t^\prime (\bm\lambda_i^{\text{\tiny OLS}}-\bm\lambda_i)
+(\bm\lambda_i^{\text{\tiny OLS}}-\bm\lambda_i)^\prime\ \mbf F_t\mbf F_t^{\prime}(\bm\lambda_i^{\text{\tiny OLS}}-\bm\lambda_i)
-2 x_{it}(\bm\lambda_i^{\text{\tiny OLS}}-\bm\lambda_i)^\prime\mbf F_t
\}\nn\\
=&\, T^{-1}\sum_{t=1}^T(x_{it}-\bm\lambda_i^\prime\mbf F_t)^2\nn\\
&+T^{-1}\sum_{t=1}^T
\{ 2\bm\lambda_i^{\text{\tiny OLS}\prime}\mbf F_t\mbf F_t^\prime (\bm\lambda_i^{\text{\tiny OLS}}-\bm\lambda_i)
+(\bm\lambda_i^{\text{\tiny OLS}}-\bm\lambda_i)^\prime\ \mbf F_t\mbf F_t^{\prime}(\bm\lambda_i^{\text{\tiny OLS}}-\bm\lambda_i)\nn\\
&-2 (\bm\lambda_i^{\text{\tiny OLS}}-\bm\lambda_i)^\prime\mbf F_t\mbf F_t^\prime\bm\lambda_i^{\text{\tiny OLS}}
-2 (\bm\lambda_i^{\text{\tiny OLS}}-\bm\lambda_i)^\prime\mbf F_t{\xi}_{it}^{\text{\tiny OLS}}
\}\nn\\
=&\, T^{-1}\sum_{t=1}^T\xi_{it}^2 + T^{-1}\sum_{t=1}^T (\bm\lambda_i^{\text{\tiny OLS}}-\bm\lambda_i)^\prime\ \mbf F_t\mbf F_t^{\prime}(\bm\lambda_i^{\text{\tiny OLS}}-\bm\lambda_i),\label{eq:homelate}
}
since by construction $\sum_{t=1}^T \mbf F_t {\xi}_{it}^{\text{\tiny OLS}}=\mbf 0_r$. Then,
\al{
\vert\sigma_i^{2\text{\tiny \upshape OLS}}-\sigma_i^2\vert \le &\, \l\vert
T^{-1}\sum_{t=1}^T (x_{it}-\bm\lambda_i^{\text{\tiny OLS}\prime}\mbf F_t)^2-
T^{-1}\sum_{t=1}^T (x_{it}-\bm\lambda_i^\prime\mbf F_t)^2
\r\vert
+
\l\vert
T^{-1}\sum_{t=1}^T (x_{it}-\bm\lambda_i^\prime\mbf F_t)^2-\E[\xi_{it}^2]
\r\vert\nn\\
\le &\, \Vert {\bm\lambda}_i^{\text{\tiny OLS}}-\bm\lambda_i\Vert^2\,  \l\Vert T^{-1}\sum_{t=1}^T  
\mbf F_t\mbf F_t^\prime\r\Vert
+ \l\vert
T^{-1}\sum_{t=1}^T\xi_{it}^2-\E[\xi_{it}^2]
\r\vert\nn\\
=&\,O_p(T^{-1/2}),\nn
}
because of part (i) and Lemma \ref{lem:consistCOV}(iv), and also by
Lemma \ref{lem:consistCOV}(i)  combined with Assumption \ref{ass:ident}(b).

For part (iv) consistency follows from the fact that $\{\mbf v_t\}$ is an independent process by Assumption \ref{ass:common}(f) and thus it is a martingale difference process \citep[Proposition 11, pp. 298-299]{Hamilton}. Part (v) follows from part (iv) and \citep[Proposition 11.2, pp. 301]{Hamilton}, since by Assumption \ref{ass:common}(g), the fourth order cumulants of $\{\mbf v_t\}$ are all finite. This completes the proof. $\Box$
%
%

\section{Lemmas necessary for proving Proposition \ref{prop:load}}


\begin{lem}\label{lem:tail}
Under Assumptions \ref{ass:common}, \ref{ass:idio}, \ref{ass:ind}, and \ref{ass:tails}
\begin{compactenum}
\item [(i)] For all $t\in\mathbb Z$, all $j=1,\ldots,r$, and all $s>0$, $\mathrm P(\vert F_{j t}\vert \ge s)\le \exp\{-K_F s^{\delta_v}\}$, for some finite positive real $K_F$ independent of $t$ and $j$;
\item [(ii)]
 for all $i\in\mathbb N$ and all $T\in\mathbb N$, 
$$
\mathrm P\l(\l\Vert T^{-1/2}\sum_{t=1}^T \mbf F_{t}\xi_{it}\r\Vert\ge s \r)\le 
r \exp\{- \kappa_{3} s^2\}+
rT \exp\{- \kappa_{4}(s\sqrt  T)^{\beta}\},
$$ 
for some finite positive reals $\kappa_3$, $\kappa_4$, and $\beta$ independent of $i$ and $T$ and such that $\frac 1\beta=\frac 1\gamma+\frac 1\delta>1$ and $\gamma=\min(\gamma_F,\gamma_\xi)$ and $\delta\in(0,\frac{\delta_v\delta_\xi}{\delta_v+\delta_\xi})$.
\end{compactenum} 
\end{lem}

\noindent\textsc{Proof.} 
For all $j=1,\ldots, r$, 
$$
\l \vert F_{jt}\r\vert\le \sum_{\ell=1}^r\l\vert \sum_{k=0}^\infty [\mbf A^k]_{j\ell} v_{\ell,t-k}\r\vert\le r\max_{\ell=1,\ldots,r}\l\vert \sum_{k=0}^\infty [\mbf A^k]_{j\ell} v_{\ell,t-k}\r\vert.
$$ 
Now, since because of Assumption \ref{ass:common}(d) the coefficients $[\mbf A^k]_{ j\ell}$ are summable over $k$, for any $\epsilon>0$ and $\eta>0$ there exists a positive integer $\bar K=\bar K(\epsilon, \eta)$ independent of $j$, $\ell$, and $t$ such that
\[
\mathrm P\l( \l\vert \sum_{k=\bar K+1}^\infty [\mbf A^k]_{j\ell} v_{\ell,t-k}\r\vert\ge \eta\r)\le \epsilon,
\]
thus we can always find $\bar K$ such that we can write
\begin{align}\label{eq:codeFF}
\l \vert F_{jt}\r\vert&\le r\max_{\ell=1,\ldots,r}\l\vert \sum_{k=0}^{\bar K} [\mbf A^k]_{j\ell} v_{\ell,t-k}\r\vert+\l\vert \sum_{k=\bar K+1}^{\infty} [\mbf A^k]_{j\ell} v_{\ell,t-k}\r\vert\nn\\
&\le r\bar K\max_{\ell=1,\ldots,r} \max_{k=0,\ldots, \bar K} \vert[\mbf A^k]_{j\ell}\vert\l\vert \sum_{k=0}^{\bar K} v_{\ell,t-k}\r\vert+o_p(1)\nn\\
&\le r\bar K C_A \l\vert \sum_{k=0}^{\bar K} v_{\ell,t-k}\r\vert+o_p(1),
\end{align}
for some positive real $C_A$ independent of $j$, where we used again Assumption \ref{ass:common}(d) to bound the coefficients.

Then,  from Assumption \ref{ass:tails}(a), \eqref{eq:codeFF}, and \citet[Corollary 4 and Section III]{maleki}, since $\{\mbf v_t\}$ is an independent process, and by Assumption \ref{ass:common} and the union bound,  for all $j=1,\ldots r$ and all $s>0$ it holds that:
\begin{align}
\mathrm P\l(\vert F_{j t}\vert \ge s\r)&\le \mathrm P\l(r\bar K C_A \l\vert \sum_{k=0}^{\bar K} v_{\ell,t-k}\r\vert\ge s\r)+o_p(1)\nn\\
&\le \exp\l\{-C^{\prime} \l(\frac{ s}{r \bar KC_A}\r)^2\r\}+\bar K \exp\l\{-C^{\prime\prime} \l(\frac{s}{r\bar K C_A} \r)^{\delta_v}\r\}\nn\\
&\le \bar K \exp\l\{-C^{\prime\prime\prime} s^{\delta_v}\r\}\le \exp\l\{-K_F s^{\delta_v}\r\},\label{eq:codeF}
\end{align}
for some finite positive reals $C^{\prime}$, $C^{\prime\prime}$, $C^{\prime\prime\prime}$, and $K_F$ independent of $t$ and $j$ and where from the second line we omitted the second term which is negligible. This proves part (i).

For all $i=1,\ldots, n$ and all $n\in\mathbb N$, consider Assumption \ref{ass:tails}(b) when $n=1$, $\bm\lambda_i=1$, and $\sigma_i^2=1$, then, because of Assumption \ref{ass:ind} and \eqref{eq:codeF}, by \citet[Lemma A2]{FLM11}, for all $j=1,\ldots r$ and all $s>0$ it holds that:
\beq\label{eq:codeFxi}
\mathrm P\l(\vert F_{j t}\xi_{it}\vert \ge s\r)\le\exp\l\{- \kappa_{F\xi} s^\delta \r\},
\eeq
for some finite positive reals $\kappa_{F\xi}$ and $\delta\in\l(0,\frac{\delta_v\delta_\xi}{\delta_v+\delta_\xi}\r)$ independent of $i$, $j$, and $t$. Moreover, because of Assumptions \ref{ass:common}(d), \ref{ass:common}(e), \ref{ass:common}(f), and \ref{ass:common}(h), by \citet[Theorem 3.1]{PT85}, $\{\mbf F_t\}$ is a strong mixing process with mixing coefficients 
\beq\label{eq:mixF}
\alpha_F(T)\le \exp\l\{-c_F T^{\gamma_F}\r\},
\eeq
for all $T\in\mathbb N$ and some finite positive reals $c_F$ and $\gamma_F$ independent of $T$ . Then, because of \eqref{eq:mixF} and Assumption \ref{ass:idio}(c), by \citet[Theorem 5.1.a]{bradley05}, we have that, for all $i\in\mathbb N$ and all $j=1,\ldots,r$, $\{F_{jt}\xi_{it}\}$ is strong mixing with mixing coefficients:
\beq\label{eq:mixFxi}
\alpha_{F\xi}(T)\le \alpha_F(T)+\alpha_\xi(T)\le \exp\l\{-c_F T^{\gamma_F}\r\}+\exp\l\{-c_\xi T^{\gamma_\xi}\r\}\le 2 \exp\l\{-c_{F\xi} T^{\gamma}\r\},
\eeq
for all $T\in\mathbb N$ and some finite positive reals $c_{F\xi}$ and $\gamma=\min(\gamma_F,\gamma_\xi)$ independent of $T$. 

Now, since for all $i\in\mathbb N$ and all $j=1,\ldots,r$, $\{F_{jt}\xi_{it}\}$ satisfies \eqref{eq:codeFxi} and \eqref{eq:mixFxi}, we can apply the results by \citet[Theorem 1]{MPR11} or equivalently by \citet[Theorem 1.4, p.31]{bosq12}, which imply that, for all $s>0$,
\beq\label{eq:rio}
\mathrm P\l(\l\vert T^{-1/2}\sum_{t=1}^T F_{jt}\xi_{it}\r\vert\ge s\r)\le 
\exp\bigg\{-c_1 s^2\bigg\}+T\exp \bigg\{-c_2 \Big(s\sqrt T\Big)^\beta  \bigg\}
\eeq
for some finite positive reals $c_1$, $c_2$, and $\beta$ independent of $i$, $j$, and $T$, where $\frac 1 \beta=\frac 1\gamma+\frac 1\delta>1$. Finally, note that
\beq\label{eq:riomax}
\l\Vert T^{-1/2}\sum_{t=1}^T \mbf F_{t}\xi_{it}\r\Vert\le \sum_{j=1}^r \l\vert T^{-1/2}\sum_{t=1}^T F_{jt}\xi_{it}\r\vert\le r\max_{j=1,\ldots, r} \l\vert T^{-1/2}\sum_{t=1}^T F_{jt}\xi_{it}\r\vert,
\eeq
and from \eqref{eq:rio} and \eqref{eq:riomax},
\begin{align}
\mathrm P\l(\l\Vert T^{-1/2}\sum_{t=1}^T \mbf F_{t}\xi_{it}\r\Vert\ge s\r)&\le \mathrm P\l(r\max_{j=1,\ldots, r}\l\vert T^{-1/2}\sum_{t=1}^T  F_{jt}\xi_{it}\r\vert\ge s\r)\nn\\
&\le  r\exp\bigg\{-c_1 \bigg(\frac s r\bigg)^2\bigg\}+rT\exp \bigg\{-c_2 \bigg(\frac{s\sqrt T}r\bigg)^\beta  \bigg\},\nn
\end{align}
by setting $\kappa_3=\frac{c_1}{r^2}$ and $\kappa_4=\frac{c_2}{r^\beta}$, we prove part (ii) and complete the proof. $\Box$
%
\begin{lem}\label{lem:funif24}
Under Assumptions \ref{ass:common} and \ref{ass:tails}, as $T\to\infty$, 
$$
\log^{-1/\delta_v} T \max_{t=1,\ldots,T} \l\Vert\mbf F_{t}\r\Vert= O_p(1).
$$
\end{lem}

\noindent
\noindent\textsc{Proof.} We have,
\beq
\max_{t=1,\ldots,T} \l\Vert\mbf F_{t}\r\Vert \le \max_{t=1,\ldots,T}\sum_{j=1}^r \vert F_{jt}\vert\le \max_{t=1,\ldots,T}\max_{j=1,\ldots, r} \vert F_{jt}\vert
=O_p(\log^{1/\delta_v} T),\nn
\eeq
by Lemma \ref{lem:tail}(i) and the union bound. This completes the proof. $\Box$
\begin{lem}\label{lem:xunif2}
Under Assumptions \ref{ass:common}, \ref{ass:idio}, \ref{ass:ind},  and \ref{ass:tails}, as $n,T\to\infty$, 
$$
\log^{-1/\delta_v} T\max_{t=1,\ldots,T} n^{-1/2} \Vert \mbf x_{nt}\Vert=O_p(1).
$$
\end{lem}

\noindent
\textsc{Proof.} 
By Assumption \ref{ass:tails}(b),  for all $s>0$, setting $\bm \lambda_i= \bm\iota_r$ and $\sigma_i^2=1$ therein, 
\beq
\mathrm P\l(n^{-1/2} \Vert {\bm \xi_{nt}}\Vert\ge s\r) \le \exp\l\{-\mathcal K_\xi s^2\r\} + n\exp\l\{- K_\xi (s\sqrt n)^{\delta_\xi}\r\}. \label{eq:caggia}
\eeq
Then, by \eqref{eq:caggia} and the union bound,  for all $s> 0$, it holds that:
\begin{align}
\mathrm P\l(\max_{t=1,\ldots,T} n^{-1/2}\Vert {\bm \xi_{nt}}\Vert\ge s\r)&\le T \mathrm P\l(n^{-1/2} \Vert {\bm \xi_{nt}}\Vert\ge s\r)\nn\\
&\le T \exp\l\{-\mathcal K_\xi s^2\r\} + nT\exp\l\{- K_\xi (s\sqrt n)^{\delta_\xi}\r\}.\label{eq:normaxiunif}
\end{align}
Thus, from Lemmas \ref{lem:lambdasqrtn} and \ref{lem:funif24}, and \eqref{eq:normaxiunif}
\begin{align}
\max_{t=1,\ldots,T}n^{-1/2}  \Vert {\mbf x_{nt}}\Vert&\le n^{-1/2}\Vert {\bm\Lambda_n}\Vert \max_{t=1,\ldots,T} \l\Vert\mbf F_{t}\r\Vert+\max_{t=1,\ldots,T}
n^{-1/2} \Vert  {\bm \xi_{nt}}\Vert\nn\\
&= O_p(\log^{1/\delta_v} T)+O_p(\sqrt{\log T})+O_p(n^{-\delta_\xi/2}\max(\log^{1/\delta_\xi} n,\log^{1/\delta_\xi} T))\nn\\
&= O_p(\log^{1/\delta_v} T).\nn
\end{align}
This completes the proof. $\Box$
\begin{lem}\label{lem:nuovo} 
Under Assumptions \ref{ass:common} and \ref{ass:idio}, for all $T\in\mathbb N$, as $n\to\infty$, 
\begin{compactenum}[(i)]
\item $\max_{t=1,\ldots, T}n \Vert\mbf P_{t|t}\Vert=O(1)$;
\item $\max_{t=1,\ldots, T}n \Vert\mbf P_{t|T}\Vert=O(1)$;
\item $\max_{t=1,\ldots, T}n \Vert\mbf P_{0,t|t}\Vert=O(1)$;
\item $\max_{t=1,\ldots, T}n \Vert\mbf P_{0,t|T}\Vert=O(1)$;
\item $\max_{t=1,\ldots, T}n^2 \Vert\mbf P_{t|T}-\mbf P_{t|t}\Vert=O(1)$;
\item $\max_{t=1,\ldots, T}n^2 \Vert\mbf P_{0,t|T}-\mbf P_{0,t|t}\Vert=O(1)$.
\end{compactenum}
\end{lem}

\noindent\textsc{Proof.} 
From \eqref{eq:up2}, using Lemma \ref{lem:wood} we have
\begin{align}
\mbf P_{t|t} &=\mbf P_{t|t-1}-\mbf P_{t|t-1}\bm\Lambda_n^\prime(\bm\Lambda_n\mbf P_{t|t-1}\bm\Lambda_n^\prime+\bm\Sigma_n^\xi)^{-1}\bm\Lambda_n\mbf P_{t|t-1}\nn\\
&=\mbf P_{t|t-1}-(\bm\Lambda_n^\prime(\bm\Sigma_n^\xi)^{-1}\bm\Lambda_n+\mbf P_{t|t-1}^{-1})^{-1}\bm\Lambda_n^\prime(\bm\Sigma_n^\xi)^{-1}\bm\Lambda_n\mbf P_{t|t-1},\label{eq:dottori}
\end{align}
indeed $\mbf P_{t|t-1}$ is positive definite and finite by Lemma \ref{lem:cazzarola}(i) and  \ref{lem:cazzarola}(ii). 
Therefore, 
since $\mbf P_{t|t-1}$ and $\mbf P_{t|t}$ are deterministic by Lemma \ref{lem:detP}(i), we have
\begin{align}
\max_{t=1,\ldots, T}n\Vert\mbf P_{t|t}\Vert &\le n\Vert \mbf P_{t|t-1}-\{(\bm\Lambda_n^\prime(\bm\Sigma_n^\xi)^{-1}\bm\Lambda_n+\mbf P_{t|t-1}^{-1})^{-1}\bm\Lambda_n^\prime(\bm\Sigma_n^\xi)^{-1}\bm\Lambda_n-\mbf I_r+\mbf I_r\}\mbf P_{t|t-1}\Vert\nn\\
 &\le n\Vert(\bm\Lambda_n^\prime(\bm\Sigma_n^\xi)^{-1}\bm\Lambda_n+\mbf P_{t|t-1}^{-1})^{-1}\bm\Lambda_n^\prime(\bm\Sigma_n^\xi)^{-1}\bm\Lambda_n-\mbf I_r\Vert\,\Vert\mbf P_{t|t-1}\Vert= O(1),
\end{align}
because of Lemma \ref{lem:LSL}(i), which is independent of $t$, and Lemma \ref{lem:cazzarola}(i). This proves part (i).

For part (ii), from \eqref{eq:KS2}, we have
\begin{align}
\Vert \mbf P_{t|T}-\mbf P_{t|t}\Vert &=\Vert \mbf P_{t|t} \mbf A^\prime \mbf P_{t+1|t}^{-1}
(\mbf P_{t+1|T}-\mbf P_{t+1|t})\mbf P_{t+1|t}^{-1} \mbf A \mbf P_{t|t}\Vert \nn\\
&\le \Vert \mbf P_{t|t} \Vert^2\, \Vert\mbf A\Vert^2\, \Vert \mbf P_{t+1|t}^{-1}\Vert^2\, \l\{\Vert \mbf P_{t+1|T}\Vert +\Vert\mbf P_{t+1|t}\Vert\r\}.\label{eq:smith}
\end{align}

Now, given that $\mbf P_{T|T}$ is obtained by the last iteration of the Kalman filter, for $t=T-1$ \eqref{eq:smith} becomes:
\begin{align}
\Vert \mbf P_{T-1|T}-\mbf P_{T-1|T-1}\Vert &\le \Vert \mbf P_{T-1|T-1} \Vert^2\, \Vert\mbf A\Vert^2\, \Vert \mbf P_{T|T-1}^{-1}\Vert^2\, \l\{\Vert \mbf P_{T|T}\Vert +\Vert\mbf P_{T|T-1}\Vert\r\}\nn\\
&\le \Vert \mbf P_{T-1|T-1} \Vert^2 \frac {M_A^2}{\underline M_P^2} \l\{\Vert \mbf P_{T|T}\Vert+M_P\r\}=O(n^{-2}),\label{eq:smith2}
\end{align}
by part (i), Assumption \ref{ass:common}(d), and Lemma \ref{lem:cazzarola}(i) and \ref{lem:cazzarola}(ii). Therefore, from part (i) and \eqref{eq:smith2}
\beq
\Vert\mbf P_{T-1|T}\Vert\le  \Vert \mbf P_{T-1|T}-\mbf P_{T-1|T-1}\Vert+ \Vert \mbf P_{T-1|T-1}\Vert = O(n^{-2})+ O(n^{-1}).\label{eq:smith4}
\eeq
For $t=T-2$, \eqref{eq:smith} becomes:
\begin{align}
\Vert \mbf P_{T-2|T}-\mbf P_{T-2|T-2}\Vert &=\Vert \mbf P_{T-2|T-2} \mbf A^\prime \mbf P_{T-1|T-2}^{-1}
(\mbf P_{T-1|T}-\mbf P_{T-1|T-2})\mbf P_{T-1|T-2}^{-1} \mbf A \mbf P_{T-2|T-2}\Vert \nn\\
&\le \Vert \mbf P_{T-2|T-2} \Vert^2\, \Vert\mbf A\Vert^2\, \Vert \mbf P_{T-1|T-2}^{-1}\Vert^2\, \l\{\Vert \mbf P_{T-1|T}\Vert +\Vert\mbf P_{T-1|T-2}\Vert\r\}\nn\\
&\le \Vert \mbf P_{T-2|T-2} \Vert^2 \frac { M_A^2}{\underline M_P^2} \l\{\Vert\mbf P_{T-1|T}\Vert+M_P\r\}=O(n^{-2}),\label{eq:smith3}
\end{align}
by part (i), Assumption \ref{ass:common}(d), Lemma \ref{lem:cazzarola}, and \eqref{eq:smith4}. Therefore, from part (i) and \eqref{eq:smith3}
\beq
\Vert\mbf P_{T-2|T}\Vert\le  \Vert \mbf P_{T-2|T}-\mbf P_{T-2|T-2}\Vert+ \Vert \mbf P_{T-2|T-2}\Vert=O(n^{-2})+O(n^{-1}).\label{eq:smith5}
\eeq
By comparing \eqref{eq:smith4} and \eqref{eq:smith5} it is clear that, the same asymptotic bound holds for all $t=T,\ldots, 1$
\beq\label{eq:smiths}
\Vert \mbf P_{t|T}-\mbf P_{t|t}\Vert \le \Vert \mbf P_{t|t} \Vert^2 \frac { M_A^2}{\underline M_P^2} \l\{\Vert\mbf P_{t+1|T}\Vert+M_P\r\}=O(n^{-2}),
\eeq
and since $\Vert\mbf P_{t+1|T}\Vert=O(n^{-1})$ then it is asymptotically negligible, thus \eqref{eq:smiths} holds uniformly in $t=T,\ldots, 1$ because of part (i), and since $\mbf P_{t|T}$ and $\mbf P_{t|t}$ are deterministic because of Lemma \ref{lem:detP}(i). It follows that, by part (i)
\[
\max_{t=1,\ldots,T}n \Vert\mbf P_{t|T}\Vert\le  \max_{t=1,\ldots,T}n\Vert \mbf P_{t|T}-\mbf P_{t|t}\Vert+\max_{t=1,\ldots,T} n\Vert \mbf P_{t|t}\Vert= O(n^{-1})+O(1)=O(1).
\]
This proves part (ii).

Parts (iii) and (iv) are proved exactly as parts (i) and (ii), respectively, but using Lemmas \ref{lem:LSL}(ii), \ref{lem:cazzarola}(iii), and \ref{lem:cazzarola}(iv) instead of Lemmas \ref{lem:LSL}(i), \ref{lem:cazzarola}(i), and \ref{lem:cazzarola}(ii). 

Part (v) is proved by \eqref{eq:smiths}. Part (vi) is proved as part (v) by repeating the same reasoning leading to the proof of part (iv). 
This completes the proof. $\Box$
\begin{lem}\label{lem:cazzarolaF}
 Under Assumptions \ref{ass:common}, \ref{ass:idio}, and \ref{ass:tails}, as $n,T\to\infty$,
\begin{compactenum}[(i)] 
\item $\Vert \mbf F_{t|t-1}\Vert=O_p(1)$, uniformly in $t$;
\item $\Vert \mbf F_{0,t|t-1}\Vert=O_p(1)$,  uniformly in $t$;
\item $\log^{-1/\delta_{v}} T\,\max_{t=1,\ldots, T}\Vert \mbf F_{0,t|t-1}\Vert=O_p(1)$.
\end{compactenum}
\end{lem}

\noindent\textsc{Proof.} Given that $\mbf F_{0|0}=\mbf 0_r$, from \eqref{eq:pred1} it follows that $\mbf F_{1|0}=\mbf 0_r$, then, from \eqref{eq:up1}
\[
\mbf F_{1|1}=\mbf P_{1|0}\bm\Lambda_n^\prime(\bm\Lambda_n\mbf P_{1|0}\bm\Lambda_n^\prime+\bm\Sigma_n^\xi)^{-1}\mbf x_{n1}=(\bm\Lambda_n^\prime(\bm\Sigma_n^\xi)^{-1}\bm\Lambda_n+\mbf P_{1|0}^{-1})^{-1}\bm\Lambda_n^\prime(\bm\Sigma_n^\xi)^{-1}\mbf x_{n1}
\]
by Lemma \ref{lem:wood} which can be applied since $\mbf P_{1|0}$ is  positive definite by Lemma \ref{lem:cazzarola}(ii). Thus, 
\begin{align}
\Vert\mbf F_{1|1}\Vert &\le \Vert (\bm\Lambda_n^\prime(\bm\Sigma_n^\xi)^{-1}\bm\Lambda_n+\mbf P_{1|0}^{-1})^{-1}\Vert \, \Vert \bm\Lambda_n\Vert\, \Vert(\bm\Sigma_n^\xi)^{-1}\Vert \, \Vert\mbf x_{n1}\Vert= O_p(1),\label{eq:pulce2}
\end{align}
by Lemmas \ref{lem:LSL2}(i), \ref{lem:lambdasqrtn}, and \ref{lem:xunif}, and Assumption \ref{ass:idio}(a). Then, from \eqref{eq:pred1}, 
\beq
\Vert\mbf F_{2|1}\Vert\le \Vert\mbf A\Vert\, \Vert\mbf F_{1|1}\Vert = O_p(1),\label{eq:pulce3}
\eeq
because of \eqref{eq:pulce2} and Assumption \ref{ass:common}(d). Therefore, from \eqref{eq:up1}, using \eqref{eq:pulce3} and the same arguments leading to \eqref{eq:pulce2}, \begin{align}
\Vert\mbf F_{2|2}\Vert &\le\Vert\mbf F_{2|1}\Vert + \Vert (\bm\Lambda_n^\prime(\bm\Sigma_n^\xi)^{-1}\bm\Lambda_n+\mbf P_{1|0}^{-1})^{-1}\Vert \, \Vert \bm\Lambda_n\Vert\, \Vert(\bm\Sigma_n^\xi)^{-1}\Vert \, \{\Vert\mbf x_{n1}\Vert+\Vert\bm\Lambda_n\Vert\,\Vert\mbf F_{2|1}\Vert \}\nn\\
&= O_p(1).\label{eq:pulce4}
\end{align}
It is then clear that \eqref{eq:pulce3} and \eqref{eq:pulce4} hold for all $t=1,\ldots, T$, and the result follows from Lemma \ref{lem:xunif}. This proves part (i).

For part (ii), the proof is identical to part (i) but using Lemma \ref{lem:LSL2}(ii) instead of Lemma \ref{lem:LSL2}(i). 
For part (iii) repeat the same steps as in part (ii) but using Lemma \ref{lem:xunif2} instead of Lemma \ref{lem:xunif}, and noticing that $\Vert \mbf x_{nt}\Vert\le \max_{t=1,\ldots, T}\Vert \mbf x_{nt}\Vert$, for all $t=1,\ldots, T$. This completes the proof. $\Box$

\begin{lem}\label{lem:KF} 
Under Assumptions \ref{ass:common}, \ref{ass:idio}, and \ref{ass:tails}, as $n,T\to\infty$, 
\begin{compactenum}[(i)]
\item $\sqrt n\Vert \mbf F_{t|t}-\mbf F_{t}\Vert=O_p(1)$, uniformly in $t$;
\item $\sqrt n\Vert \mbf F_{0,t|t}-\mbf F_{t}\Vert=O_p(1)$, uniformly in $t$;
\item $\log^{-1/\delta_{v}} T \sqrt n\,\max_{t=1,\ldots, T}\Vert \mbf F_{0,t|t}-\mbf F_t\Vert=O_p(1)$.
\end{compactenum}
\end{lem}

\noindent
\textsc{Proof.} Since $\mbf P_{t|t-1}$ is positive definite because of Lemma \ref{lem:cazzarola}(ii), we can use the Woodbury formula in Lemma \ref{lem:wood}
so that, for any $t=1,\ldots, T$, the Kalman filter estimator defined in \eqref{eq:up1} can be written as:
\begin{align}
\mbf F_{t|t}=&\, \mbf F_{t|t-1} + \mbf P_{t|t-1}\bm\Lambda_n^\prime(\bm\Lambda_n\mbf P_{t|t-1}\bm\Lambda_n^\prime+\bm\Sigma_n^\xi)^{-1}(\mbf x_{nt}-\bm\Lambda_n\mbf F_{t|t-1})\nn\\
=&\,\mbf F_{t|t-1} +(\bm\Lambda_n^\prime(\bm\Sigma_n^\xi)^{-1}\bm\Lambda_n+\mbf P_{t|t-1}^{-1})^{-1}\bm\Lambda_n^\prime(\bm\Sigma_n^\xi)^{-1}(\mbf x_{nt}-\bm\Lambda_n\mbf F_{t|t-1})\nn\\
=&\,\mbf F_{t|t-1} +(\bm\Lambda_n^\prime(\bm\Sigma_n^\xi)^{-1}\bm\Lambda_n+\mbf P_{t|t-1}^{-1})^{-1}\bm\Lambda_n^\prime(\bm\Sigma_n^\xi)^{-1}(\bm\Lambda_n\mbf F_t+\bm\xi_{nt}-\bm\Lambda_n\mbf F_{t|t-1})\nn\\
=&\,\mbf F_{t|t-1} +(\bm\Lambda_n^\prime(\bm\Sigma_n^\xi)^{-1}\bm\Lambda_n+\mbf P_{t|t-1}^{-1})^{-1}\bm\Lambda_n^\prime(\bm\Sigma_n^\xi)^{-1}\bm\Lambda_n\mbf F_t\nn\\
&-(\bm\Lambda_n^\prime(\bm\Sigma_n^\xi)^{-1}\bm\Lambda_n+\mbf P_{t|t-1}^{-1})^{-1}\bm\Lambda_n^\prime(\bm\Sigma_n^\xi)^{-1}\bm\Lambda_n\mbf F_{t|t-1}\nn\\
&+(\bm\Lambda_n^\prime(\bm\Sigma_n^\xi)^{-1}\bm\Lambda_n+\mbf P_{t|t-1}^{-1})^{-1}\bm\Lambda_n^\prime(\bm\Sigma_n^\xi)^{-1}\bm\xi_{nt},\label{eq:zotran}
\end{align}
where in the last step we used the definition of $\mbf x_{nt}$ in \eqref{eq:SDFM1R}.
Then, 
\al{
\Vert \mbf F_{t|t-1}&-(\bm\Lambda_n^\prime(\bm\Sigma_n^\xi)^{-1}\bm\Lambda_n+\mbf P_{t|t-1}^{-1})^{-1}\bm\Lambda_n^\prime(\bm\Sigma_n^\xi)^{-1}\bm\Lambda_n\mbf F_{t|t-1}\Vert\nn\\
&\le \Vert \mbf F_{t|t-1}\Vert \, \Vert\mbf I_r-(\bm\Lambda_n^\prime(\bm\Sigma_n^\xi)^{-1}\bm\Lambda_n+\mbf P_{t|t-1}^{-1})^{-1}\bm\Lambda_n^\prime(\bm\Sigma_n^\xi)^{-1}\bm\Lambda_n\Vert= O_p(n^{-1}),\label{eq:krapfen}
}
by Lemma \ref{lem:cazzarolaF}(i) and Lemma \ref{lem:LSL}(i), which can be applied by Lemma \ref{lem:cazzarola}(i). Similarly,
\al{
\Vert&(\bm\Lambda_n^\prime(\bm\Sigma_n^\xi)^{-1}\bm\Lambda_n+\mbf P_{t|t-1}^{-1})^{-1}\bm\Lambda_n^\prime(\bm\Sigma_n^\xi)^{-1}\bm\Lambda_n\mbf F_t-\mbf F_t\Vert\nn\\
&\le
\Vert(\bm\Lambda_n^\prime(\bm\Sigma_n^\xi)^{-1}\bm\Lambda_n+\mbf P_{t|t-1}^{-1})^{-1}\bm\Lambda_n^\prime(\bm\Sigma_n^\xi)^{-1}\bm\Lambda_n-\mbf I_r\Vert\, \Vert\mbf F_t\Vert= O_p(n^{-1}),\label{eq:krapfen2}
}
by the same arguments leading to \eqref{eq:krapfen} and since $\Vert\mbf F_t\Vert = O_p(1)$ uniformly in $t=1,\ldots, T$, because $\E[\Vert\mbf F_t\Vert^2]=\sum_{j=1}^r \E[F_{jt}^2]=\text{tr}(\bm\Gamma^F)=r$ by Assumption \ref{ass:ident}(b). From \eqref{eq:krapfen} and \eqref{eq:krapfen2} it follows that
\begin{align}
\Vert \mbf F_{t|t}- \mbf F_t\Vert
&\le\Vert (\bm\Lambda_n^\prime(\bm\Sigma_n^\xi)^{-1}\bm\Lambda_n)^{-1}\bm\Lambda_n^\prime(\bm\Sigma_n^\xi)^{-1}\bm\xi_{nt}\Vert+ O_p(n^{-1})\nn\\
&\le \Vert  n(\bm\Lambda_n^\prime(\bm\Sigma_n^\xi)^{-1}\bm\Lambda_n)^{-1}\Vert\, \Vert n^{-1}\bm\Lambda_n^\prime(\bm\Sigma_n^\xi)^{-1}\bm\xi_{nt}\Vert+ O_p(n^{-1})\nn\\
&=O_p(n^{-1/2}),\label{eq:c18}
\end{align}
by Lemma \ref{lem:LSL2}(iii) and Lemma \ref{lem:LSXi}(i). This proves part (i).

Part (ii) is proved in the same way, but using Lemmas \ref{lem:cazzarola}(iii), \ref{lem:cazzarola}(iv), \ref{lem:cazzarolaF}(ii), \ref{lem:LSL}(ii), \ref{lem:LSL2}(iv), and \ref{lem:LSXi}(ii), instead of Lemmas  \ref{lem:cazzarola}(i), \ref{lem:cazzarola}(ii), \ref{lem:cazzarolaF}(i), \ref{lem:LSL}(i), \ref{lem:LSL2}(iii), and \ref{lem:LSXi}(i). 
For part (iii) repeat the same steps as in part (ii) 
but using Lemma \ref{lem:cazzarolaF}(iii) instead of Lemma \ref{lem:cazzarolaF}(ii), and since $\max_{t=1,\ldots, T}\Vert \mbf F_t\Vert = O_p({\log^{1/\delta_v} T})$ by Lemma \ref{lem:funif24}. 
This completes the proof. $\Box$

\begin{lem}\label{lem:KSKF}
Under Assumptions \ref{ass:common} and \ref{ass:idio}, as $n\to\infty$, 
\begin{compactenum}[(i)]
\item $n\Vert \mbf F_{t|T}-\mbf F_{t|t}\Vert=O_p(1)$, uniformly in $t$;
\item $n\Vert \mbf F_{0,t|T}-\mbf F_{0,t|t}\Vert=O_p(1)$, uniformly in $t$.
\end{compactenum}
\end{lem}

\noindent
\noindent\textsc{Proof.} For any $t=1,\ldots,T$, using \eqref{eq:KS1}
\begin{align}
\Vert \mbf F_{t|T}-\mbf F_{t|t}\Vert\le \Vert\mbf P_{t|t}\Vert\, \Vert \mbf A^\prime\Vert\,\Vert\mbf P_{t+1|t}^{-1}\Vert\, \Vert(\mbf F_{t+1|T}-\mbf A\mbf F_{t|t})\Vert.\label{eq:pavia}
\end{align}
Now, given that $\mbf F_{T|T}$ is obtained by the last iteration of the Kalman filter, for $t=T-1$ \eqref{eq:pavia} becomes:
\begin{align}
\Vert& \mbf F_{T-1|T}-\mbf F_{T-1|T-1}\Vert \le \Vert \mbf P_{T-1|T-1} \Vert\, \Vert\mbf A\Vert \, \Vert \mbf P_{T|T-1}^{-1}\Vert\, \l\{\Vert \mbf F_{T|T}\Vert +\Vert \mbf A\Vert\,\Vert\mbf F_{T-1|T-1}\Vert\r\}\nn\\
&\le \Vert \mbf P_{T-1|T-1} \Vert\,\frac{M_A}{\underline M_P}\l\{
\Vert \mbf F_{T|T}-\mbf F_T\Vert+\Vert \mbf F_{T}\Vert+ M_A \l[\Vert \mbf F_{T-1|T-1}-\mbf F_{T-1}\Vert+ \Vert \mbf F_{T-1}\Vert\r]
\r\}\nn\\
&= O_p(n^{-1}),\label{eq:pavia2}
\end{align}
by Assumption \ref{ass:common}(d), and Lemmas \ref{lem:cazzarola}(ii), \ref{lem:nuovo}(i), and \ref{lem:KF}(i), and since $\Vert\mbf F_t\Vert = O_p(1)$ uniformly in $t=1,\ldots, T$, because $\E[\Vert\mbf F_t\Vert^2]=\sum_{j=1}^r \E[F_{jt}^2]=\text{tr}(\bm\Gamma^F)=r$ by Assumption \ref{ass:ident}(b). 
For $t=T-2$ \eqref{eq:pavia} becomes:
\al{
\Vert& \mbf F_{T-2|T}-\mbf F_{T-2|T-2}\Vert\le \Vert\mbf P_{T-2|T-2}\Vert\, \Vert \mbf A^\prime\Vert\,\Vert\mbf P_{T-1|T-2}^{-1}\Vert\, \Vert(\mbf F_{T-1|T}-\mbf A\mbf F_{T-2|T-2})\Vert\nn\\
\le&\, 
\Vert \mbf P_{T-1|T-1} \Vert\,\frac{M_A}{\underline M_P}
\l\{
\Vert\mbf F_{T-1|T}-\mbf F_{T-1|T-1}\Vert+\Vert \mbf F_{T-1|T-1}-\mbf F_{T-1}\Vert +\Vert \mbf F_{T-1}\Vert\r.\nn\\
&\l.+ M_A \l[\Vert \mbf F_{T-2|T-2}-\mbf F_{T-2}\Vert+ \Vert \mbf F_{T-2}\Vert\r]
\r\}\nn\\
=&\,O_p(n^{-1}),
\label{eq:pavia4}
}
by \eqref{eq:pavia2}, Assumption \ref{ass:common}(d), and Lemmas \ref{lem:cazzarola}(ii), \ref{lem:nuovo}(i), and \ref{lem:KF}(i), and since $\Vert\mbf F_t\Vert = O_p(1)$ uniformly in $t=1,\ldots, T$.  By comparing  \eqref{eq:pavia2} and \eqref{eq:pavia4} it is clear that, the same asymptotic bound holds uniformly in $t=T,\ldots, 1$, i.e., from \eqref{eq:pavia} we get
\al{
\Vert \mbf F_{t|T}-\mbf F_{t|t}\Vert\le&\, \Vert\mbf P_{t|t}\Vert\, \frac{M_A}{\underline M_P} \l\{\Vert\mbf F_{t+1|T}-\mbf F_{t+1|t+1}\Vert+\Vert\mbf F_{t+1|t+1}-\mbf F_{t+1}\Vert+\Vert\mbf F_{t+1}\Vert\r.\nn\\
&\l.+M_A\l[\Vert\mbf F_{t|t}-\mbf F_t\Vert+\Vert\mbf F_t\Vert\r]\r\}=O_p(n^{-1}).
\label{eq:E7iidet}
}
This proves part (i). 

Part (ii) is proved in the same way but using Lemmas \ref{lem:cazzarola}(iv), \ref{lem:nuovo}(iii), and \ref{lem:KF}(ii), instead of Lemmas \ref{lem:cazzarola}(ii), \ref{lem:nuovo}(i), and \ref{lem:KF}(i). 
This completes the proof. $\Box$


\begin{lem}\label{lem:KFGLS}
Under Assumptions \ref{ass:common}, \ref{ass:idio}, \ref{ass:ind}, and \ref{ass:tails}, as $n\to\infty$, 
\begin{compactenum}[(i)]
\item $n\Vert \mbf F_{t|t}-\mbf F_t^{\text{\tiny \upshape WLS}}\Vert = O_p(1)$, uniformly in $t$;
\item $n\Vert \mbf F_{0,t|t}-\mbf F_t^{\text{\tiny \upshape GLS}}\Vert = O_p(1)$, uniformly in $t$;
\item $n^{2}\max_{t=1,\ldots,T}\Vert \mbf P_{t|t}-(\bm\Lambda_n^\prime(\bm\Sigma_n^\xi)^{-1}\bm\Lambda_n)^{-1}\Vert = O_p(1)$;
\item $n^{2}\max_{t=1,\ldots,T}\Vert \mbf P_{0,t|t}-(\bm\Lambda_n^\prime(\bm\Gamma_n^\xi)^{-1}\bm\Lambda_n)^{-1}\Vert = O_p(1)$;
\item $n^{2}\max_{t=1,\ldots,T}\Vert \mbf P_{t|T}-(\bm\Lambda_n^\prime(\bm\Sigma_n^\xi)^{-1}\bm\Lambda_n)^{-1}\Vert = O_p(1)$;
\item $n^{2}\max_{t=1,\ldots,T}\Vert \mbf P_{0,t|T}-(\bm\Lambda_n^\prime(\bm\Gamma_n^\xi)^{-1}\bm\Lambda_n)^{-1}\Vert = O_p(1)$;
\end{compactenum} 
where $\mbf F_t^{\text{\tiny \upshape WLS}} = (\bm\Lambda_n^\prime(\bm\Sigma_n^\xi)^{-1}\bm\Lambda_n)^{-1}\bm\Lambda_n^\prime(\bm\Sigma_n^\xi)^{-1}\mbf x_t$ and 
$\mbf F_t^{\text{\tiny \upshape GLS}} = (\bm\Lambda_n^\prime(\bm\Gamma_n^\xi)^{-1}\bm\Lambda_n)^{-1}\bm\Lambda_n^\prime(\bm\Gamma_n^\xi)^{-1}\mbf x_t$. 
\end{lem}

\noindent
\noindent\textsc{Proof.} For part (i), from \eqref{eq:up1} we have
\al{
\mbf F_{t|t}-\mbf F_t^{\text{\tiny \upshape WLS}}=&\,(\bm\Lambda_n^\prime(\bm\Sigma_n^\xi)^{-1}\bm\Lambda_n+\mbf P_{t|t-1})^{-1}\bm\Lambda_n^\prime(\bm\Sigma_n^\xi)^{-1}\mbf x_t -\mbf F_t^{\text{\tiny \upshape WLS}}\nn\\
&+ \mbf F_{t|t-1} -(\bm\Lambda_n^\prime(\bm\Sigma_n^\xi)^{-1}\bm\Lambda_n+\mbf P_{t|t-1})^{-1}\bm\Lambda_n^\prime(\bm\Sigma_n^\xi)^{-1}\bm\Lambda_n\mbf F_{t|t-1}.\label{eq:crest}
}
Then,
\al{
\Vert&(\bm\Lambda_n^\prime(\bm\Sigma_n^\xi)^{-1}\bm\Lambda_n+\mbf P_{t|t-1})^{-1}\bm\Lambda_n^\prime(\bm\Sigma_n^\xi)^{-1}\mbf x_t -\mbf F_t^{\text{\tiny \upshape WLS}}\Vert \label{eq:crest2}\\
&\le\Vert (\bm\Lambda_n^\prime(\bm\Sigma_n^\xi)^{-1}\bm\Lambda_n+\mbf P_{t|t-1})^{-1}- (\bm\Lambda_n^\prime(\bm\Sigma_n^\xi)^{-1}\bm\Lambda_n)^{-1}\Vert\, \Vert \bm\Lambda_n\Vert\, \Vert(\bm\Sigma_n^\xi)^{-1}\Vert \,\Vert \mbf x_t \Vert = O_p(n^{-1}),\nn
}
by Lemmas  \ref{lem:LSL}(iii), \ref{lem:lambdasqrtn}, and \ref{lem:xunif}, and Assumption \ref{ass:idio}(a). And,
\al{
\Vert &\mbf F_{t|t-1} -(\bm\Lambda_n^\prime(\bm\Sigma_n^\xi)^{-1}\bm\Lambda_n+\mbf P_{t|t-1})^{-1}\bm\Lambda_n^\prime(\bm\Sigma_n^\xi)^{-1}\bm\Lambda_n\mbf F_{t|t-1}\Vert \nn\\
&\le \Vert \mbf I_r -(\bm\Lambda_n^\prime(\bm\Sigma_n^\xi)^{-1}\bm\Lambda_n+\mbf P_{t|t-1})^{-1}\bm\Lambda_n^\prime(\bm\Sigma_n^\xi)^{-1}\bm\Lambda_n\Vert\,\Vert \mbf F_{t|t-1}\Vert \nn\\
&= O_p(n^{-1}), \label{eq:crest3}
}
by \eqref{eq:krapfen} in the proof of Lemma \ref{lem:KF}(i). By substituting \eqref{eq:crest2} and \eqref{eq:crest3} into \eqref{eq:crest}, we prove part (i).

Part (ii) is proved as part (i), but using Lemmas \ref{lem:LSL}(iv) and \ref{lem:KF}(ii) instead of Lemmas \ref{lem:LSL}(iii) and \ref{lem:KF}(i). 

For part (iii), from \eqref{eq:up2}, using the same steps leading to \eqref{PtThathat} in the proof of Lemma \ref{lem:cazzarolahat00}, but when using the true parameters, it holds that:
\begin{align}
{\mbf P}_{t|t}
=&\,({\bm\Lambda}_n^{\prime}({\bm\Sigma}_n^{\xi})^{-1}{\bm\Lambda}_n)^{-1}\nn\\
&- {\mbf P}_{t|t-1}^{}
(({\bm\Lambda}_n^{\prime}({\bm\Sigma}_n^{\xi})^{-1}{\bm\Lambda}_n)^{-1}+ {\mbf P}_{t|t-1})^{-1}
({\bm\Lambda}_n^{\prime}({\bm\Sigma}_n^{\xi})^{-1}{\bm\Lambda}_n)^{-1}({\mbf P}_{t|t-1})^{-1}({\bm\Lambda}_n^{\prime}({\bm\Sigma}_n^{\xi})^{-1}{\bm\Lambda}_n)^{-1}.
\label{PtThathat_bernoulli}
\end{align}
Notice that all inverses in \eqref{PtThathat_bernoulli} are well defined because of Lemmas \ref{lem:LSL2}(iii), \ref{lem:cazzarola}(i), and \ref{lem:cazzarola}(ii) and Assumption \ref{ass:idio}(a).
Therefore, from \eqref{PtThathat_bernoulli}
\al{
n^2\max_{t=1,\ldots, T} \Vert {\mbf P}_{t|t}- ({\bm\Lambda}_n^{\prime}({\bm\Sigma}_n^{\xi})^{-1}{\bm\Lambda}_n)^{-1}\Vert
&\le \max_{t=1,\ldots, T} \Vert {\mbf P}_{t|t-1}\Vert \, \max_{t=1,\ldots, T} \Vert ({\mbf P}_{t|t-1})^{-1}\Vert\,
n^2\Vert ({\bm\Lambda}_n^{\prime}({\bm\Sigma}_n^{\xi})^{-1}{\bm\Lambda}_n)^{-1}\Vert^2\nn\\
&\cdot \Vert (({\bm\Lambda}_n^{\prime}({\bm\Sigma}_n^{\xi})^{-1}{\bm\Lambda}_n)^{-1}+ {\mbf P}_{t|t-1})^{-1}\Vert\nn\\
=&\, O(1),\label{eq:duffel}
}
 because of Lemmas \ref{lem:LSL2}(iii), \ref{lem:cazzarola}(i), and \ref{lem:cazzarola}(ii), and  since $\Vert (({\bm\Lambda}_n^{\prime}({\bm\Sigma}_n^{\xi})^{-1}{\bm\Lambda}_n)^{-1}+ {\mbf P}_{t|t-1})^{-1}\Vert=O(1)$, because of the same arguments leading to \eqref{eq:bochum} in the proof of Lemma \ref{lem:cazzarolahat00}, which this time hold by  
 Lemmas \ref{lem:LSL2}(iii), \ref{lem:cazzarola}(i), and \ref{lem:cazzarola}(ii), instead of Lemmas \ref{lem:est0_LAST}(iii), \ref{lem:cazzarolahat0}(i), and \ref{lem:cazzarolahat0}(ii).
 
 Part (iv) is proved as part (iii), but using  Lemmas \ref{lem:LSL2}(iv), and Assumption \ref{ass:idio}(f) instead of  Lemmas \ref{lem:LSL2}(iii) and Assumption \ref{ass:idio}(a). 
 
 For part (v),
 \[
 n^{2}\max_{t=1,\ldots,T}\Vert \mbf P_{t|T}-(\bm\Lambda_n^\prime(\bm\Sigma_n^\xi)^{-1}\bm\Lambda_n)^{-1}\Vert \le
 n^{2}\max_{t=1,\ldots,T}\Vert \mbf P_{t|T}-P_{t|t}\Vert +
 n^{2}\max_{t=1,\ldots,T}\Vert \mbf P_{t|t}-(\bm\Lambda_n^\prime(\bm\Sigma_n^\xi)^{-1}\bm\Lambda_n)^{-1}\Vert = O(1),
 \]
 because of part (iii) and Lemma \ref{lem:nuovo}(v). Part (vi) is proved  as part (v), but using part (iv) and Lemma \ref{lem:nuovo}(vi). This completes the proof. $\Box$

\begin{lem}\label{lem:2max} 
Let $\ell_0(\bm X_{nT};\underline{\bm \phi}_n)$ be the log-likelihood obtained when $\mbf A=\mbf 0_{r\times r}$ and $\bm\Gamma^v=\mbf I_r$. 
Then, under Assumptions \ref{ass:common}, \ref{ass:idio}, \ref{ass:ind},  \ref{ass:linear}, \ref{ass:tails}, and \ref{ass:ident}, as $n,T\to\infty$,
 \[
n\log^{-2/\delta_{v}} T \,\sup_{\underline{\bm\varphi}_n\in\mathcal O_n}(nT)^{-1}\l\vert\ell(\bm X_{nT};\underline{\bm\phi}_n,\underline{\bm\theta}) -\ell_0(\bm X_{nT};\underline{\bm\phi}_n)\r\vert = 
O_p(1).
\]
\end{lem}

\noindent\textsc{Proof.}
Throughout we consider generic values of the parameters such that $\underline{\bm\varphi}_n\in\mathcal O_n$ where $\mathcal O_n=\{\mathcal O_{\lambda_i}^n\cap \mathcal E_{\Lambda_n}\}\times\{\mathcal O_{\sigma_i^2}^n\cap\mathcal E_{\Gamma^\xi_n}\} \times \mathcal O_{\mathcal A}\times \mathcal O_{\Gamma^v}$ as defined in Section \ref{enzuccio}. Thus the elements of  $\underline{\bm\varphi}_n$ satisfy Assumptions \ref{ass:common}(a), \ref{ass:common}(d), \ref{ass:common}(e), \ref{ass:idio}(a), \ref{ass:idio}(b), and \ref{ass:idio}(f).

First, recall that the log-likelihood \eqref{eq:LL0true} can be written as
\begin{align}
\ell(\bm X_{nT};\underline{\bm\phi}_n,\underline{\bm\theta})=&\, \ell(\bm X_{nT}|\bm F_T;\underline{\bm\phi}_n)+\ell(\bm F_{T};\underline{\bm\theta})-\ell(\bm F_T|\bm X_{nT};\underline{\bm\phi}_n,\underline{\bm\theta})\nn\\
=&\,-\frac T2\log\det(\underline{\bm\Sigma}_n^\xi)-\frac 12\sum_{t=1}^T (\mbf x_{nt}-\underline{\bm\Lambda}_n\mbf F_t)^\prime(\underline{\bm\Sigma}_n^\xi)^{-1}(\mbf x_{nt}-\underline{\bm\Lambda}_n\mbf F_t)\nn\\
&-\frac T2 \log \det(\underline{\bm\Gamma}^v)-\frac 12\sum_{t=1}^T(\mbf F_t-\underline{\mbf A}\mbf F_{t-1})^\prime(\underline{\bm\Gamma}^v)^{-1}(\mbf F_t-\underline{\mbf A}\mbf F_{t-1})\label{giannone}\\
&+\frac 12\sum_{t=1}^T\log \det(\bm{\mathcal P}_{0,t|T}(\underline{\bm\phi}_n,\underline{\bm\theta}))\nn\\
&+\frac 12\sum_{t=1}^T(\mbf F_t-\bm{\mathcal F}_{0,t|T}(\underline{\bm\phi}_n,\underline{\bm\theta}))^\prime(\bm{\mathcal P}_{0,t|T}(\underline{\bm\phi}_n,\underline{\bm\theta}))^{-1}(\mbf F_t-\bm{\mathcal F}_{0,t|T}(\underline{\bm\phi}_n,\underline{\bm\theta})),\nn
\end{align}
where we used the fact that $\mbf F_0=\mbf 0_r$ by Assumption \ref{ass:common}(i) and we used the definitions:
\begin{align}
&\bm{\mathcal F}_{0,t|T}(\underline{\bm\phi}_n,\underline{\bm\theta})=\E_{\underline{\varphi}_n}[\mbf F_t|\bm X_{nT}]\equiv \underline{\bm{\mathcal F}}_{0,t|T},\nn\\
&\bm{\mathcal P}_{0,t|T}(\underline{\bm\phi}_n,\underline{\bm\theta})=\E_{\underline{\varphi}_n}[(\mbf F_t-\bm{\mathcal F}_{0,t|T}(\underline{\bm\phi}_n,\underline{\bm\theta}))(\mbf F_t-\bm{\mathcal F}_{0,t|T}(\underline{\bm\phi}_n,\underline{\bm\theta}))^\prime|\bm X_{nT}]\equiv\underline{\bm{\mathcal P}}_{0,t|T}.\label{doz}
\end{align}
Now, since \eqref{giannone} holds for any $\mbf F_t$ we can always choose $\mbf F_t=\underline{\bm{\mathcal F}}_{0,t|T}$ for all $t=1,\ldots, T$, so that
\begin{align}
\ell(\bm X_{nT};\underline{\bm\phi}_n,\underline{\bm\theta})=&\,-\frac T2\log\det(\underline{\bm\Sigma}_n^\xi)-\frac 12\sum_{t=1}^T (\mbf x_{nt}-\underline{\bm\Lambda}_n\underline{\bm{\mathcal F}}_{0,t|T})^\prime(\underline{\bm\Sigma}_n^\xi)^{-1}(\mbf x_{nt}-\underline{\bm\Lambda}_n\underline{\bm{\mathcal F}}_{0,t|T})\nn\\
&-\frac T2 \log \det(\underline{\bm\Gamma}^v)-\frac 12\sum_{t=1}^T(\underline{\bm{\mathcal F}}_{0,t|T}-\underline{\mbf A}\,\underline{\bm{\mathcal F}}_{0,t-1|T})^\prime(\underline{\bm\Gamma}^v)^{-1}(\underline{\bm{\mathcal F}}_{0,t|T}-\underline{\mbf A}\,\underline{\bm{\mathcal F}}_{0,t-1|T})\nn\\
&+\frac 12\sum_{t=1}^T\log\det (\underline{\bm{\mathcal P}}_{0,t|T}).\label{giannone2}
\end{align}

Second,  consider the log-likelihood \eqref{eq:LL0true} when the autocorrelation of the factors is not accounted for, i.e., when $\mbf A=\mbf 0_{r\times r}$ and $\bm\Gamma^v=\mbf I_r$, 
\beq\label{eq:LL0iidEM}
\ell_0(\bm X_{nT};\underline{\bm \phi}_n)=-\frac T2 \log\det\l(\underline{\bm \Lambda}_{n}\underline{\bm \Lambda}_{n}^\prime + \underline{\bm\Sigma}_n^\xi\r)-\frac 12\sum_{t=1}^T \l[ \mbf x_{nt}^\prime \l(\underline{\bm \Lambda}_{n}\underline{\bm \Lambda}_{n}^\prime + \underline{\bm\Sigma}_n^\xi\r)^{-1}\mbf x_{nt} \r],
\eeq
where we are imposing Assumption \ref{ass:ident}(b), so that we can set $\underline{\bm\Gamma}^F=\mbf I_r$ in the log-likelihood. Clearly, \eqref{eq:LL0iidEM} can also be written as
\begin{align}
\ell_0(\bm X_{nT};\underline{\bm\phi}_n)=&\, \ell_0(\bm X_{nT}|\bm F_T;\underline{\bm\phi}_n)+\ell_0(\bm F_{T})-\ell_0(\bm F_T|\bm X_{nT};\underline{\bm\phi}_n)\nn\\
=&\,-\frac T2\log\det(\underline{\bm\Sigma}_n^\xi)-\frac 12\sum_{t=1}^T (\mbf x_{nt}-\underline{\bm\Lambda}_n\mbf F_t)^\prime(\underline{\bm\Sigma}_n^\xi)^{-1}(\mbf x_{nt}-\underline{\bm\Lambda}_n\mbf F_t)\nn\\
&-\frac 12\sum_{t=1}^T \mbf F_t^\prime\mbf F_t+\frac 12\sum_{t=1}^T\log \det(\bm{\mathcal V}_{0,t|T}(\underline{\bm\phi}_n))\label{giannone33}\\
&+\frac 12\sum_{t=1}^T(\mbf F_t-\bm{\mathcal G}_{0,t|T}(\underline{\bm\phi}_n))^\prime(\bm{\mathcal V}_{0,t|T}(\underline{\bm\phi}_n))^{-1}(\mbf F_t-\bm{\mathcal G}_{0,t|T}(\underline{\bm\phi}_n)),\nn
\end{align}
where
\begin{align}
&\bm{\mathcal G}_{0,t|T}(\underline{\bm\phi}_n)=\E_{\underline{\phi}_n}[\mbf F_t|\bm X_{nT}]\equiv \underline{\bm{\mathcal G}}_{0,t|T},\nn\\
&\bm{\mathcal V}_{0,t|T}(\underline{\bm\phi}_n)=\E_{\underline{\phi}_n}[(\mbf F_t-\bm{\mathcal G}_{0,t|T}(\underline{\bm\phi}_n))(\mbf F_t-\bm{\mathcal G}_{0,t|T}(\underline{\bm\phi}_n))^\prime|\bm X_{nT}]\equiv\underline{\bm{\mathcal V}}_{0,t|T}.\label{doz2}
\end{align}
Now, since \eqref{giannone33} holds for any $\mbf F_t$ we can always choose $\mbf F_t=\underline{\bm{\mathcal G}}_{0,t|T}$ for all $t=1,\ldots, T$, so that
\begin{align}
\ell_0(\bm X_{nT};\underline{\bm\phi}_n)
=&\,-\frac T2\log\det(\underline{\bm\Sigma}_n^\xi)-\frac 12\sum_{t=1}^T (\mbf x_{nt}-\underline{\bm\Lambda}_n\underline{\bm{\mathcal G}}_{0,t|T})^\prime(\underline{\bm\Sigma}_n^\xi)^{-1}(\mbf x_{nt}-\underline{\bm\Lambda}_n\underline{\bm{\mathcal G}}_{0,t|T})\nn\\
&-\frac 12\sum_{t=1}^T \underline{\bm{\mathcal G}}_{0,t|T}^\prime\underline{\bm{\mathcal G}}_{0,t|T}+\frac 12\sum_{t=1}^T\log \det(\underline{\bm{\mathcal V}}_{0,t|T})\label{giannone3}.
\end{align}

Under Assumption \ref{ass:linear} the conditional mean, $\underline{\bm{\mathcal F}}_{0,t|T}$  in \eqref{doz} is a linear function of $\bm X_{nT}$, so it can be obtained by linear projection, and, therefore it is given by the correctly specified Kalman smoother, i.e., using as parameters $\underline{\bm\phi}_n$ and $\underline{\bm\theta}$ and when replacing  $\underline{\bm\Sigma}_n^\xi$ with $\underline{\bm\Gamma}_n^\xi$, thus,
\begin{align}
\underline{\bm{\mathcal F}}_{0,t|T}=\underline{\mbf F}_{0,t|T}.
\label{eq:F0P0}
\end{align}
Likewise, under Assumption \ref{ass:linear}, $\underline{\bm{\mathcal G}}_{0,t|T}$ 
in \eqref{doz2} is  also linear, however, since this is the conditional mean for the case in which no dynamics for the factors is specified, it is given by the simpler linear projection (using Lemma \ref{lem:wood})
\begin{align}
\underline{\bm{\mathcal G}}_{0,t|T}&=\underline{\bm\Lambda}_n^\prime(\underline{\bm\Lambda}_n\underline{\bm\Lambda}_n^\prime+\underline{\bm\Gamma}_n^\xi)^{-1}\mbf x_{nt}= ( \underline{\bm\Lambda}_n^\prime(\underline{\bm\Gamma}_n^\xi)^{-1}\underline{\bm\Lambda}_n+\mbf I_r)^{-1} \underline{\bm\Lambda}_n^\prime(\underline{\bm\Gamma}_n^\xi)^{-1}\mbf x_{nt}=\underline{\mbf F}_{0,t}^{\text{\tiny REG}}.\label{eq:pinnagialla}
\end{align}

Now, consider the generalized least squares estimator of the factors computed for generic values of the parameters $\underline{\bm\phi}_n$:
\beq
\underline{\mbf F}_t^{\text{\tiny GLS}}=( \underline{\bm\Lambda}_n^\prime(\underline{\bm\Gamma}_n^\xi)^{-1}\underline{\bm\Lambda}_n)^{-1} \underline{\bm\Lambda}_n^\prime(\underline{\bm\Gamma}_n^\xi)^{-1}\mbf x_{nt},
\eeq
Then, since we restrict to $\underline{\bm\varphi}_n\in\mathcal O_n$,  it follows that
\al{
\Vert\underline{\mbf F}_{0,t}^{\text{\tiny REG}}-\underline{\mbf F}_t^{\text{\tiny GLS}}\Vert 
&\le 
n \Vert ( \underline{\bm\Lambda}_n^\prime(\underline{\bm\Gamma}_n^\xi)^{-1}\underline{\bm\Lambda}_n+\mbf I_r)^{-1} 
 -
 ( \underline{\bm\Lambda}_n^\prime(\underline{\bm\Gamma}_n^\xi)^{-1}\underline{\bm\Lambda}_n)^{-1}\Vert\, n^{-1/2}\Vert \underline{\bm\Lambda}_n^\prime(\underline{\bm\Gamma}_n^\xi)^{-1}\Vert \, n^{-1/2}\Vert\mbf x_{nt}\Vert\nn\\
&=  O(n^{-1}) O_p(1),\label{eq:reg_wls}
}
by Lemmas  \ref{lem:LSL2}(viii), \ref{lem:LSL}(iv), and \ref{lem:xunif}.

Moreover, by Lemmas \ref{lem:KSKF}(ii) and \ref{lem:KFGLS}(ii) 
\beq\label{eq:lemma9}
\Vert\underline{\mbf F}_{0,t|T}-\underline{\mbf F}_t^{\text{\tiny GLS}}\Vert\le \Vert\underline{\mbf F}_{0,t|T}-\underline{\mbf F}_{0,t|t}\Vert+\Vert\underline{\mbf F}_{0,t|t}-\underline{\mbf F}_t^{\text{\tiny GLS}}\Vert = O(n^{-1}) O_p(1).
\eeq
In particular, the bound in \eqref{eq:lemma9} is a product of a stochastic and a non-stochastic term because of \eqref{eq:E7iidet} in the proof of Lemma \ref{lem:KSKF}, and \eqref{eq:crest2} and \eqref{eq:crest3} in the proof of Lemma \ref{lem:KFGLS}.

Therefore, from \eqref{eq:reg_wls} and \eqref{eq:lemma9} (see also \citealp[]{baili16})
\beq
\Vert\underline{\mbf F}_{0,t|T}-\underline{\mbf F}_{0,t}^{\text{\tiny REG}}\Vert = O(n^{-1})O_p(1).\label{eq:reg_KS}
\eeq
The results in \eqref{eq:reg_wls}-\eqref{eq:reg_KS} depend on $t$ only through the $O_p(1)$ term which in turn is just function of $\mbf x_{nt}$ and do not depend on the choice of the parameters as long as they belong to $\mathcal O_n$, as assumed. Since the parameters are deterministic and using Lemma \ref{lem:xunif2}, we then have
\beq\label{eq:fufix}
\sup_{\underline{\bm\varphi}_n\in\mathcal O_n}\;\max_{t=1,\ldots, T}\Vert\underline{\mbf F}_{0,t|T}-\underline{\mbf F}_{0,t}^{\text{\tiny REG}}\Vert = 
O(n^{-1})O_p(\log^{1/\delta_{v}} T).
\eeq

Then, consider the conditional covariance in \eqref{doz2} and letting $\underline{\mbf K}_n=( \underline{\bm\Lambda}_n^\prime(\underline{\bm\Gamma}_n^\xi)^{-1}\underline{\bm\Lambda}_n+\mbf I_r)^{-1} \underline{\bm\Lambda}_n^\prime(\underline{\bm\Gamma}_n^\xi)^{-1}$, we have
\al{
\underline{\bm{\mathcal V}}_{0,t|T}=&\, 
(\mbf I_r- \underline{\mbf K}_n\bm\Lambda_n)
\E_{\underline{\phi}_n}[\mbf F_t\mbf F_t^\prime|\bm X_{nT}]
(\mbf I_r- \underline{\mbf K}_n\bm\Lambda_n)^\prime+ \underline{\mbf K}_n\E_{\underline{\phi}_n}[\bm\xi_{nt}\bm\xi_{nt}^\prime|\bm X_{nT}]\underline{\mbf K}_n^\prime,\nn
}
where we used also Lemma \ref{lem:fidio}. Moreover, 
\al{
\E[\underline{\bm{\mathcal V}}_{0,t|T}] = (\mbf I_r- \underline{\mbf K}_n\bm\Lambda_n)
\E[\E_{\underline{\phi}_n}[\mbf F_t\mbf F_t^\prime|\bm X_{nT}]]
(\mbf I_r- \underline{\mbf K}_n\bm\Lambda_n)^\prime+ \underline{\mbf K}_n\E[\E_{\underline{\phi}_n}[\bm\xi_{nt}\bm\xi_{nt}^\prime|\bm X_{nT}]]\underline{\mbf K}_n^\prime,\nn
}
is finite and positive definite, 
since  $\Vert \underline{\mbf K}_n\bm\Lambda_n\Vert =O(1)$ by Lemmas \ref{lem:lambdasqrtn}, \ref{lem:LSL2}(vi), and \ref{lem:LSL2}(viii), and because\linebreak $\E[\sup_{\underline{\bm\varphi}_n\in\mathcal O_n}\;\max_{t=1,\ldots, T}\E_{\underline{\phi}_n}[\mbf F_t\mbf F_t^\prime|\bm X_{nT}]]$ and $\E[\sup_{\underline{\bm\varphi}_n\in\mathcal O_n}\;\max_{t=1,\ldots, T}\E_{\underline{\phi}_n}[\bm\xi_{nt}\bm\xi_{nt}^\prime|\bm X_{nT}]]$ are both finite and positive definite by Assumptions \ref{ass:common}(b) and \ref{ass:idio}(f), and Lemma \ref{lem:Gxi}(v).
Therefore, by Markov's inequality
\al{
\sup_{\underline{\bm\varphi}_n\in\mathcal O_n\setminus \{\bm\varphi_n\}}\;\max_{t=1,\ldots, T}\Vert \underline{\bm{\mathcal V}}_{0,t|T}\Vert = O_p(1), \quad \sup_{\underline{\bm\varphi}_n\in\mathcal O_n\setminus \{\bm\varphi_n\}}\;\max_{t=1,\ldots, T}\Vert (\underline{\bm{\mathcal V}}_{0,t|T})^{-1}\Vert = O_p(1).\label{eq:fufix2}
}
This bound is tighter when $\underline{\bm\varphi}_n=\bm\varphi_n$. Indeed, in that case from \eqref{eq:up2} and \eqref{eq:KS2}, setting $\mbf A=\mbf 0_{r\times r}$ therein, we have 
\[
{\bm{\mathcal V}}_{0,t|T}(\bm\phi_n)= \mbf I_r-\mbf I_r \bm\Lambda_n^\prime(\bm\Lambda_n\bm\Lambda_n^\prime+\bm\Gamma_n^\xi)^{-1}\bm\Lambda_n=
\mbf I_r-(\bm\Lambda_n^\prime(\bm\Gamma_n^\xi)^{-1}\bm\Lambda_n+\mbf I_r)^{-1}\bm\Lambda_n^\prime(\bm\Gamma_n^\xi)^{-1}\bm\Lambda_n
.
\]
where we also used Assumption \ref{ass:ident}(b) for which $\bm\Gamma^F=\mbf I_r$. Therefore,  ${\bm{\mathcal V}}_{0,t|T}(\bm\phi_n)$ is deterministic and independent of $t$, and such that
\al{
\max_{t=1,\ldots, T}\Vert {\bm{\mathcal V}}_{0,t|T}(\bm\phi_n)\Vert =O(n^{-1}),\qquad \max_{t=1,\ldots, T}\Vert ({\bm{\mathcal V}}_{0,t|T}(\bm\phi_n))^{-1}\Vert = O(n), \label{altre2}
} 
because by Lemma \ref{lem:LSL}(ii) which holds since $\bm\Gamma_n^\xi$ is positive definite by Assumption \ref{ass:idio}(f). Furthermore, following the same steps leading to \eqref{eq:duffel} in the proof of Lemma \ref{lem:KFGLS} we have
\al{
\max_{t=1,\ldots, T} \Vert {\bm{\mathcal V}}_{0,t|T}(\bm\phi_n)- ({\bm\Lambda}_n^{\prime}({\bm\Gamma}_n^{\xi})^{-1}{\bm\Lambda}_n)^{-1}\Vert = O(n^{-2}).\label{eq:duffel2}
}

%
 
Turning to the conditional covariance in \eqref{doz}, because of \eqref{eq:F0P0} we can write
\al{
\underline{\bm{\mathcal P}}_{0,t|T}=&\,\E_{\underline{\varphi}_n}[(\mbf F_t-\underline{\mbf F}_{0,t|T}+\underline {\mbf F}_{0,t}^{\text{\tiny REG}}-\underline {\mbf F}_{0,t}^{\text{\tiny REG}})(\mbf F_t-\underline{\mbf F}_{0,t|T}+\underline {\mbf F}_{0,t}^{\text{\tiny REG}}-\underline {\mbf F}_{0,t}^{\text{\tiny REG}})^\prime|\bm X_{nT}]\nn\\
=&\,\underline{\bm{\mathcal V}}_{0,t|T}+\E_{\underline{\varphi}_n}[(\underline {\mbf F}_{0,t}^{\text{\tiny REG}}-\underline{\mbf F}_{0,t|T})(\underline {\mbf F}_{0,t}^{\text{\tiny REG}}-\underline{\mbf F}_{0,t|T})^\prime|\bm X_{nT}]\nn\\
&+\E_{\underline{\varphi}_n}[(\underline {\mbf F}_{0,t}^{\text{\tiny REG}}-\underline{\mbf F}_{0,t|T})(\mbf F_t-\underline {\mbf F}_{0,t}^{\text{\tiny REG}})^\prime|\bm X_{nT}]+\E_{\underline{\varphi}_n}[(\mbf F_t-\underline {\mbf F}_{0,t}^{\text{\tiny REG}})(\underline {\mbf F}_{0,t}^{\text{\tiny REG}}-\underline{\mbf F}_{0,t|T})^\prime|\bm X_{nT}].\label{eq:bielefeld}
}
From \eqref{eq:fufix} and \eqref{eq:bielefeld} we have
\al{
\sup_{\underline{\bm\varphi}_n\in\mathcal O_n\setminus \{\bm\varphi_n\}}\;\max_{t=1,\ldots, T}\Vert \underline{\bm{\mathcal P}}_{0,t|T}-\underline{\bm{\mathcal V}}_{0,t|T}\Vert  = O_p(n^{-1}\log^{1/\delta_v}T).\label{eq:fufix3}
}
This bound is tighter when  $\underline{\bm\varphi}_n=\bm\varphi_n$. Indeed, we have ${\bm{\mathcal P}}_{0,t|T}(\bm\varphi_n)=\mbf P_{0,t|T}$, which is deterministic by Lemma \ref{lem:detP}(iii), and 
\al{
\Vert {\bm{\mathcal P}}_{0,t|T}(\bm\varphi_n)-{\bm{\mathcal V}}_{0,t|T}(\bm\varphi_n)\Vert
&=\Vert \mbf P_{0,t|T}-{\bm{\mathcal V}}_{0,t|T}(\bm\varphi_n)\Vert\nn\\
&\le
\Vert \mbf P_{0,t|T}-(\bm\Lambda_n^\prime(\bm\Gamma_n^\xi)^{-1}\bm\Lambda_n)^{-1}\Vert
+
\Vert (\bm\Lambda_n^\prime(\bm\Gamma_n^\xi)^{-1}\bm\Lambda_n)^{-1}-{\bm{\mathcal V}}_{0,t|T}(\bm\varphi_n)\Vert\nn\\
&=O(n^{-2}),\label{altre22}
} 
by \eqref{eq:duffel2} and Lemma \ref{lem:KFGLS}(vi).


\color{black}

Now, consider
\begin{align}
&\ell(\bm X_{nT};\underline{\bm\phi}_n,\underline{\bm\theta}) -\ell_0(\bm X_{nT};\underline{\bm\phi}_n)\nn\\
=&\,
-\frac 12\sum_{t=1}^T
\bigg\{ 
\l(\underline{\mbf F}_{0,t|T}-\underline{\mbf F}_{0,t}^{\text{\tiny REG}}\r)^\prime
\underline{\bm\Lambda}_n^\prime(\underline{\bm\Sigma}_n^\xi)^{-1}\underline{\bm\Lambda}_n
\l(\underline{\mbf F}_{0,t|T}-\underline{\mbf F}_{0,t}^{\text{\tiny REG}}\r)\nn\\
&+2\l(\underline{\mbf F}_{0,t|T}-\underline{\mbf F}_{0,t}^{\text{\tiny REG}}\r)^\prime
\underline{\bm\Lambda}_n^\prime(\underline{\bm\Sigma}_n^\xi)^{-1}\underline{\bm\Lambda}_n
\underline{\mbf F}_{0,t}^{\text{\tiny REG}}
+2\mbf x_{nt}^\prime(\underline{\bm\Sigma}_n^\xi)^{-1}\underline{\bm\Lambda}_n
\l(\underline{\mbf F}_{0,t|T}-\underline{\mbf F}_{0,t}^{\text{\tiny REG}}\r)\bigg\}-\frac T2 \log \det(\underline{\bm\Gamma}^v)\nn\\
&-\frac 12\sum_{t=1}^T(\underline{\mbf F}_{0,t|T}-\underline{\mbf A}\,\underline{\mbf F}_{0,t-1|T})^\prime
(\underline{\bm\Gamma}^v)^{-1}
(\underline{\mbf F}_{0,t|T}-\underline{\mbf A}\,\underline{\mbf F}_{0,t-1|T})+\frac 12\sum_{t=1}^T \underline{\mbf F}_{0,t}^{\text{\tiny REG}\prime}\underline{\mbf F}_{0,t}^{\text{\tiny REG}}\nn\\
&+\frac 12\sum_{t=1}^T\log\det (\underline{\bm{\mathcal P}}_{0,t|T})
- \frac 12\sum_{t=1}^T\log \det(\underline{\bm{\mathcal V}}_{0,t|T}).\label{eq:fufix4}
\end{align}
Consider all terms on the rhs of \eqref{eq:fufix4}. First, by \eqref{eq:fufix}, for all $\underline{\bm{\varphi}}_n\in\mathcal O_n$,
\begin{align}
&\bigg\vert -\frac 12\sum_{t=1}^T
\bigg\{ 
\l(\underline{\mbf F}_{0,t|T}-\underline{\mbf F}_{0,t}^{\text{\tiny REG}}\r)^\prime
\underline{\bm\Lambda}_n^\prime(\underline{\bm\Sigma}_n^\xi)^{-1}\underline{\bm\Lambda}_n
\l(\underline{\mbf F}_{0,t|T}-\underline{\mbf F}_{0,t}^{\text{\tiny REG}}\r)
+2\l(\underline{\mbf F}_{0,t|T}-\underline{\mbf F}_{0,t}^{\text{\tiny REG}}\r)^\prime
\underline{\bm\Lambda}_n^\prime(\underline{\bm\Sigma}_n^\xi)^{-1}\underline{\bm\Lambda}_n
\underline{\mbf F}_{0,t}^{\text{\tiny REG}}
\nn\\
&+2\mbf x_{nt}^\prime
\underline{\bm\Lambda}_n^\prime(\underline{\bm\Sigma}_n^\xi)^{-1}\underline{\bm\Lambda}_n
\l(\underline{\mbf F}_{0,t|T}-\underline{\mbf F}_{0,t}^{\text{\tiny REG}}\r)\bigg\}
\bigg\vert\nn\\
\le &\, 
T \max_{t=1,\ldots, T}\Vert\underline{\mbf F}_{0,t|T}-\underline{\mbf F}_{0,t}^{\text{\tiny REG}}\Vert
\,
\Vert\underline{\bm\Lambda}_n\Vert^2 \,\Vert(\underline{\bm\Sigma}_n^\xi)^{-1}\Vert\,
\max_{t=1,\ldots, T} \Vert\underline{\mbf F}_{0,t}^{\text{\tiny REG}}\Vert\nn\\
&+
T \max_{t=1,\ldots, T}\Vert\underline{\mbf F}_{0,t|T}-\underline{\mbf F}_{0,t}^{\text{\tiny REG}}\Vert
\,
\Vert(\underline{\bm\Sigma}_n^\xi)^{-1}\Vert\, \Vert\underline{\bm\Lambda}_n\Vert\,
\max_{t=1,\ldots, T}\Vert{\mbf x}_{nt}\Vert \nn\\ 
&+T\max_{t=1,\ldots, T}\Vert\underline{\mbf F}_{0,t|T}-\underline{\mbf F}_{0,t}^{\text{\tiny REG}}\Vert^2
\,
\Vert\underline{\bm\Lambda}_n\Vert^2\, \Vert(\underline{\bm\Sigma}_n^\xi)^{-1}\Vert\max_{t=1,\ldots, T}\Vert{\mbf x}_{nt}\Vert\nn\\
=&\,O_p(T\log^{2/\delta_v} T),\label{eq:fufix5}
\end{align}
where we used Lemmas \ref{lem:lambdasqrtn} and \ref{lem:xunif2}, which implies also that $\Vert\underline{\mbf F}_{0,t}^{\text{\tiny REG}}\Vert=O_p(\log^{1/\delta_v}T)$, and
 Assumption \ref{ass:idio}(a). 

Second, for all $\underline{\bm{\varphi}}_n\in\mathcal O_n$,
\begin{align}
&\bigg\vert -\frac T2 \log \det(\underline{\bm\Gamma}^v)-\frac 12\sum_{t=1}^T(\underline{\mbf F}_{0,t|T}-\underline{\mbf A}\,\underline{\mbf F}_{0,t-1|T})^\prime
(\underline{\bm\Gamma}^v)^{-1}
(\underline{\mbf F}_{0,t|T}-\underline{\mbf A}\,\underline{\mbf F}_{0,t-1|T})
+\frac 12\sum_{t=1}^T \underline{\mbf F}_{0,t}^{\text{\tiny REG}\prime}\underline{\mbf F}_{0,t}^{\text{\tiny REG}}\bigg\vert\nn\\
\le&\, \l\vert \frac T2 \log \det(\underline{\bm\Gamma}^v)\r\vert + \l\vert \frac 12\sum_{t=1}^T(\underline{\mbf F}_{0,t|T}-\underline{\mbf A}\,\underline{\mbf F}_{0,t-1|T})^\prime
(\underline{\bm\Gamma}^v)^{-1}
(\underline{\mbf F}_{0,t|T}-\underline{\mbf A}\,\underline{\mbf F}_{0,t-1|T})\r\vert+\l\vert \frac 12\sum_{t=1}^T \underline{\mbf F}_{0,t}^{\text{\tiny REG}\prime}\underline{\mbf F}_{0,t}^{\text{\tiny REG}}\r\vert\nn\\
=&\,O_p(T\log^{2/\delta_v} T),\label{eq:fufix6}
\end{align}
because of  \eqref{eq:reg_KS}, Assumptions \ref{ass:common}(d) and \ref{ass:common}(e),
and Lemma \ref{lem:xunif2} jointly with the same arguments leading to \eqref{eq:fufix5}.

Third, by \eqref{eq:fufix2} and \eqref{eq:fufix3}, and \citet[Theorem 1]{MK04}, which is Weyl's inequality 
\begin{align}
\bigg\vert& \frac 12\sum_{t=1}^T\log\det (\underline{\bm{\mathcal P}}_{0,t|T})
-\frac 12\sum_{t=1}^T \log \det(\underline{\bm{\mathcal V}}_{0,t|T})\bigg\vert=
\bigg\vert \frac 12\sum_{t=1}^T \log \l(\det \{\underline{\bm{\mathcal P}}_{0,t|T} (\underline{\bm{\mathcal V}}_{0,t|T})^{-1}\} \r)
\bigg\vert 
\nn\\
= &\, \bigg\vert\frac 12\sum_{t=1}^T \sum_{j=1}^r\log\l(
\l\{ \nu^{(j)}( \underline{\bm {\mathcal P}}_{0,t|T}) \r\}
 \l\{ \nu^{(j)}(\underline{\bm {\mathcal V}}_{0,t|T})\r\}^{-1}\r)\bigg\vert\nn\\
\le &\, \bigg\vert\frac 12\sum_{t=1}^T \sum_{j=1}^r\log\l(
\l\{ \nu^{(1)}( \underline{\bm{\mathcal P}}_{0,t|T}-\underline{\bm{\mathcal V}}_{0,t|T})+\nu^{(j)}(\underline{\bm{\mathcal V}}_{0,t|T}) \r\}
 \l\{ \nu^{(j)}(\underline{\bm{\mathcal V}}_{0,t|T})\r\}^{-1}\r)\bigg\vert\nn\\
 \le &\, \bigg\vert\frac 12\sum_{t=1}^T \sum_{j=1}^r\log\l(
1+\l\{ n\Vert \underline{\bm{\mathcal P}}_{0,t|T}-\underline{\bm{\mathcal V}}_{0,t|T}\Vert \r\}
 \l\{ n\nu^{(j)}(\underline{\bm{\mathcal V}}_{0,t|T})\r\}^{-1}\r)\bigg\vert\nn\\
 \le&\, \max_{t=1,\ldots, T}\bigg\vert\frac T2 \l\{ n\Vert \underline{\bm{\mathcal P}}_{0,t|T}-\underline{\bm{\mathcal V}}_{0,t|T}\Vert\r\}
 \l\{ n\nu^{(j)}(\underline{\bm{\mathcal V}}_{0,t|T})\r\}^{-1}\bigg\vert +o(Tn^{-1})\nn\\
 =&\, O_p(Tn^{-1}\log^{1/\delta_v} T),
\label{eq:fufix7}
\end{align}
where in the second last line we took into account also \eqref{altre2} and \eqref{altre22}, hence, \eqref{eq:fufix7} holds for all $\underline{\bm\varphi}_n\in\mathcal O_n$.

Summing up, by noticing that \eqref{eq:fufix5}, \eqref{eq:fufix6}, and  \eqref{eq:fufix7} hold for all $\underline{\bm\varphi}_n\in\mathcal O_n$, from  \eqref{eq:fufix4} we have:
\beq
\!\!\!\!\!\sup_{\underline{\bm\varphi}_n\in\mathcal O_n}(nT)^{-1}\l\vert\ell(\bm X_{nT};\underline{\bm\phi}_n,\underline{\bm\theta}) -\ell_0(\bm X_{nT};\underline{\bm\phi}_n)\r\vert = 
O_p(n^{-1}\log^{2/\delta_{v}} T).\nn
\eeq
This completes the proof. $\Box$


\begin{lem}\label{lem:dagstar}
Let $\wh{\bm\phi}_n^*=(\text{\upshape vec}(\wh{\bm\Lambda}_n^*)^\prime\; \wh{\sigma}^{2*}_{1}\cdots \wh{\sigma}^{2*}_{n})^\prime$ be the vector of QML estimators of the entries of ${\bm\phi}_n$ maximizing
 $\ell(\bm X_{nT};\underline{\bm \phi}_n,\underline{\bm\theta})$ defined in \eqref{eq:LL0true}, 
and let $\wh{\bm\phi}_n^\dag=(\text{\upshape vec}(\wh{\bm\Lambda}_n^\dag)^\prime\; \wh{\sigma}^{2\dag}_{1}\cdots \wh{\sigma}^{2\dag}_{n})^\prime$ be the vector of QML estimators of the entries of ${\bm\phi}_n$ maximizing
 $\ell_0(\bm X_{nT};\underline{\bm \phi}_n)$ defined in \eqref{eq:LL0iidEM} in the proof of Lemma \ref{lem:2max}. Then, under Assumptions \ref{ass:common}, \ref{ass:idio}, \ref{ass:ind}, \ref{ass:linear}, \ref{ass:tails}, and \ref{ass:ident}, as $n,T\to\infty$,
\begin{compactenum}[(i)]
\item $n\log^{-2/\delta_{v}} T \max_{i=1,\ldots, n}\Vert \wh{\bm\lambda}_i^*-\wh{\bm\lambda}_i^\dag\Vert = O_p(1)$;
\item $n\log^{-2/\delta_{v}} T n^{-1/2}\Vert \wh{\bm\Lambda}_n^*-\wh{\bm\Lambda}_n^\dag\Vert = O_p(1)$;
\item $n\log^{-2/\delta_{v}} T \max_{i=1,\ldots, n}\vert \wh{\sigma}_i^{2*}-\wh{\sigma}_i^{2\dag}\vert =O_p(1)$.
\end{compactenum}
\end{lem}

\noindent
\textsc{Proof.} From Lemma \ref{lem:2max} we have
 \al{
 (nT)^{-1}&\l\vert\sup_{\underline{\bm\phi}_n\in \mathcal O_n} \ell(\bm X_{nT};\underline{\bm\phi}_n,\underline{\bm\theta}) - \sup_{\underline{\bm\phi}_n\in \mathcal O_n}\ell_0(\bm X_{nT};\underline{\bm\phi}_n) \r\vert\nn\\
&\le (nT)^{-1}\sup_{\underline{\bm\phi}_n\in \mathcal O_n} \l\vert\ell(\bm X_{nT};\underline{\bm\phi}_n,\underline{\bm\theta}) -\ell_0(\bm X_{nT};\underline{\bm\phi}_n)\r\vert\nn\\
& = O_p(n^{-1}\log^{2/\delta_{v}} T). \label{eq:consiglio17}
}
Therefore, by continuity of the log-likelihoods and \eqref{eq:consiglio17}, 
 we have
\begin{align}
\max_{i=1,\ldots, n}\l(
\ba{c}
\Vert \wh{\bm\lambda}_i^*-\wh{\bm\lambda}_i^\dag\Vert \\
\vert \wh{\sigma}_i^{2*}-\wh{\sigma}_i^{2\dag}\vert
\ea
\r)
&=\max_{i=1,\ldots, n}\bigg\Vert\arg\!\!\!\!\!\!\!\!\!\max_{(\underline{\bm\lambda}_i^\prime\,\underline{\sigma}_i^2)^\prime\in \mathcal O_n} \ell(\bm X_{nT};\underline{\bm\phi}_n,\underline{\bm\theta}) - \arg\!\!\!\!\!\!\!\!\!\max_{(\underline{\bm\lambda}_i^\prime\,\underline{\sigma}_i^2)^\prime\in \mathcal O_n} \ell_0(\bm X_{nT};\underline{\bm\phi}_n) \bigg\Vert\nn\\
&=
\max_{i=1,\ldots, n}\bigg\Vert\arg\!\!\!\!\!\!\!\!\!\max_{(\underline{\bm\lambda}_i^\prime\,\underline{\sigma}_i^2)^\prime\in \mathcal O_1} (nT)^{-1}\ell(\bm X_{nT};\underline{\bm\phi}_n,\underline{\bm\theta}) - \arg\!\!\!\!\!\!\!\!\!\max_{(\underline{\bm\lambda}_i^\prime\,\underline{\sigma}_i^2)^\prime\in \mathcal O_n} (nT)^{-1}\ell_0(\bm X_{nT};\underline{\bm\phi}_n) \bigg\Vert\nn\\
& = O_p(n^{-1}\log^{2/\delta_{v}} T).\nn
\end{align}
This proves parts (i) and (iii), while part (ii) is a direct consequence of part (i). This completes the proof. $\Box$

\begin{lem}\label{lem:starols} 
Under Assumptions \ref{ass:common}, \ref{ass:idio}, \ref{ass:ind}, \ref{ass:linear}, \ref{ass:tails}, and \ref{ass:ident}, as $n,T\to\infty$,
\begin{compactenum}
\item [(i)] 
$\min(n \log^{-2/\delta_v} T, \sqrt{nT}, T)\, \Vert \wh{\bm\lambda}_i^{*} -\bm\lambda_i^{\text{\tiny \upshape OLS}}\Vert=O_p(1)$, uniformly in $i$;
\item [(ii)] $\min(n\log^{-2/\delta_v} T, \sqrt{nT}, T)\,
n^{-1/2} \Vert \wh{\bm\Lambda}_n^{*} -\bm\Lambda_n^{\text{\tiny \upshape OLS}}\Vert= O_p(1)$;
\item [(iii)] $\min(n\log^{-2/\delta_v} T, \sqrt{nT}, T)\,\vert \wh{\sigma}_i^{2*}-\sigma_i^{2\text{\tiny \upshape OLS}}\vert=O_p(1)$, uniformly in $i$.
\end{compactenum}
\end{lem}

\noindent\textsc{Proof.} For part (i) 
\al{
\Vert \wh{\bm\lambda}_i^{*} -\bm\lambda_i^{\text{\tiny \upshape OLS}}\Vert \le &\, \Vert \wh{\bm\lambda}_i^{*} -\bm\lambda_i^{\dag}\Vert+
\Vert \wh{\bm\lambda}_i^{\dag} -\bm\lambda_i^{(0)}\Vert+\Vert \wh{\bm\lambda}_i^{(0)} -\bm\lambda_i^{\text{\tiny \upshape OLS}}\Vert,\label{eq:5minuti}
}
where $\wh{\bm\lambda}_i^{(0)}$ is the PC estimator defined in Appendix \ref{app:prest}. Then, from \citet[Theorem 3]{MBPCAQML} 
\begin{align}
&\Vert\wh{\bm\lambda}_i^\dag-\wh{\bm\lambda}_i^{(0)} \Vert = O_p(n^{-1}),\label{eq:chiave2}
\end{align}
while from \citet[Corollary 1 and Proposition B.3]{MBPCAQML} 
\begin{align}\label{eq:chiave3}
&\Vert\wh{\bm\lambda}_i^{(0)}-{\bm\lambda}_i^{\text{\tiny OLS}}\Vert = O_p(\max(n^{-1},n^{-1/2}T^{-1/2},T^{-1}) ).
\end{align}
Part (i) then follows by substituting \eqref{eq:chiave2}, \eqref{eq:chiave3}, and Lemma \ref{lem:dagstar}(i) into \eqref{eq:5minuti}.

For part (ii), from \citet[Theorem 3]{MBPCAQML}
\begin{align}
&n^{-1/2}\Vert\wh{\bm\Lambda}_n^\dag-\wh{\bm\Lambda}_n^{(0)} \Vert = O_p(n^{-1}),\label{eq:chiave2vec}
\end{align}
while from \citet[Corollary 1 and Proposition B.3]{MBPCAQML} 
\begin{align}\label{eq:chiave3vec}
&n^{-1/2}\Vert\wh{\bm\Lambda}_n^{(0)}-{\bm\Lambda}_n^{\text{\tiny OLS}}\Vert = O_p(\max(n^{-1},n^{-1/2}T^{-1/2},T^{-1}) ).
\end{align}
Then, part (ii) is proved analogously to part (i), using \eqref{eq:chiave2vec}, \eqref{eq:chiave3vec}, and Lemma \ref{lem:dagstar}(ii). 

For part (iii), let $\wh{\sigma}_i^{2\ddag}=T^{-1}\sum_{t=1}^T (x_{it}-\wh{\bm\lambda}_i^{\dag\prime}\mbf F_t)^2$. Then, consider
\al{
\vert \wh{\sigma}_i^{2*}-\sigma_i^{2\text{\tiny OLS}}\vert \le \vert \wh{\sigma}_i^{2*}-\wh{\sigma}_i^{2\dag}\vert+\vert \wh{\sigma}_i^{2\dag}-\wh{\sigma}_i^{2\ddag}\vert
+\vert \wh{\sigma}_i^{2\ddag}-\sigma_i^{2\text{\tiny OLS}}\vert.\label{eq:sigmastarOLS}
}
From \citet[Theorem S.2  and eq. (S.33) in the online supplement]{baili16}  and noticing that the estimator of the idiosyncratic variances is unaffected by the chosen identifying constraints, we have
\al{
\vert \wh{\sigma}_i^{2\dag}-\wh{\sigma}_i^{2\ddag}\vert=O_p(\max(n^{-1},n^{-1/2}T^{-1/2},T^{-1})).\label{eq:sigmadagddag}
}
Moreover,
\al{
\vert \wh{\sigma}_i^{2\ddag}-\sigma_i^{2\text{\tiny OLS}}\vert = &\,
\l\vert T^{-1}\sum_{t=1}^T (x_{it}-\wh{\bm\lambda}_i^{\dag\prime}\mbf F_t)^2-
T^{-1}\sum_{t=1}^T (x_{it}-{\bm\lambda}_i^{\text{\tiny OLS}\prime}\mbf F_t)^2\r\vert\nn\\
\le&\, \l\vert T^{-1}\sum_{t=1}^T \l\{\wh{\bm\lambda}_i^{\dag\prime}\mbf F_t\mbf F_t^\prime\wh{\bm\lambda}_i^{\dag}
-{\bm\lambda}_i^{\text{\tiny OLS}\prime}\mbf F_t\mbf F_t^\prime{\bm\lambda}_i^{\text{\tiny OLS}}\r\}\r\vert+2\l\vert T^{-1}\sum_{t=1}^T  \l\{\wh{\bm\lambda}_i^{\dag\prime}\mbf F_t x_{it}-{\bm\lambda}_i^{\text{\tiny OLS}\prime}\mbf F_tx_{it} \r\}\r\vert\nn\\
\le&\, \Vert \wh{\bm\lambda}_i^{\dag}-{\bm\lambda}_i^{\text{\tiny OLS}}\Vert^2 \, \l\Vert T^{-1}\sum_{t=1}^T \mbf F_t\mbf F_{t}^\prime \r\Vert
+ 2 \Vert \wh{\bm\lambda}_i^{\dag}-{\bm\lambda}_i^{\text{\tiny OLS}}\Vert \, \l\Vert T^{-1}\sum_{t=1}^T \mbf F_t\mbf F_{t}^\prime \r\Vert\,\Vert {\bm\lambda}_i^{\text{\tiny OLS}}\Vert\nn\\
&+2\Vert \wh{\bm\lambda}_i^{\dag}-{\bm\lambda}_i^{\text{\tiny OLS}}\Vert \,\l\Vert T^{-1}\sum_{t=1}^T \mbf F_tx_{it} \r\Vert\nn\\
=&\, O_p(\max(n^{-1},n^{-1/2}T^{-1/2},T^{-1})),\label{eq:sigmaddagols}
}
by \eqref{eq:chiave2}, \eqref{eq:chiave3}, Lemma \ref{lem:consistCOV}(i)  combined with
Assumption \ref{ass:ident}(b), Lemma \ref{lem:consistCOV}(ii), and since $\Vert {\bm\lambda}_i^{\text{\tiny OLS}}\Vert\le \Vert {\bm\lambda}_i^{\text{\tiny OLS}}-\bm\lambda_i\Vert+\Vert {\bm\lambda}_i\Vert=O_p(1)$ by  Assumption \ref{ass:common}(a) and Lemma \ref{lem:olsT12}(i). Part (iii) then follows by substituting \eqref{eq:sigmadagddag}, \eqref{eq:sigmaddagols}, and Lemma \ref{lem:dagstar}(iii) into \eqref{eq:sigmastarOLS}. 
 This completes the proof. $\Box$

\begin{lem}\label{lem:starolsVAR} 
Let $\wh{\bm\theta}^*=(\text{\upshape vec}(\wh{\mbf A}^*)^\prime\;\text{\upshape vech}(\wh{\bm\Gamma}^{v*})^\prime)^\prime$ be the QML estimator of the entries of ${\bm\theta}$ maximizing
 $\ell(\bm X_{nT};\underline{\bm \phi}_n,\underline{\bm\theta})$ defined in \eqref{eq:LL0true}, then,
under Assumptions \ref{ass:common}, \ref{ass:idio}, \ref{ass:ind}, \ref{ass:linear}, \ref{ass:tails}, and \ref{ass:ident}, as $n,T\to\infty$,
\begin{compactenum}
\item [(i)] $n\log^{-2/\delta_v} T\,\Vert \wh{\mbf A}^*-\mbf A^{\text{\tiny \upshape OLS}}\Vert = O_p(1)$;
\item [(ii)] $n\log^{-2/\delta_v} T\,\Vert \wh{\bm \Gamma}^{v*}-\bm\Gamma^{v\text{\tiny \upshape OLS}}\Vert = O_p(1)$.
\end{compactenum}
\end{lem}

\noindent\textsc{Proof.} Throughout we consider generic values of the parameters such that $\underline{\bm\varphi}_n\in\mathcal O_n$ where $\mathcal O_n=\{\mathcal O_{\lambda_i}^n\cap \mathcal E_{\Lambda_n}\}\times\{\mathcal O_{\sigma_i^2}^n\cap\mathcal E_{\Gamma^\xi_n}\} \times \mathcal O_{\mathcal A}\times \mathcal O_{\Gamma^v}$ as defined in Section \ref{enzuccio}. Thus the elements of  $\underline{\bm\varphi}_n$ satisfy Assumptions \ref{ass:common}(a), \ref{ass:common}(d), \ref{ass:common}(e), \ref{ass:idio}(a), \ref{ass:idio}(b), and \ref{ass:idio}(f).

The log-likelihood  depends on $\underline{\mbf A}$ and $\underline{\bm\Gamma}^v$ only through $\ell(\bm F_T;\underline{\bm\theta})$  and $\ell(\bm F_{T}|\bm X_{nT};\underline{\bm\phi}_n,\underline{\bm\theta})$. Let us consider both log-likelihoods separately. First, 
since by Assumption \ref{ass:common}(i) $\mbf F_0=\mbf 0_r$, we have
\begin{align}
\ell(\bm F_T;\underline{\bm\theta}) &= -\frac{T}{2}\log\det(\underline{\bm\Gamma}^v)-\frac 12\sum_{t=1}^T(\mbf F_t-\underline{\mbf A}\mbf F_{t-1})^\prime(\underline{\bm\Gamma}^v)^{-1}(\mbf F_t-\underline{\mbf A}\mbf F_{t-1}),\nn
\end{align}
which is clearly maximized by ${\bm\theta}^{\text{\tiny \upshape OLS}}$.

Second, (see also \eqref{giannone}, \eqref{doz}, and \eqref{eq:F0P0} in the proof of Lemma \ref{lem:2max})
\al{
\ell(\bm F_{T}|\bm X_{nT};\underline{\bm\phi}_n,\underline{\bm\theta})=&\,-\frac 12\sum_{t=1}^T\log \det(\underline{\bm{\mathcal P}}_{0,t|T})-\frac 12\sum_{t=1}^T(\mbf F_t-\underline{\mbf F}_{0,t|T})^\prime
(\underline{\bm{\mathcal P}}_{0,t|T})^{-1}(\mbf F_t-\underline{\mbf F}_{0,t|T}),\nn
}
Now, by \eqref{eq:fufix}, \eqref{eq:fufix3}, and \eqref{altre22} in the proof of Lemma \ref{lem:2max}, we have that
\al{
\sup_{\underline{\bm\varphi}_n\in\mathcal O_n}\bigg\vert \ell(\bm F_{T}|\bm X_{nT};\underline{\bm\phi}_n,\underline{\bm\theta})-\ell_0(\bm F_{T}|\bm X_{nT};\underline{\bm\phi}_n)\bigg\vert = O_p(n^{-1}\log^{2/\delta_v} T).\label{eq:dgrculo}
}
The proof of parts (i) and (ii) follows from \eqref{eq:dgrculo} and by continuity of the log-likelihood and since $\ell_0(\bm F_{T}|\bm X_{nT};\underline{\bm\phi}_n)$ does not depend on $\underline{\bm\theta}$. This completes the proof. $\Box$

\begin{lem}\label{lem:starT12} 
Under Assumptions \ref{ass:common}, \ref{ass:idio}, \ref{ass:ind}, \ref{ass:linear}, \ref{ass:tails}, and \ref{ass:ident}, as $n,T\to\infty$,
\begin{compactenum}
\item [(i)] 
$\min(n \log^{-2/\delta_v} T, \sqrt T)\, \Vert \wh{\bm\lambda}_i^{*} -\bm\lambda_i\Vert=O_p(1)$, uniformly in $i$;
\item [(ii)] 
$\min(n\log^{-2/\delta_v} T, \sqrt{T})\, n^{-1/2} \Vert \wh{\bm\Lambda}_n^{*} -\bm\Lambda_n\Vert= O_p(1)$;
\item [(iii)] $\min(n\log^{-2/\delta_v} T, \sqrt{T})\,\vert \wh{\sigma}_i^{2*}-\sigma_i^2\vert=O_p(1)$, uniformly in $i$;
\item [(iv)] $\min(n\log^{-2/\delta_v} T,\sqrt T)\,\Vert \wh{\mbf A}^*-\mbf A\Vert = O_p(1)$;
\item [(v)] $\min(n\log^{-2/\delta_v} T,\sqrt T)\,\Vert \wh{\bm \Gamma}^{v*}-\bm\Gamma^v\Vert = O_p(1)$.
\end{compactenum}
\end{lem}

\noindent\textsc{Proof.} The proof follows directly from Lemmas \ref{lem:olsT12}, \ref{lem:starols}, and \ref{lem:starolsVAR}. $\Box$

\begin{lem}\label{lem:sigmaunif}
Under Assumptions \ref{ass:common}, \ref{ass:idio}, \ref{ass:ind}, \ref{ass:linear}, \ref{ass:tails}, and \ref{ass:ident}, as $n,T\to\infty$, 
\begin{compactenum}
\item [(i)] $\min(\sqrt T \log^{-1/2} n, n\log^{-2/\delta_v} T) \,\max_{i=1,\ldots,n} \Vert \wh{\bm\lambda}_i^*-\bm\lambda_i\Vert=O_p(1)$;
\item [(ii)] $\min(\sqrt T \log^{-1/2} n, n\log^{-2/\delta_v} T) \,\max_{i=1,\ldots,n} \vert \wh{\sigma}_i^{2*}-\sigma_i^2\vert=O_p(1)$;
\item [(iii)] $\min(\sqrt T \log^{-1/2} n, n\log^{-2/\delta_v} T) \, \Vert \wh{\bm\Sigma}_n^{\xi*}-\bm\Sigma_n^\xi\Vert = O_p(1)$;
\item [(iv)] $\Vert (\wh{\bm\Sigma}_n^{\xi*})^{-1}\Vert = O_p(1)$;
\item [(v)] $\min(\sqrt T \log^{-1/2} n, n\log^{-2/\delta_v} T)\, \Vert (\wh{\bm\Sigma}_n^{\xi*})^{-1}- ({\bm\Sigma}_n^{\xi})^{-1}\Vert = O_p(1)$.
\end{compactenum}
\end{lem}

\noindent
\noindent\textsc{Proof.} For part (i) consider
\al{
\max_{i=1,\ldots, n}\Vert \wh{\bm\lambda}_i^{*} -\bm\lambda_i\Vert \le &\, 
\max_{i=1,\ldots, n}\Vert \wh{\bm\lambda}_i^{*} -\wh{\bm\lambda}_i^{(0)}\Vert+\max_{i=1,\ldots, n}\Vert \wh{\bm\lambda}_i^{(0)} -\bm\lambda_i\Vert.\label{eq:unifload0}
}
First, note that  \eqref{eq:chiave2} in the proof of Lemma \ref{lem:starols} holds for all $i$, this is seen from the proof of \citet[Theorem 3]{MBPCAQML}, thus, jointly with Lemma \ref{lem:dagstar}(i), 
\al{
\max_{i=1,\ldots, n} \Vert \wh{\bm\lambda}_i^*-\wh{\bm\lambda}_i^{(0)}\Vert&\le 
\max_{i=1,\ldots, n} \Vert \wh{\bm\lambda}_i^*-\wh{\bm\lambda}_i^{\dag}\Vert+
\max_{i=1,\ldots, n} \Vert \wh{\bm\lambda}_i^{\dag}-\wh{\bm\lambda}_i^{(0)}\Vert=O_p(n^{-1}\log^{2/\delta_v} T).\label{eq:unifload2}
}
Second, from \citet[equation (A.5) in the proof of Theorem 1 in the supplementary material ]{MBPCAQML}, when imposing Assumption \ref{ass:ident}(b), 
\al{
 \wh{\bm\lambda}_i^{(0)}-\bm\lambda_i=&\, \l\{(nT)^{-1}\bm\lambda_i^\prime\sum_{t=1}^T\sum_{j=1}^n\mbf F_t\xi_{jt}\bm\lambda_j^\prime  n (\wh{\bm\Lambda}_n^{(0)\prime}\wh{\bm\Lambda}_n^{(0)})^{-1} \r\}
 +\l\{T^{-1}\sum_{t=1}^T \mbf F_{t}\xi_{it} (\bm\Lambda_n^\prime\bm\Lambda_n) (\wh{\bm\Lambda}_n^{(0)\prime}\wh{\bm\Lambda}_n^{(0)})^{-1}\r\}\nn\\
 &+\l\{(nT)^{-1}\sum_{t=1}^T\sum_{j=1}^n\xi_{it}\xi_{jt}\bm\lambda_j^\prime n (\wh{\bm\Lambda}_n^{(0)\prime}\wh{\bm\Lambda}_n^{(0)})^{-1}
 \r\}\nn\\
& +\l\{(nT)^{-1}\bm\lambda_i^\prime\sum_{t=1}^T\sum_{j=1}^n\mbf F_t\xi_{jt}(\wh{\bm\lambda}_j^{(0)}-\bm\lambda_j )^\prime  n (\wh{\bm\Lambda}_n^{(0)\prime}\wh{\bm\Lambda}_n^{(0)})^{-1}
 \r\}\nn\\
 &+\l\{
 T^{-1}\sum_{t=1}^T \mbf F_{t}\xi_{it} \bm\Lambda_n^\prime(\wh{\bm\Lambda}_n^{(0)}-\bm\Lambda_n)(\wh{\bm\Lambda}_n^{(0)\prime}\wh{\bm\Lambda}_n^{(0)})^{-1}
 \r\}\nn\\
 &+\l\{
 (nT)^{-1}\sum_{t=1}^T\sum_{j=1}^n\xi_{it}\xi_{jt}(\wh{\bm\lambda}_j^{(0)}-\bm\lambda_j )^\prime n(\wh{\bm\Lambda}_n^{(0)\prime}\wh{\bm\Lambda}_n^{(0)})^{-1}
 \r\}\nn\\
=&\, 1.a+1.b+1.c+1.d+1.e+1.f, \;\text{ say,}\nn
}
where 
\al{
&\max_{i=1,\ldots,n}\Vert 1.a\Vert  = O_p(n^{-1/2}T^{-1/2}),
\nn\\
&\max_{i=1,\ldots,n}\Vert 1.b\Vert  =O_p(T^{-1/2}\sqrt{\log n}),\label{eq:cambioposto}\\
&\max_{i=1,\ldots,n}\Vert 1.c\Vert =O_p(\max(n^{-1}, n^{-1/2}T^{-1/2}\sqrt{\log n})),\nn
\nn
}
which follows by using
\bit
\item for $1.a$ Assumption \ref{ass:common}(a), Lemma \ref{lem:est0LOAD}(ii), and \citet[Proposition B.3(a)]{MBPCAQML}; 
\item for $1.b$ Assumption \ref{ass:common}(a), Lemma \ref{lem:est0LOAD}(ii), and Lemma \ref{lem:tail}(ii) and the union bound;
\item for $1.c$ Assumption \ref{ass:common}(a), Lemma \ref{lem:est0LOAD}(ii), \ref{lem:tail}(ii) and the union bound, and \citet[Proposition B.3(c)]{MBPCAQML}, where we also used the fact that $\max_{i=1,\ldots, n} \vert (nT)^{-1}\sum_{t=1}^T\sum_{j=1}^n \E[\xi_{it}\xi_{jt}]\vert = O(n^{-1})$ by Assumption \ref{ass:idio}(b); 
\eit
while $1.d$, $1.e$, and $1.f$  are dominated by and $1.a$, $1.b$, and $1.c$, respectively, because of Lemma \ref{lem:est0LOAD}(i).
It follows that
\al{
\max_{i=1,\ldots,n}\Vert \wh{\bm\lambda}_i^{(0)}-\bm\lambda_i\Vert = O_p(\max(n^{-1},T^{-1/2}\sqrt{\log n}))\label{eq:unifload3}.
}
By substituting  \eqref{eq:unifload2} and  \eqref{eq:unifload3} into  \eqref{eq:unifload0}, we have
\begin{align}
\max_{i=1,\ldots, n} \Vert \wh{\bm\lambda}_i^*-\bm\lambda_i\Vert
&= O_p(\max(n^{-1},T^{-1/2}\sqrt{\log n}))+O_p(n^{-1}\log^{2/\delta_v} T).\nn
\end{align}
This proves part (i).

For part (ii), first,  consider 
\al{
\max_{i=1,\ldots,n} \vert \wh{\sigma}_i^{2*}-\sigma_i^{2} \vert \le&\, \max_{i=1,\ldots,n}  \vert \wh{\sigma}_i^{2*}-\wh{\sigma}_i^{2\dag}\vert
+\max_{i=1,\ldots,n} \vert \wh{\sigma}_i^{2\dag}-\wh{\sigma}_i^{2\ddag}\vert\nn\\
&+\max_{i=1,\ldots,n} \vert \wh{\sigma}_i^{2\ddag}-\sigma_i^{2\text{\tiny OLS}}\vert+\max_{i=1,\ldots,n} \vert \sigma_i^{2\text{\tiny OLS}}-\sigma_i^2\vert.
\label{eq:italoc}
}
From \eqref{eq:sigmadagddag} in the proof of Lemma \ref{lem:starols}
\beq
\max_{i=1,\ldots,n} \vert \wh{\sigma}_i^{2\dag}-\wh{\sigma}_i^{2\ddag}\vert = O_p(\max(n^{-1},T^{-1}\sqrt{\log n},n^{-1/2}T^{-1/2} \sqrt{\log n})). \label{eq:bidi}
\eeq
Indeed from \citet[equation (S.25) in the online supplement]{baili16} we see that \eqref{eq:bidi} is decomposed into the sum of 13 terms, and all depend on $i$ only through $\bm\lambda_i$ which is such that $\max_{i=1,\ldots, n} \Vert\bm\lambda_i\Vert\le M_\lambda$ by Assumption \ref{ass:common}(a), with the exceptions of the following terms
\al{
A_1&=\max_{i=1,\ldots, n} \l\vert (\wh{\bm\lambda}_i^{\dag}-\bm\lambda_i)^\prime \l\{T^{-1}\sum_{t=1}^T \mbf F_t\mbf F_t^\prime\r\}(\wh{\bm\lambda}_i^{\dag}-\bm\lambda_i) \r\vert=O_p(\max(T^{-1} \log n, n^{-2}\log^{4/\delta_v} T)),\nn\\
A_2&=\max_{i=1,\ldots, n} 2\l\vert \bm\lambda_i^\prime \wh{\bm\Lambda}_n^{\dag\prime}(\wh{\bm\Sigma}_n^{\xi\dag})^{-1}(\wh{\bm\Lambda}_n^{\dag}-\bm\Lambda_n)T^{-1}\sum_{t=1}^T \mbf F_t\xi_{it}\r\vert=O_p((n^{-1/2}T^{-1/2} \sqrt{\log n}),\nn\\
A_3&=\max_{i=1,\ldots, n} 2\l\vert \bm\lambda_i^\prime \wh{\bm\Lambda}_n^{\dag\prime}(\wh{\bm\Sigma}_n^{\xi\dag})^{-1}\l\{T^{-1}\sum_{t=1}^T \bm\xi_{nt}\r\}\l\{T^{-1}\sum_{t=1}^T \xi_{it}\r\}
\r\vert=O_p(T^{-1} \sqrt{\log n}),\nn
}
where we used: 
\bit
\item for $A_1$ part (i) and Lemmas \ref{lem:dagstar}(i) and \ref{lem:consistCOV}(i) jointly with Assumption \ref{ass:ident}(b);
\item for $A_2$ Lemma \ref{lem:tail}(ii) and the union bound,
plus \citet[Lemma B.4 and Corollary A.1 in the online supplement]{baili12};
\item for $A_3$ the same arguments used for $A_2$ plus \citet[Lemma S.10 in the online supplement]{baili16}.
\eit

Second, notice that, by Lemmas \ref{lem:consistCOV}(i) and \ref{lem:est0LOAD}(ii),  term $1.b$ in \eqref{eq:cambioposto} is such that
\al{
 1.b  &= ({\bm\lambda}_i^{\text{\tiny OLS}}-\bm\lambda_i )+ o_p(T^{-1/2}\sqrt{\log n}).\label{eq:bosco}
}
Thus, from the above arguments and \eqref{eq:bosco},
\al{
\max_{i=1,\ldots, n} \Vert \wh{\bm\lambda}_i^{(0)}-{\bm\lambda}_i^{\text{\tiny OLS}} \Vert = O_p(\max(n^{-1},n^{-1/2}T^{-1/2}\sqrt{\log n})).\label{eq:bosco2}
}
and, therefore, from \eqref{eq:unifload2} and \eqref{eq:bosco2}
\al{
\max_{i=1,\ldots, n} \Vert \wh{\bm\lambda}_i^{\dag}-{\bm\lambda}_i^{\text{\tiny OLS}} \Vert \le&\, \max_{i=1,\ldots, n} \Vert \wh{\bm\lambda}_i^{\dag}-\wh{\bm\lambda}_i^{(0)} \Vert 
+\max_{i=1,\ldots, n} \Vert \wh{\bm\lambda}_i^{(0)}-{\bm\lambda}_i^{\text{\tiny OLS}} \Vert\nn\\
 =&\, O_p(\max(n^{-1}\log^{2/\delta_v}T,n^{-1/2}T^{-1/2}\sqrt{\log n})).\label{eq:bosco3}
}

Then, from \eqref{eq:sigmaddagols} in the proof of Lemma \ref{lem:starols} 
\al{
\max_{i=1,\ldots,n} \vert \wh{\sigma}_i^{2\ddag}-\sigma_i^{2\text{\tiny OLS}}\vert\le&\,
\max_{i=1,\ldots,n} \Vert \wh{\bm\lambda}_i^{\dag}-{\bm\lambda}_i^{\text{\tiny OLS}}\Vert^2 \, \l\Vert T^{-1}\sum_{t=1}^T \mbf F_t\mbf F_{t}^\prime \r\Vert\nn\\
&+ 2 \max_{i=1,\ldots,n} \Vert \wh{\bm\lambda}_i^{\dag}-{\bm\lambda}_i^{\text{\tiny OLS}}\Vert \, \l\Vert T^{-1}\sum_{t=1}^T \mbf F_t\mbf F_{t}^\prime \r\Vert\, \max_{i=1,\ldots,n}\Vert {\bm\lambda}_i^{\text{\tiny OLS}}\Vert\nn\\
&+2\max_{i=1,\ldots,n} \Vert \wh{\bm\lambda}_i^{\dag}-{\bm\lambda}_i^{\text{\tiny OLS}}\Vert \,\l\Vert T^{-1}\sum_{t=1}^T \mbf F_t\mbf F_{t}^\prime \r\Vert\, \max_{i=1,\ldots,n}\Vert\bm\lambda_i\Vert\nn\\
&+2\max_{i=1,\ldots,n} \Vert \wh{\bm\lambda}_i^{\dag}-{\bm\lambda}_i^{\text{\tiny OLS}}\Vert \,\max_{i=1,\ldots,n}\l\Vert T^{-1}\sum_{t=1}^T \mbf F_t\xi_{it} \r\Vert\nn\\
=&\, O_p(\max(n^{-1}\log^{2/\delta_v}T,n^{-1/2}T^{-1/2}\sqrt{\log n})),\label{eq:italoc2}
}
by \eqref{eq:cambioposto}, \eqref{eq:bosco3},  Assumption \ref{ass:common}(a), Lemma \ref{lem:consistCOV}(i) combined with Assumption \ref{ass:ident}(b), and  since 
$\max_{i=1,\ldots,n}\Vert {\bm\lambda}_i^{\text{\tiny OLS}}\Vert\le \max_{i=1,\ldots,n}\Vert {\bm\lambda}_i^{\text{\tiny OLS}}-\bm\lambda_i\Vert+
\max_{i=1,\ldots,n}\Vert {\bm\lambda}_i\Vert= O_p(T^{-1/2}\sqrt{\log n})+O(1)$ again by \eqref{eq:cambioposto} and Assumption \ref{ass:common}(a).

Finally,
\al{
\max_{i=1,\ldots,n} \vert \sigma_i^{2\text{\tiny OLS}}-\sigma_i^2\vert\le &\,
\max_{i=1,\ldots,n}\l\vert T^{-1}\sum_{t=1}^T x_{it}^2-\E[x_{it}^2]\r\vert
+ \max_{i=1,\ldots,n} \l\vert \bm\lambda_i^{\text{\tiny OLS}\prime}T^{-1}\sum_{t=1}^T\mbf F_t\mbf F_t^\prime \bm\lambda_i^{\text{\tiny OLS}}-\bm\lambda_i^\prime\bm\lambda_i
\r\vert\nn\\
&+2 \max_{i=1,\ldots,n} \l\Vert\bm\lambda_i^{\text{\tiny OLS}\prime}T^{-1}\sum_{t=1}^T\mbf F_t\mbf F_t^\prime -\bm\lambda_i^\prime\r\Vert\,
 \max_{i=1,\ldots,n}\Vert \bm\lambda_i\Vert\nn\\
&+2 \max_{i=1,\ldots,n} \Vert \bm\lambda_i^{\text{\tiny OLS}}\Vert \max_{i=1,\ldots,n} \l\Vert T^{-1}\sum_{t=1}^T\mbf F_t\xi_{it} \r\Vert\nn\\
=&\, O_p(T^{-1/2}\sqrt{\log n}).
\label{eq:italomerda}
}
Indeed, for the first term on the rhs of \eqref{eq:italomerda} we have
\al{
\max_{i=1,\ldots,n}\l\vert T^{-1}\sum_{t=1}^T x_{it}^2-\E[x_{it}^2]\r\vert\le&\, 
\max_{i=1,\ldots,n} \Vert\bm\lambda_i\Vert^2 \l\Vert T^{-1}\sum_{t=1}^T \mbf F_t\mbf F_t^\prime-\mbf I_r\r\Vert+
\max_{i=1,\ldots,n}\l\vert T^{-1}\sum_{t=1}^T \xi_{it}^2-\E[\xi_{it}^2]\r\vert\nn\\
&+
2\max_{i=1,\ldots,n}\Vert\bm\lambda_i\Vert \max_{i=1,\ldots,n} \l\Vert T^{-1}\sum_{t=1}^T \mbf F_t\xi_{it}\r\Vert\nn\\
=&\,O_p(T^{-1/2}\sqrt{\log n}),\label{eq:MINKIEVA}
}
by Assumption \ref{ass:common}(a), Lemma \ref{lem:consistCOV}(i) for the first term, \citet[Lemma 3.ii]{BCF18} and \citet[Lemmas A3 and B1]{FLM11}, which we can apply because of Assumption \ref{ass:tails}(b) for the second term, and Lemma \ref{lem:tail}(ii) and the union bound for the third term. For all other on the rhs of \eqref{eq:italomerda} we just need to use  Assumption \ref{ass:common}(a), Lemma \ref{lem:consistCOV}(i) and the fact that 
$\max_{i=1,\ldots,n}\Vert {\bm\lambda}_i^{\text{\tiny OLS}}\Vert\le \max_{i=1,\ldots,n}\Vert {\bm\lambda}_i^{\text{\tiny OLS}}-\bm\lambda_i\Vert+
\max_{i=1,\ldots,n}\Vert {\bm\lambda}_i\Vert= O_p(T^{-1/2}\sqrt{\log n})+O(1)$ again by \eqref{eq:cambioposto} and Assumption \ref{ass:common}(a).

By using Lemma \ref{lem:dagstar}(iii), \eqref{eq:bidi}, \eqref{eq:italoc2}, and \eqref{eq:italomerda} into \eqref{eq:italoc} we have
\al{
\max_{i=1,\ldots,n} \vert \wh{\sigma}_i^{2*}-\sigma_i^{2} \vert=&\,O_p(n^{-1}\log^{2/\delta_v}T)+O_p(\max(n^{-1},T^{-1}\sqrt{\log n},n^{-1/2}T^{-1/2} \sqrt{\log n}))\nn\\
&+O_p(\max(n^{-1}\log^{2/\delta_v}T,n^{-1/2}T^{-1/2}\sqrt{\log n}))+O_p(T^{-1/2}\sqrt{\log n})\nn\\
=&\, O_p(\max(n^{-1}\log^{2/\delta_v}T, T^{-1/2}\sqrt{\log n})),\nn
}
which proves part (ii).

Part (iii) immediately follows from part (ii), indeed
\al{
\Vert \wh{\bm\Sigma}_n^{\xi*}-\bm\Sigma_n^\xi\Vert \le \max_{i=1,\ldots, n} \vert  \wh{\sigma}_i^{2*}-\sigma_i^{2} \vert=O_p(\max(n^{-1}\log^{2/\delta_v}T, T^{-1/2}\sqrt{\log n})).\nn
}

For part (iv) we have
\al{
\Vert (\wh{\bm\Sigma}_n^{\xi*})^{-1}\Vert = \l\{\min_{i=1,\ldots, n} \wh{\sigma}_i^{2*}\r\}^{-1}\le C_\xi +O_p(\max(n^{-1}\log^{2/\delta_v}T, T^{-1/2}\sqrt{\log n})),\nn
}
because of part (ii) and Assumption \ref{ass:idio}(a).

To conclude, for part (v) we have
\al{
\Vert (\wh{\bm\Sigma}_n^{\xi*})^{-1}-(\bm\Sigma_n^\xi)^{-1}\Vert \le \Vert (\wh{\bm\Sigma}_n^{\xi*})^{-1}\Vert \, \Vert \wh{\bm\Sigma}_n^{\xi*}-\bm\Sigma_n^\xi\Vert\,
\Vert ({\bm\Sigma}_n^{\xi})^{-1}\Vert = O_p(\max(n^{-1}\log^{2/\delta_v}T, T^{-1/2}\sqrt{\log n})),\nn
}
by parts (iii), (iv), and Assumption \ref{ass:idio}(a). This completes the proof. $\Box$

\begin{lem}\label{lem:eststar_LAST}
Under Assumptions \ref{ass:common}, \ref{ass:idio}, \ref{ass:ind}, \ref{ass:linear}, \ref{ass:tails}, and \ref{ass:ident}, as $n,T\to\infty$:
\begin{compactenum}
\item [(i)] $\min(n\log^{-2/\delta_v}T,\sqrt T \log^{-1/2}n
)\,n^{-1}\Vert\wh{\bm\Lambda}_n^{*\prime}(\wh{\bm \Sigma}_n^{\xi*})^{-1}\wh{\bm\Lambda}_n^{*}-\bm\Lambda_n^\prime(\bm\Sigma_n^\xi)^{-1}\bm\Lambda_n\Vert = O_p(1)$;
\item [(ii)] $\min(n\log^{-2/\delta_v}T,\sqrt T \log^{-1/2}n
)\,n^{-1/2}\Vert\wh{\bm\Lambda}_n^{*\prime}(\wh{\bm \Sigma}_n^{\xi*})^{-1}-\bm\Lambda_n^\prime(\bm\Sigma_n^\xi)^{-1}\Vert = O_p(1)$;
\item [(iii)] $n\Vert(\wh{\bm\Lambda}_n^{*\prime}(\wh{\bm \Sigma}_n^{\xi*})^{-1}\wh{\bm\Lambda}_n^{*})^{-1}\Vert = O_p(1)$;
\item [(iv)] $\min(n\log^{-2/\delta_v}T,\sqrt T \log^{-1/2}n
)\,n\Vert(\wh{\bm\Lambda}_n^{*\prime}(\wh{\bm \Sigma}_n^{\xi*})^{-1}\wh{\bm\Lambda}_n^{*})^{-1}-(\bm\Lambda_n^\prime(\bm\Sigma_n^\xi)^{-1}\bm\Lambda_n)^{-1}\Vert = O_p(1)$;
\item [(v)] $\omega_{n,T\delta_v}\,\sqrt n\Vert(\wh{\bm\Lambda}_n^{*\prime}(\wh{\bm \Sigma}_n^{\xi*})^{-1}\wh{\bm\Lambda}_n^{*})^{-1}
\wh{\bm\Lambda}_n^{*\prime}(\wh{\bm \Sigma}_n^{\xi*})^{-1}
-(\bm\Lambda_n^\prime(\bm\Sigma_n^\xi)^{-1}\bm\Lambda_n)^{-1}\bm\Lambda_n^\prime(\bm\Sigma_n^\xi)^{-1}\Vert = O_p(1)$,\\ with $\omega_{n,T\delta_v}=\min(n\log^{-2/\delta_v}T,\sqrt T \log^{-1/2}n
).$
\end{compactenum}
\end{lem}

\noindent\textsc{Proof.} 
For part (i) we have
\al{
n^{-1}\Vert\wh{\bm\Lambda}_n^{*\prime}(\wh{\bm \Sigma}_n^{\xi*})^{-1}\wh{\bm\Lambda}_n^{*}-\bm\Lambda_n^\prime(\bm\Sigma_n^\xi)^{-1}\bm\Lambda_n\Vert 
\le &\,
2 n^{-1}\Vert \{\wh{\bm\Lambda}_n^{*}-\bm\Lambda_n\}^\prime({\bm\Sigma}_n^\xi)^{-1}{\bm\Lambda}_n\Vert\nn\\
  &+n^{-1}\Vert  {\bm\Lambda}_n^\prime \{(\wh{\bm\Sigma}_n^{\xi*})^{-1} -({\bm\Sigma}_n^\xi)^{-1} \} {\bm\Lambda}_n\Vert\nn\\
 &+2n^{-1}\Vert
 \{\wh{\bm\Lambda}_n^{*}-\bm\Lambda_n\}^\prime \{(\wh{\bm\Sigma}_n^{\xi*})^{-1} -({\bm\Sigma}_n^\xi)^{-1} \} {\bm\Lambda}_n
 \Vert\nn\\
  &+n^{-1}\Vert
 \{\wh{\bm\Lambda}_n^{*}-\bm\Lambda_n\}^\prime \{(\wh{\bm\Sigma}_n^{\xi*})^{-1} -({\bm\Sigma}_n^\xi)^{-1} \} \{\wh{\bm\Lambda}_n^{*}-\bm\Lambda_n\}\nn
 \Vert\nn\\
 \le &\, 2 n^{-1/2}\Vert \wh{\bm\Lambda}_n^{*}-\bm\Lambda_n\Vert \, \Vert({\bm\Sigma}_n^\xi)^{-1}\Vert \, n^{-1/2}\Vert {\bm\Lambda}_n\Vert\nn\\
  &+\Vert \{(\wh{\bm\Sigma}_n^{\xi*})^{-1} -({\bm\Sigma}_n^\xi)^{-1} \} \Vert \, n^{-1}\Vert {\bm\Lambda}_n\Vert^2\nn\\
&+2n^{-1/2}\Vert
 \wh{\bm\Lambda}_n^{*}-\bm\Lambda_n \Vert \Vert \{(\wh{\bm\Sigma}_n^{\xi*})^{-1} -({\bm\Sigma}_n^\xi)^{-1} \}\Vert\, n^{-1/2}\Vert {\bm\Lambda}_n
 \Vert\nn\\
  &+n^{-1}\Vert \wh{\bm\Lambda}_n^{*}- {\bm\Lambda}_n\Vert^2 \Vert (\wh{\bm\Sigma}_n^{\xi*})^{-1} -({\bm\Sigma}_n^\xi)^{-1}\Vert\nn\\
  =&\, O_p(\max(n^{-1}\log^{2/\delta_v}T,T^{-1/2}\sqrt{\log n})),
\label{eq:renollopazzo}
}
by Assumptions \ref{ass:common}(a), \ref{ass:idio}(a), and Lemmas \ref{lem:starT12}(ii) and \ref{lem:sigmaunif}(v). This proves part (i).

Part (ii) is proved in the same way as part (i).

For part (iii), by part (ii) and \citet[Theorem 1]{MK04} which is Weyl's inequality, we have
\al{
n^{-1}\vert \nu^{(r)}&
(\wh{\bm\Lambda}_n^{*\prime}(\wh{\bm \Sigma}_n^{\xi*})^{-1}\wh{\bm\Lambda}_n^{*})
-
\nu^{(r)}({\bm\Lambda}_n^{\prime}({\bm \Sigma}_n^{\xi})^{-1}){\bm\Lambda}_n)\vert\nn\\
\le&\, 
 n^{-1}\Vert \wh{\bm\Lambda}_n^{*\prime}(\wh{\bm \Sigma}_n^{\xi*})^{-1}\wh{\bm\Lambda}_n^{*}
- {\bm\Lambda}_n^{\prime}({\bm \Sigma}_n^{\xi})^{-1}{\bm\Lambda}_n\Vert \nn\\
=&\, O_p(\max(n^{-1}\log^{2/\delta_v}T,T^{-1/2}\sqrt{\log n})).\label{eq:CRUCIALstar}
}
Moreover (note that $x-y\ge -|x-y|$ for any $x,y\in\mathbb R$), 
\al{
\det ( \wh{\bm\Lambda}_n^{*\prime}(\wh{\bm \Sigma}_n^{\xi*})^{-1}\wh{\bm\Lambda}_n^{*}) =&\,\prod_{j=1}^r \nu^{(j)}
(\wh{\bm\Lambda}_n^{*\prime}(\wh{\bm \Sigma}_n^{\xi*})^{-1}\wh{\bm\Lambda}_n^{*})\ge \l\{ \nu^{(r)}
(\wh{\bm\Lambda}_n^{*\prime}(\wh{\bm \Sigma}_n^{\xi*})^{-1}\wh{\bm\Lambda}_n^{*}) \r\}^r\nn\\
\ge&\, \l\{ \nu^{(r)}
({\bm\Lambda}_n^{\prime}({\bm \Sigma}_n^{\xi})^{-1}{\bm\Lambda}_n) -
\vert\nu^{(r)}
(\wh{\bm\Lambda}_n^{*\prime}(\wh{\bm \Sigma}_n^{\xi*})^{-1}\wh{\bm\Lambda}_n^{*})
-
\nu^{(r)}({\bm\Lambda}_n^{\prime}({\bm \Sigma}_n^{\xi})^{-1}){\bm\Lambda}_n)\vert \r\}^r,\nn
}
thus, by Lemma \ref{lem:LSL2}(iv), which implies $\lim_{n\to\infty}n^{-1}\nu^{(r)}
({\bm\Lambda}_n^{\prime}({\bm \Sigma}_n^{\xi})^{-1}{\bm\Lambda}_n)>0$,
 and \eqref{eq:CRUCIALstar} it follows that, with probability tending to one as $n,T\to\infty$, we have $\det (n^{-1} \wh{\bm\Lambda}_n^{*\prime}(\wh{\bm \Sigma}_n^{\xi*})^{-1}\wh{\bm\Lambda}_n^{*})>0$, or, equivalently $n^{-1}\wh{\bm\Lambda}_n^{*\prime}(\wh{\bm \Sigma}_n^{\xi*})^{-1}\wh{\bm\Lambda}_n^{*}$ is positive definite, 
 i.e. $n\Vert(\wh{\bm\Lambda}_n^{*\prime}(\wh{\bm \Sigma}_n^{\xi*})^{-1}\wh{\bm\Lambda}_n^{*})^{-1}\Vert=O_p(1)$.
  This proves part (iii).

For part (iv), we have
\al{
n\Vert&(\wh{\bm\Lambda}_n^{*\prime}(\wh{\bm \Sigma}_n^{\xi*})^{-1}\wh{\bm\Lambda}_n^{*})^{-1}-(\bm\Lambda_n^\prime(\bm\Sigma_n^\xi)^{-1}\bm\Lambda_n)^{-1}\Vert\nn\\
&\le
n\Vert(\wh{\bm\Lambda}_n^{*\prime}(\wh{\bm \Sigma}_n^{\xi*})^{-1}\wh{\bm\Lambda}_n^{*})^{-1}\Vert\,
n^{-1}\Vert \wh{\bm\Lambda}_n^{*\prime}(\wh{\bm \Sigma}_n^{\xi*})^{-1}\wh{\bm\Lambda}_n^{*}-\bm\Lambda_n^\prime(\bm\Sigma_n^\xi)^{-1}\bm\Lambda_n\Vert\,
n\Vert (\bm\Lambda_n^\prime(\bm\Sigma_n^\xi)^{-1}\bm\Lambda_n)^{-1}\Vert\nn\\
&=O_p(\max(n^{-1}\log^{2/\delta_v}T,T^{-1/2}\sqrt{\log n})),\nn
}
because of parts (i) and (iii) and Lemma \ref{lem:LSL2}(iii). 

Part (v) follows directly from parts (ii) and (iv). This completes the proof. $\Box$

\begin{lem}\label{lem:PPOnstar}
 Under Assumptions \ref{ass:common}, \ref{ass:idio}, \ref{ass:ind}, \ref{ass:linear}, \ref{ass:tails}, and \ref{ass:ident}, as $n,T\to\infty$:
\begin{compactenum}[(i)]
\item $\max_{t=1,\ldots, T}\Vert \mbf P^{*}_{t|t-1}\Vert=O_p(1)$; 
\item $\max_{t=1,\ldots, T}\Vert(\mbf P^{*}_{t|t-1})^{-1}\Vert =O_p(1)$;
\item $\max_{t=1,\ldots, T}n\Vert \mbf P^{*}_{t|t}\Vert=O_p(1)$;
\item $\max_{t=1,\ldots, T}n\Vert \mbf P^{*}_{t|T}\Vert=O_p(1)$.
\end{compactenum}
\end{lem}

\noindent
\textsc{Proof.} 
For part (i), 
\al{
\max_{t=1,\ldots, T}\Vert {\mbf P}^{*}_{t|t-1}\Vert&\le \max_{t=1,\ldots, T}\Vert{\mbf P}_{t|t-1}\Vert+ \max_{t=1,\ldots, T}\Vert{\mbf P}_{t|t-1}^{*}-{\mbf P}_{t|t-1}\Vert\nn\\
&= O_p(1) + O_p(\max(n^{-1}\log^{2/\delta_v} T,T^{-1/2}\sqrt{\log n})),\nn
}
by Lemma \ref{lem:cazzarola}(i) and since the second term on the rhs depends only on the estimation error of $\wh{\mbf A}^{*}$, $\wh{\bm\Gamma}^{v*}$, $n^{-1/2}\wh{\bm\Lambda}_n^{*}$, 
$n^{-1/2}\wh{\bm\Lambda}_n^{*\prime}(\wh{\bm\Sigma}_n^{\xi*})^{-1}$
and $n^{-1}(\wh{\bm\Lambda}_n^{*\prime}
(\wh{\bm\Sigma}_n^{\xi*})^{-1}
\wh{\bm\Lambda}_n^{*})^{-1}$, which are all bounded by Lemmas \ref{lem:starT12}(ii), \ref{lem:starT12}(iv), \ref{lem:starT12}(v), \ref{lem:eststar_LAST}(ii), and \ref{lem:eststar_LAST}(iv). This proves part (i).

Part (ii) is proved in the same way as part (i) but using Lemma  \ref{lem:cazzarola}(ii).

For part (iii), from \eqref{eq:up2} using the same steps leading to \eqref{PtThathat} in the proof of Lemma \ref{lem:cazzarolahat00} but when using as parameters $\wh{\bm\varphi}_n^*$, we have
\begin{align}
{\mbf P}_{t|t}^{*}=&\, {\mbf P}_{t|t-1}^{*}-{\mbf P}_{t|t-1}^{*}\wh{\bm\Lambda}_n^{*\prime}
(\wh{\bm\Lambda}_n^{*}{\mbf P}_{t|t-1}^{*}\wh{\bm\Lambda}_n^{*\prime}+\wh{\bm\Sigma}_n^{\xi*})^{-1}\wh{\bm\Lambda}_n^{*}{\mbf P}_{t|t-1}^{*}\nn\\
=&\,(\wh{\bm\Lambda}_n^{*\prime}(\wh{\bm\Sigma}_n^{\xi*})^{-1}\wh{\bm\Lambda}_n^{*})^{-1}\label{PtThathatstar}\\
&- {\mbf P}_{t|t-1}^{*}
((\wh{\bm\Lambda}_n^{*\prime}(\wh{\bm\Sigma}_n^{\xi*})^{-1}\wh{\bm\Lambda}_n^{*})^{-1}+ {\mbf P}_{t|t-1}^{*})^{-1}
(\wh{\bm\Lambda}_n^{*\prime}(\wh{\bm\Sigma}_n^{\xi*})^{-1}\wh{\bm\Lambda}_n^{*})^{-1}({\mbf P}_{t|t-1}^{*})^{-1}(\wh{\bm\Lambda}_n^{*\prime}(\wh{\bm\Sigma}_n^{\xi*})^{-1}\wh{\bm\Lambda}_n^{*})^{-1}.\nn
\end{align}
Notice that all inverses in \eqref{PtThathatstar} are well defined because of part (ii) and Lemmas \ref{lem:eststar_LAST}(iii) and  \ref{lem:sigmaunif}(iv).

Therefore, from \eqref{PtThathatstar}
\al{
\max_{t=1,\ldots, T} n\Vert {\mbf P}_{t|t}^{*}\Vert \le &\, 
 n \Vert (\wh{\bm\Lambda}_n^{*\prime}(\wh{\bm\Sigma}_n^{\xi*})^{-1}\wh{\bm\Lambda}_n^{*})^{-1}\Vert\nn\\
&+ \max_{t=1,\ldots, T} \Vert {\mbf P}_{t|t-1}^{*}\Vert \, \max_{t=1,\ldots, T} \Vert ({\mbf P}_{t|t-1}^{*})^{-1}\Vert\,
n\Vert (\wh{\bm\Lambda}_n^{*\prime}(\wh{\bm\Sigma}_n^{\xi*})^{-1}\wh{\bm\Lambda}_n^{*})^{-1}\Vert^2\nn\\
&\cdot \Vert ((\wh{\bm\Lambda}_n^{*\prime}(\wh{\bm\Sigma}_n^{\xi*})^{-1}\wh{\bm\Lambda}_n^{*})^{-1}+ {\mbf P}_{t|t-1}^{*})^{-1}\Vert\nn\\
=&\, O_p(1) + O_p(n^{-1}),\nn
}
 because of parts (i) and (ii), Lemma \ref{lem:eststar_LAST}(iii), and  since, by \citet[Theorem 1]{MK04} which is Weyl's inequality,
 \al{
 \Vert &((\wh{\bm\Lambda}_n^{*\prime}(\wh{\bm\Sigma}_n^{\xi*})^{-1}\wh{\bm\Lambda}_n^{*})^{-1}+ {\mbf P}_{t|t-1}^{*})^{-1}\Vert =
 \l\{\nu^{(r)}( (\wh{\bm\Lambda}_n^{*\prime}(\wh{\bm\Sigma}_n^{\xi*})^{-1}\wh{\bm\Lambda}_n^{*})^{-1}+ {\mbf P}_{t|t-1}^{*})  \r\}^{-1}\nn\\
 \le &\,  \l\{\nu^{(r)}( (\wh{\bm\Lambda}_n^{*\prime}(\wh{\bm\Sigma}_n^{\xi*})^{-1}\wh{\bm\Lambda}_n^{*})^{-1})+ \nu^{(r)}( {\mbf P}_{t|t-1}^{*})  \r\}^{-1}\nn\\
 = &\,  \l\{
 \l[\nu^{(1)} (\wh{\bm\Lambda}_n^{*\prime}(\wh{\bm\Sigma}_n^{\xi*})^{-1}\wh{\bm\Lambda}_n^{*})\r]^{-1}
 + \nu^{(r)}( {\mbf P}_{t|t-1}^{*})  \r\}^{-1}\nn\\
  = &\,  \l\{
 \l[\nu^{(1)} (\wh{\bm\Lambda}_n^{*\prime}(\wh{\bm\Sigma}_n^{\xi*})^{-1}\wh{\bm\Lambda}_n^{*})\nu^{(r)}( {\mbf P}_{t|t-1}^{*})\r]^{-1}
 + 1  \r\}^{-1} \l\{\nu^{(r)}( {\mbf P}_{t|t-1}^{*})\r\}^{-1}\nn\\
 =&\, \l\{1-\l[\nu^{(1)} (\wh{\bm\Lambda}_n^{*\prime}(\wh{\bm\Sigma}_n^{\xi*})^{-1}\wh{\bm\Lambda}_n^{*})\nu^{(r)}( {\mbf P}_{t|t-1}^{*})\r]^{-1}\r\}\l\{\nu^{(r)}( {\mbf P}_{t|t-1}^{*})\r\}^{-1}+O_p(n^{-2})\nn\\
 =&\, O_p(1),\nn
 }
again by parts (i) and (ii) and Lemma \ref{lem:eststar_LAST}(iii). This proves part (iii).

For part (iv), from \eqref{eq:KS2}, we get
\al{
\Vert \mbf P_{t|T}^{*}-\mbf P_{t|t}^{*}\Vert\le&\, \Vert \mbf P_{t|t}^{*}\Vert^2\, \Vert \wh{\mbf A}^{*}\Vert^2\,  
\Vert(\mbf P_{t+1|t}^{*})^{-1}\Vert^2
\{\Vert\mbf P_{t+1|T}^{*}\Vert+\Vert\mbf P_{t+1|t}^{*}\Vert\}.\label{eq:pivastar}
}
Start with $t=T-1$, then from \eqref{eq:pivastar},
\al{
\Vert \mbf P_{T-1|T}^{*}-\mbf P_{T-1|T-1}^{*}\Vert\le&\, \Vert \mbf P_{T-1|T-1}^{*}\Vert^2\, \Vert \wh{\mbf A}^{*}\Vert^2\,  
\Vert(\mbf P_{T|T-1}^{*})^{-1}\Vert^2
\{\Vert\mbf P_{T|T}^{*}\Vert+\Vert\mbf P_{T|T-1}^{*}\Vert\}\nn\\
 =&\, O_p(n^{-2}).\label{eq:pivastar2}
}
by parts (i), (ii), and (iii),  and since $\Vert\wh{\mbf A}^{*}\Vert\le \Vert{\mbf A}\Vert+\Vert\wh{\mbf A}^{*}-\mbf A\Vert =O_p(1)$, by Assumption \ref{ass:common}(d) and Lemma \ref{lem:starT12}(iv). From \eqref{eq:pivastar2} 
it follows that 
\al{
\Vert \mbf P^{*}_{T-1|T}\Vert\le \Vert \mbf P^{*}_{T-1|T-1}\Vert+\Vert \mbf P^{*}_{T-1|T}-\mbf P^{*}_{T-1|T-1}\Vert = O_p(n^{-1})+O_p(n^{-2}).\label{eq:pivastar3}
}
Thus, at $t=T-2$, from \eqref{eq:pivastar} and \eqref{eq:pivastar3}, 
\al{
\Vert \mbf P_{T-2|T}^{*}-\mbf P_{T-2|T-2}^{*}\Vert\le&\, \Vert \mbf P_{T-2|T-2}^{*}\Vert^2\, \Vert \wh{\mbf A}^{*}\Vert^2\,  
\Vert(\mbf P_{T-1|T-2}^{*})^{-1}\Vert^2
\{\Vert\mbf P_{T-1|T}^{*}\Vert+\Vert\mbf P_{T-1|T-2}^{*}\Vert\}\nn\\
 =&\, O_p(n^{-2}).\label{eq:pivastar4}
}
From \eqref{eq:pivastar4} it follows that 
\al{
\Vert \mbf P^{*}_{T-2|T}\Vert\le \Vert \mbf P^{*}_{T-2|T-2}\Vert+\Vert \mbf P^{*}_{T-2|T}-\mbf P^{*}_{T-2|T-2}\Vert = O_p(n^{-1})+O_p(n^{-2}).\label{eq:pivastar5}
}
Since all the bounds in \eqref{eq:pivastar2}-\eqref{eq:pivastar5} are the same for all $t$, from part (i) and \eqref{eq:pivastar} we have
\al{
\max_{t=1,\ldots, T}n\Vert \mbf P^{*}_{t|T}\Vert\le \max_{t=1,\ldots, T}n\Vert \mbf P^{*}_{t|t}\Vert+\max_{t=1,\ldots, T}n\Vert \mbf P^{*}_{t|T}-\mbf P^{*}_{t|t}\Vert = O_p(1)+O_p(n^{-1}).\nn
}
This proves part (iv) and completes the proof. $\Box$

\begin{lem}\label{lem:FFOnstar}
 Under Assumptions \ref{ass:common}, \ref{ass:idio}, \ref{ass:ind}, \ref{ass:linear}, \ref{ass:tails}, and \ref{ass:ident}, as $n,T\to\infty$:
\begin{compactenum}[(i)] 
\item for all $s=0,\ldots, T$, $\Vert \mbf F^{*}_{t|s}\Vert=O_p(1)$, uniformly in $t\le s$;
\item $n\Vert \mbf F^{*}_{t|T}-\mbf F^{*}_{t|t}\Vert=O_p(1)$, uniformly in $t$;
\item $n \Vert {\mbf F}_{t|t}^{*}-\wh{\mbf F}_{t}^{\text{\tiny \upshape {WLS}}*}\Vert=O_p(1)$, uniformly in  $t$;
\end{compactenum}
where $\wh{\mbf F}_{t}^{\text{\tiny \upshape {WLS}}*}=(\wh{\bm\Lambda}_n^{*\prime}(\wh{\bm\Sigma}_n^{\xi*})^{-1}\wh{\bm\Lambda}_n^{*})^{-1}\wh{\bm\Lambda}_n^{*\prime}(\wh{\bm\Sigma}_n^{\xi*})^{-1}\mbf x_{nt}.$
\end{lem}
\noindent
\textsc{Proof.} 
 The proof of part (i) follows the same steps as the proof of Lemma \ref{lem:FFO1} but when using Lemmas \ref{lem:starT12}, \ref{lem:sigmaunif}, and \ref{lem:eststar_LAST}, instead of Lemmas \ref{lem:est0LOAD}, \ref{lem:est0VAR}, \ref{lem:est0}, and \ref{lem:est0_LAST}.

For part (ii), from \eqref{eq:KS1} and \eqref{eq:pred1}
\begin{align}
\Vert{\mbf F}_{t|T}^{*}- {\mbf F}_{t|t}^{*}\Vert&\le \Vert{\mbf P}_{t|t}^{*}\Vert \, \Vert\wh{\mbf A}^{*}\Vert\, \Vert({\mbf P}_{t+1|t}^{*})^{-1}\Vert\, \{\Vert{\mbf F}_{t+1|T}^{*}\Vert+\Vert{\mbf F}_{t+1|t}^{*}\Vert\}\nn\\
&\le \Vert{\mbf P}_{t|t}^{*}\Vert \, \Vert\wh{\mbf A}^{*}\Vert\, \Vert({\mbf P}_{t+1|t}^{*})^{-1}\Vert\, \{\Vert{\mbf F}_{t+1|T}^{*}\Vert+\Vert\wh{\mbf A}^{*}\Vert\,\Vert{\mbf F}_{t|t}^{*}\Vert\}\nn\\
&= O_p(n^{-1}),\nn
\end{align}
by part (i) (when $s=T$ and $s=t$),
Lemmas 
\ref{lem:PPOnstar}(ii) and
\ref{lem:PPOnstar}(iii), and since $\Vert\wh{\mbf A}^{*}\Vert\le \Vert{\mbf A}\Vert+\Vert\wh{\mbf A}^{*}-\mbf A\Vert =O_p(1)$, by Assumption \ref{ass:common}(d) and Lemma \ref{lem:starT12}(iv). This proves part (ii).

For part (iii), from \eqref{eq:up1} and \eqref{eq:pred1}, by using Lemma \ref{lem:wood} (see also \eqref{eq:KFhathat} in the proof of Lemma \ref{lem:gennaio22bis} for more details)
\begin{align}
{\mbf F}_{t|t}^{*}&={\mbf F}_{t|t-1}^{*}+ {\mbf P}_{t|t-1}^{*}\wh{\bm\Lambda}_n^{*\prime}(\wh{\bm\Lambda}_n^{*}{\mbf P}_{t|t-1}^{*}\wh{\bm\Lambda}_n^{*\prime}+\wh{\bm\Sigma}_n^{\xi*})^{-1}(\mbf x_{nt}-\wh{\bm\Lambda}_n^{*}{\mbf F}_{t|t-1}^{*})\nn\\
=&\,  
(\wh{\bm\Lambda}_n^{*\prime}(\wh{\bm\Sigma}_n^{\xi*})^{-1}\wh{\bm\Lambda}_n^{*})^{-1}
\wh{\bm\Lambda}_n^{*\prime}(\wh{\bm\Sigma}_n^{\xi*})^{-1}\mbf x_{nt}\nn\\
&+\l\{(\wh{\bm\Lambda}_n^{*\prime}(\wh{\bm\Sigma}_n^{\xi*})^{-1}\wh{\bm\Lambda}_n^{*}+({\mbf P}_{t|t-1}^{*})^{-1})^{-1}
-
(\wh{\bm\Lambda}_n^{*\prime}(\wh{\bm\Sigma}_n^{\xi*})^{-1}\wh{\bm\Lambda}_n^{*})^{-1}
\r\}
\wh{\bm\Lambda}_n^{*\prime}(\wh{\bm\Sigma}_n^{\xi*})^{-1}\mbf x_{nt}\nn\\
&+\l\{\mbf I_r-
(\wh{\bm\Lambda}_n^{*\prime}(\wh{\bm\Sigma}_n^{\xi*})^{-1}\wh{\bm\Lambda}_n^{*}+(\wh{\mbf P}_{t|t-1}^{*})^{-1})^{-1}\wh{\bm\Lambda}_n^{*\prime}(\wh{\bm\Sigma}_n^{\xi*})^{-1}\wh{\bm\Lambda}_n^{*}\r\}\wh{\mbf A}^{*}{\mbf F}_{t-1|t-1}^{*}.
\label{eq:KFhathatstar}
\end{align}
Notice that the inverses in \eqref{eq:KFhathatstar} are all well defined by Lemmas \ref{lem:sigmaunif}(iv), \ref{lem:eststar_LAST}(iii), and \ref{lem:PPOnstar}(ii).

Now, by Lemma \ref{lem:denom}
\al{
\Vert & (\wh{\bm\Lambda}_n^{*\prime}(\wh{\bm\Sigma}_n^{\xi*})^{-1}\wh{\bm\Lambda}_n^{*}+({\mbf P}_{t|t-1}^{*})^{-1})^{-1}\wh{\bm\Lambda}_n^{*\prime}(\wh{\bm\Sigma}_n^{\xi*})^{-1}\wh{\bm\Lambda}_n^{*}-\mbf I_r\Vert
= O_p(n^{-1}).\label{eq:otitestar}
}
Furthermore, by Lemmas \ref{lem:LSL}(iii) and \ref{lem:eststar_LAST}(iv)
\al{
\Vert (\wh{\bm\Lambda}_n^{*\prime}(\wh{\bm\Sigma}_n^{\xi*})^{-1}\wh{\bm\Lambda}_n^{*}+({\mbf P}_{t|t-1}^{*})^{-1})^{-1}
-
(\wh{\bm\Lambda}_n^{*\prime}(\wh{\bm\Sigma}_n^{\xi*})^{-1}\wh{\bm\Lambda}_n^{*})^{-1}\Vert =O(n^{-2}),\label{eq:ennemeno2star}
}
and by Lemmas \ref{lem:LSL2}(vii) and \ref{lem:eststar_LAST}(ii), 
\al{
\Vert \wh{\bm\Lambda}_n^{*\prime}(\wh{\bm\Sigma}_n^{\xi*})^{-1}\Vert = O_p(\sqrt n).\label{eq:ennemeno05star}
}
Indeed, we can apply Lemmas \ref{lem:denom} and \ref{lem:LSL}(iii), since $\Vert({\mbf P}_{t|t-1}^{*})^{-1}\Vert=O_p(1)$ by Lemma \ref{lem:PPOnstar}(ii),
$\Vert(\wh{\bm\Sigma}_{n}^{\xi*})^{-1}\Vert=O_p(1)$
by Lemma \ref{lem:sigmaunif}(iv), and, by Lemmas \ref{lem:lambdasqrtn} and \ref{lem:starT12}(ii)  we have
\al{
n^{-1}\Vert \wh{\bm\Lambda}_n^{*\prime}\wh{\bm \Lambda}_n^{*}-{\bm\Lambda}_n^\prime{\bm \Lambda}_n\Vert\le&\, 2n^{-1} \Vert{\bm \Lambda}_n^\prime(\wh{\bm \Lambda}_n^{*}-{\bm \Lambda}_n)\Vert + n^{-1}\Vert 
(\wh{\bm \Lambda}_n^{*}-{\bm \Lambda}_n)^\prime(\wh{\bm \Lambda}_n^{*}-{\bm \Lambda}_n)\Vert\nn\\
\le &\, 2n^{-1/2}\Vert  {\bm \Lambda}_n\Vert \, n^{-1/2} \Vert \wh{\bm \Lambda}_n^{*}-{\bm \Lambda}_n\Vert + n^{-1}\Vert \wh{\bm \Lambda}_n^{*}-{\bm \Lambda}_n\Vert^2\nn\\
=&\,O_p(\max(n^{-1}\log^{2/\delta_v}T,{T}^{-1/2})).\nn
}
which, by Weyl's inequality \citep[Theorem 1]{MK04}, implies 
\al{
n^{-1}\vert \nu^{(j)}(\wh{\bm \Lambda}_n^{*\prime}\wh{\bm \Lambda}_n^{*})- \nu^{(j)}({\bm \Lambda}_n^\prime{\bm \Lambda}_n)\vert \le&\, n^{-1}\Vert \wh{\bm\Lambda}_n^{*\prime}\wh{\bm \Lambda}_n^{*}-{\bm\Lambda}_n^\prime{\bm \Lambda}_n\Vert=O_p(\max(n^{-1}\log^{2/\delta_v}T,{T}^{-1/2})),\nn
}
and, therefore, for $j=1,\ldots, r,$,
$$
\underline C_{j}\le \text{p-}\lim\!\!\!\!\!\!\inf_{n,T\to\infty} n^{-1} \nu^{(j)}(\wh{\bm \Lambda}_n^{*\prime}\wh{\bm \Lambda}_n^{*}) \le
\text{p-}\lim\!\!\!\!\!\sup_{n,T\to\infty} n^{-1} \nu^{(j)}(\wh{\bm\Lambda}_n^{*\prime}\wh{\bm \Lambda}_n^{*}) \le \overline C_j.
$$ 

By using \eqref{eq:otitestar}, \eqref{eq:ennemeno2star}, \eqref{eq:ennemeno05star} into \eqref{eq:KFhathatstar}:
\al{
\Vert& {\mbf F}_{t|t}^{*}-(\wh{\bm\Lambda}_n^{*\prime}(\wh{\bm\Sigma}_n^{\xi*})^{-1}\wh{\bm\Lambda}_n^{*})^{-1}\wh{\bm\Lambda}_n^{*\prime}(\wh{\bm\Sigma}_n^{\xi*})^{-1}\mbf x_{nt}\Vert\nn\\
 \le &\, n
\Vert(\wh{\bm\Lambda}_n^{*\prime}(\wh{\bm\Sigma}_n^{\xi*})^{-1}\wh{\bm\Lambda}_n^{*}+({\mbf P}_{t|t-1}^{*})^{-1})^{-1}
-
(\wh{\bm\Lambda}_n^{*\prime}(\wh{\bm\Sigma}_n^{\xi*})^{-1}\wh{\bm\Lambda}_n^{*})^{-1}
\Vert\,n^{-1/2}\Vert
\wh{\bm\Lambda}_n^{*\prime}(\wh{\bm\Sigma}_n^{\xi*})^{-1}\Vert\,n^{-1/2}\Vert \mbf x_{nt}\Vert \nn\\
&+ \Vert \mbf I_r-
(\wh{\bm\Lambda}_n^{*\prime}(\wh{\bm\Sigma}_n^{\xi*})^{-1}\wh{\bm\Lambda}_n^{*}+(\wh{\mbf P}_{t|t-1}^{*})^{-1})^{-1}\wh{\bm\Lambda}_n^{*\prime}(\wh{\bm\Sigma}_n^{\xi*})^{-1}\wh{\bm\Lambda}_n^{*}\Vert
\,\Vert \wh{\mbf A}^{*}\Vert\, \Vert {\mbf F}_{t-1|t-1}^{*}\Vert\nn\\
=&\, O_p(n^{-1}),\label{eq:girandolastar}
}
by part (i) (when $s=t-1$), Lemma \ref{lem:xunif}, and since $\Vert\wh{\mbf A}^{*}\Vert\le \Vert{\mbf A}\Vert+\Vert\wh{\mbf A}^{*}-\mbf A\Vert =O_p(1)$, by Assumption \ref{ass:common}(d) and Lemma \ref{lem:starT12}(i). This completes the proof. $\Box$

\begin{lem}\label{lem:FFO1sumstar}
Under Assumptions \ref{ass:common}, \ref{ass:idio}, \ref{ass:ind}, \ref{ass:linear}, \ref{ass:tails}, and \ref{ass:ident}, as $n,T\to\infty$, for $s=t$ and $s=T$:
\begin{compactenum}
\item [(i)] $\Vert T^{-1}\sum_{t=1}^T  \mbf F^{*}_{t|s}\mbf F_t^\prime\Vert=O_p(1)$;
\item [(i)] $\Vert T^{-1}\sum_{t=1}^T  \mbf F^{*}_{t|s}\mbf F_t^\prime\bm\lambda_i\Vert=O_p(1)$, uniformly in $i$;
\item [(ii)] $\Vert T^{-1}\sum_{t=1}^T  \mbf F^{*}_{t|s}\xi_{it}\Vert=O_p(1)$, uniformly in $i$.
\end{compactenum}
\end{lem}

\noindent
\textsc{Proof.} 
First notice that, for all $k=t-T,\ldots, t-1$,
\al{
\l\Vert n^{-1/2} T^{-1}\sum_{t=1}^T  \mbf x_{n,t-k}\mbf F_t^\prime\r\Vert=O_p(1), \label{eq:campanestar}
}
by Lemma \ref{lem:consistCOV}. The proof of part (i) follows by iterating either forward or backwards since both
$\Vert T^{-1}\sum_{t=1}^T  \mbf F_{t|t}^{*} \mbf F_t^\prime\Vert$ and $\Vert  T^{-1}\sum_{t=1}^T  \mbf F_{t|T}^{*} \mbf F_t^\prime\Vert$
are functions of \eqref{eq:campanestar} because of Lemma \ref{lem:FFOnstar}.

Part (ii) follows from part (i) and Assumption \ref{ass:common}(a). Part (iii) follows by substituting $\mbf F_t$ with $\xi_{it}$ in \eqref{eq:campanestar} and then by applying Lemma \ref{lem:consistCOV}(ii). This completes the proof. $\Box$

\begin{lem}\label{lem:milanmilan}
Under Assumptions \ref{ass:common}, \ref{ass:idio}, \ref{ass:ind}, \ref{ass:linear}, \ref{ass:tails}, and \ref{ass:ident}, as $n,T\to\infty$:
$$
\min(
n\log^{-2/\delta_v}T,
\sqrt{nT},
T{\log^{-1/2} n})\,\l\Vert T^{-1} \sum_{t=1}^T \{\mbf F_{t|T}^{*}\mbf F_{t|T}^{*\prime}+\mbf P_{t|T}^{*}\}-T^{-1} \sum_{t=1}^T  \mbf F_{t}\mbf F_{t}^{\prime}\r\Vert=O_p(1).
$$
\end{lem}

\noindent
\textsc{Proof.} 
Start with
\al{
\l\Vert T^{-1} \sum_{t=1}^T \mbf F_{t|T}^{*}\mbf F_{t|T}^{*\prime}-T^{-1} \sum_{t=1}^T  \mbf F_{t}\mbf F_{t}^{\prime}\r\Vert \le &\,
2\l\Vert T^{-1} \sum_{t=1}^T (\mbf F_{t|T}^{*}-\mbf F_t)\mbf F_{t}^{\prime}\r\Vert\nn\\
&+ \l\Vert T^{-1} \sum_{t=1}^T (\mbf F_{t|T}^{*}-\mbf F_t)(\mbf F_{t|T}^{*}-\mbf F_t)^{\prime}\r\Vert,\label{eq:vaccaboia1star}
}
and notice that if the first term on the rhs is $o_p(1)$ then the second term is dominated by the first one. So let us consider the first term on the rhs of \eqref{eq:vaccaboia1star}:
\al{
\l\Vert T^{-1} \sum_{t=1}^T (\mbf F_{t|T}^{*}-\mbf F_{t})\mbf F_t^\prime \r\Vert\le&\, 
\l\Vert T^{-1} \sum_{t=1}^T (\mbf F_{t|T}^{*}-\mbf F_{t|t}^{*})\mbf F_t^\prime \r\Vert\nn\\
&+\l\Vert T^{-1} \sum_{t=1}^T (\mbf F_{t|t}^{*}-\wh{\mbf F}_{t}^{\text{\tiny WLS}*})\mbf F_t^\prime \r\Vert\nn\\
&+\l\Vert T^{-1} \sum_{t=1}^T (\wh{\mbf F}_{t}^{\text{\tiny WLS}*}-\mbf F_t)\mbf F_t^\prime \r\Vert\nn\\
=&\, I^*+II^*+III^*, \;\text{ say.}\label{eq:vaccastar}
}
Let us consider each term in \eqref{eq:vaccastar}. First,
\al{
I^*\le&\,\max_{t=1,\ldots, T} \Vert{\mbf P}_{t|t}^{*}\Vert \, \Vert\wh{\mbf A}^{*}\Vert\, \max_{t=1,\ldots, T}\Vert({\mbf P}_{t+1|t}^{*})^{-1}\Vert\nn\\
&\cdot\l\{\l\Vert T^{-1}\sum_{t=1}^T {\mbf F}_{t+1|T}^{*}\mbf F_t^\prime\r\Vert
+\Vert\wh{\mbf A}^{*}\Vert\,\l\Vert T^{-1}\sum_{t=1}^T {\mbf F}_{t+1|t+1}^{*}\mbf F_t^\prime\r\Vert\r\}\nn\\
=&\, O_p(n^{-1}),\label{eq:Istar}
}
by Lemmas \ref{lem:PPOnstar}(ii), \ref{lem:PPOnstar}(iii), and \ref{lem:FFO1sumstar}, and since $\Vert\wh{\mbf A}^{*}\Vert\le \Vert{\mbf A}\Vert+\Vert\wh{\mbf A}^{*}-\mbf A\Vert =O_p(1)$, by Assumption \ref{ass:common}(d) and Lemma \ref{lem:starT12}(iv).  Second, from \eqref{eq:KFhathatstar} and \eqref{eq:girandolastar} in the proof of Lemma \ref{lem:FFOnstar}
\al{
II^*\le &\, O_p(n^{-1}) \l\{\l\Vert n^{-1/2}T^{-1}\sum_{t=1}^T \mbf x_{nt}\mbf F_t^\prime\r\Vert+\Vert\wh{\mbf A}^{*}\Vert\,\l\Vert T^{-1}\sum_{t=1}^T \mbf F_{t-1|t-1}^{*}\mbf F_t^\prime  \r\Vert\r\}\nn\\
&\, O_p(n^{-1}) \l\{n^{-1/2}\Vert \bm\Lambda_n\vert\, \l\Vert T^{-1}\sum_{t=1}^T \mbf F_t\mbf F_t^\prime\r\Vert+
\l\Vert n^{-1/2}T^{-1}\sum_{t=1}^T \bm\xi_{nt}\mbf F_t^\prime\r\Vert
+\Vert\wh{\mbf A}^{*}\Vert\,\l\Vert T^{-1}\sum_{t=1}^T \mbf F_{t-1|t-1}^{*}\mbf F_t^\prime  \r\Vert\r\}\nn\\
=&\, O_p(n^{-1}),\label{eq:IIstar}
}
because of Lemmas \ref{lem:lambdasqrtn}, \ref{lem:consistCOV}(i), combined with Assumption \ref{ass:ident}(b), \ref{lem:consistCOV}(iii), and \ref{lem:FFO1sumstar} and since $\Vert\wh{\mbf A}^{*}\Vert=O_p(1)$ by Assumption \ref{ass:common}(d) and Lemma \ref{lem:starT12}(iv).

Finally, let us consider the last term in \eqref{eq:vaccastar}. First, notice that we can write
\al{
\wh{\mbf F}_t^{\text{\tiny WLS}*}-\mbf F_t =&\, (\wh{\bm\Lambda}_n^{*\prime}(\wh{\bm\Sigma}_n^{\xi*})^{-1}\wh{\bm\Lambda}_n^{*})^{-1}\wh{\bm\Lambda}_n^{*\prime}(\wh{\bm\Sigma}_n^{\xi*})^{-1}(\bm\Lambda_n-\wh{\bm\Lambda}_n^{*})\mbf F_t+ (\wh{\bm\Lambda}_n^{*\prime}(\wh{\bm\Sigma}_n^{\xi*})^{-1}\wh{\bm\Lambda}_n^{*})^{-1}\wh{\bm\Lambda}_n^{*\prime}(\wh{\bm\Sigma}_n^{\xi*})^{-1}\bm\xi_{nt},\nn
}
which implies 
\al{
III^* \le &\, \l\Vert
T^{-1}\sum_{t=1}^T (\wh{\bm\Lambda}_n^{*\prime}(\wh{\bm\Sigma}_n^{\xi*})^{-1}\wh{\bm\Lambda}_n^{*})^{-1}\wh{\bm\Lambda}_n^{*\prime}(\wh{\bm\Sigma}_n^{\xi*})^{-1}(\bm\Lambda_n-\wh{\bm\Lambda}_n^{*})\mbf F_t\mbf F_t^\prime
\r\Vert\nn\\
&+ \l\Vert T^{-1}\sum_{t=1}^T (\wh{\bm\Lambda}_n^{*\prime}(\wh{\bm\Sigma}_n^{\xi*})^{-1}\wh{\bm\Lambda}_n^{*})^{-1}\wh{\bm\Lambda}_n^{*\prime}(\wh{\bm\Sigma}_n^{\xi*})^{-1}\bm\xi_{nt}\mbf F_t^\prime
\r\Vert\nn\\
\le&\, 
n\Vert(\bm\Lambda_n^\prime(\bm\Sigma_n^\xi)^{-1}\bm\Lambda_n)^{-1}\Vert \, n^{-1}\Vert \bm\Lambda_n^\prime(\bm\Sigma_n^\xi)^{-1} (\bm\Lambda_n-\wh{\bm\Lambda}_n^{*}) \Vert \l\Vert T^{-1}\sum_{t=1}^T \mbf F_t\mbf F_t^\prime\r\Vert\nn\\ 
&+\sqrt n\Vert
(\wh{\bm\Lambda}_n^{*\prime}(\wh{\bm\Sigma}_n^{\xi*})^{-1}\wh{\bm\Lambda}_n^{*})^{-1}\wh{\bm\Lambda}_n^{*\prime}(\wh{\bm\Sigma}_n^{\xi*})^{-1}-({\bm\Lambda}_n^{\prime}({\bm\Sigma}_n^{\xi})^{-1}{\bm\Lambda}_n)^{-1}{\bm\Lambda}_n^{\prime}({\bm\Sigma}_n^{\xi})^{-1}\Vert\nn\\
&\cdot n^{-1/2}\Vert \bm\Lambda_n-\wh{\bm\Lambda}_n^{*}\Vert \, \l\Vert T^{-1}\sum_{t=1}^T \mbf F_t\mbf F_t^\prime\r\Vert\nn\\ 
&+n\Vert(\bm\Lambda_n^\prime(\bm\Sigma_n^\xi)^{-1}\bm\Lambda_n)^{-1}\Vert \,n^{-1}\l\Vert T^{-1}\sum_{t=1}^T \bm\Lambda_n^\prime(\bm\Sigma_n^\xi)^{-1}\bm\xi_{nt}\mbf F_t^\prime \r\Vert\nn\\
&+ \sqrt n \Vert
(\wh{\bm\Lambda}_n^{*\prime}(\wh{\bm\Sigma}_n^{\xi*})^{-1}\wh{\bm\Lambda}_n^{*})^{-1}\wh{\bm\Lambda}_n^{*\prime}(\wh{\bm\Sigma}_n^{\xi*})^{-1}-({\bm\Lambda}_n^{\prime}({\bm\Sigma}_n^{\xi})^{-1}{\bm\Lambda}_n)^{-1}{\bm\Lambda}_n^{\prime}({\bm\Sigma}_n^{\xi})^{-1}\Vert\nn\\
&\cdot n^{-1/2}\l\Vert T^{-1}\sum_{t=1}^T\bm\xi_{nt}\mbf F_t^\prime\r\Vert\nn\\ 
=&\, III_a^*+III_b^*+III_c^*+III_d^*, \; \text{ say.}\label{eq:IIIstar}
}
Then, 
\al{
III_a^*\le&\, n\Vert(\bm\Lambda_n^\prime(\bm\Sigma_n^\xi)^{-1}\bm\Lambda_n)^{-1}\Vert \, n^{-1}\Vert \bm\Lambda_n^\prime(\bm\Sigma_n^\xi)^{-1} (\bm\Lambda_n-\wh{\bm\Lambda}_n^{\text{\tiny OLS}}) \Vert \l\Vert T^{-1}\sum_{t=1}^T \mbf F_t\mbf F_t^\prime\r\Vert\nn\\
&+n\Vert(\bm\Lambda_n^\prime(\bm\Sigma_n^\xi)^{-1}\bm\Lambda_n)^{-1}\Vert \, n^{-1/2}\Vert \bm\Lambda_n\Vert\,\Vert (\bm\Sigma_n^\xi)^{-1} \Vert\, n^{-1/2}\Vert\bm\Lambda_n^{\text{\tiny OLS}}-\wh{\bm\Lambda}_n^{*} \Vert \l\Vert T^{-1}\sum_{t=1}^T \mbf F_t\mbf F_t^\prime\r\Vert\nn\\
=&\,O_p(n^{-1/2}T^{-1/2})+
O_p(\max(
n^{-1}\log^{2/\delta_v}T,
n^{-1/2}T^{-1/2},
T^{-1}
)),
\label{eq:IIIastar}
}
where we used: for the first term on the rhs 
Lemmas \ref{lem:LSL2}(iii), \ref{lem:consistCOV}(i), combined with Assumption \ref{ass:ident}(b), and \ref{lem:COVFF0}(iv) since
$
n^{-1}\Vert \bm\Lambda_n^\prime(\bm\Sigma_n^\xi)^{-1} (\bm\Lambda_n-\wh{\bm\Lambda}_n^{\text{\tiny OLS}}) \Vert = n^{-1}T^{-1} \Vert \sum_{t=1}^T\bm\Lambda_n^\prime(\bm\Sigma_n^\xi)^{-1} \bm\xi_{nt}{\mbf F}_{t}^\prime\Vert$, while for the second term on the rhs we used
Assumption \ref{ass:idio}(a), and 
Lemmas 
\ref{lem:lambdasqrtn},
\ref{lem:LSL2}(iii), 
\ref{lem:consistCOV}(i), combined with Assumption \ref{ass:ident}(b), and \ref{lem:starols}(ii).
Moreover,
\al{
III_b^* = O_p(\max(
n^{-2}\log^{4/\delta_v} T,
n^{-1}T^{-1/2}\log^{2/\delta_v} T\sqrt{\log n},
T^{-1}\sqrt{\log n}
)),\label{eq:IIIbstar}
}
by  Lemmas \ref{lem:consistCOV}(i), combined with Assumption \ref{ass:ident}(b), \ref{lem:starT12}(ii), and \ref{lem:eststar_LAST}(v), 
\al{
III_c^* = O_p(n^{-1/2}T^{-1/2}),\label{eq:IIIcstar}
}
by Lemmas \ref{lem:LSL2}(iii) and \ref{lem:COVFF0}(iv), 
and
\al{
III_d^* = O_p(\max(
n^{-1}T^{-1/2}\log^{2/\delta_v} T
,
T^{-1}\sqrt{\log n})),\label{eq:IIIdstar}
}
by  Lemmas \ref{lem:consistCOV}(iii) and \ref{lem:eststar_LAST}(v). By substituting \eqref{eq:IIIastar}, \eqref{eq:IIIbstar}, \eqref{eq:IIIcstar}, and \eqref{eq:IIIdstar} into \eqref{eq:IIIstar} we have
\al{
III^* =O_p(\max(
n^{-1}\log^{2/\delta_v}T,
n^{-1/2}T^{-1/2},
T^{-1}\sqrt{\log n}
)).\label{eq:IIIstarbis}
}
Combining \eqref{eq:Istar}, \eqref{eq:IIstar}, and \eqref{eq:IIIstarbis} we have
\al{
\l\Vert T^{-1} \sum_{t=1}^T (\mbf F_{t|T}^{*}-\mbf F_t)\mbf F_{t}^{\prime}\r\Vert=O_p(\max(
n^{-1}\log^{2/\delta_v}T,
n^{-1/2}T^{-1/2},
T^{-1}\sqrt{\log n}
)),\label{eq:bella1star}
}
which once substituted into \eqref{eq:vaccaboia1star}, jointly with Lemma \ref{lem:PPOnstar}(iv) give
\al{
&\l\Vert T^{-1} \sum_{t=1}^T \{\mbf F_{t|T}^{*}\mbf F_{t|T}^{*\prime}+\mbf P_{t|T}^{*}\}-T^{-1} \sum_{t=1}^T  \mbf F_{t}\mbf F_{t}^{\prime}\r\Vert\nn\\
\le &\, 
\l\Vert T^{-1} \sum_{t=1}^T \mbf F_{t|T}^{*}\mbf F_{t|T}^{*\prime}-T^{-1} \sum_{t=1}^T  \mbf F_{t}\mbf F_{t}^{\prime}\r\Vert +\max_{t=1,\ldots, T} \Vert \mbf P_{t|T}^{*}\Vert \nn\\
=&\, O_p(\max(
n^{-1}\log^{2/\delta_v}T,
n^{-1/2}T^{-1/2},
T^{-1}\sqrt{\log n}
)) + O_p(n^{-1}),\nn
}
which completes the proof. $\Box$

\begin{lem} \label{lem:existEM}
Under Assumptions \ref{ass:common}, \ref{ass:idio}, and \ref{ass:ident}, the EM estimators of the parameters $\wh{\bm\varphi}_n\equiv\wh{\bm\varphi}_n^{(k+1)}$ defined in Section \ref{sec:Mstep}, exist and are unique, for any $k\ge0$.
\end{lem}

\noindent
\noindent\textsc{Proof.} 
Let $\mathcal O_i=\mathcal O_{\lambda_i}\times \mathcal O_{\sigma_i^2}\subset \mathbb R^{r+1}$, and 
$\mathcal O_r=\mathcal O_{\mathcal A}\times \mathcal O_{\Gamma^v}\subset \mathbb R^{r^2+r(r+1)/2}$ 
with $\mathcal O_{\lambda_i}$, $ \mathcal O_{\sigma_i^2}$, $\mathcal O_{\mathcal A}$, and $\mathcal O_{\Gamma^v}$ defined in Section \ref{enzuccio}.

At any iteration $k\ge 0$ the M-step requires solving the $n+1$ finite dimensional maximizations: 
\beq\label{O1}
(\wh{\bm\lambda}_i^{(k+1)}, \wh{\sigma}^{2(k+1)}_{i})=\arg\!\!\!\!\!\!\!\!\max_{(\underline{\bm\lambda}_i,\underline{\sigma}^2_{i})\in \mathcal O_i}  \E_{\wh{\varphi}_n^{(k)}}[\ell_i(\mbf x_{i1}\ldots \mbf x_{iT}|\bm F_T;\underline{\bm\lambda}_i,\underline{\sigma}^2_{i})],\qquad i=1,\ldots, n,
\eeq
and
\beq\label{O2}
(\wh{\mbf A}^{(k+1)}, \wh{\bm\Gamma}^{v(k+1)})=\arg\!\!\!\!\!\!\!\max_{(\underline{\mbf A},\underline{\bm\Gamma}^v)\in\mathcal O_r}\E_{\wh{\varphi}_n^{(k)}}[\ell(\bm F_T;\underline{\mbf A},\underline{\bm\Gamma}^v)],
\eeq
where 
 \beq\label{eq:LLbayes_exp_bis1App}
 \ell_i(\mbf x_{i1}\ldots \mbf x_{iT}|\bm F_T;\underline{\bm\lambda}_i,\underline{\sigma}^2_{i})
 =-\frac T2 \log(\underline{\sigma}^2_{i})-\frac 12\sum_{t=1}^T\frac{(x_{it}-\underline{\bm\lambda}_i^\prime\mbf F_t)^2}{\underline{\sigma}^2_{i}},
 \eeq
 and
\beq\label{eq:LLbayes_exp_bis2App}
\ell(\bm F_T;\underline{\mbf A},\underline{\bm\Gamma}^v)=-\frac {T}2 \log\det (\underline {\bm\Gamma}^v)-\frac 12\!\sum_{t=1}^T\l(\mbf F_t-\underline{\mbf A}\mbf F_{t-1}\r)^\prime(\underline {\bm\Gamma}^v)^{-1}\l(\mbf F_t-\underline{\mbf A}\mbf F_{t-1}\r).
\eeq
Now, the log-likelihoods \eqref{eq:LLbayes_exp_bis1App} and \eqref{eq:LLbayes_exp_bis2App} to be maximized are continuous and differentiable in $\mathcal O_i$ and $\mathcal O_r$, which are a compact sets by Assumptions \ref{ass:common}(a), \ref{ass:common}(d), \ref{ass:common}(e), and \ref{ass:idio}(a). Moreover, the log-likelihoods are  concave in their arguments, and \eqref{O1} and \eqref{O2} have a closed form expressions given in \eqref{eq:param1}-\eqref{eq:param3} and \eqref{eq:param4}-\eqref{eq:paramGv}, respectively. Last, notice that the true values of the parameters are fully identified by Assumption \ref{ass:ident}.

Therefore, by, e.g., \citet[Property 7.11 p.182]{GM95}, 
$\wh{\bm\lambda}_i^{(k+1)}$ and $\wh{\sigma}^{2(k+1)}_{i}$, $i=1,\ldots, n$, and  $\wh{\mbf A}^{(k+1)}$ and $\wh{\bm\Gamma}^{v(k+1)}$ exist and are unique for any $k\ge 0$. In particular, this result holds for $k=k^*$, i.e., for the EM estimators, and for any $n\in\mathbb N$ since \eqref{O1} can be solved separately for any $i$. This completes the proof. $\Box$

\begin{lem}\label{lem:localmax}
Under Assumptions \ref{ass:common}, \ref{ass:idio}, \ref{ass:linear}, and \ref{ass:ident}, $\ell(\bm X_{nT};\underline{\bm\varphi}_n)$ has a a local maximum denoted as\linebreak
$\wh{\bm\varphi}_n^{**}=(\wh{\bm\lambda}_1^{**\prime}\cdots \wh{\bm\lambda}_n^{**\prime}\; {\sigma}^{2**}_{1}\cdots {\sigma}^{2**}_{n}\; \text{\upshape vec}(\wh{\mbf A})^{**\prime},\text{\upshape vech}(\wh{\bm\Gamma})^{v**\prime})$, such that
\begin{compactenum}[(i)]
\item for all $i=1,\ldots, n$, $\lim_{k\to\infty}\Vert \wh{\bm\lambda}_{i}^{(k)} - \wh{\bm\lambda}_i^{**}\Vert=0$; 
\item for all $i=1,\ldots, n$, $\lim_{k\to\infty}\vert\wh{\sigma}_{i}^{2(k)} - \wh{\sigma}_i^{2**}\vert=0$; 
\item $\lim_{k\to\infty}\Vert \wh{\mbf A}^{(k)} - \wh{\mbf A}^{**}\Vert=0$;
\item $\lim_{k\to\infty}\Vert \wh{\bm\Gamma}^{v(k)} - \wh{\bm \Gamma}^{v**}\Vert=0$;
\end{compactenum}
where all convergences are monotonic.
\end{lem}

\noindent
\noindent\textsc{Proof.}
First notice that any maximum of $\ell(\bm X_{nT};\underline{\bm\varphi}_n)$ is also a maximum of $\tilde{\ell}(\bm X_{nT};\underline{\bm\varphi}_n)\equiv(nT)^{-1}\ell(\bm X_{nT};\underline{\bm\varphi}_n)$. Now, the starting point of the EM is such that $\tilde \ell(\bm X_{nT};\wh{\bm\varphi}_n^{(0)})>-\infty$ and since $\tilde\ell(\bm X_{nT};\underline{\bm\varphi}_n)$ is continuous and differentiable in the interior of $\mathcal O_n$, then $\{\tilde \ell(\bm X_{nT};\wh{\bm\varphi}_n^{(k)})\}_{k\ge 0}$ is bounded from above for any $n,T\in\mathbb N$.

Consider the definitions:
\begin{align}
\tilde \ell(\bm X_{nT};\underline{\bm\varphi}_n)&=\E_{\wh{\varphi}_{n}^{(k)}}[\tilde \ell(\bm X_{nT}|\bm F_T;\underline{\bm\varphi}_n)+\tilde \ell(\bm F_T;\underline{\bm\varphi}_n)|\bm X_{nT}]
-\E_{\wh{\varphi}_{n}^{(k)}}[\tilde \ell(\bm F_T|\bm X_{nT};\underline{\bm\varphi}_n)|\bm X_{nT}] \nn\\
&=\tilde {\mathcal Q}(\underline{\bm\varphi}_n,\wh{\bm\varphi}_{n}^{(k)})-\tilde {\mathcal H}(\underline{\bm\varphi}_n,\wh{\bm\varphi}_{n}^{(k)}).\label{eq:LLX2}
\end{align}
Then, for any $\underline{\bm\varphi}_n$ and any $n$ and $T$,
\begin{align}
\tilde {\mathcal H}(\underline{\bm\varphi}_{n};\wh{\bm\varphi}_n^{(k)})- \tilde {\mathcal H}(\wh{\bm\varphi}_n^{(k)};\wh{\bm\varphi}_n^{(k)})&=\E_{\wh{\bm\varphi}_n^{(k)}}[\tilde\ell(\bm F_T|\bm X_{nT};\underline{\bm\varphi}_{n})-\tilde\ell(\bm F_T|\bm X_{nT};\wh{\bm\varphi}_n^{(k)})|\bm X_{nT} ]\nn\\
&=\frac 1{nT} \E_{\wh{\bm\varphi}_n^{(k)}}\l[ \l.\log \frac{f(\bm F_T|\bm X_{nT};\underline{\bm\varphi}_{n})}{f(\bm F_T|\bm X_{nT};\wh{\bm\varphi}_n^{(k)})}\r|\bm X_{nT}\r]\nn\\
&\le\frac 1{nT} \log \E_{\wh{\bm\varphi}_n^{(k)}}\l[\l.\frac{f(\bm F_T|\bm X_{nT};\underline{\bm\varphi}_{n})}{f(\bm F_T|\bm X_{nT};\wh{\bm\varphi}_n^{(k)})}\r|\bm X_{nT}\r]\nn\\
&= \frac 1{nT} \log \int_{\mathbb R^{rT}} \frac{f(\bm F_T|\bm X_{nT};\underline{\bm\varphi}_{n})}{f(\bm F_T|\bm X_{nT};\wh{\bm\varphi}_n^{(k)})}f(\bm F_T|\bm X_{nT};\wh{\bm\varphi}_n^{(k)})\mathrm d \bm F_T\nn\\
&= \frac 1{nT} \log \int_{\mathbb R^{rT}} f(\bm F_T|\bm X_{nT};\underline{\bm\varphi}_{n}) \mathrm d \bm F_T = \frac 1{nT}\log(1)=0,\nn
\end{align}
by Jensen's inequality. Hence, we have (see also \citealp[Lemma 1]{DLR77})
\begin{align}
\tilde {\mathcal H}(\wh{\bm\varphi}_{n}^{(k+1)};\wh{\bm\varphi}_n^{(k)})&\leq \tilde {\mathcal H}(\wh{\bm\varphi}_n^{(k)};\wh{\bm\varphi}_n^{(k)}).\label{eq:HHH}
\end{align}
Therefore, from \eqref{eq:LLX2} and \eqref{eq:HHH}, for any $k$,
\beq\label{eq:convLLQQ}
\tilde \ell(\bm{X}_{nT};\wh{\bm\varphi}_{n}^{(k+1)})-\tilde  \ell(\bm{X}_{nT};\wh{\bm\varphi}_{n}^{(k)})\geq\tilde {\mathcal Q}(\wh{\bm\varphi}_{n}^{(k+1)};\wh{\bm\varphi}_{n}^{(k)})- \tilde {\mathcal Q}(\wh{\bm\varphi}_n^{(k)};\wh{\bm\varphi}_n^{(k)})\geq 0,
\eeq
where the last inequality holds by definition of the M-step. This shows that the log-likelihood $\ell(\bm{X}_{nT};\wh{\bm\varphi}_{n}^{(k)})$ increases monotonically as $k$ increases. 

Given that 
$\tilde {\mathcal Q}(\underline{\bm\varphi}_{n};\wh{\bm\varphi}_{n}^{(k)})$ has a unique maximum by Lemma \ref{lem:existEM},
and since for any $\underline{\bm\varphi}_n\in\mathcal O_n$ and $\overline{\bm\varphi}_n\in\mathcal O_n$, the function $\tilde{\mathcal Q}(\underline{\bm\varphi}_n;\overline{\bm\varphi}_n)$ is continuous in $\underline{\bm\varphi}_n$ and $\overline{\bm\varphi}_n$ and the components of the gradient vector $\nabla_{\underline{\bm\varphi}_n}\tilde{\mathcal Q}(\underline{\bm\varphi}_n;\overline{\bm\varphi}_n)$ are continuous in $\underline{\bm{\varphi}}_n$, from \citet[Theorem 3]{wu83}
\beq\label{eq:convLLQQ2}
\lim_{k\to\infty}\tilde\ell(\bm{X}_{nT}; \wh{\bm\varphi}_{n}^{(k)})=\tilde\ell(\bm{X}_{nT};\wh{\bm\varphi}_n^{**}),
\eeq
where the convergence is monotonic and $\wh{\bm\varphi}_n^{**}$ is a local maximum of $\tilde\ell(\bm{X}_{nT};\underline{\bm\varphi}_n)$.

Now notice that the solution of the M-step is such that $\wh{\bm\lambda}_i^{(k+1)}$ and $\wh{\mbf A}^{(k+1)}$ do not depend on other parameters at the same iteration, while $\wh{\sigma}_i^{2(k+1)}$ depends only on  $\wh{\bm\lambda}_i^{(k+1)}$, and $\wh{\bm\Gamma}^{v(k+1)}$ depends only on $\wh{\mbf A}^{(k+1)}$.
 
Since we are considering a Gaussian log-likelihood, from \citet[Condition 1]{wu83} holds, we have that, for any $i=1,\ldots, n$,
\al{
&\tilde{\mathcal Q}(\wh{\bm{\lambda}}_i^{(k+1)};\wh{\bm\varphi}_n^{(k)})-\tilde{\mathcal Q}(\wh{\bm{\lambda}}_i^{(k)};\wh{\bm\varphi}_n^{(k)}) \ge M_{\mathcal Q} \Vert \wh{\bm{\lambda}}_i^{(k+1)}-\wh{\bm\lambda}_i^{(k)}\Vert^2,\label{eq:convLLQQ3}\\
&\tilde{\mathcal Q}(\wh{\bm{\lambda}}_i^{(k+1)}, \wh{\sigma}_i^{2(k+1)};\wh{\bm\varphi}_n^{(k)})-\tilde{\mathcal Q}(\wh{\bm{\lambda}}_i^{(k)}, \wh{\sigma}_i^{2(k)};\wh{\bm\varphi}_n^{(k)}) \ge M_{\mathcal Q^\prime} \Vert (\wh{\bm{\lambda}}_i^{(k+1)\prime}\;\wh{\sigma}_i^{2(k+1)})-(\wh{\bm\lambda}_i^{(k)\prime}\;\wh{\sigma}_i^{2(k)})\Vert^2,\nn
}
for some finite positive reals $ M_{\mathcal Q}$ and $ M_{\mathcal Q^\prime}$, independent of $i$, and  where we use the short-hand notations:
\al{
&\tilde{\mathcal Q}(\wh{\bm{\lambda}}_i^{(k)};\wh{\bm\varphi}_n^{(k)})=\tilde{\mathcal Q}(\underline{\bm\lambda}_1,\ldots, \wh{\bm{\lambda}}_i^{(k)},\ldots,\underline{\bm\lambda}_n, \underline{\sigma}_1^2,\ldots,\underline{\sigma}_n^2,\underline{\mbf A}, \underline{\bm\Gamma}^v;\wh{\bm\varphi}_n^{(k)}),\nn\\
&\tilde{\mathcal Q}(\wh{\bm{\lambda}}_i^{(k)}, \wh{\sigma}_i^{2(k)};\wh{\bm\varphi}_n^{(k)})=\tilde{\mathcal Q}(\underline{\bm\lambda}_1,\ldots, \wh{\bm{\lambda}}_i^{(k)},\ldots,\underline{\bm\lambda}_n, \underline{\sigma}_1^2,\ldots,\wh{\sigma}_i^{2(k)},\ldots,\underline{\sigma}_n^2,\underline{\mbf A}, \underline{\bm\Gamma}^v;\wh{\bm\varphi}_n^{(k)}).\nn
}
And, similarly, 
\al{
&\tilde{\mathcal Q}(\wh{\mbf{A}}^{(k+1)};\wh{\bm\varphi}_n^{(k)})-\tilde{\mathcal Q}(\wh{\mbf{A}}^{(k)};\wh{\bm\varphi}_n^{(k)}) \ge M_{A^\prime} \Vert \wh{\mbf{A}}^{(k+1)}-\wh{\mbf A}^{(k)}\Vert^2,\label{eq:convLLQQ3A}\\
&\tilde{\mathcal Q}(\wh{\mbf{A}}^{(k+1)}, \wh{\bm{\Gamma}}^{v(k+1)};\wh{\bm\varphi}_n^{(k)})-\tilde{\mathcal Q}(\wh{\mbf{A}}^{(k)}, \wh{\bm{\Gamma}}^{v(k)};\wh{\bm\varphi}_n^{(k)}) \ge M_{v^\prime} \Vert \text{vec}(\wh{\mbf{A}}^{(k+1)}\; \wh{\bm{\Gamma}}^{v(k+1)})-\text{vec}(\wh{\mbf A}^{(k)};\wh{\bm \Gamma}^{v(k)})\Vert^2,\nn
}
for some finite positive reals $M_{A^\prime}$ and $M_{v^\prime}$.

Therefore, from \eqref{eq:convLLQQ}, \eqref{eq:convLLQQ2}, \eqref{eq:convLLQQ3}, and \eqref{eq:convLLQQ3A}, we have
\al{
&\lim_{k\to\infty}\Vert (\wh{\bm{\lambda}}_i^{(k+1)\prime}\;\wh{\sigma}_i^{2(k+1)})-(\wh{\bm\lambda}_i^{(k)\prime}\;\wh{\sigma}_i^{2(k)})\Vert=0,\nn\\
&\lim_{k\to\infty} \Vert \text{vec}(\wh{\mbf{A}}^{(k+1)}\; \wh{\bm{\Gamma}}^{v(k+1)})-\text{vec}(\wh{\mbf A}^{(k)};\wh{\bm \Gamma}^{v(k)})\Vert=0,\nn
}
which are sufficient conditions for applying the result in \citet[Theorem 6]{wu83}, i.e., such that
\begin{align}
&\lim_{k\to\infty} \Vert(\wh{\bm{\lambda}}_i^{(k)\prime}\;\wh{\sigma}_i^{2(k)})-  (\wh{\bm{\lambda}}_i^{**\prime}\;\wh{\sigma}_i^{2**})\Vert=0,\nn\\
&\lim_{k\to\infty} \Vert\text{vec}(\wh{\mbf{A}}^{(k)}\; \wh{\bm{\Gamma}}^{v(k)})-  \text{vec}(\wh{\mbf{A}}^{**}\; \wh{\bm{\Gamma}}^{v**})\Vert=0.\nn
\end{align}
This completes the proof.
$\Box$

\begin{lem}\label{lem:localglobal}
Under Assumptions \ref{ass:common}, \ref{ass:idio}, \ref{ass:ind}, \ref{ass:linear}, \ref{ass:tails}, and \ref{ass:ident}, as $n,T\to\infty$:
\begin{compactenum}
\item [(i)] $\min(n\log^{-2/\delta_v}, \sqrt{nT}\log^{-1/2}n,T)\,\max_{i=1,\ldots,n}\Vert\wh{\bm\lambda}_i^{**}-\wh{\bm\lambda}_i^{*}\Vert= O_p(1)$;
\item [(ii)] $\min(n\log^{-2/\delta_v}, \sqrt{nT}\log^{-1/2}n,T)\,\max_{i=1,\ldots,n} \vert\wh{\sigma}^{2**}_{i}-\wh{\sigma}^{2*}_{i}\vert=O_p(1)$;
\item [(iii)] $n\log^{-2/\delta_v}\,\Vert \wh{\mbf A}^{**}-\wh{\mbf A}^{*}\Vert=O_p(1)$; 
\item [(iv)] $n\log^{-2/\delta_v}\,\Vert \wh{\bm\Gamma}^{v**}-\wh{\bm\Gamma}^{v*}\Vert=O_p(1)$.
\end{compactenum}
\end{lem}

\noindent
\noindent\textsc{Proof.} 
First, notice that running the EM algorithm using $\ell(\bm X_{nT};\underline{\bm\varphi}_n)$ is equivalent to running the EM using 
$\ell(\bm X_{nT},\bm F_T;\underline{\bm\varphi}_n)$, and such EM will converge to a local maximum of $\ell(\bm X_{nT},\bm F_T;\underline{\bm\varphi}_n)$ because Lemma \ref{lem:localmax} would hold also in this case. Moreover,  since $\ell(\bm X_{nT},\bm F_T;\underline{\bm\varphi}_n)$ has clearly a unique maximum, then $\wh{\bm\varphi}_n^{**}$ is a local maximum of  $\ell(\bm X_{nT};\underline{\bm\varphi}_n)$ but also the unique global  maximum of $\ell(\bm X_{nT},\bm F_T;\underline{\bm\varphi}_n)$.

Consider the global QML estimator $\wh{\bm\varphi}_n^{*}=(\wh{\bm\phi}_n^{*\prime}\;\wh{\bm\theta}^{*\prime})^\prime$, where we defined $\wh{\bm\phi}_n^*=(\text{vec}(\wh{\bm\Lambda}_n^*)^\prime\; \wh{\sigma}_1^{2*}\cdots \wh{\sigma}_n^{2*})^\prime$ and $\wh{\bm\theta}^*=(\text{vec}(\wh{\mbf A}^*)^\prime\;\text{vech}(\wh{\bm\Gamma}^{v*})^\prime)$ maximizing $\ell(\bm X_{nT};\underline{\bm\varphi}_n)$.
The elements of $\wh{\bm\varphi}_n^{*}$ satisfy Lemma \ref{lem:starT12}. 

Now,  by definition the components of the gradient of $\ell(\bm X_{nT};\underline{\bm\varphi}_n)$ computed in $\wh{\bm\varphi}_n^{*}$ are such that (notice that $\ell(\bm F_T;\underline{\bm\varphi}_n)$ does not depend on $\underline{\bm\lambda}_i$)
\al{
\mbf 0_r &=  \nabla_{\underline{\bm\lambda}_i}\ell(\bm X_{nT};\underline{\bm\varphi}_n)\vert_{\underline{\bm\varphi}_n=\wh{\bm\varphi}_n^{*}}\nonumber\\
&= \nabla_{\underline{\bm\lambda}_i}\ell(\bm X_{nT},\bm F_T;\underline{\bm\varphi}_n)\vert_{\underline{\bm\varphi}_n=\wh{\bm\varphi}_n^{*}}- \nabla_{\underline{\bm\lambda}_i}\ell(\bm F_T|\bm X_{nT};\underline{\bm\varphi}_n)\vert_{\underline{\bm\varphi}_n=\wh{\bm\varphi}_n^{*}}\nn\\
&= \sum_{t=1}^T\mbf F_t (x_{it}-\mbf F_t^\prime \wh{\bm\lambda}_i^{*})- \nabla_{\underline{\bm\lambda}_i}\ell(\bm F_T|\bm X_{nT};\underline{\bm\varphi}_n)\vert_{\underline{\bm\varphi}_n=\wh{\bm\varphi}_n^{*}}\nn\\
&= \sum_{t=1}^T\mbf F_t (\mbf F_t^\prime\bm\lambda_i+\xi_{it}-\mbf F_t^\prime \wh{\bm\lambda}_i^{*})- \nabla_{\underline{\bm\lambda}_i}\ell(\bm F_T|\bm X_{nT};\underline{\bm\varphi}_n)\vert_{\underline{\bm\varphi}_n=\wh{\bm\varphi}_n^{*}},\label{eq:perbene}
}
and \eqref{eq:perbene} is equivalent to
\al{
(\wh{\bm\lambda}_i^{*}-\bm\lambda_i) =&\,\l(T^{-1}\sum_{t=1}^T\mbf F_t\mbf F_t^\prime\r)^{-1}\l(T^{-1}\sum_{t=1}^T\mbf F_t\xi_{it}\r)\nn\\
&+\l(T^{-1}\sum_{t=1}^T\mbf F_t\mbf F_t^\prime\r)^{-1} T^{-1}\nabla_{\underline{\bm\lambda}_i}\ell(\bm F_T|\bm X_{nT};\underline{\bm\varphi}_n)\vert_{\underline{\bm\varphi}_n=\wh{\bm\varphi}_n^{*}}\nn\\
=&\, (\bm\lambda_i-\bm\lambda_i^{\text{\tiny OLS}})+\l(T^{-1}\sum_{t=1}^T\mbf F_t\mbf F_t^\prime\r)^{-1} T^{-1}\nabla_{\underline{\bm\lambda}_i}\ell(\bm F_T|\bm X_{nT};\underline{\bm\varphi}_n)\vert_{\underline{\bm\varphi}_n=\wh{\bm\varphi}_n^{*}}.
\label{eq:perbene2}
}
By comparing \eqref{eq:perbene2} with  Lemma \ref{lem:sigmaunif}(i) (see in particular \eqref{eq:unifload2} and \eqref{eq:cambioposto} in its proof), 
and by Lemma \ref{lem:frida}, we have
\begin{align}
\max_{i=1,\ldots,n}\l\Vert
T^{-1} \nabla_{\underline{\bm\lambda}_i}\ell(\bm F_T|\bm X_{nT};\underline{\bm\varphi}_n)\vert_{\underline{\bm\varphi}_n
=\wh{\bm\varphi}_n^{*}}\r\Vert
=O_p(\max(n^{-1}\log^{2/\delta_v}T,n^{-1/2}T^{-1/2}\sqrt{\log n},T^{-1})),
\nn
\end{align}
which, from \eqref{eq:perbene}, implies
\begin{align}
\max_{i=1,\ldots,n}\l\Vert
T^{-1} \nabla_{\underline{\bm\lambda}_i}\ell(\bm X_{nT}, \bm F_T;\underline{\bm\varphi}_n)\vert_{\underline{\bm\varphi}_n
=\wh{\bm\varphi}_n^{*}}\r\Vert&=\max_{i=1,\ldots,n}\l\Vert
T^{-1} \nabla_{\underline{\bm\lambda}_i}\ell(\bm F_T|\bm X_{nT};\underline{\bm\varphi}_n)\vert_{\underline{\bm\varphi}_n
=\wh{\bm\varphi}_n^{*}}\r\Vert\label{eq:smallH}\\
&=O_p(\max(n^{-1}\log^{2/\delta_v}T,n^{-1/2}T^{-1/2}\sqrt{\log n},T^{-1})).\nn
\end{align}
Then, since $\nabla_{\underline{\bm\lambda}_i}\ell(\bm X_{nT}, \bm F_T;\underline{\bm\varphi}_n)$ is linear in $\underline{\bm\lambda}_i$, there exists a positive real $M_S^\prime$ such that
\begin{align}
\Vert T^{-1} \nabla_{\underline{\bm\lambda}_i}\ell(\bm X_{nT}, \bm F_T;\underline{\bm\varphi}_n)\vert_{\underline{\bm\varphi}_n=\wh{\bm\varphi}_n^{*}}\Vert \ge{ M_{S}^\prime} \Vert\wh{\bm\lambda}_i^*-\wh{\bm\lambda}_i^{**}\Vert,\label{eq:MSprime}
\end{align}
since, by definition, $T^{-1}\nabla_{\underline{\bm\lambda}_i}\ell(\bm X_{nT}, \bm F_T;\underline{\bm\varphi}_n)\vert_{\underline{\bm\varphi}_n=\wh{\bm\varphi}_n^{**}}=\mbf 0_r$. Therefore, from \eqref{eq:smallH} and \eqref{eq:MSprime}
\[
\max_{i=1,\ldots,n}\Vert\wh{\bm\lambda}_i^*-\wh{\bm\lambda}_i^{**}\Vert= O_p(\max(n^{-1}\log^{2/\delta_v}T,n^{-1/2}T^{-1/2}\sqrt{\log n},T^{-1})),
\]
which proves part (i).

Similarly, using Lemma \eqref{lem:sigmaunif}(ii) (see in particular \eqref{eq:italoc}, \eqref{eq:bidi}, and \eqref{eq:italoc2} in its proof), we can show that 
\beq\label{eq:bolle}
\max_{i=1,\ldots, n}\vert T^{-1} \nabla_{\underline{\sigma}_i^2}
\ell(\bm X_{nT}, \bm F_T;\underline{\bm\varphi}_n)
\vert_{\underline{\bm\varphi}_n=\wh{\bm\varphi}_n^{*}}\vert
=O_p(\max(n^{-1}\log^{2/\delta_v}T,n^{-1/2}T^{-1/2}\sqrt{\log n},T^{-1})),
\eeq
and
\begin{align}
\max_{i=1,\ldots, n}\vert T^{-1} \nabla_{\underline{\sigma}_i^2}
\ell(\bm X_{nT}, \bm F_T;\underline{\bm\varphi}_n)
\vert_{\underline{\bm\varphi}_n=\wh{\bm\varphi}_n^{*}}\vert
\ge M_{G}^\prime \vert\wh{\sigma}^{2**}_{i}-\wh{\sigma}^{2*}_{i}\vert.\label{eq:bolle2}
\end{align}
By using \eqref{eq:bolle} into \eqref{eq:bolle2}, we prove part (ii). 

Parts (iii) and (iv) follow in a similar way from Lemma \ref{lem:starT12}(iv) and \ref{lem:starT12}(v), respectively. This completes the proof. $\Box$

\begin{lem}\label{lem:convEM1}
Consider the EM estimator $\wh{\bm\lambda}_{i}\equiv \wh{\bm\lambda}_{i}^{(k+1)}$, for any $k\ge 0$.
Under Assumptions \ref{ass:common}, \ref{ass:idio}, \ref{ass:ind}, \ref{ass:linear}, \ref{ass:tails}, and \ref{ass:ident}, as $n,T\to\infty$, 
$$\min(
n^{2}\log^{-4/\delta_v}T,
nT\log^{-1/\delta_v}T{\log^{-1/2} n},
T^{3/2}{\log^{-1/2} n}
)\,
\Vert\wh{\bm\lambda}_{i}-\wh{\bm\lambda}_i^{**} \Vert = O_p(1),$$ 
uniformly in $i$. 
\end{lem}

\noindent\textsc{Proof.} First, notice that for any $\underline{\bm\varphi}_n$ and any $n$ and $T$
\al{
(nT)^{-1} \nabla_{\underline{\bm\varphi}_n} \ell(\bm X_{nT};\underline{\bm\varphi}_n)&=(nT)^{-1}\frac {\nabla_{\underline{\bm\varphi}_n}f(\bm X_{nT};\underline{\bm\varphi}_n)}{ f(\bm X_{nT};\underline{\bm\varphi}_n)} = (nT)^{-1} \int_{\mathbb R^{rT}} \frac { \nabla_{\underline{\bm\varphi}_n}f(\bm X_{nT},\bm F_{T};\underline{\bm\varphi}_n)}{ f(\bm X_{nT};\underline{\bm\varphi}_n)}\mathrm d\bm F_T\nn\\
&= (nT)^{-1}\int_{\mathbb R^{rT}}  \frac{\nabla_{\underline{\bm\varphi}_n}f(\bm X_{nT},\bm F_{T};\underline{\bm\varphi}_n)}{f(\bm X_{nT},\bm F_{T};\underline{\bm\varphi}_n)}\frac{f(\bm X_{nT},\bm F_{T};\underline{\bm\varphi}_n)}{ f(\bm X_{nT};\underline{\bm\varphi}_n)}\mathrm d\bm F_T\nn\\
&=(nT)^{-1}\int_{\mathbb R^{rT}}  \nabla_{\underline{\bm\varphi}_n}\ell(\bm X_{nT},\bm F_{T};\underline{\bm\varphi}_n)\frac{f(\bm X_{nT},\bm F_{T};\underline{\bm\varphi}_n)}{ f(\bm X_{nT};\underline{\bm\varphi}_n)}\mathrm d\bm F_T\nn\\
&= (nT)^{-1} \int_{\mathbb R^{rT}}  \nabla_{\underline{\bm\varphi}_n}\ell(\bm X_{nT},\bm F_{T};\underline{\bm\varphi}_n)f(\bm F_{T}|\bm X_{nT};\underline{\bm\varphi}_n)\mathrm d\bm F_T\nn\\
&=(nT)^{-1} \E_{\underline{\bm\varphi}_n}[ \nabla_{\underline{\bm\varphi}_n}\ell(\bm X_{nT},\bm F_{T};\underline{\bm\varphi}_n)|\bm X_{nT}]\nn\\
&=(nT)^{-1} \nabla_{\underline{\bm\varphi}_n}\mathcal Q(\underline{\bm\varphi}_n^\prime;\underline{\bm\varphi}_n)\vert_{\underline{\bm\varphi}_n^\prime=\underline{\bm\varphi}_n}.\label{pasticciscore}
}
Second, for any $k\ge 0$, by definition of the M-step, by \eqref{pasticciscore} 
and by a Taylor expansion about $\wh{\bm\lambda}_i^{(k)}$,
\al{
\mbf 0_r&= \nabla_{\underline{\bm\lambda}_i}\mathcal Q(\underline{\bm\varphi}_n;\wh{\bm\varphi}_n^{(k)})\vert_{\underline{\bm\varphi}_n=\wh{\bm\varphi}_n^{(k+1)}}\label{taylorQscore}\\
=&\, \nabla_{\underline{\bm\lambda}_i}\mathcal Q(\underline{\bm\varphi}_n;\wh{\bm\varphi}_n^{(k)})\vert_{\underline{\bm\varphi}_n=\wh{\bm\varphi}_n^{(k)}} 
+\nabla_{\underline{\bm\lambda}_i\underline{\bm\lambda}_i^\prime}\mathcal Q(\underline{\bm\varphi}_n;\wh{\bm\varphi}_n^{(k)})\vert_{\underline{\bm\varphi}_n=\wh{\bm\varphi}_n^{(k)}} (\wh{\bm\lambda}_i^{(k+1)}-\wh{\bm\lambda}_i^{(k)})+O(\Vert\wh{\bm\lambda}_i^{(k+1)}-\wh{\bm\lambda}_i^{(k)}\Vert^2)\nn\\
=&\, \nabla_{\underline{\bm\lambda}_i} \ell(\bm X_{nT};\underline{\bm\varphi}_n)\vert_{\underline{\bm\varphi}_n=\wh{\bm\varphi}_n^{(k)}}+ \nabla_{\underline{\bm\lambda}_i\underline{\bm\lambda}_i^\prime}\mathcal Q(\underline{\bm\varphi}_n;\wh{\bm\varphi}_n^{(k)})\vert_{\underline{\bm\varphi}_n=\wh{\bm\varphi}_n^{(k)}} (\wh{\bm\lambda}_i^{(k+1)}-\wh{\bm\lambda}_i^{(k)})+O(\Vert\wh{\bm\lambda}_i^{(k+1)}-\wh{\bm\lambda}_i^{(k)}\Vert^2)\nn\\
=&\, \nabla_{\underline{\bm\lambda}_i} \ell(\bm X_{nT};\underline{\bm\varphi}_n)\vert_{\underline{\bm\varphi}_n=\wh{\bm\varphi}_n^{(k)}}+ \nabla_{\underline{\bm\lambda}_i\underline{\bm\lambda}_i^\prime}\mathcal Q(\underline{\bm\varphi}_n;\wh{\bm\varphi}_n^{(k)})\vert_{\underline{\bm\varphi}_n=\wh{\bm\varphi}_n^{(k)}} (\wh{\bm\lambda}_i^{(k+1)}-\wh{\bm\lambda}_i^{(k)}),\nn
} 
since the third derivative of $\mathcal Q(\underline{\bm\varphi}_n;\wh{\bm\varphi}_n^{(k)})$ with respect to $\underline{\bm\lambda}_i$ is zero since the second derivative does not depend on $\underline{\bm\lambda}_i$ because $\mathcal Q(\underline{\bm\varphi}_n;\wh{\bm\varphi}_n^{(k)})$ is quadratic in $\underline{\bm\lambda}_i$.
And from \eqref{taylorQscore} it follows that
\al{
T^{-1}\nabla_{\underline{\bm\lambda}_i} \ell(\bm X_{nT};\underline{\bm\varphi}_n)\vert_{\underline{\bm\varphi}_n=\wh{\bm\varphi}_n^{(k)}}=&\, -T^{-1}\nabla_{\underline{\bm\lambda}_i\underline{\bm\lambda}_i^\prime}\mathcal Q(\underline{\bm\varphi}_n;\wh{\bm\varphi}_n^{(k)})\vert_{\underline{\bm\varphi}_n=\wh{\bm\varphi}_n^{(k)}} (\wh{\bm\lambda}_i^{(k+1)}-\wh{\bm\lambda}_i^{(k)}).
\label{taylorQscore2}
}
Third, by a Taylor expansion about $\wh{\bm\lambda}_i^{(k)}$, by definition of $\wh{\bm\lambda}_i^{**}$ which is a local maximum of $\ell(\bm X_{nT};\underline{\bm\varphi}_n)$, we have
\al{
\mbf 0_r &=T^{-1} \nabla_{\underline{\bm\lambda}_i} \ell(\bm X_{nT};\underline{\bm\varphi}_n)\vert_{\underline{\bm\varphi}_n=\wh{\bm\varphi}_n^{**}}\nn\\
=&\,T^{-1} \nabla_{\underline{\bm\lambda}_i} \ell(\bm X_{nT};\underline{\bm\varphi}_n)\vert_{\underline{\bm\varphi}_n=\wh{\bm\varphi}_n^{(k)}}+T^{-1}
 \nabla_{\underline{\bm\lambda}_i\underline{\bm\lambda}_i^\prime}\ell(\bm X_{nT};\underline{\bm\varphi}_n)\vert_{\underline{\bm\varphi}_n=\wh{\bm\varphi}_n^{(k)}}
(\wh{\bm\lambda}_i^{**}-\wh{\bm\lambda}_i^{(k)})\label{taylorQscore3}\\
&+ \frac 12 ((\wh{\bm\lambda}_i^{**}-\wh{\bm\lambda}_i^{(k)})^\prime\otimes \mbf I_r) 
\l\{T^{-1}\nabla_{\underline{\bm\lambda}_i} \text{vec}( \nabla_{\underline{\bm\lambda}_i\underline{\bm\lambda}_i^\prime}\ell(\bm X_{nT};\underline{\bm\varphi}_n))\vert_{\underline{\bm\varphi}_n=\check{\bm\varphi}_n}\r\} (\wh{\bm\lambda}_i^{**}-\wh{\bm\lambda}_i^{(k)})\nn\\
=&\, T^{-1} \nabla_{\underline{\bm\lambda}_i} \ell(\bm X_{nT};\underline{\bm\varphi}_n)\vert_{\underline{\bm\varphi}_n=\wh{\bm\varphi}_n^{(k)}}+T^{-1}
 \nabla_{\underline{\bm\lambda}_i\underline{\bm\lambda}_i^\prime}\ell(\bm X_{nT};\underline{\bm\varphi}_n)\vert_{\underline{\bm\varphi}_n=\wh{\bm\varphi}_n^{(k)}}
(\wh{\bm\lambda}_i^{**}-\wh{\bm\lambda}_i^{(k)})\nn\\
&+O_p\l(\Vert\wh{\bm\lambda}_i^{**}-\wh{\bm\lambda}_i^{(k)} \Vert^2n^{-1}\log^{2/\delta_v}T\r),\nn
}
where $\check{\bm\varphi}_n$ is such that $n^{-1/2} \Vert \wh{\bm\varphi}_n^{**}-\check{\bm\varphi}_n\Vert \le n^{-1/2} \Vert \wh{\bm\varphi}_n^{**}-\wh{\bm\varphi}_n^{(k+1)}\Vert$ and $n^{-1/2} \Vert \wh{\bm\varphi}_n^{(k+1)}-\check{\bm\varphi}_n\Vert \le n^{-1/2} \Vert \wh{\bm\varphi}_n^{**}-\wh{\bm\varphi}_n^{(k+1)}\Vert$. In particular, 
 the last term in \eqref{taylorQscore3} follows from Lemma \ref{lem:2max} and \citet[eq. (A.58) in the proof of Theorem 5]{MBPCAQML}, which imply
\al{
\sup_{\underline{\bm\varphi}_n\in \mathcal O_n} &\vert \nabla_{\underline{\bm\lambda}_i} \text{vec}( \nabla_{\underline{\bm\lambda}_i\underline{\bm\lambda}_i^\prime}\ell(\bm X_{nT};\underline{\bm\varphi}_n))\vert \nn\\
\le&\, 
\sup_{\underline{\bm\varphi}_n\in \mathcal O_n} \vert\nabla_{\underline{\bm\lambda}_i} \text{vec}( \nabla_{\underline{\bm\lambda}_i\underline{\bm\lambda}_i^\prime}\ell(\bm X_{nT};\underline{\bm\varphi}_n))
-\nabla_{\underline{\bm\lambda}_i} \text{vec}( \nabla_{\underline{\bm\lambda}_i\underline{\bm\lambda}_i^\prime}\ell_0(\bm X_{nT};\underline{\bm\phi}_n))\vert\nn\\
&+\sup_{\underline{\bm\varphi}_n\in\mathcal O_n} \vert\nabla_{\underline{\bm\lambda}_i} \text{vec}( \nabla_{\underline{\bm\lambda}_i\underline{\bm\lambda}_i^\prime}\ell_0(\bm X_{nT};\underline{\bm\phi}_n))\vert\nn\\
=&\, O_p(n^{-1}T \log^{2/\delta_v}T)+O(n^{-1}T)= O_p(n^{-1}T \log^{2/\delta_v}T).\label{eq:derivataterza}
}

Define the following $r\times r$ matrices:
\al{
\bm{\mathcal I}_{\text{c}}(\underline{\bm\lambda}_i) &=-\nabla_{\underline{\bm\lambda}_i\underline{\bm\lambda}_i^\prime}\mathcal Q(\underline{\bm\varphi}_n;\underline{\bm\varphi}_n)=- \E_{\underline{\varphi}_n}[\nabla_{\underline{\bm\lambda}_i\underline{\bm\lambda}_i^\prime} \ell(\bm{X}_{nT},\bm{F}_T;\underline{\bm\varphi}_n)|\bm{X}_{nT}]\nn\\
\bm I(\underline{\bm\lambda}_i) &=-\nabla_{\underline{\bm\lambda}_i\underline{\bm\lambda}_i^\prime}\ell(\bm X_{nT};\underline{\bm\varphi}_n).\nn
}
By substituting \eqref{taylorQscore2} into \eqref{taylorQscore3} and rearranging
\al{
(\wh{\bm\lambda}_i^{**}-\wh{\bm\lambda}_i^{(k)})=&\, - \l\{\nabla_{\underline{\bm\lambda}_i\underline{\bm\lambda}_i^\prime}\ell(\bm X_{nT};\underline{\bm\varphi}_n)\vert_{\underline{\bm\varphi}_n=\wh{\bm\varphi}_n^{(k)}}\r\}^{-1} \l\{\nabla_{\underline{\bm\lambda}_i} \ell(\bm X_{nT};\underline{\bm\varphi}_n)\vert_{\underline{\bm\varphi}_n=\wh{\bm\varphi}_n^{(k)}}\r\}\nn\\
&+O_p\l(\Vert\wh{\bm\lambda}_i^{**}-\wh{\bm\lambda}_i^{(k)} \Vert^2n^{-1}\log^{2/\delta_v}T\r)\nn\\
=&\,\l\{\nabla_{\underline{\bm\lambda}_i\underline{\bm\lambda}_i^\prime}\ell(\bm X_{nT};\underline{\bm\varphi}_n)\vert_{\underline{\bm\varphi}_n=\wh{\bm\varphi}_n^{(k)}}\r\}^{-1}
\l\{\nabla_{\underline{\bm\lambda}_i\underline{\bm\lambda}_i^\prime}\mathcal Q(\underline{\bm\varphi}_n;\wh{\bm\varphi}_n^{(k)})\vert_{\underline{\bm\varphi}_n=\wh{\bm\varphi}_n^{(k)}}\r\} (\wh{\bm\lambda}_i^{(k+1)}-\wh{\bm\lambda}_i^{(k)})\nn\\
&+O_p\l(\Vert\wh{\bm\lambda}_i^{**}-\wh{\bm\lambda}_i^{(k)} \Vert^2n^{-1}\log^{2/\delta_v}T\r)\nn\\
=&\, \l\{\bm I(\wh{\bm\lambda}_i^{(k)})  \r\}^{-1} \l\{\bm{\mathcal I}_{\text{c}}(\wh{\bm\lambda}_i^{(k)})\r\} (\wh{\bm\lambda}_i^{(k+1)}-\wh{\bm\lambda}_i^{**}+\wh{\bm\lambda}_i^{**}-\wh{\bm\lambda}_i^{(k)})\nn\\
&+O_p\l(\Vert\wh{\bm\lambda}_i^{**}-\wh{\bm\lambda}_i^{(k)} \Vert^2n^{-1}\log^{2/\delta_v}T\r).\label{eq:sundberg}
}
Moreover, by a Taylor approximation about $\wh{\bm\lambda}_i^{**}$,
\al{
&T^{-1}\bm{\mathcal I}_{\text{c}}(\wh{\bm\lambda}_i^{(k)})=T^{-1}\bm{\mathcal I}_{\text{c}}(\wh{\bm\lambda}_i^{**}),\nn\\
&T^{-1}\bm I(\wh{\bm\lambda}_i^{(k)}) =T^{-1} \bm I(\wh{\bm\lambda}_i^{**})+O_p\l(\Vert\wh{\bm\lambda}_i^{**}-\wh{\bm\lambda}_i^{(k)} \Vert n^{-1}\log^{2/\delta_v}T\r),\label{eq:peano}
}
where the first relation follows from the fact that, as noticed above, $\bm{\mathcal I}_{\text{c}}(\underline{\bm\lambda}_i)$ does not depend on $\underline{\bm\lambda}_i$, while the second follows again from \eqref{eq:derivataterza}.

Therefore, from \eqref{eq:sundberg} and \eqref{eq:peano}
\al{
(\wh{\bm\lambda}_i^{(k+1)}-\wh{\bm\lambda}_i^{**})=&\,\l(\mbf I_r -  \l\{\bm{\mathcal I}_{\text{c}}(\wh{\bm\lambda}_i^{(k)})\r\}^{-1}\l\{\bm I(\wh{\bm\lambda}_i^{(k)})  \r\} \r)(\wh{\bm\lambda}_i^{(k)}-\wh{\bm\lambda}_i^{**})\nn\\
&+O_p\l(\Vert\wh{\bm\lambda}_i^{**}-\wh{\bm\lambda}_i^{(k)} \Vert^2 n^{-1}\log^{2/\delta_v}T\r)\nn\\
=&\,\l(\mbf I_r -  \l\{\bm{\mathcal I}_{\text{c}}(\wh{\bm\lambda}_i^{**})\r\}^{-1}\l\{\bm I(\wh{\bm\lambda}_i^{**})  \r\} \r)(\wh{\bm\lambda}_i^{(k)}-\wh{\bm\lambda}_i^{**})\nn\\
&+O_p\l(\Vert\wh{\bm\lambda}_i^{**}-\wh{\bm\lambda}_i^{(k)} \Vert n^{-1}\log^{2/\delta_v}T\r),\label{eq:sundberg0}
}
see also \citet{sundberg74,sundberg76}, \citet{DLR77}, \citet{MR94}, \citet[Chapter 3.9, pp. 99-103]{MLT07}, and \citet[Chapter 8]{sundberg2019statistical}. 

Let $\bm R(\wh{\bm\lambda}_i^{**})= \mbf I_r -  \{\bm{\mathcal I}_{\text{c}}(\wh{\bm\lambda}_i^{**})\}^{-1}\{\bm I(\wh{\bm\lambda}_i^{**})\}$, then, by setting $k=k^*$ in \eqref{eq:sundberg0}, since $\wh{\bm\lambda}_i\equiv \wh{\bm\lambda}_i^{(k^*+1)}$, we have
\al{
(\wh{\bm\lambda}_i-\wh{\bm\lambda}_i^{**})=&\,\bm R(\wh{\bm\lambda}_i^{**})(\wh{\bm\lambda}_i^{(k^*)}-\wh{\bm\lambda}_i^{**})+O_p\l(\Vert\wh{\bm\lambda}_i^{**}-\wh{\bm\lambda}_i^{(k)} \Vert n^{-1}\log^{2/\delta_v}T\r).\label{taylorIII}
}
Hence, from \eqref{taylorIII}, by iterating backwards we get
\begin{align}
\Vert\wh{\bm\lambda}_i-\wh{\bm\lambda}_i^{**}\Vert \le&\,  
\Vert\wh{\bm\lambda}_i^{(k^*)}-\wh{\bm\lambda}_i^{**}\Vert \ \Vert \bm{ R}(\wh{\bm\lambda}_i^{**})\Vert + O_p\l(\Vert\wh{\bm\lambda}_i^{(k^*)}-\wh{\bm\lambda}_i^{**} \Vert n^{-1}\log^{2/\delta_v}T\r)\nn\\
\le&\, \l\{\Vert\wh{\bm\lambda}_i^{(k^*-1)}-\wh{\bm\lambda}_i^{**}\Vert \ \Vert \bm{ R}(\wh{\bm\lambda}_i^{**})\Vert + O_p\l(\Vert\wh{\bm\lambda}_i^{(k^*)}-\wh{\bm\lambda}_i^{**} \Vert n^{-1}\log^{2/\delta_v}T\r)
\r\} \ \Vert \bm{ R}(\wh{\bm\lambda}_i^{**})\Vert\nn\\
& + O_p\l(\Vert\wh{\bm\lambda}_i^{(k^*)}-\wh{\bm\lambda}_i^{**} \Vert n^{-1}\log^{2/\delta_v}T\r)\label{taylorIIoad}\\
\le&\,  \Vert\wh{\bm\lambda}_i^{(0)}-\wh{\bm\lambda}_i^{**}\Vert \ \Vert \bm{ R}(\wh{\bm\lambda}_i^{**})\Vert^{k^*+1} + \l\{\sum_{j=0}^{k^*}   \Vert \bm{ R}(\wh{\bm\lambda}_i^{**})\Vert^{j}  \r\} O_p\l(\Vert\wh{\bm\lambda}_i^{(k^*)}-\wh{\bm\lambda}_i^{**} \Vert n^{-1}\log^{2/\delta_v}T\r)\nn\\
\le&\,  \Vert\wh{\bm\lambda}_i^{(0)}-\wh{\bm\lambda}_i^{**}\Vert \ \Vert \bm{ R}(\wh{\bm\lambda}_i^{**})\Vert^{k^*+1} + \l\{\sum_{j=0}^{k^*}   \Vert \bm{ R}(\wh{\bm\lambda}_i^{**})\Vert^{j}  \r\} O_p\l(\Vert\wh{\bm\lambda}_i^{(0)}-\wh{\bm\lambda}_i^{**} \Vert n^{-1}\log^{2/\delta_v}T\r),\nn
\end{align}
where in the second and third step we used Lemma \ref{lem:localmax}, according to which $\Vert\wh{\bm\lambda}_i^{(k+1)}-\wh{\bm\lambda}_i^{**} \Vert\le \Vert\wh{\bm\lambda}_i^{(k)}-\wh{\bm\lambda}_i^{**} \Vert$ for any $k\ge 0$.

Let us consider separately the terms on the rhs of \eqref{taylorIIoad}. Because of Lemma \ref{lem:localglobal}(i),
\beq
\max_{i=1,\ldots,n}\Vert\bm{ R}(\wh{\bm\lambda}_i^{**}) - \bm{ R}(\wh{\bm\lambda}_i^{*})\Vert = O_p(\max(n^{-1}\log^{2/\delta_v}T,n^{-1/2}T^{-1/2}\sqrt{\log n},T^{-1})),\label{eq:approxR}
\eeq
since $\bm{ R}(\underline{\bm\lambda}_i)$ is continuous and differentiable in $\underline{\bm\lambda}_i$.
Consider the two matrices in $\bm{ R}(\wh{\bm\lambda}_i^{*})$. First, from \eqref{eq:LLbayes_exp_bis1} we can easily see that
\al{
T^{-1}\bm{\mathcal I}_{\text c} (\wh{\bm\lambda}_i^{*}) =&\, -T^{-1} \E_{\wh{\varphi}_n^{*}} [\nabla_{\underline{\bm\lambda}_i\underline{\bm\lambda}_i^\prime} \ell(\bm{X}_{nT},\bm{F}_T;\underline{\bm\varphi}_n)\vert_{\underline{\bm\varphi}_n=\wh{\bm\varphi}_n^*}|\bm{X}_{nT}]\nn\\
=&\,-T^{-1} \E_{\wh{\varphi}_n^{*}} [\nabla_{\underline{\bm\lambda}_i\underline{\bm\lambda}_i^\prime} \ell(\bm{X}_{nT}|\bm{F}_T;\underline{\bm\varphi}_n)\vert_{\underline{\bm\varphi}_n=\wh{\bm\varphi}_n^*}|\bm{X}_{nT}]\nn\\
=&\, T^{-1}(\wh{\sigma}_i^{2*})^{-1}\sum_{t=1}^T \E_{\wh{\varphi}_n^*}[\mbf F_t\mbf F_t^\prime|\bm X_{nT}]\nn\\
=&\, T^{-1}(\wh{\sigma}_i^{2*})^{-1}\sum_{t=1}^T \E_{\wh{\varphi}_n^*}[\mbf F_t|\bm X_{nT}]\E_{\wh{\varphi}_n^*}[\mbf F_t^\prime|\bm X_{nT}]\nn\\
&+ T^{-1}(\wh{\sigma}_i^{2*})^{-1}\sum_{t=1}^T \E_{\wh{\varphi}_n^*}[(\mbf F_t-\E_{\wh{\varphi}_n^*}[\mbf F_t|\bm X_{nT}])(\mbf F_t-\E_{\wh{\varphi}_n^*}[\mbf F_t|\bm X_{nT}])^\prime|\bm X_{nT}]\nn\\
=&\, T^{-1}(\wh{\sigma}_i^{2*})^{-1}\sum_{t=1}^T \l\{\mbf F_{t|T}^*\mbf F_{t|T}^{*\prime}+\mbf P_{t|T}^*\r\},\label{eq:erbette}
}
where in the last step we used Assumption \ref{ass:linear} which implies $\mbf F_{t|T}^*=\E_{\wh{\varphi}_n^*}[\mbf F_t|\bm X_{nT}]$ for all $t=1,\ldots, T$. Also note that for  \eqref{eq:HHH} and \eqref{eq:convLLQQ} in the proof of Lemma \ref{lem:localmax} to hold expectations have to be taken with respect to the same distribution as the one used to compute the log-likelihood. Hence, given that we maximize a mis-specified log-likelihood with diagonal idiosyncratic covariance, in the last step of \eqref{eq:erbette} we use $\mbf F_{t|T}^*$ and $\mbf P_{t|T}^*$, i.e., which are the outputs of the Kalman smoother implemented using $\wh{\bm\Sigma}_n^{\xi*}$.


Moreover, by Lemma \ref{lem:milanmilan}
\begin{align}
\l\Vert T^{-1}\sum_{t=1}^T\l\{\mbf F_{t|T}^*\mbf F_{t|T}^{*\prime}+\mbf P^*_{t|T}-\mbf F_t\mbf F_t^\prime\r\}\r\Vert
=&\,O_p(\max(
n^{-1}\log^{2/\delta_v}T,
n^{-1/2}T^{-1/2},
T^{-1}\sqrt{\log n}
)),\nn
\end{align}
which, once substituted into \eqref{eq:erbette}, gives
\al{
\max_{i=1,\ldots, n} &T^{-1}\l\Vert \bm{\mathcal I}_{\text c} (\wh{\bm\lambda}_i^{*}) -(\wh{\sigma}_i^{2*})^{-1}\sum_{t=1}^T\mbf F_t\mbf F_t^\prime\r\Vert\nn\\
&= 
O_p(\max(
n^{-1}\log^{2/\delta_v}T,
n^{-1/2}T^{-1/2},
T^{-1}\sqrt{\log n}
)).\label{eq:pezzotto}
}
This also shows that $T^{-1}\bm{\mathcal I}_{\text c} (\wh{\bm\lambda}_i^{*})$ is finite and positive definite, as $n,T\to\infty$, since $T^{-1}\sum_{t=1}^T \mbf F_t\mbf F_t^\prime$ is finite and positive definite with probability tending to one as $n,T\to\infty$, because of Lemma \ref{lem:consistCOV}(i) combined with Assumption \ref{ass:ident}(b), and since $\wh\sigma_i^{2*}$ is finite a positive for all $i=1,\ldots, n$ because of Lemma \ref{lem:sigmaunif}(ii) and Assumption \ref{ass:idio}(a).

Second, from Lemma \ref{lem:2max}
it follows that
\al{
\max_{i=1,\ldots, n}T^{-1}\l\Vert \nabla_{\underline{\bm\lambda}_i \underline{\bm\lambda}_i^\prime}\ell(\bm X_{nT};\underline{\bm\varphi}_n)\vert_{\underline{\bm\varphi}_n=\wh{\bm\varphi}_n^*}\r.&-\l.
 \nabla_{\underline{\bm\lambda}_i \underline{\bm\lambda}_i^\prime}\ell_0(\bm X_{nT};\underline{\bm\phi}_n)\vert_{\underline{\bm\phi}_n=\wh{\bm\phi}_n^*}\r\Vert\nn\\
 =&\, \max_{i=1,\ldots, n}T^{-1}\l\Vert \bm I(\wh{\bm\lambda}_i^*)-
 \nabla_{\underline{\bm\lambda}_i \underline{\bm\lambda}_i^\prime}\ell_0(\bm X_{nT};\underline{\bm\phi}_n)\vert_{\underline{\bm\phi}_n=\wh{\bm\phi}_n^*}\r\Vert\nn\\
  =&\, O_p(n^{-1}\log^{2/\delta_{v}}T).\label{eq:breitung4}
}
Furthermore, by a Taylor expansion about $\wh{\bm\lambda}_i^\dag$
\al{
\max_{i=1,\ldots, n}T^{-1}
\l\Vert\r.&\nabla_{\underline{\bm\lambda}_i \underline{\bm\lambda}_i^\prime}\ell_0(\bm X_{nT};\underline{\bm\phi}_n)\vert_{\underline{\bm\phi}_n=\wh{\bm\phi}_n^*}-\l.
\nabla_{\underline{\bm\lambda}_i \underline{\bm\lambda}_i^\prime}\ell_0(\bm X_{nT};\underline{\bm\phi}_n)\vert_{\underline{\bm\phi}_n=\wh{\bm\phi}_n^\dag}\r\Vert\nn\\
\le&\, \max_{i=1,\ldots, n}
T^{-1}\l\Vert
\nabla_{\underline{\bm\lambda}_i}\text{vec} \l(\nabla_{\underline{\bm\lambda}_i \underline{\bm\lambda}_i^\prime}\ell_0(\bm X_{nT};\underline{\bm\phi}_n) \r)\vert_{\underline{\bm\phi}_n=\wh{\bm\phi}_n^\dag}
\r\Vert\, \Vert\wh{\bm\lambda}_i^\dag-\wh{\bm\lambda}_i^* \Vert + O_p(\max_{i=1,\ldots, n}\Vert\wh{\bm\lambda}_i^\dag-\wh{\bm\lambda}_i^* \Vert^2)\nn\\
=&\, O_p(n^{-2}\log^{4/\delta_{v}}T),\label{eq:breitung5}
}
by Lemma \ref{lem:dagstar}(i) and since $\Vert
\nabla_{\underline{\bm\lambda}_i}\text{vec} (\nabla_{\underline{\bm\lambda}_i \underline{\bm\lambda}_i^\prime}\ell_0(\bm X_{nT};\underline{\bm\phi}_n) )
\Vert=O_p(n^{-1}T\log^{2/\delta_v}T)$ for all $\underline{\bm\phi}_n$, by \eqref{eq:derivataterza}.

Third, from \citet[Theorem 5]{MBPCAQML}:
\al{
\max_{i=1,\ldots, n}T^{-1}&\l\Vert 
\nabla_{\underline{\bm\lambda}_i \underline{\bm\lambda}_i^\prime}\ell_0(\bm X_{nT};\underline{\bm\phi}_n)\vert_{\underline{\bm\phi}_n=\wh{\bm\phi}_n^\dag}-
\nabla_{\underline{\bm\lambda}_i \underline{\bm\lambda}_i^\prime}\ell_0(\bm X_{nT}|\bm F_T;\underline{\bm\phi}_n)\vert_{\underline{\bm\phi}_n=\wh{\bm\phi}_n^\dag}
\r\Vert\nn\\
=&\,\max_{i=1,\ldots, n}T^{-1}\l\Vert 
\nabla_{\underline{\bm\lambda}_i \underline{\bm\lambda}_i^\prime}\ell_0(\bm X_{nT};\underline{\bm\phi}_n)\vert_{\underline{\bm\phi}_n=\wh{\bm\phi}_n^\dag}-\l(-(\wh\sigma_i^{2\dag})^{-1}\sum_{t=1}^T \mbf F_t\mbf F_t^\prime\r)
\r\Vert\nn\\
=&\, O_p(\max(n^{-1}, n^{-1/2}T^{-1/2}))
.\label{eq:44gatti}
}
Last, from Lemma \ref{lem:dagstar}(iii):
\beq\label{eq:UHT}
\max_{i=1,\ldots, n}\vert \wh{\sigma}_i^{2\dag}-\wh{\sigma}_i^{2*}\vert =O_p(n^{-1}\log^{2/\delta_v}T).
\eeq
By combining \eqref{eq:breitung4}, \eqref{eq:breitung5}, \eqref{eq:44gatti}, and \eqref{eq:UHT}, we have
\al{
\max_{i=1,\ldots, n}&T^{-1}\l\Vert \bm I(\wh{\bm\lambda}_i^*)-(\wh\sigma_i^{2*})^{-1}\sum_{t=1}^T \mbf F_t\mbf F_t^\prime\r\Vert=O_p(\max(n^{-1}\log^{2/\delta_v}T, n^{-1/2}T^{-1/2}).\label{eq:pezznove}
}
This is also shows that $T^{-1}\bm I(\wh{\bm\lambda}_i^*)$ is finite and positive definite, as $n,T\to\infty$, since $T^{-1}\sum_{t=1}^T \mbf F_t\mbf F_t^\prime$ is finite and positive definite with probability tending to one as $n,T\to\infty$, because of Lemma \ref{lem:consistCOV}(i) combined with Assumption \ref{ass:ident}(b), and since $\wh\sigma_i^{2*}$ is finite a positive for all $i=1,\ldots, n$ because of Lemma \ref{lem:sigmaunif}(ii) and Assumption \ref{ass:idio}(a).

Therefore, by using \eqref{eq:approxR}, \eqref{eq:pezzotto}, and \eqref{eq:pezznove}, and since, as remarked above $T^{-1}\bm{\mathcal I}_{\text c} (\wh{\bm\lambda}_i^{*})$
is finite and positive definite, we have
\al{
\max_{i=1,\ldots, n}\Vert\bm R(\wh{\bm\lambda}_i^{**})\Vert&=\Vert \mbf I_r -  \{\bm{\mathcal I}_{\text{c}}(\wh{\bm\lambda}_i^{*})\}^{-1}\{\bm I(\wh{\bm\lambda}_i^{*})\}\Vert\nn\\
=&\, \max_{i=1,\ldots, n} \Vert
\{\bm{\mathcal I}_{\text{c}}(\wh{\bm\lambda}_i^{*})^{-1}\}\{ \bm{\mathcal I}_{\text{c}}(\wh{\bm\lambda}_i^{*}) -\bm I(\wh{\bm\lambda}_i^{*}) \}
\Vert\nn\\
\le&\, \max_{i=1,\ldots, n} T \Vert
\{\bm{\mathcal I}_{\text{c}}(\wh{\bm\lambda}_i^{*})^{-1}\}\Vert \, T^{-1}\Vert \bm{\mathcal I}_{\text{c}}(\wh{\bm\lambda}_i^{*}) -\bm I(\wh{\bm\lambda}_i^{*}) 
\Vert\nn\\
\le&\, \max_{i=1,\ldots, n} T \Vert
\{\bm{\mathcal I}_{\text{c}}(\wh{\bm\lambda}_i^{*})^{-1}\}\Vert \,
T^{-1}\l\Vert \bm{\mathcal I}_{\text{c}}(\wh{\bm\lambda}_i^{*}) - (\wh{\sigma}_i^{2*})^{-1}\sum_{t=1}^T\mbf F_t\mbf F_t^\prime
\r\Vert\nn\\
&+\max_{i=1,\ldots, n} T \Vert
\{\bm{\mathcal I}_{\text{c}}(\wh{\bm\lambda}_i^{*})^{-1}\}\Vert \,\l\Vert
(\wh{\sigma}_i^{2*})^{-1}\sum_{t=1}^T\mbf F_t\mbf F_t^\prime
-\bm I(\wh{\bm\lambda}_i^{*}) 
\r\Vert
\nn\\
=&\, 
O_p(\max (n^{-1}\log^{2/\delta_v}T, n^{-1/2}T^{-1/2}\sqrt{\log n},T^{-1}\sqrt{\log n}).\label{eq:RRR2}
}

Moreover,
\begin{align}
\Vert \wh{\bm\lambda}_i^{(0)}- \wh{\bm\lambda}_i^{**}\Vert\le&\, \Vert \wh{\bm\lambda}_i^{(0)}- \wh{\bm\lambda}_i^{*}\Vert+\Vert \wh{\bm\lambda}_i^{*}- \wh{\bm\lambda}_i^{**}\Vert\nn\\
\le&\, \Vert \wh{\bm\lambda}_i^{(0)}- {\bm\lambda}_i\Vert+\Vert {\bm\lambda}_i- \wh{\bm\lambda}_i^{*}\Vert+\Vert \wh{\bm\lambda}_i^{*}- \wh{\bm\lambda}_i^{**}\Vert\nn\\
=&\,O_p(\max(n^{-1},T^{-1/2})+O_p(\max(n^{-1}\log^{2/\delta_v}T,T^{-1/2}))\nn\\
&+O_p(\max(n^{-1}\log^{2/\delta_v}T,n^{-1/2}T^{-1/2}\sqrt{\log n},T^{-1}))\nn\\
=&\,O_p(\max(n^{-1}\log^{2/\delta_v}T,T^{-1/2}),\label{eq:PCA0star}
\end{align}
because of Lemmas \ref{lem:est0LOAD}(i), \ref{lem:starT12}(i), and \ref{lem:localglobal}(i), respectively. 

By substituting  \eqref{eq:RRR2} and \eqref{eq:PCA0star}  into \eqref{taylorIIoad}, we have
\begin{align}
\Vert\wh{\bm\lambda}_{i}-\wh{\bm\lambda}_i^{**}\Vert\le&\,  \Vert\wh{\bm\lambda}_i^{(0)}-\wh{\bm\lambda}_i^{**}\Vert \ \Vert \bm{ R}(\wh{\bm\lambda}_i^{**})\Vert^{k^*+1} + o_p(\max(n^{-1}\log^{2/\delta_v}T,T^{-1/2})) \nn\\
=&\,O_p(\max(n^{-1}\log^{2/\delta_v}T,T^{-1/2}))\nn\\
&\cdot\l\{O_p(\max (n^{-1}\log^{2/\delta_v}T, n^{-1/2}T^{-1/2}\sqrt{\log n},T^{-1}\sqrt{\log n})\r\}^{k^*+1}\nn\\
&+o_p(\max(n^{-1}\log^{2/\delta_v}T,T^{-1/2}))\label{convEMML3}\\
=&\, O_p(\max(
n^{-2}\log^{4/\delta_v}T,
n^{-1}T^{-1}\log^{1/\delta_v}T\sqrt{\log n},
T^{-3/2}\sqrt{\log n}
)),\nn
\end{align}
uniformly in $i$ since the rhs of \eqref{eq:RRR2} does not depend on $i$. 
This completes the proof. $\Box$

\begin{lem}\label{lem:convEMiAH}
Consider the EM estimators $\wh{\sigma}_{i}^2\equiv \wh{\sigma}_{i}^{2(k+1)}$, 
$\wh{\mbf A}\equiv \wh{\mbf A}^{(k+1)}$, and $\wh{\bm\Gamma}^v\equiv\wh{\bm\Gamma}^{v(k+1)}$,
for any $k\ge 0$.
Under Assumptions \ref{ass:common}, \ref{ass:idio}, \ref{ass:ind}, \ref{ass:linear}, \ref{ass:tails}, and \ref{ass:ident}, as $n,T\to\infty$: 
\begin{compactenum}
\item [(i)] $\min(n\log^{-2/\delta_v}T,\sqrt T)\, \vert\wh{\sigma}_i^2-\wh{\sigma}_i^{2**}\vert=O_p(1);$
\item [(ii)] $\min(n\log^{-2/\delta_v}T,\sqrt T)\, \Vert\wh{\mbf A}-\wh{\mbf A}^{**}\Vert=O_p(1);$ 
\item [(iii)] $\min(n\log^{-2/\delta_v}T,\sqrt T)\, \Vert\wh{\bm\Gamma}^v-\wh{\bm\Gamma}^{v**}\Vert=O_p(1)$. 
\end{compactenum}
\end{lem}

\noindent
\textsc{Proof.} 
Let,  
\al{
 { R}(\wh{\sigma}_i^{2**})&= 
 1-\l\{ \E_{\underline{\varphi}_n}[\nabla_{\underline{\sigma}_i^2\underline{\sigma}_i^2} \ell(\bm{X}_{nT},\bm{F}_T;\underline{\bm\varphi}_n)
 \vert_{\underline{\sigma}_i^2=\wh{\sigma}_i^{2**}}
 |\bm{X}_{nT}]
 \r\}^{-1}
 \l\{ \nabla_{\underline{\sigma}_i^2\underline{\sigma}_i^2} \ell(\bm{X}_{nT};\underline{\bm\varphi}_n)
 \vert_{\underline{\sigma}_i^2=\wh{\sigma}_i^{2**}}
 \r\}.\nn
}
Then, following the same steps leading to \eqref{taylorIIoad} in the proof of Lemma \ref{lem:convEM1}, we obtain
 \begin{align}
\vert\wh{\sigma}_i^2-\wh{\sigma}_i^{2**}\vert \le&\, \vert\wh{\sigma}_i^{2(0)}-\wh{\sigma}_i^{2**}\vert \ \vert { R}(\wh{\sigma}_i^{2**})\vert^{k^*+1} + C\vert\wh{\sigma}_i^{2(0)}-\wh{\sigma}_i^{2**}\vert\nn\\
\le&\, \l\{
 \vert\wh{\sigma}_i^{2(0)}-{\sigma}_i^{2}\vert + \vert{\sigma}_i^{2}-\wh{\sigma}_i^{2*}\vert+\vert\wh{\sigma}_i^{2*}-\wh{\sigma}_i^{2**}\vert
\r\} \vert { R}(\wh{\sigma}_i^{2**})\vert^{k^*+1} + C\vert\wh{\sigma}_i^{2(0)}-\wh{\sigma}_i^{2**}\vert\nn\\
=&\,O_p(\max(n^{-1},T^{-1/2})+O_p(\max(n^{-1}\log^{2/\delta_v}T,T^{-1/2})\nn\\
&+ O_p(\max(n^{-1}\log^{2/\delta_v}T,n^{-1/2}T^{-1/2}\sqrt{\log n},T^{-1}))\nn\\
=&\,O_p(\max(n^{-1}\log^{2/\delta_v}T,T^{-1/2})),\label{eq:starwars}
\end{align}
for some finite positive real $C$ independent of $i$, by Lemmas \ref{lem:est0}(i), 
\ref{lem:starT12}(iii), and \ref{lem:localglobal}(ii),  
and since $\vert { R}(\wh{\sigma}_i^{2**})\vert=O_p(1)$. This proves part (i).

For part (ii),  following again the same steps leading to \eqref{taylorIIoad} in the proof of Lemma \ref{lem:convEM1}, we obtain
\al{
\Vert\wh{\mbf A}-\wh{\mbf A}^{**}\Vert\le&\,  \Vert\wh{\mbf A}^{(0)}-\wh{\mbf A}^{**}\Vert \ \Vert {\bm R}(\wh{\mbf A}^{**})\Vert^{k^*+1} + C\Vert\wh{\mbf A}^{(0)}-\wh{\mbf A}^{**}\Vert\nn\\
\le&\,\l\{
\Vert\wh{\mbf A}^{(0)}-{\mbf A}\Vert+\Vert{\mbf A}-\wh{\mbf A}^{*}\Vert+\Vert\wh{\mbf A}^{*}-\wh{\mbf A}^{**}\Vert\r\} \nn\\
&\cdot \Vert {\bm R}(\wh{\mbf A}^{**})\Vert^{k^*+1} + C\Vert\wh{\mbf A}^{(0)}-\wh{\mbf A}^{**}\Vert\nn\\
=&\,O_p(n^{-1},T^{-1/2})+O_p(\min(n^{-1}\log^{2/\delta_v}T, T^{-1/2}))+O_p(n^{-1}\log^{2/\delta_v}T)\nn\\
=&\,O_p(\max(n^{-1}\log^{2/\delta_v}T,T^{-1/2})),\nn
}
for some finite positive real $C$, by Lemmas \ref{lem:est0VAR}(i), \ref{lem:starT12}(iv), \ref{lem:localglobal}(iii), and since $\Vert {\bm R}(\wh{\mbf A}^{**})\Vert=O_p(1)$. This proves part (ii).

Part (iii) is proved in the same way as part (ii) but using Lemmas \ref{lem:est0VAR}(ii), \ref{lem:starT12}(v), and \ref{lem:localglobal}(iv). This completes the proof. $\Box$


\begin{lem}\label{lem:convEM1unif}
Consider the EM estimators $\wh{\bm\lambda}_{i}\equiv \wh{\bm\lambda}_{i}^{(k+1)}$ and $\wh{\sigma}_{i}^2\equiv \wh{\sigma}_{i}^{2(k+1)}$, for any $k\ge 0$.
Under Assumptions \ref{ass:common}, \ref{ass:idio}, \ref{ass:ind}, \ref{ass:linear}, \ref{ass:tails}, and \ref{ass:ident}, as $n,T\to\infty$: 
\begin{compactenum}
\item [(i)]  $\min(
n^{2}\log^{-4/\delta_v}T,
nT\log^{-1/\delta_v}T{\log^{-1/2} n},
T^{3/2}{\log^{-1} n}
)\,\max_{i=1,\ldots,n}\Vert\wh{\bm\lambda}_{i}-\wh{\bm\lambda}_i^{**}\Vert=O_p(1)$;
\item [(ii)] $\min(n\log^{-2/\delta_v}T,\sqrt{T}{\log^{-1/2} n})\, \max_{i=1,\ldots,n}\vert \wh{\sigma}_i^2-\wh{\sigma}_i^{2**}\vert=O_p(1)$.
\end{compactenum}
\end{lem}

\noindent
\noindent\textsc{Proof.} First notice that from \eqref{eq:unifload3} in the proof of Lemma \ref{lem:sigmaunif}
\al{
\max_{i=1,\ldots,n}\Vert\wh{\bm\lambda}_{i}^{(0)}-{\bm\lambda}_i\Vert=O_p(\max(n^{-1},T^{-1/2}\sqrt{\log n})),\label{eq:bicoccaarci}
}
while from \eqref{eq:laringite} in the proof of Lemma \ref{lem:est0} it is clear that
\al{
\max_{i=1,\ldots,n}\vert \wh{\sigma}_i^{2(0)}-{\sigma}_i^{2}\vert=O_p(\max(n^{-1},T^{-1/2}\sqrt{\log n})),\label{eq:bicoccaarci2}
}
because of \eqref{eq:bicoccaarci}.

Then,
\begin{align}
\max_{i=1,\ldots,n}\Vert \wh{\bm\lambda}_i^{(0)}- \wh{\bm\lambda}_i^{**}\Vert
\le&\, \max_{i=1,\ldots,n}\Vert \wh{\bm\lambda}_i^{(0)}- {\bm\lambda}_i\Vert+\max_{i=1,\ldots,n}\Vert {\bm\lambda}_i- \wh{\bm\lambda}_i^{*}\Vert+\max_{i=1,\ldots,n}\Vert \wh{\bm\lambda}_i^{*}- \wh{\bm\lambda}_i^{**}\Vert\nn\\
=&\,O_p(\max(n^{-1},T^{-1/2}\sqrt{\log n})+O_p(\max(n^{-1}\log^{2/\delta_v}T,T^{-1/2}\sqrt{\log n}))\nn\\
&+O_p(\max(n^{-1}\log^{2/\delta_v}T,n^{-1/2}T^{-1/2}\sqrt{\log n},T^{-1}))\nn\\
=&\,O_p(\max(n^{-1}\log^{2/\delta_v}T,T^{-1/2}\sqrt{\log n}),\label{eq:PCA0starunif}
\end{align}
because of \eqref{eq:bicoccaarci} and Lemmas \ref{lem:sigmaunif}(i), \ref{lem:localglobal}(i).

Then, from \eqref{convEMML3} in the proof of Lemma \ref{lem:convEM1}
\al{
\max_{i=1,\ldots,n}\Vert\wh{\bm\lambda}_{i}-\wh{\bm\lambda}_i^{**}\Vert\le&\, \max_{i=1,\ldots,n} \Vert\wh{\bm\lambda}_i^{(0)}-\wh{\bm\lambda}_i^{**}\Vert \ \max_{i=1,\ldots,n}\Vert \bm{ R}(\wh{\bm\lambda}_i^{**})\Vert^{k^*+1}\nn\\
& + o_p(\max(n^{-1}\log^{2/\delta_v}T,T^{-1/2}\sqrt{\log n})) \nn\\
=&\, O_p(\max(n^{-1}\log^{2/\delta_v}T,T^{-1/2}\sqrt{\log n}))\nn\\
&\cdot\l\{O_p(\max (n^{-1}\log^{2/\delta_v}T, n^{-1/2}T^{-1/2}\sqrt{\log n},T^{-1}\sqrt{\log n})\r\}^{k^*+1}\nn\\
&+o_p(\max(n^{-1}\log^{2/\delta_v}T,T^{-1/2}\sqrt{\log n})\nn\\
=&\, O_p(\max(
n^{-2}\log^{4/\delta_v}T,
n^{-1}T^{-1}\log^{1/\delta_v}T\sqrt{\log n},
T^{-3/2}{\log n}
)),\label{eq:viennetta}
}
because of \eqref{eq:PCA0starunif} and \eqref{eq:RRR2} in the proof of Lemma \ref{lem:convEM1}. This proves part (i).

For part (ii), from \eqref{eq:starwars} in the proof of Lemma \ref{lem:convEMiAH}
\begin{align}
\max_{i=1,\ldots,n}\vert\wh{\sigma}_i^2-\wh{\sigma}_i^{2**}\vert 
\le&\, \l\{
 \max_{i=1,\ldots,n}\vert\wh{\sigma}_i^{2(0)}-{\sigma}_i^{2}\vert +\max_{i=1,\ldots,n} \vert{\sigma}_i^{2}-\wh{\sigma}_i^{2*}\vert+\max_{i=1,\ldots,n}\vert\wh{\sigma}_i^{2*}-\wh{\sigma}_i^{2**}\vert
\r\}\nn\\
&\cdot \max_{i=1,\ldots,n} \vert { R}(\wh{\sigma}_i^{2**})\vert^{k^*+1} + C \max_{i=1,\ldots,n}\vert\wh{\sigma}_i^{2(0)}-\wh{\sigma}_i^{2**}\vert\nn\\
=&\,O_p(\max(n^{-1},T^{-1/2}\sqrt{\log n})+O_p(\max(n^{-1}\log^{2/\delta_v}T,T^{-1/2}\sqrt{\log n})\nn\\
&+ O_p(\max(n^{-1}\log^{2/\delta_v}T,n^{-1/2}T^{-1/2}\sqrt{\log n},T^{-1}))\nn\\
=&\,O_p(\max(n^{-1}\log^{2/\delta_v}T,T^{-1/2}\sqrt{\log n})),\nn
\end{align}
because of \eqref{eq:bicoccaarci2}, and Lemmas \ref{lem:sigmaunif}(ii), \ref{lem:localglobal}(ii). This proves part (ii) and completes the proof. $\Box$

%

\begin{lem}\label{lem:convEM1muzzo}
Consider the EM algorithm initialized with any deterministic loadings $\check{\bm\Lambda}_n^{(0)}=(\check{\bm\lambda}_1^{(0)}\cdots\check{\bm\lambda}_n^{(0)})^\prime$ such that 
$\text{\upshape vec}(\check{\bm\Lambda}_n^{(0)})\in\{\mathcal O_{\lambda_i}^n\cap \mathcal E_{\Lambda_n}\}$ as defined in Section \ref{enzuccio}.
Then, under Assumptions \ref{ass:common}, \ref{ass:idio}, \ref{ass:ind}, \ref{ass:linear}, \ref{ass:tails}, and \ref{ass:ident}, as $n,T\to\infty$, 
\begin{compactenum}[(i)]
\item 
$\min(
n\log^{-2/\delta_v}T,
\sqrt{nT}{\log^{-1/2} n},
T{\log^{-1/2} n}
)\,
\Vert\wh{{\bm\lambda}}_{i}-\wh{\bm\lambda}_i^{**} \Vert = O_p(1),$
uniformly in $i$.
\item $\min(
n\log^{-2/\delta_v}T,
\sqrt{nT}{\log^{-1/2} n},
T{\log^{-1/2} n}
)\,\max_{i=1,\ldots,n}
\Vert\wh{{\bm\lambda}}_{i}-\wh{\bm\lambda}_i^{**} \Vert = O_p(1).$
\end{compactenum}
\end{lem}

\noindent\textsc{Proof.} From \eqref{taylorIIoad} in the proof of Lemma \ref{lem:convEM1},
\begin{align}
\Vert{\wh{\bm\lambda}}_{i}-\wh{\bm\lambda}_i^{**}\Vert\le&\, 
\Vert\check{\bm\lambda}_i^{(0)}-\wh{\bm\lambda}_i^{**}\Vert \ \Vert \bm{ R}(\wh{\bm\lambda}_i^{**})\Vert^{k^*+1} + \l\{\sum_{j=0}^{k^*}   \Vert \bm{ R}(\wh{\bm\lambda}_i^{**})\Vert^{j}  \r\} O_p\l(\Vert\check{\bm\lambda}_i^{(0)}-\wh{\bm\lambda}_i^{**} \Vert n^{-1}\log^{2/\delta_v}T\r)\nn\\
\le&\, 
\Vert\check{\bm\lambda}_i^{(0)}-\wh{\bm\lambda}_i^{**}\Vert \ \Vert \bm{ R}(\wh{\bm\lambda}_i^{**})\Vert^{k^*+1} +  O_p\l(\Vert\check{\bm\lambda}_i^{(0)}-\wh{\bm\lambda}_i^{**} \Vert n^{-1}\log^{2/\delta_v}T\r)\nn\\
=&\, \l\{O_p(\max (n^{-1}\log^{2/\delta_v}T, n^{-1/2}T^{-1/2}\sqrt{\log n},T^{-1}\sqrt{\log n})\r\}^{k^*+1}+ O_p( n^{-1}\log^{2/\delta_v}T)\nn\\
=&\, O_p(\max (n^{-1}\log^{2/\delta_v}T, n^{-1/2}T^{-1/2}\sqrt{\log n},T^{-1}\sqrt{\log n}),\nn
\end{align}
by \eqref{eq:RRR2} in the proof of Lemma \ref{lem:convEM1}, and since $\Vert \check{\bm\lambda}_i^{(0)}-\wh{\bm\lambda}_i^{**}\Vert \le 
\Vert \check{\bm\lambda}_i^{(0)}\Vert + \Vert \wh{\bm\lambda}_i^{**}-\wh{\bm\lambda}_i^{*}\Vert+\Vert \wh{\bm\lambda}_i^{*}-\bm\lambda_i\Vert \le M_\lambda O(1)$, by definition of the initial estimator, Assumption \ref{ass:common}(a), and Lemmas \ref{lem:starT12}(i) and \ref{lem:localglobal}(i). This proves part (i).

Part (ii) is proved in the same way but starting from \eqref{eq:viennetta} in the proof of Lemma \ref{lem:convEM1unif} and by noting that $\max_{i=1,\ldots, n}\Vert \check{\bm\lambda}_i^{(0)}\Vert\le M_\lambda$ by the same arguments as before. This completes the proof. $\Box$

%

\section{Lemmas necessary for proving Proposition \ref{prop:factors}}

\begin{lem}\label{lem:sigmaunifhat}
Consider the EM estimators $\wh{\bm\lambda}_{i}\equiv \wh{\bm\lambda}_{i}^{(k+1)}$, $\wh{\sigma}_{i}^2\equiv \wh{\sigma}_{i}^{2(k+1)}$, 
and $\wh{\bm\Sigma}_n^\xi\equiv\wh{\bm\Sigma}_n^{\xi(k+1)}$,
for any $k\ge 0$.
Under Assumptions \ref{ass:common}, \ref{ass:idio}, \ref{ass:ind}, \ref{ass:linear}, \ref{ass:tails}, and \ref{ass:ident}, as $n,T\to\infty$, 
\begin{compactenum}
\item [(i)] $\min(\sqrt T \log^{-1/2} n, n\log^{-2/\delta_v} T) \,\max_{i=1,\ldots,n} \Vert \wh{\bm\lambda}_i-\bm\lambda_i\Vert=O_p(1)$;
\item [(ii)] $\min(\sqrt T \log^{-1/2} n, n\log^{-2/\delta_v} T) \,\max_{i=1,\ldots,n} \vert \wh{\sigma}_i^{2}-\sigma_i^2\vert=O_p(1)$;
\item [(iii)] $\min(\sqrt T \log^{-1/2} n, n\log^{-2/\delta_v} T) \, \Vert \wh{\bm\Sigma}_n^{\xi}-\bm\Sigma_n^\xi\Vert = O_p(1)$;
\item [(iv)] $\Vert (\wh{\bm\Sigma}_n^{\xi})^{-1}\Vert = O_p(1)$;
\item [(v)] $\min(\sqrt T \log^{-1/2} n, n\log^{-2/\delta_v} T)\, \Vert (\wh{\bm\Sigma}_n^{\xi})^{-1}- ({\bm\Sigma}_n^{\xi})^{-1}\Vert = O_p(1)$.
\end{compactenum}
\end{lem}

\noindent
\noindent\textsc{Proof.}  For part (i)
\begin{align}
\max_{i=1,\ldots,n}\Vert \wh{\bm\lambda}_{i}-{\bm\lambda}_{i}\Vert \le&\, \max_{i=1,\ldots,n} \Vert\wh{\bm\lambda}_{i}-\wh{\bm\lambda}_{i}^{**}\Vert+\max_{i=1,\ldots,n} \Vert\wh{\bm\lambda}_{i}^{**}-\wh{\bm\lambda}_{i}^{*}\Vert+\max_{i=1,\ldots,n}\Vert\wh{\bm\lambda}_{i}^{*}-\bm\lambda_i\Vert \nn\\
=&\,O_p(\max(
n^{-2}\log^{4/\delta_v}T,
n^{-1}T^{-1}\log^{1/\delta_v}T\sqrt{\log n},
T^{-3/2}{\log n}
))\nn\\
&+O_p(\max(n^{-1}\log^{2/\delta_v}T,n^{-1/2}T^{-1/2}\sqrt{\log n},T^{-1}))\nn\\
&+ O_p(\max(n^{-1},T^{-1/2}\sqrt{\log n}))+O_p(n^{-1}\log^{2/\delta_v} T)\nn\\
=&\, O_p(\max(n^{-1}\log^{2/\delta_v}T,T^{-1/2}\sqrt{\log n}) ),\nn
\end{align}
by Lemmas \ref{lem:sigmaunif}(i), \ref{lem:localglobal}(i), and \ref{lem:convEM1unif}(i). This proves part (i).

Part (ii) is proved as part (i) but using Lemmas \ref{lem:sigmaunif}(ii), \ref{lem:localglobal}(ii), and \ref{lem:convEM1unif}(ii).

Part (iii) immediately follows from part (ii), indeed
\al{
\Vert \wh{\bm\Sigma}_n^{\xi}-\bm\Sigma_n^\xi\Vert \le \max_{i=1,\ldots, n} \vert  \wh{\sigma}_i^{2}-\sigma_i^{2} \vert=O_p(\max(n^{-1}\log^{2/\delta_v}T, T^{-1/2}\sqrt{\log n})).\nn
}

For part (iv) we have
\al{
\Vert (\wh{\bm\Sigma}_n^{\xi})^{-1}\Vert = \l\{\min_{i=1,\ldots, n} \wh{\sigma}_i^{2}\r\}^{-1}\le C_\xi +O_p(\max(n^{-1}\log^{2/\delta_v}T, T^{-1/2}\sqrt{\log n})),\nn
}
because of part (ii) and Assumption \ref{ass:idio}(a).

To conclude, for part (v) we have
\al{
\Vert (\wh{\bm\Sigma}_n^{\xi})^{-1}-(\bm\Sigma_n^\xi)^{-1}\Vert \le \Vert (\wh{\bm\Sigma}_n^{\xi})^{-1}\Vert \, \Vert \wh{\bm\Sigma}_n^{\xi}-\bm\Sigma_n^\xi\Vert\,
\Vert ({\bm\Sigma}_n^{\xi})^{-1}\Vert = O_p(\max(n^{-1}\log^{2/\delta_v}T, T^{-1/2}\sqrt{\log n})),\nn
}
by parts (iii), (iv), and Assumption \ref{ass:idio}(a). This completes the proof. $\Box$

\begin{lem}\label{lem:eststar_LASTHAT}
Consider the EM estimators $\wh{\bm\lambda}_{i}\equiv \wh{\bm\lambda}_{i}^{(k+1)}$, 
$\wh{\bm\Lambda}_{n}\equiv \wh{\bm\Lambda}_{n}^{(k+1)}$,
$\wh{\sigma}_{i}^2\equiv \wh{\sigma}_{i}^{2(k+1)}$, 
and $\wh{\bm\Sigma}_n^\xi\equiv\wh{\bm\Sigma}_n^{\xi(k+1)}$,
for any $k\ge 0$.
Under Assumptions \ref{ass:common}, \ref{ass:idio}, \ref{ass:ind}, \ref{ass:linear}, \ref{ass:tails}, and \ref{ass:ident}, as $n,T\to\infty$:
\begin{compactenum}
\item [(i)] $\min(n\log^{-2/\delta_v}T,\sqrt T \log^{-1/2}n
)\,n^{-1}\Vert\wh{\bm\Lambda}_n^{\prime}(\wh{\bm \Sigma}_n^{\xi})^{-1}\wh{\bm\Lambda}_n-\bm\Lambda_n^\prime(\bm\Sigma_n^\xi)^{-1}\bm\Lambda_n\Vert = O_p(1)$;
\item [(ii)] $\min(n\log^{-2/\delta_v}T,\sqrt T \log^{-1/2}n
)\,n^{-1/2}\Vert\wh{\bm\Lambda}_n^{\prime}(\wh{\bm \Sigma}_n^{\xi})^{-1}-\bm\Lambda_n^\prime(\bm\Sigma_n^\xi)^{-1}\Vert = O_p(1)$;
\item [(iii)] $n\Vert(\wh{\bm\Lambda}_n^{\prime}(\wh{\bm \Sigma}_n^{\xi})^{-1}\wh{\bm\Lambda}_n)^{-1}\Vert = O_p(1)$;
\item [(iv)] $\min(n\log^{-2/\delta_v}T,\sqrt T \log^{-1/2}n
)\,n\Vert(\wh{\bm\Lambda}_n^{\prime}(\wh{\bm \Sigma}_n^{\xi})^{-1}\wh{\bm\Lambda}_n)^{-1}-(\bm\Lambda_n^\prime(\bm\Sigma_n^\xi)^{-1}\bm\Lambda_n)^{-1}\Vert = O_p(1)$;
\item [(v)] $\omega_{n,T,\delta_v}\,\sqrt n\Vert(\wh{\bm\Lambda}_n^{\prime}(\wh{\bm \Sigma}_n^{\xi})^{-1}\wh{\bm\Lambda}_n)^{-1}
\wh{\bm\Lambda}_n^{\prime}(\wh{\bm \Sigma}_n^{\xi})^{-1}
-(\bm\Lambda_n^\prime(\bm\Sigma_n^\xi)^{-1}\bm\Lambda_n)^{-1}\bm\Lambda_n^\prime(\bm\Sigma_n^\xi)^{-1}\Vert = O_p(1)$,\\ with $\omega_{n,T,\delta_v}=\min(n\log^{-2/\delta_v}T,\sqrt T \log^{-1/2}n
).$
\end{compactenum}
\end{lem}

\noindent\textsc{Proof.} The proof is the same as the proof of Lemma \ref{lem:eststar_LAST} but using Proposition \ref{prop:load}(a) (which does not require this lemma to be proved) and Lemma \ref{lem:sigmaunifhat} instead of Lemmas \ref{lem:starT12} and \ref{lem:sigmaunif}. $\Box$

\begin{lem}\label{lem:PPOnhat}
Consider the MSEs estimators $\wh{\mbf P}_{t|t-1}\equiv  {\mbf P}_{t|t-1}^{(k+1)}$, $\wh{\mbf P}_{t|t}\equiv  {\mbf P}_{t|t}^{(k+1)}$, and $\wh{\mbf P}_{t|T}\equiv  {\mbf P}_{t|T}^{(k+1)}$ derived from the Kalman filter and smoother and obtained using the EM estimators $\wh{\bm\lambda}_{i}\equiv \wh{\bm\lambda}_{i}^{(k+1)}$, 
$\wh{\bm\Lambda}_{n}\equiv \wh{\bm\Lambda}_{n}^{(k+1)}$,
$\wh{\sigma}_{i}^2\equiv \wh{\sigma}_{i}^{2(k+1)}$, 
and $\wh{\bm\Sigma}_n^\xi\equiv\wh{\bm\Sigma}_n^{\xi(k+1)}$,
for any $k\ge 0$. Under Assumptions \ref{ass:common}, \ref{ass:idio}, \ref{ass:ind}, \ref{ass:linear}, \ref{ass:tails}, and \ref{ass:ident}, as $n,T\to\infty$:
\begin{compactenum}[(i)]
\item $\max_{t=1,\ldots, T}\Vert \wh{\mbf P}_{t|t-1}\Vert=O_p(1)$; 
\item $\max_{t=1,\ldots, T}\Vert(\wh{\mbf P}_{t|t-1})^{-1}\Vert =O_p(1)$;
\item $\max_{t=1,\ldots, T}n\Vert \wh{\mbf P}_{t|t}\Vert=O_p(1)$;
\item $\max_{t=1,\ldots, T}n\Vert \wh{\mbf P}_{t|T}\Vert=O_p(1)$.
\end{compactenum}
\end{lem}

\noindent
\textsc{Proof.} The proof is the same as the proof of Lemma  \ref{lem:PPOnstar}  but using Proposition \ref{prop:load}(a) (which does not require this lemma to be proved) and Lemmas \ref{lem:sigmaunifhat} and \ref{lem:eststar_LASTHAT} instead of Lemmas \ref{lem:starT12}, \ref{lem:sigmaunif}, and \ref{lem:eststar_LAST}. $\Box$

\begin{lem}\label{lem:FFOnhat}
Consider the Kalman filter and smoother estimators $\wh{\mbf F}_{t|t-1}\equiv  {\mbf F}_{t|t-1}^{(k+1)}$, $\wh{\mbf F}_{t|t}\equiv  {\mbf F}_{t|t}^{(k+1)}$, and $\wh{\mbf F}_{t|T}\equiv  {\mbf F}_{t|T}^{(k+1)}$ obtained using the EM estimators $\wh{\bm\lambda}_{i}\equiv \wh{\bm\lambda}_{i}^{(k+1)}$, 
$\wh{\bm\Lambda}_{n}\equiv \wh{\bm\Lambda}_{n}^{(k+1)}$,
$\wh{\sigma}_{i}^2\equiv \wh{\sigma}_{i}^{2(k+1)}$, 
and $\wh{\bm\Sigma}_n^\xi\equiv\wh{\bm\Sigma}_n^{\xi(k+1)}$,
for any $k\ge 0$.
 Under Assumptions \ref{ass:common}, \ref{ass:idio}, \ref{ass:ind}, \ref{ass:linear}, \ref{ass:tails}, and \ref{ass:ident}, as $n,T\to\infty$:
\begin{compactenum}[(i)]
\item for all $s=0,\ldots, T$, $\Vert \wh{\mbf F}_{t|s}\Vert=O_p(1)$, uniformly in $t\le s$;
\item $n\Vert \wh{\mbf F}_{t|T}-\wh{\mbf F}_{t|t}\Vert=O_p(1)$, uniformly in $t$;
\item $n \Vert \wh{\mbf F}_{t|t}-\wh{\mbf F}_{t}^{\text{\tiny \upshape {WLS}}}\Vert=O_p(1)$, uniformly in  $t$;
\end{compactenum}
where $\wh{\mbf F}_{t}^{\text{\tiny \upshape {WLS}}}=(\wh{\bm\Lambda}_n^{\prime}(\wh{\bm\Sigma}_n^{\xi})^{-1}\wh{\bm\Lambda}_n)^{-1}\wh{\bm\Lambda}_n^{\prime}(\wh{\bm\Sigma}_n^{\xi})^{-1}\mbf x_{nt}.$
\end{lem}

\noindent
\textsc{Proof.} The proof is the same as the proof of Lemma \ref{lem:FFOnstar} but using 
Proposition \ref{prop:load}(a) (which does not require this lemma to be proved) and Lemmas \ref{lem:sigmaunifhat}, \ref{lem:eststar_LASTHAT}, and \ref{lem:PPOnhat} instead of 
 Lemmas \ref{lem:starT12}, \ref{lem:sigmaunif}, \ref{lem:eststar_LAST}, and \ref{lem:PPOn}. $\Box$

\begin{lem}\label{lem:FFOnhatunif}
Consider the Kalman filter and smoother estimators $\wh{\mbf F}_{t|t-1}\equiv  {\mbf F}_{t|t-1}^{(k+1)}$, $\wh{\mbf F}_{t|t}\equiv  {\mbf F}_{t|t}^{(k+1)}$, and $\wh{\mbf F}_{t|T}\equiv  {\mbf F}_{t|T}^{(k+1)}$ obtained using the EM estimators $\wh{\bm\lambda}_{i}\equiv \wh{\bm\lambda}_{i}^{(k+1)}$, 
$\wh{\bm\Lambda}_{n}\equiv \wh{\bm\Lambda}_{n}^{(k+1)}$,
$\wh{\sigma}_{i}^2\equiv \wh{\sigma}_{i}^{2(k+1)}$, 
and $\wh{\bm\Sigma}_n^\xi\equiv\wh{\bm\Sigma}_n^{\xi(k+1)}$,
for any $k\ge 0$.
 Under Assumptions \ref{ass:common}, \ref{ass:idio}, \ref{ass:ind}, \ref{ass:linear}, \ref{ass:tails}, and \ref{ass:ident}, as $n,T\to\infty$:
\begin{compactenum}[(i)]
\item for all $s=0,\ldots, T$, $\log^{-1/\delta_v}T\, \max_{t=1,\ldots, s}\Vert \wh{\mbf F}_{t|s}\Vert=O_p(1)$;
\item $\log^{-1/\delta_v}T\, \max_{t=1,\ldots, T}\,n\Vert \wh{\mbf F}_{t|T}-\wh{\mbf F}_{t|t}\Vert=O_p(1)$; 
\item $\log^{-1/\delta_v}T\, \max_{t=1,\ldots, T}\,n \Vert \wh{\mbf F}_{t|t}-\wh{\mbf F}_{t}^{\text{\tiny \upshape {WLS}}}\Vert=O_p(1)$;
\end{compactenum}
where $\wh{\mbf F}_{t}^{\text{\tiny \upshape {WLS}}}=(\wh{\bm\Lambda}_n^{\prime}(\wh{\bm\Sigma}_n^{\xi})^{-1}\wh{\bm\Lambda}_n)^{-1}\wh{\bm\Lambda}_n^{\prime}(\wh{\bm\Sigma}_n^{\xi})^{-1}\mbf x_{nt}.$
\end{lem}

\noindent
\textsc{Proof.} From \eqref{eq:ponti} in the proof of Lemma \ref{lem:FFO1}, but when computed using the EM estimator of the parameters, we see that the only modification is that we need to use $\max_{t=1,\ldots, T}n^{-1/2}\Vert\mbf x_t\Vert=O_p()$ from Lemma \ref{lem:xunif2}. This proves part (i).

Part (ii) follows from part (i) and \eqref{eq:ponti2} in the proof of Lemma \ref{lem:FFOn}, but when computed using the EM estimator of the parameters.

Part (iii) follows from \eqref{eq:girandola} in the proof of Lemma \ref{lem:gennaio22bis}, but when computed using the EM estimator of the parameters, and using Lemma \ref{lem:xunif2} again
and the fact that $\wh{\mbf F}_{t-1|t-1}$ is a weighted average of $\mbf x_{n1},\ldots,\mbf x_{n,t-1}$. This completes the proof. $\Box$

\section{Lemmas necessary for proving Proposition \ref{prop:altri}}
\begin{lem}\label{lem:est0_24}
Under Assumptions \ref{ass:common}, \ref{ass:idio}, \ref{ass:ind}, \ref{ass:tails}, and \ref{ass:ident}, as $n,T\to\infty$,
 $$\min(n,\sqrt T\log^{-1/2} n)\,\max_{i=1,\ldots, n}\vert\wh{\sigma}_{i}^{(0)2}-\sigma_i^2\vert = O_p(1).$$ 
 \end{lem}

\noindent\textsc{Proof.} From \eqref{eq:laringite} in the proof of Lemma \ref{lem:est0}, we have
\al{
\max_{i=1,\ldots, n}\vert\wh{\sigma}_i^{2(0)}- \sigma_i^2\vert\le &\, \max_{i=1,\ldots, n}\l\vert T^{-1}\sum_{t=1}^T  x_{it}^2-\E[ x_{it}^2] \r\vert+ 
\max_{i=1,\ldots, n}\l \vert
\wh{\bm\lambda}_i^{(0)\prime}\wh{\bm\lambda}_i^{(0)}- \bm\lambda_i^\prime\bm\lambda_i\r\vert\nn\\
&+2\max_{i=1,\ldots, n}\Vert \bm\lambda_i\Vert\, \max_{i=1,\ldots, n} \Vert \wh{\bm\lambda}_i^{(0)}  -\bm\lambda_i\Vert\nn\\
&+2\max_{i=1,\ldots, n}\Vert \bm\lambda_i\Vert\, \l\Vert T^{-1}\sum_{t=1}^T\mbf F_t {\mbf F}_t^\prime-\mbf I_r\r\Vert\, \l\{\max_{i=1,\ldots, n}\Vert\wh{\bm\lambda}_i^{(0)}-\bm\lambda_i\Vert+\max_{i=1,\ldots, n}\Vert\bm\lambda_i\Vert \r\} \nn\\
&+2\max_{i=1,\ldots, n} \Vert \bm\lambda_i\Vert\, \l\Vert T^{-1}\sum_{t=1}^T \mbf F_t (\wt{\mbf F}_t-\mbf F_t)^\prime\r\Vert\,\l\{\max_{i=1,\ldots, n} \Vert \wh{\bm\lambda}_i^{(0)}-\bm\lambda_i\Vert +\max_{i=1,\ldots, n}\Vert\bm\lambda_i\Vert\r\}\nn\\
& +2\max_{i=1,\ldots, n} \l\Vert T^{-1}\sum_{t=1}^T \xi_{it}{\mbf F}_t^\prime\r\Vert\,\max_{i=1,\ldots, n}\Vert{\bm\lambda}_i \Vert\nn\\
&+2\max_{i=1,\ldots, n}\l\Vert T^{-1}\sum_{t=1}^T \xi_{it}{\mbf F}_t^\prime \r\Vert\,\max_{i=1,\ldots, n}\Vert\wh{\bm\lambda}_i^{(0)}-\bm\lambda_i  \Vert\nn\\
& +2\max_{i=1,\ldots, n} \l\Vert T^{-1}\sum_{t=1}^T \xi_{it}(\wt{\mbf F}_t-\mbf F_t)^\prime\r\Vert \, \l\{\max_{i=1,\ldots, n} \Vert \wh{\bm\lambda}_i^{(0)} -\bm\lambda_i\Vert+\max_{i=1,\ldots, n} \Vert\bm\lambda_i\Vert\r\}\nn\\
=&\,O_p(\max(n^{-1},T^{-1/2}\sqrt{\log n})),\nn
}
by  Assumption \ref{ass:common}(a), Lemma \ref{lem:tail}(ii) and the union bound, 
and \eqref{eq:unifload3} and \eqref{eq:MINKIEVA} in the proof of Lemma \ref{lem:sigmaunif}. This completes the proof. $\Box$

\begin{lem}\label{lem:ultimo}
Under Assumptions \ref{ass:common}, \ref{ass:idio}, \ref{ass:ind}, \ref{ass:tails}, and \ref{ass:ident}, as $n,T\to\infty$, for any $k\ge 0$,
\begin{compactenum}
\item [(i)] $\min(n\log^{-2/\delta_v}T,\sqrt{nT}{\log^{-1/2} n}, T{\log^{-1/2} n}),\nn\, \Vert T^{-1}\sum_{t=1}^T(\mbf F_{t|T}^{(k)}-\mbf F_t)\mbf F_t^\prime\Vert=O_p(1)$;
\item [(ii)] $\min(n\log^{-2/\delta_v}T,\sqrt{nT}{\log^{-1/2} n}, T{\log^{-1/2} n})\, \Vert T^{-1}\sum_{t=1}^T(\mbf F_{t|T}^{(k)}-\mbf F_t) \xi_{it}\Vert=O_p(1)$, uniformly in $i$.
\end{compactenum}
 \end{lem}
 
 \noindent\textsc{Proof.} Throughout, let $\bm y_t=\mbf F_t$ or $\bm y_t = \xi_{it}$.
 Then, for any $k\ge 0$ we have to consider
\al{
\l\Vert T^{-1} \sum_{t=1}^T (\mbf F_{t|T}^{(k)}-\mbf F_{t})\bm y_t^\prime \r\Vert\le&\, 
\l\Vert T^{-1} \sum_{t=1}^T (\mbf F_{t|T}^{(k)}-\mbf F_{t|t}^{(k)})\bm y_t^\prime \r\Vert\nn\\
&+\l\Vert T^{-1} \sum_{t=1}^T (\mbf F_{t|t}^{(k)}-\wh{\mbf F}_{t}^{\text{\tiny WLS}(k)})\bm y_t^\prime \r\Vert\nn\\
&+\l\Vert T^{-1} \sum_{t=1}^T (\wh{\mbf F}_{t}^{\text{\tiny WLS}(k)}-\mbf F_t)\bm y_t^\prime \r\Vert\nn\\
=&\, \mathcal I+\mathcal {II}+\mathcal {III}, \;\text{ say.}
\label{eq:tucano}
}
where $\wh{\mbf F}_{t}^{\text{\tiny \upshape {WLS}}(k)}=(\wh{\bm\Lambda}_n^{(k)\prime}(\wh{\bm\Sigma}_n^{\xi(k)})^{-1}\wh{\bm\Lambda}_n^{(k)})^{-1}\wh{\bm\Lambda}_n^{(k)\prime}(\wh{\bm\Sigma}_n^{\xi(k)})^{-1}\mbf x_{nt}.$
Let us consider each term in \eqref{eq:tucano}. First, when $k=0$ we have $\mathcal I=O_p(n^{-1})$ by \eqref{eq:I00} in the proof of Proposition \ref{prop:plainvanilla} (which does not require this lemma to be proved), while, when $k\ge 1$
\al{
\mathcal I\le&\,\max_{t=1,\ldots, T} \Vert \mbf F_{t|T}^{(k)}-\mbf F_{t|t}^{(k)}\Vert \,\l\Vert T^{-1}\sum_{t=1}^T \bm y_t\r\Vert = O_p(n^{-1}\log^{1/\delta_v}T),\label{eq:I00tuc}
}
by Lemma \ref{lem:FFOnhatunif}(ii) and since 
\beq
\E\l[\l\Vert T^{-1}\sum_{t=1}^T \mbf F_t\r\Vert^2\r]\le r T^{-2} \max_{j=1,\ldots, r} \sum_{t,s=1}^T \vert \E[F_{jt}F_{js}]\vert \le 1,\label{bollini1}
\eeq
by Cauchy-Schwarz inequality and Assumption \ref{ass:ident}(b), and also
\beq
\E\l[\l\vert T^{-1}\sum_{t=1}^T \xi_{it}\r\vert^2\r]\le T^{-2} \sum_{t,s=1}^T \vert \E[\xi_{it}\xi_{is}]\vert \le T^{-1}M_3,\label{bollini2}
\eeq
by Lemma \ref{lem:Gxi}(iii).
Second, when $k=0$ we have $\mathcal {II}=O_p(n^{-1})$ by \eqref{eq:II00} in the proof of Proposition \ref{prop:plainvanilla} (which does not require this lemma to be proved), while, when $k\ge 1$
\al{
\mathcal {II}\le&\,\max_{t=1,\ldots, T} \Vert \mbf F_{t|t}^{(k)}-\wh{\mbf F}_{t}^{\text{\tiny WLS}(k)}\Vert \,\l\Vert T^{-1}\sum_{t=1}^T \bm y_t\r\Vert = O_p(n^{-1}\log^{1/\delta_v}T),\label{eq:II00tuc}
}
by Lemma \ref{lem:FFOnhatunif}(iii), \eqref{bollini1}, and \eqref{bollini2}.

Third, 
\al{
\mathcal {III}
\le&\, 
\Vert(\bm\Lambda_n^\prime(\bm\Sigma_n^\xi)^{-1}\bm\Lambda_n)^{-1}\Vert \, \Vert \bm\Lambda_n^\prime(\bm\Sigma_n^\xi)^{-1} (\bm\Lambda_n-\wh{\bm\Lambda}_n^{(k)}) \Vert \l\Vert T^{-1}\sum_{t=1}^T \mbf F_t\bm y_t^\prime\r\Vert\nn\\ 
&+\Vert
n(\wh{\bm\Lambda}_n^{(k)\prime}(\wh{\bm\Sigma}_n^{\xi(k)})^{-1}\wh{\bm\Lambda}_n^{(k)})^{-1}n^{-1/2}\wh{\bm\Lambda}_n^{(k)\prime}(\wh{\bm\Sigma}_n^{\xi(k)})^{-1}-n({\bm\Lambda}_n^{\prime}({\bm\Sigma}_n^{\xi})^{-1}{\bm\Lambda}_n)^{-1}n^{-1/2}{\bm\Lambda}_n^{\prime}({\bm\Sigma}_n^{\xi})^{-1}\Vert \nn\\
&\cdot n^{-1/2}\Vert \bm\Lambda_n-\wh{\bm\Lambda}_n^{(k)}\Vert \, \l\Vert T^{-1}\sum_{t=1}^T \mbf F_t\bm y_t^\prime\r\Vert\nn\\ 
&+n\Vert(\wh{\bm\Lambda}_n^{(k)\prime}(\wh{\bm\Sigma}_n^{\xi(k)})^{-1}\wh{\bm\Lambda}_n^{(k)})^{-1}\Vert \,n^{-1}\l\Vert T^{-1}\sum_{t=1}^T \bm\Lambda_n^\prime(\bm\Sigma_n^\xi)^{-1}\bm\xi_{nt}\bm y_t^\prime \r\Vert\nn\\
&+n \Vert
(\wh{\bm\Lambda}_n^{(k)\prime}(\wh{\bm\Sigma}_n^{\xi(k)})^{-1} \wh{\bm\Lambda}_n^{(k)})^{-1}\Vert\, n^{-1}\l\Vert T^{-1}\sum_{t=1}^T \{\wh{\bm\Lambda}_n^{(k)\prime}(\wh{\bm\Sigma}_n^{\xi(k)})^{-1}-{\bm\Lambda}_n^{\prime}({\bm\Sigma}_n^{\xi})^{-1}\} \bm\xi_{nt}\bm y_t^\prime\r\Vert\nn\\ 
=&\, \mathcal{III}_a+\mathcal{III}_b+\mathcal{III}_c+\mathcal{III}_d,\;\text{say.} \label{eq:III0tuc}
}
Then, when $k=0$, $\mathcal{III}_a=III_a=O_p(n^{-1/2}T^{-1/2})$, $\mathcal{III}_b=III_b=O_p(\max(n^{-2},T^{-1}))$, and $\mathcal{III}_c=III_c=O_p(n^{-1/2}T^{-1/2})$, because of  \eqref{eq:IIIa}, \eqref{eq:IIIb}, \eqref{eq:IIIc1}, and \eqref{eq:IIIc2}  in the proof of Proposition \ref{prop:plainvanilla}  (which does not require this lemma to be proved).
While, if $k\ge 1$,
\al{
\mathcal{III}_a=O_p(\max(n^{-1}\log^{2/\delta_v}T,n^{-1/2}T^{-1/2}\sqrt{\log n},T^{-1})),\label{eq:III0tucA}
}
by \eqref{eq:ultimopezzoAAhat}, \eqref{eq:ultimopezzoAA1hat}, and \eqref{eq:ultimopezzoAA0hat} in the proof of Proposition \ref{prop:factors}  (which does not require this lemma to be proved), and since
$\Vert T^{-1}\sum_{t=1}^T \mbf F_t\mbf F_t^\prime\Vert = O_p(1)$ by Lemma \ref{lem:consistCOV}(i) and Assumption \ref{ass:ident}(b), and  $\Vert T^{-1}\sum_{t=1}^T \mbf F_t\xi_{it}\Vert = O_p(T^{-1/2})$ by Lemma \ref{lem:consistCOV}(ii). These last arguments, and \eqref{eq:ultimopezzoBBhat} in the proof of Proposition \ref{prop:factors}  (which does not require this lemma to be proved)
imply also that, if $k\ge 1$,
\al{
\mathcal{III}_b=O_p(\max(n^{-2}\log^{4/\delta_v}T,T^{-1} \sqrt{\log n}, n^{-1}T^{-1/2} \log^{2/\delta_v}T\sqrt{\log n})).\label{eq:III0tucB}
}
Moreover, if $\bm y_t=\mbf F_t$ and $k\ge 1$, then
\al{
\mathcal{III}_c=&\, n\Vert(\wh{\bm\Lambda}_n^{(k)\prime}(\wh{\bm\Sigma}_n^{\xi(k)})^{-1}\wh{\bm\Lambda}_n^{(k)})^{-1}\Vert \,n^{-1}\l\Vert T^{-1}\sum_{t=1}^T \bm\Lambda_n^\prime(\bm\Sigma_n^\xi)^{-1}\bm\xi_{nt}\mbf F_t^\prime \r\Vert\nn\\
=&\, O_p(n^{-1/2}T^{-1/2}),\label{eq:III0tucC1}
}
by Lemmas \ref{lem:eststar_LASTHAT}(iii) and \ref{lem:COVFF0}(iv). While if $\bm y_t=\xi_{it}$ and $k\ge 1$, then
\al{
\mathcal{III}_c=&\, n\Vert(\wh{\bm\Lambda}_n^{(k)\prime}(\wh{\bm\Sigma}_n^{\xi(k)})^{-1}\wh{\bm\Lambda}_n^{(k)})^{-1}\Vert \,n^{-1}\l\Vert T^{-1}\sum_{t=1}^T \bm\Lambda_n^\prime(\bm\Sigma_n^\xi)^{-1}\bm\xi_{nt}\xi_{it} \r\Vert\nn\\
=&\, O_p(n^{-1/2}T^{-1/2}),\label{eq:III0tucC2}
}
by Lemmas \ref{lem:eststar_LASTHAT}(iii) and \ref{lem:COVFF0}(v).

Similarly, if $\bm y_t=\mbf F_t$, we have:
\al{
\mathcal {III}_d =&\, 
n \Vert(\wh{\bm\Lambda}_n^{(k)\prime}(\wh{\bm\Sigma}_n^{\xi(k)})^{-1} \wh{\bm\Lambda}_n^{(k)})^{-1}\Vert\, 
n^{-1}\l\Vert T^{-1}\sum_{t=1}^T \{\wh{\bm\Lambda}_n^{(k)\prime}(\wh{\bm\Sigma}_n^{\xi(k)})^{-1}-{\bm\Lambda}_n^{\prime}({\bm\Sigma}_n^{\xi})^{-1}\} \bm\xi_{nt}\mbf F_t^\prime\r\Vert\nn\\
\le &\,n \Vert(\wh{\bm\Lambda}_n^{(k)\prime}(\wh{\bm\Sigma}_n^{\xi(k)})^{-1} \wh{\bm\Lambda}_n^{(k)})^{-1}\Vert\, 
n^{-1/2} \Vert \wh{\bm\Lambda}_n^{(k)\prime}(\wh{\bm\Sigma}_n^{\xi(k)})^{-1}-{\bm\Lambda}_n^{\prime}({\bm\Sigma}_n^{\xi})^{-1}\Vert\, n^{-1/2}\l\Vert  T^{-1}\sum_{t=1}^T\bm\xi_{nt}\mbf F_t^\prime\r\Vert\nn\\
=&\, \l\{\ba{lll}
 O_p(\max(n^{-1}T^{-1/2},T^{-1})), &\text{ if }& k=0,\\
 O_p(\max(n^{-1}T^{-1/2}\log^{2/\delta_v}T ,T^{-1}\sqrt{\log n})),&\text{ if }& k\ge 1,
 \ea
\r.\label{eq:III0tucD1}
}
where, if $k=0$, we used Lemmas \ref{lem:est0_LAST}(ii), \ref{lem:est0_LAST}(iii), and \ref{lem:consistCOV}(iii), while, if $k\ge 1$, we used Lemmas \ref{lem:eststar_LASTHAT}(ii), \ref{lem:eststar_LASTHAT}(iii), and \ref{lem:consistCOV}(iii).

Finally, if $\bm y_t=\xi_{it}$, we have:
\al{
\mathcal{III}_d =&\, 
n \Vert(\wh{\bm\Lambda}_n^{(k)\prime}(\wh{\bm\Sigma}_n^{\xi(k)})^{-1} \wh{\bm\Lambda}_n^{(k)})^{-1}\Vert\, 
n^{-1}\l\Vert T^{-1}\sum_{t=1}^T \{\wh{\bm\Lambda}_n^{(k)\prime}(\wh{\bm\Sigma}_n^{\xi(k)})^{-1}-{\bm\Lambda}_n^{\prime}({\bm\Sigma}_n^{\xi})^{-1}\} \bm\xi_{nt}\xi_{it}\r\Vert\nn\\
\le&\, n \Vert(\wh{\bm\Lambda}_n^{(k)\prime}(\wh{\bm\Sigma}_n^{\xi(k)})^{-1} \wh{\bm\Lambda}_n^{(k)})^{-1}\Vert\, 
\bigg\{n^{-1}\bigg\Vert T^{-1}\sum_{t=1}^T (\wh{\bm\Lambda}_n^{(k)}-{\bm\Lambda}_n)^{\prime}({\bm\Sigma}_n^{\xi})^{-1} \bm\xi_{nt}\xi_{it}\bigg\Vert\nn\\
&+n^{-1}\bigg\Vert T^{-1}\sum_{t=1}^T{\bm\Lambda}_n^\prime 
\{(\wh{\bm\Sigma}_n^{\xi(k)})^{-1}-({\bm\Sigma}_n^{\xi})^{-1}\}\bm\xi_{nt}\xi_{it}\bigg\Vert\nn\\
&+ n^{-1}\bigg\Vert T^{-1}\sum_{t=1}^T (\wh{\bm\Lambda}_n^{(k)}-{\bm\Lambda}_n)^{\prime}\{(\wh{\bm\Sigma}_n^{\xi(k)})^{-1}-({\bm\Sigma}_n^{\xi})^{-1}\}\bm\xi_{nt}\xi_{it}\bigg\Vert
\bigg\}\nn\\
=&\, n \Vert(\wh{\bm\Lambda}_n^{(k)\prime}(\wh{\bm\Sigma}_n^{\xi(k)})^{-1} \wh{\bm\Lambda}_n^{(k)})^{-1}\Vert\, 
\{
\mathcal{III}_{d.1}+\mathcal{III}_{d.2}+\mathcal{III}_{d.3}
\},\;\text{say.}\label{eq:III0tucD2temp}
}
Now,
\al{
\mathcal{III}_{d.1}=&\, n^{-1}\l\Vert T^{-1}\sum_{t=1}^T ({\bm\Lambda}_n^{\text{\tiny OLS}}-{\bm\Lambda}_n)^{\prime}({\bm\Sigma}_n^{\xi})^{-1} \bm\xi_{nt}\xi_{it}\r\Vert
\nn\\
&+n^{-1/2} \Vert \wh{\bm\Lambda}_n^{(k)}-{\bm\Lambda}_n^{\text{\tiny OLS}}\Vert \,\Vert ({\bm\Sigma}_n^{\xi})^{-1}\Vert \,  n^{-1/2}
\l\Vert T^{-1}\sum_{t=1}^T
 \bm\xi_{nt}\xi_{it}
\r\Vert\nn\\
=&\, \mathcal{III}_{d.1.1}+\mathcal{III}_{d.1.2}, \text{say.}
\label{eq:III0tucD21temp}
}
And,
\al{
\mathcal{III}_{d.1.1}=&\, n^{-1}\l\Vert T^{-1}\sum_{t=1}^T \l(T^{-1}\sum_{s=1}^T \mbf F_s\mbf F_s^\prime\r)^{-1}\l(
T^{-1}\sum_{s=1}^T \mbf F_s\bm\xi_{ns}^\prime\r)({\bm\Sigma}_n^{\xi})^{-1} \bm\xi_{nt}\xi_{it}\r\Vert\nn\\
\le &\, \l\Vert  \l(T^{-1}\sum_{s=1}^T \mbf F_s\mbf F_s^\prime\r)^{-1}\r\Vert\, n^{-1}\l\Vert T^{-2}\sum_{s,t=1}^T \mbf F_s\bm\xi_{ns}^\prime({\bm\Sigma}_n^{\xi})^{-1} \bm\xi_{nt}\xi_{it}\r\Vert\nn\\
\le &\, \l\Vert  \l(T^{-1}\sum_{s=1}^T \mbf F_s\mbf F_s^\prime\r)^{-1}\r\Vert\, n^{-1}\l\Vert T^{-2}\sum_{s,t=1}^T \mbf F_s\bm\xi_{ns}^\prime({\bm\Sigma}_n^{\xi})^{-1}\bm\xi_{nt} \r\Vert\, \max_{t=1,\ldots, T} \vert \xi_{it}\vert\nn\\
=&\, O_p(n^{-1/2}T^{-3/2}\log^{1/\delta_v}T),
\label{eq:III0tucD211}
}
by  
\eqref{eq:normaxiunif} in the proof of Lemma \ref{lem:xunif2}, and Lemmas  \ref{lem:COVFF0}(vi) and \ref{lem:frida}.

\al{
\mathcal{III}_{d.1.2}=&\,\l\{
\ba{lll}
O_p(\max(n^{-1},n^{-1/2}T^{-1/2},T^{-1})), &\text{ if }& k=0,\\
O_p(\max(n^{-1}\log^{2/\delta_v}T,n^{-1/2}T^{-1/2}\sqrt{\log n},T^{-1})),&\text{ if }& k\ge 1,
\ea
\r.
\label{eq:III0tucD212}
}
where we used Assumption \ref{ass:idio}(a), the last relation in \eqref{eq:varibounds} in the proof of Proposition \ref{prop:plainvanilla} (which does not require this lemma to be proved), and, if $k=0$, we used also \citet[Corollary 1 and Proposition B.3]{MBPCAQML}, or, if $k\ge 1$, we used \eqref{eq:reggioemilia} in the proof of Proposition \ref{prop:factors} (which does not require this lemma to be proved).

By using \eqref{eq:III0tucD211} and \eqref{eq:III0tucD212} into \eqref{eq:III0tucD21temp}, we have
\al{
\mathcal{III}_{d.1}=&\,O_p(n^{-1/2}T^{-3/2}\log^{1/\delta_v}T)\nn\\
&+\l\{
\ba{lll}
O_p(\max(n^{-1},n^{-1/2}T^{-1/2},T^{-1})), &\text{ if }& k=0,\\
O_p(\max(n^{-1}\log^{2/\delta_v}T,n^{-1/2}T^{-1/2}\sqrt{\log n},T^{-1})),&\text{ if }& k\ge 1,
\ea
\r.
\label{eq:III0tucD21}
}
Furthermore,
\al{
\mathcal{III}_{d.2}=&\, n^{-1}\bigg\Vert T^{-1}\sum_{t=1}^T\sum_{j=1}^n {\bm\lambda}_j 
\{(\wh{\sigma}_j^{2(k)})^{-1}-({\sigma}_j^{2})^{-1}\}\xi_{jt}\xi_{it}\bigg\Vert\nn\\
\le &\, C_\xi^2 \max_{j=1,\ldots, n} \vert\wh{\sigma}_j^{2(k)}-{\sigma}_j^{(k)} \vert 
n^{-1}\bigg\Vert T^{-1}\sum_{t=1}^T\sum_{j=1}^n {\bm\lambda}_j \xi_{jt}\xi_{it}\bigg\Vert\nn\\
\le &\, C_\xi^2 \max_{j=1,\ldots, n} \vert\wh{\sigma}_j^{2(k)}-{\sigma}_j^{(k)} \vert 
n^{-1}\bigg\Vert T^{-1}\sum_{t=1}^T {\bm\Lambda}_n^\prime \bm\xi_{nt}\xi_{it}\bigg\Vert\nn\\
=&\, O_p(n^{-1/2}T^{-1/2})\cdot \l\{\ba{lll}
O_p(\max(n^{-1},T^{-1/2}\sqrt{\log n})), &\text{ if }& k=0,\\
O_p(\max(n^{-1}\log^{2/\delta_v}T,T^{-1/2}\sqrt{\log n})), &\text{ if }& k\ge 1,
\ea
\r.
\label{eq:III0tucD22}
}
where, if $k=0$, we used Lemmas \ref{lem:est0_24} and \ref{lem:COVFF0}(ii), while, if $k\ge 1$, we used Lemmas \ref{lem:sigmaunifhat}(ii) and \ref{lem:COVFF0}(ii). Term $\mathcal{III}_{d.3}$ is dominated by $\mathcal{III}_{d.1}$.

By using Lemmas \ref{lem:eststar_LASTHAT}(iii) or \ref{lem:est0_LAST}(iii), together with \eqref{eq:III0tucD21} and \eqref{eq:III0tucD22} into \eqref{eq:III0tucD2temp}, we have that if $\bm y_t=\xi_{it}$, then
\al{
\mathcal{III}_d =&\, O_p(n^{-1/2}T^{-3/2}\log^{1/\delta_v}T)+O_p(\max(n^{-1}\log^{2/\delta_v}T,n^{-1/2}T^{-1/2}\sqrt{\log n},T^{-1}))\nn\\
&+O_p(\max(n^{-3/2}T^{-1/2}\log^{2/\delta_v}T,n^{-1/2}T^{-1}\sqrt{\log n})).
\label{eq:III0tucD2}
}

By substituting \eqref{eq:III0tucA}, \eqref{eq:III0tucB}, \eqref{eq:III0tucC1}, \eqref{eq:III0tucC2}, \eqref{eq:III0tucD1}, and \eqref{eq:III0tucD2} into \eqref{eq:III0tuc}, if $\bm y_t=\mbf F_t$, we have
\al{
\mathcal {III}=&\, O_p(\max(n^{-1}\log^{2/\delta_v}T,n^{-1/2}T^{-1/2}\sqrt{\log n},T^{-1}))\nn\\
&+O_p(\max(n^{-2}\log^{4/\delta_v}T,T^{-1} \sqrt{\log n}, n^{-1}T^{-1/2} \log^{2/\delta_v}T\sqrt{\log n}))+ O_p(n^{-1/2}T^{-1/2})\nn\\
&+ O_p(\max(n^{-1}T^{-1/2}\log^{2/\delta_v}T ,T^{-1}\sqrt{\log n})),\label{eq:matisse}
}
while, if $\bm y_t=\xi_{it}$, we have
\al{
\mathcal {III}=&\, O_p(\max(n^{-1}\log^{2/\delta_v}T,n^{-1/2}T^{-1/2}\sqrt{\log n},T^{-1}))\nn\\
&+O_p(\max(n^{-2}\log^{4/\delta_v}T,T^{-1} \sqrt{\log n}, n^{-1}T^{-1/2} \log^{2/\delta_v}T\sqrt{\log n}))+ O_p(n^{-1/2}T^{-1/2})\nn\\
&+O_p(n^{-1/2}T^{-3/2}\log^{1/\delta_v}T)+O_p(\max(n^{-1}\log^{2/\delta_v}T,n^{-1/2}T^{-1/2}\sqrt{\log n},T^{-1}))\nn\\
&+O_p(\max(n^{-3/2}T^{-1/2}\log^{2/\delta_v}T,n^{-1/2}T^{-1}\sqrt{\log n})).\label{eq:matisse2}
}

To conclude by substituting \eqref{eq:I00tuc}, \eqref{eq:II00tuc}, and either \eqref{eq:matisse} or \eqref{eq:matisse2} into \eqref{eq:tucano}, we have
\al{
\l\Vert T^{-1} \sum_{t=1}^T (\mbf F_{t|T}^{(k)}-\mbf F_{t})\mbf F_t^\prime \r\Vert=&\,O_p(\max(n^{-1}\log^{2/\delta_v}T, n^{-1/2}T^{-1/2}\sqrt{\log n}, T^{-1}\sqrt{\log n})),\nn
}
which proves part (i), and 
\al{
\l\Vert T^{-1} \sum_{t=1}^T (\mbf F_{t|T}^{(k)}-\mbf F_{t})\xi_{it} \r\Vert=&\,O_p(\max(n^{-1}\log^{2/\delta_v}T, n^{-1/2}T^{-1/2}\sqrt{\log n}, T^{-1}\sqrt{\log n})),\nn
}
which proves part (ii) and completes the proof. $\Box$

\section{Lemmas necessary for proving Proposition \ref{prop:eff}}

\begin{lem}\label{lem:PCAF}
Consider the initial estimator of the factors $\wt{\mbf F}_t$ defined in Section \ref{app:prest}, then, under Assumptions \ref{ass:common}, \ref{ass:idio}, \ref{ass:ind}, and \ref{ass:ident},
as $n,T\to\infty$, if $\sqrt n/T\to 0$, 
    \[
    \sqrt n (\wt{\mbf F}_{t}-\mbf F_t)\stackrel{d}{\to}\mathcal N(\mbf 0_r,\bm{\mathcal W}_t^{\text{\tiny \upshape PC}}),
    \]
    for any given $t=1,\ldots, T$, where 
    $$
    \bm{\mathcal W}_t^{\text{\tiny \upshape PC}}=(\bm\Sigma_{\Lambda})^{-1} \l(\lim_{n\to\infty} n^{-1}\sum_{i,j=1}^n {\E[\xi_{it}\xi_{jt}] \bm\lambda_i\bm\lambda_j  } \r)(\bm\Sigma_{\Lambda})^{-1}
    $$ 
    with $\bm\Sigma_{\Lambda}=\lim_{n\to\infty} n^{-1}\sum_{i=1}^n \bm\lambda_i \bm\lambda_i^\prime$;
\end{lem}

\noindent
\textsc{Proof.}
By definition of the pre-estimator in Section \ref{app:prest}:
\al{
\wt{\mbf F}_t-\mbf F_t=&\,(\wh{\bm\Lambda}_n^{(0)\prime}\wh{\bm\Lambda}_n^{(0)})^{-1}\wh{\bm\Lambda}_n^{(0)\prime}\mbf x_{nt}-\mbf F_t\nn\\
=&\, 
n(\mbf M_n^\chi)^{-1}
\l\{n^{-1}\wh{\bm\Lambda}_n^{(0)\prime}(\bm\Lambda_n-\wh{\bm\Lambda}_n^{(0)}) \mbf F_t+n^{-1}(\wh{\bm\Lambda}_n^{(0)}-\bm\Lambda_n)^\prime\bm\xi_{nt}
+n^{-1}\bm\Lambda_n^\prime\bm \xi_{nt}
\r\}\label{eq:panzona}\\
&+
n\l\{(\wh{\mbf M}_n^x)^{-1}-(\mbf M_n^\chi)^{-1}\r\}
\l\{n^{-1}\wh{\bm\Lambda}_n^{(0)\prime}(\bm\Lambda_n-\wh{\bm\Lambda}_n^{(0)}) \mbf F_t+n^{-1}(\wh{\bm\Lambda}_n^{(0)}-\bm\Lambda_n)^\prime\bm\xi_{nt}
+n^{-1}\bm\Lambda_n^\prime\bm \xi_{nt}
\r\}.\nn
}
Then,
\al{
n \Vert(\mbf M_n^\chi)^{-1}\Vert = n\{\mu_{nr}^\chi \}^{-1}\ge \underline C_r,\label{eq:panza1}
}
by Lemma \ref{lem:Gxi}(iv). Moreover, by Lemmas, \ref{lem:Gxi}(v), \ref{lem:fidio} and \ref{lem:consistCOV}(vii), and \citet[Theorem 1]{MK04}, which is Weyl's inequality, for all $j=1,\ldots, r$,
\al{
n^{-1}\vert \wh{\mu}_j^x-{\mu}_j^\chi\vert \le n^{-1}\l\Vert T^{-1}\sum_{t=1}^T \mbf x_t\mbf x_t^\prime-\bm\Gamma^x\r\Vert+n^{-1}\Vert\bm\Gamma^\xi\Vert =O_p(\max(n^{-1},T^{-1/2})),\label{eq:panza2}
}
which, jointly with \eqref{eq:panza1}, implies
\al{
\det(n^{-1} \wh{\mbf M}_n^x) =\prod_{j=1}^r n^{-1}\wh{\mu}_j^x \ge \{ n^{-1}\wh{\mu}_r^x\}^r\ge  \{ n^{-1}{\mu}_r^\chi-n^{-1}\vert \wh{\mu}_r^x-{\mu}_r^\chi\vert \}^r>0,\nn
}
thus, 
\al{
n \Vert(\wh{\mbf M}_n^x)^{-1}\Vert=O_p(1).\label{eq:panza3}
}
From \eqref{eq:panza1}, \eqref{eq:panza2}, and \eqref{eq:panza3}
\al{
n\l\Vert(\wh{\mbf M}_n^x)^{-1}-(\mbf M_n^\chi)^{-1}\r\Vert
&\le n \Vert(\mbf M_n^\chi)^{-1}\Vert \,
n^{-1}\l\Vert\wh{\mbf M}_n^x-\mbf M_n^\chi\r\Vert
\,
n \Vert(\wh{\mbf M}_n^x)^{-1}\Vert=O_p(\max(n^{-1},T^{-1/2})).\label{eq:BUM}
}
Furthermore, 
\al{
n^{-1}\Vert \wh{\bm\Lambda}_n^{(0)\prime}(\bm\Lambda_n-\wh{\bm\Lambda}_n^{(0)}) \mbf F_t\Vert \le&\,
n^{-1}\Vert {\bm\Lambda}_n^{\prime}(\bm\Lambda_n-\wh{\bm\Lambda}_n^{(0)})\Vert\, \Vert \mbf F_t\Vert +n^{-1}\Vert\bm\Lambda_n-\wh{\bm\Lambda}_n^{(0)}\Vert^2\, \Vert \mbf F_t\Vert \nn\\
\le&\, 
n^{-1}\Vert {\bm\Lambda}_n^{\prime}(\bm\Lambda_n-\wh{\bm\Lambda}_n^{\text{\tiny OLS}})\Vert\, \Vert \mbf F_t\Vert +
n^{-1/2}\Vert {\bm\Lambda}_n\Vert\, n^{-1/2}\Vert\wh{\bm\Lambda}_n^{(0)}-\bm\Lambda_n^{\text{\tiny OLS}}\Vert\, \Vert \mbf F_t\Vert \nn\\
&+n^{-1}\Vert\bm\Lambda_n-\wh{\bm\Lambda}_n^{(0)}\Vert^2\, \Vert \mbf F_t\Vert \nn\\
\le&\,
n^{-1}\l\Vert 
T^{-1}\sum_{t=1}^T {\bm\Lambda}_n^\prime\bm\xi_{nt}\mbf F_t^\prime\r\Vert\,
\l\Vert \l(T^{-1}\sum_{t=1}^T \mbf F_t\mbf F_t^\prime\r)^{-1}
\r\Vert
\Vert \mbf F_t\Vert\nn\\
 &+n^{-1/2}\Vert {\bm\Lambda}_n\Vert\, n^{-1/2}\Vert\wh{\bm\Lambda}_n^{(0)}-\bm\Lambda_n^{\text{\tiny OLS}}\Vert\, \Vert \mbf F_t\Vert +n^{-1}\Vert\bm\Lambda_n-\wh{\bm\Lambda}_n^{(0)}\Vert^2\, \Vert \mbf F_t\Vert \nn\\
=&\, O_p(n^{-1/2}T^{-1/2})+O_p(\max(n^{-1},n^{-1/2}T^{-1/2}))+O_p(\max(n^{-2},T^{-1}))\nn\\
=&\, O_p(\max(n^{-1},n^{-1/2}T^{-1/2},T^{-1})),\label{eq:BUM2}
}
by \citet[Corollary 1]{MBPCAQML}, 
Lemmas \ref{lem:lambdasqrtn}, \ref{lem:COVFF0}(i), \ref{lem:frida}, and \ref{lem:est0LOAD}(b), 
and since $\Vert \mbf F_t\Vert = O_p(1)$ because $\E[F_{jt}^2]=1$, $j=1,\ldots, r$, by Assumption \ref{ass:ident}(b). 

Last,
\al{
n^{-1}\Vert(\wh{\bm\Lambda}_n^{(0)}-\bm\Lambda_n)^\prime\bm\xi_{nt}\Vert\le&\, 
n^{-1}\Vert(\bm\Lambda_n^{\text{\tiny OLS}}-{\bm\Lambda}_n)^\prime\bm\xi_{nt}\Vert
+n^{-1/2}\Vert\wh{\bm\Lambda}_n^{(0)}-\bm\Lambda_n^{\text{\tiny OLS}}\Vert\, n^{-1/2}\Vert\bm\xi_{nt}\Vert\nn\\
\le&\, n^{-1}\l\Vert T^{-1}\sum_{t=1}^T \mbf F_t\bm\xi_{nt}^\prime\bm\xi_{nt} \r\Vert\,
\l\Vert \l(T^{-1}\sum_{t=1}^T \mbf F_t\mbf F_t^\prime\r)^{-1}
\r\Vert\nn\\
&+n^{-1/2}\Vert\wh{\bm\Lambda}_n^{(0)}-\bm\Lambda_n^{\text{\tiny OLS}}\Vert\, n^{-1/2}\Vert\bm\xi_{nt}\Vert\nn\\
=&\, O_p(n^{-1/2}T^{-1/2})+O_p(\max(n^{-1},n^{-1/2}T^{-1/2}))\nn\\
=&\, O_p(\max(n^{-1},n^{-1/2}T^{-1/2})),
\label{eq:BUM3}
}
by \citet[Corollary 1]{MBPCAQML}, Lemmas \ref{lem:COVFF0}(vi) (setting $\bm\Sigma_n^\xi=\mbf I_n$ and using $\bm\xi_{nt}$ in place of $\bm{\mathcal E}_{nT}$ therein), 
\ref{lem:frida}, and \ref{lem:est0LOAD}(b), 
and since $n^{-1/2}\Vert\bm\xi_{nt}\Vert=O_p(1)$ because $\sum_{i=1}^n\sigma_i^2\le n C_\xi$ by Assumption \ref{ass:idio}(a).

By substituting \eqref{eq:BUM}, \eqref{eq:BUM2},  and \eqref{eq:BUM3} into \eqref{eq:panzona}, it follows that if $\sqrt n/T\to 0$, as $n,T\to\infty$, 
\al{
\sqrt n (\wt{\mbf F}_t-\mbf F_t) = n(\mbf M_n^\chi)^{-1}(n^{-1/2}\bm\Lambda_n^\prime\bm\xi_{nt})
+o_p(1).\label{eq:sviluppoFPC}
}
and since $\lim_{n\to\infty} n^{-1}\mbf M_n^\chi= \bm\Sigma_\Lambda$ by Assumption \ref{ass:ident}(b) which is positive definite, by Assumption \ref{ass:idio}(e) (when setting $\sigma_i^2=1$ therein) and Slutsky's Theorem we complete the proof. $\Box$

\section{Derivation of the Kalman filter MSE}\label{app:foche}

\begin{lem}\label{lem:control} 
Under Assumption \ref{ass:common}, the DFM \eqref{eq:SDFM1R}-\eqref{eq:SDFM2R} is both stabilizable and detectable, for all $n\in\mathbb N$.
\end{lem}

\noindent\textsc{Proof.} We use the definitions in \citet[Appendix C, p. 341-342]{AM79}. The DFM \eqref{eq:SDFM1R}-\eqref{eq:SDFM2R} is a linear systems with $r$ states. A linear system is stabilizable if its unstable states are controllable and all uncontrollable states are stable, and it is detectable if its unstable states are observable and all unobservable states are stable.

First, by factorizing $\bm \Gamma^{v}=\mbf H\mbf H^\prime$ for some $\mbf H$ having full-column rank, we see that $\text{rk}[\mbf H\; (\mbf A\mbf H)\cdots (\mbf A^{(r-1)}\mbf H)]=r$, thus the linear system is controllable, by Assumption \ref{ass:common}(e). Moreover, there are no unstable states, since because of  Assumption \ref{ass:common}(d) all eigenvalues of ${\mbf A}$ are smaller than one in absolute value. This implies that the model is stabilizable. 

Second, since by Assumption \ref{ass:common}(a) for any given $n\in\mathbb N$ there exists at least one $i=1,\ldots,n$ such that $\Vert{\bm\lambda}_i\Vert \ge m_\lambda$, then $\mbox{rk}({\bm \Lambda}_n)\ge 1$, while, $\mbox{rk}({\bm \Lambda}_n)=r$ only for $n>N$, then, for a given $n$, there might be $(r-1)$ unobservable states, however, as already noticed, they are all stable by Assumption \ref{ass:common}(d). Thus the model is detectable. This completes the proof.  $\Box$

\begin{lem}\label{lem:steady1} 
Under Assumptions \ref{ass:common} and \ref{ass:idio}, the matrix $\mbf P_{t|t-1}$ has a steady-state denoted as ${\mbf P}=\lim_{t\to\infty} \mbf P_{t|t-1}$.
\end{lem}

\noindent\textsc{Proof.} 
First, given that with our initialization $\mbf P_{0|0}=\bm\Gamma^F$, then it is positive definite by Assumption \ref{ass:common}(b), therefore also $\mbf P_{1|0}$ is positive definite (see also \eqref{eq:pred2}). Second,  as proved in Lemma \ref{lem:control}, the linear system defining the DFM \eqref{eq:SDFM1R}-\eqref{eq:SDFM2R} is stabilizable and detectable.
Therefore,  because of \citet[Theorem 4.1]{CGS84}, as $t\to\infty$, ${\mbf P}_{t|t-1}$ converges to a steady-state $\mbf P$ exponentially fast (see Lemma \ref{lem:steady_rate} below for the rate of convergence), which is a solution of the algebraic Riccati equation (ARE) derived from \eqref{eq:riccati}
\al{
\mbf P=&\, \mbf A\mbf P\mbf A^\prime-\mbf A\mbf P\bm\Lambda_n^\prime(\bm\Lambda_n\mbf P\bm\Lambda_n^\prime+ \bm\Sigma_n^\xi)^{-1}\bm\Lambda_n\mbf P\mbf A^\prime+\bm\Gamma^v.\nn
}
Moreover, since Lemmas \ref{lem:cazzarola}(i)  and \ref{lem:cazzarola}(ii) hold for all $T\in\mathbb N$, we have 
\beq\label{eq:PrincipeGatto}
\Vert \mbf P\Vert =O(1),\qquad \Vert \mbf P^{-1}\Vert =O(1).
\eeq 
This completes the proof. $\Box$

\begin{lem}\label{lem:steady1pi} 
Under Assumptions \ref{ass:common} and \ref{ass:idio}, the matrix $\bm\Pi_{t|t-1}$ has a steady-state denoted as ${\bm \Pi}=\lim_{t\to\infty} \bm\Pi_{t|t-1}$.
\end{lem}

\noindent\textsc{Proof.} The existence of ${\bm\Pi}$ follows from the same arguments used in Lemma \ref{lem:steady1}. Moreover, from \citet[Section 2.2]{harvey2009computing} we have that ${\bm\Pi}$ must satisfy
\al{
\bm \Pi=&\,\mbf A\bm\Pi\mbf A^\prime+\mbf A\mbf P\bm\Lambda_n^\prime(\bm\Lambda_n\mbf P\bm\Lambda_n^\prime+ \bm\Sigma_n^\xi)^{-1}
(\bm\Lambda_n\bm\Pi\bm\Lambda_n^\prime+ \bm\Gamma_n^\xi)
(\bm\Lambda_n\mbf P\bm\Lambda_n^\prime+ \bm\Sigma_n^\xi)^{-1}\bm\Lambda_n\mbf P\mbf A^\prime\nn\\
&-\mbf A\mbf P\bm\Lambda_n^\prime(\bm\Lambda_n\mbf P\bm\Lambda_n^\prime+ \bm\Sigma_n^\xi)^{-1}\bm\Lambda_n\bm\Pi \mbf A^\prime
-\mbf A\bm\Pi \bm\Lambda_n^\prime(\bm\Lambda_n\mbf P\bm\Lambda_n^\prime+ \bm\Sigma_n^\xi)^{-1}\bm\Lambda_n\mbf P\mbf A^\prime+\bm\Gamma^v.\nn
}
Moreover, by the same arguments in Lemmas \ref{lem:cazzarola}(i) and \ref{lem:cazzarola}(ii) it holds that $\max_{t=1,\ldots, T}\Vert \bm \Pi_{t|t-1}\Vert=O(1)$ and also $\max_{t=1,\ldots, T}\Vert  (\bm\Pi_{t|t-1})^{-1}\Vert=O(1)$, for all $T\in\mathbb N$, thus 
\beq
\Vert \bm \Pi\Vert=O(1),\qquad \Vert \bm \Pi^{-1}\Vert=O(1).\label{eq:PrincipeCane}
\eeq 
This completes the proof. $\Box$
\begin{lem}\label{lem:steady_rate} 
Under Assumptions \ref{ass:common}, \ref{ass:idio}, and \ref{ass:ident}, if $\Vert\mbf P_{0|0}\Vert=O(n^\gamma)$ for some $\gamma>0$, then 
$n\max_{t=\bar t,\ldots, T}\Vert\mbf P_{t|t-1}-{\mbf P} \Vert=o(1)$ and
$n\max_{t=\bar t,\ldots, T}\Vert\bm \Pi_{t|t-1}-{\bm \Pi} \Vert=o(1)$, where $\bar t = \lceil 2 +\gamma/2\rceil $.
\end{lem}

\noindent\textsc{Proof.} 
Let $\bm\Psi_{1,1}=\mbf I_r$ and, for $t=2,\ldots,T$, 
\beq\label{Psit1}
\bm\Psi_{t,1} =\prod_{s=1}^{t-1}[{\mbf A}-{\mbf A}{\mbf P}_{s|s-1}{\bm \Lambda}_n^\prime ({\bm \Lambda}_n{\mbf P}_{s|s-1}{\bm \Lambda}_n^\prime+\bm\Sigma_n^\xi)^{-1}{\bm \Lambda}_n].
\eeq
Then, from \citet[Chapter 4.4, pp. 76-81]{AM79}, we have, for $t=1,\ldots,T$,
\al{
{\mbf P}_{t|t-1} -{\mbf P}  =&\,  \{{\mbf A}-{\mbf A}{\mbf P}{\bm \Lambda}_n^\prime ({\bm \Lambda}_n{\mbf P}{\bm \Lambda}_n^\prime+\bm\Sigma_n^\xi)^{-1}{\bm \Lambda}_n\}^{t-1} ({\mbf P}_{1|0} -{\mbf P}) \bm\Psi_{t,1}\nn
\\
=&\,{\mbf A}^{t-1}\{\mbf I_r-{\mbf P}{\bm \Lambda}_n^\prime ({\bm \Lambda}_n{\mbf P}{\bm \Lambda}_n^\prime+\bm\Sigma_n^\xi)^{-1}{\bm \Lambda}_n\}^{t-1} ({\mbf P}_{1|0} -{\mbf P}) \bm\Psi_{t,1}\nn\\
=&\, {\mbf A}^{t-1}\{
\mbf I_r- (\bm\Lambda_n^\prime(\bm\Sigma_n^\xi)^{-1}{\bm \Lambda}_n+\mbf P^{-1})^{-1}{\bm \Lambda}_n^\prime(\bm\Sigma_n^\xi)^{-1}{\bm \Lambda}_n
\}^{t-1}({\mbf P}_{1|0} -{\mbf P}) \bm\Psi_{t,1},\label{eq:convergence_PP0}
}
where we used Lemma \ref{lem:wood}. 

Now, for $t=2,\ldots, T$, from \eqref{Psit1} and using again Lemma \ref{lem:wood}, there exists a $\bar n$ such that for all $n\ge \bar n$
\begin{align}
\label{Psit1_cont}
\Vert \bm\Psi_{t,1}\Vert&\le \Vert{\mbf A}\Vert^{t-1} \prod_{s=1}^{t-1}\Vert \mbf I_r- (\bm\Lambda_n^\prime(\bm\Sigma_n^\xi)^{-1}{\bm \Lambda}_n+\mbf P_{s|s-1}^{-1})^{-1}{\bm \Lambda}_n^\prime(\bm\Sigma_n^\xi)^{-1}{\bm \Lambda}_n\Vert \le M_0 n^{-(t-1)},
\end{align}
for some finite positive real $M_0$ independent of $n$ and $t$, 
because of Lemmas \ref{lem:LSL}(i) and \ref{lem:cazzarola}(ii), and Assumption \ref{ass:common}(d). Moreover, 
\beq\label{eq:P00Steady}
\Vert {\mbf P}_{1|0} -{\mbf P}\Vert = \Vert \mbf A{\mbf P}_{0|0} -{\mbf P}\Vert\le \Vert \mbf A\Vert \, \Vert {\mbf P}_{0|0}\Vert +\Vert \mbf A\Vert \, \Vert {\mbf P}\Vert\le 2 \Vert \mbf A\Vert \, \Vert {\mbf P}_{0|0}\Vert,
\eeq
since because of Lemma \ref{lem:detP}, $\Vert \mbf P\Vert \le \Vert\mbf P_{t|t-1}\Vert\le \Vert \mbf P_{1|0}\Vert\le \mbf P_{0|0}\Vert$.

Therefore,  for $t=2,\ldots, T$, from \eqref{eq:convergence_PP0}, \eqref{Psit1_cont}, and \eqref{eq:P00Steady},  there exists a $\bar n$ such that for all $n\ge \bar n$
\al{
\Vert{\mbf P}_{t|t-1} -{\mbf P} \Vert \le&\, 2 \Vert {\mbf A}\Vert^{t}\Vert \mbf I_r- ({\bm \Lambda}_n^{\prime}(\bm\Sigma_n^\xi)^{-1}{\bm \Lambda}_n+{\mbf P}^{-1})^{-1}
{\bm \Lambda}_n^{\prime}(\bm\Sigma_n^\xi)^{-1}{\bm \Lambda}_n\}\Vert^{t-1}\, \Vert {\mbf P}_{0|0} \Vert\, \Vert \bm\Psi_{t,1}\Vert\nn\\
\le&\, M_1 n^{-2(t-1)} \Vert {\mbf P}_{0|0} \Vert,\label{eq:convergence_PP}
}
for some finite positive real $M_1$ independent of $n$ and $t$, because of Lemmas \ref{lem:LSL}(i) and \ref{lem:cazzarola}(ii), and Assumption \ref{ass:common}(d). 

Now, define $\bar t$ the first point in time such that the rhs of \eqref{eq:convergence_PP} is $o(n^{-1})$, i.e., such that
\beq\nn
n^{-2(\bar t-1)} n^{1+\epsilon}  \le M_2 \Vert{\mbf P}_{0|0}\Vert^{-1},
\eeq
for any $\epsilon>0$ and some finite positive real $M_2$ independent of $n$ and $t$, or equivalently, by letting $K=\log M_2$,  $\bar t$ is such that:
\beq\label{bartt2}
\bar t \ge (3+\epsilon)/2+{\log \Vert{\mbf P}_{0|0}\Vert}/(2\log n)-K/(2\log n).
\eeq
Clearly \eqref{bartt2} is always satisfied if we set $\bar t = \lceil 2+{\log \Vert{\mbf P}_{0|0}\Vert}/(2\log n)\rceil$. Letting now $\Vert{\mbf P}_{0|0}\Vert=O(n^\gamma)$ for some $\gamma>0$, then  we have $\bar t=\lceil 2+\gamma/2\rceil$ satisfies \eqref{bartt2} and therefore we have at least $\max_{t=\bar t,\ldots, T}\Vert{\mbf P}_{t|t-1} -{\mbf P} \Vert= O(n^{-1})$.
Notice that if $\gamma=0$, i.e., $\Vert{\mbf P}_{0|0}\Vert$ is finite, then $\bar t=2$ and in this case from \eqref{eq:convergence_PP} we have an even tighter bound as we get $\max_{t=2,\ldots, T}\Vert{\mbf P}_{t|t-1} -{\mbf P} \Vert= O(n^{-2})$. 

The proof for $\bm\Pi_{t|t-1}$ can be done analogously but using the recursions in \citet{harvey2009computing}. This completes the proof. $\Box$\\

\begin{lem}\label{lem:foche} 
Under Assumptions \ref{ass:common}, \ref{ass:idio}, and \ref{ass:ident}, 
\begin{compactenum}[(i)]
\item $n\max_{t=\bar t,\ldots, T}\Vert\bm{\mathcal P}_t\Vert = O(1)$;
\item $n\max_{t=1,\ldots, T}\Vert\bm{\Pi}_{t|t}\Vert = O(1)$;
\item $n \max_{\bar t,\ldots, T} \Vert \bm \Pi_{t|t}-\bm{\mathcal P}_t\Vert = o(1)$,
\item $\max_{t=\bar t,\ldots, T}\Vert n\bm{\mathcal P}_t-\bm{\mathcal W}_t\Vert = o(1)$;
\end{compactenum}
where 
$ \bm \Pi_{t|t}$ is defined in \eqref{eq:dellefoche}, $ \bm {\mathcal P}_{t}$ is obtained from $ \bm \Pi_{t|t}$ when replacing $\mbf P_{t|t-1}$ and $\bm \Pi_{t|t-1}$ with their steady states defined in Lemmas \ref{lem:steady1} and \ref{lem:steady1pi}, respectively, 
and
$\bm{\mathcal W}_t$ is defined in Proposition \ref{prop:factors}.
\end{lem}

\noindent\textsc{Proof.} We have
\al{
\bm{\mathcal P}_t =&\, \bm \Pi +  \mbf P \bm\Lambda_n^\prime ( \bm\Lambda_n\mbf P\bm\Lambda_n^\prime+\bm\Sigma_n^\xi)^{-1} \bm\Lambda_n\bm\Pi \bm\Lambda_n^\prime ( \bm\Lambda_n\mbf P\bm\Lambda_n^\prime+\bm\Sigma_n^\xi)^{-1}\bm\Lambda_n \mbf P\nn\\
&- \mbf P \bm\Lambda_n^\prime ( \bm\Lambda_n\mbf P\bm\Lambda_n^\prime+\bm\Sigma_n^\xi)^{-1}\bm\Lambda_n \bm \Pi
- \bm \Pi \bm\Lambda_n^\prime ( \bm\Lambda_n\mbf P\bm\Lambda_n^\prime+\bm\Sigma_n^\xi)^{-1}\bm\Lambda_n\mbf P\nn\\
&+  \mbf P \bm\Lambda_n^\prime ( \bm\Lambda_n\mbf P\bm\Lambda_n^\prime+\bm\Sigma_n^\xi)^{-1} \bm\Gamma_n^\xi ( \bm\Lambda_n\mbf P\bm\Lambda_n^\prime+\bm\Sigma_n^\xi)^{-1}\bm\Lambda_n \mbf P.\label{eq:gattopardo0}
}
Hereafter, for simplicity of notation let $\bm H=({\bm\Lambda}_n^{\prime}({\bm\Sigma}_n^{\xi})^{-1}{\bm\Lambda}_n)^{-1}$.
Using twice Lemma \ref{lem:taylorinv} we have
\al{
\mbf P(\bm H+\mbf P)^{-1} 
&= \mbf P\l\{\mbf P^{-1}- (\bm H+\mbf P)^{-1}\bm H \mbf P^{-1} \r\}\nn\\
&=\mbf I_r- \mbf P\l\{\mbf P^{-1}-(\bm H+\mbf P)^{-1}\bm H \mbf P^{-1} 
\r\}\bm H \mbf P^{-1} \nn\\
&= \mbf I_r-\bm H \mbf P^{-1}+\mbf P(\bm H+\mbf P)^{-1}\bm H \mbf P^{-1} \bm H \mbf P^{-1}\nn\\
&= \mbf I_r-\bm H \mbf P^{-1} +\mbf C,\; \text{ say.} \label{eq:taylor2nd}
}
Then, by Lemma \ref{lem:pioveadir8} and \eqref{eq:taylor2nd}
\al{
\bm \Pi- \mbf P \bm\Lambda_n^\prime ( \bm\Lambda_n\mbf P\bm\Lambda_n^\prime+\bm\Sigma_n^\xi)^{-1}\bm\Lambda_n \bm \Pi=&\,
\l\{\mbf I_r- \mbf P(\bm H+\mbf P)^{-1}\r\}\bm \Pi
= \bm H\mbf P^{-1} \bm \Pi -\mbf C \bm \Pi. \label{eq:gattopardo1}
}
Similarly, again by Lemma \ref{lem:pioveadir8} and \eqref{eq:taylor2nd}
\al{
 \mbf P  \bm\Lambda_n^\prime& ( \bm\Lambda_n\mbf P\bm\Lambda_n^\prime+\bm\Sigma_n^\xi)^{-1} \bm\Lambda_n\bm\Pi \bm\Lambda_n^\prime ( \bm\Lambda_n\mbf P\bm\Lambda_n^\prime+\bm\Sigma_n^\xi)^{-1}\bm\Lambda_n \mbf P- \bm \Pi \bm\Lambda_n^\prime ( \bm\Lambda_n\mbf P\bm\Lambda_n^\prime+\bm\Sigma_n^\xi)^{-1}\bm\Lambda_n\mbf P\nn\\
 =&\,\mbf P(\bm H+\mbf P)^{-1} \bm \Pi (\bm H+\mbf P)^{-1} \mbf P-\bm \Pi \mbf P (\bm H+\mbf P)^{-1}\nn\\
  =&\,
  \l\{
   \mbf I_r-\bm H \mbf P^{-1} +\mbf C
  \r\}
  \bm \Pi 
  \l\{
   \mbf I_r- \mbf P^{-1} \bm H +\mbf C
  \r\}
    -\bm \Pi  \l\{
   \mbf I_r- \mbf P^{-1} \bm H +\mbf C
  \r\}\nn\\
  =&\, \bm \Pi -\bm \Pi \mbf P^{-1}\bm H-\bm H\mbf P^{-1} \bm \Pi +\bm H\mbf P^{-1}\bm\Pi\mbf P^{-1}\bm H-\bm \Pi +\bm \Pi \mbf P^{-1}\bm H\nn\\
  &+\mbf C\bm \Pi+\bm \Pi\mbf C -\bm H\mbf P^{-1}\bm \Pi\mbf C-\mbf C\bm\Pi \mbf P^{-1}\bm H+\mbf C\bm\Pi\mbf C-\bm \Pi\mbf C\nn\\
  =&\,  -\bm H\mbf P^{-1} \bm \Pi +\bm H\mbf P^{-1}\bm\Pi\mbf P^{-1}\bm H +\mbf C\bm \Pi -\bm H\mbf P^{-1}\bm \Pi\mbf C-\mbf C\bm\Pi \mbf P^{-1}\bm H+\mbf C\bm\Pi\mbf C.
  \label{eq:gattopardo2}
}
By substituting \eqref{eq:gattopardo1} and \eqref{eq:gattopardo2} into \eqref{eq:gattopardo0}:
\al{
\bm{\mathcal P}_t =&\,  \mbf P \bm\Lambda_n^\prime ( \bm\Lambda_n\mbf P\bm\Lambda_n^\prime+\bm\Sigma_n^\xi)^{-1} \bm\Gamma_n^\xi ( \bm\Lambda_n\mbf P\bm\Lambda_n^\prime+\bm\Sigma_n^\xi)^{-1}\bm\Lambda_n \mbf P\nn\\
&+\bm H\mbf P^{-1} \bm \Pi -\mbf C \bm \Pi -\bm H\mbf P^{-1} \bm \Pi +\bm H\mbf P^{-1}\bm\Pi\mbf P^{-1}\bm H +\mbf C\bm \Pi -\bm H\mbf P^{-1}\mbf C\bm \Pi-\mbf C\bm\Pi \mbf P^{-1}\bm H+\mbf C\bm\Pi\mbf C\nn\\
=&\,  \mbf P \bm\Lambda_n^\prime ( \bm\Lambda_n\mbf P\bm\Lambda_n^\prime+\bm\Sigma_n^\xi)^{-1} \bm\Gamma_n^\xi ( \bm\Lambda_n\mbf P\bm\Lambda_n^\prime+\bm\Sigma_n^\xi)^{-1}\bm\Lambda_n \mbf P\nn\\
&-\mbf C \bm \Pi +\bm H\mbf P^{-1}\bm\Pi\mbf P^{-1}\bm H +\mbf C\bm \Pi -\bm H\mbf P^{-1}\bm \Pi\mbf C-\mbf C\bm\Pi \mbf P^{-1}\bm H+\mbf C\bm\Pi\mbf C\nn\\
=&\,  ( \bm\Lambda_n^\prime(\bm\Sigma_n^\xi)^{-1}\bm\Lambda_n+\mbf P^{-1})^{-1}
\bm\Lambda_n^\prime(\bm\Sigma_n^\xi)^{-1}
 \bm\Gamma_n^\xi 
 (\bm\Sigma_n^\xi)^{-1}\bm\Lambda_n( \bm\Lambda_n^\prime(\bm\Sigma_n^\xi)^{-1}\bm\Lambda_n+\mbf P^{-1})^{-1}\nn\\
& +\bm H\mbf P^{-1}\bm\Pi\mbf P^{-1}\bm H  -\bm H\mbf P^{-1}\bm \Pi\mbf C-\mbf C\bm\Pi \mbf P^{-1}\bm H+\mbf C\bm\Pi\mbf C.\label{eq:gattopardo3}
}
where in the last step we used Lemma \ref{lem:wood}. 
Moreover,
\al{
\Vert \bm H\Vert =&\,\Vert ({\bm\Lambda}_n^{\prime}({\bm\Sigma}_n^{\xi})^{-1}{\bm\Lambda}_n)^{-1}\Vert= O(n^{-1}), \label{eq:gattopardo4bis}\\
\Vert\mbf C\Vert =&\,\Vert \mbf P(({\bm\Lambda}_n^{\prime}({\bm\Sigma}_n^{\xi})^{-1}{\bm\Lambda}_n)^{-1}+\mbf P)^{-1}({\bm\Lambda}_n^{\prime}({\bm\Sigma}_n^{\xi})^{-1}{\bm\Lambda}_n)^{-1} \mbf P^{-1}({\bm\Lambda}_n^{\prime}({\bm\Sigma}_n^{\xi})^{-1}{\bm\Lambda}_n)^{-1} \mbf P^{-1}\Vert\nn\\
&\le \Vert \mbf P\Vert \, \Vert (\mbf P)^{-1}\Vert^2\Vert (({\bm\Lambda}_n^{\prime}({\bm\Sigma}_n^{\xi})^{-1}{\bm\Lambda}_n)^{-1}+\mbf P)^{-1}\Vert\, \Vert ({\bm\Lambda}_n^{\prime}({\bm\Sigma}_n^{\xi})^{-1}{\bm\Lambda}_n)^{-1}\Vert^2= O(n^{-2}),.\label{eq:gattopardo4}
}
because of Lemma \ref{lem:LSL2}(iii), \eqref{eq:PrincipeGatto} in the proof of Lemma \ref{lem:steady1},
and  since, by \citet[Theorem 1]{MK04} which is Weyl's inequality,
 \al{
 \Vert &(({\bm\Lambda}_n^{\prime}({\bm\Sigma}_n^{\xi})^{-1}{\bm\Lambda}_n)^{-1}+ {\mbf P})^{-1}\Vert =
 \l\{\nu^{(r)}( ({\bm\Lambda}_n^{\prime}({\bm\Sigma}_n^{\xi})^{-1}{\bm\Lambda}_n)^{-1}+ {\mbf P})  \r\}^{-1}\nn\\
 \le &\,  \l\{\nu^{(r)}( ({\bm\Lambda}_n^\prime({\bm\Sigma}_n^{\xi})^{-1}{\bm\Lambda}_n)^{-1})+ \nu^{(r)}( {\mbf P})  \r\}^{-1}\nn\\
 = &\,  \l\{
 \l[\nu^{(1)} ({\bm\Lambda}_n^{\prime}({\bm\Sigma}_n^{\xi})^{-1}{\bm\Lambda}_n)\r]^{-1}
 + \nu^{(r)}( {\mbf P})  \r\}^{-1}\nn\\
  = &\,  \l\{
 \l[\nu^{(1)} ({\bm\Lambda}_n^{\prime}({\bm\Sigma}_n^{\xi})^{-1}{\bm\Lambda}_n)\nu^{(r)}( {\mbf P})\r]^{-1}
 + 1  \r\}^{-1} \l\{\nu^{(r)}( {\mbf P})\r\}^{-1}\nn\\
 =&\, \l\{1-\l[\nu^{(1)} ({\bm\Lambda}_n^{\prime}({\bm\Sigma}_n^{\xi})^{-1}{\bm\Lambda}_n)\nu^{(r)}( {\mbf P})\r]^{-1}\r\}\l\{\nu^{(r)}( {\mbf P})\r\}^{-1}+O_p(n^{-2})\nn\\
 =&\, O_p(1),\nn
 }
 again by Lemma \ref{lem:LSL2}(iii) and \eqref{eq:PrincipeGatto} in the proof of Lemma \ref{lem:steady1}. Therefore, from \eqref{eq:gattopardo3}, \eqref{eq:gattopardo4bis}, and \eqref{eq:gattopardo4}
\al{
n\Vert \bm{\mathcal P}_t\Vert \le&\, n
\Vert ( \bm\Lambda_n^\prime(\bm\Sigma_n^\xi)^{-1}\bm\Lambda_n+\mbf P^{-1})^{-1}\Vert^2
\Vert \bm\Lambda_n^\prime(\bm\Sigma_n^\xi)^{-1}\Vert^2
\Vert  \bm\Gamma_n^\xi  \Vert \nn\\
&+n \Vert \bm H\Vert^2\Vert(\mbf P)^{-1}\Vert^2 \Vert\bm\Pi\Vert+2n \Vert\bm H\Vert \, \Vert(\mbf P)^{-1}\Vert \, \Vert \bm\Pi\Vert\, \Vert\mbf C\Vert
+n \Vert\mbf C\Vert^2\Vert\bm\Pi\Vert\nn\\
=&\, O(1)+ O(n^{-1})+ O(n^{-2}) + O(n^{-3}),\label{eq:calamaretto}
}
where we used also \eqref{eq:PrincipeGatto} and \eqref{eq:PrincipeCane} in the proofs of Lemmas \ref{lem:steady1} and \ref{lem:steady1pi}, respectively. Since $\bm{\mathcal P}_t$ does not depend on $t$ we prove part (i). Part (ii) is proved analogously by substituting in part (i) $\mbf P$ and $\bm \Pi$ with $\mbf P_{t|t-1}$ and $\bm \Pi_{t|t-1}$ and using Lemma \ref{lem:cazzarola}(i) and \ref{lem:cazzarola}(ii) instead of \eqref{eq:PrincipeGatto}  and \eqref{eq:PrincipeCane}. 

Part (iii) follows directly from Lemmas \ref{lem:steady1}, \ref{lem:steady1pi}, and \ref{lem:steady_rate}, and parts (i) and (ii).

Turning to part (iv), consider the first term on the rhs of \eqref{eq:gattopardo3}, using Lemma \ref{lem:wood} we have
\al{
\mbf P \bm\Lambda_n^\prime& ( \bm\Lambda_n\mbf P\bm\Lambda_n^\prime+\bm\Sigma_n^\xi)^{-1} \bm\Gamma_n^\xi ( \bm\Lambda_n\mbf P\bm\Lambda_n^\prime+\bm\Sigma_n^\xi)^{-1}\bm\Lambda_n \mbf P= (\bm H^{-1}+\mbf P^{-1} )^{-1} \bm\Lambda_n^\prime(\bm\Sigma_n^\xi)^{-1}
 \bm\Gamma_n^\xi 
 (\bm\Sigma_n^\xi)^{-1}\bm\Lambda_n (\bm H^{-1}+\mbf P^{-1} )^{-1} \nn\\
 =&\,
 \l\{
 (\bm H^{-1}+\mbf P^{-1} )^{-1} \bm\Lambda_n^\prime(\bm\Sigma_n^\xi)^{-1}-
 \bm H \bm\Lambda_n^\prime(\bm\Sigma_n^\xi)^{-1}
 +\bm H \bm\Lambda_n^\prime(\bm\Sigma_n^\xi)^{-1}
 \r\}
 \bm\Gamma_n^\xi\nn\\
&\cdot \l\{
 (\bm\Sigma_n^\xi)^{-1}\bm\Lambda_n(\bm H^{-1}+\mbf P^{-1} )^{-1}
 -
  (\bm\Sigma_n^\xi)^{-1}\bm\Lambda_n\bm H
  +  (\bm\Sigma_n^\xi)^{-1}\bm\Lambda_n\bm H
 \r\}\nn\\
 =&\,  
 \bm H \bm\Lambda_n^\prime(\bm\Sigma_n^\xi)^{-1} \bm\Gamma_n^\xi(\bm\Sigma_n^\xi)^{-1}\bm\Lambda_n\bm H
 \nn\\
  &+
  \l\{
 (\bm H^{-1}+\mbf P^{-1} )^{-1} \bm\Lambda_n^\prime(\bm\Sigma_n^\xi)^{-1}-
 \bm H \bm\Lambda_n^\prime(\bm\Sigma_n^\xi)^{-1}\r\} \bm\Gamma_n^\xi(\bm\Sigma_n^\xi)^{-1}\bm\Lambda_n\bm H\nn\\
 &+\bm H \bm\Lambda_n^\prime(\bm\Sigma_n^\xi)^{-1} \bm\Gamma_n^\xi \l\{
 (\bm\Sigma_n^\xi)^{-1}\bm\Lambda_n(\bm H^{-1}+\mbf P^{-1} )^{-1}
 -
  (\bm\Sigma_n^\xi)^{-1}\bm\Lambda_n\bm H\r\}\nn\\
 &+ \l\{
 (\bm H^{-1}+\mbf P^{-1} )^{-1} \bm\Lambda_n^\prime(\bm\Sigma_n^\xi)^{-1}-
 \bm H \bm\Lambda_n^\prime(\bm\Sigma_n^\xi)^{-1}\r\} \bm\Gamma_n^\xi
 \l\{
 (\bm\Sigma_n^\xi)^{-1}\bm\Lambda_n(\bm H^{-1}+\mbf P^{-1} )^{-1}
 -
  (\bm\Sigma_n^\xi)^{-1}\bm\Lambda_n\bm H\r\}\nn\\
  =&\, \bm H \bm\Lambda_n^\prime(\bm\Sigma_n^\xi)^{-1} \bm\Gamma_n^\xi(\bm\Sigma_n^\xi)^{-1}\bm\Lambda_n\bm H + \bm{ \mathcal A} +\bm{ \mathcal A}^\prime +\bm{ \mathcal B}, \; \text{say.}\label{eq:AABedge}
}
Then, 
\al{
n \Vert\bm{\mathcal A}\Vert &\le 
n \Vert  ({\bm\Lambda}_n^{\prime}({\bm\Sigma}_n^{\xi})^{-1}{\bm\Lambda}_n+\mbf P^{-1} )^{-1} \bm\Lambda_n^\prime(\bm\Sigma_n^\xi)^{-1}-
({\bm\Lambda}_n^{\prime}({\bm\Sigma}_n^{\xi})^{-1}{\bm\Lambda}_n)^{-1} \bm\Lambda_n^\prime(\bm\Sigma_n^\xi)^{-1}\Vert\, 
\Vert \bm H\Vert\, \Vert\bm\Lambda_n^\prime(\bm\Sigma_n^\xi)^{-1}\Vert\, \Vert\bm\Gamma_n^\xi\Vert\nn\\
&= O(n^{-3/2}),\label{eq:piccolo}
}
because of Lemmas \ref{lem:Gxi}(v), \ref{lem:LSL2}(vii), and \ref{lem:LSL}(iii), and \eqref{eq:gattopardo4bis}. Similarly,
\al{
n \Vert\bm{\mathcal B}\Vert &\le n \Vert  ({\bm\Lambda}_n^{\prime}({\bm\Sigma}_n^{\xi})^{-1}{\bm\Lambda}_n+\mbf P^{-1} )^{-1} \bm\Lambda_n^\prime(\bm\Sigma_n^\xi)^{-1}-
({\bm\Lambda}_n^{\prime}({\bm\Sigma}_n^{\xi})^{-1}{\bm\Lambda}_n)^{-1} \bm\Lambda_n^\prime(\bm\Sigma_n^\xi)^{-1}\Vert^2\Vert\bm\Gamma_n^\xi\Vert\nn\\
&= O(n^{-3}),\label{eq:piccolo2}
}
because of Lemmas \ref{lem:Gxi}(v) and \ref{lem:LSL}(iii). Furthermore, 
\al{
n \Vert  \bm H \bm\Lambda_n^\prime(\bm\Sigma_n^\xi)^{-1} \bm\Gamma_n^\xi(\bm\Sigma_n^\xi)^{-1}\bm\Lambda_n\bm H\Vert &\le 
n \Vert  \bm H \Vert^2 \Vert \bm\Lambda_n^\prime(\bm\Sigma_n^\xi)^{-1}\Vert^2 \Vert \bm\Gamma_n^\xi\Vert =O(1),\label{eq:piccolo3}
}
because of Lemmas \ref{lem:Gxi}(v),  \ref{lem:LSL2}(vii),  and \eqref{eq:gattopardo4bis}.

Thus, by using \eqref{eq:piccolo}, \eqref{eq:piccolo2}, and \eqref{eq:piccolo3}, from \eqref{eq:AABedge}, we have
\al{
n \Vert & \mbf P \bm\Lambda_n^\prime ( \bm\Lambda_n\mbf P\bm\Lambda_n^\prime+\bm\Sigma_n^\xi)^{-1} \bm\Gamma_n^\xi ( \bm\Lambda_n\mbf P\bm\Lambda_n^\prime+\bm\Sigma_n^\xi)^{-1}\bm\Lambda_n \mbf P-  \bm H \bm\Lambda_n^\prime(\bm\Sigma_n^\xi)^{-1} \bm\Gamma_n^\xi(\bm\Sigma_n^\xi)^{-1}\bm\Lambda_n\bm H
 \Vert\nn\\
 & \le 2 n\Vert\bm{\mathcal A}\Vert + n \Vert\bm{\mathcal B}\Vert = O(n^{-3/2}).\label{eq:calamaro}
}
And, by using \eqref{eq:calamaro} and \eqref{eq:calamaretto} into \eqref{eq:gattopardo3} we have
\al{
n\Vert \bm{\mathcal P}_t -   \bm H \bm\Lambda_n^\prime(\bm\Sigma_n^\xi)^{-1} \bm\Gamma_n^\xi(\bm\Sigma_n^\xi)^{-1}\bm\Lambda_n\bm H\Vert\le&\,
n\Vert
\mbf P \bm\Lambda_n^\prime ( \bm\Lambda_n\mbf P\bm\Lambda_n^\prime+\bm\Sigma_n^\xi)^{-1} \bm\Gamma_n^\xi ( \bm\Lambda_n\mbf P\bm\Lambda_n^\prime+\bm\Sigma_n^\xi)^{-1}\bm\Lambda_n \mbf P
-
 \bm H \bm\Lambda_n^\prime(\bm\Sigma_n^\xi)^{-1} \bm\Gamma_n^\xi(\bm\Sigma_n^\xi)^{-1}\bm\Lambda_n\bm H
 \Vert
\nn\\
&+n \Vert \bm H\Vert^2\Vert(\mbf P)^{-1}\Vert^2 \Vert\bm\Pi\Vert+2n \Vert\bm H\Vert \, \Vert(\mbf P)^{-1}\Vert \, \Vert \bm\Pi\Vert\, \Vert\mbf C\Vert
+n \Vert\mbf C\Vert^2\Vert\bm\Pi\Vert\nn\\
=&\, O(n^{-1}).\label{eq:statale}
}
Finally, notice that, by definition:
\al{
\lim_{n\to\infty} n \bm H \bm\Lambda_n^\prime(\bm\Sigma_n^\xi)^{-1} \bm\Gamma_n^\xi(\bm\Sigma_n^\xi)^{-1}\bm\Lambda_n \bm H  & =
\lim_{n\to\infty}  n \bm H  n^{-1}\bm\Lambda_n^\prime(\bm\Sigma_n^\xi)^{-1} \bm\Gamma_n^\xi(\bm\Sigma_n^\xi)^{-1}\bm\Lambda_n n \bm H\nn\\
&  =  (\bm\Sigma_{\Lambda\Sigma\Lambda})^{-1} 
\l\{\lim_{n\to\infty} n^{-1}\bm\Lambda_n^\prime(\bm\Sigma_n^\xi)^{-1} \bm\Gamma_n^\xi(\bm\Sigma_n^\xi)^{-1}\bm\Lambda_n\r\}
(\bm\Sigma_{\Lambda\Sigma\Lambda})^{-1}=\bm{\mathcal W}_t.\label{eq:HLPW}
%
}
Therefore, from \eqref{eq:statale} and \eqref{eq:HLPW}
\al{
\Vert n \bm{\mathcal P}_t -  \bm{\mathcal W}_t \Vert &\le  \Vert n \bm{\mathcal P}_t -  n \bm H \bm\Lambda_n^\prime(\bm\Sigma_n^\xi)^{-1} \bm\Gamma_n^\xi(\bm\Sigma_n^\xi)^{-1}\bm\Lambda_n\bm H\Vert 
+  \Vert n \bm H \bm\Lambda_n^\prime(\bm\Sigma_n^\xi)^{-1} \bm\Gamma_n^\xi(\bm\Sigma_n^\xi)^{-1}\bm\Lambda_n\bm H - \bm{\mathcal W}_t \Vert=O(n^{-1})+o(1),\nn
}
and since $\bm{\mathcal P}_t$ and $\bm{\mathcal W}_t$ do  not depend on $t$ we complete the proof. $\Box$

\end{appendix}

\end{document}